\documentclass[12pt,a4paper, twoside]{amsbook}

\usepackage[utf8]{inputenc}
\usepackage[T1]{fontenc}

\usepackage{amscd}
\usepackage{amsmath,amssymb,amsthm,mathrsfs}
\usepackage[margin=2cm]{geometry}

\usepackage{fancyhdr}
\fancypagestyle{plain}{%
	\fancyhf{}
	\fancyhead[LO,RE]{\thepage}
}

\usepackage{thmtools}

\usepackage{pgf,tikz} 
\usetikzlibrary{matrix,calc,positioning,intersections,through,angles,patterns}
\usetikzlibrary{arrows.meta}

\usetikzlibrary{external}
\tikzset{wqgvert/.style={draw,black,rectangle,fill=white,minimum size=4.5pt,inner sep=0pt}  } 
\tikzset{bqgvert/.style={draw,black,rectangle,fill=black,minimum size=4.5pt,inner sep=0pt}  } 
\tikzset{qvert/.style={draw,black,circle,fill=gray,minimum size=5pt,inner sep=0pt}  } 
\tikzset{bvert/.style={draw,circle,fill=black,minimum size=5pt,inner sep=0pt}  } 
\tikzset{wvert/.style={draw,circle,fill=white,minimum size=5pt,inner sep=0pt}  } 
\tikzset{nvert/.style={}  } 
\tikzset{hvert/.style={draw,circle,minimum size=10pt,inner sep=0pt,line width=0.8pt,gray!75}  } 
\tikzset{hqvert/.style={draw,black,circle,pattern color=gray,pattern=north west lines,minimum size=7pt,inner sep=0pt}  } 
\tikzset{wfoc/.style={draw,rectangle,fill=white,minimum size=4pt,inner sep=0pt}  } 
\tikzset{bfoc/.style={draw,rectangle,fill=black,minimum size=4pt,inner sep=0pt}  } 
\tikzset{gfoc/.style={draw,rectangle,fill=gray,minimum size=4pt,inner sep=0pt}  } 

\usepackage[hidelinks,pdftex,pdfauthor={Niklas Christoph Affolter},pdftitle={Discrete Differential Geometry and Cluster Algebras via TCD maps},pdfkeywords={PhD thesis, discrete differential geometry, cluster algebras, TCD maps, discrete integrable systems, exactly solvable models}, pdflang={english}, bookmarksopen={true}]{hyperref} 
\usepackage{caption} 

\usepackage{parskip}
\usepackage{enumerate}
\usepackage{enumitem}
\newcounter{descriptcount}

\usepackage{calc}

\usepackage{ifthen}

\usepackage{bm} 

\usepackage{cancel} 

\newcommand{\eq}{\Leftrightarrow}

\newcommand*\circled[1]{\raisebox{.5pt}{\textcircled{\raisebox{-.9pt} {#1}}}} 

\newcommand{\qedabovehere}{\par \vspace{-2\baselineskip}\qedhere} 

\DeclareMathSymbol{\shortminus}{\mathbin}{AMSa}{"39}

\DeclareMathOperator{\cro}{cr} 
\DeclareMathOperator{\spa}{span} 
\DeclareMathOperator{\mr}{mr} 

\renewcommand{\Im}{\mathrm{Im}}

\newcommand{\proj}[2]{\pi_{#1 \rightarrow #2}} 

\newcommand{\tcd}{\mathcal T} 
\newcommand{\qg}{{\mathfrak Q}} 
\newcommand{\qui}{{\mathcal Q}} 
\newcommand{\lrp}{{\mathcal{L}\kern-0.15em\mathcal{R}}} 
\newcommand{\pb}{{\mathcal G}} 
\newcommand{\pbm}{\mathcal G^{\bm\shortminus}} 
\newcommand{\lio}{\mathcal O} 
\newcommand{\liom}{\mathcal O^{\bm\shortminus}} 
\newcommand{\lip}{\mathcal P} 
\newcommand{\lipm}{\mathcal P^{\bm\shortminus}} 
\newcommand{\lima}{M} 
\newcommand{\limam}{M^{\bm\shortminus}} 
\newcommand{\fl}{{\mathbf{F}}} 
\newcommand{\lat}{{\mathcal{L}}} 
\newcommand{\fp}{\digamma} 
\newcommand{\enm}[2]{{\mathcal S_{#2}^{#1}}} 

\newcommand{\tgra}{{G_{\mathrm T}}} 
\newcommand{\tgrad}{{G_{\mathrm D}}} 
\newcommand{\tg}{\mathfrak T} 
\newcommand{\tcdm}{{\mathcal T_{\bm{\circlearrowright}}}} 
\newcommand{\tcdp}{{\mathcal T_{\bm{\circlearrowleft}}}} 
\newcommand{\vrc}{R} 
\newcommand{\alt}{\mathscr A} 
\newcommand{\altg}{\mathscr G} 
\newcommand{\plu}{\mathscr{P}} 

\newcommand{\tcdarc}[3]{  
	\draw[-]
		(#1) edge[mid arrow] ++(0,#3*0.6)
		(#1) -- ++(0,#3*0.7) to[out=90,in=180] ++(#3*0.3,#3*0.3) -- ($(#2)+(-#3*0.3,#3)$) to[out=0,in=90] ++(#3*0.3,-#3*0.3) -- (#2)
		(#2) edge[mid rarrow] ++(0,#3*0.6)
	;
}
\newcommand{\tcdarcr}[3]{  
	\draw[-]
		(#1) edge[mid rarrow] ++(0,#3*0.6)
		(#1) -- ++(0,#3*0.7) to[out=90,in=180] ++(#3*0.3,#3*0.3) -- ($(#2)+(-#3*0.3,#3)$) to[out=0,in=90] ++(#3*0.3,-#3*0.3) -- (#2)
		(#2) edge[mid arrow] ++(0,#3*0.6)
	;
}
\newcommand{\tcdcross}[5]{  
	\draw[-]
		(#5) edge[mid rarrow, out=210,in=90] (#1) edge[mid arrow, out=330,in=90] (#2) edge[out=30,in=270] (#3) edge[out=150,in=270] (#4)
	;
}
\newcommand{\tcdcrossr}[5]{  
	\draw[-]
		(#5) edge[mid arrow, out=210,in=90] (#1) edge[mid rarrow, out=330,in=90] (#2) edge[out=30,in=270] (#3) edge[out=150,in=270] (#4)
	;
}

\DeclareMathOperator{\pro }{\mathbf{Pro}} 
\DeclareMathOperator{\aff }{\mathbf{Aff}} 

\newcommand{\pent}{{\mathbf T}} 

\newcommand{\del}{\partial}

\newcommand{\an}{{\mathcal A}} 
\newcommand{\ban}{{\bar{ \mathcal A}}} 

\newcommand{\qn}{\mathbf{QN}}
\newcommand{\dm}{\mathbf{DM}}
\newcommand{\lc}{\mathbf{LC}}
\newcommand{\cub}{\mathrm{Cub}}
\newcommand{\hex}{\mathrm{Hex}}

\newcommand{\N}{\mathbb{N}} 
\newcommand{\Z}{\mathbb{Z}} 
\newcommand{\R}{\mathbb{R}} 
\newcommand{\RP}{\mathbb{R}\mathrm{P}} 
\newcommand{\C}{\mathbb{C}} 
\newcommand{\Cx}{\mathbb C_{\times}}

\newcommand{\p}[1]{\bm{#1}}

\renewcommand{\S}{\mathbb{S}}

\newcommand{\CP}{\mathbb{C}\mathsf{P}} 

\newcommand{\f}[2]{\frac{#1}{#2}} 

\newcommand{\ve}[1]{\begin{pmatrix}#1\end{pmatrix}}

\newcommand{\menge}[2][]{ 
	  \if\relax\detokenize{#1}\relax
		  \left\{#2\right\}
	  \else
		  \left\{#1\ | \ #2\right\}
	  \fi}

\def\centerarc[#1](#2)(#3:#4:#5)
{ \draw[#1] ($(#2)+({#5*cos(#3)},{#5*sin(#3)})$) arc (#3:#4:#5); }

\usetikzlibrary{decorations.pathreplacing,decorations.markings}
\tikzset{ 
	mid arrow/.style={postaction={decorate,decoration={
				markings,
				mark=at position  0.5*\pgfdecoratedpathlength+2.5pt with {\arrow[line width=0.4 mm]{latex}}
	}}},
	mid rarrow/.style={postaction={decorate,decoration={
				markings,
				mark=at position 0.5*\pgfdecoratedpathlength+2.5pt with {\arrow[line width=0.4 mm]{latex reversed}}
	}}},
}
\tikzset{ 
	orient/.style={thick,postaction={decorate,decoration={
				markings,
				mark=at position 0.55*\pgfdecoratedpathlength+2pt with {\arrow{Triangle}}
	}}},
	rorient/.style={thick,postaction={decorate,decoration={
				markings,
				mark=at position 0.45*\pgfdecoratedpathlength+1pt with {\arrow{Triangle[reversed]}}
	}}},
}
\tikzset{ 
	orientb/.style={thick,postaction={decorate,decoration={
				markings,
				mark=at position 0.5*\pgfdecoratedpathlength+1pt with {\arrow{Classical TikZ Rightarrow}},
				mark=at position 0.5*\pgfdecoratedpathlength-1pt with {\arrow{Classical TikZ Rightarrow}},				
	}}},
	rorientb/.style={thick,postaction={decorate,decoration={
				markings,
				mark=at position 0.5*\pgfdecoratedpathlength+1pt with {\arrow{Classical TikZ Rightarrow[reversed]}},
				mark=at position 0.5*\pgfdecoratedpathlength-1pt with {\arrow{Classical TikZ Rightarrow[reversed]}}				
	}}},
}

\tikzset{circle through 3 points/.style n args={3}{%
		insert path={let    \p1=($(#1)!0.5!(#2)$),
			\p2=($(#1)!0.5!(#3)$),
			\p3=($(#1)!0.5!(#2)!1!-90:(#2)$),
			\p4=($(#1)!0.5!(#3)!1!90:(#3)$),
			\p5=(intersection of \p1--\p3 and \p2--\p4)
			in },
		at={(\p5)},
		circle through= {(#1)}
}}

\newcommand{\tcoo}[3]{\coordinate ({v#1#2#3}) at ($#1*(e1)+#2*(e2)+#3*(e3)$);}

\newcommand{\tztriangle}[3]{
\begin{tikzpicture}[baseline={([yshift=-.7ex]current bounding box.center)},scale=.8] 
\tikzstyle{bvert}=[draw,circle,fill=black,minimum size=2.5pt,inner sep=0pt]

\node[bvert] (v1)  at (0,0) {};
\node[bvert] (v2) at (1,0) {};
\node[bvert] (v3) at (.5,0.875) {};

\ifnum#1=1 \draw[-] (v2) -- (v3);\else \draw[-,gray!50] (v2) -- (v3); \fi
\ifnum#2=1 \draw[-] (v1) -- (v3);\else \draw[-,gray!50] (v1) -- (v3); \fi
\ifnum#3=1 \draw[-] (v1) -- (v2);\else \draw[-,gray!50] (v1) -- (v2); \fi
\end{tikzpicture}}

\newcommand{\tzstar}[3]{
\begin{tikzpicture}[baseline={([yshift=-.7ex]current bounding box.center)},scale=.8] 
\tikzstyle{bvert}=[draw,circle,fill=black,minimum size=2.5pt,inner sep=0pt]

\node[bvert] (v1)  at (0,0) {};
\node[bvert] (v2) at (1,0) {};
\node[bvert] (v3) at (.5,0.875) {};
\node[bvert] (v4) at (.5,0.35) {};

\ifnum#1=1 \draw[-] (v4) -- (v1); \else \draw[-,gray!50] (v4) -- (v1); \fi
\ifnum#2=1 \draw[-] (v4) -- (v2); \else \draw[-,gray!50] (v4) -- (v2); \fi
\ifnum#3=1 \draw[-] (v4) -- (v3); \else \draw[-,gray!50] (v4) -- (v3); \fi

\end{tikzpicture}}

\newcommand{\tzspider}[3][0]{
\begin{tikzpicture}[baseline={([yshift=-.7ex]current bounding box.center)},scale=1,rotate=#1] 
\tikzstyle{bvert}=[draw,circle,fill=black,minimum size=2.5pt,inner sep=0pt]

\node[bvert] (v1)  at (0,0) {};
\node[bvert] (v2) at (1,0) {};
\node[bvert] (v3) at (1,1) {};
\node[bvert] (v4) at (0,1) {};
\node[bvert] (w1) at (.25,.25) {};
\node[bvert] (w3) at (.75,.75) {};

\ifnum#2=1 \draw[-] (w1) -- (v1);\else \draw[-,gray!50] (w1) -- (v1); \fi
\ifnum#2=2 \draw[-] (w1) -- (v2);\else \draw[-,gray!50] (w1) -- (v2);  \fi
\ifnum#2=3 \draw[-] (w1) -- (v4);\else \draw[-,gray!50] (w1) -- (v4);  \fi
\ifnum#3=1 \draw[-] (w3) -- (v3);\else \draw[-,gray!50] (w3) -- (v3);  \fi
\ifnum#3=2 \draw[-] (w3) -- (v4);\else \draw[-,gray!50] (w3) -- (v4);  \fi
\ifnum#3=3 \draw[-] (w3) -- (v2);\else \draw[-,gray!50] (w3) -- (v2);  \fi
\end{tikzpicture}} 

\numberwithin{equation}{chapter}

\declaretheorem[style=definition,name=Theorem,qed={\tiny$\blacksquare$},numberwithin=chapter]{theorem}
\declaretheorem[style=definition,name=Corollary,sibling=theorem,qed={\tiny$\blacksquare$}]{corollary}
\declaretheorem[style=definition,name=Definition,sibling=theorem,qed={\tiny$\blacksquare$}]{definition}
\declaretheorem[style=definition,name=Lemma,sibling=theorem,qed={\tiny$\blacksquare$}]{lemma}

\declaretheorem[style=definition,name=Example,sibling=theorem,qed={\tiny$\blacksquare$}]{example}
\declaretheorem[style=definition,name=Remark,sibling=theorem,qed={\tiny$\blacksquare$}]{remark}
\declaretheorem[style=definition,name=Question,sibling=theorem,qed={\tiny$\blacksquare$}]{question}
\declaretheorem[style=definition,name=Conjecture,sibling=theorem,qed={\tiny$\blacksquare$}]{conjecture}

\renewcommand{\emptyset}{\mbox{\O}}

\begin{document}

\frontmatter

\title{Discrete Differential Geometry and Cluster Algebras via TCD maps}

\author{Niklas Christoph Affolter}
\address{TU Berlin, Institute of Mathematics, Strasse des 17. Juni~136,
	10623 Berlin, Germany}
\email{affolter@posteo.net}

\begin{titlepage}
	\begin{center}
		\LARGE \bfseries Discrete Differential Geometry \\and Cluster Algebras via TCD maps
	\end{center}
	\vspace{2cm}
	\begin{center}
		\large vorgelegt von\\M. Sc.\\ Niklas Christoph Affolter\\ORCID: 0000-0003-3936-2060
	\end{center}
	\vspace{0.7cm}
	\begin{center}
		\large an der Fakultät II -- Mathematik und Naturwissenschaften\\der Technischen Universität Berlin\\
		zur Erlangung des akademischen Grades
	\end{center}
	\vspace{.7cm}
	\begin{center}
		\large Doktor der Naturwissenschaften\\Dr. rer. nat.
	\end{center}
	\vspace{.7cm}
	\begin{center}
		\large genehmigte Dissertation
	\end{center}
	\vspace{3cm}
	{
		\large \ Promotionsausschuss:\vspace{.5cm}\\ 
		\begin{tabular}{lll}				
			Vorsitz: & Prof. Dr. Stefan Felsner & (TU Berlin)\\
			Gutachter: & Prof. Dr. Boris Springborn & (TU Berlin)\\ 
			Gutachter: & Prof. Dr. Yuri B. Suris & (TU Berlin)\\ 
			Gutachter: & Prof. Dr. Christian Müller & (TU Wien)\\
			Gutachter: & Prof. Dr. Richard Kenyon & (Yale U)\\
		\end{tabular}
		\vspace{.5cm}\\
		\hspace*{0mm} Tag der wissenschaftlichen Aussprache: 27. März 2023
	}
	\vspace{3cm}
	\begin{center}
		\large Berlin 2023
	\end{center}

\end{titlepage}

\label{deg:}
\label{dof:}

\setcounter{tocdepth}{1}
\tableofcontents

\mainmatter

\chapter{Introduction} \label{sec:introduction}

\section{Outline}

The main goal of this thesis is to find and investigate occurrences of cluster algebras in objects of discrete differential geometry (DDG). During the investigation it became apparent that there are multiple ways in which cluster algebras occur for multiple objects of DDG. In order to understand the different cluster algebras systematically, it is practical to introduce the common framework of triple crossing diagram maps (TCD maps). For TCD maps the multiple occurrences of cluster algebras can be defined and related systematically. It turns out the framework is exhaustive, in the sense that it covers a long list of examples and indeed every discrete 3D-integrable system that is defined in geometric terms that we are aware of. In particular, all the known examples of DDG, discrete integrable systems and embeddings associated to exactly solvable models that feature cluster algebras are included.

For this reason, we begin with a brief informal introduction to TCD maps in Section \ref{sec:introtcdmaps}. We then give a summary of our approach and our findings in Section \ref{sec:introselected}. Expanding on that, we list the full results in Section \ref{sec:introfull}. Subsequently, we provide a list of examples included in our framework in Section \ref{sec:examplelist}, both to show the applicability of the framework and to guide the reader who is interested in particular examples. We finish the introduction with open questions and directions of further research in Section \ref{sec:introfuture}.

\section{TCD maps}\label{sec:introtcdmaps}

\begin{figure}
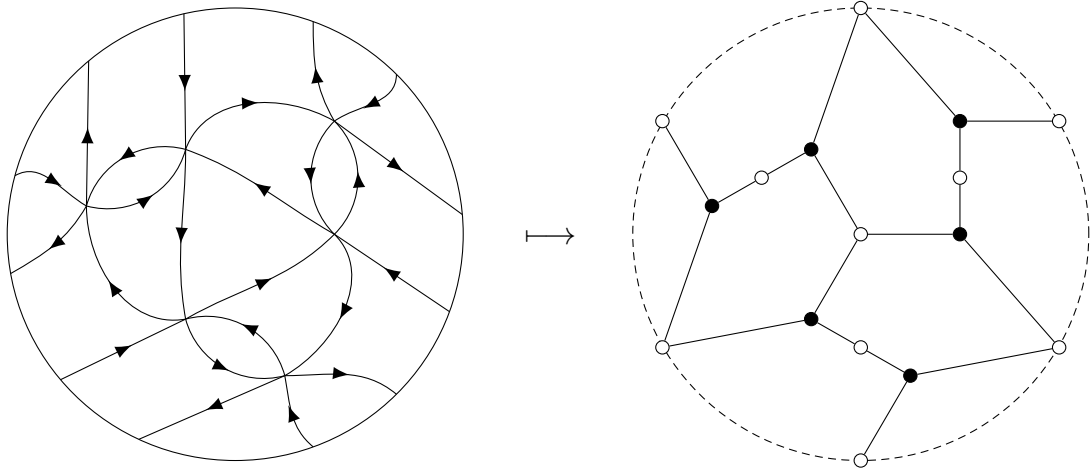

	%
	\caption{An example of a TCD $\tcd$ and the corresponding bipartite graph $\pb$.}
	\label{fig:firsttcdexample}
\end{figure}

This is a short and informal explanation of the most important objects with regard to TCD maps. A formal introduction to TCD maps follows in Section \ref{sec:tcdmap}. A \emph{triple crossing diagram (TCD) $\tcd$} is a set of oriented curves in the disc such that every curve either begins and ends at the boundary or is a closed loop, and such that locally every intersection of curves looks like this:

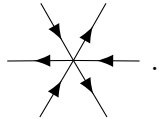
\begin{center}
	\begin{tikzpicture}[scale=1.5, baseline={([yshift=-.7ex]current bounding box.center)}]
		\coordinate (v1) at (90:1);
		\coordinate (v2) at (210:1);
		\coordinate (v3) at (330:1);			
		\coordinate (e12) at ($(v1)!.33!(v2)$);
		\coordinate (e13) at ($(v1)!.33!(v3)$);
		\coordinate (e23) at ($(v2)!.33!(v3)$);
		\coordinate (e21) at ($(v2)!.33!(v1)$);
		\coordinate (e32) at ($(v3)!.33!(v2)$);
		\coordinate (e31) at ($(v3)!.33!(v1)$);
		\coordinate (o) at (0,0);
		
		\draw[-]
			(o) edge[mid rarrow] (e12) edge[mid rarrow] (e23) edge[mid rarrow] (e31)
			(o) edge[mid arrow] (e21) edge[mid arrow] (e32) edge[mid arrow] (e13)
		;
	\end{tikzpicture}
	.
\end{center}
See Figure \ref{fig:firsttcdexample} for an example of a TCD. Therefore the only intersection points are \emph{triple intersection points} and every face of $\tcd$ is consistently oriented either clockwise or counterclockwise. Triple crossing diagrams were introduced by Dylan Thurston \cite{thurstontriple}, we give a formal definition and an introduction in Section \ref{sec:tcds}. 

To each TCD $\tcd$ we also associate a \emph{planar bipartite graph $\pb$}. The white vertices of $\pb$ are the counterclockwise oriented faces of $\tcd$, the black vertices are the intersection points of $\tcd$. There is an edge $(w,b) \in E(\pb)$ if the corresponding counterclockwise face and intersection point are incident in $\tcd$. Therefore at every intersection point, the local configuration looks like this:
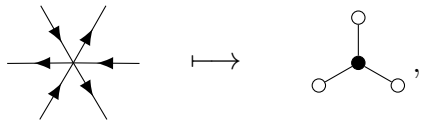
\begin{center}
\hspace{1mm},
\end{center}
where on the left is the intersection point in $\tcd$ and on the right is the corresponding local piece of $\pb$. 

A \emph{TCD map} is a map $T:\tcdp \rightarrow \CP^n$ defined on the counterclockwise faces $\tcdp$ of a TCD $\tcd$, or equivalently on the white vertices of $\pb$, such that at every black vertex of $\pb$ the images of the three incident white vertices are contained in a projective line in complex projective space $\CP^n$. We also make the genericity assumption that no two white vertices adjacent to the same black vertex coincide. 

By considering local changes of combinatorics of TCDs, we introduce \emph{dynamics} to TCD maps. The local moves of TCDs are called \emph{2-2 moves} and come in the two types

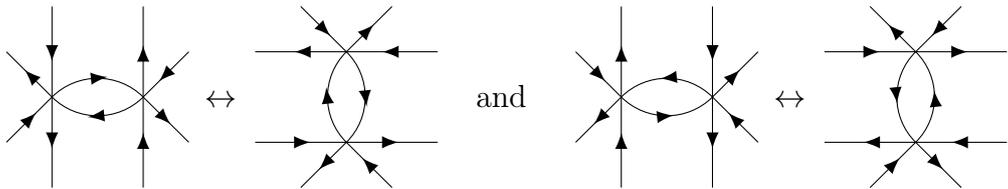
\begin{center}
	%
.
\end{center}
Note that the left 2-2 move is centered at a clockwise face and the right 2-2 move is centered at a counterclockwise face. In terms of the bipartite graph $\pb$ the two moves are

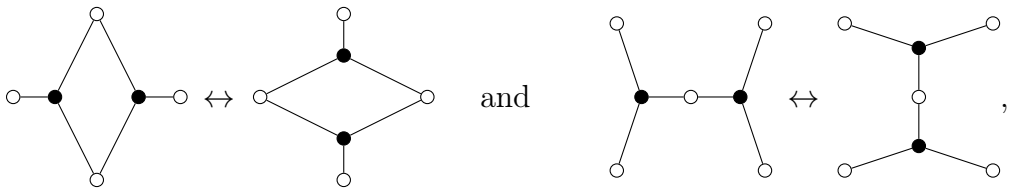
\begin{center}
	\hspace{-1mm}
	%
,

\end{center}

and are called the \emph{spider move} and the \emph{resplit} respectively. The geometric data of a TCD map changes as

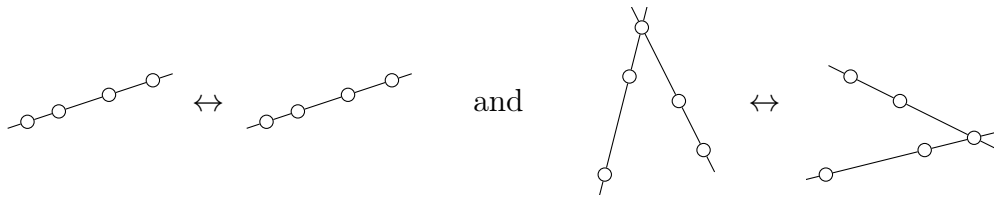
\begin{center}
	\hspace{3mm}
	%
\hspace{1mm},
	
\end{center}
under the spider move and resplit respectively. Thus, we observe that the image of a TCD map does not change under a spider move, because the white vertex set of $\pb$ and the incidence requirements are not changed by the spider move. Contrarily, when performing the resplit the geometric data of a TCD map does change, in the sense that we have to change the image of the white vertex in the center. In the generic case, the new point is defined by the incidence requirements as the intersection of the two lines that correspond to the two adjacent black vertices after the resplit. If the dimension of the projective space is one and incidence geometry is therefore unavailable, the resplit is still well defined via the \emph{dSKP equation}, see Section \ref{sec:tcddskp}.

Other important building blocks of the TCD framework are the \emph{vector-relation configurations}, see Section \ref{sec:vrc}, where we introduce homogeneous lifts of TCD maps to $\C^{n+1}$. The homogeneous lifts allow us to introduce linear relations which determine \emph{edge weights} on $\pb$. These edge weights can be used to define projective cluster variables, see Section \ref{sec:projcluster}. In a specific gauge, the edge weights define the \emph{affine cluster variables}, see Section \ref{sec:affcluster}. Furthermore, the edge weights define complex valued partition functions of \emph{almost perfect matchings}, see Section \ref{def:almostperfect}, which relates TCD maps to statistical mechanics. Another important tool are \emph{labeling induced orientations (li-orientations)}, see Section \ref{sec:sweeps}, which allow us to construct TCD maps via boundary data and the projective cluster variables. Moreover in the generic setup, li-orientations allow us to show that TCD maps are parametrized by their projective cluster variables up to projective transformations, see Section~\ref{sec:fromprojinvariants}.

\section{Summary}\label{sec:introselected}

The initial goal of this research project was to find cluster algebras in DDG examples. More specifically, we were looking for a way to attribute quivers and coefficient type cluster variables to DDG maps such that the dynamics of the DDG maps correspond to mutations of the quiver and cluster variables. Initially, we found several new but different cluster structures for different examples (Miquel dynamics, Q-nets, Darboux maps). The cluster structures differed in the sense that some had variables that are invariant under projective transformations while other cluster structures had variables that are invariant under affine transformations. We then realized there is a unified description for these examples that simultaneously gives rise to both types of cluster structures in every example. The unified description is via the framework of \emph{TCD maps}, as we outlined in Section \ref{sec:introtcdmaps} and as we will formally introduce in Section \ref{sec:tcdmap}. The projective cluster structure is defined in Section \ref{sec:projcluster} and the affine cluster structure in Section \ref{sec:affcluster}. 

Subsequently, we realized that even examples not only from DDG, but also from discrete integrable systems and statistical mechanics fit in the TCD map framework. Indeed, the previously known cluster algebra structures of T-graphs \cite{kenyonsheffield}, pentagram map \cite{glickpentagram}, dSKP lattices \cite{absoctahedron} and ideal hyperbolic triangulations via shear coordinates \cite{pennerteichmuller} can all be recovered as special cases. We introduce the examples during the course of the thesis while we develop the theory of TCD maps. We give a list of all examples covered in Section \ref{sec:examplelist}. We consider the list to be almost exhaustive, in the sense that we are currently not aware of any geometric discretely integrable 3D-system that is not covered. 

We also discovered that TCD maps allow a unified description of other interesting properties besides cluster algebras. For example, one can show that TCD maps are \emph{multi-dimensionally consistent} via the flip-graph and Desargues' theorem, see Section \ref{sec:tcdconsistency}. Thus we obtain a unified proof of multi-dimensional consistency for all examples. We can also show that dynamics of TCD maps correspond to propagation of the dSKP equation in $A_N$ lattices, see Section \ref{sec:tcddskp}. Moreover, we show that the projective cluster variables together with boundary data actually characterize a TCD map up to projective transformation. This is a typical result in many DDG examples, that we managed to translate to the whole TCD map framework. In fact, we can define what the \emph{maximal dimension} of a TCD map is, and in maximal dimension a TCD map is characterized by its projective cluster variables up to projective transformations without any boundary data, see Section \ref{sec:fromprojinvariants}.

Of course, once we are able to capture multiple examples in the TCD framework, it is natural to ask \emph{how} they relate inside the framework. It turns out an important ingredient is to understand how to take a \emph{section} $\sigma_H(T)$ of a TCD map $T$, by intersecting the lines represented in a TCD map with a hyperplane $H$ such that the result is again a TCD map, see Section \ref{sec:sections}. For example, we can then show that the section of a line complex is a Q-net, the section of a Q-net is a Darboux map and the section of a Darboux map is a line compound. But apart from relating the geometry of examples, we also find that sections relate different cluster algebras to each other. Specifically, we can show that the projective cluster structure of the section $\sigma_H(T)$ with the hyperplane at infinity $H$ is the affine cluster structure of the TCD map $T$, see Section \ref{sec:sectioncluster}. Moreover, the combinatorics of taking a section is an operation that relates TCDs with different endpoint matchings, but does not cut or glue the strands themselves. We think this operation is of combinatorial interest in itself even without the geometric and algebraic interpretation. Another interesting projective operation is to take a projective dual of a TCD map, see Section \ref{sec:projduality}. We do not go into detail here, but on the TCD level we define a strand preserving operation that also relates to the TCD with all strand orientations reversed. Moreover, we are identify the cluster structures of primal and dual TCD maps after inverting all cluster variables, see Section \ref{sec:projduality}.

\begin{figure}
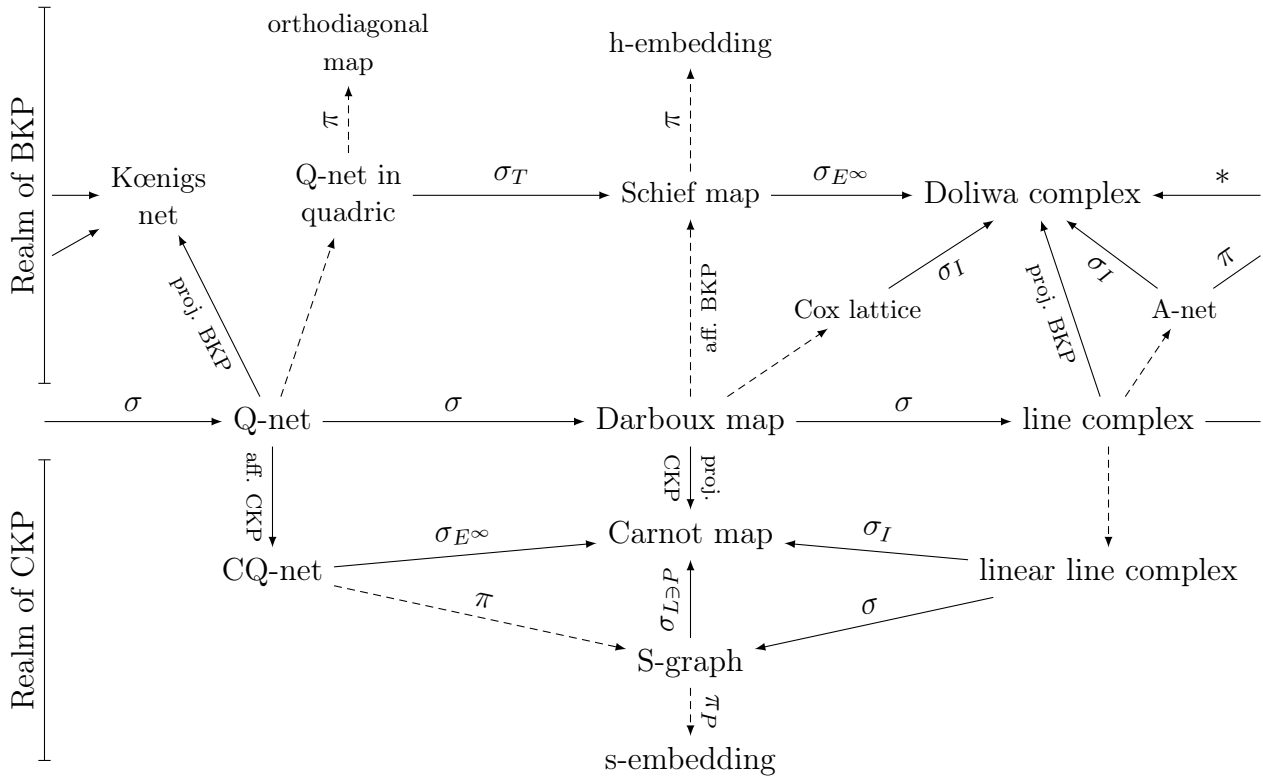


	\caption{Relations between TCD maps for $\Z^3$ combinatorics. A dashed arrow from A to B indicates that B is a special case of A. Sections with respect to subspaces $E$ are denoted by $\sigma_E$, projections by $\pi$.}
	\label{fig:bkpckprelations}
\end{figure}

While we mostly focus on the projective cluster variables $X$ and affine cluster variables $Y$ of TCD maps which are both of coefficient type, we also find the other type of cluster variables for TCD maps which we call $\tau$-variables, see Section \ref{sec:tautcd}. Moreover, we attribute global invariants to TCD maps that are combinations of dimer or almost perfect matching partition functions $Z$ with complex weights, see Section \ref{sec:dimerinvariants}. There are two well known subvarities for the face weights of the dimer model, the resistor subvariety and the Ising subvariety that give rise to the spanning tree model and the Ising model as reductions of the dimer model. Therefore, the ability to associate the dimer model to TCD maps raises the question of how these two subvarieties occur in TCD maps. We are able to show that the cases in which the dimer model reduces to the spanning tree or Ising model correspond to known reductions of TCD maps. In fact, we show that every known example that we are aware of that features the BKP equation (Equation \eqref{eq:bkp}) has, as a TCD map, some of its cluster variables in the resistor subvariety. Analogously, every example that features the CKP equation (Equation \eqref{eq:ckp}) has some of its cluster variables in the Ising subvariety. The viewpoint via the subvarities has a significant advantage. Previously it was only possible to state whether a map in \emph{some} way features the BKP (respectively CKP) equation, we can now \emph{distinguish} whether a TCD map is projectively BKP (CKP) or affinely BKP (CKP) with respect to some subspace. We can also \emph{define} what the set of all TCD maps is that is projectively BKP (CKP) or affinely BKP (CKP). For example, we can show that K{\oe}nigs nets \cite{bsmoutard,doliwatnets} are \emph{exactly} those Q-nets that are projectively BKP, see Section \ref{sec:konigs}. On the other hand, Schief maps \cite{schieflattice} are only a subset of those Darboux maps that are affinely BKP, see Section \ref{sec:schiefmaps}. There is a considerable number, one could say a zoo, of examples that feature the BKP or the CKP equation in the literature, but their interrelations have rarely been considered. The ability to recognize the BKP and CKP equation in a canonical manner together with an improved understanding of sections and projective dual allow us to find order in the zoo, as we illustrate in Figure \ref{fig:bkpckprelations} and work out in detail in Chapter \ref{cha:subvar} and Chapter \ref{cha:bilinear}.

We give a full list of results in Section \ref{sec:introfull}, but let us mention a few more results that may be of particular interest. First, in Section \ref{sec:cptemb} we illustrate how to associate two different positive cluster structures to Miquel dynamics on circle patterns. The affine cluster structure was found already \cite{amiquel}, and independently discovered by Kenyon, Lam, Ramassamy and Russkikh \cite{kenyonlam}, whereas the projective cluster structure is new. Secondly, we relate the BMS variables \cite{bmscircular} that were introduced for circular Q-nets to the affine cluster structure of circular Q-nets viewed as TCD maps, and this relation also relates the cluster Poisson algebra and the quantization, see Section \ref{sec:circularnets}. Thirdly, we can also explain the special case of projective flag configurations arising in Fock and Goncharovs higher Teichmüller theory \cite{fghighertm} via TCD maps geometrically in a very  direct way, and retrieve their cluster variables as special case of the (dual) projective cluster variables.

Finally, let us mention two ``soft results''. First of all, we think one big advantage of the TCD map framework is that it connects results and questions from different math communities. We especially think of the DDG community and the already closely connected discrete integrable systems community, the dynamical systems community around the pentagram map that is intertwined with the totally positive Grassmannian community, as well as the exactly solvable models community. These relations are not necessarily surprising, but we think we have made some progress in making them more concrete and universal. We also think or hope that the TCD map framework could enable a lot more result being transferred between the communities in the future. 

The second soft result we want to mention is that the high amount of structure in TCD maps, in our experience, makes many results easy to guess. Especially with respect to finding BKP and CKP structures, the first reasonable guess surprisingly often turns out to be true. This is because due to the TCD map framework, we have a clear understanding of the combinatorics and geometric invariance. For example, due to work of Bobenko and Schief on linear line complexes \cite{bobenkoschieflinecomplexes}, we knew that linear line complexes are related to the CKP equation. But looking at the combinatorics of line complexes in $\RP^3$, the only way that the CKP equation can canonically appear is via sections with a line. On the other hand, the distinguished object in a linear line complex is a null-polarity and the only distinguished lines in a null-polarity are isotropic lines. Thus the only reasonable guess is that the section of a linear line complex with an isotropic line is CKP. And indeed, this turned out to be true, see Section \ref{sec:linearlinecomplexes}. Similar reasoning can be applied in many other examples. 

\section{List of results}\label{sec:introfull}

\begin{description}[font=\bfseries R\stepcounter{descriptcount}\thedescriptcount)~,align=left]
	\item[TCD maps] The notion of a \emph{triple crossing diagram} was introduced by Dylan Thurston \cite{thurstontriple}. We give a short review of the definitions and known results in Section \ref{sec:tcds}. We then define the fundamental object of this thesis, the so called \emph{TCD maps} in Section \ref{sec:tcdmap}. We also already gave an informal introduction in Section \ref{sec:introtcdmaps}. We discuss the two basic local moves, the \emph{resplit} and the \emph{spider move} in terms of combinatorics and geometry. As a first lemma, we can relate the resplit to the \emph{dSKP equation} (see \cite{ksclifford}).
	\item[Multi-dimensional consistency] First we show that TCD maps naturally define maps from the $A_n$ lattice to projective spaces via flips. Then we rephrase the notion of multi-dimensional consistency on $A_n$ known from discrete integrable systems \cite{absoctahedron} into consistency along cycles in the flip graph of TCDs. A recent result \cite{bwtriple} shows that there are three types of cycles that generate all cycles of the flip graph. This allows us to prove multi-dimensional consistency of TCD maps in general (Section \ref{sec:tcdconsistency}). The consistency along the non-trivial cycles in terms of geometry are proven using Desargues' theorem. Multi-dimensional consistency has been already known for most examples we consider (see \cite{ddgbook} for a guide), although in every example specific arguments were used. Therefore it is an added benefit of the TCD map framework that we now have one universal proof. Moreover, for t-embeddings and s-embeddings we believe the multi-dimensional consistency to be a new result.
	\item[Cluster structures]Finding cluster structures for Q-nets and other objects of DDG was one of the main motivations of this thesis. Indeed, we find two cluster structures associated to each TCD map. The \emph{projective cluster structure} (Section \ref{sec:projcluster}) is invariant under projective transformations, the \emph{affine cluster structure} (Section \ref{sec:affcluster}) is invariant under affine transformations. Both structures are new for almost all examples considered. The exceptions are that the projective cluster structure specializes for the pentagram map (Section \ref{sec:pentagram}) to findings of Glick \cite{glickpentagram}, and in the special case of dSKP lattices to the Y-system found by Adler, Bobenko and Suris \cite{absoctahedron}. Moreover, the affine cluster structure is new for t-embeddings (Section \ref{sec:cptemb}), but was also found independently \cite{kenyonlam, amiquel}. The affine cluster structure also specializes to findings of Kenyon and Sheffield  \cite{kenyonsheffield} for T-graphs (Section \ref{sec:tgraphcluster}), and specializes to recent results for h- and s-embeddings (Sections \ref{sec:harmonicemb}, \ref{sec:sembeddings}) by Kenyon, Lam, Ramassamy and Russkikh \cite{kenyonlam}. 
	 
	\item[Projections]We show that TCD maps are compatible with projections, in the sense that for a generic central projection dynamics remain well-defined. We explain how to take advantage of this result to define cube-flips in Q-nets, Darboux maps and line complexes in all dimensions including dimension one (Section \ref{sec:degendm}). This result is new except for Darboux maps \cite{schieflattice}, and in the three-dimensional case reproduces the idea of \emph{fundamental line complexes} investigated by Bobenko and Schief \cite{bobenkoschieflinecomplexes}. In all these cases we explain the occurrence of multi-ratio equations via certain TCD strands. Indeed, the combinatorics of the equations are derived from zig-zag paths while the algebra is derived from a lemma related to Menelaus theorem. We also obtain the new results that (Section \ref{sec:harmonicemb}) h-embeddings are projections of Darboux maps, that (Section \ref{sec:sembeddings}) s-embeddings \cite{chelkaksembeddings} are projections of Q-nets and that (Section \ref{sec:cptemb}) Miquel dynamics \cite{ramassamymiquel} are a projection of Laplace-Darboux dynamics \cite{doliwalaplace}.
	
	\item[Sections]We show how to construct a section $\sigma_H(T)$ of a TCD map with respect to a generic hyperplane $H$ (Section \ref{sec:sections}). We do this construction a priori on the level of the graph $\pb$, to obtain the geometry we require. Surprisingly, on the level of the triple crossing diagram this operation is local and preserves the strands. We suspect this operation is of combinatorial interest in itself, as it yields a way to compare TCDs with different endpoint matchings. We also show that this operation preserves minimality. On the algebraic level, we show that in fact the affine cluster structure of $T$ coincides with the projective cluster structure of $\sigma_H(T)$ (\ref{sec:sectioncluster}). Because we can also iterate sections, this shows that we can actually associate a whole sequence of cluster structures to a TCD map. Regarding examples of TCD maps associated to $\Z^3$, it was already known that the section of a Q-net is a Darboux map \cite{ddgbook}. We extend this observation and show that the section of a Darboux map is a line complex, and the section of a line complex is a Q-net. On the level of quivers for stepped surfaces, this leads to the cyclic sequence of sections hexahedral $\mapsto$ cuboctahedral $\mapsto$ cuboctahedral $\mapsto$ hexahedral again.
	
	\item[Sweeps]Recall that we consider TCDs to be closely related to the $A_n$ lattice. The analogue of TCDs for the $\Z^n$ lattice are \emph{pseudoline arrangements} \cite{felsnerbook}. In pseudoline arrangments there is the notion of \emph{sweeps}. In statistical mechanics this technique is informally called ``pulling through the strand'' \cite{baxterplanar}. In discrete integrable systems this corresponds to the propagation of a new solution in an additional dimension, also sometimes called a discrete Bäcklund transform. We study the TCDs for which sweeps are possible to perform with a single strand (Section \ref{sec:sweeps}), and call these TCDs \emph{sweepable TCDs}. In these cases we show that depending on the projective dimensions, sweeps produce either sections of TCD maps or \emph{extensions} of TCD maps (Section \ref{sec:tcdextensions}). As tools we introduce \emph{labeling induced orientations}, very practical acyclic orientations that are very reminiscent of \emph{perfect orientations} that occur in the study of totally positive Grassmannians \cite{postgrass}, but are also analogues of \emph{extremal matchings} as studied in dimer theory \cite{broomheaddimer}. Furthermore, we show that TCD maps defined on minimal TCDs can always be incrementally supplemented by marked points at the boundary so that the TCD map is defined on a sweepable TCD. Each addition of a marked point corresponds on the combinatorial level to steps in the \emph{circular Bruhat order} \cite{postgrass,williamsgrass} as well as the the BCFW-bridge decomposition \cite{scattamp}.
	
	\item[Projective duality]We introduce a definition of the projective dual $T^\star$ of a TCD map $T$ in Section \ref{sec:projduality}. To make the definition of the projective dual unique, we need to look at \emph{flags of TCD maps}. Essentially a flag of TCD maps is a sequence of TCD maps $(T_{k})_{1\leq k \leq n}$ related by sections, that is $T_{k-1} = \sigma(T_k)$. Then to each flag of TCD maps there is a unique projective dual flag of TCD maps. The first important ingredient is that we consider the operation of a \emph{line dual} $\eta(\tcd)$ on a TCD $\tcd$, or the accompanying planar bipartite graph $\pb$. While the section $\sigma(\tcd)$ captures the relations of lines of $T$ via common planes, the line dual $\eta(\tcd)$ captures the relations of lines of $T$ via common intersection points. We also introduce the operation $\iota(\tcd)$ that just reverses the orientation of all strands of $\tcd$. Note that $\sigma(\tcd)$ and $\eta(\tcd)$ involve choices, while $\iota(\tcd)$ is defined without choices. The details are contained in Lemma \ref{lem:iotaetasigma}, informally the important part is that there are choices such that the three equations
	\begin{align}
		\iota \circ \sigma \circ \iota \circ \sigma(\tcd) = \tcd,\quad \sigma \circ \sigma(\tcd) = \iota \circ \eta (\tcd), \quad \iota \circ \sigma (\tcd) = \sigma \circ \eta(\tcd)
	\end{align}
	hold. In fact we show which choices in these equations uniquely determine the remaining choices. This allows us to show that there are unique TCDs such that the TCDs $\tcd_k$ of the primal and the TCDs $\tcd_k^\star$ of the dual satisfy
	\begin{align}
		\sigma(\tcd^\star_{k}) = \tcd^\star_{k-1}, \quad \iota(\tcd^\star_k) = \tcd_{n-k},  \quad \eta(\tcd^\star_k) = \tcd_{n-k+2},
	\end{align}
	for all reasonable $k$. On the geometric side of things, we keep track of certain subspaces $U_k: \tcdp_k \rightarrow \CP^n$ when taking the sections that define the primal flag, see Definition \ref{def:subspacemap}. The subspace map of the primal and dual flag satisfy $U^\star_k(w^\star) = (U_{n-k+1}(w))^\perp$ for all reasonable $k$. Moreover, in Section \ref{sec:projclusterduality} we find that we cannot just relate the geometry and combinatorics of primal and dual flags, but also the cluster structures. More explicitly, the cluster structures satisfy
	\begin{align}
		\aff_{E_{k-1}}(T_k) = \rho(\aff_{E_{n-k}^\star}(T_{n-k+1}^\star)) \quad \mbox{and} \quad \pro(T_{k}) = \rho(\pro(T^\star_{n-k-1}))
	\end{align}
	for all reasonable $k$, where $\rho$ reverses all arrows of the quiver and maps every cluster variable $X_v$ to $X_v^{-1}$.

	\item[Focal nets and Laplace transforms]
	The study of discrete focal and discrete Laplace transforms of Q-nets was initiated by Doliwa \cite{doliwalaplace}. We show that focal transforms both on $\Z^2$ (Section \ref{sec:laplacedarboux}) and stepped surfaces (Section \ref{sec:focalnets}) correspond to 2-2 moves on TCD maps. Due to our investigation of sections, we essentially consider Q-nets, Darboux maps and line complexes on equal footing. Thus we propose a new extension of the definitions of focal transforms to not only Q-nets and line complexes but also Darboux maps. This extension very naturally operates on the $\Z^3/(1,1,1)$ lattice of maps. On this lattice the Laplace transform naturally operates on three disjoint $A_2$ sublattices, one for each type of map.
	\item[Uniqueness from invariants]
	In DDG the typical approach is to associate some case-specific quantities to example systems and then show that together with some case-specific boundary data the system is uniquely determined. In the framework of TCD maps, we generalize and strengthen these results. We show that if a TCD map assumes its maximal dimension and is defined on a sweepable TCD, the TCD map is uniquely defined up to projective transformations by its projective cluster variables (Section \ref{sec:fromprojinvariants}). We want to stress that there is no free boundary data in this case. If in addition the projective cluster variables are strictly positive, we can also show existence. If the TCD is not sweepable or not in maximal dimension, we show uniqueness given the projective cluster variables together with boundary data. Indeed, the necessary boundary data corresponds to the minimal elements of the li-orientation. Note that typical DDG examples are sweepable but not in maximal dimension. Thus one can interpret the additional choice of boundary data in these cases as the choice of projection from a unique map in maximal dimension.
	
	\item[Dimer models]In Section \ref{sec:dimers} we give a very short introduction into the dimer model. We explain how both the projective as well as the affine cluster structures are naturally associated to dimer models (with generically non-positive weights) defined on $\pb$ and $\pbm$ respectively. By introducing \emph{almost perfect matchings} we can associate projective invariants to a TCD map that are also invariants under moves (Section \ref{sec:dimerinvariants}). We do not give a full exposition, but this is certainly an interesting direction for further research. In particular, there are clear parallels to ideas employed in the study of totally positive Grassmannians as introduced by Postnikov \cite{postgrass}.
	
	\item[Subvarieties from statistical mechanics]We recall how two other models of statistical mechanics, the spanning tree model \cite{kirchhoff} and the Ising model \cite{lenzising} can be viewed as reductions \cite{temperleybijection, dubedattrick} of the dimer model. Each of the two reductions requires a particular subclass of combinatorics, that is associated to quad-graphs. Additionally, the cluster variables are in certain algebraic subvarieties, namely the resistor subvariety (Section \ref{sec:spanningtrees}) and the Ising subvariety (Section \ref{sec:ising}). This raises the question whether it is possible to identify those TCD maps, that have the right combinatorics and projective cluster variables in the resistor and Ising subvariety respectively. The TCD framework enables us to pose this question in the first place and provides an affirmative answer. We show that there are two possibilities for TCD maps with projective variables in the resistor subvariety. In the Q-net case (Section \ref{sec:konigs}) these are exactly the K{\oe}nigs nets \cite{bsmoutard, doliwatnets}. In the line compound case (\ref{sec:Doliwacomplex}) these are exactly the \emph{Doliwa compounds}. The latter are maps that we have newly defined, as line compounds that are characterized by a local multi-ratio condition. They are related by a focal transform to certain Q-nets, which are also called K{\oe}nigs nets by Doliwa \cite{doliwakoenigs}, but are not to be confused with the K{\oe}nigs nets as mentioned before (or as in Definition \ref{def:koenigs}). We also show that various maps defined in the literature, especially diagonal intersection nets \cite{ddgbook}, Schief maps \cite{schieflattice} and reciprocal figures \cite{ksreciprocal} are related to Doliwa compounds and K{\oe}nigs nets within the framework of TCD maps. For the Ising subvariety the only candidate case due to combinatorics is the case of Darboux maps (Section \ref{sec:darbouxckp}). We show that the TCD maps with projective variables in the Ising subvariety are exactly the \emph{Carnot maps} (also called CKP Darboux maps in the literature \cite{schieflattice}). We show that Carnot maps are in fact sections of what we call \emph{CQ-nets}, which were introduced as C-quadrilateral lattices by Doliwa \cite{doliwacqnet}. Additionally, we show that Carnot maps are also sections of S-graphs \cite{chelkaksgraphs}, and that S-graphs can be obtained as projections of certain CQ-nets.
	
	\item[Q-nets in quadrics and related maps]
	Consider $\CP^3$ with a distinguished symmetric bilinear form $b:\C^4\times \C^4\rightarrow \C$ on the homogeneous coordinates $\C^4$. The (projectivization of the) zeroset of such a bilinear form is called a \emph{quadric}. Q-nets with points in a quadric \cite{doliwaqnetsinquadrics} are an important object in DDG, as one can study subgeometries of projective geometry (for example Möbius geometry) by distinguishing a quadric in projective space \cite{bsorganizing}. Let us call a line $\ell \subset \CP^3$ isotropic if for any two points $p_1,p_2\in \ell$ holds that $b(\hat p_1, \hat p_2) = 0$, where $\hat p_1,\hat p_2\in \C^4$ are arbitrary homogeneous lifts of $p_1,p_2$. The fact that we can capture Q-nets via TCD maps and relate those in turn to dimer cluster structures, together with our investigations into sections allow us to formulate the following novel result (Section \ref{sec:quadricqnet}): The section $\sigma_\ell(q)$ with respect to an isotropic line $\ell$ of a any Q-net $q$ with all points contained in a quadric has projective invariants in the resistor subvariety. Equivalently, $\sigma_\ell(q)$ is a Doliwa compound. In the process of proving that statement, we obtain several other interesting results. For example we find that several maps in $\CP^2$ have affine variables with respect to certain points that are in the resistor subvariety. These maps include Q-nets inscribed in conics, but also the so called \emph{2-conical nets} that include circle patterns. Moreover, we also show that if the quadric in $\CP^3$ is in fact a cone, we recover Schief maps \cite{schieflattice} (which Schief calls BKP maps) as sections $\sigma_H(q)$, where $q$ is a Q-net with points in the cone and $H$ is a tangent plane. 
	
	\item[Linear line complexes and related maps]
	Instead of a symmetric bilinear form we can also consider an anti-symmetric bilinear form $\omega: \C^4\times \C^4\rightarrow \C$. All points are isotropic with respect to an anti-symmetric bilinear form, but not all lines are. A line complex such that all its lines are isotropic with respect to a bilinear form is called a \emph{linear line complex} \cite{bobenkoschieflinecomplexes}. Note that Bobenko and Schief introduced linear line complexes as sections of the Plücker quadric, for the relation to anti-symmetric forms see \cite{semplekneebone}. We provide an another new result: The section $\sigma_\ell(l)$ of a linear line complex $l$ with respect to an isotropic line $\ell$ has projective invariants in the Ising subvariety. Equivalently, $\sigma_\ell(l)$ is a Carnot map \cite{schieflattice}. Similarly to the case of Q-nets with points in quadrics, the path to the proof illuminates new relations between several other maps in the literature. In particular, for any hyperplane $H$ the section $\sigma_H(l)$ is an S-graph \cite{chelkaksgraphs} with respect to the unique polar point $P$ of $H$. This also leads to relations to s-embeddings and CQ-nets.
	
	\item[A-nets and Cox lattices]A-nets were introduced by Sauer \cite{saueranet} and Cox lattices by King and Schief \cite{kscox}. The Plücker lift $\hat a$ of an A-net $a$ was first considered by Doliwa \cite{doliwaanetspluecker}, and he discovered that the Plücker lift of an A-net is an isotropic line complex (Section \ref{sec:anets}). We newly apply the same idea to Cox lattices and show that the Plücker lift $\hat c$ of a Cox lattice $c$ is in fact an isotropic Darboux map (Section \ref{sec:cox}). We also apply our insights on the canonical appearances of BKP structures to both Plücker lifts. In the case of A-nets we find that sections $\sigma_H(\hat a)$ with isotropic planes $H$ are Doliwa compounds, or equivalently are projectively BKP line complexes. We then go on to show that $\sigma_H(\hat a)$ is actually the projective dual of the projection $\pi(a)$ of $a$ to a plane. Because of our results on the cluster variables of projective duals of TCD maps, we see that $\pi(a)$ is actually a K{\oe}nigs net. This reproduces a result noted in the DDG book \cite[Exercise 2.29]{ddgbook}, with the advantage that it is now not a peculiarity but a canonical deduction. In the case of Cox lattices we find not one but two BKP structures associated to sections of $\hat c$. It is unclear if one of the two BKP structures that we discovered coincides with the BKP structure found by King and Schief, and we leave this as an open question. We also apply our understanding of general focal transforms to show that the focal transform of the Plücker lift of a Cox lattice is the Plücker lift of an A-net.
	
	\item[Quadrirational Yang-Baxter maps]Quadrirational Yang-Baxter maps \cite{absyangbaxter} were introduced together with an interpretation in terms of pencils of conics. We observe that quadrirational Yang-Baxter maps are a reduction of Darboux maps (Section \ref{sec:quadrirationalybmaps}). We apply the techniques of TCD maps and sections, to show that certain sections of quadrirational Yang-Baxter maps are Doliwa complexes and therefore feature the BKP equation. We also show quadrirational Yang-Baxter maps themselves arise naturally as sections of what we call \emph{alternating generator Q-nets}. These Q-nets are inscribed in a quadric intersection curve, and along strips of quads the edge-lines of the Q-net are generators of a fixed quadric of the pencil of quadrics that intersects in the quadric intersection curve.
		
	\item[Circular Q-nets and BMS-variables]Circular Q-nets \cite{bobenkocircular,cdscircular} are Q-nets such that every quad is inscribed in a circle. Bazhanov, Mangazeev and Sergeev \cite{bmscircular} found variables (the \emph{BMS variables}) for circular Q-nets, that feature an ultra-local Poisson bracket and quantization. We show in Section \ref{sec:circularnets}, that the affine cluster variables of circular Q-nets can be factorized into the BMS variables. We show this factorization also relates the Poisson brackets of BMS to the canonical Poisson bracket of the affine cluster algebra, and an analogous statement for the quantization. Interestingly, it appears the factorization exists also for non-circular Q-nets, although it is not clear if this factorization is invariant under the cube-flip for general Q-nets.
		
	\item[Embeddings from statistical mechanics]In Chapter \ref{cha:cpone} we study relations between \emph{circle patterns} and statistical mechanics and investigate geometric and algebraic properties of embeddings known from statistical mechanics, in particular \emph{t-embeddings}, \emph{h-embeddings} and \emph{s-embedding}. First, we show that the centers of circle patterns under the \emph{Miquel move} satisfy the dSKP equation. This was shown by the author \cite{amiquel} and independently by Kenyon, Lam, Russkikh and Ramassamy \cite{kenyonlam}. In fact, the cluster structure of t-embeddings coincides with the affine cluster structure of t-embeddings. We also introduce maps to $\CP^1$ that we call \emph{u-embeddings}. Unlike t-embeddings, they do not have real affine cluster variables but real projective cluster variables. We then show how to associate a u-embedding to a circle pattern with $\Z^2$ combinatorics. As a result, we obtain a second cluster structure associated to circle patterns that is projectively invariant (in fact, Möbius invariant). We then show that the t-embedding and the u-embedding associated to a cluster structure are actually related by a sweep. In other words, they can be understood as one map from a subset of $A_4$ such that this combined map satisfies the dSKP equation on $A_4$. We then proceed to study h-embeddings, which are related to harmonic embeddings introduced by Tutte \cite{tutteembedding}. We also study s-embeddings, introduced by Chelkak \cite{chelkaksembeddings}. A cluster structure for both h- and s-embeddings was introduced by Kenyon et al. \cite{kenyonlam}. We show that these cluster structures are the canonical affine cluster structures if we view h-embeddings as Darboux maps and s-embeddings as Q-nets. We also prove that h-embeddings are special cases of Schief maps, while s-embeddings are a special case of a type of map that we newly introduce and call \emph{fixed focal point map}. We also show that the \emph{orthodiagonal maps} that accompany harmonic embeddings can be viewed as Q-nets and that the projective cluster structure of orthodiagonal maps coincides with the affine cluster structure of h-embeddings and vice versa. Moreover, in a lift to higher projective dimension the h-embedding can be understood as the section of the orthodiagonal map, which in turn is inscribed in a quadric. We also find another peculiarity, namely that the projective cluster variables of h-embeddings satisfy equations that are very close (but not the same) to the equations defining the Ising subvariety. In the case of s-embeddings our TCD map framework allows to see s-embeddings as certain projections of \emph{S-graphs}, which in turn relates s-embeddings to CQ-nets and linear line complexes. These observations then lead to the observation that s-embeddings are in fact affine CKP in any chart of $\CP^1$, not just in the affine chart in which the s-embedding is given.

	\item[$\tau$-variables]In most of the thesis we only consider cluster variables that are projectively invariant and of coefficient type in the sense of cluster algebras. In Chapter \ref{cha:tau} we show how to associate the other type of cluster variables, which we call $\tau$-variables to TCD maps. The $\tau$-variables obey the \emph{discrete KP equation}, also known as the \emph{octahedron recurrence}. We restrict ourselves to a generic setup, although we explain informally how to generalize the definition to less generic setups. The $\tau$-variables are defined via determinants of certain hyperplanes associated to the strands of the TCD. These hyperplanes are actually points in the projective dual flag. Because we employ determinants, the $\tau$-variables are only defined up to choice of homogeneous lifts of the hyperplanes in dual space. This is the typical gauge freedom in solutions of the dKP equation.
	
	\item[Projective flag configurations]In Section \ref{sec:fgmoduli} we investigate flags of projective subspaces that are attached to vertices of a triangulation. This is a special case of systems studied by Fock and Goncharov \cite{fghighertm, fgtwoflags, gwebs}. We show that projective flag configurations can be described as TCD maps and that the $X$-variables employed by Fock and Goncharov coincide with the projective $X$-variables of TCD maps. We conjecture that the $\tau$-variables of TCD map also coincide with the $\tau$-variables of Fock and Goncharov. We also show that flips in the triangulation correspond to sequences of 2-2 moves.	
	
\end{description}

\section{List of examples of TCD maps}\label{sec:examplelist}

We give a list of all examples of objects and dynamics that we cast in the TCD map framework. This should help the interested reader to navigate to her or his example of choice. Of course, the list also illustrates the wide applicability of the TCD map framework. We sort the examples into categories depending on the type of specialization of a TCD map they are.

The first group of examples can be cast as TCD maps with particular constraints on the combinatorics.

\setcounter{descriptcount}{0}
\begin{enumerate}[font=\bfseries E\stepcounter{descriptcount}\thedescriptcount)~,align=left]
	\item[Q-nets:] Q-nets with $\Z^2$ combinatorics were first considered by Sauer \cite{sauerqnet}, Q-nets on $\Z^N$ by Doliwa and Santini \cite{doliwasantiniqnet} as quadrilateral lattices. We show in Section \ref{sec:qnets} that Q-nets defined on quad-graphs correspond to TCD maps on TCDs with particular combinatorics. The propagation of initial data via cube flips corresponds to a sequence of 2-2 moves. We show in Section \ref{sec:laplacedarboux} that Laplace-Darboux dynamics \cite{doliwalaplace} for Q-nets corresponds to a global sequence of 2-2 moves. We show in Lemma \ref{lem:laplaceinvprojvar} that the Laplace invariants \cite{doliwalaplace} of a Q-net coincide with the projective cluster variables of the corresponding TCD map.
	\item[Darboux maps:] Darboux maps were introduced by Schief \cite{schieflattice}. We show in Section~\ref{sec:darbouxmap} that Darboux maps defined on quad-graphs correspond to TCD maps on TCDs with particular combinatorics. The propagation of initial data via cube flips corresponds to a sequence of 2-2 moves.
	\item[Line complexes:] Line complexes on $\Z^N$ were introduced by Bobenko and Schief \cite{bobenkoschieflinecomplexes} as discrete fundamental line complexes, building on work on discrete line congruences by Doliwa, Santini and Mañas \cite{dsmlinecongruence}. We show in Section \ref{sec:linecomplex} that line complexes defined on quad-graphs correspond to TCD maps on TCDs with particular combinatorics. The propagation of initial data via cube flips corresponds to a sequence of 2-2 moves.
	\item[Line compounds:] We newly introduce line compounds as a different generalization of line complexes from $\Z^3$ to $\Z^N$ with $N>3$ in Section \ref{sec:lcco}. We show that line compounds defined on quad-graphs correspond to TCD maps with particular combinatorics. The propagation of initial data via cube flips corresponds to a sequence of 2-2 moves.
	\item[Desargues maps:] Desargues maps were introduced by Doliwa \cite{doliwadesargues}. We show in Section \ref{sec:desargues} how Desargues maps are TCD maps, and that the propagation of initial data corresponds to a sequence of 2-2 moves that ``pulls through a strand''.
	\item[Ideal hyperbolic triangulations:] We show in Section \ref{sec:extriangulation} how ideal hyperbolic triangulations are TCD maps and how shear coordinates \cite{pennerteichmuller} of the ideal hyperbolic triangulation correspond to the projective cluster variables of the corresponding TCD map. We also show that every edge flip in the ideal hyperbolic triangulation corresponds to a spider move.
	\item[Projective flag configurations:] Projective flag configurations are a special case of objects considered by Fock and Goncharov in their study of higher Teichmüller theory \cite{fghighertm}. We show in Section \ref{sec:fgmoduli} how projective flag configurations are TCD maps. We also show that the TCD and the associated bipartite graph $\pb$ already occur in work of Fock and Goncharov \cite{fghighertm,gwebs}. We show that the $\mathcal X$ variables of Fock-Goncharov coincide with the projective cluster variables of the corresponding dual TCD map, and we conjecture that the $\mathcal A$ variables coincide with the $\tau$-variables of the corresponding dual TCD map. We also show how edge flips in projective flag configurations correspond to sequences of 2-2 moves.
\end{enumerate}

The following examples allow for general TCD combinatorics but make constraints on the positivity of the cluster variables and require certain embeddedness properties.
\begin{enumerate}[resume,font=\bfseries E\stepcounter{descriptcount}\thedescriptcount)~,align=left]
	\item[T-graphs:] T-graphs were introduced by Kenyon and Sheffield \cite{kenyonsheffield}. We show in Section \ref{sec:tgraphcluster} how T-graphs are TCD maps. We show that the dimer face weights for T-graphs introduced by Kenyon and Sheffield coincide with the affine cluster variables of the corresponding TCD map. We show that resplits in the corresponding TCD map correspond to moves on the T-graph. We also explain how the bipartite graphs appearing in T-graphs relate to the bipartite graphs appearing in TCD maps.
	\item[Circle patterns and t-embeddings:]A t-embedding \cite{clrtembeddings} consists of the centers of a circle pattern. T-embeddings are also known as conical nets \cite{muellerconical} or Coulomb gauge \cite{kenyonlam}. In Section \ref{sec:cptemb} we show how t-embeddings are TCD maps. We show that Miquel dynamics \cite{ramassamymiquel} correspond to a global sequence of 2-2 moves, this was published in a previous paper by the author \cite{amiquel} and found independently and simultaneously by Kenyon, Lam, Russkikh and Ramassamy \cite{kenyonlam}.	We show that the affine cluster variables of the corresponding TCD map coincide with the dimer face weights introduced in \cite{kenyonlam} and equivalently with the star-ratios introduced in \cite{amiquel}.
	\item[Circle patterns and u-embeddings:]We define a new type of map called a u-embedding that consists of certain intersection points of circle patterns in Section \ref{sec:cptemb}. We show that Miquel dynamics correspond to a global sequence of 2-2 moves. U-embeddings are also interesting because they feature real projective cluster variables, unlike t-embeddings that feature real affine cluster variables. We also explain additional relations between t- and u-embeddings.
\end{enumerate}

The following maps are particular cases of Q-nets, Darboux maps or line compounds and thus TCD maps with special combinatorics. Additionally these maps can be characterized as being precisely the TCD maps with particular combinatorics that have projective cluster variables in a either the resistor or the Ising subvariety.

\begin{enumerate}[resume,font=\bfseries E\stepcounter{descriptcount}\thedescriptcount)~,align=left]
	\item[K{\oe}nigs nets:] K{\oe}nigs nets are due to Bobenko and Suris \cite{bsmoutard} and independently due to Doliwa \cite{doliwatnets}, the latter calls them B-quadrilateral lattices. In Section~\ref{sec:konigs} we show that K{\oe}nigs nets are exactly the Q-nets with projective cluster variables in the resistor subvariety.
	\item[Carnot maps:] Schief introduced Carnot maps as a consistent reduction of Darboux maps \cite{schieflattice}, which he calls CKP maps. In Section \ref{sec:darbouxckp} we show that Carnot maps are exactly the Darboux maps with projective cluster variables in the Ising subvariety.
	\item[Doliwa compounds:] We introduce Doliwa compounds ourselves in Section \ref{sec:Doliwacomplex}. They are related to a different (and non-equivalent) definition of K{\oe}nigs nets due to Doliwa \cite{doliwakoenigs}. We show that Doliwa compounds are exactly the line compounds with projective cluster variables in the resistor subvariety.
	\item [CQ-nets:] Doliwa introduced CQ-nets as C-quadrilateral lattices \cite{doliwacqnet}, a reduction of Q-nets that features the CKP equation in an affine gauge. In Section \ref{sec:cqnets} we show that CQ-nets are exactly the Q-nets with affine cluster variables in the Ising subvariety. Equivalently, CQ-nets are exactly the Q-nets such that the section with the hyperplane at infinity is a Carnot map.
\end{enumerate}

The following maps are special cases of Q-nets, Darboux maps, line complexes or line compounds and thus TCD maps with particular combinatorics. Additionally these maps are TCD maps that have cluster variables in a particular subvariety, but this does not fully characterize the maps.

\begin{enumerate}[resume,font=\bfseries E\stepcounter{descriptcount}\thedescriptcount)~,align=left]
	\item[Q-nets in quadrics:] Q-nets with points in a quadric were studied by Doliwa \cite{doliwaqnetsinquadrics}. We show in Section \ref{sec:quadricqnet} that the section of such a Q-net with an isotropic line is a Doliwa compound. Equivalently, the affine cluster variables with respect to the isotropic line are in the resistor subvariety. We also show that a tangent section of a Q-net with points in a cone is a Schief map.
	\item[Linear line complexes:] Linear line complexes were introduced by Bobenko and Schief \cite{bobenkoschieflinecomplexes} as discrete fundamental linear line complexes. We show in Section \ref{sec:linearlinecomplexes} that the section of a linear line complex with an isotropic line is a Carnot map. Equivalently, the affine cluster variables with respect to the isotropic line are in the Ising subvariety. We also show that the section of a linear line complex with a hyperplane is an S-graph with respect to the polar point of the hyper plane.
	\item[Schief maps:] Schief introduced Schief maps as BKP maps \cite{schieflattice}. We show in Section~\ref{sec:pardm} that Schief maps are a special case of Darboux maps, such that the section with the hyperplane at infinity is a Doliwa compound. Therefore, the affine cluster variables of a Schief map are in the resistor subvariety.
	\item[Reciprocal figures:] King and Schief \cite{kstetraoctacubo} investigated reciprocal figures and detected two different BKP structures. In Section \ref{sec:reciprocal} we show that reciprocal figures are Doliwa compounds and that both BKP structures are equivalent to the observation that the projective cluster variables of a lift and of a section of the reciprocal figures are in the resistor subvariety.
\end{enumerate}

The following maps have certain lifts that are special cases of Darboux maps or line complexes and thus TCD maps with particular combinatorics. Additionally these maps are TCD maps that have cluster variables in a particular subvariety.

\begin{enumerate}[resume,font=\bfseries E\stepcounter{descriptcount}\thedescriptcount)~,align=left]
	\item[A-nets:] A-nets were introduced by Sauer \cite{saueranet}. The Plücker lift of an A-net was first considered by Doliwa \cite{doliwaanetspluecker}. Doliwa also showed that the Plücker lift of an A-net is an isotropic line complex. We show in Section \ref{sec:anets} that this implies that the section of the Plücker lift of an A-net with an isotropic plane of the Plücker quadric is a Doliwa compound. Equivalently, the affine cluster variables with respect ot the isotropic plane of the Plücker lift of an A-net are in the resistor subvariety. We also give a new proof for the fact that the projection of an A-net is a K{\oe}nigs net, and show that the BKP structures of the section of the Plücker lift coincides with the BKP structure of the projection of the A-net. 
	\item[Cox lattices:] Cox lattices were introduced by King and Schief \cite{kscox}. We introduce the Plücker-lift of a Cox map in Section \ref{sec:cox}, and show that the Plücker lifts are isotropic Darboux maps. We show that the section of the Plücker lift of a Cox map with a 3-space tangent to the Plücker quadric is a K{\oe}nigs net. Additionally, we show that the section of the Plücker lift of a Cox map with an isotropic line of the Plücker quadric is a Doliwa compound. Accordingly, the affine cluster variables of the Plücker lift of a Cox lattice with respect to the tangent space as well as to the isotropic line are in the resistor subvariety. We also explain that a hyperplane section of the Plücker lift of a Cox lattice is the Plücker lift of an A-net.
	\item[Anti-fundamental line-circle complexes:]Bo\-ben\-ko and Schief introduced anti-fun\-da\-men\-tal line-circle com\-plex\-es \cite{bobenkoschiefcirclecomplexes}, as a reduction of linear line complexes. In Section \ref{sec:antifundamental} we show that the Blaschke lift of an anti-fundamental line-circle complex is a Cox lattice. Via the Plücker lift of the Blaschke lift, an anti-fundamental line-circle complex is accompanied by two BKP structures.
\end{enumerate}

The following maps are special cases of Q-nets or Darboux maps and thus TCD maps with particular combinatorics. Additionally these maps are TCD maps that have cluster variables in a particular subvariety and satisfy a positivity constraint.

\begin{enumerate}[resume,font=\bfseries E\stepcounter{descriptcount}\thedescriptcount)~,align=left]
	\item[Harmonic embeddings:] Harmonic embeddings were first considered by Tutte \cite{tutteembedding} and are also known as \emph{Tutte embedding}. To each harmonic embedding one can associate an h-embedding \cite{kenyonlam}. We show in Section \ref{sec:harmonicemb} that h-embeddings are a special case of Schief maps, and thus also a special case of Darboux maps. As a consequence of being a special case of Schief maps the affine cluster variables of h-embeddings are in the resistor subvariety, which coincides with a result of \cite{kenyonlam}.
	\item[Orthodiagonal maps:] Orthodiagonal maps naturally accompany h-embeddings. We show in Section \ref{sec:harmonicemb} that while h-embeddings are tangent sections of a Q-net with points in a conic, orthodiagonal maps are the corresponding stereographic projections. We also show that orthodiagonal maps are u-embeddings and that the projective (resp. affine) cluster variables of an orthodiagonal map coincide with the affine (projective) cluster variables of the corresponding h-embedding.
	\item[S-embeddings:] S-embeddings were introduced by Chelkak \cite{chelkaksembeddings}. We show in Section \ref{sec:sembeddings} that s-embeddings are a special case of fixed focal maps, a special case of Q-nets that we introduce ourselves. We show that fixed focal maps and therefore s-embeddings have affine cluster variables in the Ising subvariety with respect to any point in $\CP^1$. The corresponding result for the point at infinity was also shown in \cite{kenyonlam}. Note that only the affine cluster variables with respect to the point at infinity are necessarily real positive.
	\item[S-graphs:] S-graphs were introduced by Chelkak \cite{chelkaksgraphs}. We show in Section \ref{sec:sgraphs} that S-graphs are Q-nets that satisfy an additional geometric constraint with respect to a distinguished point at infinity. We also show that the section of an S-graph with a line through the distinguished point is a Carnot map. Therefore the affine cluster variables of an S-graph are in the Ising subvariety with respect to any such line. We also show the projection of an S-graph from the distinguished point is an s-embedding. We also show that S-graphs occur as sections of linear line complexes, and as projections of certain CQ-nets.
\end{enumerate}

Finally, let us mention two examples that do not fall into one of the previous categories.
\begin{enumerate}[resume,font=\bfseries E\stepcounter{descriptcount}\thedescriptcount)~,align=left]
	\item[Quadrirational Yang-Baxter maps:] Quadrirational Yang-Baxter maps were introduced by Adler, Bobenko and Suris \cite{absyangbaxter}. We observe in Section \ref{sec:quadrirationalybmaps} that quadrirational Yang-Baxter maps are a reduction of Darboux maps. Moreover, generically there are four distinguished points for a quadrirational Yang-Baxter map. We show that the section of a quadrirational Yang-Baxter map with respect to a line through any two of the distinguished points is a Doliwa compound. Therefore, the affine cluster variables of a quadrirational Yang-Baxter map with respect to such a line are in the resistor subvariety. We also show that quadrirational Yang-Baxter maps are sections of certain Q-nets inscribed in the intersection curve of two quadrics. 
	\item[Pentagram map:] The pentagram map as a discrete evolution of polygons was introduced by Schwartz \cite{schwartz}. In  Section \ref{sec:pentagram} we explain how the polygons can be considered to be doubly periodic Q-nets on which the pentagram map acts as Laplace-Darboux dynamics. This identification has already been sketched by Schief in a notorious talk \cite{schieftalk}, but has not been published. We also show that the cluster variables introduced by Glick \cite{glickpentagram} coincide with the projective cluster variables that we can associate to polygons via the doubly periodic Q-nets.
	\item[Circular Q-nets]Circular Q-nets \cite{bobenkocircular,cdscircular} are Q-nets such that every quad is inscribed in a circle. In Section \ref{sec:circularnets} we explain how the affine cluster variables can be expressed as a product of the BMS-variables (see \cite{bmscircular}), and how the canonical Poisson bracket of the cluster algebra corresponds to the ultra-local BMS Poisson bracket.
\end{enumerate}

Let us also note that there are other periodic reductions of TCD maps that are of interest, but beyond the scope of this thesis. For example Konopelchenko and Schief discovered \cite{ksclifford} how \textbf{Schramm circle packings} \cite{schramm} and \textbf{discrete holomorphic functions} \cite{bpdisosurfaces} can be understood as periodic reductions of the dSKP equation. There is also  \textbf{polygon recutting} \cite{adlerrecutting}, for which Izosimov \cite{izosimovrecutting} recently found a cluster structure. In upcoming and joint work with Melotti and de Tilière \cite{amdtdskp} we expand on Schief's and Izosimov's results, and additionally show how \textbf{integrable cross-ratio systems} \cite{bmsconformal} and \textbf{circle intersection dynamics} \cite{gdevron} are periodic reductions of the dSKP equation and thus TCD maps. We also focus more on periodic TCD maps in a collaboration with George and Ramassamy \cite{agrcrdyn} with special attention given to discrete holomorphic functions. Other periodic maps are various notions of \textbf{higher pentagram maps}, we claim that the projective quivers for these maps as TCD maps can be found in work of Glick and Pylyavskyy on \textbf{Y-meshes} \cite{gpymesh}, without providing further details here.

\section{Open questions and future directions}\label{sec:introfuture}

The following list is long but we think each item is of significant interest.

\setcounter{descriptcount}{0}
\begin{enumerate}[font=\bfseries Q\stepcounter{descriptcount}\thedescriptcount)~,align=left]
	\item[Positivity and discrete surfaces:]We have encountered some TCD maps with a cluster structure with strictly positive cluster variables, namely t-embeddings, {T-graphs} and their reductions, h-embeddings, s-embeddings and S-graphs. It would be interesting to understand the geometric interpretation of line complexes, Q-nets, Darboux maps and line compounds having a positive cluster structure. The geometric interpretation is not completely obvious, for example Q-nets with positive projective cluster variables have only quads that are non-convex.
	
	\item[Continuous limits:]How do we take continuous limits of TCD maps and what are the limits? For some of the DDG examples limits are well understood (see \cite[Chapter~5]{ddgbook}). For example the limit of a Q-net is a surface in conjugate parametrization, the limit of an A-net is a surface in asymptotic parametrization. For general TCD maps however, it is not clear how to refine the combinatorics while taking the limit or what the interpretation of the strands in some sort of limit surface is. Note that with our conventions, the limit of a TCD map corresponding to a Q-net is actually three disjoint surfaces: the surfaces that the Q-net converges to and the two focal surfaces.

	\item[Dimers and positive surfaces in the limit:]Of course, the limits of dimer partition functions and associated objects (like the height function) have been studied extensively in the literature (see for example \cite{ckpdimers}). Thus the question is, can we give a projective geometric interpretation to the limits of the partition functions in the case of a limit of a TCD map? For example, assuming we have a Q-net that has a positive cluster structure and a limit of the Q-net that is a surface in conjugate parametrization: What is the meaning of the limit height function of the dimer model for the limit surface? Note that to some extent, this has been studied for S-graphs in statistical mechanics, and we have shown in this thesis that S-graphs are reductions of Q-nets. However, S-graphs were only studied in $\RP^2$, which makes it non-obvious how to interpret them as surfaces.
	
	\item[Grassmannians:]Although this was not the original intent, many of the methods we employ are very similar to the methods employed in the study of the moduli space of (totally positive) Grassmannians \cite{postgrass}. In the end, we strongly believe that in an appropriate sense, TCD maps are the dual story to Postnikov's ideas. However, what is the precise correspondence? Is there some geometric meaning to the Grassmannian that corresponds to a TCD map? As we often study for example Q-nets far away from their maximal dimension, how does this translate to Grassmannians? Also every specific reduction of a TCD map should therefore correspond to specific reductions of Grassmannians, so what is the Grassmannian of an A-net? A K{\oe}nigs net? And so on? There is the converse question as well, for example there are the so called orthogonal Grassmannians that relate to the Ising model \cite{gpisingorthogonal}, how does this relate to TCD maps?

	\item[Grassmannians and limits:]The questions above can also be applied to the limits. Note that a possible process of taking a limit of a Q-net is to use domains of definition that become larger and larger subsets of $\Z^2$. This means that the maximal dimension of the Q-net increases while taking the limit. What are the implications for the Grassmannian? Can we or do we need to consider infinite dimensional Grassmannians?

	\item[Dimer configurations:]Although the dimer face weights as well as the dimer partition functions have an interpretation in terms of TCD maps, there is no obvious direct interpretation of the meaning of the dimer configurations. 
	
	\item[On the torus:]An important next step is to understand TCD maps, such that the underlying TCD is defined on the torus not in the disc. The importance is because we claim that many well-studied systems are actually TCD maps on the torus, like the pentagram map (see Section \ref{sec:pentagram}). In \cite{agrcrdyn} we show that cross-ratio dynamics can also be understood as dynamics on TCD maps on the torus. We have also outlined some more general results in \cite{agrcrdyn}, but it would certainly be helpful to have a complete and detailed understanding of TCD maps on the torus. While we do not know of any obvious geometric examples on other surfaces, it would also be interesting to understand TCD maps on other surfaces. On the dimer side, interesting research has already been done, see for example \cite{crdimerspins}.
	
	\item[$Z$-variables, $X$-variables and Muller-Speyer:] Can we use the results of Muller-Speyer \cite{mstwist} to interpret (alternating ratios of) the partition functions $Z$ of a TCD map $T$ as the $X$ variables of the dual TCD map $T^\star$? How do the results of Section \ref{sec:projclusterduality} about the coincidence of cluster variables in the primal and dual translate to the viewpoint of Muller-Speyer?
	
	\item[The cluster variables of all sections:] Consider a TCD map and all its sections, see also Remark \ref{rem:envelopingalgebra}. We observed that except for a finite number of very small TCD maps, the number of all projective cluster variables in all sections is larger than the dimension of the space of TCD maps and choices of sections. Therefore, these cluster variables are not all independent. What are the relations?

	\item[Reductions of flag configurations:]Projective flag configurations are a special case of the theory of Fock-Goncharov moduli-spaces \cite{fghighertm} for the projective groups. However, it is possible to consider certain reductions of projective flag configurations, for example flags that are in particular position with respect to bilinear forms. These reductions are reflected in subvarieties of the cluster variables, which is part of ongoing research. Are these reductions of projective flag configurations related to the more general Fock-Goncharov theory?
	
	\item[De/genericity:]We claim that the TCD map that are not flip-generic can be recognized as the TCD maps that have projective cluster variables not in a union of certain subvarieties. For example, the TCD map with four points on a line and cross-ratio 1 is not flip-generic, which is reflected in the fact that the unique projective cluster variable is -1. How do these subvarieties look in general? We also claim that non-flip-genericity is reflected in some of the almost perfect matching partition functions being 0. Can this be formalized? By the previous claims TCD maps that are not $1$-generic should be recognizable via their affine cluster variables. Can we also recognize TCD maps that are not $k$-generic in the projective cluster variables? Note that whether a TCD map attains maximal dimension can not be read off the projective cluster variables because the projective cluster variables are invariant under projections. In general, it would be interesting to precisely understand the (positive and general) moduli-space of TCD maps beyond the non-generic case. We assume the latter question is in some sense dual to the theory of totally positive Grassmannians as proposed by Postnikov \cite{postgrass}.
	
	\item[Non-generic projections/sections:]What happens when we apply non-generic projections and sections to generic TCD maps? Should we adapt the combinatorics? When doing dynamics in terms of local moves, are there cases when it is possible to pass through the generated singularities? Should one consider the lift of the projection, that is the original map as a sort of blow-up of the projected TCD map and what degrees of freedom are there?
	
	\item[Non-minimal TCDs and subvarieties:]We claim that the cluster variables of TCD maps defined on non-minimal TCDs are in certain subvarieties. For example, if we consider ideal hyperbolic triangulations with interior vertices in the triangulation and we label the $m$ projective cluster variables (shear coordinates in this case) around an interior vertex by $X_1,X_2,\dots,X_m$, then the well-known closing conditions are
	\begin{align}
		X_1 \cdot X_2 \cdots X_m = 1 \quad \mbox{ and }\quad X_1 + X_1\cdot X_2 +  \dots + X_1\cdot X_2 \cdots X_m = 0.
	\end{align}
	We conjecture that any closed loop in a TCD contributes two such equations for the corresponding TCD map. What about other violations of non-minimality? This is research in progress.
	
	\item[Bilinear forms and reductions:] Is it an accident that for TCD maps related to symmetric bilinear forms (Q-nets in quadrics, see Section \ref{sec:quadricqnet}) the BKP equation appears and for maps related to anti-symmetric bilinear forms (Linear line complexes, see Section \ref{sec:linearlinecomplexes}) the CKP equation appears? We have shown that these particular TCD maps feature BKP respectively CKP when one considers sections with isotropic lines. What exactly are all particular TCD maps that feature BKP or CKP in a section with a line? The latter question is already somewhat clear from our remarks to the corresponding theorems. But this question can also be asked for many of the other maps. For example, any section with a line through the distinguished point of an S-graph is a Carnot map, but does the converse hold? In general, what about higher-dimensions, for example Q-nets with vertices in a quadric in $\CP^4$ instead of $\CP^3$?
	
	\item[Missing BKP/CKP identifications:]Bobenko and Schief found the CKP equation accompanying linear line complexes \cite{bobenkoschieflinecomplexes} via M-systems. We did so as well via sections with an isotropic lines (Section \ref{sec:linearlinecomplexes}), but we have not worked out if and how these two occurrences of CKP coincide. The structure that Bobenko and Schief found is formulated in terms of the $\tau$-variables, our results of Section \ref{sec:tautcd} on $\tau$-variables of TCD maps may be useful. Note that because there is a 3-parameter family of isotropic lines, we found a whole 3-parameter family of occurrences of CKP, how does this relate to the results of Bobenko and Schief? Similarly, King and Schief found the BKP equation in Cox lattices via lifts and the choice of a distinguished bilinear form \cite{kscox}. We found two occurrences of families of BKP equations (Section \ref{sec:cox}), but via sections with isotropic lines as well as tangent 3-spaces in the Plücker lift of Cox lattices. How do these findings relate?

	\item[Bilinear forms and TCD maps with general combinatorics]The resistor and Ising subvarities occur for TCD maps with certain admissible combinatorics. Ultimately, this stems from the two bijection tricks between the dimer model and the spanning tree and the Ising model respectively \cite{temperleybijection, dubedattrick}. However, it is clear that one can do local dimer moves away from the admissible combinatorics, then return to admissible combinatorics and land in the subvarities again, due to the consistency of the system. Therefore, the corresponding subvarities should also exist for non-admissible combinatorics. This becomes especially relevant as we believe that the same holds for the geometric properties of the corresponding TCD maps. For example, we believe it is possible to recognize whether a TCD map is flip equivalent to a Q-net in a quadric via certain incidence and polarity lemmas. This is related to not yet published results on discrete surface patches in specific parametrizations by Fairley \cite{fairleythesis}.
	
	\item[Non-commutative TCD maps:]It is rather clear that the definitions of TCD maps and VRC can be generalized to relations between projective Grassmannians, much in the spirit of \cite{absgrassmannian}. Note that a priori this is not the same as the moduli-space theory of Grassmannians as done by Postnikov. In fact, we believe such a non-commutative theory has interesting geometric features that vanish in the commutative case. This could also provide a clue to non-commutative cluster structures, although this is a difficult field.
	
	\item[TCD maps and the Darboux system:]We have ignored one very interesting approach to Q-nets via the so called \emph{Darboux system} \cite{bkdarboux, doliwasantiniqnet}. This system can be related by taking advantage of the specific combinatorics of Q-nets and affine gauge as well as additional specific choices on the relations. We believe the Darboux system or at least parts thereof (including the rotation coefficients) can also be found generally in TCD maps, which would essentially constitute a $A_N$-lattice version of the Darboux system.

	\item[TCD maps and M-systems:]M-systems make an appearance in work by Bobenko and Schief \cite{bobenkoschieflinecomplexes} specifically for the case of line complexes, but are clearly intended for much more general purposes. Can we relate M-systems and TCD maps, or is there an $A_N$-lattice version of M-systems? Note that M-systems are based on minor-relations of certain matrices, possibly indicating a very geometric relation to TCD maps and cluster algebras.

	\item[Signotopes for TCDs:]There is a higher dimensional generalization of pseudoline arrangements called \emph{signotopes} \cite{fwsignotopes}. We would like to know whether an analogous higher dimensional generalization exists for triple crossing diagrams. In particular, the theory of signotopes comes with a lot of structure including a full understanding of the flip-graph. Thus there is a chance that a higher dimensional generalization of TCDs might help to find a higher dimensional integrable system as well.
	
	\item[Poisson-Lie groups:]Can we relate results of Izosimov \cite{izosimovpoissonlie} on relations between Poisson-Lie groups and pentagram maps to TCD maps on the torus, possibly to TCD maps in some sense in general?
	
	\item[Section combinatorics:]During the completion of this thesis it was brought to the author's attention that operations similar to how we define sections (see Section \ref{sec:sections}) have appeared in the recent literature \cite{galashincritical, psbwamplituhedron}. What are the relations to our work?
	
	\item[Aztec diamond theorem for TCD maps:]We generalize the Aztec diamond theorem \cite{speyerdimers} from dKP to dSKP in upcoming joint work with Melotti and de Tilière \cite{amdtdskp}. This provides another interesting relation between partition functions and TCD maps, although we only cover the case of combinatorics corresponding to $A_3$. Generalizing this to all combinatorics is still work in progress.
	
	\item[Complex TCD maps via real geometry:]We have shown how t-embeddings, h-embeddings and s-embeddings are special cases of TCD maps in affine charts of $\CP^1$ that have positive affine cluster variables. We have related these maps in many ways to TCD maps in general and to specific examples, especially from DDG. However, we do not think that we have given appropriate attention to the geometric meaning of the reality constraint on the affine cluster variables of these TCD maps in $\CP^1$. We intend to give a more enlightening explanation in upcoming work with Chelkak.
	
	\item[Integrable 2D systems:]We have focused almost exclusively on 3D integrable systems in this thesis. It would be interesting to understand what 2D systems fit into the TCD map framework as well. Indeed, there are several 2D systems that can be viewed as periodic reductions of TCD maps, because they are periodic reductions of the dSKP equation and thus periodic reductions of TCD maps, see \cite{amdtdskpb} for an exposition of examples. However, this covers only one partial example (Q1 with $\delta = 0$) of the ABS classification of 2D integrable quad systems \cite{absquads}. What about all the other examples of the ABS classification? What about face-centered quad equations \cite{kelsfcc}? The question can also be posed the other way around: is it possible to classify reasonable 2D reductions of TCD maps? In the scalar case, the latter question is asking for a classification of 2D systems on the $A_N$ lattice.

	\item[Subdivision:]Is it possible to subdivide TCD maps in a projectively invariant manner? How to subdivide the combinatorics? It is a notoriously different problem to subdivide discrete integrable systems, as it is difficult to subdivide ``close'' to the original map while preserving the ``rigidity''.

	\item[Lelieuvre normal fields:]One interesting aspect of DDG that we did not yet relate to TCD maps are the so called Lelieuvre normal fields \cite{kplelieuvre}, which are associated to A-nets. In the case that there is a relation, it would be interesting to see if one can transfer some results to other TCD map examples.
	
	\item[Tropical geometry:]Is there a reasonable theory of TCD maps in projective tropical geometry (possibly as in \cite{rgsttropical})? Are there other relations to tropical geometry, possibly via cluster algebras or Fock-Goncharov moduli-spaces? What about relations to ultra-discrete systems?
	
	\item[Finite geometry:]Although we do not know if this question is of importance, but what about considering TCD maps for projective spaces over finite fields? In this case, there will only be a finite number of TCD maps for every TCD, can we enumerate them?
	
	\item[Scattering theory:]Relations between positive Grassmannians and scattering amplitudes have been developed in \cite{scattamp}. Indeed, \cite[Section 5]{scattamp} already gives some thought to point configurations similar to TCD maps. It would be interesting to understand this in more depth, can we understand the geometric meaning of the scattering amplitudes in some way, especially in the case of specific examples like Q-nets? 

	\item[Random walks:]For T-graphs Kenyon and Sheffield \cite{kenyonsheffield} also introduced certain martingale random walks using the edge-weights of the affine cluster structure. Can we generalize this procedure to arbitrary TCD maps with positive affine cluster variables? Note that T-graphs already allow for general combinatorics, but the particular embedding properties of T-graphs are used to choose the direction of the random walks. Furthermore, is it also possible to understand these random walks for positive projective cluster structures? It would be interesting to understand the geometric meaning of these walks on DDG examples. Also, as many DDG examples generically feature negative cluster variables, is it possible to give an interpretation of the random walks even in this case?
	
	\item[Random TCD maps:]When Dylan Thurston introduced TCDs \cite{thurstontriple} he motivated his ideas by showing that TCDs are a generalization of dimer configurations on subgrids of $\Z^2$, so called \emph{domino tilings}. More specifically, to each domino tiling we associate a TCD by replacing every domino
	
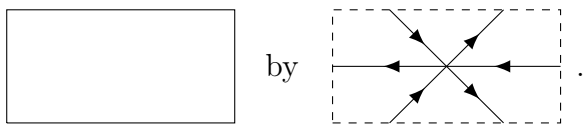
\begin{center}
		\begin{tikzpicture}[scale=1.5, baseline={([yshift=-.7ex]current bounding box.center)}]
			\draw[-]
				(0,0) -- (2,0) -- (2,1) -- (0,1) -- (0,0)
			;
		\end{tikzpicture}\hspace{4mm}by\hspace{3mm}
		\begin{tikzpicture}[scale=1.5, baseline={([yshift=-.7ex]current bounding box.center)}]
			\draw[dashed]
				(0,0) -- (2,0) -- (2,1) -- (0,1) -- (0,0)
			;
			\draw[-]
				(1,0.5) edge[mid arrow] (0, 0.5) edge[mid arrow] (1.5, 0) edge[mid arrow] (1.5, 1)
				(1,0.5) edge[mid rarrow] (2, 0.5) edge[mid rarrow] (0.5, 0) edge[mid rarrow] (0.5, 1)
			;
		\end{tikzpicture}\hspace{2mm}.
	\end{center}
	Thurston illustrated that \emph{domino flips} correspond to 2-2 moves. On the other hand, the sampling of dimer configurations is extensively studied (see Section \ref{sec:dimers}). Is it possible to extend the sampling idea from dimers to TCDs? Are there weights for such a sampling that relate to TCD maps?
	
	\item[Double nets:]There is an interesting new development in DDG to describe discretizations via coupled pairs of particular nets \cite{bsstconfocala,bsstconfocalb,techterthesis}. The integrability of these coupled nets is still opaque so far. In joint work with Techter we aim to solve the integrability issue of coupled nets with the help of insights from TCD maps.

	\item[Double random currents:]In \cite{dlrandomcurrent} a new reduction of planar dimer models (related to the Ising model) is introduced. Does this reduction correspond to something in geometry via TCD maps? Are there other reductions of the dimer model and how do they relate to geometry? 
	
	\item[Schubert dynamics:] The following system was proposed by Glick \cite{gdevron} due to some interesting numeric evidence. Consider four given lines in $\CP^3$, then in a generic situation there are exactly two lines that intersect the four given lines. This idea can be used to defined octahedral-type dynamics (as in the case of Miquel or Laplace-Darboux dynamics). It is not clear and we have not tested whether these dynamics are discretely integrable and it is not clear if they fit in some way into the TCD map framework, and we leave this as an open question.
\end{enumerate}

\section{Publications}
The integrability result of t-embeddings and Miquel dynamics has already appeared in a preprint \cite{amiquel} of the author (and also independently \cite{kenyonlam}). The first projective cluster structures for Q-nets and Darboux maps have appeared in a joint preprint with Glick, Pylyavskyy and Ramassamy \cite{vrc}. The vector-relation configurations are also introduced in \cite{vrc}, together with a result on uniqueness from the boundary that is not part of the thesis. We do not know exactly in what format the results of this thesis will be published, but there will be a joint publication with Glick and Ramassamy containing parts of the thesis. Related work not included in the thesis is contained in a preprint with George and Ramassamy \cite{agrcrdyn}. More related work not included in the thesis is available in two preprints on joint work with de Tilière and Melotti \cite{amdtdskp, amdtdskpb}.

\section{Acknowledgments}
Special thanks to my supervisor Boris Springborn, Stefan Felsner, my older math-brother Ananth Sridhar and my math-sparring partner Jan Techter. Additional thanks go to Alexander Bobenko, Ulrike Bücking, Dmitry Chelkak, Albert Chern, Alexander Fairley, Terrence George, Roman Gietz, Max Glick, Richard Kenyon, Hana Dal Poz Kouřimská, Wai Yeung Lam, Paul Melotti, Christian Müller, Sanjay Ramassamy, Isabella Retter, Thilo Rörig, Andy Sageman-Furnas, Wolfgang Schief, Lara Skuppin, Nina Smeenk, Yuri Suris for countless discussions and help with proof-reading.

\chapter{TCD maps and discrete integrability}\label{cha:tcd}

\section{Triple crossing diagrams (TCDs)}\label{sec:tcds}

What we call a \emph{triple crossing diagram} was introduced by Dylan Thurston \cite{thurstontriple} as \emph{triple point diagrams} and also as \emph{triple diagrams}. Another name given to such a diagram is Thurston diagram \cite{fmloop}. The name triple crossing diagram has also appeared in the literature before \cite{bocklandtdimerabc,bwtriple}. We choose to stick with the term triple crossing diagram because it minimizes the potential for confusion. We will now give several basic definitions and theorems from Thurston's initial work before relating triple crossing diagrams to vector relation configurations in the next section.

\begin{definition}[\cite{thurstontriple}]
	A \emph{triple crossing diagram (TCD)} $\tcd$ is a collection of oriented closed intervals and circles immersed smoothly into a disk. The image of a connected component is a \emph{strand}; it is either an \emph{arc} (the image of an interval) or a \emph{loop} (the image of a circle). The immersions are required to satisfy:
	\begin{enumerate}
		\item three strands cross at each point of intersection,
		\item the endpoints of arcs are distinct points on the boundary of the disk, and no other points map to the boundary,
		\item the orientations on the strands induce consistent orientations on the complementary regions.
	\end{enumerate}
	Triple crossing diagrams are considered up to homotopy among such diagrams.
\end{definition}
This makes them essentially combinatorial objects, 6-valent graphs with some extra structure. Given a triple crossing diagram $\tcd$, we may reverse the orientation of all strands to obtain a new triple crossing diagram $\tcd'$. For the TCD maps that we will introduce in Section \ref{sec:tcdmap} the orientation of the strands will matter, thus we really consider $\tcd$ and $\tcd'$ as two fundamentally different diagrams. Also note that the strands at the boundary have to alternate between in- and out-endpoints. This is because along each face of the diagram including the faces at the boundary the orientations of the strands are consistent. 

\begin{definition}\label{def:labeledtcd}
	Let $\tcd$ be a TCD with $n$ strands of which $m$ are arcs, not loops. We call $\tcd$ a \emph{labeled TCD} when we also fix a labeling of its $n$ strands by the numbers $1,2,\dots,n$ such that
	\begin{enumerate}
		\item no two strands carry the same label,
		\item the labels $1,2,\dots,m$ appear in counterclockwise order at the in-endpoints.\qedhere
	\end{enumerate}	
\end{definition}

\begin{definition}\label{def:endpointmatching}
	Let $\tcd$ be a labeled triple crossing diagram with $m$ in-endpoints. Label the out-endpoints by numbers $1,2,\dots, m$ in counterclockwise order such that 1 is the first out-endpoint after the in-endpoint 1 in counterclockwise direction. Let $C_\tcd$ be the permutation such that $C_\tcd(i) = j$ if there is a strand from in-endpoint $i$ to out-endpoint $j$. We call $C_\tcd$ the \emph{endpoint matching}. For labeled TCDs we consider two endpoint matchings to be equivalent if they differ by conjugation with a cyclic permutation. We say strand $i$ has \emph{length} $k\in \Z_m$ if $C_\tcd(i) = i+k$. We say a triple crossing diagram with $n$ strands has endpoint matching $\enm nk$ if every strand has length $k$. If a triple crossing diagram $\tcd$ has endpoint matching $\enm nk$ for some $k\in \N$ we say $\tcd$ is \emph{balanced}.
\end{definition}

\begin{figure}
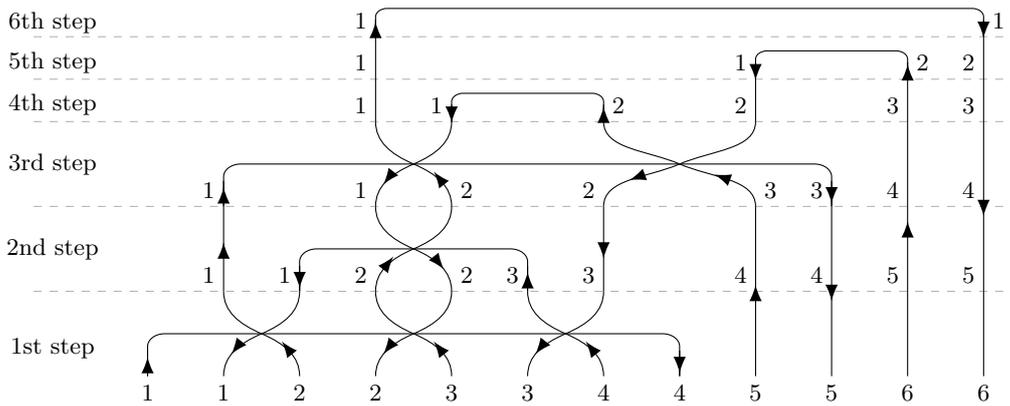

\hspace{-0.5cm}
\vspace{-4mm}
\caption{A half-plane drawing of a standard diagram for the endpoint matching $\enm 63$.}
\label{fig:standardtcd}
\end{figure}

An interesting question is whether there are global restrictions on the endpoint matching. Let us speak of a \emph{half-plane drawing} of a TCD when we draw the in- and out-endpoints on the real axis in ascending order and the whole diagram is contained in the upper half-plane $\{z\in \C \ : \ \Im(z) \geq 0\}$, see Figure \ref{fig:standardtcd}.

\begin{definition}[\cite{thurstontriple}]\label{def:standardtcd}
	Let $\tcd$ be a labeled TCD with $n$ strands and without loops. Let $C_\tcd$ be the endpoint matching of $\tcd$. A \emph{standard diagram} is one that is constructed recursively (see Figure \ref{fig:standardtcd}).  Choose an interval $[k,C_\tcd(k)]$, that is minimal with respect to the partial inclusion order of intervals. Let TCD $\tcd'$ be a standard diagram with $n-1$ strands and endpoint matching
	\begin{align}
		 C_{\tcd'}(i) = \begin{cases}
		 	C_\tcd(i+1) - 1 & \mbox{for } C_\tcd(i+1) > C(k),\\
		 	C_\tcd(i+1) & \mbox{for } C_\tcd(i+1) < C(k).
		 \end{cases}
	\end{align}       
	We construct $\tcd$ from $\tcd'$ by adding a strand that runs along the boundary, such that we introduce $k-1$ new crossings, each crossing involving the new strand and strands $j$ and $ C_{\tcd'}^{-1}(j)$ of $\tcd'$ for $1\leq j < C(k)$. As $\tcd'$ has $n-1$ strands, the recursion ends after $n$ steps.	
\end{definition}

Note that we do indeed consider intervals, therefore $[k,C_\tcd(k)] = [C_\tcd(k),k]$. The interval that we consider in the first step in Figure \ref{fig:standardtcd} is $[1,4]$. In the second step we see one strand less and the interval is $[3,1]$. After that the intervals are $[1,3],[2,1],[2,1],[1,1]$. 

\begin{theorem}[\cite{thurstontriple}]
	In a disk with $2n$ endpoints on the boundary, all $n!$ endpoint matchings are achievable by some triple crossing diagram without loops.
\end{theorem}
\proof{Follows from Definition \ref{def:standardtcd}, where we gave an algorithm to construct a TCD for a given endpoint matching.\qed}

\begin{definition}[\cite{thurstontriple}]
	A \emph{connected triple crossing diagram} is a diagram in which the image of the immersed curves together with the boundary of the disk is connected. Equivalently, it is a diagram in which each complementary region to the image is a disk. 
\end{definition}

For example, any diagram without loops is connected.

\begin{figure}
	\centering
	\begin{tikzpicture}[baseline={([yshift=-.7ex]current bounding box.center)},scale=.75,rotate=90]
		\coordinate (ws) at (-2,-1);
		\coordinate (wn) at (-2,1);
		\coordinate (es) at (2,-1);
		\coordinate (en) at (2,1);
		\coordinate (nw) at (-1,2);
		\coordinate (ne) at (1,2);
		\coordinate (sw) at (-1,-2);
		\coordinate (se) at (1,-2);
		\coordinate (n) at (0,1);
		\coordinate (s) at (0,-1);
		
		\draw[-]
			(n) edge[bend right=45,mid rarrow] (s) edge[mid rarrow] (nw) edge[mid arrow] (ne) edge[mid rarrow] (en) edge[mid arrow] (wn)
			(s) edge[bend right=45,mid rarrow] (n) edge[mid arrow] (sw) edge[mid rarrow] (se) edge[mid arrow] (es) edge[mid rarrow] (ws)
		;		
	\end{tikzpicture}\hspace{.5cm}$\longleftrightarrow$\hspace{.5cm}
	\begin{tikzpicture}[scale=.75,baseline={([yshift=-.7ex]current bounding box.center)}]
		\coordinate (ws) at (-2,-1);
		\coordinate (wn) at (-2,1);
		\coordinate (es) at (2,-1);
		\coordinate (en) at (2,1);
		\coordinate (nw) at (-1,2);
		\coordinate (ne) at (1,2);
		\coordinate (sw) at (-1,-2);
		\coordinate (se) at (1,-2);
		\coordinate (n) at (0,1);
		\coordinate (s) at (0,-1);
		
		\draw[-]
			(n) edge[bend right=45,mid rarrow] (s) edge[mid rarrow] (nw) edge[mid arrow] (ne) edge[mid rarrow] (en) edge[mid arrow] (wn)
			(s) edge[bend right=45,mid rarrow] (n) edge[mid arrow] (sw) edge[mid rarrow] (se) edge[mid arrow] (es) edge[mid rarrow] (ws)
		;		
	\end{tikzpicture}
	\caption{The 2-2 move in a TCD.}\label{fig:twotwo}
\end{figure}
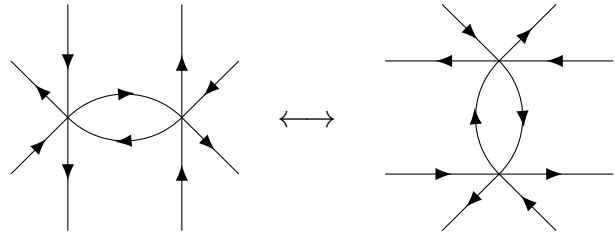

\begin{definition}[\cite{thurstontriple}]
	A \emph{2-2 move} in a TCD is the local rearrangement of strands at a bigon that preserves the endpoint matching as depicted in Figure \ref{fig:twotwo}.
\end{definition}
The 2-2 move will prove to be of central importance for the geometric dynamics that we study. Note that the orientation of the bigon is the same before and after the move. As we distinguish orientations, we will later consider the 2-2 move at a counterclockwise and the 2-2 move at a clockwise oriented bigon as two different types of 2-2 moves. One will correspond to a change in local geometry, while the other one will correspond to a change of parametrization (see Definition \ref{def:tcdmapmoves}).

\begin{definition}[\cite{thurstontriple}]\label{def:minimaltcd}
	A \emph{minimal TCD} is a connected diagram with no more triple points than any other triple crossing diagram with the same endpoint matching.
\end{definition}
Throughout the thesis, we will only consider minimal diagrams. One particularly useful property of minimal diagrams is the next theorem.

\begin{theorem}[\cite{thurstontriple}]\label{th:tcdflipsconnected}
	Any two minimal TCDs with the same matching on the endpoints are related by a sequence of 2-2 moves.
\end{theorem}

The following theorem states that one can recognize a non-minimal TCD by the absence of forbidden configurations.
\begin{theorem}[\cite{thurstontriple}]\label{th:mintcdforbidden}
	A connected TCD is minimal if and only if it has no loop, no strand which intersects itself (self intersection) or pairs of strands which intersect at two points $x$~and~$y$, with both strands oriented from $x$ to $y$ (parallel intersection).
\end{theorem}

\begin{lemma}[\cite{thurstontriple}]
	A standard diagram is minimal.
\end{lemma}

\begin{figure}
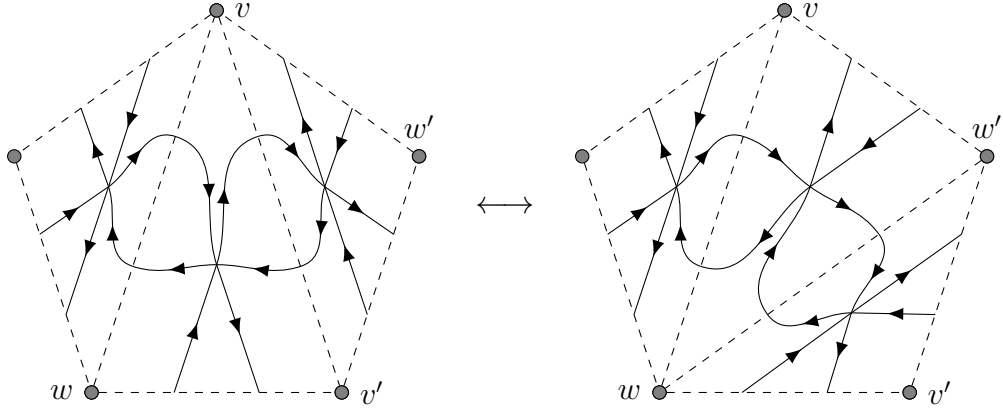


\caption{Triple crossing diagram $\tcd_{\tg}$ corresponding to a triangulation $\tg$ (dashed) of the pentagon and the effect of an edge-flip.}
\label{fig:tcdpentagon}
\end{figure}

\section{Triangulations}\label{sec:tcdtriangulations}

Before we proceed to geometry, let us give some combinatorial examples of triple crossing diagrams that can be associated to two reoccurring objects: triangulations and quad-graphs. We begin with triangulations.

\begin{definition}
	A \emph{disc triangulation}, or \emph{triangulation} for short, is a closed disc that is a CW-complex whose faces (2-cells) are triangles which are glued edge-to-edge. We denote the set of vertices (0-cells), edges (1-cells) and faces of a triangulation $\tg$ by $V,E,F$. 
\end{definition}

In general, given a triangulation that is a disc we call the vertices on the boundary of the disc \emph{boundary vertices}. Edges contained in the boundary of the disc are called \emph{boundary edges} and faces that contain a boundary edge are called \emph{boundary faces}.

\begin{definition}
	In a triangulation $\tg$, consider the local situation where two different triangles $t = (v,w,v')$ and $t'=(v',w',v)$ share exactly one edge $e = (v,v')$. Then the \emph{edge flip} at $e$ yields a new triangulation $\tg'$ which coincides with $\tg$ except that we remove the two triangles $t,t'$ and replace them with the two triangles $(w,v',w')$ and $(w',v,w)$, see also Figure~\ref{fig:tcdpentagon}.
\end{definition}

In other words, we flipped the edge $(v,v')$ to become the edge $(w,w')$. Now we associate a TCD to a given triangulation and show that an edge flip corresponds to a 2-2 move.

\begin{samepage}
\begin{definition}\label{ex:tcdtriangulation}
	Let $\tg$ be a triangulation. Into every triangle glue the following piece \\
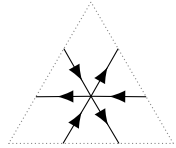
	\begin{center}\vspace{-5mm}
		\begin{tikzpicture}[scale=1.25]
			\coordinate (v1) at (90:1);
			\coordinate (v2) at (210:1);
			\coordinate (v3) at (330:1);			
			\draw[gray,-,densely dotted]
				(v1) -- (v2) -- (v3) -- (v1)
			;
			\coordinate (e12) at ($(v1)!.33!(v2)$);
			\coordinate (e13) at ($(v1)!.33!(v3)$);
			\coordinate (e23) at ($(v2)!.33!(v3)$);
			\coordinate (e21) at ($(v2)!.33!(v1)$);
			\coordinate (e32) at ($(v3)!.33!(v2)$);
			\coordinate (e31) at ($(v3)!.33!(v1)$);
			\coordinate (o) at (0,0);
			
			\draw[-]
				(o) edge[mid rarrow] (e12) edge[mid rarrow] (e23) edge[mid rarrow] (e31)
				(o) edge[mid arrow] (e21) edge[mid arrow] (e32) edge[mid arrow] (e13)
			;
		\end{tikzpicture}
	\end{center}
	of a TCD such that the pieces of strands of neighbouring triangles meet along edges. We say the resulting TCD $\tcd_{\tg}$ is the TCD corresponding to $\tg$.
\end{definition}
\end{samepage}

By construction, for each triangle of $\tg$ there is exactly one triple intersection point in $\tcd_{\tg}$. Moreover, we observe that to each vertex $v$ of $\tg$ corresponds a strand in $\tcd_{\tg}$ that is contained in all triangles incident to $v$ and passes clockwise through all the crossing points of the incident triangles. Therefore, if there is a vertex $v$ in $\tg$ that is not a boundary vertex, then the corresponding strand in $\tcd_{\tg}$ is a loop and $\tcd_{\tg}$ is not a minimal TCD. If a triangulation has $n$ boundary vertices and no interior vertices, we call it a \emph{triangulation of the $n$-gon} (see Figure \ref{fig:tcdpentagon}). The corresponding TCD is always a minimal TCD with endpoint matching $\enm n1$. Vice versa, every TCD endpoint matching $\enm n1$ corresponds to the triangulation of an $n$-gon. 

\begin{lemma}
	Every edge-flip in a triangulation $\tg$ corresponds to a 2-2 move at a clockwise bigon in  $\tcd_{\tg}$.
\end{lemma}
\proof{See Figure \ref{fig:tcdpentagon}.\qed }

Thus, as a consequence of Theorem \ref{th:tcdflipsconnected} the (edge-)flip-graph of the triangulations of an $n$-gon is connected as well.

We will also consider some geometry and algebra associated to triangulations in Section \ref{sec:extriangulation}. 

\section{Quad-graphs} \label{sec:quadgraphs}

Let us now turn to the other example that we want to discuss in this section, namely a class of TCDs associated to quad-graphs.

\begin{definition}\label{def:quadgraph}
	A \emph{quad-graph} $\qg$ is a closed disc that is a CW-complex whose faces (2-cells) are quads which are glued edge-to-edge. We will denote the set of vertices (0-cells), edge (1-cells) and faces of a quad-graph $\qg$ by $V,E,F$. 
\end{definition}

By definition every quad-graph $\qg$ is bipartite and we generally assume that we have already chosen a partition of $V$ into white and black vertices. 

\begin{figure}
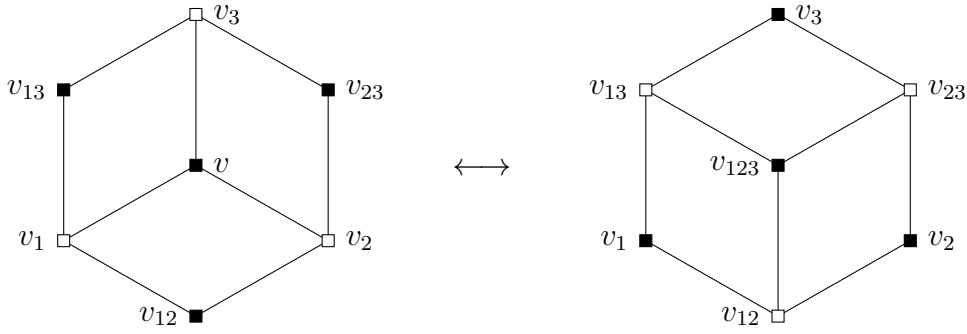


	\caption{The cube flip in a quad-graph.}	
	\label{fig:cubeflipandtcd}	
\end{figure}

\begin{definition}\label{def:cubeflip}
	In a quad-graph $\qg$, consider the local situation where three quads 
	\begin{align}
		q^{12} = (v,v_1,v_{12},v_2),\ q^{23} = (v,v_2,v_{23},v_3) \mbox{ and } q^{13} = (v,v_3,v_{13},v_1)
	\end{align}
	are glued together around a common vertex $v$. Then the \emph{cube flip} at $v$ yields a new quad-graph $\qg'$ which coincides with $\qg$ except that we replace the three quads $q^{12},q^{23},q^{13}$ by the three quads
	\begin{align}
		q^{12}_3 = (v_3,v_{13},v_{123},v_{23}),\ q^{23}_1 = (v_1,v_{12},v_{123},v_{13}) \mbox{ and } q^{13}_2 = (v_2,v_{23},v_{123},v_{12}),
	\end{align}
	where we also replaced the vertex $v$ with the vertex $v_{123}$, see also Figure \ref{fig:cubeflipandtcd}.
\end{definition}

In other words, we flipped the ``back view'' of the cube to the ``front view'' of the cube. Now we associate a TCD to a given quad-graph and show that the cube flip corresponds to a sequence of eleven 2-2 moves.

\begin{samepage}
\begin{definition}\label{ex:tcdquad}\label{def:tcdquad}
	Let $\qg$ be a quad-graph. Into every quad glue the following piece \\
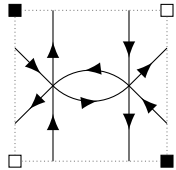
	\begin{center}\vspace{-5mm}

	\end{center}
	of a TCD such that the pieces of strands of neighbouring quads meet along edges. We say the resulting TCD $\tcd_{\qg}$ is the TCD corresponding to $\qg$.
\end{definition}
\end{samepage}

Quad-graphs and the corresponding TCDs will be an important recurring theme throughout this thesis. In particular we can use quad-graphs to describe Q-nets, Darboux maps and line complexes.

\begin{figure}
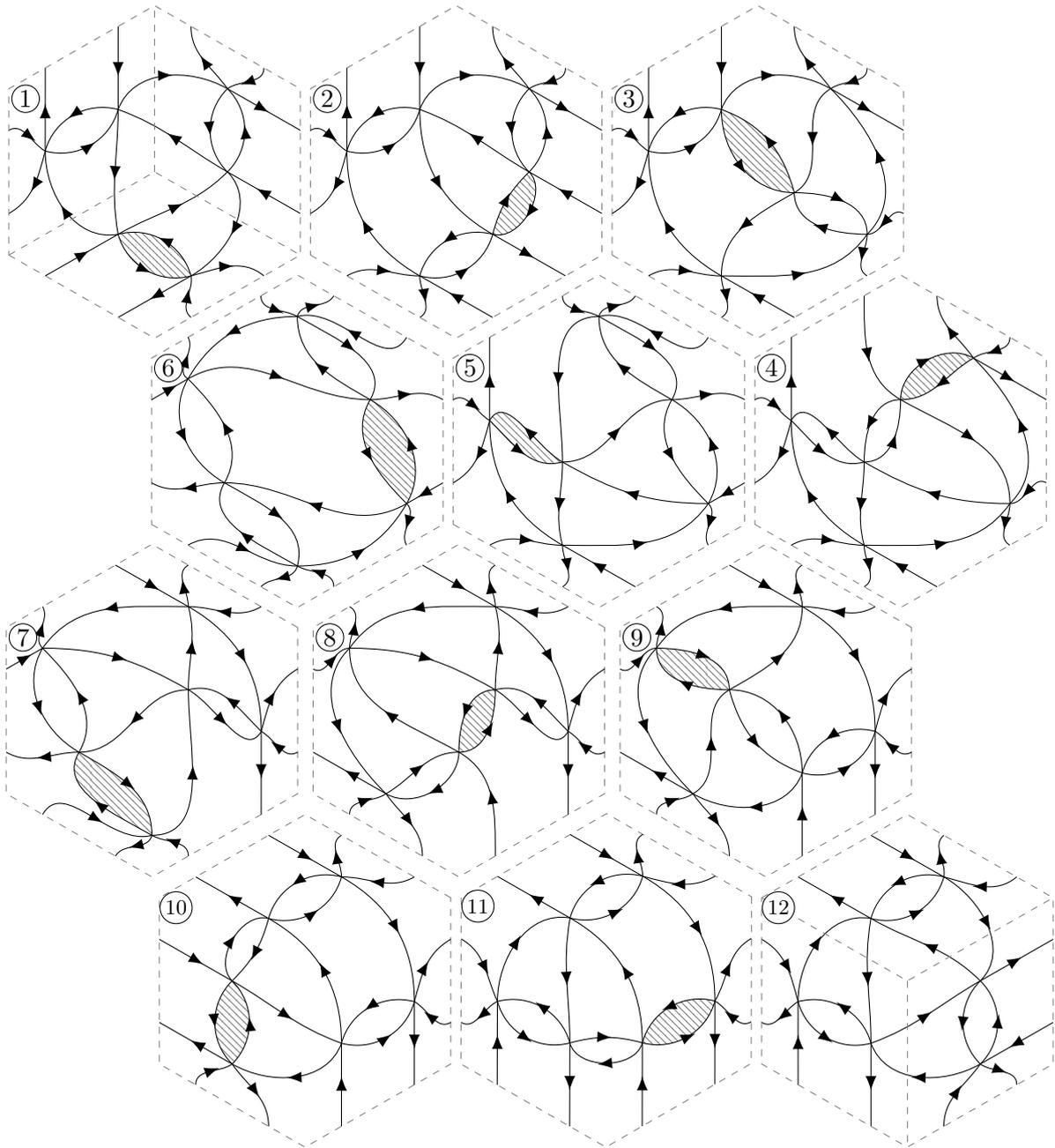

	\hspace{-2.2cm}
	%

	\caption{Eleven 2-2 moves induce the cube flip.}

	\label{fig:tcdcubeflip}
\end{figure}

\begin{lemma}\label{lem:cubefliptcd}
	A cube-flip in a quad-graph $\qg$ corresponds to a sequence of eleven 2-2 moves in $\tcd_{\qg}$.
\end{lemma}
\proof{The sequence is given in Figure \ref{fig:tcdcubeflip}.\qed }

The sequence of 2-2 moves that induces the cube-flip is not unique. Moreover, we could also allow that some of the quads are filled with the TCD of Definition \ref{ex:tcdquad}, but rotated about ninety degrees. In this case there are sequences with less then eleven 2-2 moves.

\begin{figure}[]
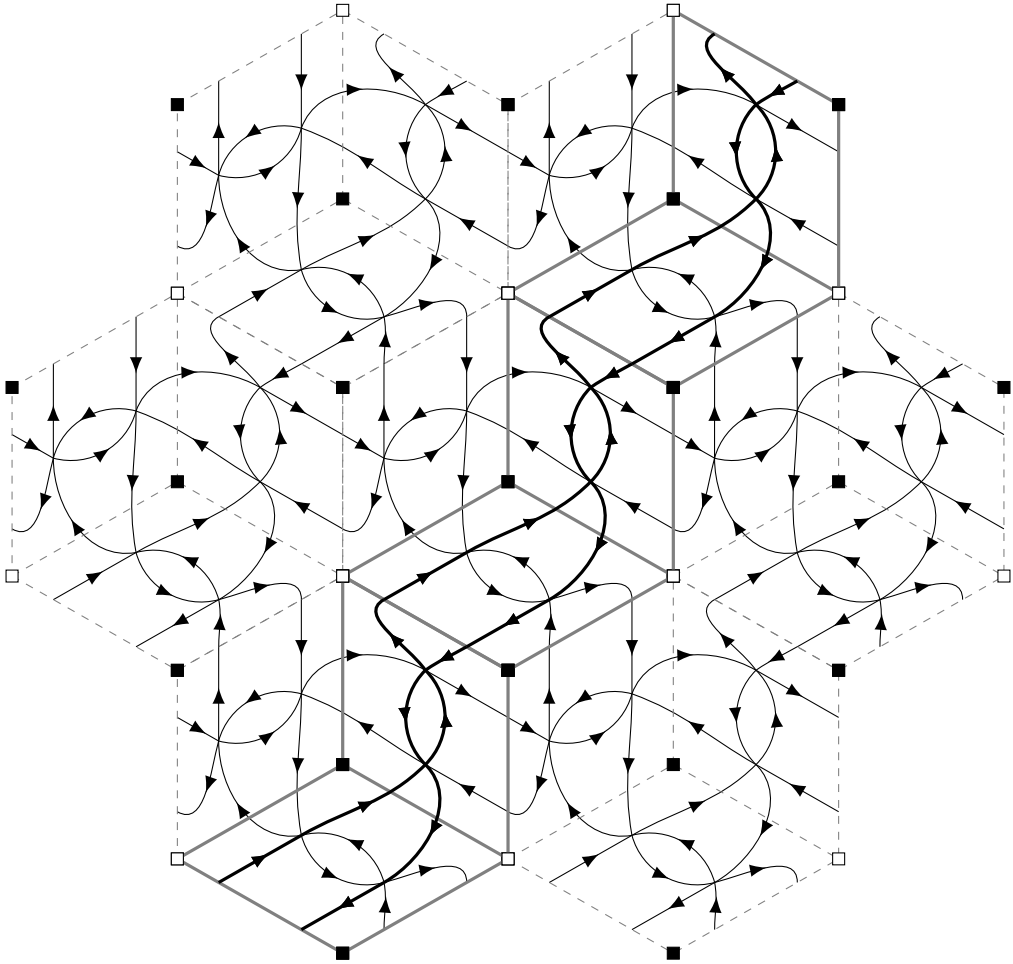


\caption{Highlighting a strip of quads and the corresponding pair of strands in the TCD of a quad-graph.}
\label{fig:quadgraphandtcd}
\end{figure}

\begin{definition}
	A \emph{strip} (of quads) in a quad-graph $\qg$ is a maximal sequence of quads $q_0,q_1,\dots,q_m$ such that for all $k, 0\leq k<m$ the quads $q_k$ and $q_{k+1}$ share an edge $e_k$ and all pairs of consecutive edges $e_k,e_{k+1}$ are disjoint, see Figure \ref{fig:quadgraphandtcd}.
\end{definition}

We also observe in Figure \ref{fig:quadgraphandtcd} that to each strip belong exactly two strands that run in opposite directions. We consider two strips to be intersecting if they share a quad. A strip intersects itself if there are indices $i,j$ with $i\neq j$ but $q_i=q_j$.

\begin{definition}\label{def:minqg}
	A quad-graph is \emph{minimal} if
	\begin{enumerate}
		\item no strip is a loop,
		\item no strip intersects itself,
		\item any two strips intersect at most once.\qedhere
	\end{enumerate}
\end{definition}

This definition is similar to the characterization of minimality in TCDs as in Theorem \ref{th:mintcdforbidden}.

\begin{lemma}\label{lem:qgmintcdmin}
	A quad-graph $\qg$ is minimal if and only if $\tcd_\qg$ is minimal.
\end{lemma}
\proof{If a strip intersects itself then the corresponding two strands intersect, and this is the only way a strand in $\tcd_\qg$ can self intersect. Now assume two different strips intersect twice and pick a strand $s$ that corresponds to the first strip. Then one of the two strands corresponding to the second strip intersects $s$ twice in a parallel manner. If no two strips intersect twice then no two strands can intersect twice, except if the two strands correspond to the same strip. In that case they do not intersect in a parallel manner.\qed
}

Thus, as a consequence of Theorem \ref{th:tcdflipsconnected} the (cube-)flip-graph of minimal quad-graphs with fixed matching of the strip endpoints is connected as well.

There is a procedure to associate a quad-graph to any planar graph. By planar graph $G$ we mean a closed disc that is a CW-complex. Thus $G=(V,E,F)$ where $F$ is the set of faces (2-cells) of $G$.

\begin{definition}\label{def:graphtoqg}
	Let $G$ be a planar graph. The corresponding quad-graph $\qg_G$ has a white vertex $w$ for every vertex $v_w$ of $G$ and a black vertex $b$ for every face $f_b$ of $G$. Moreover, there is an edge in $\qg_G$ between two vertices $w,b$ of $\qg_G$ if $v_w$ is incident to $f_b$ in $G$.
\end{definition}

\clearpage

\section{Vector-relation configurations}\label{sec:vrc}
We continue by defining vector-relation configurations. Throughout the thesis, we will mostly use a subset of vector-relation configurations associated to triple crossing diagrams. However, for some calculations and figures it will prove to be useful to have the general notion of a vector-relation configuration available. 

\begin{definition}[\cite{vrc}]\label{def:vrc}
	Let $G$ be a planar bipartite graph with vertex set $B \mathbin{\dot{\cup}} W$. For $b \in B$ let $N(b) \subset W$ denote its set of neighbors. A \emph{vector-relation configuration (VRC)} $\vrc$ consists of choices of
	\begin{enumerate}
		\item a vector $R_w\in \C^n \setminus\{0\}$ for each $w\in W$ and
		\item an edge weight $\mu_e \in \C \setminus\{0\}$ for each $e \in E$,
	\end{enumerate}
	such that for each $b \in B$ the vectors  $\{R_w : w \in N(b)\}$ satisfy the linear relation 
	\begin{equation}
		R_b:\ \sum_{w\in N(b)} \mu_{bw} R_w = 0.\qedhere
	\end{equation}
\end{definition}

Note that in particular for each black vertex $b\in B(G)$ the vectors $R_{w_1},R_{w_2},\dots R_{w_{d_w}}$ must be linearly dependent, where $w_1,w_2,\dots, w_{d_w} = N(b)$. 

We also consider \emph{gauge transformations} of a VRC $R$. Gauge transformations, come in two kinds:
\begin{enumerate}
	\item Scaling by a factor $\lambda \in \C\setminus\{0\}$ at a white vertex, that is $R_w \mapsto \lambda R_w$ and $\mu_e \mapsto \lambda^{-1} \mu_e$ for all $e \sim w$.
	\item Scaling by a factor $\lambda \in \C\setminus\{0\}$ at a black vertex, that is $R_b \mapsto \lambda R_b$ which corresponds to $\mu_e \mapsto \lambda \mu_e$ for all $e \sim b$.
\end{enumerate}
In the remainder we will mostly consider the vectors of a VRC as the homogeneous lifts of some points in $\CP^n$, thus a scaling at a white vertex only changes the lift but not the point and a scaling at a black vertex does not change the projective subspace in which the adjacent points in $\CP^n$ are contained. Thus the use of the term \emph{gauge} is justified.

It is an interesting question what the invariants of a VRC modulo gauge and projective transformations are. However, we postpone an introduction of these quantities to Chapter~\ref{cha:clusters}, where they will play an integral part in the definition of cluster structures associated to vector-relation configurations.

\begin{figure}
	\begin{tikzpicture}[baseline={([yshift=-.7ex]current bounding box.center)},scale=2.4] 
	\tikzstyle{bvert}=[draw,circle,fill=black,minimum size=5pt,inner sep=0pt]
	\tikzstyle{wvert}=[draw,circle,fill=white,minimum size=5pt,inner sep=0pt]
	\tikzstyle{nvert}=[]
	
	\node[wvert] (v1)  at (-.87,-.5) {};
	\node[wvert] (v2) at (.87,-.5) {};
	\node[wvert] (v3) at (0,1) {};
	\node[wvert] (v12)  at (0,-1) {};
	\node[wvert] (v23) at (.87,.5) {};
	\node[wvert] (v13) at (-.87,.5) {};
	\node[wvert] (v) at (0,0) {};
	
	\node[bvert] (f12) at ($(v)!.5!(v12)$) {};
	\node[bvert] (f23) at ($(v)!.5!(v23)$) {};
	\node[bvert] (f13) at ($(v)!.5!(v13)$) {};
	
	\coordinate (e1) at ($(v)!.5!(v1)$) {};
	\coordinate (e2) at ($(v)!.5!(v2)$) {};
	\coordinate (e3) at ($(v)!.5!(v3)$) {};	
	\coordinate (e12) at ($(v2)!.5!(v12)$) {};
	\coordinate (e13) at ($(v3)!.5!(v13)$) {};
	\coordinate (e21) at ($(v1)!.5!(v12)$) {};
	\coordinate (e23) at ($(v3)!.5!(v23)$) {};
	\coordinate (e31) at ($(v1)!.5!(v13)$) {};
	\coordinate (e32) at ($(v2)!.5!(v23)$) {};
	
	\draw[-] (v1) -- (f12) -- (v2) -- (f23) -- (v3) -- (f13) -- (v1);
	\draw[-] (v12) -- (f12) -- (v) -- (f23) -- (v23) (v) -- (f13) -- (v13);
	\draw[-,gray,dashed] (v1) -- (v) -- (v2) (v) -- (v3);
	\draw[-,gray,dashed] (v1) -- (v12) -- (v2) -- (v23) -- (v3) -- (v13) -- (v1);
	\end{tikzpicture}
	\caption{The graph $G$ of a VRC of a Q-net defined on three quads (dashed).}
	\label{fig:vrcqnet}
\end{figure}
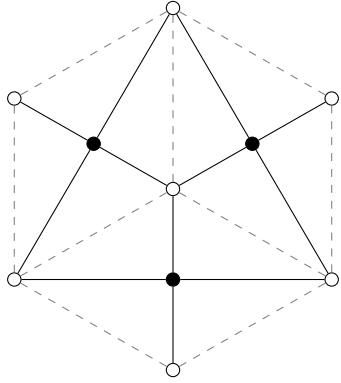

\begin{example}\label{ex:vrcquad}
	A \emph{Q-net} is a map from the vertices of a quad-graph $\qg$ (see Definition \ref{def:quadgraph}) to $\CP^n$ such that the image of each quad is contained in a 2-plane. Since $\qg$ is planar and bipartite, any homogeneous lift of a Q-net to $\C^{n+1}$ can be represented as a VRC. See Figure \ref{fig:vrcqnet} for an example of the respective graph $G$.
\end{example}

There is one particular useful gauge, both for calculations in general and for the affine cluster structures that we define in Section \ref{sec:affcluster}.

\begin{definition}\label{def:affinegauge}
	Consider a VRC $\vrc$. The edge-weights $\mu$ of $\vrc$ are in \emph{affine gauge} if 
	\begin{align}
		\sum_{w \sim b} \mu_{bw} = 0
	\end{align}
	holds for every black vertex $b\in B$.
\end{definition}
Let us explain why we call this an affine gauge. Assume each vector $R_w \in \C^{n+1}$ of $\vrc$ is the homogeneous lift of a point $T_w$ in $\CP^{n}$, and assume that none of the vectors have a zero in their last component. Then we can apply gauge transformations at each white vertex $w$ of $G$ such that the last component, that is the $(n+1)$-th component of $R_w$ is 1. In this gauge the edge-weights have to satisfy the condition of Definition \ref{def:affinegauge}. Moreover, in this gauge the $n$-tuple of coordinates $1,\dots,n$ are called \emph{affine coordinates} of $T_w$ in textbook projective geometry.

Assume we have a VRC in an affine gauge. Then at a black vertex of degree three, we can always perform a gauge transform at the black vertex such that the three incident edge weights are $1,\mu$ and $-(1+\mu)$. Later on, we will mostly work in setups in which all black vertices are of degree three, so this is a useful trick for calculations.

Besides gauge transformations we are also interested in a set of combinatorial transformations.

\begin{figure}
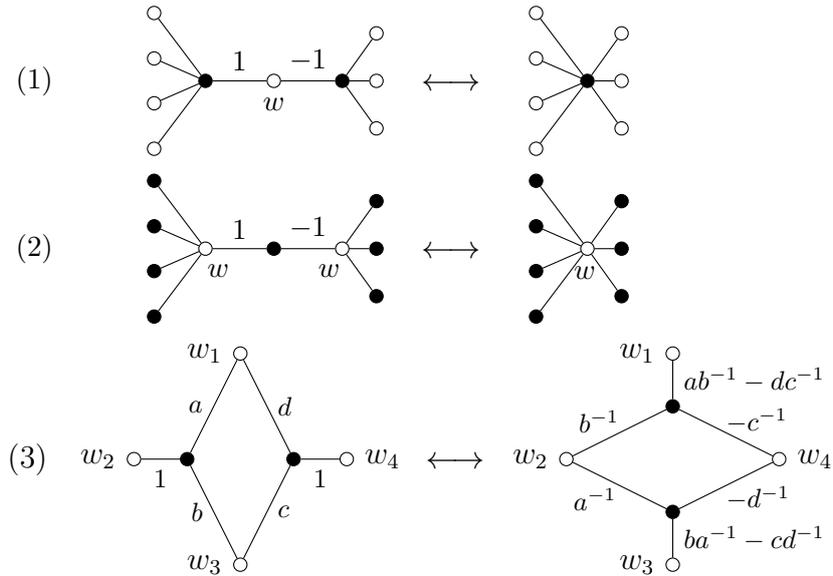

\vspace{-2mm}
	\caption{Local moves: (1) contraction at  a white vertex, (2) contraction at a black vertex and (3) the spider move.} 
	\label{fig:vrclocalmoves}
\end{figure}

\begin{definition}\label{def:vrclocaltrafos}
	We define three \emph{local transformations} that locally alter a VRC $\vrc$. The three local transformations (and their inverses) are called
	\begin{enumerate}
		\item the \emph{contraction} (resp.\ \emph{split}) of a white (black) vertex,
		\item the \emph{contraction} (resp.\ \emph{split}) of a black (white) vertex and
		\item the \emph{spider move}.
	\end{enumerate}
	The changes of combinatorics and edge-weights of these transformations are shown in Figure \ref{fig:vrclocalmoves}. The edges that carry no label in Figure \ref{fig:vrclocalmoves} do not change their weight in the respective move.
\end{definition}
It is a straightforward calculation to check that the vector $R_w$ in the split is indeed well-defined from the edge weights. In the case of the spider move another quick calculation shows that the relations on the left imply the relations on the right and vice versa. Note that if the combinatorics are as in Figure \ref{fig:vrclocalmoves}, then we can always apply gauge transformations such that the edge weights are as on the left of Figure \ref{fig:vrclocalmoves} as well. Therefore whether a local transformations is possible or not depends purely on the combinatorics, not on the vectors. Moreover, note that the edge weights for the spider move as shown in Figure \ref{fig:vrclocalmoves} are not symmetric with respect to a rotation by 180 degree, but they are symmetric up to gauge.

Let us also say a few words about the geometric meaning of the local transformations. Assume each vector $R_w \in \C^{n+1}$ is the homogeneous lift of a point in $T_w \in \CP^n$. We consider the moves:
\begin{enumerate}
	\item On the left in Figure \ref{fig:vrclocalmoves} there are two black vertices and we label their neighbours $\{w_1,\dots,w_k,w\}$ and $\{w'_1,\dots,w'_l,w\}$. Assume that $n \geq l+k-2$ and that all the points are in general position. Then we can remove $w$ and join the two black vertices, because $\{T_{w_1},\dots,T_{w_k},T_{w'_1},\dots,T_{w'_l}\}$ span an $(l+k-2)$-dimensional projective space. Vice versa, given a black vertex with neighbours $\{w_1,\dots,w_k,w'_1,\dots,w'_l\}$ we can introduce the point $T_w = \spa \{T_{w_1},\dots,T_{w_k}\} \cap \spa \{T_{w'_1},\dots,T_{w'_l}\}$.
	\item If we have a black vertex incident to two white vertices, then the corresponding two projective points have to be identical. We can thus contract the black vertex and replace the two white vertices with one, also corresponding to the same point. We can also split any white vertex by adding a two valent black vertex.
	\item Consider two black vertices with neighbours $\{w_1,w_2,w_4\}$ and $\{w_2,w_3,w_4\}$. Then the projective images $T_{w_1}, T_{w_2}, T_{w_3}, T_{w_4}$ have to be on a common projective line. This configuration can clearly also be captured by two black vertices with neighbours $\{w_1,w_2,w_3\}$ and $\{w_3,w_4,w_1\}$. 
\end{enumerate}
There are other known local transformations at a quad that we did not list here. However, they can all be obtained by a combination of the three local transformations above.

\section{TCD maps} \label{sec:tcdmap}

The goal of this section is to introduce geometric maps associated to triple crossing diagrams. They encode the linear relations of point configurations, and allow us to treat a large number of examples of maps from discrete differential geometry, discrete integrable systems and statistical mechanics in a unified framework. They will also allow us in Chapter \ref{cha:clusters} to introduce two canonically associated cluster structures, a projective and an affine structure. The canonical existence of these two structures explains the occurrence of the two types of cluster structures associated to examples in the literature.  Moreover, the combinatorial framework of TCD maps is flexible enough to allow us to study the meaning of the projection, section and projective dual of a TCD map in terms of geometry and combinatorics. Additionally, we obtain an understanding of the algebraic and combinatorial identities of the cluster structures under these operations.

\begin{figure}
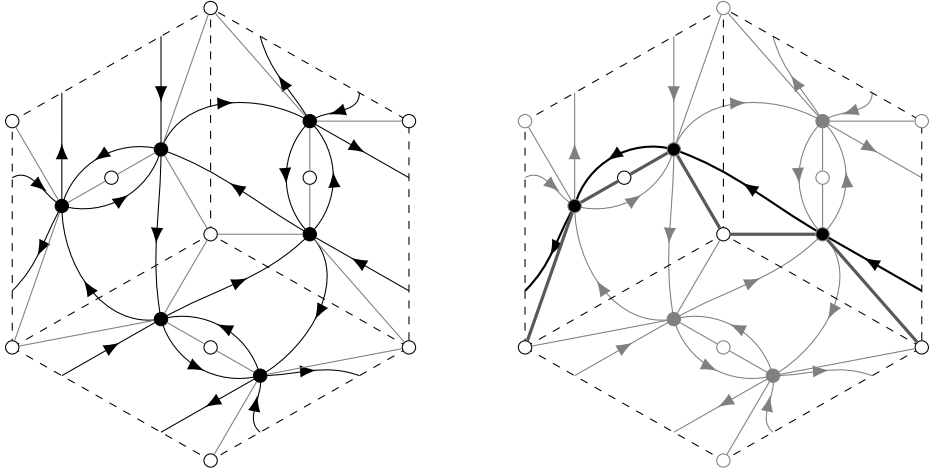

	
	\caption{Left: A graph $\pb$ (gray) of three planar quads (dashed) and the associated TCD (black). Right: We highlighted one strand of $\tcd$ as well as the corresponding zig-zag path in $\pb$.}
	\label{fig:vrcvstcd}
\end{figure}

\begin{definition}
	Let $\tcd$ be a TCD. Then let $\tcdp$ denote the set of faces, including the boundary faces, of $\tcd$ that are oriented counterclockwise and let $\tcdm$ denote the set of faces, including the boundary faces, of $\tcd$ that are oriented clockwise.
\end{definition}


\begin{definition}
	A \emph{TCD map $T: \tcdp \rightarrow \CP^n$} assigns a point in $\CP^n$ to each counterclockwise oriented face of $\tcd$, such that at every triple crossing any two of the three incident points span a line that contains the third point. To every TCD map we can choose an \emph{associated vector-relation configuration (VRC) $\vrc$}, which has a white vertex for every counterclockwise face of $\tcd$ that is mapped to a lift in $\C^{n+1}\setminus\{0\}$ of the corresponding point of $T$, as well as a black vertex for every triple crossing point. There is an edge from a black to a white vertex in $\vrc$ if the corresponding crossing point and face are incident in $\tcd$. A \emph{boundary vertex} is a white vertex that corresponds to a counterclockwise boundary face of $\tcd$.
\end{definition}

In a sense, TCD maps capture the most fundamental invariant among points in projective spaces, namely that three points are on a common line. To some extent, this explains the usefulness of considering triple crossings. Why the planar structure is so beneficial is less obvious. One reason is that throughout this thesis we consider 3D-systems (or reductions thereof), which are systems that are defined by two dimensional Cauchy data. Thus it is less surprising that the essential data is associated to planar combinatorics.

\begin{definition}\label{def:pb}
	Let $\tcd$ be a TCD. The \emph{associated bipartite planar graph $\pb$} has
	\begin{enumerate}
		\item a white vertex $w$ for each counterclockwise face (including boundary faces) of $\tcd$,
		\item a black vertex $b$ for each triple crossing of $\tcd$,
		\item an edge $e=(w,b)$ for each pair of counterclockwise face and adjacent triple crossing.\qedhere
	\end{enumerate}
\end{definition}

Whenever we work with a TCD $\tcd$, we assume that $\pb$ is the associated bipartite planar graph. To distinguish, we generally denote other graphs by $G$. Note that by definition the set of white vertices $W$ of $\pb$ is in bijection with $\tcdp$. Also by definition, every black vertex of $\pb$ has degree three. In fact, it is easy to see that any bipartite planar graph $G$ with only white boundary vertices and only black vertices of degree three defines a triple crossing diagram.

\begin{definition}
	Let $T: \tcdp \rightarrow \CP^n$ be a TCD map. An \emph{associated vector-relation configuration} $\vrc$ is a VRC $\vrc:\pb \rightarrow \C^{n+1}$, such that $R(w)$ is a homogeneous lift of $T(w)$ for every $w\in W$.
\end{definition}
Whenever we work with a TCD map $T: \tcdp \rightarrow \CP^n$ we assume that we have also fixed an associated VRC $\vrc:\pb \rightarrow \C^{n+1}$. Of course, this means that we have to take care whether statements and proofs are independent of the choice of homogeneous lifts in $\vrc$. In practice however, this does not present any difficulty. 

To any VRC $\vrc: W(G)\rightarrow \C^n$, we can also associate a TCD map $T: \tcdp \rightarrow \CP^{n-1}$ for some TCD $\tcd$. Too see this, note that with the local contraction and split operations explained in Definition \ref{def:vrclocaltrafos}, we can transform any bipartite planar graph $G$ into a graph $\pb$, such that all its black vertices have degree three. In the process, new white vertices and thus new vectors may appear, but up to gauge, these are uniquely determined by the edge-weights of $\vrc$. Finally, the corresponding TCD map $T$ simply maps each white vertex of $\pb$ to the projectivization of the vector $R(w)$. Due to the process of splitting black vertices the associated TCD map for a given VRC is not unique. However, two different TCD maps associated to the same VRC are always related by a sequence of local moves. We explain local moves below in Definition \ref{def:tcdmapmoves}.

We consider the TCD maps to be the fundamental objects instead of VRCs, as we are primarily interested in objects and maps in projective geometry. Moreover, most of the interesting quantities we study later on, like the cluster structures and partition functions are considered modulo gauge transformations of VRCs, which is another reason why it is natural to consider TCD maps as the fundamental objects. Another advantage of TCD maps is that in the associated VRC, all the edge-weights are determined from the vectors up to gauge transformations, unlike in general VRCs. Thus all the information about a TCD map is completely contained in the points of the TCD map and no information is hidden in the edge-weights.

In general, we attempt to rely as little as possible on calculations and prefer to rely on incidence theorems and a flexible collection of lemmas. However, on the occasions where we do need to use calculations, the use of VRCs is most practical.

\begin{example}
	We looked at the VRC of planar quads before, see Example \ref{ex:vrcquad} and Figure \ref{fig:vrcqnet}. In that case the planarity of each quad was represented by a degree four black vertex. We obtain a triple crossing diagram if we split each of those black vertices, see Figure~\ref{fig:vrcvstcd}. The extra points that we created this way are the intersection points of opposite lines of a planar quad. These extra points are called the \emph{focal points}, as we will explain in more detail in Section~\ref{sec:qnets}.
\end{example}

Let us also give a short explanation of how to recognize the strands of a TCD in the associated bipartite graph $\pb$. In fact, a generalization of these strands exists for arbitrary planar bipartite graphs.

\begin{definition}\label{def:zigzagtcd}
	Let $G$ be a planar bipartite graph. A \emph{zig-zag} path in $G$ is a path $(w_1,b_1,w_2,b_2,\dots,b_{n-1},w_n)$ that 
	\begin{enumerate}
		\item turns maximally left at black vertices,
		\item turns maximally right at white vertices,
		\item begins and ends at the boundary.\qedhere
	\end{enumerate}
\end{definition}

Clearly, if $G$ is actually the associated bipartite graph $\pb$ of a TCD $\tcd$, then the zig-zag paths of $\pb$ are in bijection with the strands of $\tcd$, see Figure \ref{fig:vrcvstcd}. Given a strand, the zig-zag path consists of all the black vertices that are on the strand and all white vertices that are just to the left of the strand.

\section{Local moves in TCD maps and the dSKP equation}\label{sec:tcddskp}

In Section \ref{sec:tcds} we looked at 2-2 moves, the local moves in triple crossing diagrams. Let us now define 2-2 moves not only on the TCD, but also on the TCD map and investigate the consequences for the associated VRC.

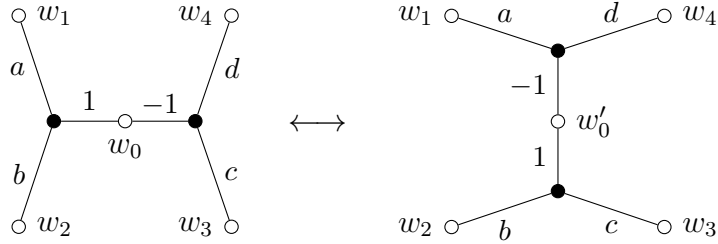
\begin{figure}
	\centering
	\begin{tikzpicture}[scale=1.4,baseline={([yshift=-.7ex]current bounding box.center)}] 
	\node[wvert,label=right:$w_1$] (w1)  at (0,0) {};
	\node[wvert,label=right:$w_2$] (w2) at (0,-2) {};
	\node[bvert] (b12)  at (.33,-1) {};
	\node[wvert,label=below:$w_0$] (x) at (1,-1) {};
	\node[bvert] (b34) at (1.66,-1) {};
	\node[wvert,label=left:$w_3$] (w3) at (2,-2) {};
	\node[wvert,label=left:$w_4$] (w4) at (2, 0) {};
	
	\path[-,font=\small] 
		(b12) edge node[left] {$a$} (w1)
		(b12) edge node[left] {$b$} (w2)
		(b12) edge node[above] {$1$} (x)
		(b34) edge node[above,inner sep=2] {$-1$} (x)
		(b34) edge node[right] {$c$} (w3)
		(b34) edge node[right] {$d$} (w4)
	;	
	\end{tikzpicture}\hspace{.5cm}$\longleftrightarrow$\hspace{.4cm}
	\begin{tikzpicture}[scale=1.4,baseline={([yshift=-.7ex]current bounding box.center)}] 
	\node[wvert,label=left:$w_1$] (w1)  at (0,0) {};
	\node[wvert,label=left:$w_2$] (w2) at (0,-2) {};
	\node[bvert] (b14)  at (1,-.33) {};
	\node[wvert,label=right:$w_0'$] (y) at (1,-1) {};
	\node[bvert] (b23) at (1,-1.66) {};
	\node[wvert,label=right:$w_3$] (w3) at (2,-2) {};
	\node[wvert,label=right:$w_4$] (w4) at (2, 0) {};
	
	\path[-,font=\small]
		(b23) edge node[below] {$b$} (w2)
		(b23) edge node[below] {$c$} (w3)
		(b23) edge node[left] {$1$} (y)
		(y) edge node[left] {$-1$} (b14)
		(b14) edge node[above] {$a$} (w1)
		(b14) edge node[above] {$d$} (w4)
	;	
	\end{tikzpicture}
	\vspace{-1mm}
	\caption{The edge-weights of the resplit.}
	\label{fig:resplit}
\end{figure}

\begin{samepage}
\begin{definition}\label{def:tcdmapmoves}
	Let $T: \tcdp \rightarrow \CP^n$ be a TCD map. A \emph{local move} replaces $T: \tcdp \rightarrow \CP^n$ with a new TCD map $\tilde T: \tilde\tcdp \rightarrow \CP^n$. The two local moves are
	\begin{enumerate}
		\item the \emph{spider move}, which is a 2-2 move at a clockwise oriented face of $\tcd$. The image points of $T$ remain unchanged, the move corresponds to the spider move in $\pb$, see $(3)$ in Figure~\ref{fig:vrclocalmoves};
		\item the \emph{resplit}, which is a 2-2 move at a counterclockwise oriented face $f$ of $\tcd$. The face $f$ corresponds to a white vertex $w_0$ of $\pb$. The point $T(w_0)$ is replaced with $\tilde T(w_0')$ by the resplit and is determined by projectivization of the vector $\tilde R(w_0')$, which in turn is determined by the edge weights of $\vrc$ as given in Figure~\ref{fig:resplit}. The other image points of $T$ remain unchanged.\qedhere
	\end{enumerate}
\end{definition}
\end{samepage}

\begin{figure}
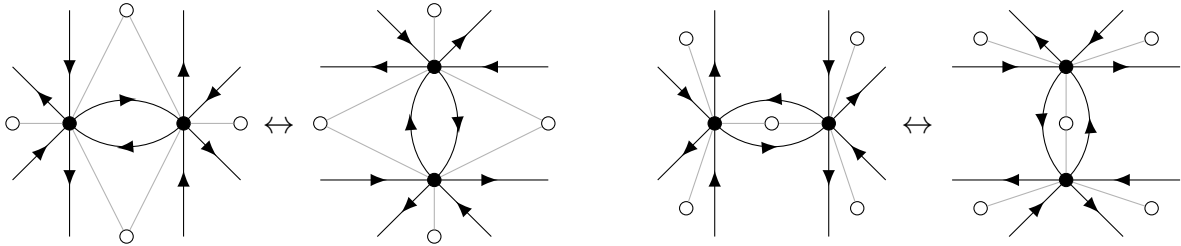

	\centering
	
	\caption{The TCD $\tcd$ and graph $\pb$ as they appear in the 2-2 move at a clockwise and a counterclockwise oriented bigon.}\label{fig:tcdtwotwo}
\end{figure}

Let us discuss the resplit first. By definition, the resplit replaces $T(w_0)$ with $\tilde T(w_0')$. We need to check that $\tilde T$ is still a TCD map, which corresponds to verifying that $T(w_0')$ does not coincide with any of the four unchanged points $T(w_1),T(w_2),T(w_3),T(w_4)$. However, as none of the edge-weights after the resplit are equal to 0, this cannot happen and $\tilde T$ is indeed a TCD map.

\begin{figure}
	\begin{tikzpicture}
		\node[wvert,label=right:$T(w_1)$] (a1) at (1.0,3.7) {};
		\node[wvert,label=right:$T(w_0)$] (a2) at (-0.6,5.4) {};
		\node[wvert,label=above:$T(w_4)$] (a3) at (-6.6,2.5) {};
		\node[wvert,label=above:$T(w_3)$] (b1) at ($(a3)!.47!(a2)$) {};
		\node[wvert,label=above:$\tilde T(w_0')$] (b2) at ($(a3)!.68!(a1)$) {};
		\node[wvert,label=right:$T(w_2)$] (b3) at (intersection of a1--a2 and b1--b2) {};
		
		\draw[-]
			(a1) -- (a2) -- (b1) -- (a3) -- (b2) -- (a1) -- (b3) -- (b2) -- (b1)
		;
	\end{tikzpicture}\vspace{-2mm}
	\caption{Menelaus' configuration.}
	\label{fig:menelaus}
\end{figure}
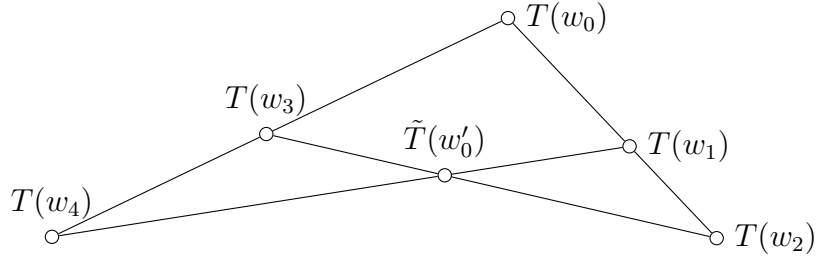

In the generic situation, the five points $T(w_0),T(w_1),T(w_2),T(w_3),T(w_4)$ involved in a resplit span a plane. In that case, we actually find that we are in the situation of Menelaus' configuration, see Figure \ref{fig:menelaus}. Indeed, $\tilde T(w'_0)$ has to be the intersection of the lines $T(w_1)T(w_4)$ and $T(w_2)T(w_3)$. Thus in the generic situation the new point $\tilde T(w'_0)$ after the resplit is uniquely defined by incidence geometry. We can also say more about the degenerate case, the case where $T(w_0),T(w_1),T(w_2),T(w_3),T(w_4)$ span a line instead of a plane, but we need to understand the concept of cross-ratios and multi-ratios first.

\begin{definition}\label{def:orientedlengthratio}
	Let $p_1,p_2,p_3$ be points on a line $\ell$ in $\CP^n$ and assume we have a fixed affine chart $\C^n$ of $\CP^n$. Then the \emph{oriented length ratio} $\lambda(p_j,p_k,p_i)$ satisfies
	\begin{align}
		(p_j-p_i) = (p_k-p_i)\lambda(p_j,p_k,p_i),
	\end{align}
	for any set of indices such that $\{i,j,k\} = \{1,2,3\}$. We also observe that the fixed chart $\C^n$ induces an affine chart $\C$ of the line $\ell$. Thus, we also write the oriented length ratio as
	\begin{align}
		\lambda(p_j,p_k,p_i) = \frac{p_j-p_i}{p_k-p_i},
	\end{align}
	where the quotient is taken in the induced chart $\C$ of $\ell$. 
\end{definition}

\begin{definition}\label{def:multiratio}
	Let $p_1,p_{1,2},p_2,p_{2,3},\dots, p_m,p_{m,1}$ be $2m$ points in $\CP^n$ such that every $p_{k,k+1}$ is on the line $p_kp_{k+1}$. Then we choose an affine chart $\C^n$ of $\CP^n$ and define the \emph{multi-ratio} of these points as
	\begin{align}
		\mr(p_1,p_{1,2},p_2,\dots, p_m,p_{m,1}) = \prod_{k=1}^m \frac{p_k - p_{k,k+1}}{p_{k,k+1} - p_{k+1}}.
	\end{align}
	In the special case of $m=2$ we call the multi-ratio the \emph{cross-ratio}, that is
	\begin{equation}
		\cro(p_1,p_{1,2},p_2,p_{2,1}) = \frac{p_1 - p_{1,2}}{p_{1,2} - p_{2}}\frac{p_2 - p_{2,1}}{p_{2,1} - p_{1}}. \qedhere
	\end{equation}
\end{definition}

Let us gather a few basic but important properties of the multi-ratio, see for example \cite{ddgbook} for proofs.
\begin{lemma}\label{lem:projinvariants}
	Let $p_1,p_{1,2},p_2,p_{2,3},\dots, p_m,p_{m,1}$ be $2m$ points in $\CP^n$ such that every $p_{k,k+1}$ is on the line $p_kp_{k+1}$. Then the multi-ratio $\mr(p_1,p_{1,2},p_2,\dots, p_m,p_{m,1})$
	\begin{enumerate}
		\item does not depend on the affine chart used to calculate it,
		\item is invariant under projective transformations and
		\item is invariant under (central) projections, if there is no $k$ such that $p_{k,k+1}$ is projected to the same point as $p_k$ or $p_{k+1}$.\qedhere
	\end{enumerate}
\end{lemma}

Note that in the thesis we adhere to the convention that \emph{projective transformations} are bijections $\CP^n\rightarrow \CP^n$. We use the term \emph{projection} short for \emph{central projection} and note that with our convention a projection is not a projective transformation.

We now show that oriented length ratios appear naturally in TCD maps.

\begin{lemma}\label{lem:edgeratios}
	Assume we have a degree three black vertex $b$ in a VRC $R$ with white neighbours $w_1,w_2,w_3$. Let $\mu_1,\mu_2,\mu_3$ be the corresponding edge-weights in the affine gauge. Then the oriented length ratio is
	\begin{align}
		\frac{R(w_j)-R(w_i)}{R(w_k)-R(w_i)} = -\frac{\mu_k}{\mu_j},
	\end{align}
	for $\{i,j,k\} = \{1,2,3\}$.
\end{lemma}
\proof{
	Since we are in the affine gauge of a projective line, we can take $R(w_1),R(w_2),R(w_3)$ to be numbers in $\C$ accompanied by the two equations
	\begin{align}
		\mu_1 R(w_1) + \mu_2 R(w_2) + \mu_3 R(w_3) &= 0,\\
		\mu_1 + \mu_2 + \mu_3  &= 0.
	\end{align}
	We insert these equations into the desired expression to obtain
	\begin{equation}
		\frac{p_j-p_i}{p_k-p_i} = \frac{\frac{\mu_i p_i+ \mu_k p_k}{\mu_i + \mu_k} - p_i}{p_k-p_i} =  \frac{\frac{-\mu_k p_i+ \mu_k p_k}{\mu_i + \mu_k}}{p_k-p_i} = -\frac{\mu_k}{\mu_j}.\qedhere
	\end{equation}
}

We are now ready to characterize the resplit in the non-generic situation via a multi-ratio equation.

\begin{lemma}\label{lem:resplitmr}
	Let $T(w_0),\tilde T(w_0'),T(w_1),T(w_2),T(w_3),T(w_4)$ be the points involved in a resplit as in Figure \ref{fig:resplit}. Then the multi-ratio equation
	\begin{align}
		\mr(T(w_1),T(w_0),T(w_2),T(w_3),\tilde T(w_0'),T(w_4)) = -1 \label{eq:resplitmr}
	\end{align}
	holds.	
\end{lemma}
\proof{
	Due to the Lemma \ref{lem:edgeratios} the multi-ratio is expressible via the edge weights. The edge weights can be read off Figure~\ref{fig:resplit} and cancel out, which yields the desired equation from the statement of the lemma.\qed
}

The multi-ratio equation from Lemma \ref{lem:resplitmr} is called the \emph{dSKP equation}, short for \emph{discrete Schwarzian Kadomtsev-Petviashvili equation} \cite{ncwqdskp, dndskp, bkkphierarchy, bkkphierarchymulti, ksclifford}. The dSKP equation is considered as a lattice equation on the $A_3$ lattice. The relation to $A_k$-lattices will be explained in Section \ref{sec:tcdconsistency}. The case of the $A_3$ lattice appears again in Section \ref{sec:laplacedarboux}. Note that Equation \eqref{eq:resplitmr} possesses a high degree of symmetry, according to the symmetries of the Menelaus' configuration (see Figure \ref{fig:menelaus}). Indeed, in the generic case, that is if the points $T(w_0),\tilde T(w_0'),T(w_1),T(w_2),T(w_3),T(w_4)$ span a plane, then Lemma \ref{lem:resplitmr} is the following well known (and ancient) incidence theorem, named after Menelaus.

\begin{theorem}[Menelaus' theorem]\label{th:menelaus}
	Let $A_1,A_2,A_3$ be the vertices of a triangle in $\CP^2$ and let $B_{12},B_{23},B_{13}$ be three points such that $B_{12}$ is on the line $A_1A_2$, $B_{23}$ is on the line $A_2A_3$ and $B_{13}$ is on the line $A_1A_3$, see Figure \ref{fig:menelaus}. Then
	\begin{align}
		\mr(A_1,B_{12},A_2,B_{23},A_3,B_{13}) = -1
	\end{align}
	if and only if $B_{12},B_{23},B_{13}$ are on a line.
\end{theorem}

The relation between Menelaus' theorem and configuration and the dSKP equation has already been studied by Ko\-no\-pel\-chen\-ko and Schief \cite{ksclifford}. Let us give a short summary. As consequence of Lemma \ref{lem:projinvariants} the projection of 6 points related by a resplit in $\CP^2$ is also a resplit in $\CP^1$, because the Lemma states that multi-ratios are invariant under projections. The converse is also true: Every resplit in $\CP^1$ is the projection of a resplit in $\CP^2$. To see this, assume we include $\CP^1$ as a line $\ell$ in $\CP^2$ and fix a point $P \in \CP^2 \setminus \ell$ from which we want to project onto $\ell$. Then we can choose the lifts of $T(w_0),T(w_1)$ resp. $T(w_4)$ on the lines $PT(w_0),PT(w_1)$ resp. $PT(w_4)$ such that the lifts are not on a line. Now, $T(w_2)$ and $T(w_3)$ are determined as intersections of lines. Moreover, the lift of $\tilde T(w_0')$ has to be the intersection point of the lines $T(w_2) T(w_3)$ and $T(w_1)T(w_4)$. This intersection point is exactly the point that is projected onto $\tilde T(w_0')$ because of the multi-ratio characterization of the resplit on a line and the invariance of the multi-ratio under projections. As a result, the construction of the lift always closes and the lift always exists.

Let us now also discuss the spider move. We observe that the spider move does not change the points of $T$ nor do new lines occur, only the combinatorics of the triple crossing diagram $\tcd$ change. This is why we think of the spider move as a \emph{reparametrization}. The geometric object itself does not change, only the parametrization of the object. However, there is a problem that can occur with the spider move. Let us denote the four white vertices involved in the spider move by $w_1,w_2,w_3,w_4$ such that the triplets $w_1,w_2,w_3$ and $w_3,w_4,w_1$ are adjacent to a common black vertex before the move and the triplets $w_2,w_3,w_4$ and $w_4,w_1,w_2$ are adjacent to a common black vertex after the move (see Figure \ref{fig:vrclocalmoves}). Then it is possible that $T(w_2) = T(w_4)$, which is not a problem for $T$. However, then we have also $\tilde T(w_2) = \tilde T(w_4)$, which is a violation of the definition of a TCD map for $\tilde T$. In these cases the spider-move is not well defined. However, in some situations, we want to assume that we can perform 2-2 moves without second thoughts.

\begin{definition}\label{def:gentcdmap}
	Let $\tcd$ be a minimal TCD and $T:\tcd \rightarrow \CP^n$ be a TCD map. We call $T$ a \emph{flip-generic} TCD map if every possible spider move is well-defined and if for any TCD map $\tilde T$ related to $T$ via a sequence of 2-2 moves, also every possible spider move is well-defined.
\end{definition}
It is important in Definition \ref{def:gentcdmap} that $\tcd$ is minimal, because in this case the flip graph is finite and well understood, as discussed in Section \ref{sec:tcds}.

\begin{remark}
	One of the main motivations of TCD maps is to study examples of discrete differential geometry (DDG). In DDG it is common not to investigate the occurrence of singularities, but simply assumes that these do not occur. It is usually plausible enough that one can simply ``wiggle'' a bit in the initial data to avoid running into singularities, and this approach is perfectly suitable to the needs of DDG. However, to some extent in this thesis we want to show that many of the singularities that are avoided in DDG are actually not a problem and sometimes even of specific interest. Section \ref{sec:nongendimension} is devoted to TCD maps that feature non-generic projective dimension. Therefore it is useful to have Definition \ref{def:gentcdmap} that singles out those singularities that we do also avoid in the TCD map framework.
\end{remark}

\section{Existence and maximal dimension of TCD maps}

One natural question that we have not answered yet is whether for a given TCD there actually is a TCD map. Of course, if the TCD map maps to $\CP^1$ then (almost) any choice of points is allowed as any three points are on a line in $\CP^1$. Thus the less trivial question is: Given $\tcd$, what is the \emph{maximal dimension} $n$ such that there is a TCD map $T:\tcd \rightarrow \CP^n$ and such that the image of $T$ spans $\CP^n$?

\begin{figure}\hspace{-0.5cm}
\vspace{-2mm}
\caption{A standard diagram and the associated graph $\pb$ for endpoint matching $\enm 63$.}
\label{fig:standardtcdmap}
\end{figure}

\begin{theorem}\label{th:maxdim}
	The maximal dimension of a minimal TCD $\tcd$ is $|W|-|B| - 1$, the difference of white and black vertices minus one in the graph $\pb$. Equivalently, the maximal dimension is the number of left moving strands in any half-plane drawing. Moreover, all points of a TCD map are included in the span of the boundary points.
\end{theorem}

\proof{We prove this for standard diagrams first. Consider the iterative construction of a standard diagram as in Definition \ref{def:standardtcd}. We assign points to $T$ as we proceed backwards through the steps. As we start our labeling of the strands with an in-endpoint on the lower left the outer face (above everything else, see Figure \ref{fig:standardtcdmap}) is counterclockwise. Thus before we begin, we add a projective point to the TCD map for the outer face. The projective dimension spanned by the points of the TCD map is therefore zero so far. In every step of the standard construction of a TCD, we add white and black vertices to the VRC while also adding points to the TCD map for every white vertex. If the added boundary-parallel strand is oriented from left to right and intersects $2m$ other strands we add $m$ black vertices for the intersections and $m$ white vertices for the new counterclockwise faces below the added strand. The new points have to be on lines determined by the points already chosen above them, and they are therefore already in the span of everything above them. If the added boundary-parallel strand is oriented from right to left and intersects $2m$ other strands we add $m$ black vertices for the intersections and $(m+1)$ white vertices for the new counterclockwise faces below the added strand. We can choose the left most new point arbitrary, so we choose it in general position. This adds one dimension to the span of all points. This proves that the maximal dimension is $|W|-|B| - 1$ for standard diagrams.

Moreover, when we add a right moving strand then the points one level above the added strand are included in the span of the new points below the strand. When we add a left moving strand, then the points one level above the added strand are included in the span of the union of the new points below the strand and the two points above left and above right of the added strand. Therefore, in each step all points are included in the span of the points of those white vertices that correspond to boundary faces of the TCD. By induction we deduce that all points of a TCD map that is constructed in this way are in the span of the boundary points. 

Now given an arbitrary TCD $\tcd$, we may first construct a TCD map as above associated to a standard diagram $\tcd'$ with the same endpoint matching as $\tcd$. Its maximal dimension is $|W|-|B|-1$. This difference does not change under 2-2 moves and there is a sequence that takes $\tcd'$ to $\tcd$. The dimension cannot change under 2-2 moves and the boundary points are not changed. Therefore the maximal dimension of $\tcd$ is $|W|-|B|-1$ as well. The endpoint matching is also unchanged under 2-2 moves. Therefore the theorem holds for any TCD.

The last argument contains a technical inaccuracy. As the TCD map $T'$ that we construct for a standard diagram $\tcd'$ is not necessarily flip-generic, it may happen that the sequence of 2-2 moves that takes $\tcd'$ to $\tcd$ involves a spider move that is not well-defined. It is plausible that we can avoid this by ``wiggling'' a little at the TCD map $T'$. A rigorous proof of the theorem is given in Section \ref{sec:fromprojinvariants}, where Lemma \ref{lem:opensetoftcdmaps} shows that the construction algorithm given in Definition \ref{def:constralgo} can construct a TCD $T$ of maximal dimension for any minimal TCD $\tcd$. However, that lemma relies on a lot of combinatorics and algebra that we will only develop in the subsequent chapters. Therefore, we chose here to give the ``plausible'' argument above instead.\qed
}

\begin{example}
	Let $T$ be a TCD map from $\tcdp$ with endpoint matching $\enm nk$. Then the maximal dimension of $T$ is $k$, because that is the number of left moving strands.
\end{example}

\section{Multi-dimensional consistency in $A_n$}\label{sec:tcdconsistency}
\begin{figure}
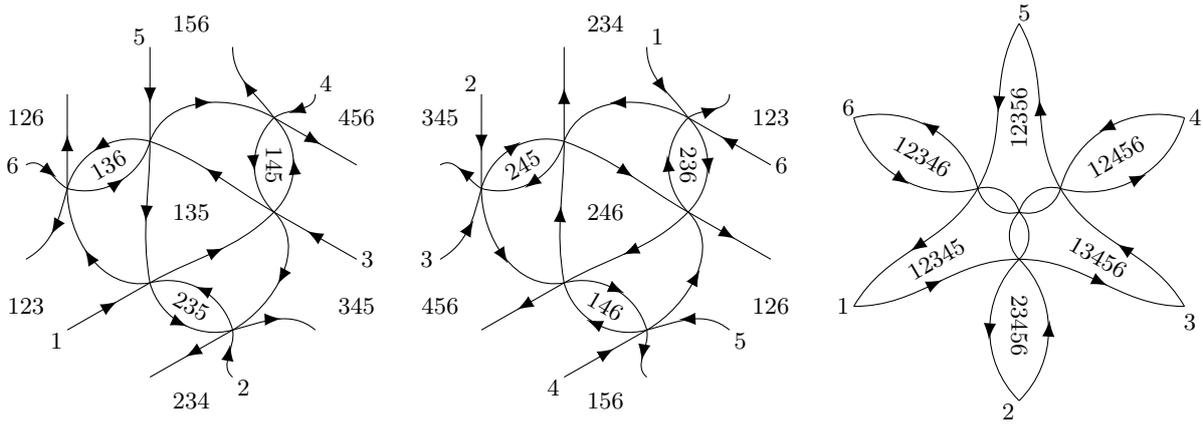


	\caption{The map $\an'$ that assigns the points of the VRC associated with a TCD to points in the $A_5$ lattice.}
	\label{fig:tcdafive}
\end{figure}

There is a natural way to embed the white vertices of $\pb$ associated to any TCD in the $A_n$ lattice, where $(n+1)$ is the number of strands in the TCD. Let
\begin{align}
	A_n =  \left\{(z_1,z_2,\dots,z_{n+1})\in\Z^{n+1}\ : \ \sum_{i=1}^{n+1}z_i = 0\right\}
\end{align}
denote the $A_n$ lattice. Thus, the $A_n$ lattice is defined as all the points in $\Z^{n+1}$ that have coordinate sum zero. We can replace the zero with any other integer and work on a shifted copy of $A_n$. This will sometimes simplify notation. In order to find the map, we can actually define a map $\an$ from all the faces of the TCD to $\Z^{n+1}$, instead of just from the counterclockwise oriented faces. We assign to each strand $s$ a unique index $i_s \in \{0,1, \dots, n\}$. We choose one counterclockwise face $f_0$ as base point and map it to the origin, i.e. $\an(f_0) = 0$. Now if we have two adjacent faces $f$ and $f'$ such that $f'$ is to the left of their common strand $s$, then we require that $\an(f')-\an(f) = e_{i_s}$. This is well defined, as around any crossing point of the TCD, we will add and subtract $i_s$ exactly once for any strand $s$ involved. Moreover, to get from one counterclockwise oriented face to another, we will have to traverse the same amount of strands from right to left as we do from left to right. Thus the coordinate sum is a constant for all white vertices. As a result all white vertices of the associated VRC are mapped to $A_n \subset \Z^{n+1}$. Because we work with minimal triple crossing diagrams, the map $\an$ only takes values in $\{0,\pm 1\}^n$. Thus we often denote the image of $\an$ simply by a set of indices, as we did in Figure \ref{fig:tcdafive}.

Alternatively, there is another map $\an'$ to a shifted copy of $A_n$ for which there is a direct construction. We assign to each face $f$ the point $\an'(f) \in \Z^{n+1}$ as follows:
\begin{align}
	\left(\an'(f)\right)_i =\begin{cases}
		1 & f \mbox{ is to the left of strand } i,\\
		0 & f \mbox{ is to the right of strand } i.\\	\end{cases} \label{eq:shiftamap}
\end{align}
Therefore the map $\an'$ indeed maps to a shifted copy of $A_n$ in $\Z^{n+1}$ and we have the relation
\begin{align}
	\an'(f) - \an'(f_0) = \an(f),
\end{align}
for all faces $f$.

Of course, on the one side one can view the map $\an$ as just the map from a subset of $A_n$ to the indices of shift-notation (which we will introduce in the beginning of Chapter \ref{cha:ddg}). The specific idea of associating $A$-type shift notation to graphs with strands already occurs in the study of cluster algebras and Grassmannians \cite{scottgrass}. A version of $\mathcal A'$ has also appeared in the literature under the name of \emph{Abel map} \cite{fockinverse, giclustermodular}. The Abel map is defined on the torus and for zig-zag paths of bipartite graphs instead of strands of a TCD. 

When looking at the triple crossing diagrams of Q-nets, Darboux maps or line complexes we observe that the strands can be partitioned into 6 (resp. 4 for $\Z^2$ Q-nets) families, such that the strands of a family do not intersect each other. This allows us to assign the same index to all strands in one family. Thus for the hexagonal cases the dimension of the A-lattice is 5, and for the square case it is 3. This view is reminiscent of the construction for isoradial graphs as combinatorial surfaces in $\Z^n$, as discussed by Bobenko, Mercat and Suris \cite{bmsconformal}.

We will take a similar view and think of TCDs as discrete combinatorial surfaces in an $A_n$ lattice, while the TCD map may be considered as a discrete geometric surface. It is natural to view the 2-2 move at a white vertex of $\pb$ (a resplit) as changing both the combinatorial surface in $A_n$ and the geometric surface in projective space. On the other hand we interpret a 2-2 move at a quad of $\pb$ (a spider move) as a reparametrization of the geometric surface because it only changes the combinatorial surface.

It is important to realize that given a TCD, we can perform a sequence of 2-2 moves (without going backwards), such that we return to the original TCD. In other words: The 2-2 flip graph has cycles. We now claim that TCD maps are consistent under 2-2 moves: if we return to the same TCD we also return to the same TCD map. To prove this, we consider the following conjecture by Thurston \cite{thurstontriple} that has recently been proven by Balitskiy and Wellman.

\begin{figure}
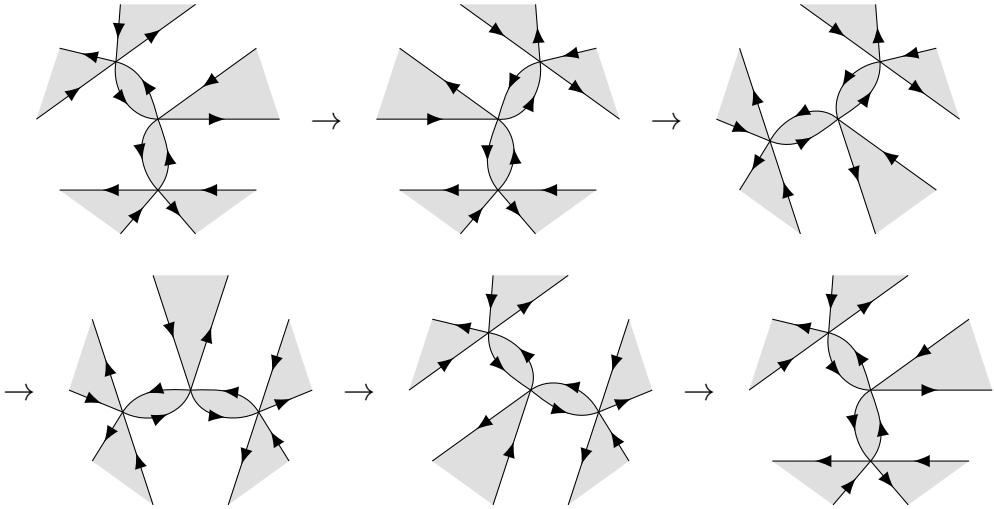

	\centering
	
	\caption{The elementary 5-cycle in the flip graph of triple crossing diagrams. Each diagram is the result of applying a 2-2 move to one of the two bigons in the previous diagram. The gray faces correspond to points in $\RP^3$. 
	}	
	\label{fig:tcdflipgraphfivecycle}
\end{figure}

\begin{figure}
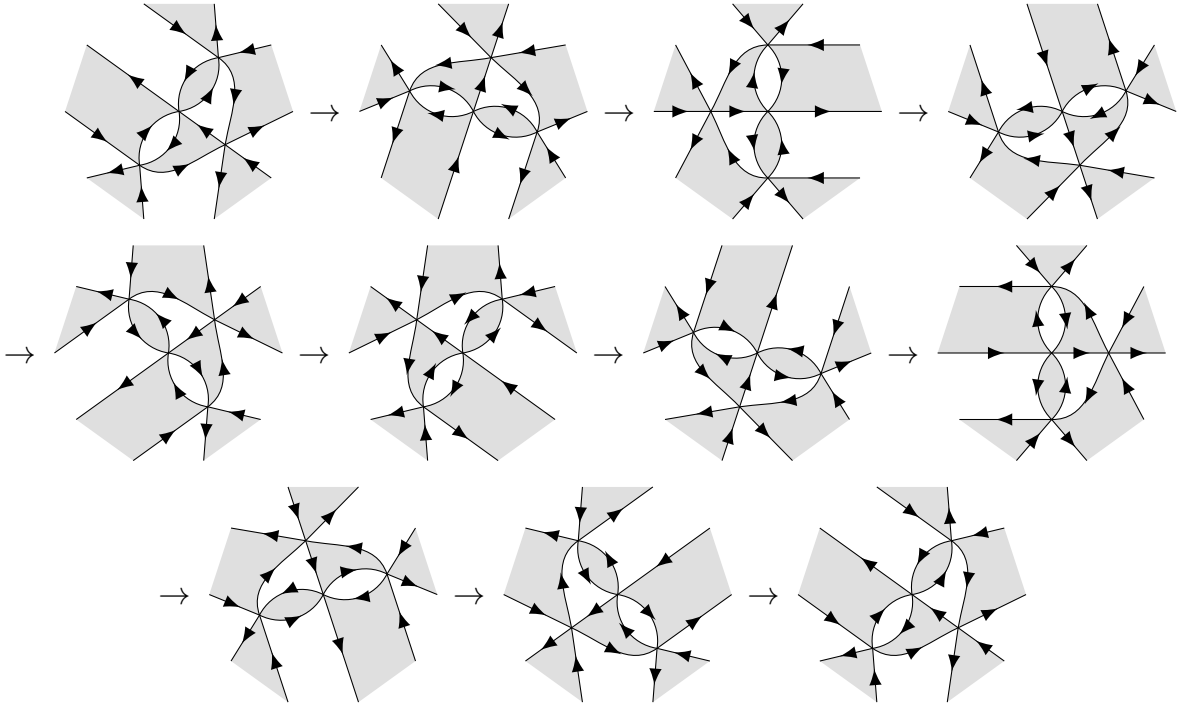

	\centering
	\hspace{8mm}
	%
			
	\caption{The elementary 10-cycle in the flip graph of triple crossing diagrams. 
	}
	\label{fig:tcdflipgraphtencycle}
\end{figure}

\begin{theorem}[\cite{bwtriple}]\label{th:tcdtwocycles}
	Let $\mathfrak F$ be the 2-complex given by the flip graph of triple crossing diagrams with fixed endpoint matching (see Definition \ref{def:endpointmatching}), with the following 2-cells glued to it:
	\begin{enumerate}
		\item a quadrilateral, wherever two flips commute because they are sufficiently far apart;
		\item a pentagon, wherever there is a 5-cycle in a subset of the diagram which	is a triple crossing diagram with endpoint matching $\enm 53$ (see Figure~\ref{fig:tcdflipgraphfivecycle}) or $\enm 51$ (Figure~\ref{fig:tcdflipgraphfivecycle} with reversed orientations);
		\item a decagon, wherever there is a 10-cycle in a subset of the diagram which is a triple crossing diagram with endpoint matching $\enm 52$ (see Figure~\ref{fig:tcdflipgraphtencycle}).
	\end{enumerate}
	Then $\mathfrak F$ is simply connected.
\end{theorem}

The statement that $\mathfrak F$ is simply connected implies that any cycle in the flip graph is generated by a combination of the three types of elementary cycles listed in the theorem. This theorem allows us to give an elementary proof of the multi-dimensional consistency of TCD maps with regard to 2-2 moves.

\begin{figure}
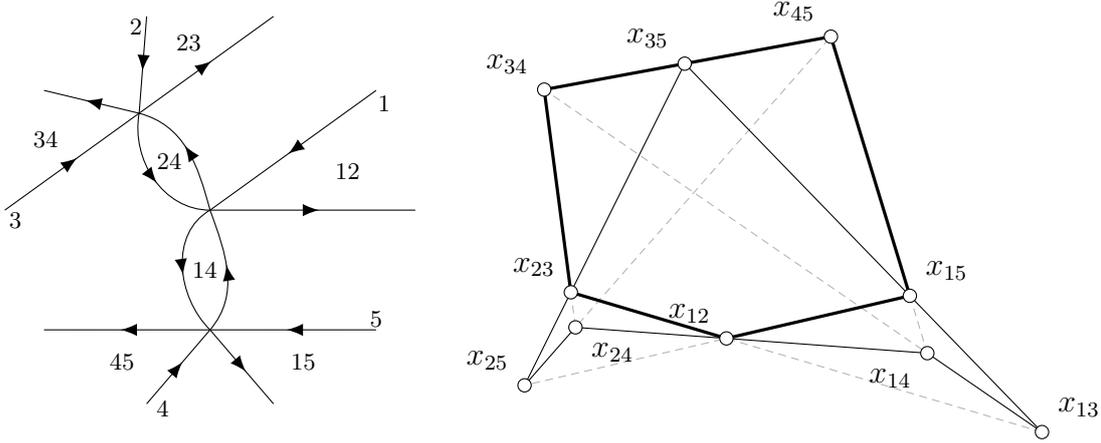



	\caption{Desargues' theorem. On the left is Cauchy data for Desargues' configuration on the right.}
	\label{fig:desargues}
\end{figure}

\begin{theorem} \label{th:tcdconsistency}
	Flip-generic TCD maps are consistent on $A_n$.
\end{theorem}
\proof{Because of Theorem \ref{th:tcdtwocycles}, it suffices to prove that this is true along the three types of fundamental cycles in $\mathfrak F$. Along the 4-cycles, it is clearly true because the recurrence is locally defined.

For the 10-cycle the configuration can span at most a space of dimension 2. Assuming it does indeed span a space of dimension 2 then the internal points are uniquely defined by the boundary points. Therefore the dSKP equation is consistent along the 10-cycles in dimension 2. In dimension 1 a direct algebraic calculation proves the result as well.

There are two 5-cycle configurations. The cycle of diagrams with endpoint matching $\enm 51$ has no internal points and is therefore trivially consistent. The interesting case is the $\enm 53$ endpoint matching. It can span a space of dimension up to 3. If it spans a space of dimension two or three then it is in fact the Desargues configuration (see Figure \ref{fig:desargues} and Theorem \ref{th:desargues}). For instance, the triangles $(x_{12},x_{13},x_{15})$ and $(x_{24}, x_{34}, x_{45})$ are in perspective with respect to $x_{14}$. Therefore the two triangles are also in perspective with respect to the line through $x_{23}$ and $x_{35}$, on which therefore also $x_{25}$ lies. Thus, the TCD map is well defined along the 10-cycle. In dimension 1 a direct algebraic calculation proves the result again. Thus the proof is complete. \qedhere
}

Note that the obstacle for non flip-generic TCD maps is that it may not be possible to flip along each of the fundamental cycles of Theorem \ref{th:tcdtwocycles}, thus it may be necessary to check along other cycles than those. We do not expect multi-dimensional consistency to fail for non flip-generic TCD maps, but we do not investigate this question further.

\begin{theorem}[Desargues' theorem]\label{th:desargues}
	Let $A_1,A_2,A_3$ and $B_1,B_2,B_3$ be the vertices of two triangles in $\CP^n, n \geq 2$. Then the following two statements are equivalent:
	\begin{enumerate}
		\item The three lines $A_1B_1,A_2B_2, A_3B_3$ intersect in a point.
		\item The three points $A_1A_2 \cap B_1B_2, A_2A_3 \cap B_2B_3, A_3A_1 \cap B_3B_1$ are on a line.\qedhere
	\end{enumerate}
\end{theorem}
The two triangles are called in perspective with respect to a point if condition (1) is satisfied and in perspective with respect to a line if condition (2) is satisfied. In Figure \ref{fig:desargues} we see that the two triangles $x_{14},x_{24},x_{34}$ and $x_{15},x_{25},x_{35}$ are in perspective with respect to $x_{45}$ and with respect to the line that contains the three points $x_{12},x_{23},x_{13}$. More generally, with respect to point $x_{ij}$ the two triangles $\{x_{ik}: k \notin \{i,j\} \}$ and $\{x_{jk}: k \notin \{i,j\} \}$ are in perspective. They are also in perspective with respect to the line $\{x_{kl} : k,l\notin \{i,j\} \}$. Alternatively, one can state that the pentagons $(x_{k,k+1})$ and $(x_{k,k+2})$ are inscribed about each other, which means that each point $x_{k,k+1}$ is on the line $x_{k,k+3}x_{k+1,k+3}$ and each point $x_{k,k+2}$ is on the line $x_{k,k+1}x_{k+1,k+2}$. 

It is not surprising that consistency is induced by Desargues' theorem. In the case of fundamental line complexes it was already noted that multi-dimensional consistency is due to Desargues' theorem \cite{bobenkoschieflinecomplexes}. Moreover, we can consider TCD maps in $\CP^1$, thus obtaining a new combinatorial proof of the multi-dimensional consistency of the dSKP equation. A (different) proof of the multi-dimensional consistency of the dSKP equation was first given by Adler, Bobenko and Suris and appeared in their classification of integrable equations of octahedral type \cite{absoctahedron}. On the other hand, King and Schief identified that the consistency of the dSKP equation is geometrically due to a conformal version of Desargues' theorem \cite{ksconformaldesargues}, which fits nicely with the interpretation we gave here.

\begin{remark}
	We will also consider so called Desargues maps in Section \ref{sec:desargues}. However, shortly after finishing the thesis it was brought to our attention that there is a newer, different definition of Desargues maps \cite[Proposition 3.1]{doliwadesarguesweyl}, that contains the original definition as a special case. This newer definition of Desargues maps is defined on $A_n$, and TCD maps can be understood as certain subsets of Desargues maps, by embedding TCD maps via the map $\an$ as explained above.
\end{remark}


\chapter{TCD maps and discrete differential geometry}\label{cha:ddg}\label{sec:ddgexamples}

Four of our main examples stem from discrete differential geometry (DDG): Q-nets, Darboux maps, line complexes and line compounds. These are projective, integrable $\Z^n$-type maps that are fundamental for DDG. In this chapter we focus on the introduction of these maps as TCD maps. We do this by associating a TCD map to Cauchy-data in each example. We show that the propagation of Cauchy-data corresponds to performing 2-2 moves in the TCD maps. After that, the integrability is simply a corollary of Theorem \ref{th:tcdconsistency}. Later in the thesis we give further arguments why these four maps should be considered on equal footing. In particular we will show they are related by taking sections (Section \ref{sec:sections}), projective duals (Section \ref{sec:projduality}) and focal transforms (Section \ref{sec:focalnets}). Furthermore, the framework of TCD maps yields a canonical algebraic description of these maps in terms of cluster variables that are projective invariants (Section \ref{sec:projcluster}). Many other systems -- like K{\oe}nigs nets, A-nets, circular nets, Cox lattices et cetera -- can be obtained as reductions of Q-nets, Darboux maps, line complexes or line compounds and are thus also TCD maps.

\section{Quad-graphs and multi-dimensional consistency in $\Z^n$}

We have already seen that TCD maps in general are discrete integrable on $A_n$ lattices in Section \ref{sec:tcdconsistency}. One of the important properties of the examples we want to introduce in this chapter however feature discrete integrability on $\Z^n$. With that goal in mind, we first explain the relation between maps defined on quad-graphs and discrete integrability on $\Z^n$, analogously to the results of Section \ref{sec:tcdconsistency}. 

\begin{figure}
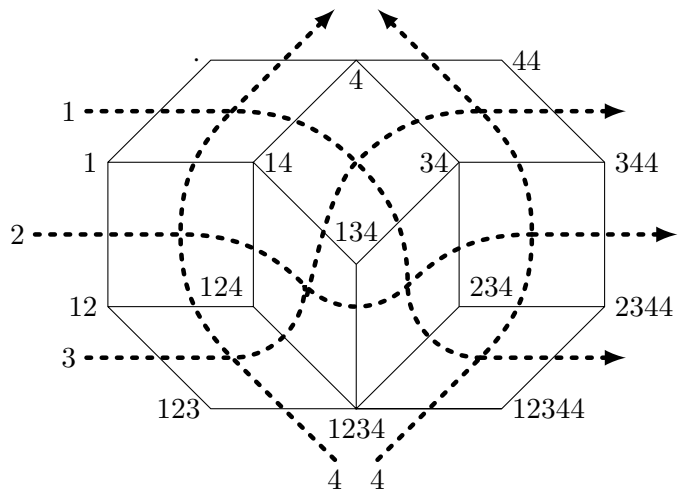


	\vspace{-4mm}
	\caption{Strips and labeling of a quad-graph.}
	\label{fig:striplabeling}
\end{figure}

We begin by introducing a way to label vertices of a minimal quad-graph $\qg$ with points of $\Z^n$, where $n$ is the number of strips of $\qg$. We also fix an arbitrary orientation of every strip of $\qg$, and assign to each strip $s$ a unique index $i_s\in\{1,2,\dots,n\}$. The labeling is a map $\mathcal Z: V(\qg) \rightarrow \Z^n$, where we write each point in $\Z^n$ by an index set as in the case of the map $\an$ in Section \ref{sec:tcdconsistency}. We choose a base-vertex $v_0\in V(\qg)$ and map it to the origin. If we have two adjacent vertices $v$ and $v'$ such that $v'$ is to the left of their common strip $s$, then we require that $\mathcal Z(v')-\mathcal Z(v) = \mathbf e_{i_s}$, where $\mathbf e_1,\mathbf e_1,\dots,\mathbf e_n$ are the unit vectors of $\Z^n$. Alternatively, there is another map $\mathcal Z'$ for which there is a direct construction. We assign to each vertex $v\in V(\qg)$ the point $\mathcal Z'(f) \in \Z^{n}$ as follows:
\begin{align}
	\left(\mathcal Z'(v)\right)_i =\begin{cases}
		1 & v \mbox{ is to the left of strip } i,\\
		0 & v \mbox{ is to the right of strip } i.\\	\end{cases}
\end{align}
The map $\mathcal Z'$ therefore satisfies the relation
\begin{align}
	\mathcal Z'(v) - \mathcal Z'(v_0) = \mathcal Z(v),
\end{align}
for any $v\in V(\qg)$ and is thus a shifted version of the map $\mathcal Z$. The image of $\mathcal Z'$ is actually contained in the hypercube $\{0,1\}^n$. However, if two strips do not intersect in $\qg$, we can assign them the same index and thus obtain labelings that do not necessarily map to the hypercube, we show an example in Figure \ref{fig:striplabeling}.

\begin{figure}
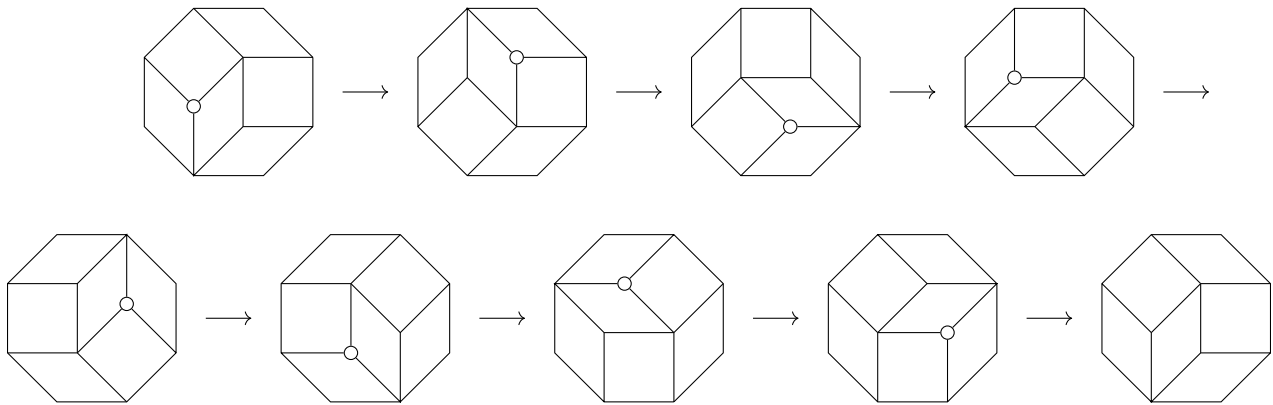


	\caption{The non-trivial cycle in the quad-flip-graph.}
	\label{fig:quadflipgraph}
\end{figure}

As explained in Section \ref{sec:quadgraphs}, the cube flip is a local move on quad-graphs. In this section we look at certain maps defined on vertices, edges or faces of quad-graphs and we introduce rules on how they change when the quad-graph undergoes a cube flip. We study whether these maps and rules are multi-dimensionally consistent. For this purpose we check that these maps are well defined along cycles in the flip graph of quad-graphs. Conveniently there is only one non-trivial cycle in the quad-flip graph, which is the 8-cycle depicted in Figure \ref{fig:quadflipgraph} (see the book by Felsner \cite{felsnerbook}, wherein the strips are called \emph{pseudolines} and the quad-graphs are called \emph{zonotopal tilings}). Thus, to check the multi-dimensional consistency of maps defined on quad-graphs it suffices to check along this 8-cycle. Note that in the DDG community multi-dimensional consistency is often checked by other means, usually by checking that the propagation of data in a 4-cube is consistent independent of the path taken. We prefer to check consistency along the 8-cycle because this yields a canonical algorithm independently of whether the map is defined on vertices, edges or faces of $\Z^n$, and because of the ease of visualization of the 8-cycle. However, that does not mean that checking along the 8-cycle is necessarily the easiest in terms of calculation. The fact that the number of strips appearing in the quad-graphs of the 8-cycle is four corresponds to the well known result \cite{ddgbook}, that a 3D-system on $\Z^3$ is multi-dimensionally consistent on $\Z^n$ if and only if it is 4D-consistent.

As we observed, all vertices of $\qg$ and quad-graphs $\qg'$ related via a sequence of flips to $\qg$ naturally live in a hypercube of $\Z^n$, if each of the $n$ strips of $\qg$ carries a different label. If some of the strips are parallel and we assign the same labels to those parallel families of strips, the vertices of $\qg$ and all $\qg'$ can be considered to be part of a \emph{brick} $B(\qg) = \prod_{i=0}^{N-1} [0,a_i] \cap \Z^N$, where the sum of all $a_i$ equals the number of strips of $\qg$ and we borrow the name brick from \cite{bmsconformal}. Thus, if a 3D-system is multi-dimensionally consistent, one may define it not only on the quad-graph $\qg$, but on all of $B(\qg)$. In fact, in the discrete integrable systems community most systems are a priori given by a definition on bricks (or all of $\Z^N$). Well-definedness from Cauchy-data and multi-dimensional consistency is proven afterwards, thus justifying the definition on bricks. We instead choose to give definitions on quad-graphs and thus Cauchy-data first, then deduce consistency from our results on TCD maps, and thus justify also working with equivalent definitions on lattices.

\begin{figure}
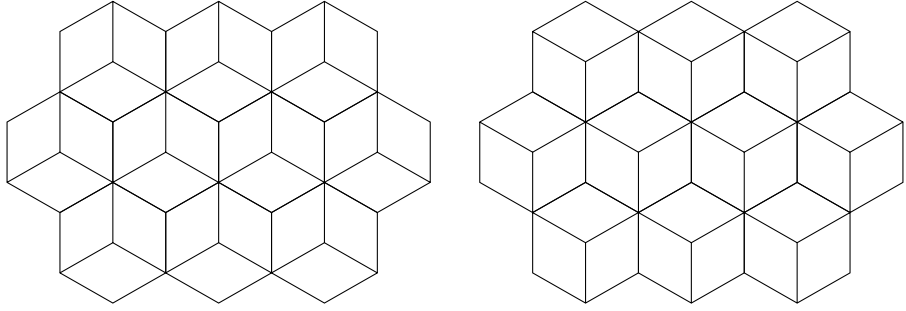


\caption{Parts of a stepped surface before and after propagation.}
\label{fig:steppedsurface}
\end{figure}

Moreover, in discrete integrable systems the quad-graphs are generally considered to be defined not only on bricks but on all of $\Z^N$. Because the lattice equations are local, most variables are local or locally integrated and multi-dimensional consistency is also a local phenomenon, this does not pose a serious problem. For example, when working with definitions on $\Z^3$, we can describe all data on $\Z^3$ by working on stepped surfaces. As a subset of $\Z^3$, a \emph{stepped surface} $U \subset \Z^3$ is for example
\begin{align}
	U = \{(z_1,z_2,z_3) \in \Z^3 \ : \ z_1 + z_2 + z_3 \in \{m-1,m,m+1\}   \},
\end{align}
for any $m\in \Z$, see also Figure \ref{fig:steppedsurface}. In this case, we consider maps that are not defined on a graph in a disc but on a whole stepped surface. If we consider a 3D system then the map restricted to a stepped surface is Cauchy-data for the map on all of $\Z^3$. This is clear because we can apply an infinite number of independent cube-flips to translate the stepped surface. By repeating this procedure we cover all of $\Z^3$.

\section{$\Z^N$ notation}\label{sec:znotation}

Before we continue, let us introduce some standard notation for the cells in $\Z^N$:
\begin{enumerate}
	\item $V(\Z^N)$ denotes the vertices (0-cells) of $\Z^N$ and thus just $\Z^N$ itself.
	\item $E(\Z^N)$ denotes the edges (1-cells) of $\Z^N$, each edge consisting of two adjacent points.
	\item $F(\Z^N)$ denotes the faces (2-cells) of $\Z^N$, each face is a quad consisting of four points.
	\item $C(\Z^N)$ denotes the cubes (3-cells) of $\Z^N$, each cube consisting of eight points.
\end{enumerate}
Moreover, we identify the dual lattice ${\Z^N}^*$ with the shifted lattice $\Z^N + \frac12 \{1\}^N$. In the case of $N=3$ this gives an identification of $V(\Z^3)$ with $C({\Z^3}^*)$, $E(\Z^3)$ with $F({\Z^3}^*)$ and vice versa.

Additionally we make use of \emph{shift notation}\label{pag:shiftnot}. Let $\mathbf e_1, \mathbf e_2, \dots, \mathbf e_N$ denote the unit vectors of $\Z^N$. If $k\in \N$ and $f$ is a function defined on the vertices of $\Z^N$ then we will abbreviate the function $n\mapsto f(n+\mathbf{e}_k)$ by $f_k$ and call this a shift in the $k$ direction. We use multiple subscripts to denote multiple shifts, $f_{kl}$ for example is a shift in $k$ and in $l$ direction. If a function $f$ is defined on the edges of $\Z^N$ then we denote $f(\{n,n_k\})$ by $f^{k}$, such that $f$ is now represented by $N$ functions $f^1, f^2,\dots, f^N$ defined on $\Z^N$. For a function $f$ on faces we similarly denote $f(\{n,n_k,n_l,n_{kl}\})$ by $f^{kl}$. We apply shift notation to functions defined on edges and faces as well. Moreover, we use a bar on subscripts to denote shifts in negative direction, that is $f_{\bar k} = f(n-\mathbf e_k)$.  Shift notation allows us to treat \emph{lattice equations} in abbreviated notation. A function $f$ on $\Z^N$ satisfies a lattice equation, if $f$ also satisfies this equation for any shift.

\begin{figure}
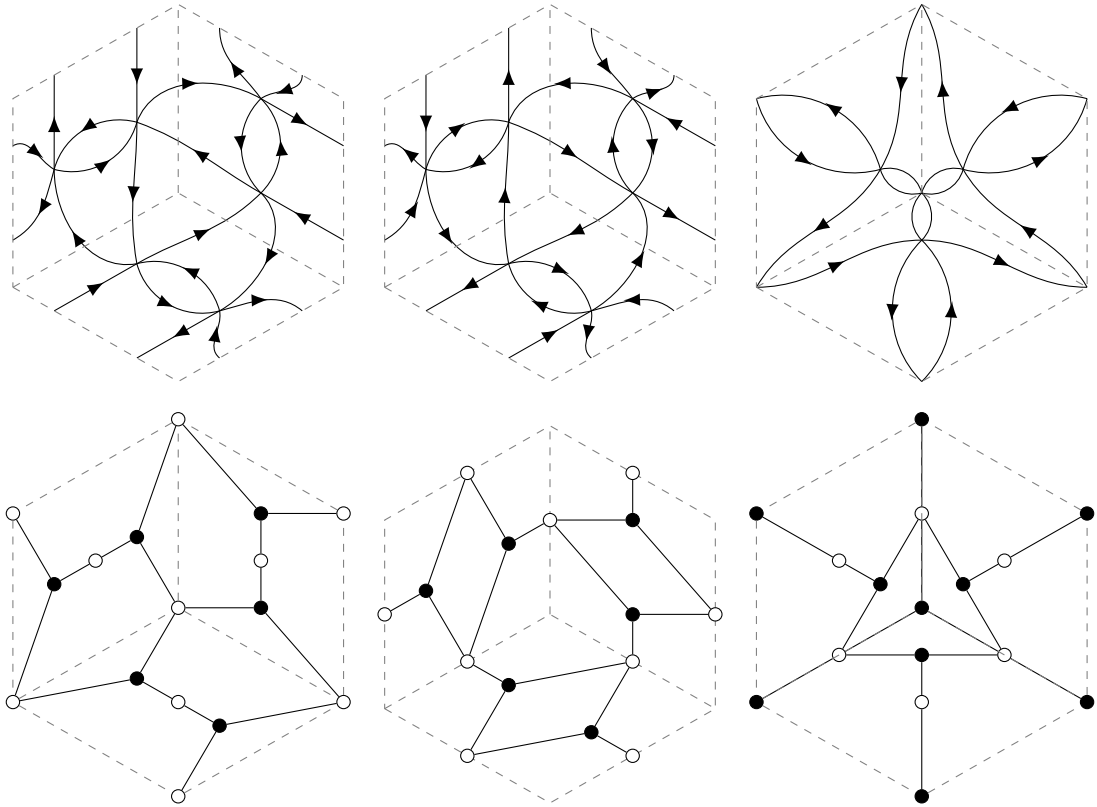

	%
	\caption{The fundamental domains of a stepped surface from left to right for: Q-net, Darboux map, line complex. The TCD in the top row and the corresponding graph $\pb$ in the bottom row.}
	\label{fig:tcdexamplesB}
\end{figure}

\section{Q-nets}\label{sec:qnets}

Q-nets with $\Z^2$ combinatorics were first considered by Sauer \cite{sauerqnet}, Q-nets on $\Z^N$ by Doliwa and Santini \cite{doliwasantiniqnet} as \emph{quadrilateral lattices}. We begin with a definition of Q-nets on quad-graphs, and give a definition of Q-nets on lattices at the end of the section.

\begin{definition}\label{def:qgqnet}
	Let $\qg$ be a quad-graph. A map $q: V(\qg) \rightarrow \CP^n$ with $n\geq 2$ is a \emph{Q-net} if the image of every quad is planar. Denote the points of a quad in $q$ which is crossed by strips $k,l$ by $q,q_k,q_l,q_{kl}$. If the four points of a quad span a plane, we also define the two \emph{focal points}
	\begin{align}
		\fp^{kl} &= qq_k \cap q_lq_{kl},\\
		\mbox{ and } \fp^{lk} &= qq_l \cap q_kq_{kl},
	\end{align}
	 as the intersections of opposite lines.
\end{definition}

The notation for the focal point is not symmetric in its two superscripts, reflecting the existence of two different focal points per quad. For a thorough introduction to Q-nets in terms of integrability and discrete differential geometry see the book by Bobenko and Suris \cite{ddgbook}.

Assume that we know the seven points $q,q_1,q_2,q_3,q_{12},q_{23},q_{13}$ of three mutually adjacent quads (the backside of a cube) in a Q-net $q$, such that these points span a 3-dimensional space. After a cube flip, the new vertex $q_{123}$ is then determined by the Q-net condition, indeed
\begin{align}
	q_{123} = \spa\{q_{12},q_{2},q_{23}\} \cap \spa\{q_{23},q_{3},q_{13}\} \cap \spa\{q_{13},q_{1},q_{12}\}. 
\end{align}
This replacement of $q$ with $q_{123}$ is the \emph{cube flip} of Q-nets. The cube flip can also be defined if the points only span a 2- or 1-dimensional space, which we explain in Section \ref{sec:nongendimension}.

Recall that we consider the vertices of a quad-graph $\qg$ to be bipartitioned into black ($\blacksquare$) and white ($\square$) vertices, and that each quad inherits an orientation from a fixed embedding of $\qg$ into a disk.

To each Q-net, we associate a TCD map. We obtain the associated graph $\pb$ by gluing the following piece of graph
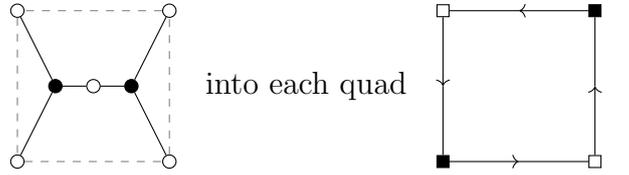
\begin{center}\vspace{-2mm}
\begin{tikzpicture}[scale=1]
	\begin{scope}[shift={(-0.2,0)}]
		\node[wvert] (v) at (0,0) {};
		\node[wvert] (v1) at (2,0) {};
		\node[wvert] (v2) at (0,2) {};
		\node[wvert] (v12) at (2,2) {};
		\node[wvert] (f21) at (1,1) {};
		\node[bvert] (b2) at (0.5,1) {};
		\node[bvert] (b21) at (1.5,1) {};
		\draw[gray, dashed]
			(v) -- (v1) -- (v12) -- (v2) -- (v)
		;	
		\draw[-]				
			(b2) edge (v) edge (v2) edge (f21)
			(b21) edge (v1) edge (v12) edge (f21)
		;
		\node[] at (3.8,1) { into each quad };
	\end{scope}
	\begin{scope}[shift={(5.4,0)}]
		\node[bqgvert] (v) at (0,0) {};
		\node[wqgvert] (v1) at (2,0) {};
		\node[wqgvert] (v2) at (0,2) {};
		\node[bqgvert] (v12) at (2,2) {};
		\coordinate (e1) at (1,0);
		\coordinate (e2) at (0,1);
		\coordinate (e12) at (1,2);
		\coordinate (e21) at (2,1);				
		\draw[black]
			(v) edge[->] (e1) edge (e2)
			(v1) edge[->] (e21) edge (e1)
			(v12) edge[->] (e12) edge (e21)
			(v2) edge[->] (e2) edge (e12)
		;			
		\node[] at (2.5,1) {.};
	\end{scope}
\end{tikzpicture}	
\end{center}

Note that one could make a different definition by rotating the glued piece of graph by 90 degrees, which would be an equally valid convention.
 
The TCD map is such that the white vertices of $\pb$ corresponding to vertices of the quads are mapped to the corresponding points of the Q-net. The white vertex inside each quad is mapped to one of the corresponding focal points, which is determined by the location of the black vertices of $\pb$. The corresponding TCD and $\pb$ for a fundamental domain of a stepped surface are shown on the left of Figure \ref{fig:tcdexamplesB}. 

In fact, we can even consider Q-nets in $\CP^1$ via the corresponding TCD maps. In this case however the focal points are not well-defined, and we may choose their position generically. This is explained in more detail in Section \ref{sec:nongendimension}.

The sequence of eleven 2-2 moves shown in Figure \ref{fig:tcdcubeflip} induces the cube flip of the Q-net. Therefore the propagation of Cauchy-data in a Q-net corresponds to sequences of 2-2 moves on the level of the TCD map. 

\begin{remark}
	Note that there is also a sequence of only ten 2-2 moves that we will discuss later (see Figure \ref{fig:tensteps}), that can also be viewed as a cube flip. That sequence however does not respect our convention on which focal point is represented by the TCD map in each quad. Of course, it is not actually necessary to explicitly give the 2-2 sequence that corresponds to the cube flip, because of Theorem \ref{th:tcdflipsconnected}. That theorem guarantees that there is such a sequence of 2-2 moves, because one easily verifies that the strand connectivity before and after the cube flip in the associated TCD map of the Q-net is the same.
\end{remark}

The multi-dimensional consistency of Q-nets is well known \cite{doliwasantiniqnet}. In our setup the consistency is an immediate corollary of the multi-dimensional consistency of TCD maps.
\begin{corollary}
	Q-nets defined on minimal quad-graphs are multi-dimensionally consistent.
\end{corollary}
\proof{It is necessary to check that Q-nets are consistent under the flip-graph cycle depicted in Figure \ref{fig:quadflipgraph}. That cycle is a sequence of cube flips which corresponds to a sequence of 2-2 moves in the associated TCD map. Moreover, TCD maps defined on a minimal TCD are consistent due to Theorem \ref{th:tcdconsistency}. \qed }

As the multi-dimensional consistency of Q-nets is now established, we know that by applying cube flips we can extend Q-nets from Cauchy-data to all of $\Z^N$. We now give a definition for Q-nets defined directly on $\Z^N$.

\begin{definition}
	Fix $n,N\in \N$ with $n\geq 3, N\geq 2$. A map $q: V(\Z^N) \rightarrow \CP^n$ is a \emph{Q-net} if the image of every quad is planar.

\end{definition}

The restriction that $n \geq 3$ is needed for $N>2$ to relate the eight vertices involved in a cube flip. The constraint $n\geq3$ can be dropped if we relate the eight vertices and twelve focal points via the cube-flip induced by 2-2 moves of the corresponding TCD map.

\section{Darboux maps}\label{sec:darbouxmap}
Our next objects of study are Darboux maps which were introduced by Schief \cite{schieflattice}. As in the case of Q-nets, we begin with a definition of Darboux maps on quad-graphs and give a lattice-based definition after we established multi-dimensional consistency.

\begin{definition}
	Let $\qg$ be a quad-graph. A map $d: E(\qg) \rightarrow \CP^n$ is a \emph{Darboux map} if the image of every quad is colinear.
\end{definition}

Assume that we know the nine points $d^1,d^1_{2},d^1_{3},d^2,d^2_{1},d^2_{3},d^3,d^3_{2},d^3_{1}$ of three mutually adjacent quads (the backside of a cube) in a Darboux map $d$, and that these points span a 2-dimensional space. After a cube-flip, the new vertices $d^{1}_{23},d^{2}_{13},d^{3}_{12}$ are then determined by the Darboux map conditions, indeed
\begin{align}
	d^1_{23} &= d^1_{2}d^2_{1} \cap d^1_{3}d^3_{1},\\
	d^2_{13} &= d^2_{1}d^1_{2} \cap d^2_{3}d^3_{2},\\
	d^3_{12} &= d^3_{1}d^1_{3} \cap d^3_{2}d^2_{3}.
\end{align}
This replacement of $d^1,d^2,d^3$ with $d^{1}_{23},d^{2}_{13},d^{3}_{12}$ is the \emph{cube flip} of Darboux maps. The cube-flip can also be defined if the points only span a 1-dimensional space, which we explain in Section~\ref{sec:nongendimension}.

To each Darboux map, we associate a TCD map. We obtain the associated graph $\pb$ by gluing the following piece of graph
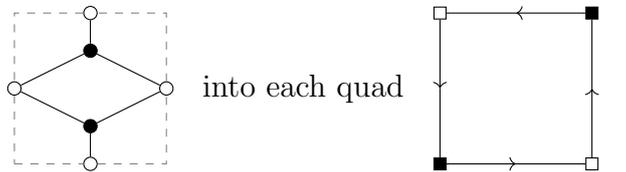
\begin{center}\vspace{-2mm}
\begin{tikzpicture}[scale=1]
	\begin{scope}[shift={(-0.2,0)}]
		\coordinate (v) at (0,0) {};
		\coordinate (v1) at (2,0) {};
		\coordinate (v2) at (0,2) {};
		\coordinate (v12) at (2,2) {};
		\node[wvert] (e1) at (1,0) {};
		\node[wvert] (e2) at (0,1) {};
		\node[wvert] (e12) at (1,2) {};
		\node[wvert] (e21) at (2,1) {};
		\node[bvert] (b1) at (1,0.5) {};
		\node[bvert] (b12) at (1,1.5) {};
		\draw[gray, dashed]
			(v) -- (e1) -- (v1) -- (e21) -- (v12) -- (e12) -- (v2) -- (e2) -- (v)
		;	
		\draw[-]				
			(b1) edge (e1) edge (e21) edge (e2)
			(b12) edge (e12) edge (e21) edge (e2)
		;
		\node[] at (3.8,1) { into each quad };
	\end{scope}
	\begin{scope}[shift={(5.4,0)}]
		\node[bqgvert] (v) at (0,0) {};
		\node[wqgvert] (v1) at (2,0) {};
		\node[wqgvert] (v2) at (0,2) {};
		\node[bqgvert] (v12) at (2,2) {};
		\coordinate (e1) at (1,0);
		\coordinate (e2) at (0,1);
		\coordinate (e12) at (1,2);
		\coordinate (e21) at (2,1);				
		\draw[black]
			(v) edge[->] (e1) edge (e2)
			(v1) edge[->] (e21) edge (e1)
			(v12) edge[->] (e12) edge (e21)
			(v2) edge[->] (e2) edge (e12)
		;			
		\node[] at (2.5,1) {.};
	\end{scope}
\end{tikzpicture}	
\end{center}

Each white vertex corresponds to an edge of $\qg$. Thus we require that the associated TCD map maps each white vertex to the corresponding point of the Darboux map. We observe that the two black vertices in each quad describe the fact that the four points of that quad are colinear. The TCD and the graph $\pb$ for a stepped surface are shown in the middle of Figure~\ref{fig:tcdexamplesB}. Note that the TCD is the same as for a Q-net except that the orientations of all strands are reversed. Thus the cube flip in Darboux maps is also captured by a sequence of 2-2 moves, as we can apply the same sequence as in the case of Q-nets.

The multi-dimensional consistency of Darboux maps was shown when they were introduced \cite{schieflattice}. In our setup the consistency is an immediate corollary of the multi-dimensional consistency of TCD maps.
\begin{corollary}
	Darboux maps  defined on minimal quad-graphs are multi-dimensionally consistent.
\end{corollary}
\proof{It is necessary to check that Darboux maps are consistent under the flip-graph cycle depicted in Figure \ref{fig:quadflipgraph}. That cycle is a sequence of cube flips which corresponds to a sequence of 2-2 moves in the associated TCD map. Moreover, TCD maps defined on minimal TCD are consistent due to Theorem \ref{th:tcdconsistency}. \qed }

As the multi-dimensional consistency of Darboux maps is now established, we give a definition for Darboux maps on $\Z^N$.

\begin{definition}
	Fix $n,N\in \N$ with $n,N \geq 2$. A map $d: E(\Z^N) \rightarrow \CP^n$ is a \emph{Darboux map} if the image of every quad is colinear.
\end{definition}

The restriction that $n \geq 2$ is needed for $N>2$ to relate the twelve vertices involved in a cube flip. The constraint $n \geq 2$ can be dropped if we relate the twelve vertices via the cube-flip induced by 2-2 moves of the corresponding TCD map.

Note that a Darboux map on $\Z^2$ is the same as a Q-net on $\Z^2$. Each quadruple of colinear points in the Darboux map corresponds to two vertices $q,q_k$ and two focal points $f^{kl}_{\bar l},f^{kl}$ that are colinear in the equivalent Q-net. 

We explain later that h-embeddings (see Section \ref{sec:harmonicemb}) are a special case of Darboux maps. We also show that certain isotropic Darboux maps correspond to Cox lattices in Section \ref{sub:cox}. Moreover, an interesting relation appears in Exercises 2.8 and 2.9 in the DDG book \cite{ddgbook}. Take a Q-net $q$ and a fixed hyperplane $H$, then intersect every line $qq_k$ with $H$. The obtained map is clearly a Darboux map. We expand on this quite a bit in Section \ref{sec:sections} on sections of TCD maps.

\section{Line complexes}\label{sec:linecomplex}
\begin{figure}
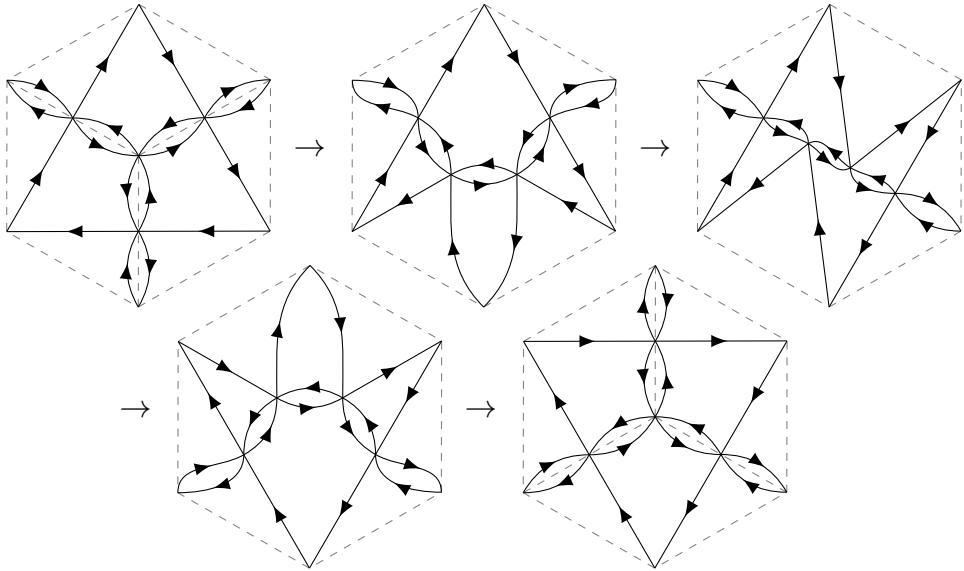

		
	
	\caption{The cube flip in a line complex. 
	}
	\label{fig:tidSTflc}
\end{figure}

Line complexes were introduced by Bobenko and Schief \cite{bobenkoschieflinecomplexes}, and their theory builds upon the study of line congruences \cite{dsmlinecongruence} by Doliwa, Santini and Mañas. We begin with a definition of line complexes on quad-graphs. We will give a definition on $\Z^n$ lattices at the end of the section.

\begin{definition} \label{def:lc}\label{def:lineslc}
	Let $\qg$ be a quad-graph. A map $l: E(\qg) \rightarrow \CP^n$ is a \emph{line complex} if for every vertex $v\in V(\qg)$ the images of the incident edges are contained in a line. A map $\ell: V(\qg) \rightarrow \{\mbox{Lines of } \CP^n \}$ yields the \emph{lines of a line complex} if adjacent lines intersect.
\end{definition}

In the literature, a line complex $l$ is defined via what we call the lines of a line complex $\ell$. However, the two definitions are equivalent. For our purposes, the definition via intersection points $l$ is more canonical and thus we primarily use this definition.

Consider the case where we already know seven lines $\ell,\ell_{1},\ell_{2},\ell_{3},\ell_{12},\ell_{13},\ell_{23}$ of three mutually adjacent quads, and that these lines span a 4-dimensional space. Then the eighth line $\ell_{123}$ has to pass through the three lines $\ell_{12},\ell_{13},\ell_{23}$. In four dimensions there is a unique line $\ell_{123}$ intersecting $\ell_{12},\ell_{13},\ell_{23}$. The replacement of $\ell$ with $\ell_{123}$ is the \emph{cube flip} of the lines of a line complex. Of course, the new line $\ell_{123}$ also defines the new intersection points $l^{1}_{23}, l^{2}_{13}, l^{3}_{12}$ and the replacement of $l^{1}, l^{2}, l^{3}$ with $l^{1}_{23}, l^{2}_{13}, l^{3}_{12}$ is the \emph{cube flip} of a line complex.

The cube flip can also be defined if the lines only span a 3-dimensional space, and in this case line complexes have been studied under the name of \emph{fundamental line complexes} \cite{bobenkoschieflinecomplexes}. Moreover, the cube flip can even be defined if the lines only span a 1-dimensional space, which we will explain in Section \ref{sec:nongendimension}.

To each line complex, we associate a TCD map. However, unlike in the cases of Q-nets and Darboux maps, our definition of a corresponding TCD does not eliminate all local choices. The graph $\pb$ of a line complex is obtained by adding one white vertex $w_e$ to $\pb$ for every edge $e$ of $\qg$. Moreover, for each vertex $v$ of $\qg$ we add a small piece of graph $\pb_v$ to $\pb$ that has endpoint matching $\enm {d_v}1$, such that the $d_v$ boundary vertices of $\pb_v$ are identified with the white vertices $w_{e_1},w_{e_2},\dots,w_{e_{d_v}}$, where $e_i$ are the edges incident to $v$ in $\qg$. To obtain a quad-wise graphical representation, we glue

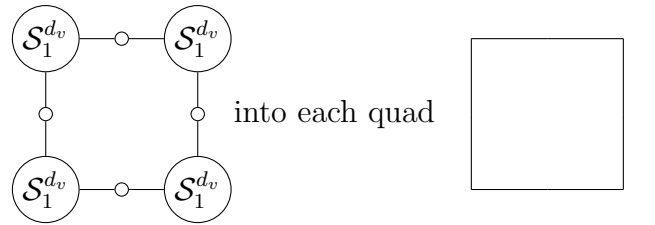
\begin{center}\vspace{-2mm}
\begin{tikzpicture}[scale=1]
	\begin{scope}[shift={(-0.2,0)}]
		\node[draw,circle,inner sep=1pt] (v) at (0,0) {$\enm {d_v}1$};
		\node[draw,circle,inner sep=1pt] (v1) at (2,0) {$\enm {d_v}1$};
		\node[draw,circle,inner sep=1pt] (v2) at (0,2) {$\enm {d_v}1$};
		\node[draw,circle,inner sep=1pt] (v12) at (2,2) {$\enm {d_v}1$};
		\node[wvert] (e1) at (1,0) {};
		\node[wvert] (e2) at (0,1) {};
		\node[wvert] (e12) at (1,2) {};
		\node[wvert] (e21) at (2,1) {};
		\draw[-]
			(v) -- (e1) -- (v1) -- (e21) -- (v12) -- (e12) -- (v2) -- (e2) -- (v)
		;				
		\node[] at (3.8,1) { into each quad };
	\end{scope}
	\begin{scope}[shift={(5.4,0)}]
		\coordinate[] (v) at (0,0) {};
		\coordinate[] (v1) at (2,0) {};
		\coordinate[] (v2) at (0,2) {};
		\coordinate[] (v12) at (2,2) {};
		\coordinate (e1) at (1,0);
		\coordinate (e2) at (0,1);
		\coordinate (e12) at (1,2);
		\coordinate (e21) at (2,1);				
		\draw[black]
			(v) edge[] (e1) edge (e2)
			(v1) edge[] (e21) edge (e1)
			(v12) edge[] (e12) edge (e21)
			(v2) edge[] (e2) edge (e12)
		;			
		\node[] at (2.5,1) {.};
	\end{scope}
\end{tikzpicture}	
\end{center}

\begin{figure}
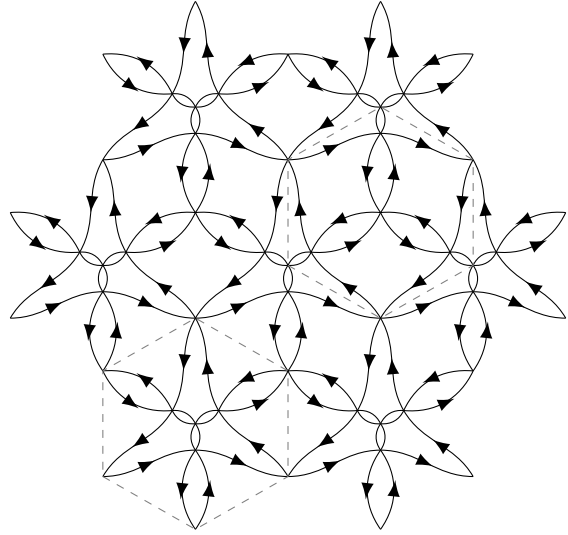


	\caption{Seven fundamental domains of the TCD of a line complex. Highlighted are two possible choices of fundamental domains that differ by a reflection about the horizontal axis.}
	\label{fig:biglctcd}
\end{figure}

The $\pb_v$ pieces encode that the intersection points on the adjacent edges of each vertex are on a line. In the case that $\qg$ is a stepped surface, one can make a highly symmetric gluing, see the right of Figure \ref{fig:tcdexamplesB}. However, there are two such possibilities as visualized in Figure \ref{fig:biglctcd}. 

As in the case of Q-nets and Darboux maps, the cube flip corresponds to a sequence of 2-2 moves in the corresponding TCD map, as depicted in Figure \ref{fig:tidSTflc}. Note that because of the different possible choices of $\pb_v$ at the vertices of $\qg$, it may be necessary to apply a finite number of spider moves (that is reparametrizations of the TCD map) at the vertices to be able to start the cube flip sequence. That is, we may need to apply spider moves at the vertices to be in the initial situation as depicted in the first diagram in Figure \ref{fig:tidSTflc}.

The multi-dimensional consistency of line complexes is well known \cite{ddgbook} and can be verified in various ways. In the setup of TCD maps however, the multi-dimensional consistency of line complexes is an elementary corollary of the multi-dimensional consistency of TCD maps.
\begin{corollary}
	Line complexes  defined on minimal quad-graphs are multi-dimensionally consistent.
\end{corollary}
\proof{It is necessary to check that line complexes are consistent under the flip-graph cycle depicted in Figure \ref{fig:quadflipgraph}. That cycle is a sequence of cube flips which corresponds to a sequence of 2-2 moves in the associated TCD map. Moreover, TCD maps defined on minimal TCD are consistent due to Theorem \ref{th:tcdconsistency}.\qed}

As the multi-dimensional consistency of line complexes is now established, we give a definition for line complexes defined on $\Z^N$.

\begin{definition}\label{def:lclattice}
	Fix $n,N\in \N$ with $n \geq 4, N \geq 2$. A map $l: E(\Z^N) \rightarrow \CP^n$ is a \emph{line complex} if for any vertex, the adjacent points are on a common line. A map $\ell: V(\Z^N) \rightarrow \{\mbox{Lines of } \CP^n\}$ is the \emph{lines of a line complex} if lines of adjacent vertices intersect.
\end{definition}

The restriction that $n \geq 4$ is needed for $N>2$ to relate the eight lines involved in a cube flip. The constraint $n \geq 4$ can be dropped if we relate the twelve intersection points via the cube-flip induced by 2-2 moves of the corresponding TCD map.

As in the case of Darboux maps, for $N=2$ a line complex is just a reformulation of a Q-net defined on $\Z^2$. Observe that in a Q-net on $\Z^2$ all the lines in one direction (for example $qq_k$ for fixed $k$) form the lines of a line complex.

\section{Line compounds}\label{sec:lcco}

We introduce a new type of map, the line compound, that has not previously appeared in the literature. We discuss later in this section why line compounds can be considered to be a different way to extend the definition of line complexes from $\Z^3$ to to $\Z^N$ for $N>3$.

\begin{definition} \label{def:colc}
	Let $\qg$ be a quad-graph. A map $l: F(\qg) \rightarrow \CP^n$ is a \emph{line compound} if for every vertex $v\in V(\qg)$ the images of the adjacent faces are contained in a $(d_v - 2)$-dimensional space $S_v$, where $d_v$ is the degree of $v$.
\end{definition}

Consider the case of a vertex $v\in V(\qg)$ of degree three. Each of the three adjacent quads $(v,v_1,v_{12},v_2), (v,v_2,v_{23},v_3), (v,v_3,v_{13},v_1)$  is mapped by a line compound to points $l^{12},l^{23},l^{13}$. Also denote the projective subspaces associated to $v$ and the other vertices of the three adjacent quads by $S,S_1,S_2,S_3,S_{12},S_{13},S_{23}$. The space $S$ is clearly a line. Denote by $Q_i$ for $i \in \{1,2,3\}$ the spaces that are spanned by the quads incident to $v_i$ except for the three quads adjacent to $v$. Because of the dimension restrictions to the spaces spanned by adjacent faces, the intersection $Q_i$ with $S$ is a point. We define the \emph{cube flip} of a line compound such that the new points after the cube flip are
\begin{align}
	l^{23}_1 = S \cap Q_1, \quad l^{13}_2 = S \cap Q_2, \quad l^{12}_3 = S \cap Q_3.
\end{align}
It is straight forward to check that the new intersection points satisfy the subspace dimension requirements. Note that by definition of the cube-flip, the six points on the faces of a cube are colinear.

To each line compound, we associate a TCD map. The corresponding TCD is not uniquely defined. Instead, we give a definition that involves some local choices. The graph $\pb$ of a line compound is obtained by adding one white vertex $w_f$ to $\pb$ for every face $f$ of $\qg$. Moreover, for each vertex $v$ of $\qg$ we add a small piece of graph $\pb_v$ to $\pb$ that has endpoint matching $\enm {d_v}{-2}$, such that the $d_v$ boundary vertices of $\pb_v$ are identified with the white vertices $w_{f_1},w_{f_2},\dots,w_{f_{d_v}}$, where $f_i$ are the faces incident to $v$ in $\qg$. To obtain a quad-wise graphical representation, we glue
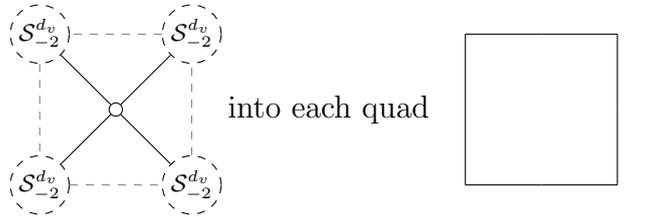
\begin{center}\vspace{-2mm}
\begin{tikzpicture}[scale=1]
	\begin{scope}[shift={(-0.2,0)}]
		\node[draw,circle,dashed,inner sep=1] (v) at (0,0) {\scriptsize$\enm {d_v}{-2}$};
		\node[draw,circle,dashed,inner sep=1] (v1) at (2,0) {\scriptsize$\enm {d_v}{-2}$};
		\node[draw,circle,dashed,inner sep=1] (v2) at (0,2) {\scriptsize$\enm {d_v}{-2}$};
		\node[draw,circle,dashed,inner sep=1] (v12) at (2,2) {\scriptsize$\enm {d_v}{-2}$};
		\node[wvert] (c) at (1,1) {};
		\draw[-]
			(c) edge (v) edge (v1) edge (v2) edge (v12)
		;
		\draw[gray, dashed]
			(v) -- (v1) -- (v12) -- (v2) -- (v)
		;				

		\node[] at (3.8,1) { into each quad };
	\end{scope}
	\begin{scope}[shift={(5.4,0)}]
		\coordinate[] (v) at (0,0) {};
		\coordinate[] (v1) at (2,0) {};
		\coordinate[] (v2) at (0,2) {};
		\coordinate[] (v12) at (2,2) {};
		\coordinate (e1) at (1,0);
		\coordinate (e2) at (0,1);
		\coordinate (e12) at (1,2);
		\coordinate (e21) at (2,1);				
		\draw[black]
			(v) edge[] (e1) edge (e2)
			(v1) edge[] (e21) edge (e1)
			(v12) edge[] (e12) edge (e21)
			(v2) edge[] (e2) edge (e12)
		;			
		\node[] at (2.5,1) {.};
	\end{scope}
\end{tikzpicture}	
\end{center}

The $\pb_v$ pieces encode that the intersection points on the adjacent faces are in the associated subspace $S_v$. In the case that $\qg$ is a stepped surface, one can make a highly symmetric gluing, which is identical to the TCD for a stepped surface of a line complex, see Figure \ref{fig:tcdexamplesB}.

The cube flip corresponds to a sequence of 2-2 moves in the corresponding TCD map. In fact the sequence is the same sequence as in the case of a line complex (see Figure  \ref{fig:tidSTflc}), only the identification of the quads in the flip is different. Note that it may be necessary to apply a finite number of resplits at the vertices to be able to start the cube flip sequence. 

As line compounds are a new object, the multi-dimensional consistency is also a new result, albeit an easy corollary of the consistency of TCD maps.

\begin{corollary}
	Line compounds defined on minimal quad-graphs are multi-dimensionally consistent.
\end{corollary}
\proof{It is necessary to check that line compounds are consistent under the flip-graph cycle depicted in Figure \ref{fig:quadflipgraph}. That cycle is a sequence of cube flips which corresponds to a sequence of 2-2 moves in the associated TCD map. Moreover, TCD maps defined on minimal TCD are consistent due to Theorem \ref{th:tcdconsistency}.\qed}

As the multi-dimensional consistency of line complexes is now established, let us now give a definition for line complexes defined on $\Z^N$.

\begin{definition}\label{def:colclattice}
	Fix $n,N\in \N$ with $n \geq 2, N \geq 3$. A map $l: F(\Z^N) \rightarrow \CP^n$ is a \emph{line compound} if for any 3-cube, the points on the six faces of that 3-cube are on a common line. A map $\ell: C(\Z^N) \rightarrow \{\mbox{Lines of } \CP^n\}$ is the \emph{lines of a line compound} if lines of face-adjacent 3-cubes intersect.
\end{definition}

As before, we can drop the restriction of the dimension of the projective space if we require that the points in a cube are related by the cube flip induces by 2-2 moves of the associated TCD map.

Note that if $\qg$ is $\Z^2$, then the definition of the line compound coincides with the definitions of Q-nets, Darboux maps and line complexes up to some identification of the combinatorics. In particular, a line compound defined on $\Z^2$ is just a Q-net defined on $(\Z^2)^*$. Moreover, a line compound defined on $\Z^3$ is simply a line complex defined on $(\Z^3)^*$. However, in $\Z^N$ for $N\geq 4$ there is no obvious identification of line compounds with line complexes. line compounds will occur when we investigate a certain reduction of line complexes (Doliwa compounds, see \ref{sec:Doliwacomplex}), that are naturally multi-dimensionally consistent as line compounds instead of line complexes. Line compounds also appear when considering sections of Darboux maps on $\Z^N$ for $N\geq 4$, see Section \ref{sec:sections}.

\section{Q-nets on $\Z^2$ and Laplace-Darboux dynamics}\label{sec:laplacedarboux}

In this section we study Laplace-Darboux dynamics \cite{doliwalaplace}, which are defined on $\Z^2$ Q-nets.

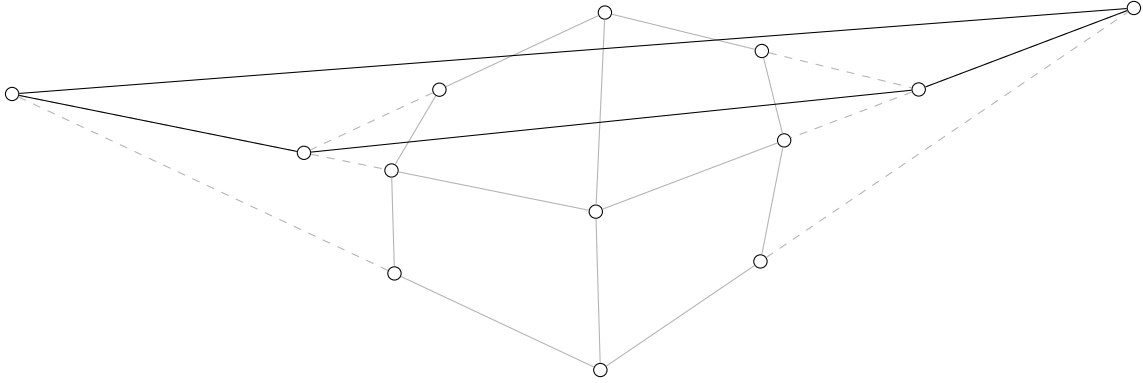
\begin{figure}
	\begin{tikzpicture}[scale=.6]
	\node[wvert] (q1) at (5.4,-9.4) {};
	\node[wvert] (q12) at (5.3,-5.9) {};
	\node[wvert] (q122) at (5.5,-1.5) {};
	\node[wvert] (f) at (-7.5,-3.3) {};
	\node[wvert] (f2) at ($(q12)!.5!(f)$) {};
	\node[wvert] (f1) at (17.1,-1.4) {};
	\node[wvert] (f12) at ($(q12)!.6!(f1)$) {};
	\node[wvert] (q) at ($(q1)!.35!(f)$) {};
	\node[wvert] (q2) at ($(q12)!.35!(f)$) {};
	\node[wvert] (q22) at ($(q122)!.55!(f2)$) {};
	\node[wvert] (q11) at ($(q1)!.3!(f1)$) {};
	\node[wvert] (q112) at ($(q12)!.35!(f1)$) {};
	\node[wvert] (q1122) at ($(q122)!.5!(f12)$) {};
	
	\draw[gray!60,-]
		(q2) -- (q) -- (q1) -- (q11) -- (q112) -- (q12) -- (q2) -- (q22) -- (q122) -- (q1122) -- (q112)
		(q1) -- (q12) -- (q122)
	;
	\draw[gray!60, dashed]
		(q) -- (f) -- (f2) -- (q2)
		(q22) -- (f2)
		(q11) -- (f1) -- (f12) -- (q112)
		(q1122) -- (f12)
	;
	\draw[-]
		(f) -- (f1) -- (f12) -- (f2) -- (f)
	;		
	\end{tikzpicture}
	\caption{Four quads of a $\Z^2$ Q-net drawn in gray and one quad of the Laplace-Darboux transform in black.}
	\label{fig:laplacedarbouxgeometry}
\end{figure}

\begin{definition}
	The \emph{Laplace transform} $\Delta^k(q)$ of a Q-net $q: V(\Z^2) \rightarrow \CP^n$ is the new Q-net $\Delta^k(q): F(\Z^2) \rightarrow \CP^n$ such that
	\begin{align}
		(\Delta^k(q))^{kl} = \fp^{kl}
	\end{align}
	holds everywhere, where $\{k,l\}=\{1,2\}$ and $\fp^{kl}$ are the focal points of $q$ as in Definition \ref{def:qgqnet}.
\end{definition}

Consider a quad of the Laplace transform $\Delta^k(q)$ as shown in Figure \ref{fig:laplacedarbouxgeometry}. We observe that
\begin{align}
	\fp^{kl}\fp^{kl}_l \cap \fp^{kl}_k\fp^{kl}_{kl} = q_{kl}\label{eq:laplacequad}
\end{align}
holds. Thus the quad is planar and therefore the Laplace transform is indeed a Q-net again, as stated in the definition. Therefore it is possible to iterate the Laplace transform. In fact, Equation \eqref{eq:laplacequad} shows that 
\begin{align}
	\Delta^1 \circ \Delta^2 = \Delta^2 \circ \Delta^1 = \mbox{id}\label{eq:laplacetrafoinverse}
\end{align}
holds.

\begin{definition}
	\emph{Laplace-Darboux dynamics} is the iteration of the Laplace transform $\Delta^1$ of a Q-net defined on $\Z^2$.
\end{definition}

Thus Laplace-Darboux dynamics generate an infinite sequence of Q-nets, subsequently related by $\Delta^1$. On the other hand, $\Delta^2$ is the reverse iteration, thus we can also think of Laplace-Darboux dynamics as generating a biinfinite sequence of Q-nets.

\begin{figure}
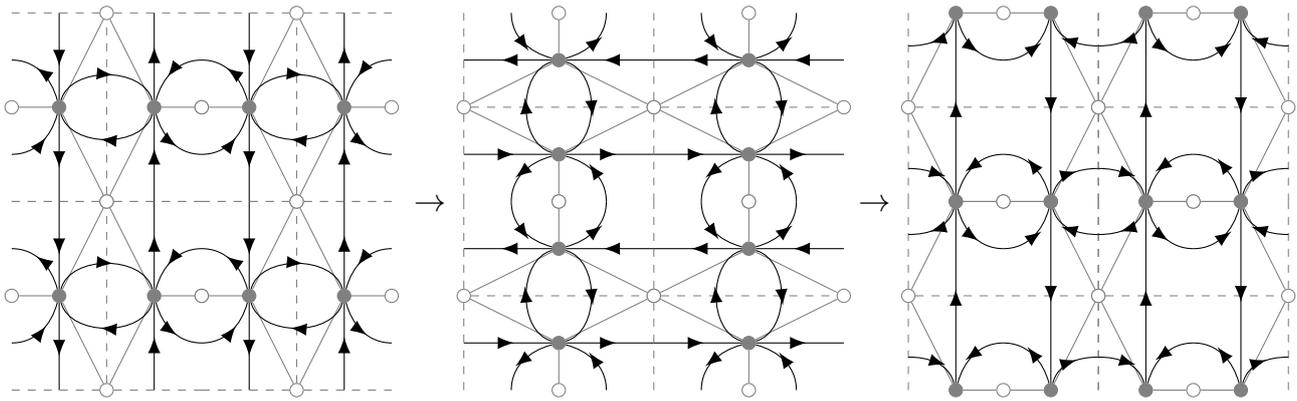
	
	
	\caption{The TCD of a Q-net under Laplace-Darboux dynamics in two steps.}
	\label{fig:tcdlaplacedarboux}
\end{figure}

From Section \ref{sec:qnets} we know that we can associate a TCD map to every Q-net. Thus it is natural to ask whether Laplace-Darboux dynamics correspond to a sequence of 2-2 moves. Indeed, the answer is yes and the proof is simply an illustration of the moves, as shown in Figure \ref{fig:tcdlaplacedarboux}. Let us give a short description of the sequence. Assume the initial TCD map captures $(q,\Delta^1(q))$, that is the points of the Q-net $q$ as well as the focal points $\fp^{12}$. On the level of the TCD one step of Laplace-Darboux dynamics consists of two steps:
\begin{enumerate}
	\item The first step consists of applying independent spider moves wherever possible. Then we can reinterpret the TCD map data as consisting of the data $(\Delta^1(q), \Delta^2 \circ \Delta^1(q))$, that is the points of $\Delta^1(q)$ as well as the focal points $\fp^{21}$ of $\Delta^1(q)$.
	\item The second step consists of applying independent resplits wherever possible. These resplits switch the focal points in all quads, such that the new TCD map data consists of $(\Delta^1(q), \Delta^1 \circ \Delta^1(q))$.
\end{enumerate}
We can then repeat these two steps to iterate Laplace-Darboux dynamics on the TCD map level.

In fact, one can also study Laplace and the related focal transforms for higher dimensional lattices $\Z^N$. Moreover, analogous definitions cannot only be given for Q-nets, but also for line complexes and Darboux maps. We study the details of these definitions and constructions in Section \ref{sec:focalnets}

The lattice generated by Laplace-Darboux dynamics also appears under the name of Menelaus-Cox lattice as a type of degenerate Cox lattice \cite{kscox}. We will take a closer look at non-degenerate Cox lattices in Section \ref{sub:cox}.

\chapter{Geometry and combinatorics of TCD maps}\label{cha:geometry}

\section{Propagation in non-general dimensions} \label{sec:nongendimension}

Maps that are defined in ambient spaces with a dimension that is not large enough for general position have been occasionally studied, and are usually defined as projections of the corresponding maps in higher dimensions \cite{ddgbook, bobenkoschieflinecomplexes}. To construct propagation in too small dimension but larger than 1, one can find case specific incidence relations. In dimension 1 this fails, but there have been approaches that characterize propagation via case specific multi-ratio equations \cite{ksclifford, schieflattice}. On the other hand, we showed how to describe the propagation of Q-nets, Darboux maps and line complexes via sequences of 2-2 moves. We did choose a definition of the 2-2 move that works in $\CP^1$ as well, see Definition \ref{def:tcdmapmoves}. Thus in the formulation via TCD maps, there is no obstacle to defining propagation in all these examples in small dimensions including dimension 1. In this way, we will show that the specific multi-ratio equations that have been studied in the literature arise naturally along strands of the TCDs. On the geometric side, these multi-ratios are closely related to the following generalization of Menelaus' theorem \cite[Theorem 9.12]{ddgbook}.

\begin{theorem}[Generalized Menelaus' theorem]\label{th:genmenelaus}
	Let $p_1,p_{1,2},p_2,\dots, p_m,p_{m,1}$ be $2m$ points in $\CP^{m-1}$ such that every point $p_{k,k+1}$ is on the line $p_kp_{k+1}$ and such that $p_1,p_2,\dots,p_m$ span $\CP^{m-1}$. Then the multi-ratio equation
	\begin{align}
		\mr(p_1,p_{12},p_2,\dots, p_m,p_{m1}) = (-1)^m
	\end{align}
	holds if and only if the $m$ points $p_{1,2},p_{2,3},\dots,p_{m,1}$ lie in an $(m-2)$ dimensional subspace.
\end{theorem}

In the case $m=3$ this clearly specializes to Menelaus' theorem (Theorem \ref{th:menelaus}). This theorem is important because, as we stated in Lemma \ref{lem:projinvariants}, multi-ratios are invariant under projective transformations and, what is important for our case, also under projection to $\CP^1$. However, the multi-ratios that occur as lattice equations require a version of Menelaus' theorem that has fewer assumptions but states only an implication, not an equivalence.

\begin{lemma}[Menelaus' lemma]\label{lem:menelaus}
	Let $p_1,p_{1,2},p_2,\dots, p_m,p_{m,1}$ be $2m$ points in $\CP^{m-1}$ such that every point $p_{k,k+1}$ is on the line $p_kp_{k+1}$, and such that $p_1,p_2,\dots,p_m$ span an $n$ dimensional space and $p_{1,2},p_{2,3},\dots,p_{m,1}$ span an $(n-1)$ dimensional space. Then the multi-ratio equation
	\begin{align}
		\mr(p_1,p_{1,2},p_2,\dots, p_m,p_{m,1}) = (-1)^m
	\end{align}
	holds.
\end{lemma}
\proof{
	Choose an affine chart where the space spanned by $p_{1,2},p_{2,3},\dots,p_{m,1}$ is at infinity. Then all the oriented length ratios cancel and only a factor $(-1)$ per edge survives.\qed
}

Another property that we observe is that the multi-ratios along zig-zag paths (see Definition \ref{def:zigzagtcd}) of a TCD map are invariant under 2-2 moves.

\begin{theorem}\label{th:zigzagratiodisc}
	Let $T$ be a TCD map and $(w_1,b_1,w_2,b_2,\dots,w_{n})$ be a zig-zag path. Denote the neighbours of each black vertex $b_k$ in counterclockwise order by $(w_k,w'_k,w_{k+1})$.
	Let $\tilde T$ be $T$ after a 2-2 move and $(\tilde w_1,\tilde b_1,\tilde w_2,\tilde b_2,\dots,\tilde w_{\tilde n})$ be the points of the same zig-zag path after the 2-2 move. Then
	\begin{align}
		&(-1)^n \mr(T(w_1),T(w'_1),T(w_2),T(w'_2),\dots,T(w_n)) \nonumber \\= &(-1)^{\tilde n} \mr(\tilde T(\tilde w_1),\tilde T(\tilde w'_1),\tilde T(\tilde w_2), \tilde T(\tilde w'_2), \dots, \tilde T(\tilde w_{\tilde n}))
	\end{align}
	holds.
\end{theorem}
\proof{We can calculate the multi-ratios in terms of edge-weights via Lemma \ref{lem:edgeratios}. The changes of edge-weights are given in Figure \ref{fig:vrclocalmoves} and Figure \ref{fig:resplit}. They are the same before and after a 2-2 move except for a change of sign if the length of the path changes.\qed
}

An interesting consequence is that we can consider two TCD maps related by 2-2 moves and glue their TCDs together to form a sphere. 

\begin{figure}
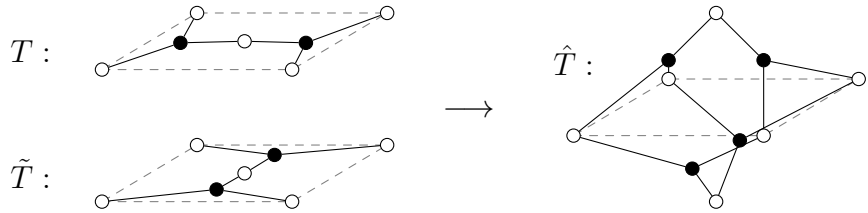



\caption{Gluing the bipartite graphs of $T$ and $\tilde T$ along the boundary (dashed) to obtain $\hat T$.}
\label{fig:glue}
\end{figure}

\begin{corollary}\label{cor:zigzagratios}
	Let $T$ and $\tilde T$ be two TCD maps related by a sequence of 2-2 moves. Denote by $\hat T$ the map that is obtained by gluing the two TCD together along their boundaries, such that we identify corresponding boundary white vertices in $\pb$ resp. $\tilde \pb$ (see Figure \ref{fig:glue}). Let $(w_1,w'_1,w_2,w'_2,\dots,w'_n,w_1)$ be the vertices on and incident to a zig-zag path (as in Theorem \ref{th:zigzagratiodisc}) of $\hat T$. Then
	\begin{align}
		\mr(\hat T(w_1),\hat T(w'_1),\hat T(w_2),\hat T(w'_2),\dots,\hat T(w'_n),\hat T(w_1)) = (-1)^n
	\end{align}	
	holds.
\end{corollary}
\proof{This is a direct consequence of Theorem \ref{th:zigzagratiodisc}, because the multi-ratios along zig-zag paths of $\hat T$ are ratios of multi-ratios of the same zig-zag paths in $T$ and $\tilde T$.\qed}

Note that $\hat T$ is technically not a TCD map because the glued TCD $\hat\tcd$ lives on a sphere and not on a disc. The advantage is that we can understand the sphere as an octahedron or cube. From that point of view, we will interpret the multi-ratio equations that occur in Corollary \ref{cor:zigzagratios} as lattice equations. We show that some of the lattice multi-ratio equations in the literature are indeed special cases of this construction.

\subsection{Q-nets and Laplace-Darboux dynamics}\label{sec:degdskp}
\begin{figure}
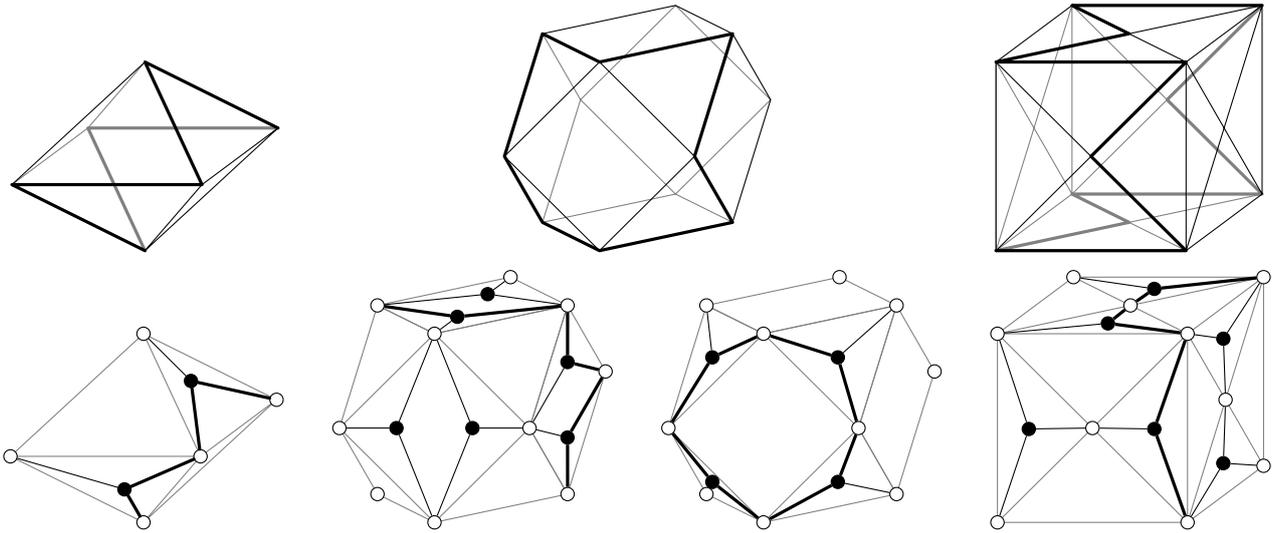


\caption{Top: The octahedron, cuboctahedron and the hexahedron with a highlighted zig-zag path each. Bottom: The corresponding (visible) parts of the TCDs of $\Z^2$ Q-net, Darboux map, line complex and $\Z^3$ Q-net.}
\label{fig:polyzigzag}
\end{figure}

The Laplace-Darboux dynamics of Q-nets with $\Z^2$ combinatorics are well defined in $\CP^2$, therefore the only non-trivial dimension is 1. Whenever we change from one focal point to another as part of the Laplace-Darboux dynamics, we are locally looking at a Menelaus configuration. If we glue the TCD maps before and after the resplit together, we obtain the octahedron (see Figures \ref{fig:glue} and \ref{fig:polyzigzag}). The dSKP equation (Equation \eqref{eq:resplitmr}) is then precisely the multi-ratio equation as occurring in Corollary \ref{cor:zigzagratios}. Lattices that satisfy the dSKP equation in dimension 1 are called dSKP lattices \cite{bkkphierarchy,bkkphierarchymulti,ksclifford}. Note that a dSKP lattice corresponds to the simultaneous projection of vertices of a Q-net and one of its focal points per quad. While all the focal points in $\CP^n$ for $n > 1$ are determined by the vertices, in $\CP^1$ we may actually choose one focal point per quad in the Cauchy data freely. The other focal point per quad is then determined via the dSKP equation. 

The relation between Q-nets and dSKP lattices will be very useful in our study of circular nets in $\CP^1$, together with other ``dimensionally degenerate'' constructions, see Section \ref{sec:cptemb}.

\subsection{Darboux maps}\label{sec:degendm}
The propagation, that is the cube-flip, of Darboux maps is well defined in $\CP^2$ and above. The case of $\CP^1$ has been studied by Schief \cite{schieflattice}. In particular, they satisfy and are characterized by the multi-ratio equation
\begin{align}
	\mr(d^3,d^1,d^3_1,d^2_1,d^3_{12},d^1_2,d^3_2,d^2) = 1 \label{eq:darbouxmr}.
\end{align}
On the one hand, this multi-ratio is an instance of Menelaus' Lemma \ref{lem:weakmenelaus}. This is because we can consider $d^3,d_1^3,d_{12}^3,d^3_2$ as a 4-gon that spans a 2-space and $d^2,d^2_1,d^1_2,d^2$ are points on the edges that only span a line. In $\CP^1$ this equation determines $d^3_{12}$ if we know the other seven points and thus defines propagation of the initial data, together with permutations of the equation that determine $d^2_{13}$ and $d^1_{23}$.

On the other hand, Equation \eqref{eq:darbouxmr} is the multi-ratio along a zig-zag path in the cuboctahedron (see Figure \ref{fig:polyzigzag}). The vertices of the cuboctahedron are the edges of a cube. Let $\hat T$ be the TCD map that we obtain by gluing the two TCD maps before and after a cube flip in a Darboux map together (see Figure \ref{fig:vrcvstcd} with reversed orientations of the strands). Then the zig-zag paths of $\hat T$ and those of the cuboctahedron coincide. We explain this phenomenon in more detail when we introduce the affine quiver of a TCD map in Section \ref{sec:affcluster}.

We will encounter $\CP^1$ Darboux maps as h-embeddings used in statistical physics in Section \ref{sec:harmonicemb}.

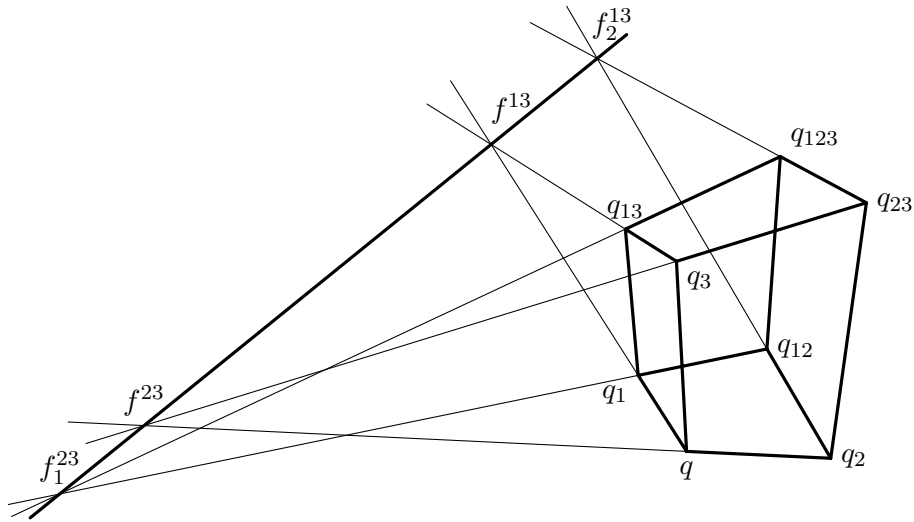
\begin{figure}
\begin{tikzpicture}[scale=.7]
	\small
	\useasboundingbox (-13,-4.3) rectangle (4.5,5.7);
	\coordinate[label=-90:$q$] (q) at (-0.25,-2.91);
	\coordinate[label=225:$q_1$] (q1) at (-1.16,-1.47);
	\coordinate[label=0:$q_2$] (q2) at (2.46,-3.04);
	\coordinate[label=330:$q_3$] (q3) at (-0.44,0.69);
	\coordinate[label=0:$q_{12}$] (q12) at (1.26,-0.97);
	\coordinate[label=90:$q_{13}$] (q13) at (-1.4,1.3);
	\coordinate[label={0:$q_{23}$}] (q23) at (3.13,1.8);
	\coordinate[label={[yshift=4,xshift=-4]85:$f^{13}$}] (f13) at (intersection of q--q1 and q3--q13);
	\coordinate[label=90:$f^{23}$] (f23) at (intersection of q--q2 and q3--q23);
	\coordinate[label={[xshift=-6,yshift=3]85:$f^{13}_2$}] (f132) at (intersection of f13--f23 and q2--q12);
	\coordinate[label={[yshift=0,xshift=-0]90:$f^{23}_1$}] (f231) at (intersection of f13--f23 and q1--q12);
	\coordinate[label=80:$q_{123}$] (q123) at (intersection of f231--q13 and f132--q23);
	\draw[-,very thick,line cap=round, line join=round]
		(q) -- (q1) -- (q12) -- (q2) -- (q)
		(q3) -- (q13) -- (q123) -- (q23) -- (q3)
		(q) -- (q3) (q1) -- (q13) (q2) -- (q23) (q12) -- (q123)
		(f132) edge[shorten >=-0.5cm, shorten <=-0.5cm] (f231)
	;	
	\draw[-,line cap=round, line join=round]
		(q) edge[shorten >=-1cm] (f23)
		(q1) edge[shorten >=-1cm] (f13) edge[shorten >=-1cm] (f231)
		(q3) edge[shorten >=-.8cm] (f23)
		(q12) edge[shorten >=-.8cm] (f132)
		(q13) edge[shorten >=-1cm] (f13) edge[shorten >=-0.7cm] (f231)
		(q123) edge[shorten >=-1cm] (f132)
	;
\end{tikzpicture}
\vspace{-2mm}
\caption{Projection of a cube in a Q-net to a plane.}
\label{fig:qnetprojection}
\end{figure}

\subsection{Q-nets}
The propagation, that is the cube-flip, of Q-nets is well defined in $\CP^3$. In $\CP^2$ a possible definition is given in the DDG book \cite[Exercise 2.2]{ddgbook}. In particular, the requirement is that the four points 
\begin{align}
	\fp^{ik}, \fp^{jk}, \fp^{ik}_j, \fp^{jk}_i
\end{align}
are colinear (see Figure \ref{fig:qnetprojection}) for any set of indices such that $\{i,j,k\}=\{1,2,3\}$. In $\CP^3$ these points are colinear as well because they all lie on the intersection line of the planes $p^{ij},  p^{ij}_k$. If we expect the Q-net in $\CP^2$ to be a projection of a Q-net in $\CP^3$, then this colinearity property is the necessary and sufficient condition. It is sufficient that this relation is satisfied for one choice of indices $\{i,j,k\}=\{1,2,3\}$. This property determines propagation in $\CP^2$.

In $\CP^1$ we can of course not give a characterization in terms of colinearity. We already know that propagation is determined uniquely via a sequence of 2-2 moves, because we can view the Q-net as a TCD map. It is unclear if there is a multi-ratio characterization in the spirit of the previous section. However, we can observe that the following multi-ratio equation
\begin{align}
	\mr(q,q_1,\fp^{12},q_{12},q_2,\fp^{31}_2,q_{23},q_{123},\fp^{12}_3,q_{13},q_3,\fp^{31}) = 1, \label{eq:mrhexahedron}
\end{align}
is satisfied. Again, this is an instance of Menelaus' lemma (Lemma \ref{th:weakmenelaus}). In $\CP^3$ and higher-dimensional spaces the points $q,\fp^{12},q_2,q_{23},\fp^{12}_3,q_3$ build a hexagon and the points $q_1,q_{12},\fp^{31}_2,q_{123},q_{13},\fp^{31}$ are in a plane. However, unlike before Equation \eqref{eq:mrhexahedron} involves 2 new points with respect to propagation. In order to propagate via the cube flip, we are looking for the four points $q_{123},\fp^{12}_3,\fp^{23}_1,\fp^{31}_2$. The permutations of Equation \eqref{eq:mrhexahedron} yield six trilinear equations for these points. By looking at the product of all the equations, we see that one of the equations is implied by the other 5. Thus we are left with 5 coupled bilinear equations for four variables. We know they have a solution but it is unclear whether the solution is unique. Numerical simulations indicate that the solution is unique, but we have not found a proof.

We also observe that Equation \eqref{eq:mrhexahedron} is the multi-ratio equation along a zig-zag path in the hexahedron (see Figure \ref{fig:polyzigzag}). The zig-zag paths of the hexahedron are exactly those that appear when we glue the TCD of a Q-net before and after the cube flip. 

We will encounter a special case of $\CP^1$ Q-nets called s-embeddings, known from statistical physics, in Section \ref{sec:sembeddings}.

\subsection{Line complexes}\label{sub:cponelc}
Line complexes are a priori only well defined in $\CP^4$ and above. Projections to $\CP^3$ have been studied by Bobenko and Schief \cite{bobenkoschieflinecomplexes}. They observed two possible equivalent characterizations of a line complex in $\CP^3$. Specifically that
\begin{enumerate}
	\item the four intersection points $l^{i},l^{i}_j,l^{i}_k,l^{i}_{jk}$ are coplanar,
	\item or the four lines $l^{i} l^{i}_{j},l^{i}_{k}l^{i}_{jk},l^{k}l^{k}_{j},l^{k}_{i}l^{k}_{ij}$ intersect in a point,
\end{enumerate}
for $\{i,j,k\} = \{1,2,3\}$. It turns out that the second property is also characterizing in $\CP^2$. Moreover, in $\CP^1$ we can characterize propagation as in the previous examples via a multi-ratio equation along a zig-zag path, which is that
\begin{align}
	\mr(l^3,l^1,l^3_1,l^2_1,l^3_{12},l^1_2,l^3_2,l^2) = 1. \label{eq:lcmr}
\end{align}
This is the same recurrence as in the case of $\CP^1$ Darboux maps. We also recognize this as an instance of the generalized Menelaus' theorem. In $\CP^4$ we can view $l^{1},l^{1}_2,l^{1}_3,l^{1}_{23}$ as a 4-gon that spans a 3-space. Then the multi-ratio equation \eqref{eq:lcmr} is equivalent to stating that $l^{1},l^{1}_2,l^{1}_3,l^{1}_{23}$ are on a plane, which is exactly characterization (1) from above. On the other hand, this multi-ratio equation can once again be viewed as induced by gluing together the TCDs from before and after the cube flip, this time for a line complex.

Note that in $\CP^1$, the combinatorics and the multi-ratio characterizations are identical for line complexes and Darboux maps. Therefore, there is no difference between line complexes and Darboux maps in $\CP^1$. One may thus view line complexes and Darboux maps as different higher dimensional generalizations of the same scalar system.

\section{Sections}\label{sec:sections}

In this section we consider what happens if we intersect a TCD map with a hyperplane (a codimension one projective subspace of the ambient space). More explicitly, we intersect all the lines corresponding to crossings in a TCD map with the hyperplane, which yields new points. These points satisfy new relations, and together they yield a new TCD map. Before we give a precise definition of a section of a TCD map, let us introduce some genericity conditions.

\begin{definition}\label{def:onegeneric}
	A TCD map $T: \tcdp \rightarrow \CP^n$ is \emph{1-generic} if for any three different white vertices $w_1,w_2,w_3$ that are on the boundary of a common face of $\pb$, the span of $T(w_1),T(w_2),T(w_3)$ is two dimensional.
\end{definition}
In a sense, the restriction that a TCD map $T$ is such that at any black vertex the images of the three adjacent vertices are different could be viewed as $T$ being 0-generic. We did incorporate the 0-genericity into the definition of TCD maps because otherwise the accompanying theory with edge-weights breaks down. We will define $k$-genericity in Definition \ref{def:kgeneric}. Note that a TCD map defined on a TCD that has maximal dimension 1 is always 1-generic, because there are no triplets of white vertices as in the definition.

\begin{definition}\label{def:generichyperplane}
	Let $T: \tcdp \rightarrow \CP^n$ be a TCD map. A hyperplane $H\subset \CP^n$ is called \emph{generic with respect to $T$} if $H$ does not contain the image $T(w)$ of any white vertex $w$ of $\pb$.
\end{definition}

The conditions for a generic hyperplane only exclude a closed codimension 1 subset of hyperplanes and thus there always exists a generic hyperplane.

\begin{figure}
	
	\caption{The sequence $\pb \rightarrow \pb^i \rightarrow \pb' = \sigma(\pb)$ or how to take a section of a TCD map.}
	\label{fig:vrcsection}
\end{figure}

\begin{definition}\label{def:section}
	Let $T: \tcdp \rightarrow \CP^n$ with  $n \geq 2$ be a 1-generic TCD map and $H$ a generic hyperplane of $\CP^n$ with respect to $T$. A \emph{section} $\sigma_H(T)$ of $T$ is a TCD map $T': \tcdp'\rightarrow H = \CP^{n-1}$, for which we give the construction below. We define the change of combinatorics in terms of $\pb$ and $\pb'$. We begin by constructing an intermediate graph $\pb^i$ starting from $\pb$ in two steps (see Figure \ref{fig:vrcsection}):
	\begin{enumerate}
		\item Add a black vertex $b_{k,k+1}$ for any two consecutive boundary vertices $w_k, w_{k+1}$ and the two edges $(b_{k,k+1}, w_k)$ and $(b_{k,k+1}, w_{k+1})$.
		\item For each face $(b_1,w_1,b_2,\dots,w_l)$ of $\pb$, triangulate the polygon $(w_1,w_2,\dots,w_l)$ and for each diagonal $(w_j,w_k)$ of the triangulation add a new black vertex $b$ as well as the two edges $(b,w_j),(b,w_k)$.
	\end{enumerate}
	Every interior face of $\pb^i$ is a quad or a hexagon. Now we construct $\pb'$ and $T'$ in the following steps:
	\begin{enumerate}
		\setcounter{enumi}{2}	
		\item Add a white vertex $w'_b$ to $\pb'$ for each black vertex $b$ of $\pb^i$. Put the corresponding point $T'(w'_b)$ at the intersection of $H$ with the line that corresponds to $b$ in $T$.
		\item For each face $f$ of $\pb^i$, add a new black vertex $b'_{f}$ to $\pb'$.
		\item Add edges between a white vertex $w_b'$ and a black vertex $b'_{f}$ in $\pb'$ whenever the corresponding black vertex $b$ and face $f$ are incident in $\pb^i$.
		\item Contract all black vertices of degree 2 in $\pb'$.\qedhere
	\end{enumerate}
\end{definition}

Let us make a few observations. If $f$ is a hexagon of $\pb^i$, then the three lines in $T$ that belong to that hexagon are in a common plane $J$, and indeed span a plane because $T$ is 1-generic. Thus in the intersection with $H$ the corresponding black vertex that we add to $\pb'$ represents the line that is the intersection of $J$ with $H$. Note that $J$ cannot be contained in $H$ because $H$ is generic with respect to $T$. If a face $f$ of $\pb^i$ is a quad, then the two black vertices of that face correspond to the same line $\ell$. The new black vertex that we add to $\pb'$ for $f$ has degree 2 and is contracted in step (6), therefore $\ell$ is represented by a single point in $\pb'$. Hence, performing 2-2 moves at faces of $\pb$ (reparametrizations) does not affect the section. Also note that the 1-genericity of $T$ guarantees that the images of no three white vertices adjacent to a black vertex in $\sigma(T)$ coincide. However, 1-genericity is a bit more than is necessary for the existence of a section. Instead, 1-genericity guarantees the existence of all sections.

The different triangulations that we can choose in step (2) lead to different sections that are related by resplits.

It is also a useful observation that each new face of $\pb'$ corresponds to a white vertex of $\pb$. 

Note that in subsequent parts of the thesis, whenever we take sections $\sigma_H(T)$ we implicitly assume that $T$ is 1-generic and $H$ is generic with respect to $T$.

We will now present another was of thinking about taking the section of a TCD $\tcd$. It is also possible to take sections of non-minimal TCD, but we restrict our attention to minimal TCD, as they are our main interest, and this restriction simplifies the exposition. We argue that we only need to perturb the strands of $\tcd$ to obtain $\sigma(\tcd)$, without the need to cut or glue any strands. We begin by introducing an alternating strand diagram  $\alt$ associated to $\tcd$. These were introduced by Postnikov \cite{postgrass} to study Grassmannians, but we forego a general definition of alternating strand diagrams as we are only interested in the particular alternating strand diagrams that appear in the next definition.

\begin{definition}\label{def:alt}
	Let $\tcd$ be a minimal TCD. Perturb every strand by a small $\epsilon > 0$ in its normal direction to the right (see Figure \ref{fig:tcdsection}), such that the triple crossings of $\tcd$ are resolved into three transversal crossings of pairs of strands, and such that no other crossings appear. Afterwards, for each strand $k$ move its in-endpoint $k$ such that it appears after out-endpoint $k$, which introduces one new crossing per strand. We call the resulting diagram the \emph{alternating strand diagram} $\alt(\tcd)$.
\end{definition}

\begin{definition}\label{def:shadowgraph}
	Let $\tcd$ be a minimal TCD and consider it as a graph. For every pair of in- and out-endpoint $k$ add a single boundary vertex $v_k$. Also add boundary edges $(v_k,v_{k+1})$ for every strand $k$ and denote the resulting \emph{shadow graph} by $\altg$ (see Figure \ref{fig:tcdsection}). Moreover, consider the superposition of $\altg$ and $\alt$. Every edge of $\altg$ is intersected exactly twice by strands of $\alt$. In each face $f$ of $\altg$ we consider the endpoint matching $C_\alt^f$ induced by the strands of $\alt$ in that face. Call the pair $(\altg,C_\alt)$ the \emph{shadow} of $\tcd$.
\end{definition}

In fact, the endpoint matching $C_\alt^f$ in a clockwise face $f$ of degree $k$ is always $\enm {k}{-2}$ and in a counterclockwise face it is always $\enm {k}{-1}$, see Figure \ref{fig:tcdsection}.

\begin{definition}\label{def:shadowtotcd}
	From a shadow $(\altg,C_\alt)$ we define a TCD $\tcd(\altg,C_\alt)$ by replacing every face $f$ in $\altg$ with a TCD that has endpoint matching $C^f_\alt$ and removing all strands of length zero.
\end{definition}
By definition, the different diagrams $\tcd(\altg,C_\alt)$ that are constructed in this way are related by a sequence of 2-2 moves.

\begin{lemma}\label{lem:sectionviatcd}
	Let $\tcd$ be a minimal TCD and $T:\tcdp \rightarrow \CP^n$ be a TCD map. Then the TCD $\sigma(\tcd)$ of any section $\sigma(T)$ corresponds to a TCD constructed from the shadow. Thus, one can make choices such that
	\begin{equation}
		\sigma(\tcd) = \tcd(\altg,C_\alt)
	\end{equation}
	holds.
\end{lemma}
\proof{
We compare the construction of a section on the level of TCDs to the construction of a section on the level of graphs in Definition \ref{def:section}. When splitting triple intersection points to obtain the alternating strand diagram $\alt(\tcd)$ we split all strands such that a new counterclockwise face emerges, this corresponds to the new white vertices that are placed in graph step (3) for black vertices that were already present in $\pb$ (and not introduced in $\pb^i$). The fact that we also swap in- and out-endpoints corresponds to the introduction of black vertices on the boundary edges in graph step (1) and their subsequent replacement by white vertices in graph step (3). Gluing new TCDs in each face of the shadow $(\altg,C_\alt)$  corresponds to the triangulation choices in graph step (2) and the subsequent introduction of white vertices in graph step (3) as well as black vertices in graph step (4). We remove strands of length 0 because they do not contribute information to the TCD map and do not appear in the TCD associated with $\sigma(\pb)$.\qed
}

\begin{figure}
	
	\caption{The sequence $\tcd \rightarrow \alt$ (black) and $\altg$ (gray) $\rightarrow\tcd(\altg,C_\alt)=\sigma(\tcd)$ or how to take a section of a TCD.}
	\label{fig:tcdsection}
\end{figure}

If we do not remove strands of length 0 from $\sigma(\tcd)$, then there is a simple formula for the endpoint matching of a section, indeed
\begin{align}
	C_{\sigma(\tcd)}(k) = C_\tcd(k) - 1
\end{align}
holds. Note that the strands of length 0 that we removed in $\sigma(\tcd)$ are strands in $\tcd$ of length 1, these have only a single white vertex to their right in $\tcd$.

\begin{corollary}\label{cor:maxdimsection}
	The maximal dimension of $\sigma(\tcd)$ is one less than the maximal dimension of $\tcd$.
\end{corollary}
\proof{
	Let $n$ be the number of strands of $\tcd$. First of all, if every strand is of length 1 then the endpoint matching is $\enm {n}{1}$ and thus corresponds to the triangulation of an $n$-gon (see Example \ref{ex:tcdtriangulation}). The TCD map belonging to an $\enm {n}{1}$ TCD is also a triangulation (see Section \ref{sec:extriangulation}) and its maximal dimension is 1. After taking a section, the TCD consists only of strands of length 0 and there is only one counterclockwise face and the maximal dimension is 0. In every other case, there is a strand of length more than 1. Choose a labeling where such a strand connects to out-endpoint $1$ and consider a half-plane drawing. In order to determine the maximal dimension of the section we use Theorem \ref{th:maxdim}, which states that the maximal dimension equals the number of left moving strands. All strands of length 1 in $\tcd$ are right moving and thus their removal does not change the maximal dimension. The strand ending at out-endpoint $1$ is necessarily left moving. In the section out-endpoint 1 gets swapped with in-endpoint $1$ and we choose a new half-plane drawing where the former out-endpoint $1$ is now out-endpoint $n$. The strand ending at the former out-endpoint $1$ is now right moving. All other left (resp. right) moving strands stay left (right) moving. Therefore the maximal dimension of $\sigma(\tcd)$ is exactly the maximal dimension of $\tcd$ minus one.\qed
}

\begin{lemma}\label{lem:sectionminimal}
	Let $\tcd$ be a minimal TCD. Then $\sigma(\tcd)$ is connected and minimal.
\end{lemma}
\proof{By construction of $\sigma(\tcd)$, no strands are cut and thus all strands remain connected to the boundary, thus $\sigma(\tcd)$ is connected. For minimality, we show first that there is no pair of strands that intersect twice in a parallel manner and then that there are no strands that self-intersect. Then the conditions of Theorem \ref{th:mintcdforbidden} are satisfied which will conclude the proof of the Lemma. 
	
Let us look at what happens when we introduce new crossings in the section $\sigma(\tcd)$. In each face $f$ of $\altg$ we rearrange the strands to yield a new TCD $\tcd^f$ with connectivity $C_{\tcd^f}(k) = (k-2)$. Observe that if two strands intersected once at the corresponding face $f'$ in $\tcd$ then they intersect also once in $\tcd^f$. If they did not intersect in $f'$ then they either do not intersect in $\tcd^f$ or they intersect twice while bounding a counterclockwise bigon.

Let $k,l$ be two strands in $\tcd$ and let $b_1 <_k b_2 <_k \dots <_k b_n$ be their intersection points ordered along the orientation of strand $k$. Because $\tcd$ is minimal, the intersection points are in reverse order along strand $l$, that is $b_1 >_l b_2 >_l \dots >_l b_n$, as else there would be a parallel bigon in $\tcd$. In $\sigma(\tcd)$ we only add new intersections points $b',b''$ in pairs such that:
\begin{align}
	b_1 <_k \dots <_k b_n\quad &\longrightarrow\quad b_1 <_k \dots <_k b' <_k b'' <_k \dots <_k b_n\\
	b_1 >_l \hspace{0.7mm} \dots >_l \hspace{0.7mm} b_n\quad &\longrightarrow\quad b_1 >_l\hspace{0.7mm} \dots >_l \hspace{0.7mm} b' >_l \hspace{0.7mm} b'' >_l \hspace{0.7mm} \dots >_l\hspace{0.7mm} b_n
\end{align}
Therefore no parallel intersections appear in $\sigma(\tcd)$.

Moreover, we argue that as $\tcd$ is minimal it is not possible that a clockwise face $f$ in $\tcd$ is bounded twice by the same strand. Assume strand $k$ bounds the face $f$ between intersection points $b_1,b_2$ and also between $b_i,b_{i+1}$ such that:
\begin{enumerate}
	\item $f$ is to the left of the strand $k$ between $b_2$ and $b_i$,
	\item the strands bounding $f$ between $b_j,b_{j+1}$ for $1 < j <i$ are distinct.
\end{enumerate}
Of course $2<i<n-1$ because else there is already a self-intersection. But then the strand that bounds the face between $b_2,b_3$ has to either intersect itself in $b_2$ or intersect strand $k$ while forming a parallel intersection. Because therefore all strands around a clockwise face in $\tcd$ are distinct, we are not introducing any self-intersections in $\sigma(\tcd)$.\qed
}

\begin{figure}
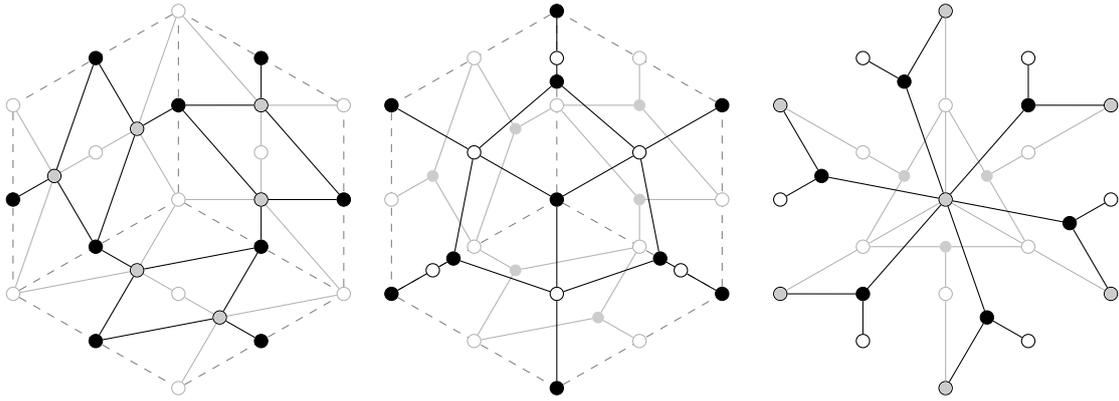


	\caption{The graph $\pb$ in gray and $\sigma(\pb)$ in black for the fundamental domains of stepped surfaces of Q-net, Darboux map and line complex.}
	\label{fig:tcdsectionexamples}
\end{figure}

\begin{figure}
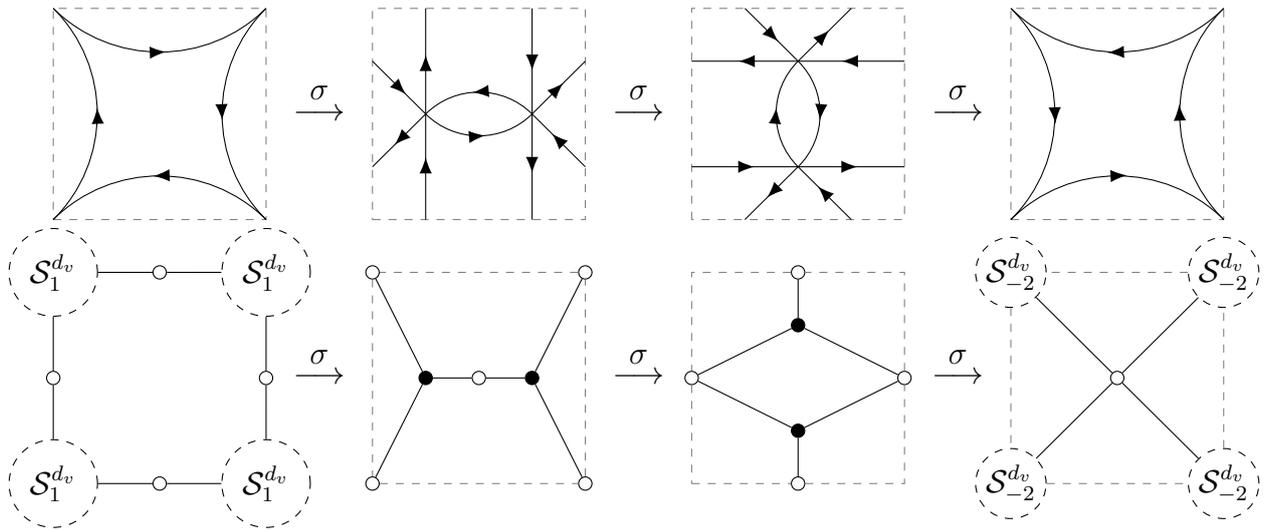

	
	\caption{The TCD and bipartite graph pieces that we glue into the quads of a line complex, a Q-net, a Darboux map, and finally a new type of map.}
	\label{fig:qgtcd}
\end{figure}

Our running examples from discrete differential geometry are actually related via sections.
\begin{lemma}\label{ex:ddgsteppedsections}
	Consider the TCD maps associated to stepped surfaces of a Q-net, Darboux map, line complex. Then sections can be chosen such that
	\begin{enumerate}
		\item the section of a Q-net is a Darboux map,
		\item the section of a Darboux map is a line complex,
		\item the section of a line complex is a Q-net.
	\end{enumerate}
	Moreover, the section of a Q-net defined on $\Z^2$ is a Q-net defined on (a translated copy of) $\Z^2$ again. If we consider maps defined on quad-graphs, then instead one can make choices such that
	\begin{enumerate}
		\item the section of a line complex is a Q-net,
		\item the section of a Q-net is a Darboux map,
		\item the section of a Darboux map is a line compound.\qedhere
	\end{enumerate}	
\end{lemma}
\proof{Follows from applying the section rules in Definition \ref{def:section}. For a visualization of the stepped surface statements, see Figure \ref{fig:tcdsectionexamples}, and for a visualization of the quad-graph statements, see Figure \ref{fig:qgtcd}.\qed\\}

We note that the first three statements above define a cyclic relation with three members. Stepped surface statement (1) is the subject of Exercises 2.8 and 2.9 in the DDG book \cite{ddgbook}. The other statements are to the best of our knowledge new, although they are all very straight-forward from the definitions.

We can also take \emph{iterated sections}, and they commute.

\begin{lemma}\label{lem:sectionscommute}
	Let $T: \tcdp \rightarrow \CP^n, n \geq 2$ be a 1-generic TCD map. Let $H_1,H_2$ be two different generic hyperplanes with respect to $T$ and let $H_{12}= H_1\cap H_2$. Let $T_1= \sigma_{H_1}(T),T_2= \sigma_{H_2}(T)$ be two sections of $T$, both defined on a TCD $\tcd'$, and choose a TCD $\tcd'' = \sigma (\tcd')$. If $T_1,T_2$ are 1-generic, $H_{12}$ is generic with respect to both $T_1$ and $T_2$, then there is a unique TCD map $T_{12}:\tcdp'' \rightarrow H_{12}$ such that
	\begin{equation}
		T_{12} = \sigma_{H_2}(T_1) = \sigma_{H_1}(T_2).\qedhere
	\end{equation}

\end{lemma}
\proof{Once the combinatorics of a section are chosen, the points of the section are determined. Thus we have to check that $\sigma_{H_2}(T_1) = \sigma_{H_1}(T_2)$. A white vertex $w$ in $\sigma_{H_2}(T_1)$ corresponds to
	\begin{align}
		H_{2}\cap H_{1} \cap H_w = H_{12}\cap H_w,
	\end{align}
	where $H_w$ is some 2-dimensional space spanned by points of $T$. The same white vertex $w$ in $\sigma_{H_1}(T_2)$ however corresponds to 
	\begin{align}
		H_{1}\cap H_{2} \cap H_w = H_{12}\cap H_w.
	\end{align}
	Therefore both, the combinatorics and the geometry agree, and the lemma is proven.\qed
}

As a consequence of Lemma \ref{lem:sectionscommute}, it makes sense to denote
\begin{align}
	\sigma_{H_{12}}(T) = \sigma_{H_2}(T_1) = \sigma_{H_1}(T_2)
\end{align}
and to call $\sigma_{H_{12}}(T)$ a codimension 2-section of $T$. More generally, for $\dim H' = n-k$, we can say that $\sigma_{H'}(T)$ is a codimension $k$-section of $T$. 

Now that we can take iterated sections in a meaningful way, let us also give a definition of ``iterated genericity''.
\begin{definition}\label{def:kgeneric}
	A TCD map $T: \tcdp \rightarrow \CP^n$ is \emph{$k$-generic} for $n>k>1$, if $T$ is 1-generic and every section $\sigma_H(T)$ is $(k-1)$-generic. Moreover, a subspace $S$ of codimension $k$ is \emph{generic} if it is generic for all sections $\sigma_H(T)$ with $H \supset S$.
\end{definition}

Let us discuss generic subspaces in $\CP^n$ first. Essentially, a generic $(n-k)$ space is a space that does not intersect any of the $(k-1)$-spaces spanned by $T$ that occur as points in a $k$-th section of $T$. In general, if $\pb$ has $m$ white vertices, the set of forbidden spaces is a subset of all $(k-1)$ spaces spanned by $k$ image points of $T$. Therefore generic $k$ spaces exist for any $k$ with $0\leq k < n-1$. It is not completely obvious whether $k$-generic TCD maps exist for any given $\tcd$. We leave the proof of existence (by examples) to the reader. We make the following suggestion: Consider a TCD $\tcd$ with endpoint matching $\enm{n}{-1}$. The corresponding graph $\pb$ consists of $n$ white vertices and no black vertices, the maximal dimension is clearly $n-1$. Consider a TCD map on $\tcd$ that attains maximal dimension. We claim it is $(n-2)$-generic.

\begin{lemma}
	Let $\tcd, \tilde \tcd$ be two minimal TCD related by a sequence of 2-2 moves. Let $T,\tilde T$ be two TCD maps defined on $\tcd, \tilde \tcd$. Then there is a sequence of 2-2 moves that takes $\sigma_H(T)$ to $\sigma_H(\tilde T)$.
\end{lemma}
\proof{
	The endpoint matchings of the TCDs $\tcd$ and $\tilde{\tcd}$ agree. 
	Therefore the endpoint matchings of the sections agree as well. By Theorem \ref{th:tcdflipsconnected}, there exists a sequence of 2-2 moves that takes the TCD of one of the sections into the other one. \qed
}

\section{Sweeps and acyclic orientations}\label{sec:sweeps}

Recall that to every $\tcd$ we associate the graph $\pb$, on which we defined the associated vector-relation configuration. In this section and later on, it turns out to be useful to have a second, related graph.
	
\begin{definition}\label{def:pbm}
	Let $\tcd$ be a TCD and $\pb$ the associated graph. The graph $\pbm$ has a black vertex $b$ for every black vertex $b$ of $\pb$, a white vertex $w$ for every face $f$ of $\pb$ and an edge $(w,b)$ whenever $f$ and $b$ are incident in $\pb$.	
\end{definition}

By definition, $\pbm$ is a planar bipartite graph and all black vertices of $\pbm$ are of degree three. Indeed, if $\iota(\tcd)$ is the TCD $\tcd$ with all strand orientations reversed, then $\pbm$ is the associated graph for $\iota(\tcd)$.

\begin{figure}
	\begin{tikzpicture}
		\coordinate[label=left:$a_1$] (a1) at (240:2);
		\coordinate (aa1) at (60:2);
		\coordinate[label=below:$a_2$] (a2) at (0:2);
		\coordinate (aa2) at (180:2);
		\coordinate[label=right:$a_3$] (a3) at (120:2);
		\coordinate (aa3) at (300:2);
		\draw[]
			(a1) edge[-latex] (aa1)
			(a2) edge[-latex] (aa2)
			(a3) edge[-latex] (aa3)
		;
		\node[bvert] (b) at (0,0) {};
		\node[wvert] (w1) at (330:1.25) {};
		\node[wvert] (w2) at (90:1.25) {};
		\node[wvert] (w3) at (210:1.25) {};
		\draw[-]
			(b) edge[orient] (w1) edge[rorient] (w2) edge[orient] (w3)
		;
	\end{tikzpicture}\hspace{1.5cm}
	\begin{tikzpicture}
		\coordinate[label=left:$a_1$] (a1) at (240:2);
		\coordinate (aa1) at (60:2);
		\coordinate[label=below:$a_2$] (a2) at (0:2);
		\coordinate (aa2) at (180:2);
		\coordinate[label=right:$a_3$] (a3) at (120:2);
		\coordinate (aa3) at (300:2);
		\draw[]
			(a1) edge[-latex] (aa1)
			(a2) edge[-latex] (aa2)
			(a3) edge[-latex] (aa3)
		;
		\node[bvert] (b) at (0,0) {};
		\node[wvert] (w13) at (270:1.25) {};
		\node[wvert] (w12) at (30:1.25) {};
		\node[wvert] (w23) at (150:1.25) {};
		\draw[-]
			(b) edge[rorient] (w12) edge[orient] (w13) edge[rorient] (w23)
		;
	\end{tikzpicture}
	\caption{The li-orientation at a black vertex for $\pb$ on the left and for $\pbm$ on the right, with $a_1<a_2<a_3$.}
	\label{fig:liorientatblack}
\end{figure}
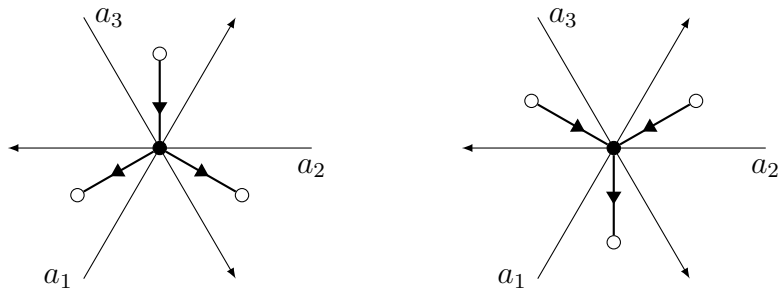

\begin{lemma}\label{lem:strandorder}
	Let $\tcd$ be a labeled minimal TCD.
	\begin{enumerate}
		\item At every intersection point $b$, there are three strands $a_1,a_2,a_3$ in counterclockwise order. If $a_1$ carries the smallest label, then $a_1 < a_2 < a_3$.
		\item At every face or boundary face $f$ there are $d_f$ strands $a_1,a_2,\dots, a_{d_f}$ in counterclockwise order. If $a_1$ carries the smallest label then $a_1 < a_2 < \dots < a_{d_f}$.\qedhere
	\end{enumerate}
\end{lemma}
\proof{
	Let us denote by $a_1',a_2',a_3'$ the part of the strands $a_1,a_2,a_3$ from the boundary to the intersection point $b$. Because of minimality, $a_1',a_2',a_3'$ can only intersect at $b$. Therefore, the first part of the lemma follows. For the second part, assume that $a'_i$ and $a'_j$ intersect before reaching the face $f$. Then the two strands $a_i$ and $a_j$ have to be non-consecutive on $f$ or the TCD is not minimal. Denote by $a''_i, a''_j$ the parts of $a_i,a_j$ from $f$ to the boundary. Then either $a''_i$ or $a''_j$ is locked in by $f, a'_i$ and $a'_j$. This is not possible in a minimal TCD and by contradiction the lemma follows.\qed	
	
}

\begin{definition}\label{def:liorient}
	Let $\tcd$ be a labeled TCD. Each edge $e$ of $\pb$ is adjacent to a unique black vertex $b$. At $b$ there are three strands $a_1< a_2 < a_3$. Orient $e$ away from $b$ if $e$ is in between $a_1$ and $a_3$ and towards $b$ in any other case. This defines the \emph{labeling induced orientation} $\lio(\tcd)$ on $\pb$, which we abbreviate by \emph{li-orientation}. Analogously, in $\pbm$ orient every edge $e$ towards $b$ if $e$ is in between $a_1$ and $a_3$ and away from $b$ else. This defines the \emph{li-orientation} $\liom(\tcd)$ on $\pbm$, see also Figure \ref{fig:liorientatblack}.
\end{definition}

In fact, the li-orientation has previously appeared for special diagrams, the so called \emph{monotone diagrams} \cite{postgrass}. For any permutation (and a labeling choice of the strands) Postnikov gives an algorithm how to construct a TCD with the given endpoint matching. This construction is accompanied by an orientation $\lio'$ of a graph that is $\pbm$ minus some edges at the boundary. It is an easy verification that $\lio'$ coincides with $\lio$.

\begin{lemma}\label{lem:lioperfect}
	Let $\tcd$ be a minimal labeled TCD and let $\lio$ (resp. $\liom$) be the li-orientation of $\pb$ (resp. $\pbm$). Then every white vertex has exactly one incoming (resp. outgoing) edge. Moreover, for each face $f$ there is exactly one incoming (resp. outgoing) edge. An incoming (resp. outgoing) edge is an edge not on the boundary of $f$ but adjacent to a black vertex $b$ on the boundary of $f$ and pointing towards (resp. away from) that black vertex $b$. Moreover, all sinks (resp. sources) are boundary vertices of $\pb$ (resp. $\pbm$). 
\end{lemma}
\proof{Direct consequence of Lemma \ref{lem:strandorder}.\qed}

By construction, every black vertex has exactly one incoming edge and every white vertex has exactly one outgoing edge. Therefore, the li-orientation is an example of a \emph{perfect orientation}, see \cite{postgrass}.

At each black vertex, a li-orientation distinguishes one of the edges from the other two. Therefore, it is possible to define corresponding li-matchings, which will also prove useful later on.

\begin{definition}
	Let $\tcd$ be a labeled TCD. Define the \emph{li-matchings $\lima$ (resp. $\limam$)} on $\pb$ (resp.~$\pbm$) as the subset of edges that point away from (resp. towards) their incident black vertex.
\end{definition}

\begin{lemma}\label{lem:limatch}
	Let $\tcd$ be a labeled minimal TCD and let $\lima$ (resp. $\limam$) be the li-matchings on $\pb$ (resp.~$\pbm$). Then every interior vertex of $\pb$ (resp.~$\pbm$) is matched exactly once. Moreover, every interior face is matched exactly once in the sense that for each face $f$ there is exactly one edge $e$ that is not a boundary edge of $f$ but incident with a black vertex on the boundary of $f$.
\end{lemma}
\proof{
	The faces, white and black vertices of $\pb$ and $\pbm$ correspond to faces of $\tcd$. Therefore the lemma follows immediately from Lemma \ref{lem:strandorder}.\qed
}

Thus li-matchings coming from minimal TCDs are perfect matchings. Note that these perfect matchings are closely related to so called \emph{extremal matchings} that have been studied for minimal bipartite graphs on the torus \cite{broomheaddimer}.

\begin{lemma}
	The li-orientations $\lio$ and $\liom$ of a labeled minimal TCD $\tcd$ are acyclic.
\end{lemma}
\proof{
	Assume $\gamma$ is a cycle of $\pb$. Consider the subgraph $G=(B\cup W,E,F)$ of $\pb$ that is bounded by the cycle $\gamma$. Also denote by $B'$ the set of black vertices on $\gamma$ such that one of the vertices adjacent to $\gamma$ is not in $G$. We know that the li-orientations induce a li-matching on $G$ such that every black vertex is matched, therefore we have $|W|=|B|$. Also, the number of edges satisfies $|E|=3|B|-|B'|$ by construction. Additionally, Lemma \ref{lem:limatch} ensures that there is exactly one matched edge for every face of $G$ and that there is exactly one face for every matched edge except for those edges adjacent to vertices of $B'$. As an equation this observation reads $|F| = |B| -|B'|$. We can now insert the three equations we found into the formula for the Euler characteristic to obtain
	\begin{align}
		0 = 1 + |E| - |F| - |B| - |W| = 1 + (3|B|-|B'|) - (|B| -|B'|) - |B| - |B| = 1,
	\end{align}
	which is a contradiction. Therefore there is no cycle of $\pb$. The argument works exactly the same for $\pbm$ and we have thus proven the Lemma.\qed
}

\begin{definition}
	Let $\lio, \liom$ be the li-orientations of a labeled minimal TCD $\tcd$. We define the \emph{li-posets} $\lip, \lipm$ as the posets $(V(\pb), \leq), (V(\pbm), \leq)$ induced by the acyclic orientations $\lio, \liom$. Explicitly, for any two vertices $v,v'$ we define $v\leq v'$ if $v=v'$ or there is a directed path in $\lio$ (resp. $\liom$) from $v'$ to $v$.
\end{definition}

\begin{figure}
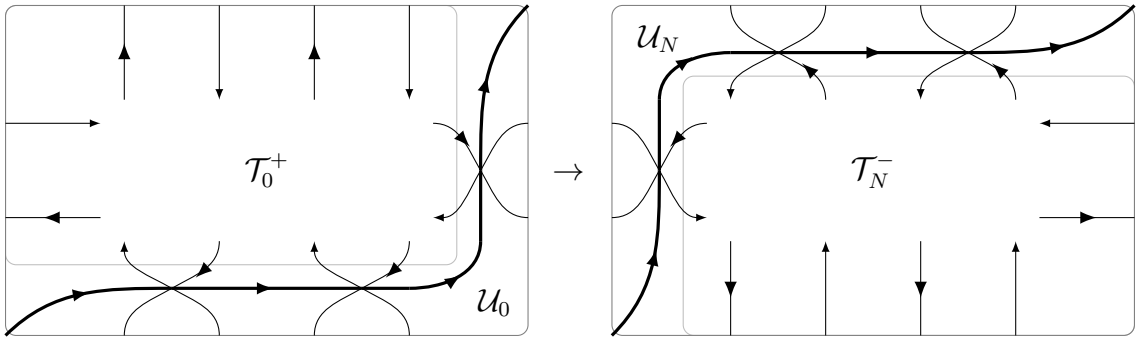


\caption{The situation before (left) and after (right) a sweep.}
\label{fig:tcdextension}
\end{figure}

\begin{definition}
	Given a TCD $\tcd$, we call $\tcd'$ a \emph{subdiagram} of $\tcd$ if there is a disk $D$ such that $\tcd' = \tcd \cap D$ is a TCD in the disk $D$.
\end{definition}

\begin{definition}
	Let $\tcd^+_0$ be a labeled minimal TCD with $n$ strands and let $k \in \N$ be such that $0 \leq k \leq n$. We define a new TCD $\mathcal U_0$ such that
	\begin{enumerate}
		\item $\mathcal U_0$ has $(n+1)$ strands,
		\item $\tcd^+_0$ is a subdiagram of $\mathcal U_0$,
		\item the additional strand of $\mathcal U_0$ not in $\tcd^+_0$ is called the \emph{sweep strand} and runs counterclockwise along the boundary of $\tcd^+_0$ in $\mathcal U_0$,
		\item the strand belonging to out-endpoint $i$ for $0 \leq i < k$ of $\tcd^+_0$ crosses the strand belonging to in-endpoint $i$ of $\tcd^+_0$ on the sweep strand, and these are the only crossings in $\mathcal U_0$ that are not in $\tcd^+_0$.
	\end{enumerate}
	See also Figure \ref{fig:tcdextension}. We say $k$ is the \emph{length of the sweep strand}.
\end{definition}

\begin{definition}
	Let $\tcd^+_0$ be a labeled minimal TCD and $\mathcal U_0$ the diagram where we added a sweep strand of length $k$. A \emph{sweep} is a sequence of TCDs $\mathcal U_0, \mathcal U_1, \mathcal U_2, \dots, \mathcal U_N$ such that
	\begin{enumerate}
		\item consecutive TCDs are related by a 2-2 move that involves the sweep strand,
		\item the number of faces to the right of the sweep strand is increasing,
		\item only boundary faces are to the left of the sweep strand in $\mathcal U_N$.
	\end{enumerate}
	We say a TCD $\tcd^+_0$ is \emph{$k$-sweepable} if a sweep with a sweep strand of length $k$ exists. We say a $\tcd^+_0$ is \emph{sweepable} if a number $k$ exists such that $\tcd^+_0$ is $k$-sweepable. We also denote by $\tcd^+_i$ (resp. $\tcd^-_i$) the subdiagram of $\mathcal U_i$ that is above (resp. below) the sweep strand, see Figure \ref{fig:tcdextension}. See Figure \ref{fig:sweepexample} for an example of a sweep.
\end{definition}

\begin{figure}


	\caption{An example of a sweep of a TCD. The oriented edges of $\lio$ are in black and of $\liom$ in gray.}
	\label{fig:sweepexample}
\end{figure}

In fact, every 2-2 moves in a sweep increases the number of faces below (that is, to the right of) the sweep strand by one. The sweep strand has to traverse every interior face of $\mathcal U_0$. The number of interior faces of $\mathcal U_0$ is just the sum of the number of interior faces of $\tcd^+_0$ and $(2k-1)$. Therefore the length $N$ of the sweep sequence is indeed 
\begin{align}
	N = |F^{\mathrm{int}}(\tcd^+_0)| + 2k-1.
\end{align}
Another trivial observation is that by definition if $\tcd^+_0$ is sweepable, then so is $\tcd^+_i$ for any other $i$. Moreover, the endpoint matching of $\mathcal U_i$ does not depend on $i$, as the diagrams for different $i$ are related by 2-2 moves. Therefore, the endpoint matching of $\tcd^-_N$ is determined by the endpoint matching of $\tcd^-_0$.

\begin{lemma} \label{lem:sweepcondi}
	Let $\tcd$ be a labeled minimal TCD with $n$ strands and endpoint-matching $\pi$, such that $\tcd$ is $k$-sweepable. Then
	\begin{enumerate}
		\item for all strands $i$ with $\pi(i) \leq k$ we have that $\pi(i) < i$ and
		\item for all strands $i$ with $\pi(i) > k$ we have that $\pi(i) > i$.\qedhere
	\end{enumerate}
\end{lemma}
\proof{If the first condition is not satisfied, then $\mathcal U_0$ is not a minimal diagram. Indeed, in this case there is an $i$ with $i < \pi(i) \leq k$. The strand $i$ can never be fully below the sweep strand. If the second condition is not satisfied, then $\mathcal U_N$ is not a minimal diagram. Indeed, in this case there is an $i$ with $i > \pi(i) \geq k$. The strand $i$ can never be passed by the sweep strand.\qed
}

In fact, we will show that the conditions of Lemma \ref{lem:sweepcondi} are not only necessary but also sufficient for the sweepability of a TCD.

\begin{figure}
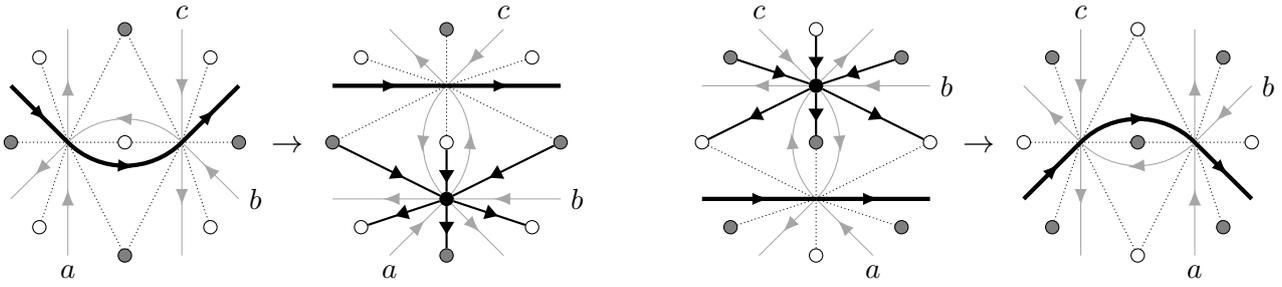

	\centering
	
	\caption{The change of the TCD and the li-orientation for the two types of 2-2 moves in a sweep. The bold strand is the sweep strand.}
	\label{fig:sweeporient}
\end{figure}

\begin{lemma}\label{lem:liomin}
	Let $\mathcal U_i$ be a TCD in a sweep. Then every clockwise 2-cell above and adjacent to the sweep strand (see Figure \ref{fig:sweeporient}) is a minimal element of $\liom(\tcd^+_i)$.
\end{lemma}
\proof{We choose the labels $a,b,c$ for the three strands involved that are not the sweep strand, see Figure \ref{fig:sweeporient}. Let us show that $a < b < c$. First of all, $c$ ends at the position of $a$. Thus Lemma \ref{lem:sweepcondi} ensures that $c > a$. Because strand $b$ crosses the strands $a$ and $c$ already at the 2-cell that we consider, it cannot cross the parts of strands $a$ and $c$ before the 2-cell. Therefore $a < b < c$ and together with the definition of the li-orientation, the 2-cell has to be a minimal element.\qed
}

Analogously, one can prove that every counterclockwise cell below is a maximal element of $\lio(\tcd^-)$, but we do not need it in the following.

\begin{lemma}\label{lem:liopreserved}
	The li-orientation $\liom(\tcd^+_i)$ is preserved along the sweep. More precisely, if $v,v'$ are vertices of $\pbm(\tcd^+_i)$ as well as of $\pbm(\tcd^+_j)$, then $v \leq v'$ with respect to $\liom(\tcd^+_i)$ if and only if $v \leq v'$ with respect to $\liom(\tcd^+_j)$.
\end{lemma}
\proof{We observe in Figure \ref{fig:sweeporient} that a 2-2 move at a counterclockwise 2-face does not affect the li-orientation $\liom$ at all. At a clockwise cell we do remove one vertex from $\pbm$ but no other relations are affected. Moreover, note that a 2-2 move at a counterclockwise 2-face removes one strand, but this strand has no crossings and is therefore irrelevant for the li-orientation. At a clockwise 2-face we split the strand $b$ into two strands $b,b'$. Therefore we have to show that the li-orientation at later crossings along strand $b$ are not affected by the split. Clearly $b' < b$ and therefore if $b$ was already the minimum at a later crossing, then this is preserved. If at a later crossing we had $x < b < y$, then one finds $x < a$ and therefore $x < b'$, because else there is again a contradiction to the minimality of the TCD. If $b$ was the maximum at a later crossing, i.e. $x < y < b$, then one finds that $x < y < a$ and therefore the li-orientation is preserved.
\qed
}

Note that Lemma \ref{lem:liomin} and Lemma \ref{lem:liopreserved} together imply that the order in which counterclockwise faces of $\tcd_0^+$ are (completely) swept is determined by $\liom$. In fact, any sweep $\mathcal U_0, \dots, \mathcal U_N$ also yields a sequence of counterclockwise faces $f_0,f_1,\dots, f_N'$. This sequence is a linear extension of $\liom(\tcd_0^+)$ in the sense of order theory.

\begin{lemma}\label{lem:sweepeyesaremin}
	Let $\tcd$ be a labeled minimal TCD that satisfies the condition of Lemma \ref{lem:sweepcondi} and such that no counterclockwise 2-face is adjacent to the sweep strand. Then the only minimal elements of $\liom(\tcd)$ are clockwise 2-faces adjacent to the sweep strand.
\end{lemma}
\proof{
	As a consequence of Lemma \ref{lem:strandorder}, every face including the boundary faces of $\tcd$ can have at most one incoming edge in $\liom$. The only faces that have only one edge are boundary faces. A clockwise boundary face $f$ has vertex degree 1 in $\pbm$ if it has face degree 3 in $\tcd$ or if it separates the TCD into two TCDs $\tcd',\tcd''$ that are only connected to each other via the boundary of the TCD. Let $\tcd'$ be the TCD that is to the right of the boundary strands of $f$. But $\tcd'$ cannot posses a right-turning strand, as else it would violate the conditions of Lemma \ref{lem:sweepcondi}. Therefore $\tcd'$ has to consist only of isolated left-turning strands, which is also not possible because then there would be a counterclockwise 2-face adjacent to the sweep strand. Analogously, one can also exclude the possibility of separating faces not adjacent to the sweep strand. Therefore the only possible sinks are indeed degree 3 faces on the boundary. We have already shown in Lemma \ref{lem:liomin} that the degree 3 faces adjacent to the sweep strand are minima. An analogous argument shows that in fact degree 3 faces not adjacent to the sweep strand are maximal elements.\qed
}

\begin{theorem}\label{th:sweepcondi}
	The conditions of Lemma \ref{lem:sweepcondi} imply $k$-sweepability.
\end{theorem}
\proof{As $\lipm(\tcd^+_i)$ is a finite poset, there are always minimal elements in $\lipm(\tcd^+_i)$. If no counterclockwise 2-face is available to continue the sweep, then Lemma \ref{lem:sweepeyesaremin} ensures that the minimal elements correspond to clockwise 2-faces where one can continue the sweep. \qed
}

\begin{lemma}\label{lem:sweepdimension}
	If a TCD $\tcd^+_0$ is $k$-sweepable then $k$ equals the maximal dimension of $\tcd^+_0$, that is $k=|W|-|B|-1$, where $W$ and $B$ are the number of white and black vertices of $\pb(\tcd^+_0)$ respectively.
\end{lemma}
\proof{
	Consider $\mathcal U_0$ with a labeling starting at the sweep strand and the li-orientation $\lio(\mathcal U_0)$. From Lemma \ref{lem:strandorder} we see that all boundary faces below the sweep strand are minima, and so is the boundary face just to the left of the beginning of the sweep strand. All other boundary faces have one outgoing edge. Therefore the corresponding li-matching $\lima(\mathcal U_0)$ matches all but $(k+1)$ of the counterclockwise faces of $\mathcal U_0$. The difference of the numbers of white and black vertices of $\pb(\mathcal U_0)$ is the same as of $\pb(\tcd^+_0)$ and therefore we obtain that $|W|-|B| = k + 1$, which proves the claim.\qed
}

As not every TCD is sweepable, it is not clear how useful sweeps are. In fact, we will show that sweeps naturally occur in many examples as a way of propagating information through TCD maps. In particular, they relate sections and projections (see Section \ref{sec:tcdextensions}) and allow a natural way of parametrizing TCD maps with given TCD via certain projective invariants (see Section \ref{sec:fromprojinvariants}). But one can also say more about sweepable TCDs from a purely combinatorial point of view.

\begin{lemma}\label{lem:balancedsweepable}
	If a minimal TCD $\tcd$ is balanced (see Definition \ref{def:endpointmatching}), that is $\tcd$ has endpoint matching $\enm{n}{k}$ for some $k,n\in\N$, it is $k$-sweepable for any choice of labeling.
\end{lemma}
\proof{As every strand in $\tcd$ has length $k$, it is clear that the conditions of Lemma \ref{lem:sweepcondi} are satisfied. Due to Theorem \ref{th:sweepcondi}, the TCD $\tcd$ is therefore $k$-sweepable.\qed}

\begin{figure}
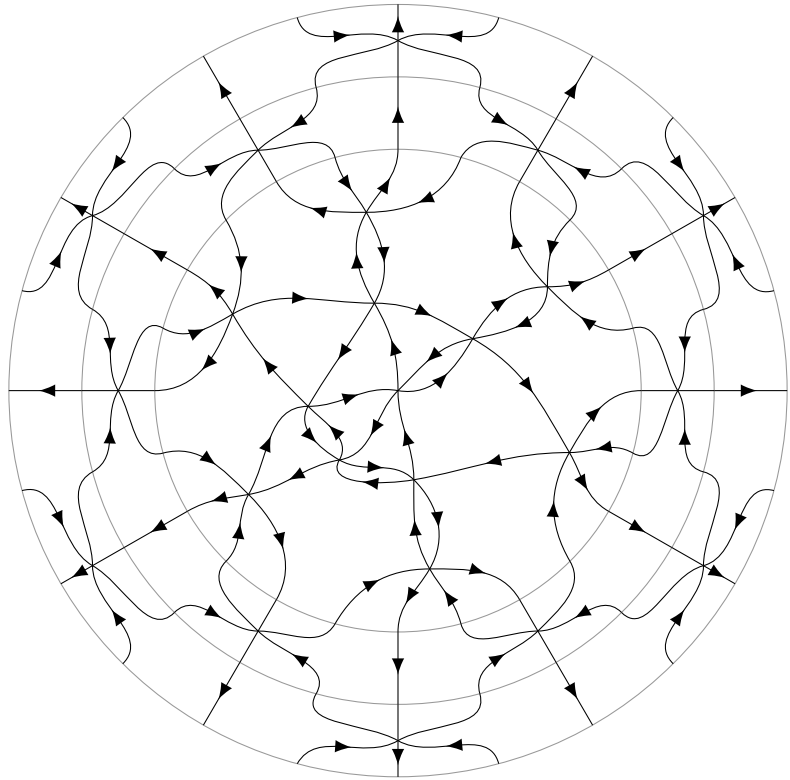


	\caption{Three TCDs $\tcd \subset \tcd' \subset \tcd''$. The innermost TCD $\tcd$ is not sweepable, $\tcd'$ is sweepable and $\tcd''$ is balanced.}
	\label{fig:bridgedecomp}
\end{figure}

\begin{lemma}\label{lem:subsweepable}
	Let $\tcd$ be a minimal TCD of maximal dimension $k$ with $n$ strands. Then there is a TCD $\tcd'$ such that
	\begin{enumerate}
		\item $\tcd$ is a subdiagram of $\tcd'$,
		\item the maximal dimension of $\tcd'$ is also $k$,
		\item the number of strands of $\tcd'$ is also $n$ and
		\item $\tcd'$ is sweepable.\qedhere
	\end{enumerate}
\end{lemma}
See also Figure \ref{fig:bridgedecomp}.
\proof{
	We show that there is a balanced TCD $\tcd''$ with $n$ strands and maximal dimension $k$ such that $\tcd$ is a subdiagram of $\tcd''$. Due to Lemma \ref{lem:balancedsweepable}, the TCD $\tcd''$ is sweepable. This suffices to prove the claim. Let us denote the length of a strand $s$ by $\ell(s) = \pi(s) - s$, where $\pi$ is the endpoint matching with the convention $0\leq \ell(s) < n$. We construct a sequence of diagrams $\tcd_0 \subset \tcd_1 \subset \dots \subset \tcd_m$ such that $\tcd_0 = \tcd$ and $\tcd_m = \tcd''$. If $\tcd_i$ is not balanced, then there is a strand $j$ of maximal length such that $\ell(j) > \ell(j+1)$. Let $\tcd_{i+1}$ be a copy of $\tcd_i$ but with an additional crossing at the boundary of $\tcd$ where the strands $j$, $j+1$ and the strand ending between the two cross. The new lengths $\ell'$ in $\tcd_{i+1}$ are therefore \begin{align}
		\ell'(j) = \ell(j+1)+1 \mbox{ and } \ell'(j+1) = \ell(j)-1.
	\end{align}
	Note that $\ell(j) \geq \ell(j+1)+2$ as else the two strands would exit the TCD at the same point. Therefore in $\tcd_{i+1}$, either the maximal length is decreased or the number of strands of maximal length is decreased. Therefore the sequence of TCDs is finite. The number of strands is also preserved. Moreover, in the bipartite graphs $\pb(\tcd_i)$ we add a black and a white vertex. Therefore the maximal dimension is also preserved. It remains to verify that the TCDs $\tcd_i$ also remain minimal. Self-intersections are not possible, because if the exiting strand is strand $j$ or $(j+1)$, then $j$ is not a strand of maximal length. Similarly, if strand $j$ and $(j+1)$ already have a crossing, then the length of strand $j$ is smaller or equal than the length of strand $(j+1)$ in contradiction to the assumptions. Therefore $\tcd_{i+1}$ is minimal as well. As the sequence is finite and only stops when all strands have equal length, the claim is proven.\qed
}

\begin{remark}
	Consider the sequence $(\tcd_i) :=\tcd_0 \subset \dots \subset \tcd_m$ of Lemma \ref{lem:subsweepable}, and assume we have an accompanying sequence of TCD maps $T_i: \tcd_i \rightarrow \CP^k$. Then the steps in the sequence $(\tcd_i)$ correspond to adding `marked' points on the boundary. Therefore in a sense, a TCD map on a non sweepable TCD is simply a TCD map that lacks a bit of information on its boundary. Secondly, there is the so called \emph{circular Bruhat order} \cite{postgrass} on permutations. The endpoint matchings of the TCDs in $(\tcd_i)$ are in fact an (upwards) chain in the circular Bruhat order. Indeed, there is a generalization of permutahedrons based on the circular Bruhat order and chains are shown to correspond to \emph{bridge-decompositions} of TCDs \cite{williamsgrass}. The sequence $(\tcd_i)$ is indeed part of such a decomposition. Note however that the sweep sequence $(\tcd^+_0, \tcd^+_1, \dots, \tcd^+_N)$ generally is not a bridge-decomposition. It is possible that the two coincide, viewing 2-2 moves of the sweep at clockwise faces as steps in the bridge-decomposition and viewing 2-2 moves at counterclockwise faces as henceforth ignoring of fixed-points. However, the bridge-decomposition does not allow all 2-2 moves at clockwise faces and sweeps do only allow steps of the bridge-decomposition at faces adjacent to the sweep strand. It is certainly an interesting question if one can develop a better understanding of sweep-decompositions analogously to bridge-decompositions. For that purpose, it may be useful to view non-sweepable diagrams as diagrams where multiple sweep strands are necessary. In that case it might be that sweep-decompositions contain bridge-decompositions as a special case. Later, we will have results specifically for TCD maps defined on sweepable TCDs, for example in Section \ref{sec:fromprojinvariants}, where we construct TCD maps from invariants. For TCD maps with non-sweepable TCD $\tcd_0$, it might be possible to extend these results, by understanding a sweep on a larger, sweepable TCD $\tcd_m \supset \tcd_0$. If we observe the sweep on $\tcd_m$ restricted to $\tcd_0$, we see that the sweep strand may not be connected in $\tcd_0$, and may also enter and leave $\tcd_0$ at several different locations along the boundary. However, the sweep still monotonously passes through $\tcd_0$ by 2-2 moves and there are just some additional operations needed at the boundary. It would be very interesting to gain a better understanding of these restricted sweeps as well as the involved orientations.
\end{remark}

Let us prepare for the consequences for geometry by looking at the combinatorics of something that we have largely ignored so far, namely $\tcd^-_N$.

\begin{figure}
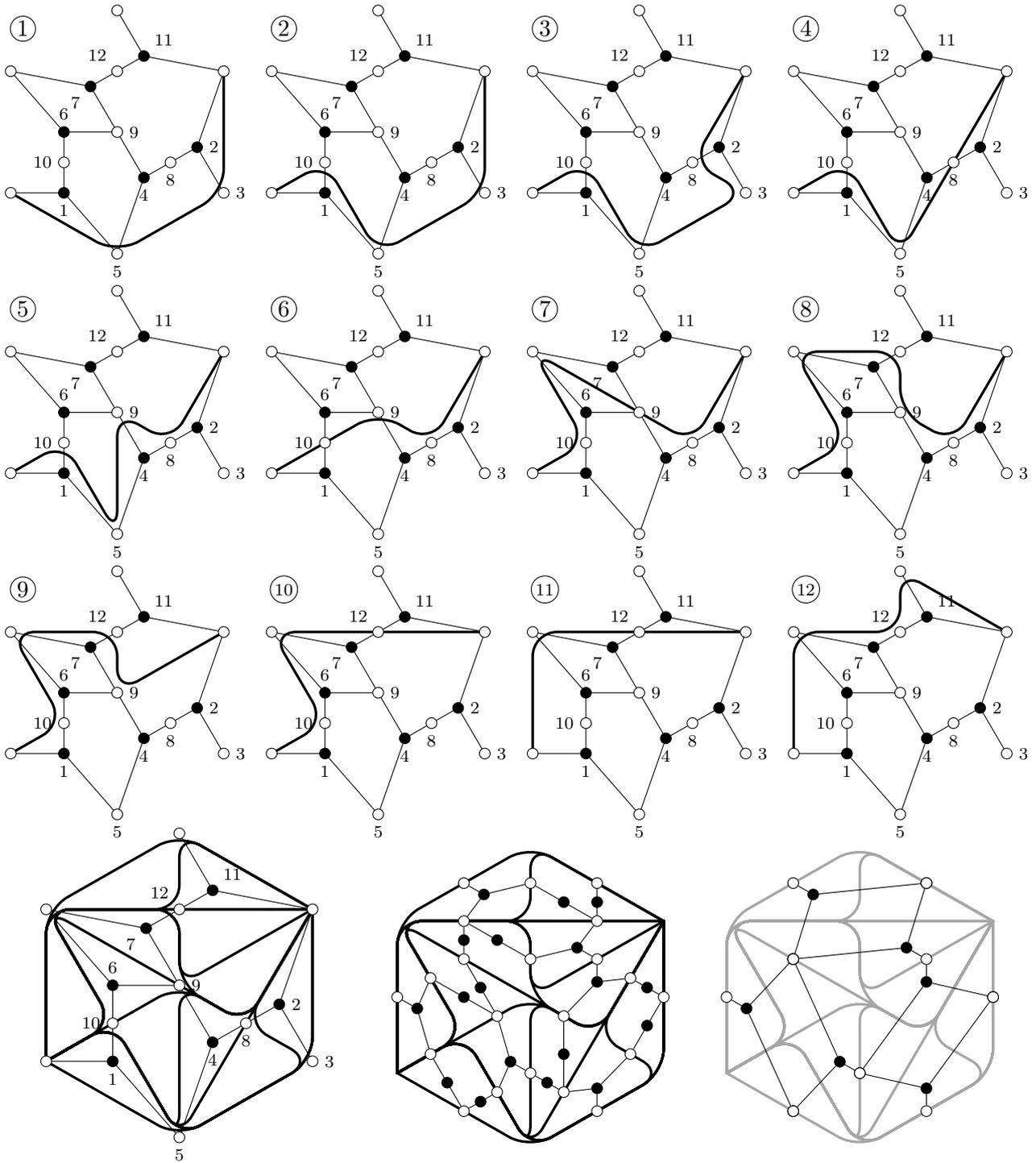

	%
	
	\caption{Top: All but the last step of a sweep. Bottom from left to right: All steps of a sweep; with replacements glued into the triangles; after contracting.}
	\label{fig:sweepsection}
\end{figure}

\begin{theorem}\label{th:sweepsection}
	Consider a sweep from $\mathcal U_0$ to $\mathcal U_N$. Then there is a choice of section such that
	\begin{equation}
		\tcd^-_N = \sigma(\tcd^+_0).\qedhere
	\end{equation}
\end{theorem}
\proof{
	Each step in the sweep erases either a black or a white vertex from $\pb(\tcd_0^+)$, thus we can label the vertices of $\pb(\tcd_0^+)$ by the sweep step in which they disappear, see Figure \ref{fig:sweepsection}. We can draw the sweep strand in $\pb(\tcd_0^+)$ as a line connecting the white vertices above and adjacent to the sweep strand. Each 2-2 move in the sweep at a clockwise face means we pull the sweep strand over the corresponding black vertex, while a 2-2 move at a counterclockwise face means we close off an open triangle adjacent to the corresponding white vertex. On the other side of the sweep strand, we insert new vertices and edges to obtain $\pb(\tcd^-_N)$. Specifically, in each triangle corresponding to a black or white vertex, we glue in
	\begin{align}
,
	\end{align}
	where the rotation of the first triangle does not matter. After gluing all such triangles together, we contract black vertices of degree 2 and obtain a new TCD $\tcd'$. Each white vertex in $\tcd'$ corresponds either to a collection of black vertices that are adjacent to each other via 4-faces of $\pb(\tcd^+_0)$, a boundary edge, or a diagonal of a face. Each black vertex corresponds to a triangle of a triangulation of a face of $\pb(\tcd^+_0)$. Therefore $\tcd'$ is indeed a section of $\tcd^+_0$.\qed
}

\section{Sweeps reduced to quad-graphs, pseudoline arrangements}\label{sec:sweepsqg}

Recall that we introduced quad-graphs in Section \ref{sec:quadgraphs}, where we show that quad-graphs can be obtained as a special case of TCDs. In particular, we observed that in a TCD that specializes to a quad-graph, each strip of the quad-graph corresponds to two strands of the TCD (see Definition \ref{def:tcdquad}). The arrangement of strips of a minimal quad-graph (see Definition \ref{def:minqg}) is also called a \emph{pseudoline arrangement} \cite{felsnerbook}. To some extent, the results of the previous section generalize results and methods used in the study of sweeps of pseudoline arrangements, see \cite{felsnerbook}. We do not need the general sweep results for pseudoline arrangements (resp. quad-graphs), therefore we only give an informal explanation of how the sweeps on TCDs reduce to sweeps of pseudoline arrangements. However, the reduction of the li-orientation $\lio$ to quad-graphs is an object that will be useful in Chapter \ref{cha:bilinear}, where we investigate certain reductions of integrable maps on quad-graphs. Therefore we will give a precise explanation of that reduction at the end of the section.

For TCDs we introduced labeled TCDs (see Definition \ref{def:labeledtcd}), and we can do something similar to quad-graphs to obtain \emph{labeled quad-graphs}. 

\begin{definition}\label{def:labeledqg}
	Let $\qg$ be a quad-graph with $n$ strips and no closed strips (loops). Let us label the endpoints of the strips in counterclockwise order by the numbers $\{0,1,\dots,2n-1\}$. We label each strip by the smaller label of the two labels assigned to the two endpoints of the strip. We call a $\qg$ with such labels a \emph{labeled quad-graph}.
\end{definition}

The \emph{endpoint matching} of a labeled quad-graph is the set of $n$ pairs of endpoint labels of the strips. We call a minimal quad-graph a \emph{balanced quad-graph} if each strip intersects each other strip exactly once. As a consequence, the endpoint matching  of a balanced labeled quad-graph consists of the $n$ pairs $\{i,i+n\}$ for $i$ such that $0\leq i < n$, and the strips are labeled by the numbers $\{0,1,\dots,n-1\}$.

Recall that we explained a procedure to obtain a TCD $\tcd_\qg$ from every quad-graph $\qg$ in Definition \ref{def:tcdquad}. We have already observed in Lemma \ref{lem:qgmintcdmin} that a quad-graph $\qg$ is minimal if and only if $\tcd_\qg$ is minimal.

\begin{lemma}
	If $\qg$ is balanced then $\tcd_\qg$ is balanced and vice versa.
\end{lemma}
\proof{
	Assume $\qg$ is balanced, and choose a labeling of $\qg$. Then there is a labeling of $\tcd_\qg$ such that $\pi(a) = a+n$ and $\pi(a+n) = a$ for all $a$ with $0\leq a<2n$, where $\pi$ is the endpoint matching of $\tcd_\qg$. Thus $\tcd_\qg$ is a balanced TCD with endpoint matching $\enm {2n}{n}$ (see Definition~\ref{def:endpointmatching}). Conversely, assume $\tcd_\qg$ is balanced and choose a labeling of $\tcd_\qg$. Because $\tcd_\qg$ is the TCD associated to a quad-graph, $\tcd_\qg$ has $2n$ strands and $\pi(a) = b$ implies $\pi(b) = a$ for all $a,b$. The only way that $\tcd_\qg$ is also balanced is if $b-a = a - b$ in $\Z_{2n}$ and therefore $b = a +n$, which in turn implies that $\qg$ is balanced.\qed
}

To perform a sweep of a marked pseudoline arrangement, or equivalently a labeled balanced quad-graph $\qg$, we add a strip that runs counterclockwise around the boundary and intersects every strip exactly once, thus obtaining a quad-graph $\qg_0$ that is also balanced. As in the case of sweeps of TCDs, the sweep itself is a sequence of quad-graphs $\qg_0,\qg_1,\dots, \qg_N$ such that in $\qg_N$ the sweep strip runs clockwise along the boundary, and consequent quad-graphs are related by cube-flips (instead of 2-2 moves), such that the number of vertices below the sweep strip is strictly increasing in every step. Such sweeps always exist, see \cite{felsnerbook}. There is one interesting difference between sweeps of TCDs and sweeps of quad-graphs. In the case of TCDs we observed that the TCD $\tcd_N^-$ after the sweep is quite different from the TCD $\tcd_0^+$ before the sweep. In fact Theorem \ref{th:sweepsection} shows that the two TCDs are related by taking a section. On the other hand, in the case of quad-graphs we have that $\qg_N^- = \qg_0^+$, that is that the combinatorics of $\qg$ are the same before and after the sweep.

\begin{remark}
	For the reader familiar with discrete integrable systems or discrete differential surface theory, the sweep captures the combinatorics of what is sometimes called a \emph{Darboux transform} or \emph{Bäcklund transform}, for example in a Q-net (see Definition \ref{def:qgqnet}).
\end{remark}

Due to Lemma \ref{lem:cubefliptcd}, we know that every cube flip in $\qg$ corresponds to a sequence of 2-2 moves in $\tcd_\qg$. Thus it is plausible that the sweep in $\qg$ actually corresponds to the composition of a sweep and a backwards sweep of $\tcd_\qg$. We do not investigate this further as it is not relevant for the remainder. What is relevant for later results is to understand the analogue of the li-orientation $\lio$ (see Definition \ref{def:liorient}) of the planar bipartite graph $\pb_\qg$ associated to $\tcd_\qg$ (see Definition \ref{def:pb}) in the case of quad-graphs.

\begin{figure}
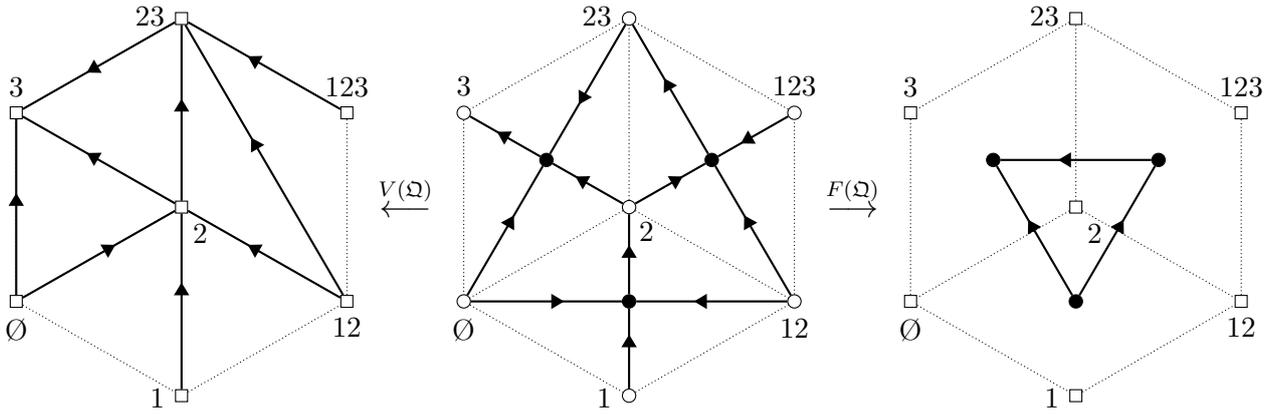


	\caption{The li-orientation $\lio(\qg^D)$ of a quad-graph in the center and the induced li-orientations on the vertices (faces) of $\qg^D$ on the left (right).}
	\label{fig:qglio}
\end{figure}

Let us denote by $\qg^*$ the graph-dual of the quad-graph $\qg$, without the outer face. Therefore $\qg^*$ has a vertex for every face of $\qg$. Let us also denote by $\qg^D$ the \emph{double graph}, which has a vertex for every vertex and every quad of $\qg$ and there is an edge between two vertices $v,v'$ of $\qg^D$ if $v$ corresponds to a vertex of $\qg$ that is incident to a quad of $\qg$. 

\begin{definition}
	Let $\qg$ be a minimal quad-graph. The \emph{li-orientation} $\lio(\qg^D)$ of $\qg^D$ is obtained from $\lio(\pb_\qg)$ by contracting all white vertices of $\pb_\qg$ that do not correspond to vertices of $\qg$. More explicitly, we obtain the li-orientation of $\qg^D$ by gluing
	
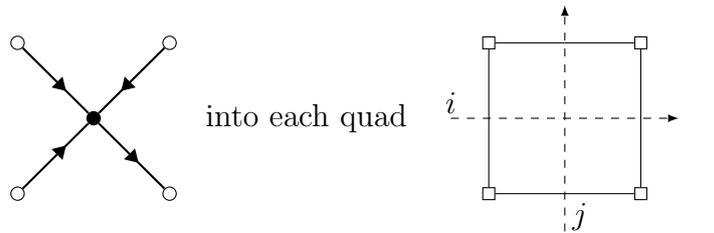
\begin{center}\vspace{-3mm}
	\begin{tikzpicture}[scale=1]
		\begin{scope}[shift={(-0.2,0)}]
			\node[wvert] (v) at (0,0) {};
			\node[wvert] (v1) at (2,0) {};
			\node[wvert] (v2) at (0,2) {};
			\node[wvert] (v12) at (2,2) {};
			\node[bvert] (f21) at (1,1) {};
			\draw[-,ultra thick]				
				(f21) edge[rorient] (v) edge[orient] (v1) edge[rorient] (v2) edge[rorient] (v12)
			;
			\node[] at (3.8,1) { into each quad };
		\end{scope}
		\begin{scope}[shift={(6,0)}]
			\node[wqgvert] (v) at (0,0) {};
			\node[wqgvert] (v1) at (2,0) {};
			\node[wqgvert] (v2) at (0,2) {};
			\node[wqgvert] (v12) at (2,2) {};
			\coordinate (e1) at (1,0);
			\coordinate (e2) at (0,1);
			\coordinate (e12) at (1,2);
			\coordinate (e21) at (2,1);				
			\draw[black]
				(v) edge (e1) edge (e2)
				(v1) edge (e21) edge (e1)
				(v12) edge (e12) edge (e21)
				(v2) edge (e2) edge (e12)
			;			
			\node[] at (3,1) {,};
			
			\draw[dashed,-latex] (-0.5,1) -- (2.5,1);
			\draw[dashed,-latex] (1,-0.5) -- (1,2.5);
			\node (i) at (-0.5,1.2) {$i$};
			\node (i) at (1.2,-0.3) {$j$};
		\end{scope}
	\end{tikzpicture}	
	\end{center}\vspace{-5mm}
	where $i,j$ are the labels of the two strands, $i < j$, and the coloring of the vertices of $\qg$ does not matter.
\end{definition}

The significance of the li-orientation $\lio(\qg^D)$ of the double of a quad-graph becomes apparent when we deal with per-quad propagation of data (see Theorem \ref{th:liftofdoublebkp} and its proof for an example).  Consider some map $g:V(\qg) \rightarrow S$ into some arbitrary set $S$. Assume there is a rule that determines the value of $g$ on the fourth vertex of any quad from the values of $g$ on the other three vertices of that quad. What is the free initial data that completely determines the system and how do we propagate from it? We simply choose the values of $g$ freely on the minimal elements of $\lio(\qg^D)$. Then we propagate the data along the li-orientation $\lio(\qg^D)$. Because 
\begin{enumerate}
	\item any black vertex of $\qg^D$ (a face of $\qg$) has exactly one outgoing edge in $\lio(\qg^D)$,
	\item every white vertex of $\qg^D$ (a vertex of $\qg$) has at most one incoming edge,
	\item and $\lio(\qg^D)$ is acyclic,
\end{enumerate}
there are no conflicts when propagating the data. Moreover, since we chose values on all minimal elements, we can propagate to the whole quad-graph $\qg$.

Note that we can also look at the locally induced li-orientations on $\qg$ as well as on $\qg^*$, see Figure \ref{fig:qglio} for an example. For the li-orientation on $\qg$, we glue the three arrows into each quad that are induced by $\lio(\qg^D)$. For the li-orientation on $\qg^*$, we glue an arrow between every pair of adjacent faces of $\qg$, the direction of which is induced by $\lio(\qg^D)$. This is well-defined because faces of $\qg^D$ cannot be oriented cyclically in $\lio(\qg^D)$, and because of properties (1) and (2) of $\lio(\qg^D)$ mentioned above.

From the viewpoint of sweeps of pseudoline arrangements, the set of sweeps is in one-to-one correspondence with the linear extensions of $\lio(\qg^*)$, by writing down the sequence of pseudoline crossings swept in the sweep. Similarly, the set of sweeps is also in one-to-one correspondence with the linear extensions (in the sense of partial order theory) of $\lio(\qg)$ (without the minimal elements), by writing down the sequence of first contacts of the sweep line with regions of the pseudoline arrangement (vertices of $\qg$) during the sweep.

Another small note is that we can also look at what the li-orientation $\liom$ induces on the quad-graph. The short answer is that we can also contract $\liom$ in each quad, that is we glue
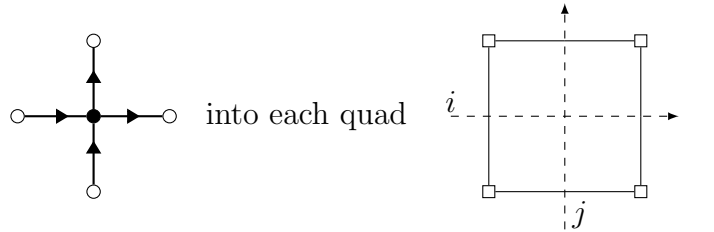
\begin{center}
\begin{tikzpicture}[scale=1]
	\begin{scope}[shift={(-0.2,0)}]
		\coordinate (v) at (0,0) {};
		\coordinate (v1) at (2,0) {};
		\coordinate (v2) at (0,2) {};
		\coordinate (v12) at (2,2) {};
		\node[bvert] (f21) at (1,1) {};
		\coordinate[wvert] (e1) at (1,0);
		\coordinate[wvert] (e2) at (0,1);
		\coordinate[wvert] (e12) at (1,2);
		\coordinate[wvert] (e21) at (2,1);				
		\draw[-,ultra thick]				
			(f21) edge[rorient] (e1) edge[rorient] (e2) edge[orient] (e12) edge[orient] (e21)
		;
		\node[] at (3.8,1) { into each quad };
	\end{scope}
	\begin{scope}[shift={(6,0)}]
		\node[wqgvert] (v) at (0,0) {};
		\node[wqgvert] (v1) at (2,0) {};
		\node[wqgvert] (v2) at (0,2) {};
		\node[wqgvert] (v12) at (2,2) {};
		\coordinate (e1) at (1,0);
		\coordinate (e2) at (0,1);
		\coordinate (e12) at (1,2);
		\coordinate (e21) at (2,1);				
		\draw[black]
			(v) edge (e1) edge (e2)
			(v1) edge (e21) edge (e1)
			(v12) edge (e12) edge (e21)
			(v2) edge (e2) edge (e12)
		;			
		\node[] at (3,1) {,};
		
		\draw[dashed,-latex] (-0.5,1) -- (2.5,1);
		\draw[dashed,-latex] (1,-0.5) -- (1,2.5);
		\node (i) at (-0.5,1.2) {$i$};
		\node (i) at (1.2,-0.3) {$j$};
	\end{scope}
\end{tikzpicture}	
\end{center}
for $i<j$ to obtain a li-orientation $\liom(\qg^+)$ of $\qg^D$. Note that $\liom(\qg^+)$ has two outgoing and two incoming edges at every black vertex, but still one incoming edge at every white vertex. Therefore this orientation can be used to propagate systems where the values of a function at any two adjacent edges of a quad determine the values at the other two edges of that quad. These two-edge to two-edge systems are called \emph{Yang-Baxter maps} in the literature. They go back to a question asked by Drinfeld \cite{drinfeldyangbaxter} on set-theoretical solutions of the quantum Yang-Baxter equation. For a survey see  \cite{veselovyangbaxter} and for a classification of quadrirational Yang-Baxter maps see \cite{absyangbaxter}, which also gives relations to projective geometry. We also briefly discuss the latter maps in Section \ref{sec:quadrirationalybmaps}.

\section{Extensions and lifts via sweeps}\label{sec:tcdextensions}

We have seen in Section \ref{sec:sweeps} that for sweepable TCDs it is possible to relate the TCD $\tcd$ itself to a section $\sigma(\tcd)$ of it. Here, we investigate the corresponding relations for TCD maps.

Recall that we introduced the notions of flip-generic in Definition \ref{def:gentcdmap}, 1-generic in Definition \ref{def:onegeneric} and generic hyperplanes in Definition \ref{def:generichyperplane}.

\begin{theorem}\label{th:tcdmapsweep}
	Let $\mathcal U_0, \mathcal U_1,\dots, \mathcal U_N$ be a sweep sequence as introduced in Section \ref{sec:sweeps}. Consider a sequence of TCD maps $U_0, U_1,\dots, U_N$ such that $U_i$ maps from the counterclockwise faces of $\mathcal U_i$ to $\CP^n$ for $0 \leq i \leq N$, and such that consecutive TCD maps are obtained by the 2-2 moves of the sweep sequence. Write $T^\pm_i$ for the restrictions of $U_i$ to $\tcd^\pm_i$ and let $H$ be the span of $T^-_0$. Assume that $U_0$ is 1-generic and flip-generic, spans $\CP^n$ and that $H$ is $(n-1)$-dimensional and generic with respect to $T^+_0$. Then there is a choice of section such that $T^-_N = \sigma_H(T^+_0)$.	
\end{theorem}

\proof{ We first claim that every point of every TCD $\tcd_i^-$ is indeed in $H$. By definition of $H$ this is true for $\tcd_0^-$. Moreover, each 2-2 move that adds a white vertex $w$ to $\pb_{i+1}^-$ of the next TCD $\tcd_{i+1}^-$ in the sweep sequence is a resplit. Each such resplit places the new point $T_{i+1}^-(w)$ on a line spanned by points of $T_i^-$. Therefore the first claim holds. Moreover, we have observed in the proof of Theorem \ref{th:sweepsection} that each new white vertex added to $\pb_{i+1}^-$ is associated to a collection of black points that belong to a line of $T_0^+$. Indeed, by construction the point in $H$ associated to the new white vertex is also on the corresponding line of $T_0^+$, and is thus the intersection of that line with $H$.\qed 
}

Therefore, if we have the data of a TCD map $T$, we can construct the section of $T$ by adding some 1D-data (in the sense of discrete integrable systems) on the boundary and propagate the information through the whole TCD map $T$ to obtain a section. Moreover, one can also choose to apply 2-2 moves to $T$ to obtain another TCD map $\tilde T$. For example, this could be two Q-nets that are related to each other by cube flips. Then one can obtain the section of $\tilde T$, which in this example is a Darboux map, via a sweep. On the other hand, because the flip graph is connected we know there is also a sequence of 2-2 moves from the section of $T$ to the section of $\tilde T$. From this point of view, the propagation of data in Darboux maps and Q-nets is compatible if they are related to each other via a section. 

Let us make a less formal remark on the reverse construction. Note that as sweeps are reversible, one can just as well start with a TCD $\tcd_N^-$ and then do the sweep backwards. If $k$ is the maximal dimension of $\tcd_N^-$, then there are $(k+2)$ points in $T_N^+$. Assume the points of $T^-_N$ span $\CP^j, j\leq k$. We can consider $\CP^j=H'$, where $H'$ is a hyperplane of $\CP^{j+1}$. Then by Theorem \ref{th:tcdmapsweep} $T^+_0$ is a ``lift'' of $T^-_N$ such that $T^-_N$ is a section of $T^+_0$. Thus it is also possible to construct lifts of a section via reverse sweeps.

If the initial TCD map $T^+_0: \tcd^+_0 \rightarrow \CP^n$ spans $\CP^n$ but $n$ is not the maximal dimension that $\tcd^+_0$ allows for, one can consider the case that the dimension of $H$ equals $n$. In order to understand this case, let us begin with a definition.

\begin{definition}\label{def:extension}
	Assume we have two TCD maps
	\begin{align}
		T^\pi: \tcdp^\pi&\rightarrow \CP^{n},\\
		T^\sigma: \tcdp^\sigma&\rightarrow \CP^{n}.
	\end{align}
	Then we call $T^\pi$ an \emph{extension} of $T^\sigma$ if there exists a map $T: \tcdp \rightarrow \CP^{n+1}$ such that
	\begin{align}
		T^\pi &= \pi (T),\\
		T^\sigma &= \sigma_E(T),
	\end{align}
	where $E = \CP^n \subset \CP^{n+1}$ and $\pi$ is a central projection onto $E$.
\end{definition}

Note that $\tcd^\pi = \tcd$ because we are not changing the combinatorics when projecting. Due to the construction, the points of $T^\sigma$ are situated on the lines of $T^\pi$. Thus one can think of the lines of $T^\pi$ as extending the points of $T^\sigma$. Before we relate sweeps to extensions, let us clarify the existence of projection-lifts in a lemma.

\begin{lemma}\label{lem:projectionlift}
	Let $\tcd$ be a minimal TCD with maximal dimension $k > 1$. Let $n<k$ and $S \subset \CP^k$ be an $n$-dimensional subspace and let $P\in \CP^k\setminus S$ be a point. Let $T: \tcdp \rightarrow S, n < k$ be a TCD map. Then there is a TCD map $T': \tcdp \rightarrow \CP^k$ such that the points of $T'$ span $\CP^k$ and $T$ is the projection $\pi_P(T')$ of $T'$ from center $P$. It is possible to choose $T'$ such that $T'(w) \neq T(w)$ for all white vertices $w$.
\end{lemma}
\proof{Fix a labeling of $\tcd$. Consider the li-orientation $\lio$ of $\pb$ and the induced li-poset $\lip$. There are $k+1$ minimal elements $\hat w_1,\hat w_2,\dots \hat w_{k+1}$ in $\lip$. For each $i$ with $1\leq i \leq k$ choose lifts $T'(\hat w_i)$ on the line $PT(\hat w_i)$, but not equal to $P$ or $T(\hat w_i)$ such that the $k+1$ points span $\CP^k$. Choose an arbitrary linear extension 
\begin{align}
	\varepsilon: W(\pb)\setminus \{\hat w_1,\hat w_2,\dots, \hat w_{k+1}\} \rightarrow \N
\end{align}
of $\lip$. We successively go upwards in $\varepsilon$ through the corresponding elements $(w_1,w_2,\dots,w_{|B|})$ to determine the lifts $T'(w_i)$. For each $w_i$ there is one black vertex $b_i$ that is matched to $w_i$ by the li-matching $\lim$, and $b_i$ has two other adjacent white vertices $w'_i, w''_i$. Because we go upwards in $\varepsilon$, the lifts $T'(w'), T'(w'')$ are already determined. Therefore the lift $T'(w)$ is also determined as the intersection of the lines $PT(w)$ and $T'(w')T'(w'')$. This intersection exists because $P,T(w),T'(w'),T'(w'')$ are all on the plane spanned by $P,T(w')$ and $T(w'')$. With this iteration we produce a map $T'$ such that $\pi_P(T') = T$. However, it is still possibl that some of the points of $T'$ coincide with points of $T$. Choose an affine chart of $\CP^n$ in which $P$ is at infinity. Then we can translate all the points of $T'$ in the direction of $P$ (a translation along parallel lines) for an arbitrary amount $t$ to obtain a map $T'_t$. Clearly $\pi_P(T_t') = T$ as well. As there are only a finite number of white vertices, there is a translation $T_t'$ of $T'$ such that $T_t'(w) \neq T(w)$ for all white vertices $w$, which proofs the Lemma. \qed

}

\begin{figure}
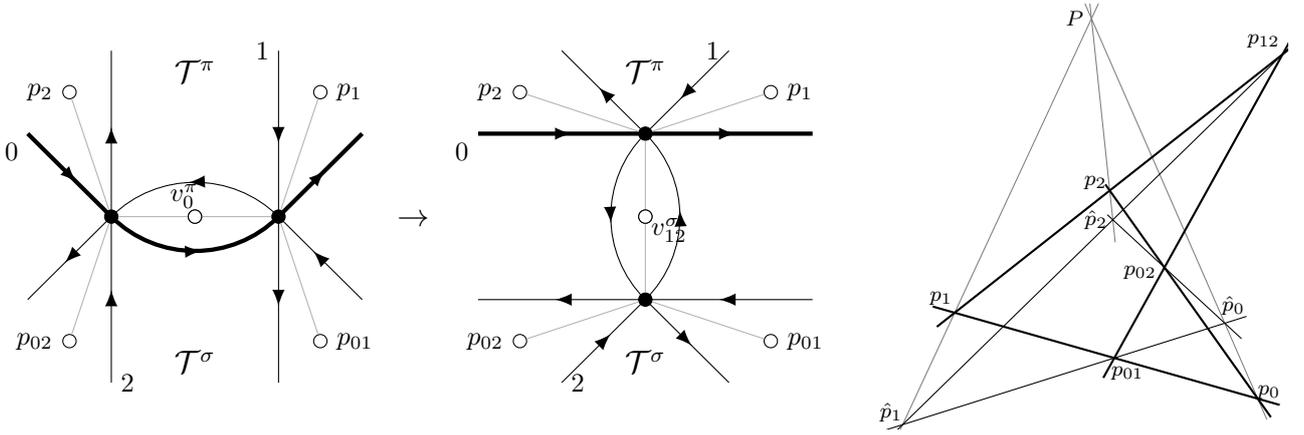

\vspace{-2mm}
\caption{A 2-2 move at a counterclockwise face in a sweep and the geometric picture of the involved points on the right.}
\label{fig:extensionsweep}
\end{figure}

\begin{theorem} \label{th:extension}
	Let $\mathcal U_0, \mathcal U_1,\dots, \mathcal U_N$ be a sweep sequence as introduced in Section \ref{sec:sweeps}. Let $n\in\N$ such that $n$ is smaller than the maximal dimension of $\mathcal U_0$. Consider a sequence of TCD maps $U_0, U_1,\dots, U_N$ such that $U_i$ maps from the counterclockwise faces of $\mathcal U_i$ to $\CP^n$ for $0 \leq i \leq N$, and such that consecutive TCD maps are obtained by the 2-2 moves of the sweep sequence. Write $T^\pm_i$ for the restrictions of $U_i$ to $\tcd^\pm_i$. Assume that $U_0$ is flip-generic, spans $\CP^n$ and that $T_0^+$ and $T_0^-$ and also both span $\CP^n$. Then $T^+_0$ is an extension of $T^-_N$.
\end{theorem}
\proof{Let us embed $\CP^n = H \subset \CP^{n+1}$ and choose a point $P\in \CP^{n+1}\setminus H$. Now, as described in Lemma \ref{lem:projectionlift}, choose a lift $\hat T^+_0:{\tcdp^+}_0\rightarrow \CP^{n+1}$ of $T^+_0$. That is a TCD map such that the image $\hat T^+_0(w)$ of each white vertex $w$ is on the line through $P$ and $T^+_0(w)$ and does not coincide with either $P$ or $T^+_0(w)$. In other words, let $\pi$ be the projection from $P$ to $H$, then we have that $\pi(\hat T^+_0) = T^+_0$. We know from Theorem \ref{th:sweepsection} that if we would start the sweep with $\hat T^+_0$ instead of $T^+_0$ then indeed $T^+_0$ would be the section $\sigma_H(\hat T^+_0)$ of $\hat T^+_0$. Now we claim that the same is true if we start with $T^+_0$. By construction, we know the claim is true on $\mathcal U_0 \setminus \tcd^+_0$. Now we show that if the claim is true at step $i$ of the sweep it is also true at step $i+1$. Consider Figure \ref{fig:extensionsweep} for labels and a picture of the geometry involved. At any counterclockwise 2-2 move in the sweep, there are two points $p_{01}, p_{02}$ that have already been swept and belong to both $T^-_i$ and $T^-_{i+1}$. Similarly, there are two points $p_1,p_2$ that belong to both $T^+_i$ and $T^+_{i+1}$. By construction, there are two points $\hat p_1,\hat p_2$ that are the lifts of $p_1,p_2$ in $\hat T^+_i$. Also by construction, the new point $p_{12}$ that we add to $T^-_{i+1}$ is the intersection of the line $p_1p_2$ with the line $p_{01}p_{02}$ and thus naturally in $H$. Moreover, we claim that $p_{12}$ is also on the line $\hat p_1\hat p_2$, if we show that then the theorem is proven by induction. There is also the point $p_0$ in $T^+_i$ that we remove and its lift $\hat p_0$. We introduce three planes. The plane $A_1$ is spanned by $p_1,p_2,\hat p_1,\hat p_2$ and $P$. The plane $A_2$ is spanned by $p_1,p_2,p_{01},p_{02}$ and $p_0$. The plane $A_3$ is spanned by $\hat p_1, \hat p_2,p_{01},p_{02}$ and $\hat p_0$, this plane exists because of our inductive assumption. But then we see that the lines $p_1p_2, \hat p_1\hat p_2$ and $p_{01}p_{02}$ are the three pairwise intersections of $A_1,A_2,A_3$ and therefore intersect in a single point $p_{12}$.\qed
}

Of course, reading the proof and looking at Figure \ref{fig:extensionsweep} we realize that this is another appearance of Desargues' theorem (see Theorem \ref{th:desargues}). Let us now look at examples of extensions in DDG that have appeared in the literature as well as at new examples that we are able to consider due to our general approach.

\begin{example}[Extending a Q-net to a line complex \cite{ddgbook}]Assume we know the vertices of a Q-net $q$ defined on $\Z^3$, which is the section $T^\sigma$ in this case. Now we are looking for the lines $\ell$ of a line complex (the projection $T^\pi$ in this case) that pass through the points of the Q-net. Let us look at a quad of the Q-net. We claim that if we already know the lines $\ell,\ell_1,\ell_2$ of the line complex associated to three of the points of the quad of a Q-net, then the fourth line $\ell_{12}$ is uniquely determined. The three lines of and the four points of $q$ together span a 3-dimensional space. In this 3-dimensional space $\ell_{12}$ has to contain the point $q_{12}$ and intersect the two lines $\ell_1,\ell_2$. As there is only one such line, $\ell_{12}$ is uniquely determined. Now given a cube in $q$, we may prescribe the lines $\ell,\ell_1,\ell_2,\ell_3$. Then the remaining lines of the cube are determined and one has to check whether $\ell_{123}$ is well-defined, because there are three different ways to construct $\ell_{123}$. Generically, $\ell_{123}$ has to be in a 4-space spanned by $\ell,\ell_1,\ell_2,\ell_3$, but also pass through the three non-intersecting lines $\ell_{12},\ell_{23},\ell_{13}$. There is only one such line, and thus $\ell_{123}$ is unique.
\end{example}

\begin{figure}
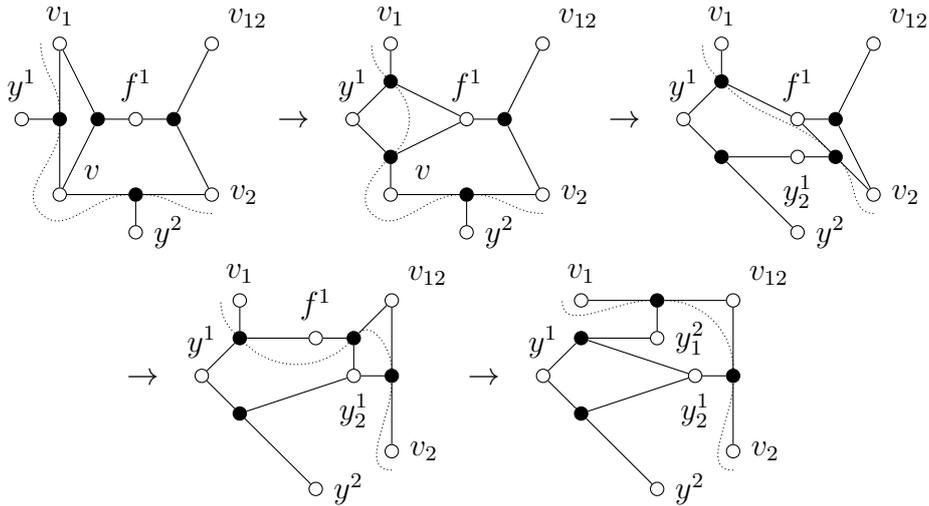


	\vspace{-2mm}
	\caption{Finding a Darboux map that is extended by a Q-net by moving the strand through the diagram.}
	\label{fig:pullingthestring}
\end{figure}

In this example we propagate the information from the three coordinate axes of $\Z^3$ to all of $\Z^3$, but we have to check consistency. On the other hand, we can consider a Q-net as given by some Cauchy data, for example a stepped surface or the coordinate planes. One readily checks that a finite TCD patch covering to a coordinate-plane-corner is sweepable. Therefore we can get Cauchy data for the extending line complex by performing a sweep through the Cauchy data of the Q-net. Then we can propagate Cauchy data of the line complex to obtain the whole line complex. Equivalently we could propagate the data of the Q-net first and then sweep. Again, the consistency of this procedure reduces to the consistency of TCD maps.

We have observed in Section \ref{sec:sections} that in the $\Z^3$-cases the section of a Q-net is a Darboux map, and of a Darboux map the section is a line complex. Together with our method for constructing extensions it is therefore also possible to construct Q-nets that extend Darboux maps, see Figure \ref{fig:pullingthestring} for an example of the combinatorics in that case. Similarly one can construct Darboux maps that extend line complexes.

\begin{figure}
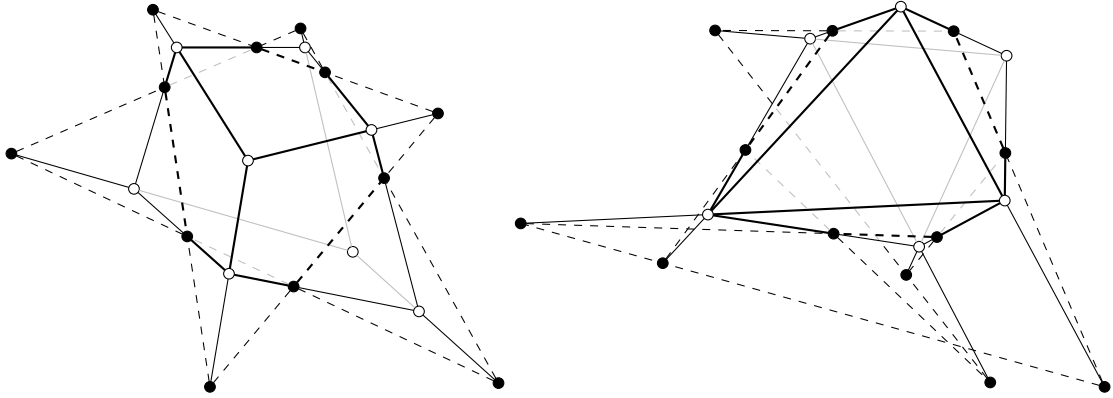

	
	\caption{Left: Drawing of a cube of a Q-net that extends a Darboux map (dashed). Right: Drawing of a cube of an octahedron of a Darboux map that extends a line complex (dashed).}
	\label{fig:extensiongeometry}
\end{figure}

Another natural place where extensions occur is when we do actual drawings of a TCD map and its section, because in this case we have to draw both the map and its section in two dimensions, see Figure \ref{fig:extensiongeometry}.

Another advantage of understanding extensions via TCD maps is that we can even define extensions in $\CP^1$, where no incidence geometry is available anymore. We will see that extensions occur quite naturally in several $\CP^1$ models and often come in the guise of a discrete connection. For example, we will show that there is a relation between the intersection points of a circle pattern and the centers of a circle pattern via an extension in Section \ref{sec:cptemb}.

\section{Desargues maps}\label{sec:desargues}

Desargues maps were introduced by Doliwa \cite{doliwadesargues}. We show that Cauchy-data of Desargues maps constitutes a TCD map. Doliwa proved that Desargues maps are multi-dimensionally consistent due to Desargues theorem.  In our setting that comes as no surprise, as the integrability of all TCD maps follows from Desargues theorem. Apart from capturing another important example of the literature, we also present Desargues maps here because they feature extensions in a quite natural way.

\begin{figure}
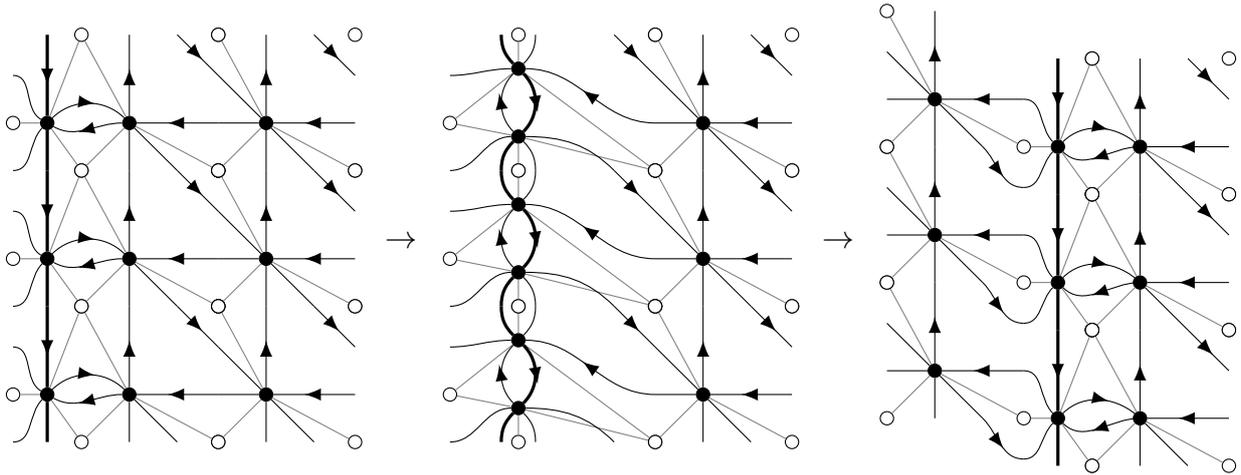


	\caption{The two steps of propagation in a Desargues map. The thick strand is the sweep strand and is pulled through from left to right}
	\label{fig:desarguesextension}
\end{figure}

\begin{definition}
	A map $D: \Z^n \rightarrow \CP^m$ with $n,m > 1$ is called a Desargues map if $D(z)$ and $D(z+\mathbf{e}_1),D(z+\mathbf{e}_2),\dots, D(z+\mathbf{e}_n)$ are colinear for every $z\in\Z$.
\end{definition}

\begin{remark}
	After finishing the thesis, it was brought to our attention that there is a newer, different definition of Desargues maps \cite[Proposition 3.1]{doliwadesarguesweyl}, that contains the original definition as a special case. We only consider Desargues maps as originally defined \cite{doliwadesargues}.
\end{remark}

Note that unlike in the other examples on $\Z^n$ that we considered, the orientations of the coordinate directions of $\Z^n$ matter. For examplem by definition the three points $D,D_1,D_2$ (in shift notation) are on a line, but in general, $D_{12},D_1,D_2$ are not on a line.

\begin{samepage}
\begin{lemma}
	Let $\qg$ be a minimal quad-graph, and identify the strips of $\qg$ with the coordinate directions of $\Z^n$, where parallel strips may correspond to the same coordinate. A Desargues map $\qg \rightarrow \CP^m$ is represented as a TCD map $T$ on the level of the graph $\pb$ by gluing 
	
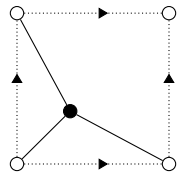
\begin{center}
		\begin{tikzpicture}[>=Triangle]
			\node[wvert] (v) at (0,0) {};
			\node[wvert] (v1) at (2,0) {};
			\node[wvert] (v2) at (0,2) {};
			\node[wvert] (v12) at (2,2) {};
			\node[bvert] (b) at (0.7,0.7) {};				
			\coordinate (d1) at (1.2,0) {};
			\coordinate (d2) at (0,1.2) {};
			\coordinate (d12) at (1.2,2) {};
			\coordinate (d21) at (2,1.2) {};
			\draw[-]
				(b) edge (v) edge (v1) edge (v2)
			;
	
			\draw[densely dotted,->]
				(v) edge (d1) edge (d2)
				(v1) edge (d21)
				(v2) edge (d12)
			;
			\draw[densely dotted]
				(d1) edge (v1)
				(d2) edge (v2)			
				(v12) edge (d12) edge (d21)
			;
		\end{tikzpicture}
	\end{center}
	into every quad. The white vertices are situated at the vertices of the quads and are mapped by $T$ to the corresponding points of the Desargues map. The arrows indicate the edges oriented corresponding to the orientations of the coordinate directions.
\end{lemma}
\end{samepage}
\proof{By construction the three points $D,D_1,D_2$ are on a line in each quad. \qed}

We do not go into more detail here, but there is a cube-flip if the involved strips are not oriented cyclically. The relation between strips and strands is more involved, as is the analysis of whether the corresponding TCD is minimal. Instead, we focus on $\Z^2$ and $\Z^3$ combinatorics. Consider two consecutive $\Z^2$ slices $D$ and $D'$ of a $\Z^3$ Desargues map. By definition, the points of $D'$ are on the lines of $D$. More specifically, the line through $D,D_1,D_2$ contains the point $D_3$. Thus $D$ is an extension of $D'$. Therefore one may ask if there is a suitable sweep and sweep strand to construct $D$ as an extension of $D'$. We recall that in the case of $\Z^2$ Q-nets a sweep strand would typically run along two of the coordinate axes. In the case of Desargues maps however, a sufficient and maximal sweep strand runs only along one coordinate axis. We have depicted the sweep strand and the two steps of the sweep in Figure \ref{fig:desarguesextension}. 

As an additional note, consider a Desargues map defined on $\Z^2$. Then there are three families of parallel strands. Thus such a Desargues map can naturally be associated to the $A_2$ lattice (see the Discussion in Section \ref{sec:tcdconsistency}). In comparison, the initial data of a Q-net defined on $\Z^2$  together with Laplace Darboux dynamics (see Section \ref{sec:laplacedarboux}), can be associated to $A_3$. The data of Q-nets, Darboux maps and line complexes defined on a stepped surface together with propagation of these maps, can be associated naturally to $A_5$.

\section{Projective duality and TCD maps}\label{sec:projduality}

\begin{figure}
		
	\caption{The sequence $\pb \rightarrow \pb^i \rightarrow \pb' = \eta(\pb)$ or how to take a line dual of a TCD map.}
	\label{fig:vrclinedual}
\end{figure}

In this section we explain a way how to construct a TCD map $T^\star:\tcdp^\star \rightarrow (\CP^n)^*$ that can be considered to be the projective dual of a given TCD map $T:\tcdp \rightarrow \CP^n$. An important role is played by the lines that are associated to the black vertices of the bipartite graphs. Before we proceed, let us clarify some conventions around projective dualization. By $(\CP^n)^*$ we denote the projective dual space of $\CP^n$. Every point $P$ in $(\CP^n)^*$ corresponds to a hyperplane $H_P$ of $\CP^n$. We denote this correspondence by writing $P=(H_P)^\perp$ as well as $H_P=(P)^\perp$. Moreover $(\CP^n)^*$ is itself isomorphic to $\CP^n$, because we also identify every $k$-space $S_k \subset \CP^n$ with the space of all hyperplanes containing $S_k$, which is the $(n-k-1)$-space $S_k^\perp \subset (\CP^n)^*$. To reduce confusion, we denote dual TCD maps with a star by $T^\star$ instead of $T^*$.

\begin{definition}\label{def:linemap}
	Let $T: \tcdp \rightarrow \CP^n$ be a TCD map. Let the \emph{associated line map} be the map $L: B(\pb) \rightarrow \{\mbox{Lines of } \CP^n\}$ such that for every black vertex $b\in B(\pb)$ holds that
	\begin{align}
		L(b) = \spa\{T(w),T(w')\},
	\end{align}
	where $w,w'$ are two different white vertices adjacent to $b$.
\end{definition}

 As the procedure for constructing the dual TCD map is somewhat involved, we find it instructive to begin with the case where the ambient dimension equals $n=2$. In $\CP^2$ we can identify the dual space $(\CP^{2})^*$ with the space of lines in $\CP^2$. Thus one may ask whether we can capture the lines occurring in $T$ and their relations as a TCD map $T^\star$. We give an affirmative answer by giving a direct operation on the associated graph $\pb$ of $T$. Before we proceed, we need the appropriate concept of genericity.

\begin{definition}
	Let $T:\tcdp \rightarrow \CP^n$ be a TCD map. Consider triplets $(w,w_1,w_2)$ of all different white vertices of $\pb$, such that
	\begin{enumerate}
		\item there is a face $f_1$ adjacent to $w$ and $w_1$,
		\item there is a face $f_2$ adjacent to $w$ and $w_2$,
		\item there is no black vertex adjacent to all three $w,w_1$ and $w_2$.
	\end{enumerate}
	 We say $T$ is \emph{cogeneric} if for any such triplet the points $T(w),T(w_1),T(w_2)$ span a plane.
\end{definition}

Note that $f_1 = f_2$ is allowed, also note that there is also a reasonable concept of $k$-cogenericity, but we do not need it.

\begin{definition}\label{def:linedual}
	Let $T: \tcdp \rightarrow \CP^2$ be a cogeneric TCD map. 
	A \emph{line dual} $\eta(T)$ of $T$ is a TCD map $T': \tcdp'\rightarrow  (\CP^{2})^*$ constructed from $T$ as follows. We define the change of combinatorics on the level of graphs $\pb$ and $\pb'$. We begin by constructing an intermediate graph $\pb^i$ starting from $\pb$ in two steps (the same first two steps as in Definition \ref{def:section}, see also Figure \ref{fig:vrclinedual}):
	\begin{enumerate}
		\item Add a black vertex $b_{k,k+1}$ for any two consecutive white boundary vertices $w_k,w_{k+1}$ and the two edges $(b_{k,k+1},w_k)$ and $(b_{k,k+1},w_{k+1})$.
		\item For each face $(b_1,w_1,b_2,\dots,w_l)$ of $\pb$, triangulate the polygon $(w_1,w_2,\dots,w_l)$ and for each diagonal $(w_j,w_k)$ of the triangulation add a new black vertex $b$ as well as the two edges $(b,w_j),(b,w_k)$.
	\end{enumerate}
	Every interior face of $\pb^i$ is a quad or a hexagon. Now we construct $\pb'$ in the following steps:
	\begin{enumerate}
		\setcounter{enumi}{2}	
		\item Contract each set of black vertices in $\pb^i$ that represent the same line to a single black vertex $b$. This is our starting diagram for $\pb'$. In other words, replace each maximal subgraph $U \subset \pb$ that corresponds to a TCD with endpoint matching $\enm {k}{1}$ for some $k\in \N$ with a black vertex $b$ that is connected with an edge to each boundary white vertex of $U$.
		\item Swap the colors of black and white vertices in $\pb'$. Thus each white vertex $w'$ was a black vertex $b$ before. Set $T'(w') = (L(b))^\perp$.
		\item Replace every black vertex $b'$ of degree $d_{b'}$ more than 2 in $\pb'$ by a graph piece, which corresponds to a TCD with endpoint matching $\enm {d_{b'}}{1}$ with $d_{b'}$ white boundary vertices, which we identify with the neighbours of $b'$.
		\item Erase all black vertices of degree 2 in $\pb'$ as well as the edges incident to the erased black vertices.\qedhere
	\end{enumerate}
\end{definition}

Steps (1) and (2) do the same as in the case of sections, they store information about the existence of lines in $T$ that were not represented as a black vertex in $\pb$. Steps (3) and (4) represent lines spanned by points in $T$ as new white vertices in $\pb'$. Step (4) also represents points in $T$ as lines in $T'$. Step (5) ensures that all former lines through a point in $T$ are represented by points that are on a line in $T'$. This is because TCDs with endpoint matching $\enm {k}{1}$ represent points on a line, as the maximal dimension of such a TCD is 1. Note that the 1-cogenericity of $T$ ensures that the resulting TCD map $T'$ is 0-generic. Black vertices of degree 2 represent white vertices in $\pb^i$ that are only on two lines in $T$ and do not add relations to $T'$. Therefore black vertices of degree 2 are erased in step (6).

Before we continue to discuss the combinatorics of the line dual, let us give some context. For a TCD map $T$ in $\CP^2$ as defined in \ref{def:linedual}, the line dual $\eta(T)$ is already the dual map $T^\star$. It clearly catches the relations of the lines of $T$ in dual space. In higher dimensions, the line dual is not itself the dual TCD map, but it is an important ingredient in the construction of the dual TCD map as we will see later in this section.

Recall that in Section \ref{sec:sections} we defined the section $\sigma(T)$ of a TCD map via an operation on $\pb$. Then, surprisingly Lemma \ref{lem:sectionviatcd} showed that this operation corresponds to an operation on the strands of $\tcd$ that does not cut or glue strands. We now show that, also surprisingly, the same holds for the line dual. Recall that in the case of sections we introduced the alternating strand diagram $\alt(\tcd)$ in Definition \ref{def:alt} and the shadow $(\altg,C_\alt)$ in Definition \ref{def:shadowgraph}.

\begin{figure}
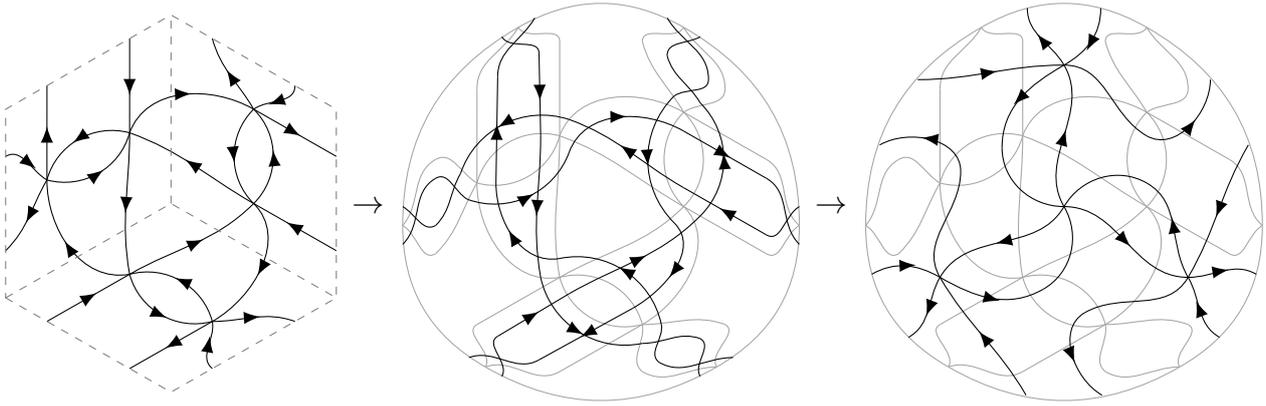

	
	\caption{The sequence $\tcd \rightarrow \alt$ (black) and $\altg$ (gray) $\rightarrow\tcd(\altg,\bar C_\alt)=\eta(\tcd)$ or how to take a line dual of a TCD.}
	\label{fig:tcdlinedual}
\end{figure}

\begin{definition}\label{def:revshadow}
	Let $\tcd$ be a minimal TCD and let $(\altg,C_\alt)$ be the corresponding shadow. Consider the endpoint matching $C^f_\alt$ in a face of $\altg$. Replace $C^f_\alt$ with the endpoint matching $\bar C^f_\alt$ defined by
	\begin{enumerate}
		\item swapping in- and out-endpoint on each edge,
		\item reversing the orientation of all strands,
		\item for every edge $e$ of $\altg$, if $e$ is an edge between $f$ and a clockwise bigon of $\alt$ then \emph{short circuit} the edge $e$. That is, if after the first two steps a strand in $f$ starting at in-endpoint $i$ ends at $e$, passes through the bigon, returns to $f$ and then leaves $f$ at out-endpoint $j$ then define $\bar C^f_\alt(i) = j$ instead. There are no endpoints on $e$ anymore. 
	\end{enumerate}
	We call the resulting pair $(\altg,\bar C_\alt)$ the \emph{reversed shadow} (see Figure \ref{fig:tcdlinedual}).
\end{definition}
In step (1) we really slide the two endpoints over each other such that in the result they are swapped. Moreover, as a result of step (3) there are no strands in clockwise bigons in the reversed shadow. Recall that in Definition \ref{def:shadowtotcd} we defined how to obtain a TCD from a shadow.
\begin{lemma}
	Let $\tcd$ be a minimal TCD and $T:\tcdp \rightarrow \CP^2$ be a TCD map. Then the TCD $\eta(\tcd)$ of the line dual $\eta(T)$ corresponds to a TCD constructed from the reversed shadow. Thus, one can make choices such that
	\begin{equation}
		\eta(\tcd) = \tcd(\altg,\bar C_\alt).\qedhere
	\end{equation}
\end{lemma}
\proof{
We compare the construction of $\eta(\tcd)$ to the graph construction of $\eta(\pb)$ in Definition \ref{def:linedual}, see Figure \ref{fig:tcdlinedual} for an example of $\eta(\tcd)$ and Figure \ref{fig:vrclinedual} for the corresponding $\eta(\pb)$. When splitting triple intersection points to obtain the alternating strand diagram $\alt(\tcd)$ we perturb all strands such that a new counterclockwise face emerges. Afterwards we reverse the orientation of the strands but also swap them on each edge of $\altg$, thus it remains a counterclockwise face that corresponds to the white vertex that we place there in graph step (4). Similarly as in the case of sections the swaps at the boundary together with reversing the orientation of the strands corresponds to introducing white vertices at the boundary. The TCDs that we glue in each face of the reversed shadow $(\altg,\bar C_\alt)$ that was oriented counterclockwise in $\tcd$ have endpoint matching $\enm {k}{1}$ for some $k\in \N$. This captures that we split the black vertices in $\pb^i$ such that the adjacent white vertices are on a line. The TCDs that we glue in each face of the reversed shadow $(\altg,\bar C_\alt)$ that was oriented clockwise in $\tcd$ have endpoint matching $\enm {k}{-2}$ for some $k\in \N$. This corresponds to the triangulation that we choose in graph step (2) and the extra points that we introduced. We remove strands of length 0 because they do not contribute information to the TCD map and do not appear in the TCD associated with $\pb$.\qed
}

If for a moment we do not remove strands of length 0, then there is a simple formula for the endpoint matching of the line dual:
\begin{align}
 	C_\tcd(k) = l \quad \Rightarrow \quad C_{\eta(\tcd)}(l-2) = k,
\end{align}
where we assume that the strand starting at endpoint $k$ in $\tcd$ and the strand ending at out-endpoint $k$ in $\eta(\tcd)$ are the same ones.

\begin{corollary}
	Let $\tcd$ be a minimal TCD with minimal length at least $2$. If the maximal dimension of $\tcd$ is $m$ then the maximal dimension of $\eta(\tcd)$ is $n-m+1$.
\end{corollary}
\proof{On the boundary we swap the strands twice, in effect taking the maximal dimension from $m$ to $(m-2)$ by Corollary \ref{cor:maxdimsection}. Additionally, we are reversing the orientations of the strands, swapping the number of left and right moving strands. Because we have to move the first out-endpoint to the very right in our half-plane drawing, the total effect of the orientation reversal is $m'\mapsto n-m'-1$. Together we have $m\mapsto n-(m-2)-1$ which proves the claim.\qed }

\begin{lemma}
	$\eta(\tcd)$ is connected and minimal.
\end{lemma}
\proof{Because we are not cutting or inserting any strands $\eta(\tcd)$ is connected. For minimality we follow the arguments of the proof of Lemma \ref{lem:sectionminimal}. First we argue that no parallel intersections exist. As in the case of a section, between any existing intersection points of two strands we are producing new intersection points only in the correct order and thus no parallel intersections can occur. The argument that the strands bounding a face of $\tcd$ are all distinct (see proof of Lemma \ref{lem:sectionminimal}) holds for counterclockwise faces as well. Thus, there are no self-intersections in $\eta(\tcd)$ and the Lemma is proven.\qed}

\begin{lemma}
	Let $\tcd$ be a minimal TCD with minimal length at least $2$. Then one can make choices in both line duals such that $\eta\circ \eta(\tcd) = \tcd$.
\end{lemma}
\proof{On the level of the $\pb$ we observe that $\eta$ replaces every set of black vertices that corresponds to a single line with a white vertex and every white vertex with a set of black vertices that corresponds to a single line. The composition $\eta \circ \eta$ then reverses these replacements. White vertices of degree 2 that are removed can be recovered by making the right choice of triangulation in $\eta(\tcd)$. White vertices that were introduced for triangulation edges in $\eta(\tcd)$ are removed again in $\eta\circ\eta(\tcd)$. On the boundary the minimal length ensures that no strands are removed. To see that $\eta\circ \eta(\tcd) = \tcd$ recall that we introduce additional points in $\eta$ for each two consecutive boundary points.\qed }

So far we have only studied line duals for TCD maps taking values in $\CP^2$. In order to obtain a definition in higher dimensions, we have to put in additional work into the combinatorics. Let us define a simple global operation on TCDs.
\begin{definition}\label{def:tcddual}
	Let $\mathcal T$ be a triple crossing diagram. The \emph{TCD dual $\iota (\mathcal T)$} is the same as $\mathcal T$ but with all strand orientations reversed.
\end{definition}
Clearly the TCD dual is an involution. The black vertices of the corresponding graph $\pb$ are unchanged, while faces and white vertices are interchanged. However, there is no obvious direct geometric way to realize $\iota(\mathcal T)$ as a TCD map. Let us consider some examples.

\begin{example}\label{ex:tcdprojdual}
	There are choices such that
	\begin{enumerate}
		\item if $\tcd$ is the TCD of a Q-net, then $\iota(\tcd)$ is the TCD of a Darboux map and vice versa;
		\item if $\tcd$ is the TCD of a line complex, then $\iota(\tcd)$ is the TCD of a line compound and vice versa;
		\item if $\tcd$ is the TCD of a Q-net, then $\eta(\tcd)$ is the TCD of a line complex and vice versa;
		\item if $\tcd$ is the TCD of a Darboux map, then $\eta(\tcd)$ is the TCD of a line compound and vice versa.\qedhere
	\end{enumerate}
\end{example}

\begin{figure}
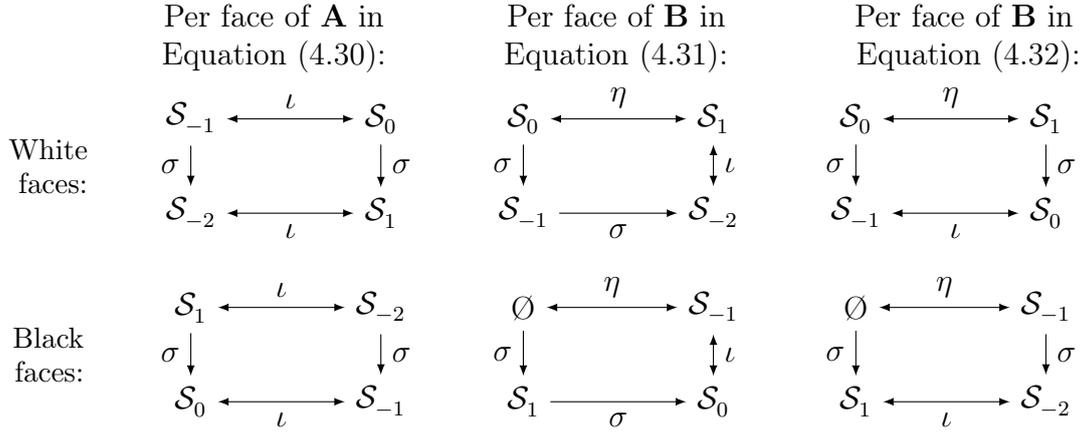

	\hspace{1.8cm}
\begin{tabular}{ccc}
	Per face of $\mathbf{A}$ in\hspace{12mm} & Per face of $\mathbf{B}$ in  & \hspace{12mm}Per face of $\mathbf{B}$ in\\
	Equation \eqref{eq:sigmaiota}:\hspace{12mm} & Equation \eqref{eq:sigmasigma}:  & \hspace{12mm}  Equation \eqref{eq:sigmaeta}:
	\end{tabular}%
	\vspace{-2mm}
	\caption{Commuting diagrams underlying Lemma \ref{lem:iotaetasigma}.}
	\label{fig:iotaetasigma}
\end{figure}

The following relations between $\iota,\eta$ and $\sigma$ are certainly interesting in their own right, but will also be essential to understanding the projective dual of a TCD map.

\begin{lemma}\label{lem:iotaetasigma}
	Let $\tcd$ be a minimal TCD with minimal length at least $2$. In order for the equation
	\begin{align}
		\iota \circ \sigma \circ \iota \circ \sigma (\mathcal T) &= \mathcal T, \label{eq:sigmaiota}
	\end{align}
	to be satisfied, the choices made in one of the two sections determine the choices in the other section. Additionally, for any choices on the left hand side of equations
	\begin{align}		
		\sigma \circ \sigma (\mathcal T) &= \iota \circ \eta (\mathcal T), \label{eq:sigmasigma}\\
		\iota \circ \sigma (\mathcal T) &= \sigma \circ \eta (\mathcal T), \label{eq:sigmaeta}
	\end{align}
	there are choices on the right hand side that satisfy the equations.
\end{lemma}
\proof{
	The proofs of all three equations involve partitioning the TCDs into the faces of auxiliary graphs and then tracking the endpoint matchings and choices per face, see Figure \ref{fig:iotaetasigma} for the underlying commuting diagrams of endpoint matchings.
	
	Let us prove Equation \eqref{eq:sigmaiota} first. We construct an auxiliary planar graph $\mathbf A$. The vertex set of $\mathbf A$ are the white vertices of $\pb$, and there is an edge $(w,w')$ in $\mathbf A$ if there is exactly one black vertex in $\pb$ adjacent to both $w$ and $w'$. Consider the restriction of $\tcd$ to the faces of $\mathbf A$. Then there are two types of faces in $\mathbf A$. We call a face black if it contains black vertices of $\pb$, in which case it corresponds to a TCD with endpoint matching $\enm{}1$. The other type of face contains no vertices of $\pb$, in which case we call it a white face and these faces correspond to endpoint matching $\enm{}{-1}$. Let $\tcd'=\sigma(\tcd)$ and consider the restriction of $\tcd'$ to the faces of $\mathbf A$. The black faces of $\mathbf A$ contain exactly one white vertex of $\tcd'$, which corresponds to endpoint matching $\enm{}{0}$. The white faces of $\mathbf A$ on the other hand contain TCDs of endpoint matching $\enm{}{-2}$, the configuration of $\tcd'$ depending on the choices made in taking the section. Let $\tcd'' = \iota(\tcd')$ and consider the restriction of $\tcd''$ to the faces of $\mathbf A$. Then the black faces of $\mathbf A$ contain TCDs with endpoint matching $\enm{}{-1}$, corresponding to faces of $\pb''$. The white faces of $\mathbf A$ contain TCDs with endpoint matching $\enm{}{1}$. Thus, the characterization via endpoint matchings of white and black faces of $\mathbf A$ with respect to $\pb''$ is reversed in comparison to $\pb$. We can therefore repeat the previous steps to obtain $\tcd''' = \sigma(\tcd'')$ and $\tcd'''' = \iota(\tcd''')$. Then the endpoint matchings in the faces of $A$ agrees with respect to $\tcd''''$ and $\tcd$. For the content to match exactly in the faces of $\mathbf A$, the choices in the second section are determined uniquely. This concludes the proof of Equation \eqref{eq:sigmaiota}.
	
	Let us turn to Equation \eqref{eq:sigmasigma}. Recall that in both Definition \ref{def:section} of the section and Definition \ref{def:linedual} of the line dual we began by constructing the intermediary graph $\pb^i$, which involved some choices. Assume we fix one choice of intermediary graph $\pb^i$. From $\pb^i$ we construct an auxiliary graph $\mathbf B$. The vertex set of $\mathbf B$ are the black vertices of $\pb^i$, and there are as many edges $(b,b')$ in $\mathbf B$ as there are white vertices in $\pb^i$ adjacent to both $b$ and $b'$. We distinguish two types of faces of $\mathbf B$, the white face faces that contain white vertices of $\pb^i$ and the black faces that contain no vertices of $\pb^i$. Note that the black faces have only degree three or degree two. Let us first consider the restriction of $\pb^i$ to $\mathbf B$. Then the white faces correspond to TCDs with endpoint matching $\enm{}{0}$ while the black faces are empty. The choices in $\sigma(\tcd)$ correspond to the choice of intermediate graph, we assume that the intermediate graph is $\pb^i$ as fixed above. Let us now consider the restriction of $\sigma(\tcd)$ to the faces of $\mathbf B$. In this restriction, the white faces of $\mathbf B$ correspond to endpoint matchings $\enm{}{-1}$ and the black faces correspond to endpoint matchings $\enm{}{1}$. Next, consider the restriction of $\sigma(\sigma(\tcd))$ to the faces of $\mathbf B$. In this restriction, the white faces of $\mathbf B$ correspond to endpoint matchings $\enm{}{-2}$ and the black faces correspond to endpoint matchings $\enm{}{0}$. The choices in $\sigma(\sigma(\tcd))$ are the choices of the configuration in the white faces of $\mathbf B$. Let us turn to the line dual $\eta(\tcd)$. We use the same intermediary graph $\pb^i$ for the construction of $\eta(\tcd)$. Consider the restriction of $\eta(\tcd)$ to $\mathbf B$. The white faces of $\mathbf B$ correspond to endpoint matchings $\enm{}{1}$ and the black faces of $\mathbf B$ correspond to endpoint matchings $\enm{}{n-1}$. The additional choices in $\eta(\tcd)$ correspond to the choice of configuration in the white faces of $\mathbf B$. We observe that in the restriction of $\iota(\eta(\tcd))$, the white faces correspond to endpoint matchings $\enm{}{n-2}$ and the black faces of $\mathbf B$ correspond to endpoint matchings $\enm{}{0}$. Therefore the endpoint matchings in the restrictions of $\sigma(\sigma(\tcd))$ and $\iota(\eta(\tcd))$ to the faces of $\mathbf B$ coincide. For the TCDs to coincide as well, we have to make the corresponding choices in the white faces of $\mathbf B$ when taking the second section in $\sigma(\sigma(\tcd))$ and when choosing the $\enm{}{1}$ matchings in $\eta(\tcd)$. This concludes the proof of Equation \eqref{eq:sigmasigma}.
	
	The proof of Equation \eqref{eq:sigmaeta} proceeds in the same manner as the proofs Equation \eqref{eq:sigmasigma}. We also employ $\mathbf B$ as auxiliary graph and then trace the choices and endpoint matchings in the faces of $\mathbf B$ constructed from a fixed choice of intermediary graph $\pb^i$. We also assume that $\pb^i$ is fixed and used for both $\eta(\tcd)$ and $\sigma(\tcd)$. The additional choices (beyond choosing $\pb^i$) in $\eta(\tcd)$ do not matter this time, because they are canceled in $\sigma(\eta(\tcd))$. The endpoint matchings per white face of $\mathbf B$ with respect to both $\sigma(\eta(\tcd))$ and $\iota(\sigma(\tcd))$ are $\enm{}{0}$, and per black face they are $\enm{}{-2}$. We only have to take care that the choices when taking the section of $\eta(\tcd)$ correspond to the 2-valent white vertices of $\pb^i$ that were lost in $\eta(\tcd)$. This concludes the proof of Equation \eqref{eq:sigmaeta}.\qed	
}

\begin{figure}
\begin{tabular}{cc|cc}
\begin{tikzpicture}[scale=.8,baseline={([yshift=-.7ex]current bounding box.center)}]
	\coordinate (ws) at (202.5:2);
	\coordinate (wn) at (157.5:2);
	\coordinate (es) at (-22.5:2);
	\coordinate (en) at (22.5:2);
	\coordinate (nw) at (112.5:2);
	\coordinate (ne) at (67.5:2);
	\coordinate (sw) at (247.5:2);
	\coordinate (se) at (292.5:2);
	\coordinate (n) at ($(wn)!.5!(en)$);
	\coordinate (s) at ($(ws)!.5!(es)$);
	\draw[dotted] (0,0) circle (2);
	\draw[-]
		(n) edge[bend right=45,mid arrow] (s) edge[mid rarrow] (ne) edge[mid arrow] (en) edge[mid rarrow] (wn)
		(s) edge[mid arrow] (se) edge[mid rarrow] (es) edge[mid arrow] (ws)
		(s) edge[bend right=45,mid arrow] (n) edge[mid rarrow] (sw)
		(n) edge[mid arrow] (nw)
	;
	\node[bvert] (n) at (n) {};
	\node[bvert] (s) at (s) {};
	\node[wvert] (vnw) at (135:2) {};
	\node[wvert] (vne) at (45:2) {};
	\node[wvert] (vsw) at (225:2) {};
	\node[wvert] (vse) at (315:2) {};
	\node[wvert] (c) at (0,0) {};
	\draw[-, gray!60]
		(n) edge (vnw) edge (vne) edge (c)
		(s) edge (vsw) edge (vse) edge (c)
	;
	\node (arr) at (270:2.5) {\hspace{2mm}$\downarrow$\hspace{.5mm}$\sigma$};
	\draw[<->] (0:2.3) -- node[above] {$\eta$} (0:3);
\end{tikzpicture} &
\hspace{-2mm}\begin{tikzpicture}[scale=.8,baseline={([yshift=-.7ex]current bounding box.center)}]
	\coordinate (ws) at (202.5:2);
	\coordinate (wn) at (157.5:2);
	\coordinate (es) at (-22.5:2);
	\coordinate (en) at (22.5:2);
	\coordinate (nw) at (112.5:2);
	\coordinate (ne) at (67.5:2);
	\coordinate (sw) at (247.5:2);
	\coordinate (se) at (292.5:2);
	\coordinate (n) at ($(wn)!.5!(en)$);
	\coordinate (s) at ($(ws)!.5!(es)$);
	\draw[dotted] (0,0) circle (2);
	\draw[-]
		(nw) edge[bend right=60, mid arrow] (ne)
		(se) edge[bend right=60, mid arrow] (sw)
		(ws) edge[bend right=60, mid arrow] (wn)
		(en) edge[bend right=60, mid arrow] (es)
	;
	\node[wvert] (vw) at (180:2) {};
	\node[wvert] (ve) at (0:2) {};
	\node[wvert] (vs) at (270:2) {};
	\node[wvert] (vn) at (90:2) {};
	\node (arr) at (270:2.5) {\hspace{2mm}$\downarrow$\hspace{.5mm}$\sigma$}; 	
\end{tikzpicture}\hspace{.4cm} & \hspace{.4cm}
\begin{tikzpicture}[scale=.8,rotate=0,baseline={([yshift=-.7ex]current bounding box.center)}]
	\coordinate (ws) at (202.5:2);
	\coordinate (wn) at (157.5:2);
	\coordinate (es) at (-22.5:2);
	\coordinate (en) at (22.5:2);
	\coordinate (nw) at (112.5:2);
	\coordinate (ne) at (67.5:2);
	\coordinate (sw) at (247.5:2);
	\coordinate (se) at (292.5:2);
	\coordinate (n) at ($(wn)!.5!(en)$);
	\coordinate (s) at ($(ws)!.5!(es)$);
	\draw[dotted] (0,0) circle (2);
	\draw[-]
		(n) edge[bend right=45,mid rarrow] (s) edge[mid arrow] (ne) edge[mid rarrow] (en) edge[mid arrow] (wn)
		(s) edge[mid rarrow] (se) edge[mid arrow] (es) edge[mid rarrow] (ws)
		(s) edge[bend right=45,mid rarrow] (n) edge[mid arrow] (sw)
		(n) edge[mid rarrow] (nw)
	;
	\node[bvert] (n) at (n) {};
	\node[bvert] (s) at (s) {};
	\node[wvert] (vw) at (180:2) {};
	\node[wvert] (ve) at (0:2) {};
	\node[wvert] (vs) at (270:2) {};
	\node[wvert] (vn) at (90:2) {};
	\draw[-, gray!60]
		(n) edge (vn) edge (vw) edge (ve)
		(s) edge (vs) edge (vw) edge (ve)
	;
	\node (arr) at (270:2.5) {\hspace{2mm}$\downarrow$\hspace{.5mm}$\sigma$}; 
	\draw[<->] (0:2.3) -- node[above] {$\eta$} (0:3);
\end{tikzpicture} &
\hspace{-2mm}\begin{tikzpicture}[scale=.8,baseline={([yshift=-.7ex]current bounding box.center)}]
	\coordinate (ws) at (202.5:2);
	\coordinate (wn) at (157.5:2);
	\coordinate (es) at (-22.5:2);
	\coordinate (en) at (22.5:2);
	\coordinate (nw) at (112.5:2);
	\coordinate (ne) at (67.5:2);
	\coordinate (sw) at (247.5:2);
	\coordinate (se) at (292.5:2);
	\coordinate (n) at ($(wn)!.5!(en)$);
	\coordinate (s) at ($(ws)!.5!(es)$);
	\draw[dotted] (0,0) circle (2);
	\draw[-]
		(nw) edge[bend right=60, mid rarrow] (ne)
		(se) edge[bend right=60, mid rarrow] (sw)
		(ws) edge[bend right=60, mid rarrow] (wn)
		(en) edge[bend right=60, mid rarrow] (es)
	;
	\node[wvert] (c) at (0,0) {};
	\node (arr) at (270:2.5) {\hspace{2mm}$\downarrow$\hspace{.5mm}$\sigma$}; 	
\end{tikzpicture}\\
\end{tabular} 
\caption{Two exceptional examples for the commuting diagram of $\sigma \circ \eta = \iota\circ \sigma$. On the right hand side we start with strands of length less than 2. Yet, if we do not remove strands of length zero, the diagram still commutes.}
\label{fig:tcddualcommutationquad}
\end{figure}
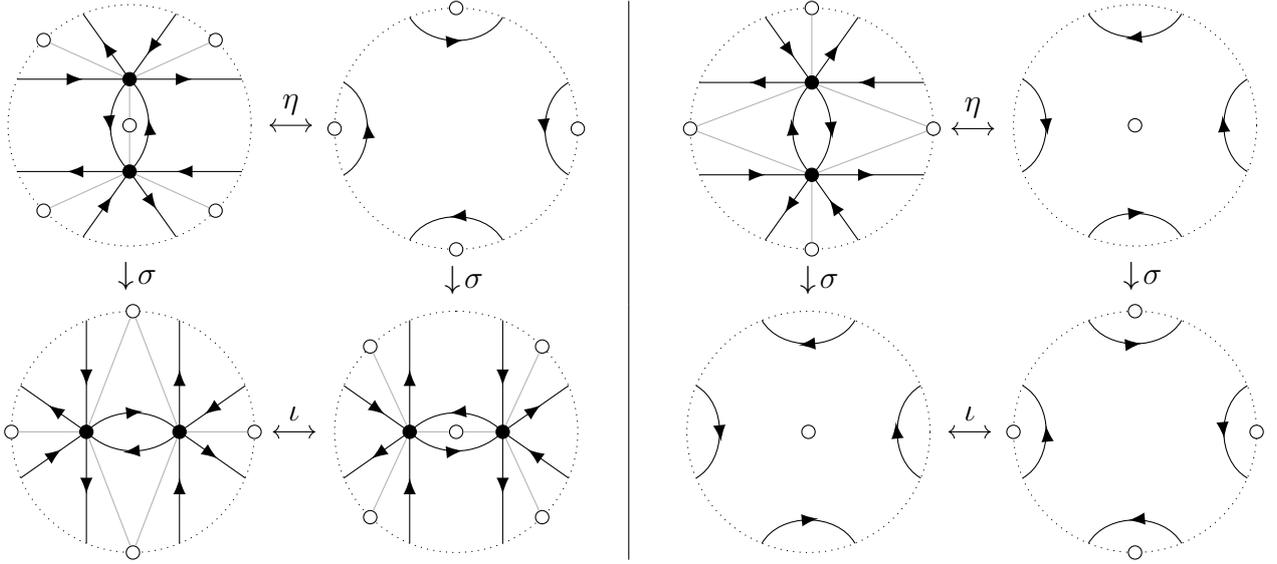

Note that there are examples (see the right hand side of Figure \ref{fig:tcddualcommutationquad}) with minimal length less than two, where Equation \eqref{eq:sigmaeta} still holds if we do not remove strands of length 0.  We expect the whole of Lemma \ref{lem:iotaetasigma} to hold without the minimal length constraint if one characterizes $\sigma$ and $\eta$ via the TCDs. However, we do not investigate this further as the geometric meaning in this case is unclear, but this may still be of interest from a purely combinatorial perspective.

In order to study the geometry of the projective dual for TCD maps in projective spaces of higher dimensions, we introduce the concept of a flag of TCD maps.
\begin{definition}\label{def:tcdflag}
	Let $n\in \N$. A \emph{flag of TCD maps} is a pair of a	collection of projective spaces $(E_k)_{0 \leq k \leq n}$ and a collection of TCD maps $(T_k)_{1 \leq k \leq n}$ such that
	\begin{enumerate}
		\item $\dim E_k = k$ for $0\leq k \leq n$,
		\item $E_k \subset E_{k+1}$ for $0\leq k \leq n-1$,
		\item $T_k: {\tcdp}_k \rightarrow E_k$ for $1\leq k \leq n$,
		\item $T_k = \sigma_{E_k}( T_{k+1})$ for $1\leq k \leq n-1$.\qedhere
	\end{enumerate}
\end{definition}

Implicitly this definition assumes that $E_k$ is generic with respect to $T_{k+1}$ for $1\leq k \leq n$, because otherwise the sections $\sigma_{E_k}(T_{k+1})$ are not defined. Not every TCD map $T_n$ is part of a flag of TCD maps, but if $T_n$ is $(n-1)$-generic, then it is, due to the definition of genericity (Definition \ref{def:kgeneric}). 

Consider a TCD $\tcd$ and its section $\tcd' = \sigma(\tcd)$ as well as the corresponding graphs $\pb,\pb'$. Recall that every black vertex of the intermediate graph $\pb^i$ is replaced by a white vertex $w_b$ in $\pb'$, see Definition \ref{def:section}. Conversely, for every white vertex $w$ in $\pb'$ there is at least one black vertex in $\pb^i$ that was replaced by $w$. For every white vertex $w$ in $\pb'$ let us denote by $b_w$ one of these black vertices of $\pb^i$.

\begin{definition}\label{def:subspacemap}
	Let $n\in \N$ and $n \geq 2$ and consider a flag of TCD maps $(T_k)_{1 \leq k \leq n}$, $(E_k)_{0 \leq k \leq n}$. For $1\leq k < n$ and each white vertex $w$ of $\pb_{k-1}$ choose two different white vertices $v_w,v'_w$ of $\pb_k$ such that $v_w$ and $v'_w$ are adjacent to $b_w$. Define the \emph{subspace maps}
	\begin{align}
		U_k: \tcdp_k \rightarrow \{(n-k)-\mbox{subspaces of } E_n \},
	\end{align}
	for $1 \leq k \leq n$ such that $U_n = T_n$ and such that $U_{k-1}(w) = \spa\{U_k(v_w),U_k(v_w')\}$ for every white vertex of $\pb_{k-1}$.
\end{definition}

Not only is $U_n = T_n$ but also $U_{n-1}(w) = L(b_w)$, where $L$ is the associated line map (see Definition \ref{def:linemap}). Also note that $U_{k-1}(w)$ is well-defined because the choice of $v_w,v'_w$ does not matter, as other choices define the same span.

\begin{figure}
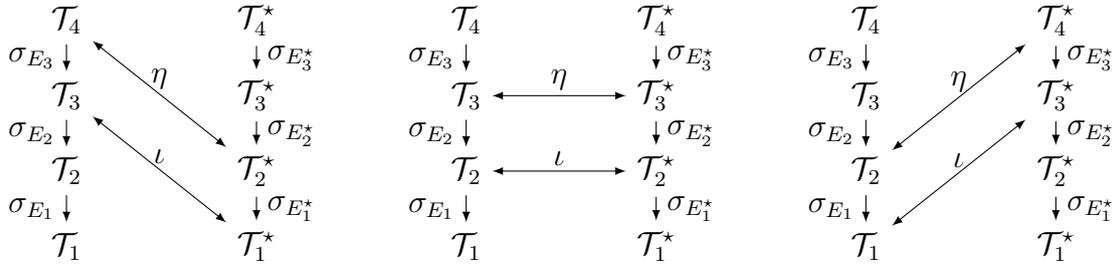


	\caption{The different ways that the TCD dual $\iota$ and the line dual $\eta$ relate TCDs of the projective dual flag of a flag of TCD maps in $\CP^4$.}
	\label{fig:tcddualcombexample}
\end{figure}

\begin{definition}\label{def:dualflag}
	Let $n\in \N$ and $n \geq 2$. A \emph{projective dual} of a flag of TCD maps $(T_k)_{1 \leq k \leq n}$, $(E_k)_{0 \leq k \leq n}$ is a flag of TCD maps $(T_k^\star)_{1 \leq k \leq n}$, $(E_k^\star)_{0 \leq k \leq n}$ such that
	\begin{enumerate}
		\item $E_{n-k-1}^\star \subset (E_n)^*$ is the projective dual $E_k^\perp$ in $E_n$ for $0\leq k < n$, \label{itm:projdualspace}
		\item the TCDs satisfy $\tcd^\star_{k+2} = \eta (\tcd_{n-k})$ for $0 \leq k < n$,\label{itm:projdualeta}
		\item the TCDs satisfy $\tcd^\star_{k+1} = \iota (\tcd_{n-k-1})$ for $0 \leq k < n$,\label{itm:projdualiota}
		\item $T_k^\star(w) = \spa \{U_{n-k+1}(w), E_{n-k-1} \}^\perp$, for every white vertex $w$ of $\pb^\star_{k}$ and $1 \leq k \leq n$. \label{itm:projdualgeometry}
	\end{enumerate}
	Here, we introduce the convention $E_{-1} = \emptyset$ and $U$ is the subspace map associated to the flag $(T_k)$ as in Definition \ref{def:subspacemap}.
\end{definition}

\begin{lemma}\label{lem:projdualexists}
	Let $n\in \N, n \geq 2$ and let $(T_k)_{1 \leq k \leq n}$ be a flag of cogeneric TCD maps and assume $E_0$ is generic with respect to $T_1$. There is a unique projective dual flag $(T^\star_k)_{1 \leq k \leq n}$,  $(E_k^\star)_{0 \leq k \leq n}$.
\end{lemma}
\proof{ 
	We consider the combinatorics first. For $1 \leq k < n$ we define $\tcd_k^\star = \iota(\tcd_{n-k})$, which is forced due to condition (\ref{itm:projdualiota}). We set $\tcd^\star_n = \eta(\tcd_2)$ with the unique choices that are forced due to Equation \eqref{eq:sigmaeta} and condition (\ref{itm:projdualeta}) which requires that $\sigma(\tcd_n^\star) = \tcd_{n-1}^\star = \iota(\tcd_1)$. Due to Equation \eqref{eq:sigmaiota}, there are indeed choices such that $\tcd_{k-1}^\star = \sigma(\tcd_k^\star)$ for $2 \leq k < n$. Moreover, Equation \eqref{eq:sigmasigma} guarantees that there are indeed choices in the line dual such that $\eta(\tcd_k) = \tcd^\star_{n-k+2}$ for $2 \leq k < n$. Thus, the combinatorial requirements of the claim are satisfied.
	
	Let us consider the geometry. The subspaces $E^\star_k$ are determined by condition (\ref{itm:projdualspace}), which requires $E^\star_k = (E_{n-k-1})^\perp$. The TCD maps $T^\star_k$ are determined by condition (\ref{itm:projdualgeometry}.) Note that the white vertex set of $\pb_{n-k+1}$ coincides with the white vertex set of $\pb_k^\star$, because $\tcd_k^\star = \sigma(\eta(\tcd_{n-k+1}))$ and the discussion in the proof of Lemma \ref{lem:iotaetasigma}. Because we define $T^\star_k$ via condition (\ref{itm:projdualgeometry}) it is necessary that $\spa \{U_{n-k+1}(w), E_{n-k-1} \}$ is a hyperplane for all white vertices $w$ of $\pb_{n-k+1}$ and $1\leq k \leq n$. Equivalently, we need that $\spa \{U_{k}(w), E_{k-2}\}$ is a hyperplane for all vertices $w$ of $\pb_{k}$ and $1\leq k \leq n$. Also equivalently, $U_{k}(w) \cap E_{k-2} = \emptyset$ for the same $k,w$ as before. For $k=n$, $U_n(w) = T_n(w)$ and because $E_{n-2} \subset E_{n-1}$ and $E_{n-1}$ is generic and therefore does not contain points of $T_n$ follows indeed $U_{n}(w) \cap E_{n-2} = \emptyset$. For $1\leq k < n$ we write	
	\begin{align}
		U_{k}(w) \cap E_{k-2} = (U_{k}(w) \cap E_{k}) \cap E_{k-2} = (T_{k}(w)) \cap E_{k-2},
	\end{align}
	which is empty because $E_{k-2}$ is generic with respect to $T_{k-1}(w)$.
	
	We also need to verify that $T^\star_k = \sigma_{E_k^\star}(T^\star_{k+1})$ for $1 \leq k < n$. Specifically, for every white vertex of $\pb^\star_k$ we need to show the identity
	\begin{align}
		\spa \{U_{n-k+1}(w), E_{n-k-1} \}^\perp = \spa \{U_{n-k}(v_w), U_{n-k}(v'_w), E_{n-k-2} \}^\perp \cap E_k^\star,
	\end{align}
	where $v_w,v'_w$ are as in Definition \ref{def:subspacemap}. By inserting the recursive definition of the subspace maps and the definition of $E^\star_k$ the equation is equivalent to
	\begin{align}
		\spa \{U_{n-k}(v'_w),U_{n-k}(v_w), E_{n-k-1} \}^\perp = \spa \{U_{n-k}(v_w), U_{n-k}(v'_w), E_{n-k-2} \}^\perp \cap E_{n-k-1}^\perp.
	\end{align}
	This equation is trivial because $E_{n-k-2} \subset  E_{n-k-1}$. 
	
	It remains to verify that $T^\star_k$ is well-defined in the sense that $T^\star_k$ is 0-generic for $1 \leq k < n$. Consider the case $k>1$ first. Let $w_1,w_2,w_3$ be three different white vertices of $\pb^\star_k$ adjacent to a common black vertex $b$. We want to show that the points $T^\star_k(w_1), T^\star_k(w_2), T^\star_k(w_3) \in E^\star_k$ are three different points. On the other hand, we know there are three black vertices $b_1,b_2,b_3$ of $\pb_{n-k+2}$ such that
	\begin{align}
		T^\star_k(w_i) = \spa\{L_{n-k+2}(b_i), E_{k-2}\}^\perp,
	\end{align}
	for $i\in\{1,2,3\}$. The line genericity of $T_{n-k+2}$ guarantees that the three lines $L_{n-k+2}(b_i)$ are three different lines, therefore the three points $T^\star_k(w_1), T^\star_k(w_2), T^\star_k(w_3)$ are different as well. Finally, consider the case $k=1$. Then there are three white vertices $w'_1,w'_2,w'_3$ of $\pb_n$ such that
	\begin{align}
		T^\star_1(w_i) = \spa\{T_n(w'_i), E_{n-2}\}^\perp,
	\end{align}
	for $i\in\{1,2,3\}$, because in general the white vertices of $\pb$ are in bijection with the white vertices of $\sigma(\eta(\pb))$. Moreover, the three vertices $w'_i$ are adjacent to a common face and thus the three points $T_n(w'_i)$ are different because we assume that $T_n$ is 1-generic.\qed
}

Given two dual flags of TCD maps, we consider $T^\star_n$ to be the projective dual of $T_n$, as $T^\star_n$ captures the geometry of the hyperplanes of $T_n$. A sufficiently generic TCD map $T$ can be extended to a flag of TCD maps such that $T_n = T$, and therefore it is possible to find a projective dual TCD map $T^\star$. However, the combinatorics of $T^\star$ depend on the choices in the sections of the flag, and therefore $T^\star$ is not unique. 

Also note that we suspect that the cogenericity assumptions on $T_k$ in Lemma \ref{lem:projdualexists} can be dropped except for $T_1$, by adapting the last argument in the proof appropriately. 

The following lemma highlights the role of the dual flag as projective dual of a given flag.

\begin{lemma}\label{lem:subspacedual}
	Let $n\in \N$ and $n \geq 2$. Consider a flag of TCD maps $(T_k)_{1 \leq k \leq n}$, $(E_k)_{0 \leq k \leq n}$ and its dual flag $(T_k^\star)_{1 \leq k \leq n}$, $(E_k^\star)_{0 \leq k \leq n}$. Let $U_k^\star$ for $1\leq k \leq n$ be the subspace maps of the dual flag. Then $U^\star_k(w^\star) = (U_{n-k+1}(w))^\perp$ for $1\leq k \leq n$ and $w$ a white vertex of $\pb_k$, where $w^\star$ is the white vertex in $\pb_{n-k+1}^\star$ corresponding to $w$.
\end{lemma}
\proof{We recall from the proof of Lemma \ref{lem:projdualexists} that the white vertices of $\pb^\star_k$ are in bijection with the white vertices of $\pb_{n-k+1}$ for $1\leq k \leq n$. We do induction over $k$ from $n$ to $1$. For $k=n$ we see that
\begin{align}
	U^\star_n(w^\star) = T^\star_n(w^\star) = \spa\{U_1(w), E_{-1}\}^\perp = U_1(w)^\perp,
\end{align}
by employing the definition of the subspace maps, the definition of the dual flag and the convention $E_{-1}=\emptyset$. We deduce the claim for $U^\star_{k-1}$ from $U^\star_k$ via
\begin{align}
	U^\star_{k-1} &= \spa\{U^\star_{k}(v^\star_w),U^\star_{k}(v^\star_w)\}=\spa\{U_{n-k+1}(v_w)^\perp,U_{n-k+1}(v'_w)^\perp\}\\
	&= (U_{n-k+1}(v_w) \cap U_{n-k+1}(v'_w))^\perp = U_{n-k+2}(w)^\perp,
\end{align}
where $v_w$ and $v'_w$ are the vertices corresponding to $v^\star_w$ and $v'^\star_w$. Note that $v^\star_w,v'^\star_w$ of $\pb_{n-k+1}$ correspond to two black vertices $b,b'$ in $\pb_{n-k+2}$ adjacent to $w$ that do not belong to a common $\enm{}{1}$ configuration. Therefore $U_{n-k+1}(v_w), U_{n-k+1}(v'_w)$ are two different $(k-1)$-dimensional spaces that by the recursive definition contain $U_{n-k+2}(w)$.\qed
}

\begin{figure}
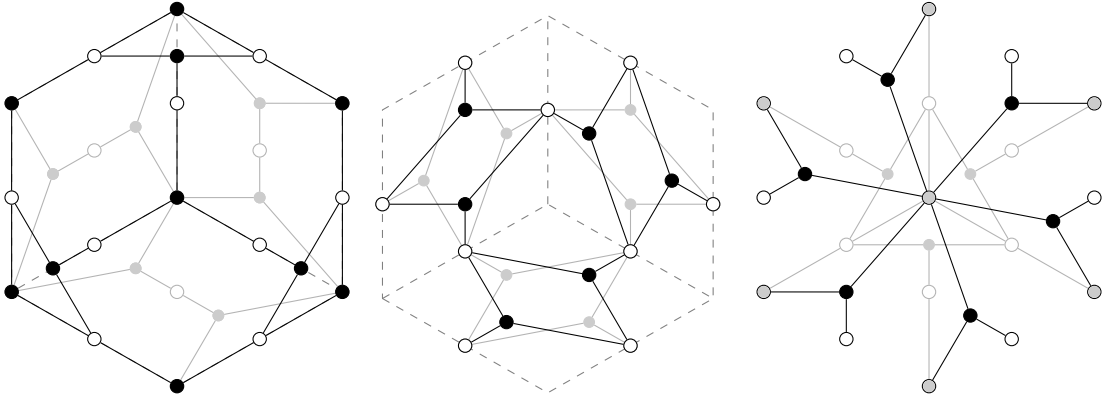


	\caption{The TCD in gray and the line dual in black for fundamental domains of Q-net, Darboux map and line complex.}
	\label{fig:tcdlinedualexamples}
\end{figure}

\begin{example}\label{ex:tcdprojdualzzz}
	In Section \ref{sec:sections} we have seen that for $\Z^3$ combinatorics, repeated sections on Q-nets yield Darboux maps, line complexes and then Q-nets again, because line complexes and line compounds coincide on $\Z^3$. Thus in the $\Z^3$ case there are flags of TCD maps consisting only of those three types of maps. Let us now think about the projective duals of those flags. A consequence of the the construction of the projective dual flag is that the projective dual $T^\star_n$ of a $\Z^3$ Q-net $T_n$ depends on the ambient dimension $n$ of the projective space $\CP^n$. In $\CP^2$ the dual is a line complex. In $\CP^3$ however the dual is a Darboux map and in $\CP^4$ it is a Q-net. More generally, let $n \in 3\Z + r$ for $r\in \Z_3$ then in the table
	\begin{align}
		\begin{tabular}{l|l|l|l}
			$r$ & Q-net &  Darboux map & line complex \\ \hline 
			0 \rule{0pt}{1em} & Darboux map & Q-net & line complex \\
			1 & Q-net & line complex & Darboux map \\
			2 & line complex & Darboux map & Q-net \\
		\end{tabular}
	\end{align}
	the entries denote the type of the projective dual $T^\star_n$ depending on the type of $T_n$ (head row).
\end{example}

Consider a projection $\proj{P}{H}$ with center $P$ to hyperplane $H$ not containing $P$. The projection of a flag of TCD maps $(T_k)_{1 \leq k \leq n}$, $(E_k)_{0 \leq k \leq n}$ is the collection $(\proj{P}{H}(T_{k}))_{1 \leq k \leq n}$, $(\proj{P}{H}(E_k))_{0 \leq k \leq n}$, which is not a flag of TCD maps in general because $\proj{P}{H}(E_n) = \proj{P}{H}(E_{n-1})$. However, if we choose $P=E_0$ then $(\proj{E_0}{H}(T_{k+1}))_{1 \leq k \leq n-1}$, $(H \cap E_{k+1})_{0 \leq k \leq n-1}$ is a flag of TCD maps. 

\begin{lemma}\label{lem:dualflagandprojections}
	Let $(T_k)_{1 \leq k \leq n}$, $(E_k)_{0 \leq k \leq n}$ be a flag of TCD maps and let $(T^\star_k)_{1 \leq k \leq n}$, $(E^\star_k)_{0 \leq k \leq n}$ be the dual flag. Let $H$ be a hyperplane not containing $E_0$. Then the projective dual flag of
	\begin{align}
		(\proj{E_0}{H}(T_{k+1}))_{1 \leq k \leq n-1},\quad (H \cap E_{k+1})_{0 \leq k \leq n-1}
	\end{align}
	is
	\begin{align}
		(T^\star_{k})_{1 \leq k \leq n-1}, \quad(E^\star_k)_{0 \leq k \leq n-1}.
	\end{align}
\end{lemma}
\proof{
	Follows immediately from Definition \ref{def:dualflag}. \qed
}

The dual flag of a projection $\proj{E_0}{H}$ of a flag can also be viewed as the section $\sigma_{E_0^\perp}$ of the dual flag.

\begin{figure}
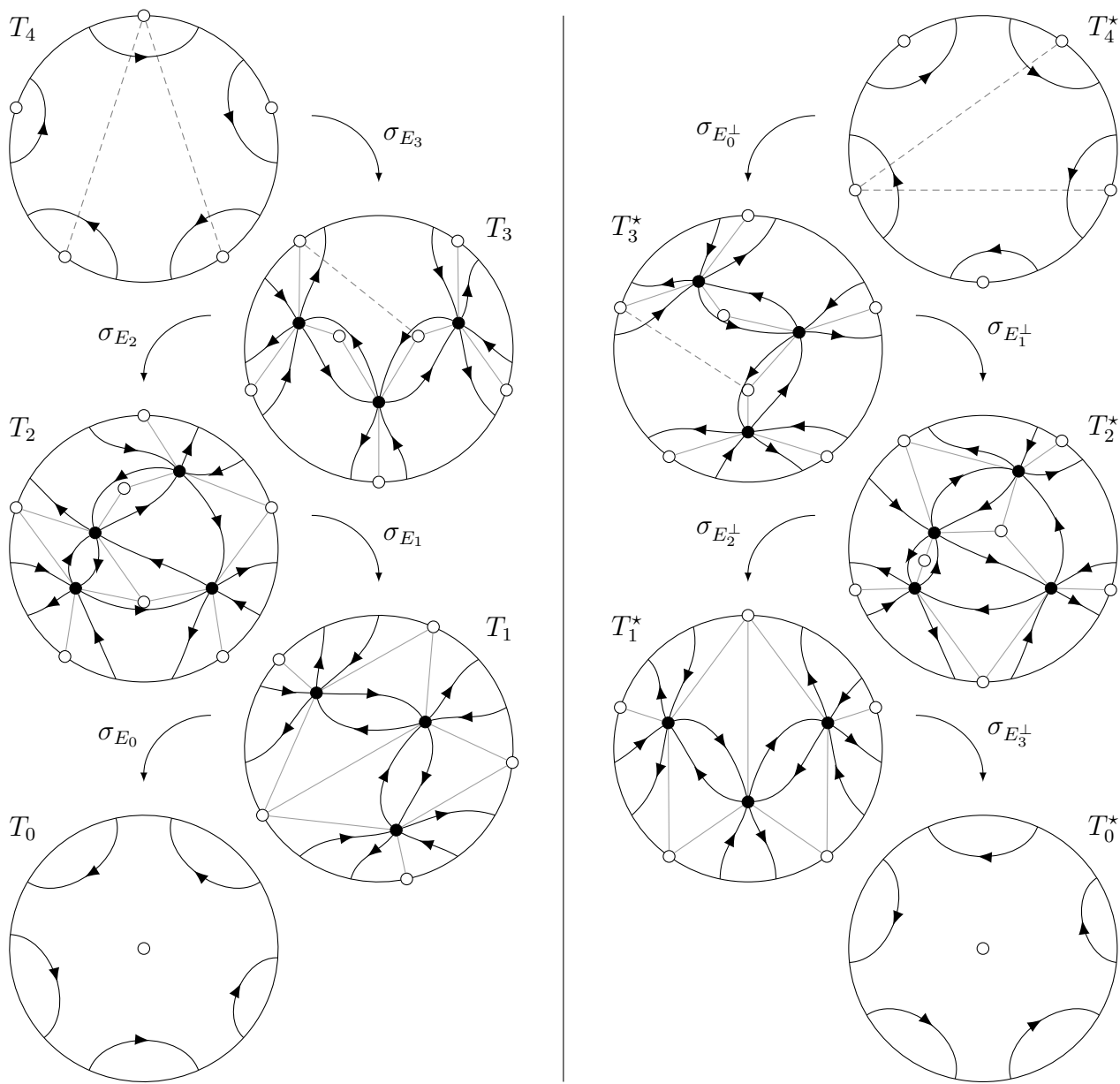


\caption{An example of a flag of TCD maps and its dual flag in $\CP^4$, where we also included $T_0$ even though it is not included in the definition of flags.}
\label{fig:balanceddualflags}
\end{figure}

\chapter{Cluster structures for TCD maps}\label{cha:clusters}

\section{Cluster algebras for planar graphs with bipartite dual}\label{sec:clusterintro}
Cluster algebras were introduced by Fomin and Zelevinsky \cite{fzclusteralgebra} in a very general algebraic and combinatorial setup. However, we will only use the theory of cluster algebras for the very particular case of planar bipartite graphs and complex numbers. Still, this restricted formulation suffices to describe local changes of combinatorics in the dimer model as well as in TCD maps. Of course, theorems that Fomin, Zelevinsky and others proved for the general cluster algebra setup still hold for the particular cases that we consider.

\begin{definition}\label{def:planquiver}
	A \emph{PDB quiver $Q$} is a planar directed graph such that its dual is bipartite. The edges are oriented such that every face is oriented either clockwise or counterclockwise.
\end{definition}
Here, PDB quiver is an abbreviation of planar dual-bipartite quiver. In the remainder of the thesis we work exclusively with PDB quivers. The edges of a quiver are called arrows.

\begin{figure}
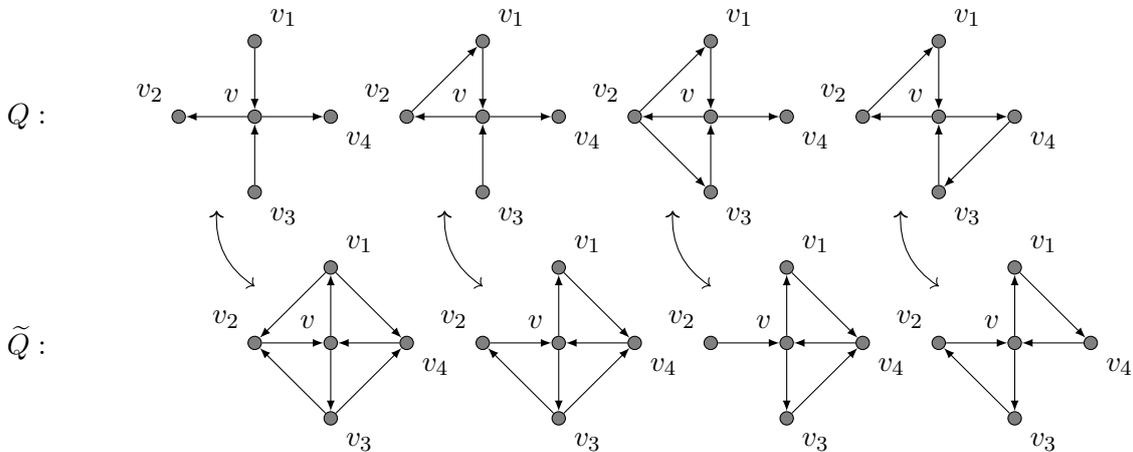


	
	\caption{Mutation at a vertex $v$ and its effect on the arrows in the neighbourhood.}
	\label{fig:quivermut}
\end{figure}

\begin{definition}
	The \emph{mutation} at a vertex $v$ of degree 4 of a PDB quiver $Q$ produces a new PDB quiver $\widetilde Q$ with the same vertex set but a different set of arrows in the neighbourhood of $v$. Let the neighbours of $v$ be $v_1,v_2,v_3,v_4$ in cyclic order and assume that the arrows incident to $v$ are $(v_1,v)$, $(v,v_2)$, $(v_3,v)$ and $(v,v_4)$. The new arrows of $\widetilde{ Q}$ are obtained via two steps from $Q$:
	\begin{enumerate}
		\item Reverse the orientation of all arrows at $v$.
		\item Add $(v_1,v_2), (v_1,v_4), (v_3,v_2), (v_3,v_4)$, possibly canceling any previously present arrows between the neighbours.\qedhere
	\end{enumerate}
\end{definition}

When we say that two arrows $(v,v')$ and $(v',v)$ cancel, we mean that if both are present in the quiver then we remove them both. It is easiest to understand the effect of a mutation on the combinatorics by drawing all the possible local configurations, as we did in Figure \ref{fig:quivermut}. In our setup we allow mutations \emph{only} at vertices of degree 4. In general cluster algebra theory it is possible to mutate at any vertex. But such mutations do not (generally) preserve the class of PDB quivers, and thus we do not allow them.

We can also capture the orientations of the arrows as an anti-symmetric form. Let us consider the free vector space $\mathcal V$ over $\R$ generated by the set of vertices of a quiver. Then let $\nu: \mathcal V \times \mathcal V\rightarrow \Z$ be the anti-symmetric bilinear form such that
\begin{align}
	\nu(v,v') = \begin{cases}
		+1 & \mbox{if there is an arrow } (v,v'),\\
		-1 & \mbox{if there is an arrow } (v',v),\\
		0 & \mbox{else}.
	\end{cases}\label{eq:nubil}
\end{align}
Then $\nu$ changes under mutation at vertex $v$ as follows:
\begin{align}
	\tilde \nu_{v_1v_2} = \begin{cases}
		-\nu_{v_1v_2} & \mbox{if } v \in \{v_1,v_2\},\\
		\nu_{v_1v_2}+ \Theta(\nu_{v_1v}) \Theta(\nu_{vv_2}) - \Theta(\nu_{v_2v}) \Theta(\nu_{vv_1}) & \mbox{else},
	\end{cases}
\end{align}
where $\Theta(x) = \max(0,x)$.

\begin{definition}\label{def:mutationx}
	The  $X$ \emph{cluster variables} for a quiver $Q$ are a set of variables $X = \{X_v\}_{v\in V}$ associated to the vertices $v\in V$ of the quiver and taking values in $\C$. When \emph{mutating} at a vertex $v$ the new variables $\tilde X$ are determined as follows:
	\begin{equation}
		\tilde X_{v'} = \begin{cases}
			X_v^{-1} & \mbox{if }v' = v,\\
			X_{v'} (1+X_v) & \mbox{if }(v,v') \mbox{ is an arrow},\\
			X_{v'} (1+X_v^{-1})^{-1} & \mbox{if }(v',v) \mbox{ is an arrow},\\
			X_{v'} & \mbox{else.}
		\end{cases}	\label{eq:mutationx}		
	\end{equation}\qedabovehere
\end{definition}
Note that this definition of the cluster variables and mutation also holds for non PDB quivers. The $X$-cluster variables are also called \emph{of coefficient type} in the literature. In combination with the quiver they are also called \emph{Y-patterns} \cite{fzcoefficients}.

A pair $(Q,X)$ of quiver and cluster variables is called a \emph{seed} in cluster algebra theory. By iterating mutations one can generate all seeds of the cluster algebra, which may be an infinite number of seeds. The whole \emph{cluster algebra} as an algebra is generated by the union of all variables of all seeds modulo the relations induced by mutation. Thus the cluster algebra can be generated from a seed. The term \emph{cluster structure} is colloquially used but not clearly defined. In the following we call a pair $(Q,X)$ a cluster structure.

In fact, different sets of variables with different mutation laws available in the literature. There is no consensus on how to call or denote the different sets of variables. The $X$-variables we have introduced here are most directly related to the geometric notions we will study in this thesis. We denote them by $X$ when they are invariant under projective transformations applied to the TCD maps that we derive them from. We denote them by $Y$ when they are invariant under affine transformations applied to the TCD maps that we derive them from. We will explain the details in the next two sections. There is one more set of cluster variables that will be useful for us, we call them the $\tau$ cluster variables. Unlike the $X$ variables, the $\tau$ variables will not be uniquely defined by the TCD maps. We adopt the letter $\tau$ because they frequently coincide with potentials introduced in the discrete integrable system community that are denoted by $\tau$.

\begin{definition}\label{def:mutationtau}
	The \emph{$\tau$ cluster variables} for a quiver $\qui$ are a set of variables $\tau = \{\tau_v\}_{v\in V}$ associated to the vertices $v\in V(\qui)$ of the quiver ands taking values in $\C$. When \emph{mutating} at a vertex $v$ of degree four with neighbours $v_k$ as in Figure \ref{fig:quivermut} the new variables $\tilde \tau$ are determined by
	\begin{equation}
		\tilde \tau_{v'} = 
		\begin{cases}
			\dfrac{\tau_{v_1}\tau_{v_3}+\tau_{v_2}\tau_{v_4}}{\tau_v} & \mbox{if }v' = v,\\
			\tau_{v'} & \mbox{else}.\hspace{5.6cm} 
		\end{cases}\label{eq:taumutation}
	\end{equation}\qedabovehere	
\end{definition}

We observe that the mutation law \eqref{eq:taumutation} is a local instance of the dKP equation, also known as the octahedron recurrence \cite{speyerdimers} or discrete Hirota equation \cite{hirotaequation}. The mutation of $\tau$-variables can also be defined for vertices with degree not equal to four, in which case $\tilde \tau_v \tau_v$ is the product of all variables at incoming arrows plus the product of all variables at outgoing arrows at $v$. The $X$ cluster variables can be calculated from the $\tau$ variables by the formula
\begin{align}
	X_v = \frac{\prod_{(v',v)} \tau_{v'}}{\prod_{ (v,v')} \tau_{v'}},\label{eq:tauandx}
\end{align}
where the products are over all incoming arrows in the numerator and outgoing arrows in the denominator. It is a straightforward exercise to check that this formula is compatible with the two mutation laws.

\begin{definition}\label{def:reciprocalcluster}
	Let $\qui$ be a PDB quiver. Construct a new quiver $\rho(\qui)$ from $\qui$ by reversing all arrows of $\qui$. We call $\rho(\qui)$ the \emph{reciprocal quiver} of $\qui$. Moreover, if $(\qui,X)$ is a cluster structure then the \emph{reciprocal cluster structure} is $(\rho(\qui), \rho(X))$, where $\rho(X)_v = X_v^{-1}$ for all vertices $v \in V(\qui)$.
\end{definition}

Clearly, $\rho \circ \rho$ is the identity. Moreover, $\rho$ is compatible with mutations in the sense of the next Lemma.

\begin{lemma}
	Let $Q$ be a quiver and $X$ its cluster variables. Let $\tilde Q$ and $\tilde X$ be the quiver and variables after a mutation at a vertex $v$. Then we have that
	\begin{align}
		\widetilde{\rho(Q)} &= \rho(\tilde Q),\\
		\widetilde{\rho(X)} &= \rho(\tilde X).\qedhere
	\end{align}
\end{lemma}
\proof{
	Both relations can be proven by straightforward calculations. \qed
}

We introduce the notion of a reciprocal cluster structure because we will encounter the reciprocal cluster structure when relating cluster structures of the projective dual to the primal in Section \ref{sec:projclusterduality}. However, we will need another operation at the boundary.

\begin{definition}\label{def:upsilon}
	Let $\qui$ be a PDB quiver and let $v_1,v_2,\dots,v_n$ be the boundary vertices of the $\qui$ in counterclockwise order. We construct the quiver $\Upsilon^+(\qui)$ from $\qui$ by adding the arrows $(v_1,v_2), (v_2,v_3),\dots, (v_n,v_1)$. Analogously, we construct the quiver $\Upsilon^-(\qui)$ from $\qui$ by adding the arrows $(v_2,v_1), (v_3,v_2),\dots, (v_1,v_n)$.
\end{definition}

Clearly $\Upsilon^\pm \circ \Upsilon^\mp$ is the identity, and $\rho \circ \Upsilon^\pm = \Upsilon^\mp \circ \rho$.

\section{Particular classes of quivers}
In this section we study the basic quivers that appear when we look at objects of discrete differential geometry. Coincidentally, these quivers occur in models of statistical physics as well. This coincidence of quivers will allow us later to establish an algebraic and combinatorial link between the two subjects via cluster algebras.

\subsection{Square grid quiver}\label{sub:squarequiver}
\begin{figure}
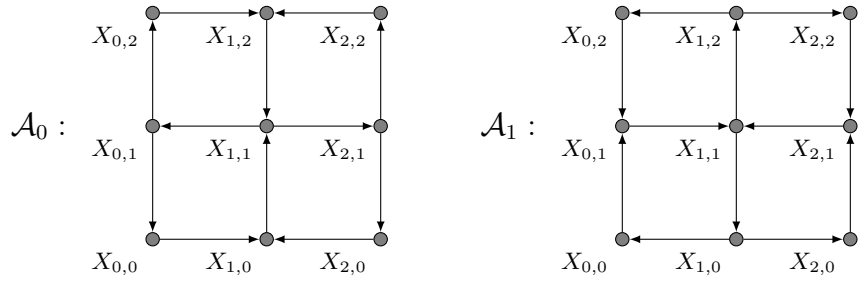
\scriptsize

	\caption{The square quivers $\mathcal A_0$ and $\mathcal A_1$ that alternate under mutation at half the vertices.}
	\label{fig:squarequivers}
\end{figure}

This quiver is simply the $\Z^2$ lattice with a consistent orientation of the arrows. The set of arrows is one of the following two (see Figure \ref{fig:squarequivers}):
\begin{align}
	\mathcal{A}_0 &= \bigcup_{i+j\in 2\Z\phantom{+0}} \{(v_{i,j},v_{i+1,j}),(v_{i,j},v_{i-1,j})),(v_{i,j+1},v_{i,j}),(v_{i,j-1},v_{i,j}) \}, \\
	\mathcal{A}_1 &= \bigcup_{i+j\in 2\Z + 1} \{(v_{i,j},v_{i+1,j}),(v_{i,j},v_{i-1,j})),(v_{i,j+1},v_{i,j}),(v_{i,j-1},v_{i,j}) \}. 	
\end{align}
One step of global dynamics consists of mutating at every other vertex of the quiver. To fix the positive direction of discrete time, we mutate at the vertices $v_{i,j}$ with $i+j\in (2\Z+k)$ whenever we are in the quiver $\mathcal A_k$. The mutation of $\mathcal A_k$ yields $\mathcal A_{1-k}$. Indeed, all arrow orientations are reversed and every diagonal that was added due to mutation at one vertex is canceled by the mutation at another vertex that is two steps away. Moreover if we mutate the quiver with arrows $\mathcal A_k$ then the cluster variables $X$ transform as:
\begin{align}
	\tilde X_{i,j} = \begin{cases}
		X_{i,j}^{-1} & \mbox{if }i+j\in 2\Z + k,\\
		\dfrac{(1+X_{i,j+1})(1+X_{i,j-1})}{(1+X_{i+1,j})(1+X_{i-1,j})}X_{i,j} & \mbox{else}.
	\end{cases}
\end{align}
This quiver will occur in the study of Laplace-Darboux dynamics and related objects like $\Z^2$ circle patterns.

\subsection{Quiver of a triangulation}

\begin{definition}
	Let $\tg$ be a triangulation. The quiver $\qui_\tg$ of $\tg$ is obtained by gluing the following piece
	\begin{center}
		%
	\end{center}
	of a quiver into every triangle, such that there is one quiver vertex per edge of $\tg$.
\end{definition}

\begin{figure}
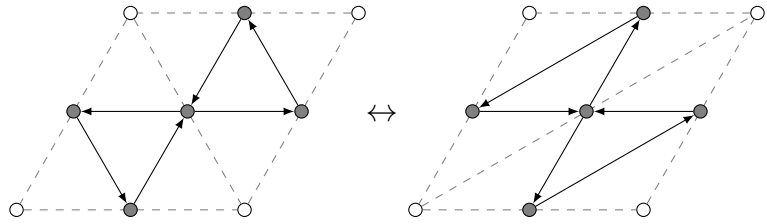

	%
	
	\caption{The edge flip in a triangulation (white vertices and dashed lines) and the corresponding mutation in a triangular quiver (arrows and gray vertices).}
	\label{fig:quiveredgeflip}
\end{figure}

\begin{lemma}
	Let $\tg$ be a triangulation and $\qui_\tg$ the corresponding quiver. The edge-flip at $e$ corresponds to mutation at the quiver vertex corresponding to $e$.
\end{lemma}
\proof{See Figure \ref{fig:quiveredgeflip}.\qed}

\subsection{Cuboctahedral quiver}\label{sus:cuboctahedral}

Consider a quad with vertices $v,v_i,v_{ij},v_j$ in counterclockwise order. Using shift notation (see Section \ref{sec:znotation}), we denote the edges in that quad by
\begin{align}
	e^i = \{v,v_i\}, \quad e^j = \{v,v_j\} , \quad e^i_j = \{v_j,v_{ij}\} , \quad e^j_i = \{v_i,v_{ij}\}.
\end{align}

\begin{definition}\label{def:cuboquiver}
	Let $\qg$ be a quad-graph. The \emph{cuboctahedral quiver} $\qui_\qg$ of $\qg$ has a vertex for every edge of $\qg$. For each quad $(v,v_1,v_{12},v_2)$ of $\qg$ we add the four arrows
	\begin{align}
		e^1 \rightarrow e^2 \rightarrow e^1_2 \rightarrow e^2_1 \rightarrow e^1
	\end{align}
	to the quiver $\qui_\qg$.
\end{definition}
If we reverse all arrows of the cuboctahedral quiver, we call the result the \emph{reversed cuboctahedral quiver}. We will also write \emph{the} cuboctahedral quiver when we mean the cuboctahedral quiver of a stepped surface.

\begin{figure}
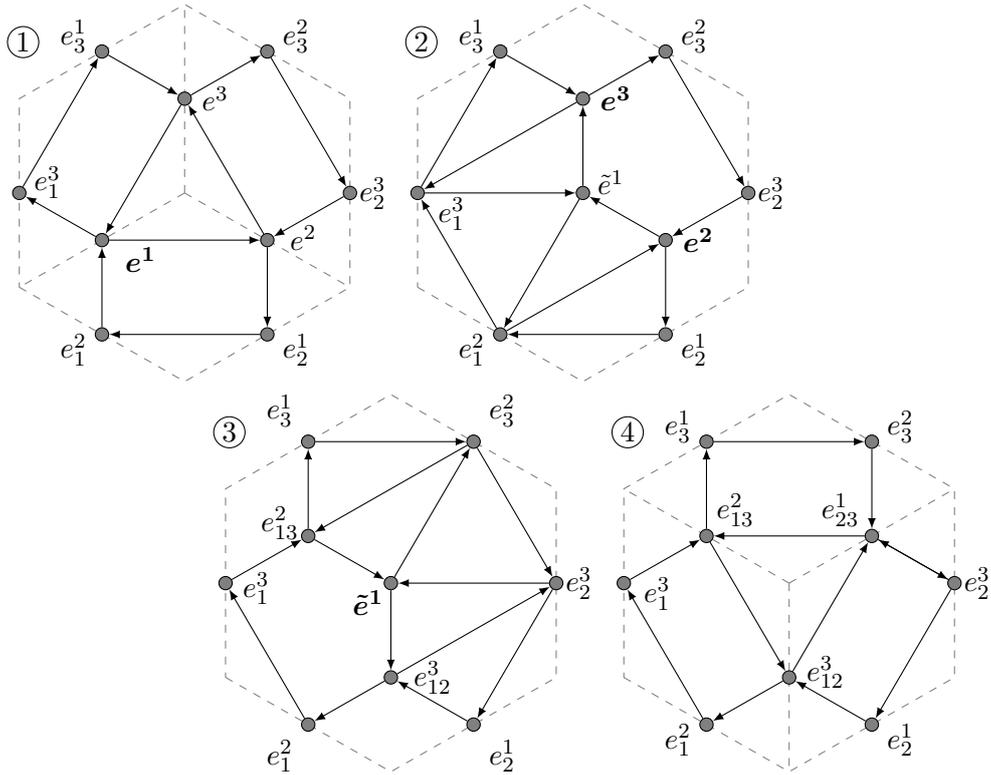
 
	\hspace{-2.7cm}\small

	\caption{The cuboctahedral quiver and the four mutations of the cube flip. First we mutate at $\bm{e^1}$, then at $\bm{e^2}$ and $\bm{e^3}$ simultaneously and finally at $\bm{\tilde e^1}$.}
	\label{fig:cuboctaquiver}
\end{figure}

\begin{definition}\label{def:cuboflip}
	Consider three incident quads in a quad-graph $\qg$ and view them as the backside of a cube, see Figure \ref{fig:cuboctaquiver}. The \emph{cuboctahedral flip} is the sequence of mutations at quiver vertices $e^1, e^3,e^2$ and again at $\tilde e^1$ in that order, where $\tilde e^1$ is $e^1$ after the first mutation. 
\end{definition}
We observe that the cuboctahedral flip in $\qui_\qg$ corresponds to the cube flip of $\qg$. Note that there are several possible sequences of mutations that correspond to the cube flip. For example in the given sequence the mutations at $e^3$ and $e^2$ can be interchanged without changing the final combinatorics.

The cuboctahedral quiver occurs in $\Z^3$ Q-nets, Darboux maps and line complexes (see Figures \ref{fig:projquivers}, \ref{fig:affquivers}) as well as the spanning tree model (see Section \ref{sub:spanningtrees}).

\subsection{Hexahedral quiver} \label{sus:hexahedral}

\begin{definition}\label{def:hexaquiver}
	Let $\qg$ be a quad-graph. The \emph{hexahedral quiver} $\qui_\qg$ of $\qg$ has a vertex for every vertex of $\qg$ and a vertex for every quad of $\qg$. For each quad $q=(v,v_1,v_{12},v_2)$ with $v$ white and in counterclockwise order we add the six arrows
	\begin{align}
		v^{12} \rightarrow v \rightarrow v_1 \rightarrow v^{12} \rightarrow v_{12} \rightarrow v_{2} \rightarrow v^{12}
	\end{align}
	to the quiver, where $v^{12}$ denotes the vertex corresponding to that quad.
\end{definition}
If we reverse all arrows of the hexahedral quiver we call the result the \emph{reversed hexahedral quiver}. The reversed hexahedral quiver is also obtained by interchanging the vertex colors. We will also write \emph{the} hexahedral quiver when we mean the hexahedral quiver of a stepped surface.

\begin{figure}
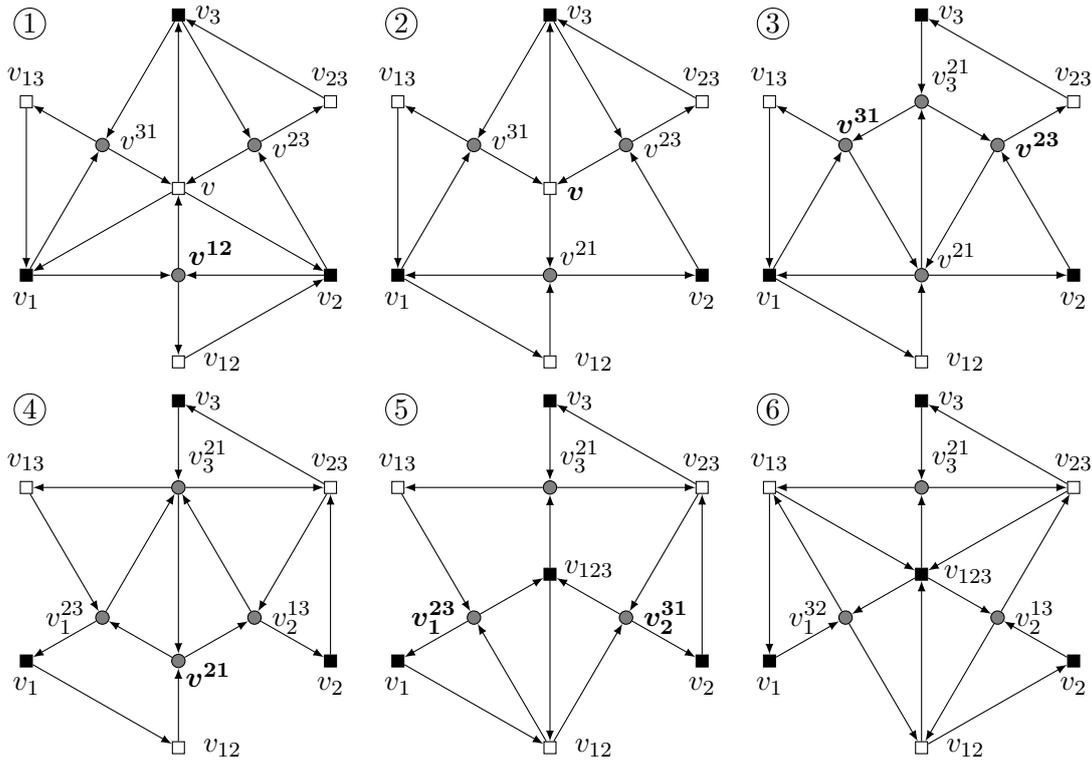

	\small

	\caption{The hexahedral quiver and the seven mutations of a cube flip. The sequence of mutations is $\bm{v^{12},v,v^{31},v^{23},v^{21},v_1^{23},v_2^{31}}$.}
	\label{fig:hexahedralquiver}
\end{figure}

\begin{definition}\label{def:hexaflip}
	Consider three pairwise incident quads in a quad-graph $\qg$ and view them as the backside of a cube, see Figure \ref{fig:hexahedralquiver}. The \emph{hexahedral flip} is the sequence of mutations that begins at quiver vertex $v^{12}$ followed by $v$, $v^{31}$, $v^{23}$, $v^{21}$ and finally at $v_1^{23}$ and $v_2^{31}$.
\end{definition} 
We observe that the hexahedral flip in $\tcd_\qg$ corresponds to the cube flip in $\qg$. Note that there are several possible sequences of mutations that correspond to the cube flip. For example in the given sequence the mutations at $v^{31}$ and $v^{23}$ can be interchanged without changing the final combinatorics.

The hexahedral cluster structure occurs both in the case of Q-nets and Darboux maps (see Figures \ref{fig:projquivers}, \ref{fig:affquivers}) as well as the Ising model (see Section \ref{sub:ising}).

\section{The projective cluster structure of a TCD map}\label{sec:projcluster}
The first cluster structure that we associate to a TCD map has a quiver which is the dual of the associated bipartite graph $\pb$. The arrows are oriented such that they turn counterclockwise around the black vertices of $\pb$. It is possible that $\pb$ has white vertices of degree 2. In this case the two arrows in the dual cancel each other. From the perspective of the TCD, the vertices of the quiver correspond to all the clockwise oriented faces of the TCD, and arrows connect faces that are incident to a common triple intersection point. We formalize this explanation in the following definition.

\begin{definition}\label{def:projclusterstructure}
	Let $T$ be a TCD map and let $\mu$ be the edge weights of an associated VRC. The vertices of the \emph{extended projective cluster structure} $\pro^\times(T)$ are located at each face $f$ of the associated graph $\pb$. The number of arrows $\nu_{ff'}$ between two faces $f,f'$ is defined such that
	\begin{align}
		\nu_{ff'} = \phantom{-} &|\{(f,w,f',b) \mbox{ counterclockwise} : b\in B, w \in W \}| \\ \nonumber
		- &|\{(f,w,f',b) \mbox{ clockwise} : b\in B, w \in W \}|.
	\end{align}
	Let $f$ be an face of $\pb$ and let $b_1,w_1,b_2,...,w_n$ be the black and white vertices bounding $f$ in counterclockwise order.	Then the projective cluster variable $X_f$ is the alternating ratio
	\begin{equation}
		X_f = (-1)^{n+1}\prod_{i=1}^n \frac{ \mu(b_i,w_i)}{\mu(w_i,b_{i+1})}.
	\end{equation}
	The \emph{projective cluster structure}  $\pro(T)$ is the restriction of the extended projective cluster structure to interior faces of $\pb$.
\end{definition}
Because we define the projective cluster variables as alternating ratios, the variables $X_f$ of interior faces $f$ do not change when we rescale the edge weights around a vertex of $\pb$. In other words, the projective cluster variables $X$ of the projective cluster structure do not depend on the choice of homogeneous lifts, that is on the choice of associated VRC $R$. The $X$-variables at the boundary of the extended projective cluster structure however do depend on the lifts and thus the choice of VRC. We will generally not deal with the extended projective cluster structure except in Section \ref{sec:sectioncluster}.

Of course, if we define a cluster structure then the next question is if there is some change of combinatorics that corresponds to a mutation of the quiver and the variables. Recall that for TCDs the local changes of combinatorics are the 2-2 moves, of which there are two types, see Section \ref{sec:tcddskp}.

\begin{figure}
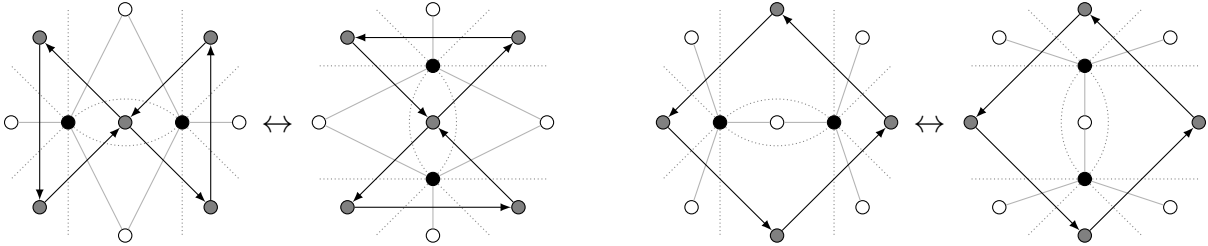

	\centering
	
	\caption{The projective quiver (gray vertices and arrows) before and after a spider move (left) as well as before and after a resplit (right).}\label{fig:twotwoquiver}
\end{figure}

\begin{lemma}
	The spider move in $\pb$ corresponds to a mutation in the projective cluster algebra, the resplit in $\pb$ leaves the projective cluster algebra invariant (cf. Figure \ref{fig:twotwoquiver}).
\end{lemma}
\proof{Looking at the edge weights in the spider move (Figure \ref{fig:vrclocalmoves}), we see that for the cluster variables $X, \tilde X$ before resp. after the spider move at the face where we mutate
	\begin{align}
		X = \frac{bd}{ac} = \tilde{X}^{-1}
	\end{align} 
	holds, as it should be. For the cluster variables of adjacent faces, we do the calculation for $X_{12}$ which is the variable of the face incident to $w_1$ and $w_2$. We obtain
	\begin{align}
		\frac{\widetilde{X}_{12}}{X_{12}} = \frac{ab^{-1}-dc^{-1}}{ab^{-1}} = 1-\frac{bd}{ac} = 1 + X,
	\end{align}
	also as it should be. The calculations for $\widetilde{X}_{23},\widetilde{X}_{34},\widetilde{X}_{14}$ are similar. Therefore the spider move acts as mutation on the cluster variables $X$. In case of the resplit (Figure \ref{fig:resplit}) one can see without calculation that the $X$ variables are unchanged for all four adjacent faces.\qed
}

As we explained in the beginning, the edge-weights of a VRC associated to a TCD map are determined up to gauge by the points of the TCD map. Thus it makes sense to also consider an expression of the $X$-variables via the points of a TCD map.

\begin{lemma}\label{lem:projclusterviadistances}
	Let $T$ be a TCD map. The alternating ratio $X_f$ at a face $f = (b_1,w_1,b_2,...,w_n)$ of $\pb$ can be expressed as the multi-ratio
	\begin{align}
		X_f &= (-1)^{n+1} \mr(T(w_1),T(v_2),T(w_2),T(v_3),\dots,T(w_n),T(v_1))  \\
			&= (-1)^{n+1}\prod_{i=1}^n \frac{T(w_{i-1})-T(v_i)}{T(v_i)-T(w_{i})},
	\end{align}	
	where $v_i,w_i,w_{i-1}$ are the three neighbours of $b_i$. The oriented length ratios are defined as in Definition \ref{def:orientedlengthratio}.
\end{lemma}
\proof{Consequence of Lemma \ref{lem:edgeratios}. Note that in Lemma \ref{lem:edgeratios} we assumed that the edge weights are given in affine gauge. However, changing gauge does not change the value of $X_f$, because the scaling factors cancel in Definition \ref{def:projclusterstructure}.\qed}

\begin{figure}
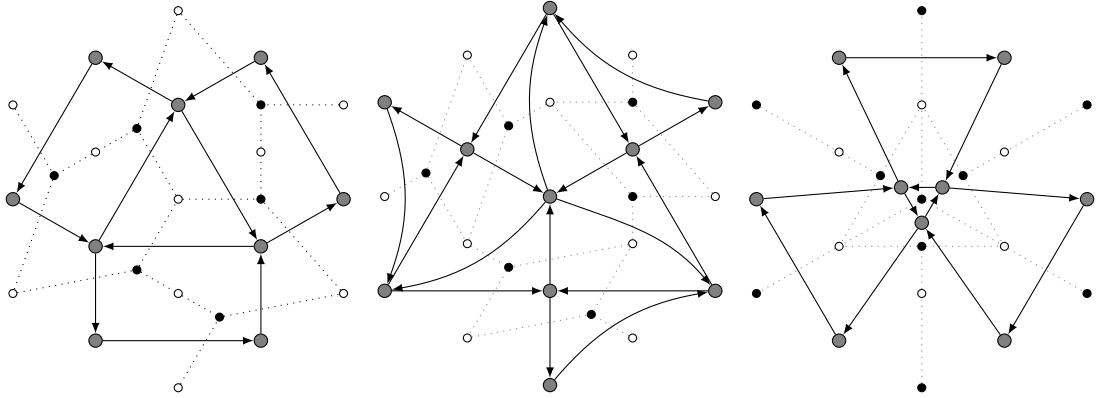


	\caption{The projective quivers (arrows and gray vertices) for a Q-net, Darboux map and line compound. The graph $\pb$ is drawn dotted. We recognize the cuboctahedral, hexahedral and again the cuboctahedral quivers.}
	\label{fig:projquivers}
\end{figure}

Recall that by Lemma \ref{lem:projinvariants} multi-ratios of polygons with additional points on the lines are invariant under projective transformations. The name projective cluster structure is therefore justified because the projective cluster variables are invariant under projective transformations. Another important point is that the combinatorics of the projective quiver determine the combinatorics of $\pb$ up to resplits. One notices for example in Figure \ref{fig:projquivers} or from the definition that all faces of the quiver that are oriented counterclockwise (white faces) contain a single white vertex. All faces that are oriented clockwise (black faces) contain $(d_f-2)$ black vertices and $(d_f-3)$ white vertices, where $d_f$ is the degree of the face. The black faces encode the information that the points in the adjacent white faces are in a $(d_f-2)$-dimensional projective subspace.

Let us consider our standard DDG examples in Figure \ref{fig:projquivers}. We observe that both $\Z^3$ Q-nets and line compounds feature a cuboctahedral quiver, while Darboux maps feature a hexahedral quiver. However, note that the orientations of the two cuboctahedral quivers are not the same, corresponding to the fact that the black vertices are in the quads of the quiver for Q-nets while the opposite is the case for line compounds. Of course, because line complexes and line compounds coincide when defined on $\Z^3$, the projective quiver of a stepped surface of a line complex is also cuboctahedral. The same is however not true for more general quad-graphs.

Recall that in Definition \ref{def:qgqnet} of a Q-net $q:V(\qg)\rightarrow \CP^n$ on a quad-graph $\qg$ we also defined two focal points per quad. On the other hand, to each edge $(v,v') = e\in E(\qg)$ we can also assign the two focal points of the adjacent quads that are on the line $q(v)q(v')$.

\begin{definition}\label{def:laplaceinvariant}
	Let $q:V(\qg)\rightarrow \CP^n$ be a flip-generic Q-net. The \emph{Laplace invariant} is a map $\Lambda: E(\qg) \rightarrow \C\setminus\{0\}$ such that
	\begin{align}
		\Lambda_e = -\cro(q(v),\fp(f),q(v'),\fp(f')),
	\end{align}
	for every edge $e=(v,v')\in E(\qg)$, where $e^*=(f,f')$ and $v,f,v',f'$ appear in counterclockwise order in $\qg$.
\end{definition}

Laplace invariants were introduced by Doliwa \cite{doliwalaplace}.

\begin{lemma}\label{lem:laplaceinvprojvar}
	Let $q:V(\qg)\rightarrow \CP^n$ be a flip-generic Q-net. For each non-boundary edge $e\in E(\qg)$ there is a unique face $f_e\in F(\pb)$ and vice versa. The Laplace invariants of $q$ and the projective cluster variables coincide, that is
	\begin{align}
		\Lambda_e = X_{f_e}
	\end{align}
	for all $e\in E(\qg)$.
\end{lemma}
\proof{
	Fix an edge $e\in E(\qg)$. We can always perform resplits at the adjacent focal points such that the face $f_e$ has degree four because $q$ is flip-generic. Recall that resplits do not change the projective invariants $X_{f_e}$. Then the claim at $e$ is equivalent to Lemma \ref{lem:projclusterviadistances}.\qed
}

Especially in Chapter \ref{cha:cpone}, on TCD maps related to statistical physics, the case of TCD maps with strictly positive $X$-variables is of particular interest. However, we can also state something interesting straight away.

\begin{lemma}
	A TCD map $T: \tcdp \rightarrow \CP^n$ defined on a minimal TCD $\tcd$ for which all $X$-variables are strictly positive real, is a flip-generic TCD map in the sense of Definition \ref{def:gentcdmap}.
\end{lemma}
\proof{Note that if the $X$-variables of $T$ are strictly positive then they are strictly positive for any TCD map $\tilde T$ related to $T$ via a sequence of 2-2 moves, because of the mutation rules explained in Definition \ref{def:mutationx}. Assume that there is a face $f$ in some 2-2 related TCD map $\tilde T$ such that a spider move cannot be performed at $f$. Then $X_f = -1$ in contradiction to the fact that all $X$-variables are strictly positive.\qed
}

\section{Ideal hyperbolic triangulations}\label{sec:extriangulation}

One of the well known occasions where cluster algebras relate to geometry is the Teichmüller theory of punctured surfaces, see the introduction to cluster algebras by Williams \cite{williamsintro}. In this section we want to briefly outline how the $X$-cluster variables of Teichmüller theory, that is the so called \emph{shear coordinates} appear as a special case of the projective invariants of TCD maps that we introduced in Section \ref{sec:projcluster}. For an introduction to Teichmüller theory, ideal hyperbolic triangulations and shear coordinates we recommend the book of Penner \cite{pennerteichmuller}. In order to avoid having to introduce the whole framework of Teichmüller theory, we compare these variables on the level of hyperbolic structures.

Consider the \emph{Poincaré disk model} of the hyperbolic plane. Here, the disk $D$ is the set of complex numbers $\{z\in \C \ | \ |z| < 1\} \subset \C \subset \CP^1$, where we view $\C$ as an affine chart of $\CP^1$. The disk is a model for the hyperbolic plane and the orientation preserving isometries of the hyperbolic plane correspond to the projective transformations of $\CP^1$ that preserve the boundary of the disk. The points on the boundary are called \emph{ideal points}.

\begin{definition}
	Let $\tg$ be a triangulation. An \emph{ideal hyperbolic triangulation} is a map $z: V(\tg) \rightarrow \partial D$ from the vertices of $\tg$ to the boundary of the Poincaré disk, such that the edges are realized as geodesics and such that no two such edge-geodesics are intersecting.
\end{definition}

Fix the triangulation $\tg$ for a moment. Then two natural questions are: when are two ideal hyperbolic triangulations $z,z'$ related by a projective transformation? And how many different ideal hyperbolic triangulations are there up to projective transformations? One way to answer these questions is via shear coordinates.

\begin{definition}
	Let $z: V(\tg) \rightarrow \partial D$ be an ideal hyperbolic triangulation. The \emph{shear coordinates} $s: E(\tg) \rightarrow \R_{>0}$ are defined such that
	\begin{align}
		s_e = -\cro(z_{v_1},z_{v_2},z_{v_3},z_{v_4}),
	\end{align}
	for every edge $e=(v_1,v_3)$ and where $v_2,v_4$ are such that $v_1,v_2,v_3$ and $v_3,v_4,v_1$ are triangles in $\tg$ and $v_1,v_2,v_3,v_4$ appear in counterclockwise order in $\tg$.
\end{definition}

Indeed, two ideal hyperbolic triangulations are related by an orientation preserving hyperbolic isometry if and only if they posses the same shear coordinates \cite{pennerteichmuller}. Moreover, for every set of shear coordinates there exists an ideal hyperbolic triangulation, which answers the two questions.

We can associate a TCD map $T_z$ to every ideal hyperbolic triangulation $z: V(\tg) \rightarrow \partial D$ by gluing
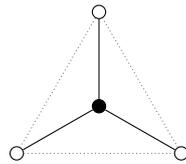
\begin{center}\vspace{-5mm}
	\begin{tikzpicture}[scale=1.25]
		\node[wvert] (v1) at (90:1) {};
		\node[wvert] (v2) at (210:1) {};
		\node[wvert] (v3) at (330:1) {};			
		\node[bvert] (b) at (0,0) {};
		\draw[gray,-,densely dotted]
			(v1) -- (v2) -- (v3) -- (v1)
		;
		
		\draw[-]
			(b) edge (v1) edge (v2) edge (v3)
		;
	\end{tikzpicture}
\end{center}
into every triangle of $\tg$. Each white vertex $w$ is associated to a vertex $v$ of $\tg$ and the TCD map $T_z$ maps $w$ to $z(v)$. We observe that the TCD necessarily has endpoint matching $\enm n1$ (see Definition \ref{def:endpointmatching}), where $n$ is the number of vertices of $\tg$. Moreover, the combinatorics of the VRC clearly constrain all vertices to be on a line, corresponding to the fact that we are considering triangulations in $\CP^1$.

\begin{lemma}
	The $X$-variables of $T_z$ coincide with the shear coordinates of $z$.
\end{lemma}
\proof{To begin with, we notice that the combinatorics are correct. That is there is a face in $\pb$ of $T_z$ for every edge of $\tg$. Moreover each face of the VRC has degree four and then the claim follows straight from Lemma \ref{lem:projclusterviadistances}.\qed
}

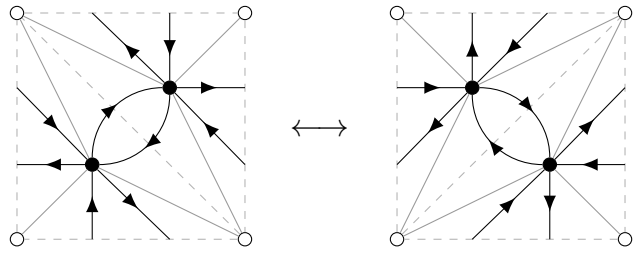
\begin{figure}
	\begin{tikzpicture}[baseline={([yshift=-.7ex]current bounding box.center)},scale=3]
		\node[wvert] (v) at (0,0) {};
		\node[wvert] (v1) at (1,0) {};
		\node[wvert] (v2) at (0,1) {};
		\node[wvert] (v12) at (1,1) {};
		\node[bvert] (b) at (.33,.33) {};
		\node[bvert] (b12) at (.67,.67) {};
		\draw[gray!60,dashed]
			(v1) -- (v) -- (v2) -- (v12) -- (v1) -- (v2)
		;
		\draw[gray!80,-]
			(b) edge (v) edge (v1) edge (v2)
			(b12) edge (v12) edge (v1) edge (v2)
		;
		\draw[-]
			(b) edge[mid rarrow] (0.33,0) edge[mid arrow] (0,0.33) edge[mid arrow] (0.67,0) edge [mid rarrow] (0,0.67) edge [bend left=45, mid arrow] (b12)
			(b12) edge[mid rarrow] (0.67,1) edge[mid arrow] (1,0.67) edge[mid arrow] (0.33,1) edge [mid rarrow] (1,0.33) edge [bend left=45, mid arrow] (b)
		;
	\end{tikzpicture}\hspace{.5cm}$\longleftrightarrow$\hspace{.4cm}
	\begin{tikzpicture}[scale=3,baseline={([yshift=-.7ex]current bounding box.center)}]
		\node[wvert] (v) at (0,0) {};
		\node[wvert] (v1) at (1,0) {};
		\node[wvert] (v2) at (0,1) {};
		\node[wvert] (v12) at (1,1) {};
		\node[bvert] (b1) at (.67,.33) {};
		\node[bvert] (b2) at (.33,.67) {};
		\draw[gray!60,dashed]
			(v) -- (v2) -- (v12) -- (v1) -- (v) -- (v12)
		;
		\draw[gray!80,-]
			(b1) edge (v) edge (v1) edge (v12)
			(b2) edge (v) edge (v2) edge (v12)
		;
		\draw[-]
			(b1) edge[mid rarrow] (1,0.33) edge[mid arrow] (0.67,0) edge[mid arrow] (1,0.67) edge [mid rarrow] (0.33,0) edge [bend left=45, mid arrow] (b2)
			(b2) edge[mid rarrow] (0,0.67) edge[mid arrow] (0.33,1) edge[mid arrow] (0,0.33) edge [mid rarrow] (0.66,1) edge [bend left=45, mid arrow] (b1)
		;
	\end{tikzpicture}
	\caption{The edge flip in a triangulation (dashed), which corresponds to a spider move in $\pb$ (black).}
	\label{fig:triaedgeflip}
\end{figure}

In the theory of ideal hyperbolic triangulations one can perform edge-flips and calculate the new shear coordinates from the old ones. Because of the previous Lemma it is no surprise that the edge-flip formulas for the shear coordinates coincide with cluster mutations. Figure \ref{fig:triaedgeflip} shows that an edge-flip in $z$ corresponds to a spider move in $T_z$.

We will study generalizations to projective $k$-flag configurations as introduced by Fock-Goncharov in Section \ref{sec:fgmoduli}. The corresponding $\tau$ cluster variables coincide with Penner coordinates in the special case of ideal hyperbolic triangulations.

Another generalization is not to restrict the image $T_z$ to $D$. In that case one can study triangulations living in $\CP^1$, although it is not clear in that case for which $X$-variables the triangle edges are non-intersecting. One can generalize even further and also consider maps $T_z$ that are defined on non-minimal TCDs, that is triangulations with interior vertices. In that case the $X$-variables are not free anymore but are constrained to certain subvarieties that we do not discuss.

\section{The pentagram map}\label{sec:pentagram}

In this section we show that the pentagram map can be viewed as acting on a doubly periodic Q-net in $\CP^2$ by Laplace-Darboux dynamics.

\begin{figure}
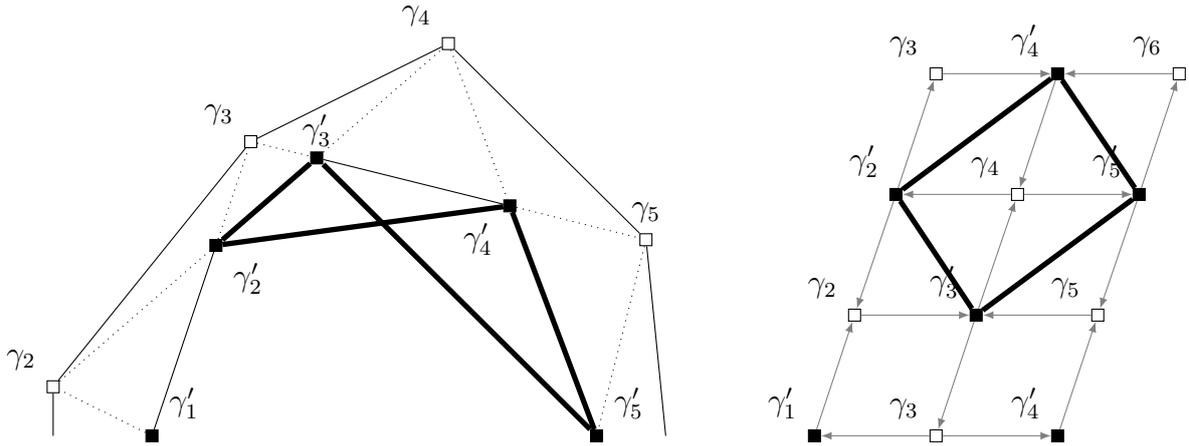


	\caption{Left: part of a polygon $\gamma$ and the polygon $\gamma'$ after applying the pentagram map. Right: The projective cluster structure associated to the pentagram map. We highlighted a quad of the associated Q-net $q_{\gamma'}$ in both pictures.}
	\label{fig:pentagram}
\end{figure}

\begin{definition}
	Let $\gamma: \Z_m \rightarrow \CP^2$ be a polygon. Define $\pent(\gamma) : \Z_m \rightarrow \CP^2$ to be the new curve (see Figure \ref{fig:pentagram}) such that
	\begin{align}
		\pent(\gamma)(k) = \gamma(k-2)\gamma(k) \cap \gamma(k-1)\gamma(k+1),
	\end{align}
	holds for all $1 \leq k \leq m$. The map $\pent: (\CP^2)^m \rightarrow (\CP^2)^m$ is called the \emph{pentagram map}.
\end{definition}
This map was introduced by Schwartz \cite{schwartz} and it was shown to be integrable \cite{ostpentagram}. It is called pentagram map because it behaves non-trivially for $m \geq 5$.

\begin{definition}
	Let $\gamma: \Z_m \rightarrow \CP^2$ be a polygon. 
	 Define the \emph{associated Q-net} $q_\gamma: V(\Z^2 / \{(m,0), (-2,1) \}) \rightarrow \CP^2$ by
	\begin{align}
		q_\gamma(k,l) &= \gamma(k + 2l),
	\end{align}
	for all $k,l\in \Z$.
\end{definition}

Clearly $q_\gamma$ is a Q-net because any four points in $\CP^2$ are on a plane.

\begin{theorem}\label{th:pentld}
	Let $\gamma: \Z_m \rightarrow \CP^2$ be a polygon and $q_\gamma$ the associated Q-net. Let $\Delta^2(q_\gamma)$ denote the Laplace transform (see Section \ref{sec:laplacedarboux}) in coordinate direction 2. Then 
	\begin{align}
		\Delta^2(q_\gamma) = q_{\pent(\gamma)}.
	\end{align}
	That is the iteration of the pentagram map corresponds to Laplace-Darboux dynamics on $q_\gamma$.
\end{theorem}
\proof{
	Every quad of $q_\gamma$ is thus made up of four points $(\gamma(k), \gamma(k+1), \gamma(k+3), \gamma(k+2))$ in counterclockwise order for $k\in \Z$. As a consequence, the two focal points of such a quad are $\pent(\gamma)(k+1)$ and $\pent^{-1}(\gamma)(k+2)$.\qed
}

This theorem can be found in different wording in a (notorious) talk given by Schief \cite{schieftalk}, but has not previously appeared in the literature. A cluster structure for the pentagram map has been found by Glick \cite{glickpentagram}.

\begin{theorem}
	The cluster structure and the cluster mutations associated to the pentagram map by Glick \cite{glickpentagram} coincides with the projective cluster structure of the pentagram map viewed as a TCD map.
\end{theorem}
\proof{Via resplits every $X$-variable of the pentagram map as a TCD map corresponds to a Laplace invariant (see Definition \ref{def:laplaceinvariant}) of the pentagram map viewed as a Q-net. The Laplace invariants are cross-ratios and coincide with Glick's variables.\qed}

\begin{remark}Except from explaining how the pentagram map is a TCD map we do not investigate the pentagram map in this thesis. However, we leave a few remarks for the interested reader.
	\begin{enumerate}
		\item In Section \ref{sec:affcluster} and Section \ref{sec:projclusterduality} we introduce several other cluster structures TCD maps. These cluster structures exist for the pentagram map as well, because the pentagram map can be viewed as a periodic reduction of a Q-net.
		\item There is follow-up work on the cluster structure of Glick relating it to directed networks \cite{gstvnetworks}.
		\item Until recently the precise relation between the cluster, Poisson and Hamiltonian structures found for the pentagram map and the dimer cluster integrable system \cite{gkdimers} was still open, although there was a sketch in lecture notes of Glick and Rupel \cite{gdnotes}. However, this was resolved in recent work by Izosimov \cite{izosimovnetworks}. In general, it would be interesting to compare our approach to the general framework in terms of Poisson-Lie groups developed by Izosimov \cite{izosimovpoissonlie}.
		\item Note that we explained how the pentagram map is a TCD map $\hat T:\tcdp \rightarrow \CP^2$ on a TCD $\tcd$ on the torus. It is also possible to consider $\hat \tcd$ on the universal cover of the torus instead. Then we can consider a TCD map $T:\hat \tcdp \rightarrow \CP^2$, such that $\hat T(w) = T(w)$ for all white vertices of $\pb$. However, in this case it is not necessary that the homogeneous lifts captured in the VRC $\hat R$ satisfy $\hat R(w) = R(w)$ for all white vertices of $\pb$. Instead, it is possible that $\hat R(w) = c_w R(w)$ for some factors $c_w\in \C$. We do not go into detail here, but it turns out that the factors $c_w$ are not arbitrary, and instead depend only on two constants $c_1,c_2$ associated to two non-trivial cycles on the torus that generate its fundamental group. This aspect should be of interest with respect to both integrable systems and statistical mechanics in the sense of dimer configurations. See Section \ref{sec:dimerinvariants} for a relation between dimer configurations and TCD maps (albeit not on the torus).
	\end{enumerate}

\end{remark}

\section{TCD maps from projective cluster variables}\label{sec:fromprojinvariants}

In Section \ref{sec:projcluster} we introduced projective invariants -- the projective cluster variables -- of TCD maps. But how do we reconstruct a TCD map from these invariants? We want to give an algorithm that constructs such a TCD step by step. 

Let us revisit some conclusions from Section \ref{sec:sweeps} that we employ in the proof. We recall that for every labeled minimal TCD we defined the li-orientation $\lio$ on the graph $\pb$. Each black vertex has only one incoming edge, while every white vertex has at most one outgoing edge. Moreover, each face $f$ has only one black boundary vertex that has an incoming edge from a white vertex not on the boundary of $f$. Assume the maximal dimension of the labeled minimal TCD $\tcd$ is $k$. Viewed as a poset, $\lio$ has $k+1$ minimal elements, and these minimal elements are all on the boundary. Also consider a linear extension $\varepsilon: W(\pb) \rightarrow \N$ of $\lio$ restricted to the white vertices. A white vertex $w$ that is not a minimal element has a unique outgoing edge $(b,w)$. The other two white vertices adjacent to $b$ therefore necessarily need to appear before $w$ in $\varepsilon$. There is also the face $f$ adjacent to $b$ but not to $w$. All white vertices on the boundary of $f$ do also appear before $w$ in $\varepsilon$.

We are now ready to define the algorithm. Note that this algorithm may fail, we will discuss details of this after the definition.

\begin{definition}\label{def:constralgo}
	Let $\tcd$ be a labeled minimal TCD, $n\in \N_{>0}$ and let $X': F(\pb) \rightarrow \Cx$. In the following we define the \emph{construction algorithm}, that either fails or produces a TCD map $T:\tcdp \rightarrow \CP^n$ such that its projective invariants satisfy $X_f = X'_f$ for all faces $f\in F(\pb)$. 
	\begin{enumerate}
		\item Choose a linear extension $\varepsilon: W(\pb) \rightarrow \N$ of $\lio$. 
		\item Choose the points associated to the minimal elements of $\lio$ arbitrarily in $\CP^n$, but such that any two minimal elements that share a black vertex span a line.
		\item Find the smallest white vertex $w$ in $\varepsilon$ that has not been placed yet. Let $b$ be the unique black vertex such that $(b,w)$ is the incoming edge at $b$, let $w',w''$ be the other two neighbours of $b$ and let $f$ be the face adjacent to $b$ but not $w$. There are now three cases:
		\begin{enumerate}
			\item If $T(w') = T(w'')$ then the algorithm fails.
			\item Else if $f$ is a boundary face then place $T(w)$ arbitrarily on the line through $T(w')$ and $T(w'')$, but not such that it coincides with one of those two points.
			\item Else there is a unique point $T(w)$ on the line through $T(w')$ and $T(w'')$ such that $X(f) = X'(f)$, choose that point.
		\end{enumerate}
		\item Repeat step (3) until all points are placed.\qedhere
	\end{enumerate}
\end{definition}

\begin{lemma}\label{lem:opensetoftcdmaps}
	Let $\tcd$ be a labeled minimal TCD with maximal dimension $k$ and let $n\in \N_{>0}$. The set of functions $X': F(\pb) \rightarrow \Cx$, such that there is a TCD map $T:\tcdp \rightarrow \CP^n$ with projective invariants $X'$ is an open an full-dimensional subset of $\C^{|F|}$. The set of TCD maps $T:\tcdp \rightarrow \CP^n$ that possess a fixed set of projective invariants is either empty or an open and full-dimensional subset of $\C^{k n + k + n}$.
\end{lemma}
\proof{
	Consider the algorithm of Definition \ref{def:constralgo}. Any TCD map can be constructed by this algorithm. Also, note that the number of minimal elements of $\lio$ equals $k+1$, while the number of boundary faces that are hit by step (3b) is $k$. Therefore the (complex) degrees of freedom of this construction are $(k+1)\times n + k$. Moreover, the only way the algorithm fails is if the point of a white vertex $w$ is placed such that it coincides with a point of a white vertex $w'$, such that $w$ and $w'$ share a black vertex that comes later in $\varepsilon$. As every white vertex is only incident to a finite number of black vertices, in each step there is only a finite number of disallowed choices. Therefore the whole set of TCD maps with for a given TCD but arbitrary invariants $X$ is an open and full dimensional subset of $\C^{k n + k + n + |F|}$. It evidently splits into the two spaces of the Lemma.\qed
}

In the special case that $\tcd$ is a sweepable TCD, the minimal elements of $\lio$ are consecutive boundary vertices and the boundary faces encountered in the construction algorithm are exactly the boundary faces between the consecutive minima. Thus the degrees of freedom in the construction are located along the boundary. The construction algorithm therefore is simply one particular way of propagating 1-dimensional, hyperbolic Cauchy-data in a 2-dimensional system. This system is discretely integrable in the sense that one can change the combinatorics via 2-2 moves but this will not affect the construction of the TCD map away from where we changed the combinatorics.

\begin{theorem}\label{th:uniquefrominvariants}
	Let $\tcd$ be a sweepable, labeled TCD of maximal dimension $k$. Let $T,T':\tcdp \rightarrow \CP^k$ be two TCD maps, such that the points of $T$ as well as the points of $T'$ span $\CP^k$. Let $X,X'$ be the projective invariants of $T,T'$. If $X=X'$ then there is a projective transformation $f$ such that $f(T') = T$. In other words, the projective invariants $X$ characterize a TCD map that attains its maximal dimension up to projective transformations.
\end{theorem}
\proof{
	Given a sequence of $k+1$ different points spanning $\CP^k$ and a hyperplane that does not contain any of those points, it is always possible to find a projective transformation mapping the points and the hyperplane to any other such sequence of points and a hyperplane. In our case the points are the $k+1$ points of $T$ (resp. $T'$) associated to minimal elements of $\lio$. Moreover, in the sweepable case the minimal elements of $\lio$ are consecutive boundary vertices $v_1,v_2,\dots, v_{k+1}$ of $\pb$. Fix an auxiliary set of numbers $X^\partial_{i} \in \Cx$ for $i$ from 1 to $k$. Denote by $f_i$ the boundary face with $v_i,v_{i+1}$ on its boundary. Introduce the points $M_1,M_2,\dots M_i$, such that each point $M_i$ is on the line $T(v_i)T(v_{i+1})$ and such that for each face $f_i$ the multi-ratio of the points along the boundary of $f_i$ together with $M_i$ is $X^\partial_i$. Define the hyperplane $H$ as the span of all $M_i$. Repeat this procedure for $T'$, with the same $X^\partial_i$ to obtain the points $M'_i$ and the hyperplane $H'$. By the initial argument of the proof, there is a projective transformation $f: \CP^k \rightarrow \CP^k$ that maps $H'$ to $H$ and the points $T'(v_i)$ to $T(v_i)$. As a result, $f$ also maps each point $M'_i$ to $M_i$. We can now verify that indeed $f(T') = T$ by following the steps of the construction algorithm, where we see that the combination of the values $X^\partial_i$ and the marked points $M_i$ defines the choices in step (3b) of the algorithm. Therefore all choices in the algorithm are the same and the theorem follows.\qed
}

To the best of our knowledge, this is a new type of result within discrete differential geometry examples. Usually, in discrete differential geometry one considers a map from a large or even infinite graph to a projective space of fixed (and small) dimension like $\RP^3$. The maximal dimension of these examples however grows with their size, therefore in this setting one is content to show that the invariants together with boundary data determine the whole map. Theorem \ref{th:uniquefrominvariants} however states that in maximal dimension, there is actually no relevant boundary data, as the map (up to projective transformations) is determined by the projective invariants. Therefore, the boundary data that appears in smaller dimensions simply captures how we project the unique map in maximal dimension to the smaller dimensional space.

Also note that in the case of non-sweepable TCDs, we expect a similar result. The degrees of freedom are clearly the same as for sweepable TCDs. However, because the minima are not consecutive it is not quite so straightforward to introduce marked points that define a hyperplane.

We now consider the case of real projective spaces and real positive invariants. In this case we strengthen our previous results to include an existence statement. 

\begin{figure}
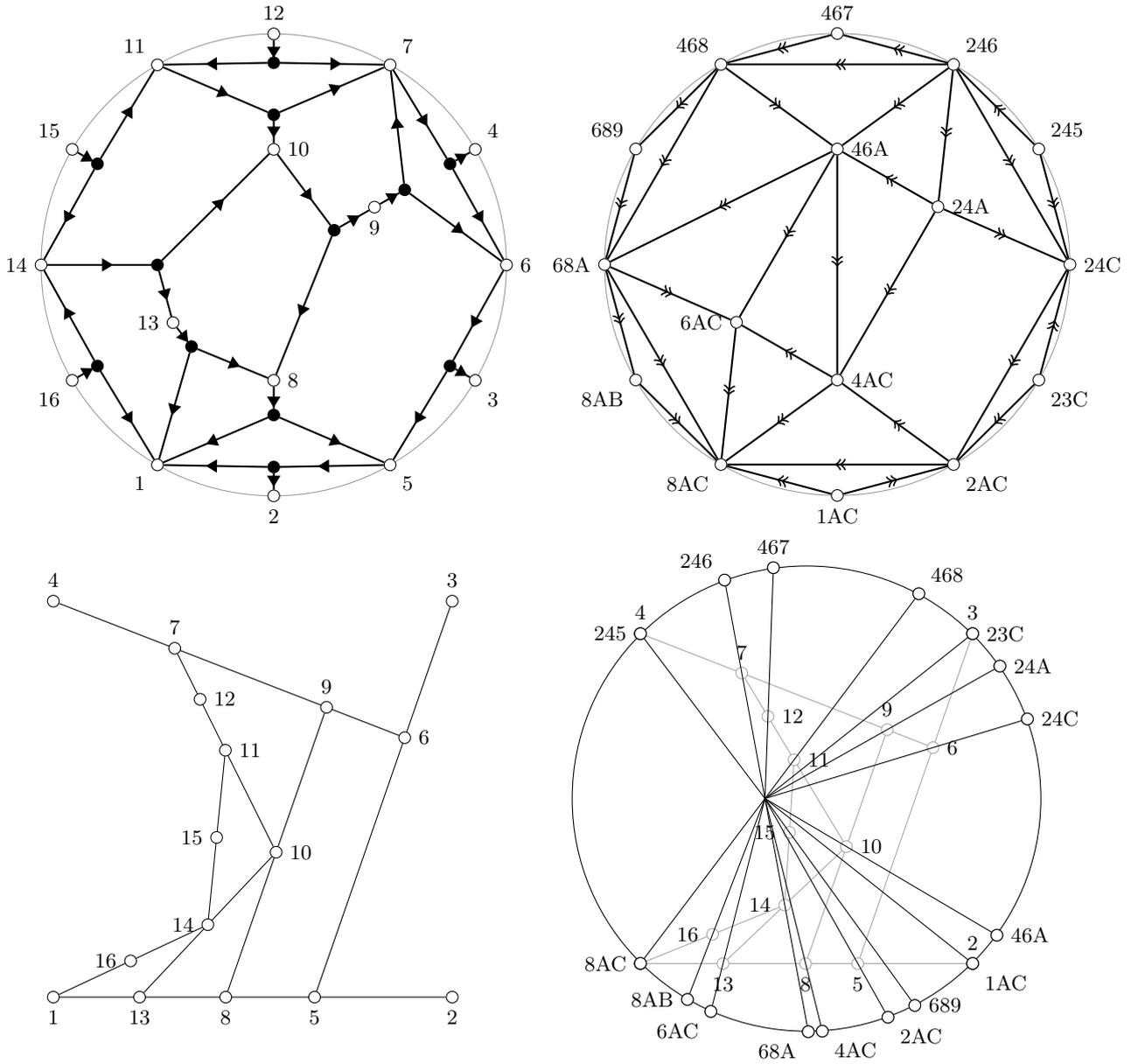

	\SMALL

	\caption{Top left: The li-orientation $\lio$ of a TCD $\tcd$ and a labeling that is a linear extension of $\lio$. Top right: The lr-orientation $G^\lrp$ and the labeling by $\mathcal{A}'$. Bottom left: A convex drawing of a TCD on $\tcd$ with positive $X$-variables. Bottom right: A central projection of the convex drawing to a circle representing $\RP^1$.}
	\label{fig:convextcd}
\end{figure}

\begin{definition}\label{def:lrposet}
	Let $\tcd$ be a labeled minimal TCD with $n$ strands. Recall that we defined the map $\mathcal A'$ from white vertices of $\pb$ to a subset of $\{1,2,\dots, n\}$ in Section \ref{sec:tcdconsistency}, where $i\in \mathcal A'(w)$ if $w$ is to the right of strand $i$. Construct an oriented graph $G^\lrp$ (see Figure \ref{fig:convextcd}) with vertex-set equivalent to the white vertices $W(\pb)$ of $\pb$. There is an edge in $G^\lrp$ from $w_1$ to $w_2$ whenever 
	\begin{enumerate}
		\item there is a black vertex $b$ adjacent to both $w_1$ and $w_2$ and
		\item $\mathcal A'(w_1) \leq \mathcal A'(w_2)$ in the lexicographic order, with $\mathcal A'(w)$ written in increasing order itself.
	\end{enumerate}
	The graph $G^\lrp$ is acyclic as it is a restriction of the complete lexicographic order on the white vertices.	We define the \emph{left-right poset} or abbreviated the \emph{lr-poset} $\lrp$ on the white vertices $W(\pb)$ of $\pb$ as the transitive closure of the acyclic orientation given by $G^\lrp$, that is $w_1 \leq w_2$ in $\lrp$ if there is an oriented path from $w_2$ to $w_1$ in $G^\lrp$.
\end{definition}

\begin{theorem}\label{th:positiveinvariants}
	Let $\tcd$ be a minimal TCD with maximal dimension $k$. Let $X': F(\pb) \rightarrow \R_{>0}$ be a positive function. Then there is a TCD map $T:\tcdp \rightarrow \RP^k$ with projective invariants $X=X'$. If $\tcd$ is sweepable, then up to projective transformations, there is a unique TCD map $T:\tcdp \rightarrow \RP^k$ that spans $\RP^k$ and has projective invariants $X=X'$.
\end{theorem}
\proof{Let us begin under the assumption that $\tcd$ is sweepable and choose a labeling of $\tcd$. Then choose a line $\ell$, an affine chart $\R$ of $\ell$ and a central projection $\pi:\RP^k\rightarrow \ell$ such that the $k+1$ images of the minimal elements of $\lio$, sorted in the order they appear on the labeled boundary of $\tcd$, appear in increasing order on the affine chart of $\ell$. This is always possible as we can map $k+1$ points in general position in $\RP^k$ to any other such points, and thus the desired position for their projections can also be achieved. Now apply the construction algorithm (see Figure \ref{fig:convextcd} for an example). In step (3b) choose the new points in between the two corresponding other points on the affine chart of $\ell$. The sign of the invariants $X$ ensures that in each step (3c), the new point is also in between the two other points. To see this, note that each black vertex on the boundary of a face $f$ contributes a minus sign to $X_f$, except the maximal black vertex in $f$. To see this, before starting the construction algorithm initialize $M = (M_1,M_2,\dots, M_{k+1})$ as the sequence of minima of $\lio$. Note that $M_i<_\lrp M_{i+1}$ for each $i$ with respect to the lr-poset (see Definition \ref{def:lrposet}). At the same time, the images satisfy the same thing on $\ell$, that is $T(M_i) < T(M_{i+1})$ for all indices $i$. In each step of the construction algorithm, update the sequence $M$ as follows: Whenever we place a new white vertex $w$ at the same black vertex as $w',w''$ (notation as in Definition \ref{def:constralgo}), add $w$ also to the sequence $M$ in between the consecutive elements $w'$ and $w''$. In this way $M$ remains sorted with respect to $\lrp$, and the images of the elements of $M$ remain sorted in $\ell$. Then remove any white vertex which has no adjacent black vertices left that have not already been considered in the construction algorithm. This also preserves the sorting of $M$ and its images. Because of this preservation of sorting during the algorithm, it is never possible that two points in $M$ coincide. The images of the vertices $w',w''$ that we consider in step (3) of the construction algorithm are always in $M$, and therefore the failure in step (3a) can never occur. In the sweepable case, we have therefore proven the existence of $T$. The uniqueness follows from Theorem \ref{th:uniquefrominvariants}. In the non-sweepable case, recall that any TCD $\tcd$ is the sub-diagram of a sweepable TCD $\tcd'$ according to Lemma \ref{lem:subsweepable}. One can therefore construct a TCD map $T'$ for $\tcd'$ in the way we explained. We can restrict $T'$ to $\tcd$ and obtain a TCD map $T$ that satisfies the assumptions of the theorem.\qed
}

\begin{remark}\label{rem:convextcd}
	Consider Theorem \ref{th:positiveinvariants} in the case of maximal dimension $k>1$. Then instead of projecting to a line $\ell$, we can project onto a 2-plane $E$. Moreover, instead of monotony on $\ell$ we use convexity in $E$. We place the minimal elements of $\lio$ on a convex curve. In the iteration we place the new points always on diagonals of a convex curve that corresponds to the iterated sequence $M$, see Figure \ref{fig:convextcd}. The two results in $\ell$ resp. $E$ are related by a central projection from a point inside the convex hull of the end-curve of the construction algorithm. In fact, viewing Figure \ref{fig:convextcd} one observes that this projection of a TCD map to $E$ are indeed maps that we have encountered before in Section \ref{sec:extgraphs}.
\end{remark}

The following corollary refers to T-graphs that we will define later in Section \ref{sec:tgraphcluster}.

\begin{corollary}
	Let $T: \tcdp \rightarrow \RP^k$ be a TCD map that attains its maximal dimension $k$ such that the projective invariants $X$ of $T$ are positive. Then there is a central projection $\pi$ onto a 2-plane $E$ in $\RP^k$ and an affine chart of $E$ in which $\pi(T)$ is a one-sided T-graph.
\end{corollary}
\proof{Direct consequence of Theorem \ref{th:positiveinvariants} and Remark \ref{rem:convextcd}. Whenever we place a white vertex $v$ in $E$, we cut off a half-space $H_v$. The other segments corresponding to black vertices adjacent to $v$ cannot enter these half-spaces. Also by construction the (open) segments are disjoint.\qed
}

However, not every T-graph has positive projective invariants. For example, one can see a T-graph in Figure \ref{fig:tgraphexample} that has several faces where the number of boundary edges with interior points is not an odd number. Still, we will be able to employ the construction algorithm in Section \ref{sec:tgraphcluster} for larger classes of T-graphs.

\begin{remark}
	Note that there is another result \cite{vrc} about uniqueness of TCD maps with given invariants and TCD. It states that a TCD map that attains maximal dimension, has given projective invariants and all of its boundary vertices prescribed, is unique if it exists. From the point of view of the results of this section, it seems overly restrictive to prescribe all of the boundary. However, what is still missing is a precise understanding of the moduli space of TCD maps that attain maximal dimension, particularly in the case that the invariants are real non-positive (or even complex) or the TCD is non-sweepable (or non-minimal). The reader is invited to check, that there are many TCDs for which any set of invariants can be uniquely realized (up to projective transformations, in maximal dimension). However, it is doubtful that this is true for all TCDs, although we did not find a counter-example.
\end{remark}

\begin{remark}\label{rem:altconstruction}
	Let us also remark that there is an \emph{alternative construction algorithm} that is closer to the construction of T-graphs (see Section \ref{sec:tgraphcluster}) introduced by Kenyon and Sheffield \cite{kenyonsheffield}. Again, we are looking for a way to construct a TCD map $T: \tcd \rightarrow \CP^n$ for some given TCD $\tcd$ and $n$ as well as prescribed projective invariants $X'$. We split the problem in two. First, we solve the \emph{gauge problem}. That is we construct edge weights for the relations of the VRC $\vrc$, such that the projective invariants $X$ read off of $\vrc$ coincide with $X'$. Consider the dual $\pb^*$ of $\pb$ and a rooted spanning tree $\mathcal S$ of $\pb^*$ that is rooted at the outer face. Assume all edges of $\mathcal S$ are oriented towards the root. For every edge $e$ such that $e^* \notin \mathcal S$ set the corresponding edge weight in $\vrc$ to 1. Now begin at a leaf $f$ of $\mathcal S$ ($f$ is a face of $\pb$) and consider the unique outgoing edge $e^* \in \mathcal S$. As we are at a leaf this is the only edge on the boundary of $f$ in $\mathcal S$, and thus there is a unique edge weight for $e^*$ such that the projective invariant $X$ at face $f$ coincides with the prescribed value $X'$. Then repeat this on $\mathcal S \setminus \{f\}$ and iterate until all edges are assigned edge weights. In each step the weight is unique because there is only one edge leaving each face. Now choose any acyclic, perfect orientation $\chi$ of $\pb$. Perfect orientation in the sense of Postnikov \cite{postgrass} meaning each black vertex has one incoming and each white vertex has one outgoing edge. Now proceed as in the construction algorithm of Definition \ref{def:constralgo}, but replace the li-orientation by $\chi$. In each step (3) where we add the point of a new white vertex, the point is now determined by the relation at the corresponding black vertex in $\vrc$. 
\end{remark}

The following lemma is also useful later on, but easy to prove here.

\begin{lemma}\label{lem:rightofstrand}
	Consider a minimal TCD $\tcd$ and a TCD map $T: \tcdp \rightarrow \CP^k$ and a strand $s$. Let $w_1,w_2,\dots, w_m$ be all the white vertices to the right of $s$ and let $U$ be the space spanned by the images of the white boundary vertices to the right of $s$. Then every point $T(w_i)$ for $1\leq i \leq m$ is in $U$.
\end{lemma}
\proof{Choose a labeling of $\tcd$ such that $s$ carries the label 1 and consider the li-orientation of $\pb$. Then at all the black vertices on $s$, the outgoing edges of the li-orientation point to the left of $s$. Therefore the minimal elements that any white vertex to the right of $s$ sees in its down-set are to the right of $s$.\qed
}

\section{The affine cluster structure of a TCD map}\label{sec:affcluster}

\begin{figure}
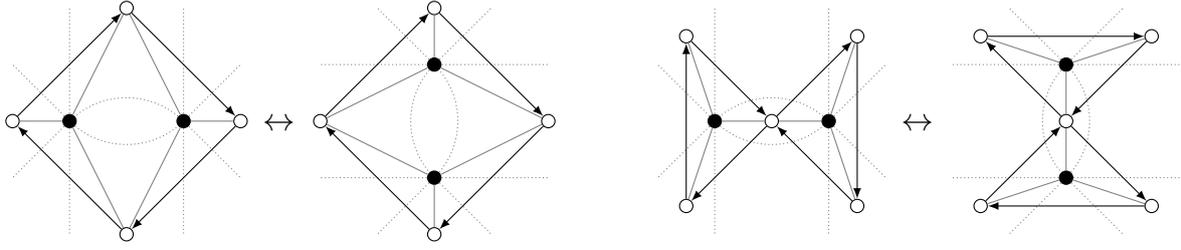

	\centering
	
	\caption{The affine quiver (white vertices and arrows) before and after a spider move (left) as well as before and after a resplit (right).}\label{fig:twotwoaffinequiver}
\end{figure}

\begin{definition} \label{def:affinecluster}
	Let $T$ be a TCD map and let $\mu$ be the edge weights of the associated VRC in affine gauge with respect to an affine chart, and denote by $H$ the hyperplane at infinity in this chart. Assume that $H$ is generic with respect to $T$. The vertices of the quiver of the \emph{extended affine cluster structure} $\aff_{H}^\times(T)$ are located at the white vertices of $\pb$. The number of arrows $\nu_{ww'}$ between two white vertices $w,w'$ is defined such that
	\begin{align}
		\nu_{ww'} = \phantom{-} &|\{(w,b,w') : w' \mbox{ clockwise after } w \mbox{ around } b\in B\}| \\-\ &|\{(w,b,w') : w' \mbox{ counterclockwise after } w \mbox{ around } b\in B\}|.
	\end{align}
	Let $w$ be a white vertex of $\pb$ and introduce labels such that the neighbourhoods in counterclockwise order are
	\begin{align}
		N(w) = (b_1,b_2,\dots,b_m)\quad \mbox{ and } \quad N(b_i) = (w,w_i,w'_i),
	\end{align}
	for all $b_i$ adjacent to $w$. Then the \emph{affine cluster variable} $Y_w$ is the alternating ratio
	\begin{align}
		Y_w = (-1)^{m+1} \prod_{i=1}^m \frac{\mu(b_i,w'_i)}{\mu(b_i,w_i)}.
	\end{align}
	We also call an affine cluster variable $Y_w$ a \emph{star-ratio}. The \emph{affine cluster structure} $\aff_{H}(T)$ is the restriction of the extended affine cluster structure to interior white vertices of $\pb$.
\end{definition}

\begin{remark}
	Let us consider a VRC in affine gauge. Recall that for a TCD $\tcd$ there is also Definition \ref{def:pbm} that introduces the graph $\pbm$. At a black vertex $b$ of $\pb$ denote the white vertices and faces adjacent to $b$ by $w_1,f_3,w_2,f_1,w_3,f_2$ in counterclockwise order. Identify each edge $(b,w_i)$ in $\pb$ with the edge $(b, f_i)$ in $\pbm$. To an edge in $\pbm$ assign the inverse of the edge weight of the corresponding edge of $\pb$ in $\vrc$. Then by definition the star-ratio $Y_w$ for $w$ in $\pb$ is $X_f$ in $\pbm$.
\end{remark}

Next, we are interested in the effect of the 2-2 moves on the quiver and cluster variables.

\begin{figure}
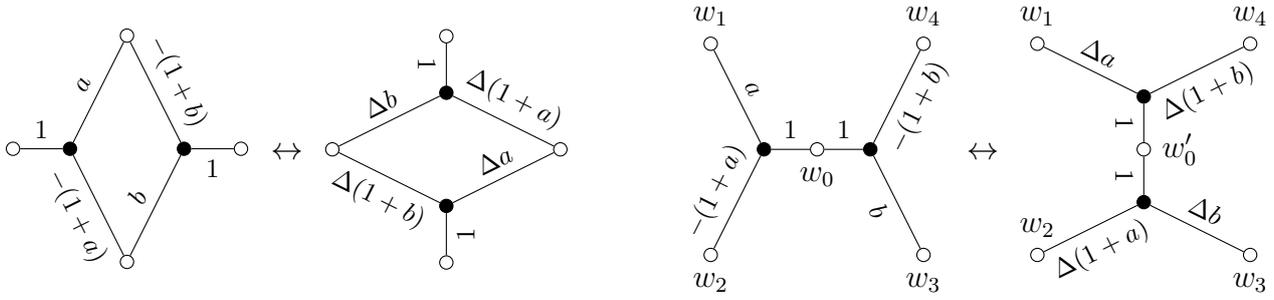

	\centering

	\caption{Change of edge weights in affine gauge, with $\Delta=-(1+a+b)^{-1}$.}\label{fig:affineweights}
\end{figure}

\begin{lemma}
	The resplit corresponds to a mutation in the affine cluster algebra, the spider move leaves the affine cluster structure invariant (cf. Figure \ref{fig:twotwoaffinequiver}).
\end{lemma}
\proof{Refer to Figure \ref{fig:affineweights} for the change of edge weights. The star-ratios are clearly unchanged by the spider move for any of the four involved vertices. For the resplit at $w_0$ we observe that
	\begin{align}
		Y_{w_0}= - (1+a^{-1})(1+b^{-1}) = \tilde{Y}_{w'_0}^{-1}
	\end{align}
	holds. This is the same expression as in the projective case and similarly we observe that 
	\begin{align}
	\frac{\tilde Y_{w_2}}{Y_{w_2}} = \frac{1}{a\Delta b} = - (a^{-1} + b^{-1} + a^{-1}b^{-1}) = 1 + Y_{w_0}.
	\end{align}
	The other calculations proceed in the same manner.\qed
}

In both the affine and the projective cluster variables the edge weights appear in pairs of edges that are incident to the same black vertex. Thus it is not surprising that we can also express the affine cluster variables via signed distance ratios.

\begin{figure}
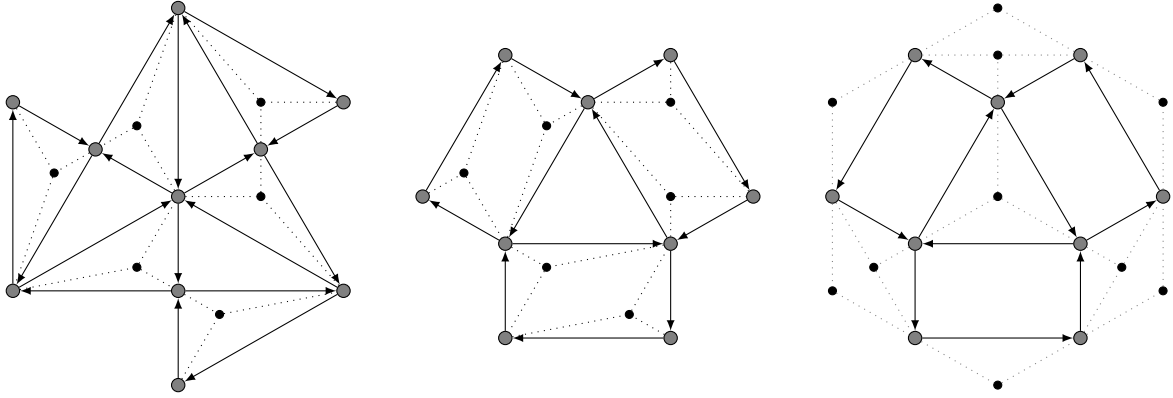


	\caption{The affine quivers (arrows and gray vertices) for Q-net, Darboux map and line compound. The graph $\pb$ is drawn dotted. We recognize the hexahedral quiver, the cuboctahedral quiver and the reversed cuboctahedral quiver.}
	\label{fig:affquivers}
\end{figure}

\begin{lemma}\label{lem:starratioviadistances}
	Let $T$ be a TCD map in affine gauge with respect to a generic hyperplane $H$. Then star-ratio at a white vertex $w$ is
	\begin{align}
		Y_w = - \prod_{i=1}^{m} \frac{T(w_i)-T(w)}{T(w'_i)-T(w)},
	\end{align}
	where the quotients are signed distance quotients in any Euclidean chart and $w_i,w_i'$ are as in Definition \ref{def:affinecluster}.
\end{lemma}
\proof{
	Direct consequence of Lemma \ref{lem:edgeratios}.
	\qed

}
The star-ratios are defined in a fixed affine gauge. The affine gauge is determined by a choice of hyperplane $H$ in $\CP^n$ as the plane at infinity for an affine chart. The affine transformations are the projective transformations that map $H$ to itself. Moreover, affine transformations do not change oriented distance ratios on a line. Thus, the name affine cluster structure and affine cluster variables are justified. As in the case of the projective cluster structure, it is possible to regain the combinatorics of $\pb$ from the affine quiver, this time up to spider moves. Clearly the vertices of the quiver correspond to the white vertices of $\pb$. As one can see in Figure \ref{fig:affquivers} or from the definition, each clockwise (black) face of the quiver corresponds to a line that is represented by $(d_f-2)$ black vertices. The corresponding images of the white vertices adjacent to a black face of the quiver are then on the line represented by the black face.

\begin{figure}
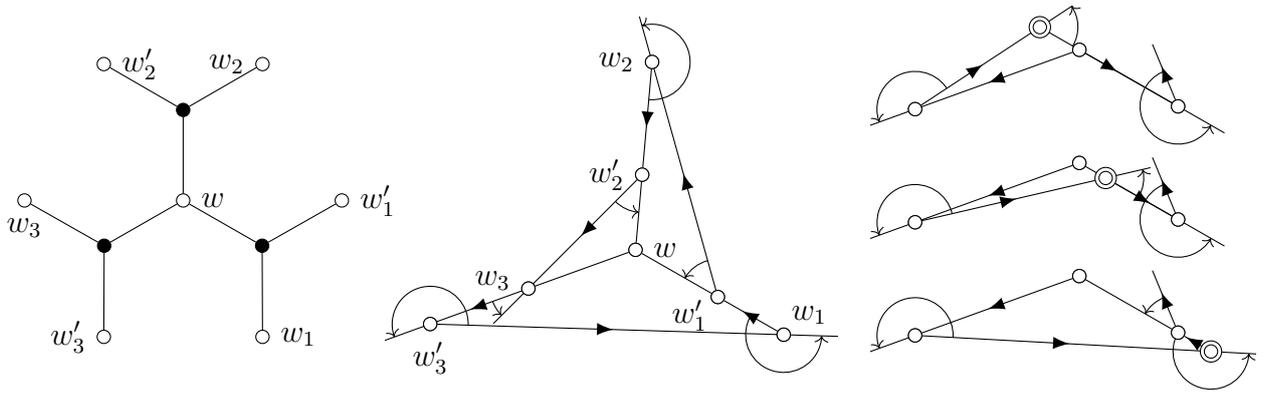


	\caption{How to express star-ratios via angles. On the right we observe that the sign of the product of sines of angles changes in accordance with the sign of the product of oriented distance ratios, the moving point is highlighted.}
	\label{fig:affviaangles}
\end{figure}

\begin{lemma}\label{lem:srviaangles}
	Let $T: \tcd \rightarrow \RP^n$ be a TCD map in affine gauge and choose an additional similarity structure on $\R^n \subset \RP^n$. With the similarity structure we are able to define the counterclockwise angles $\alpha_k$ from direction $T(w_k)\rightarrow T(w'_k)$ to direction $T(w'_{k-1})\rightarrow T(w_k)$, as well as the angles $\beta_k$ from direction $T(w'_k)\rightarrow T(w_{k+1})$ to direction $T(w_{k})\rightarrow T(w'_k)$. Then the star-ratio $Y_w$ at a white vertex $w$ is equal to
	\begin{align}
		Y_w &= (-1)^{m+1} \prod_{i=1}^m \frac{\sin \beta_k }{\sin \alpha_k}.
	\end{align}	
	Thus it is possible to express star-ratios via oriented angles.
\end{lemma}
\proof{
	This is a consequence of Lemma \ref{lem:starratioviadistances} and the law of sines. In particular, the law of sines yields that
	\begin{align}
		\left|\frac{\sin \beta_{k-1}}{\sin \alpha_k}\right| = \left|\frac{w_k - w}{w'_{k-1} - w}\right|
	\end{align}
	holds. To verify that the signs are correct, we begin with the case in which all oriented distance ratios are positive. In this case, every black vertex contributes a minus sign to the angle expression, confirming the claim (see also Figure \ref{fig:affviaangles}). On the other hand, whenever we move a point along a line that corresponds to a black vertex until the corresponding oriented distance ratio becomes negative, then also exactly one of the signs in the angle expression changes (see Figure \ref{fig:affviaangles}).\qed
}

This lemma yields a way to express the cluster variables assigned to TCD maps via angles. It therefore raises the possibility that integrable systems that are described by angle coordinates actually relate to the affine cluster structure of some TCD map. In fact, we use this viewpoint in Section \ref{sub:circularq} to introduce an angle-based cluster structure for circular Q-nets, which relates to the treatment by Bazhanov, Mangazeev and Sergeev \cite{bmscircular}.

\begin{remark}
	The projective and affine quivers come from the triple crossing diagram in similar manner. In particular the affine quiver of a diagram $\tcd$ is just the projective quiver of the TCD dual $\iota(\tcd)$ (see Definition \ref{def:tcddual}). Given a TCD map $T:\tcdp\rightarrow \CP^n$ one may therefore wonder if there is actually a second TCD map $T':\tcdp'\rightarrow \CP^n$ such that the projective and affine cluster structures are reversed in $T'$. In order to deal with this question we need to investigate the behaviour of the cluster structures under sections and projective duality first, but we will shed some light onto this question in Section \ref{sec:perfdual}.
\end{remark}

Let us consider our standard DDG examples, see Figure \ref{fig:affquivers}. This time we notice that both Darboux maps and line compounds feature a cuboctahedral quiver, while Q-nets possess a hexahedral quiver. As in the case of the projective quivers, the orientations of the two cuboctahedral quivers are not the same. Note that in the case of stepped surfaces, the affine quiver of a line compound is also cuboctahedral, but for general quad-graphs this is not the case.

\begin{remark}
	We gave a sequence of eleven 2-2 moves that corresponds to the cube flip of a Q-net in Figure \ref{fig:tcdcubeflip}. However, if we flip the cube such that after the flip the role of black and white vertices of the quad-graph is interchanged, then there is a shorter sequence of only ten 2-2 moves, see Figure \ref{fig:tensteps}. Four of these 2-2 moves are spider-moves and correspond to the four mutations of the cuboctahedral flip, see Definition \ref{def:cuboflip}. Six of these 2-2 moves are resplits and correspond to six mutations in the hexahedral quiver. This is a shorter version of the hexahedral flip, see Definition \ref{def:hexaflip}, where the flip also interchanges the role of black and white vertices of the quad-graph. Thus, the TCD flip decomposes into a combination of a cuboctahedral and a hexahedral flip.
\end{remark}

\newpage
\thispagestyle{empty}
	
\begin{figure}[p]
	\vspace{-3mm}
	\hspace{-2.2cm}

	\vspace{-3mm}
	\caption{Sequence of $\tcd$, $\pb$, proj. and affine quivers of a cube flip in a Q-net.}
	\label{fig:tensteps}
\end{figure}


\chapter{Cluster structures and projective operations}\label{cha:clustergeometry}

\section{Cluster variables and projections of TCD maps}
As we have discussed before, the projective cluster variables are invariant under projective transformations as well as under central projections. The affine cluster variables on the other hand are in general not invariant under central projections, but there is an exception. Before we proceed, let us make another genericity definition to enable us to be precise about projections.

\begin{definition}
	Let $T:\tcdp \rightarrow \CP^n$ be a TCD map. A point $P\in \CP^n$ is \emph{cogeneric} if there is no black vertex $b$ of $\pb$ such that $P\in L(b)$.
\end{definition}

\begin{theorem}\label{th:projectionandaffine}
	Let $T: \mathcal T \rightarrow \CP^n$ with $n>1$  be a TCD map, let $H$ be a hyperplane and let $P \in \CP^n\setminus H$ be a cogeneric point. Let $\pi: \CP^n \rightarrow H$ be the projection with center $P$ to $H$. Then
	\begin{align}
		\pro(\pi(T)) &= \pro(T)
	\end{align}
	holds. Moreover, let $E$ be a hyperplane that contains $P$, such that $E$ is 1-generic with respect to $T$ and $E\cap H$ is 1-generic with respect to $\pi(T)$. If $T$ and $\pi(T)$ are 1-generic then 
	\begin{align}
		\aff_{H\cap E}(\pi(T)) &= \aff_E (T)
	\end{align}
	holds as well.
\end{theorem}
\proof{The cogenericity of $P$ guarantees the 0-genericity and thus well-definedness of $\pi(T)$. The first equation in the theorem is a consequence of the invariance of multi-ratios and thus the projective cluster variables under central projection, see Lemmas \ref{lem:projinvariants} and \ref{lem:projclusterviadistances}. For the second equation, consider affine coordinates of $\CP^n$ for which $E$ is the hyperplane at infinity. Then the projection $\pi$ is in fact a parallel projection. Remember that we can think of the star-ratios as given by the distance ratios as in Lemma \ref{lem:starratioviadistances}. But as the distance ratios of three points on a line are invariant under parallel projection, so are the star-ratios and the second equation is proven as well.\qed 
}

\section{Cluster variables and sections of TCD maps}\label{sec:sectioncluster}
In Section \ref{sec:sections} we studied the geometry and combinatorics of sections of TCD map. Now we want to build on these results and study the behaviour of the cluster variables of sections. In particular we will show that the affine cluster structure of a TCD map can be identified with the projective cluster structure of a section of that TCD map.

\begin{lemma} \label{lem:srviacr}
	Let $T$ be a 1-generic TCD map and $w$ a white vertex of $\pb$. Use the same labeling for the vertices close to $w$ as in Definition \ref{def:affinecluster}. Let $H$ be a generic hyperplane and let $Y$ be the affine cluster variables of $\aff_H(T)$. Let $\ell_k = T(w)T(w_k)$ be the lines around $w$ and denote by $S_k = \ell_k \cap H$ a fourth point on $\ell_k$. Then we can express the star-ratio $Y_w$ with respect to $H$ as
\begin{equation}
	Y_w = \prod_k \cro (T(w),T(w_k),S_k,T(w'_k)).\qedhere
\end{equation}
\end{lemma}
\proof{Choose affine coordinates with $H$ at infinity. Then we have that the distances $|T(w_k)-S_k|$ and $|T(w'_k)-S_k|$ are infinite and thus
\begin{align}
	\cro (T(w),T(w_k),S_k,T(w'_k)) = \frac{T(w_k)-T(w)}{T(w'_k)-T(w)}
\end{align}
holds. Together with Lemma \ref{lem:starratioviadistances} the proof is complete.\qed
}

\begin{figure}
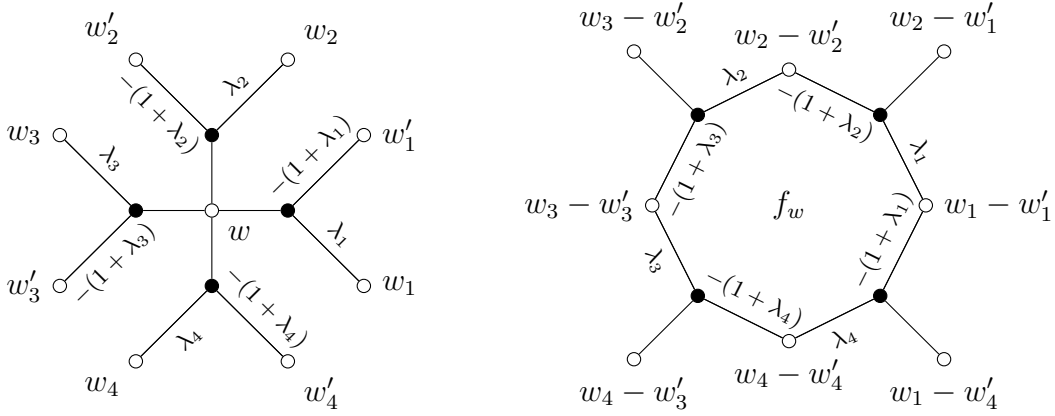

	
	\caption{The edge-weights in $R$ and $\sigma_H(R)$ as employed in the proof of Theorem \ref{th:affprojcluster}. Edges without label carry weight 1.}
	\label{fig:sectionweights}	
\end{figure}

We observe that the additional points occurring in the lemma are actually the points of the section $\sigma_H(T)$. Thus we are lead to the following theorem.
\begin{theorem}\label{th:affprojcluster}
	Let $\tcd$ be a minimal TCD and $T: \tcdp \rightarrow \CP^n,\ n>1$ be a 1-generic and flip-generic TCD map and $H$ a generic hyperplane. Then there is the identity
	\begin{equation}
		\pro(\sigma_H(T)) = \aff_H(T).\qedhere
	\end{equation}
\end{theorem}
\proof{
	We have already noticed before (in Section \ref{sec:sections}) that the white vertices of $\pb$ are in bijection with the faces of $\sigma_H(\pb)$, and the interior white vertices are in bijection with the interior faces as well. This is of course necessary for the theorem to hold on a combinatorial level and allows us to prove the theorem per white vertex $w$ of $\pb$. Let $m$ be the degree of $w$. Denote by $f_w$ the face in $\sigma(\pb)$ that corresponds to $w$ in $\pb$. We are allowed to perform spider moves in $T$ and resplits in $\sigma_H(T)$ without affecting the statement of the theorem.  Thus, and due to minimality of $\tcd$, we can assume that we performed a sequence of spider moves in $T$ such that there is no degree four face incident to $w$. Note that 1-genericity is not affected by spider moves. In $\sigma_H(T)$ there is a white vertex $\bar w_k$ for each line $T(w_k)T(w'_k)$ and we also choose white vertices $\bar w'_k$ for each line $T(w'_k)T(w_{k+1})$. Note that we can apply resplits in $\sigma_H(T)$ that change the choice of section because all sections are well defined due to the 1-genericity of $T$. Let us denote the VRC of $T$ by $R$, see \ref{def:vrc}. Assume we choose affine coordinates of $\CP^n$ with $H$ at infinity, then the white vertices $\bar w_k, \bar w_k'$ in $\sigma_H(R)$ are represented by
	\begin{align}
		\sigma_H(R)(\bar w_k)=R(w_k)-R(w'_k) \mbox{ and } \sigma_H(R)(\bar w'_k) =R(w_{k+1}) - R(w'_k),
	\end{align} 
	 for $k=1,2,\dots ,m$. Choose the edge weights of $R$ as in Figure \ref{fig:sectionweights} on the left. They correspond to the relations
	\begin{align}
		R(w) + \lambda_k R(w_k) = (1+\lambda_k) R(w'_k).
	\end{align}
	On the other hand in the section $\sigma_H(R)$ we claim we can locally choose edge weights as in Figure \ref{fig:sectionweights} on the right. The corresponding relations in $\sigma_H(T)$ are
	\begin{align}
		\lambda_k(R(w_k)-R(w'_k)) + (R(w_{k+1})-R(w'_k)) = (1+\lambda_{k+1})(R(w_{k+1})-R(w'_{k+1})),
	\end{align}
	$k=1,2,\dots, m$. The relations in $\sigma_H(R)$ are satisfied because they are equivalent to the difference of two relations in $R$. Now we can read off the edge weights of $\sigma_H(R)$ and $R$, that indeed the multi-ratio around $f_w$, is the star-ratio at $w$ in $T$ and thus the theorem is proven.\qed
}

There is an alternative proof involving $\deg(w)$ copies of Menelaus' configurations. Also note that on the combinatorial level, there is a relation between the extended projective and affine cluster structures too. Let $\qui^\times$ be the extended affine quiver of $T$ and let $\qui^\times_\sigma$ be the extended projective quiver of $\sigma(T)$. Then in fact
\begin{align}
	\qui^\times_\sigma = \Upsilon^+(\qui^\times)
\end{align}
holds, where $\qui^\times_\sigma$ changes the arrows at the boundary, see Definition \ref{def:upsilon}. On the geometric level however, it does not appear to be possible to identify the extended cluster structures unless one fixes some additional data, as the projective cluster variables of boundary faces of $\sigma(T)$ are not invariant under affine transformations of $\CP^n$ with $H$ at infinity.

\begin{lemma}\label{th:sectionsectioncluster}
	Let $\tcd$ be a minimal TCD $T: \tcdp \rightarrow \CP^n,\ n>1$ be a 2-generic TCD map and $E,E'$ be two generic hyperplanes. Consider two sections $\sigma_E(T), \sigma_{E'}(T)$ with TCDs such that $\sigma_E(\tcd) = \sigma_{E'}(\tcd)$. Assume $H=E\cap E'$ is 1-generic with respect to both sections and that $\sigma_{E'}(\sigma_E(\tcd)) = \sigma_{E}(\sigma_{E'}(\tcd))$. If $\dim H > 0$, then
	\begin{align}
		\pro(\sigma_{H}(\sigma_E(T))) = \pro(\sigma_{H}(\sigma_{E'}(T)))
	\end{align}
	holds. Moreover, for $\dim H \geq 0$ holds
	\begin{equation}
		\aff_H(\sigma_E(T)) = \aff_{H}(\sigma_{E'}(T)).\qedhere
	\end{equation}
\end{lemma}
\proof{The first statement is a trivial consequence of Lemma \ref{lem:sectionscommute}, which shows that $\sigma_{E'}(\sigma_E(T)) = \sigma_{E}(\sigma_{E'}(T))$. The second follows from the first one because of Theorem \ref{th:affprojcluster} in all cases except if $n = 2$, which implies that $H$ is a point and we thus cannot identify the affine cluster structures $\aff_H$ with the projective cluster structure $\pro(\sigma_H)$. However, if the maximal dimension of $\sigma(\tcd)$ is larger than 1, we can use Lemma \ref{lem:projectionlift} to construct a lift $\hat T$ of $\sigma(T)$ to a space of dimension more than 1. Because of Theorem \ref{th:projectionandaffine} we can identify the affine cluster structure of $\sigma_H(T)$ with the projective cluster structure of the corresponding section of $\hat T$ and then apply the same arguments as in the other cases. If the maximal dimension of $\sigma(\tcd)$ is 1 then the affine quiver has no arrows and there is nothing to show.\qed
	
}

Under some genericity assumptions, we have now established a way to explain every affine cluster structure via a projective cluster structure, except in $\CP^1$. Of course, this is only true for TCD maps that do not attain maximal dimension. In the maximal dimension case the affine cluster structure of a TCD map in $\CP^1$ is actually trivial. Moreover, we will give a projective interpretation for sections with 1-dimensional subspaces in Section \ref{sec:projclusterduality}, where we investigate cluster structures of the projective dual. Note that the affine cluster structures in the complex line are of particular interest. Indeed, in Chapter \ref{cha:cpone} we show that it is exactly the affine cluster structures in $\C$, which have been used to define embeddings of statistical models in recent research.

\begin{remark}\label{rem:envelopingalgebra}
	Consider a TCD map $T$ defined on a sweepable TCD $\tcd$. Then there is a way to relate projective and affine cluster variables via sweeps, as discussed in Section \ref{sec:tcdextensions}. More precisely, there is a TCD $\mathcal U$ that contains $\tcd$ as a subdiagram and a TCD $\tilde{\mathcal U}$ that contains $\sigma(\tcd)$ as a subdiagram, such that the two TCDs $\mathcal U$ and $\tilde{\mathcal U}$ are related by a sequence of mutations. Additionally, there are also two TCD maps $U,\tilde U$ defined on $\mathcal U, \tilde{\mathcal U}$ such that their respective restrictions are $T$ and $\sigma(T)$. Therefore both the projective cluster variables $X(T)$ as well as the affine cluster variables $Y(T) = X(\sigma(T))$ are expressible as a function of all the projective cluster variables of the TCD map $U$ defined on $\mathcal U$. Of course, if the section $\sigma(\tcd)$ is sweepable again, this process is repeatable. Let us consider a fundamental example, the unique TCD $\tcd^n_{n-1}$ that has endpoint matching $\enm n{-1}$. Then, there is a TCD map $T^n_{n-1}: \tcd \rightarrow \CP^{n-1}$ that represents $n$ points in general position. Now we can add a sweep strand of length $n-1$, and then add a sweep strand of length $n-2$ and so on until the last strand of length 1. Denote the TCD with all the strands added $\mathcal U^\cup$, and denote the $i$-th section of $\tcd^n_{n-1}$ by $\tcd^n_{n-i-1}$. Then for each $i$ with $0\leq i \leq n$ there is a mutation of $\mathcal U^\cup$ such that $\tcd^n_{n-1}$ is a subdiagram of $\mathcal U^\cup$. Therefore all the projective variables that occur in all the sections are functions of the projective variables of the TCD map defined on $\mathcal U^\cup$. One can make counting arguments to show that the number of variables on all of $\mathcal U^\cup$ is $\binom{n-1}{2}$. However, the total number of all variables appearing in all the sections is $\binom{n-1}{3}$. Thus, for $n>5$ there are more variables appearing in the sections than there are variables in the encompassing TCD $\mathcal U^\cup$. Therefore there must be additional relations between the cluster variables appearing in the sections. It would be very interesting to have an understanding of these relations. This remains an open question.
\end{remark}

\section{T-graphs}\label{sec:tgraphcluster}\label{sec:extgraphs}

T-graphs were introduced by Kenyon and Sheffield \cite{kenyonsheffield} to study relations between bipartite planar graphs, planar Markov chains and tilings with convex polygons. In this section we want to show that if we look at certain generic T-graphs, we can capture the incidence relations of T-graphs as TCD maps and also that the known cluster structure on a T-graph \cite{kenyonsheffield} coincides with the affine cluster structure of the corresponding TCD map.

\begin{definition}
	Let $D\subset \R^2$ be homeomorphic to a closed disc. Let $L=(L_1,L_2,\dots,L_m)$ be a collection of disjoint, open line segments such that $\partial D \cup ( \cup_i L_i)$ is connected and closed and such that $\cup_i L_i$ is contained in $D$. We define the \emph{T-graph $\tgra=(V,E,F)$}, where 
	\begin{enumerate}
		\item $V$ is the set of endpoints of line segments,
		\item $E$ contains $(v,v')$ if $v,v'$ appear consecutively on a line segment and
		\item $F$ is the set of discs in $D\setminus (\cup_i L_i)$.
	\end{enumerate}
	Vertices in $\partial D$ are called \emph{boundary vertices}, all other vertices are \emph{interior vertices}. For each interior vertex $v$ denote the segment that contains $v$ in its interior by $L_v$.
\end{definition}

\begin{figure}
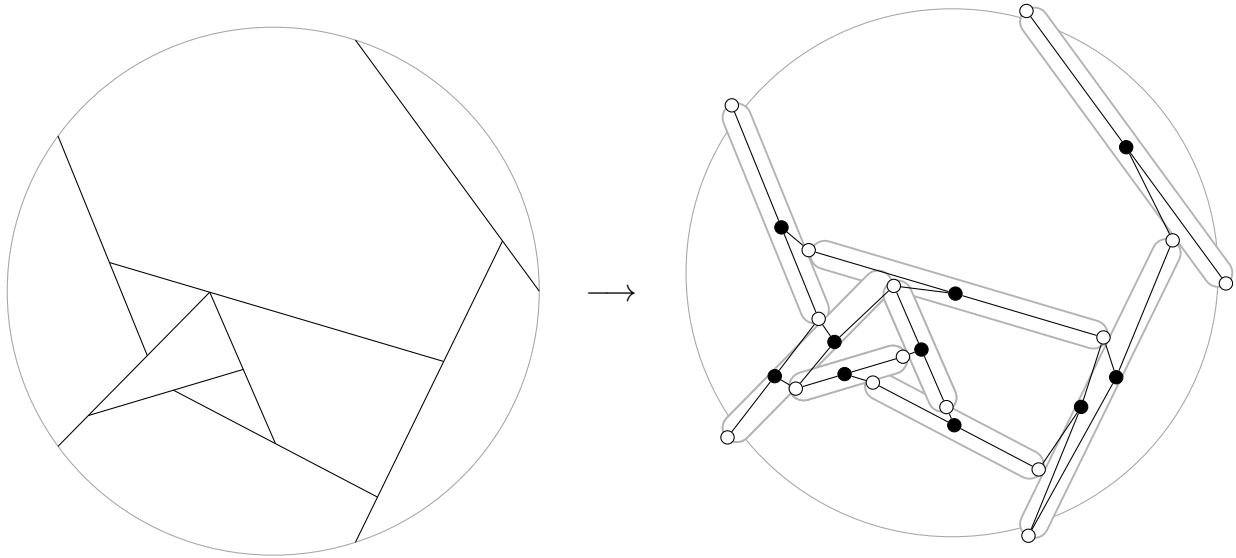

	
	\caption{A one-sided T-graph on the left. Each segment $L_v$ is expanded to a disc $C_v$ (gray), which gives rise to the graph $\pb$ (black and white).}
	\label{fig:tgraphexample}
\end{figure}

If we add a condition to the definition then a T-graph can be represented as a TCD map.

\begin{definition}
	Let $\tgra$ be a T-graph. We call $\tgra$ a \emph{one-sided T-graph} if we can assign a half-space $H_v \subset \R^2$ to each interior vertex $v$ of $\tgra$, such that $L_v \subset H_v$ and for every segment $L_i$ and every interior vertex $v'\in \partial L_i$ holds that $L_i\subset H_{v'}$.
\end{definition}

Thus if we look at a interior vertex $v$ in a one-sided T-graph, the other line segments ending at $v$ have to be on only one side of $L_v$, not both. Also note that $H_v$ is actually uniquely determined for all interior vertices $v$. We can think of one-sidedness as a genericity constraint, because a vertex that is not one-sided can be seen as a limiting case of several one-sided vertices.

Let us explain how to associate a TCD map $T$ to a one-sided T-graph $\tgra$. We begin by removing all line segments from $\tgra$ that do not contain a vertex. If there are interior vertices left that are not on the boundary of any line segment we remove them as well. Then we thicken each remaining line segment $L_v$ such that it becomes a disc $C_v$, see Figure \ref{fig:tgraphexample}. For each vertex of $\tgra$ add a white vertex $w_v$ to $\pb$, such that the new vertex $w_v$ is at the closest point to $v$ that is in $\partial C_v \cap H_v$.  On the boundary of each disc $C_v$ there are now $j_v$ white vertices. We add $j_v - 2$ black vertices inside $C_v$ and add edges, such that we obtain a graph corresponding to a TCD with $\enm {j_v}{1}$ endpoint matching in each disc $C_v$. The TCD map $T$ simply maps every white vertex $w_v$ to $v\in \R^2$. The graph pieces in each disc $C_v$ capture the fact that the corresponding white vertices are on a line in $\R^2$. The condition that $\tgra$ is one-sided guarantees that we can position the white vertices on the boundaries of each disc $C_v$ such that we can actually glue the $\enm {j_v}{1}$ pieces into $C_v$. If $\tgra$ is not one-sided, we may still capture the line-incidence of vertices of $\tgra$ by adding enough black vertices to $\pb$, but $\pb$ would not be a planar graph anymore.

\begin{figure}
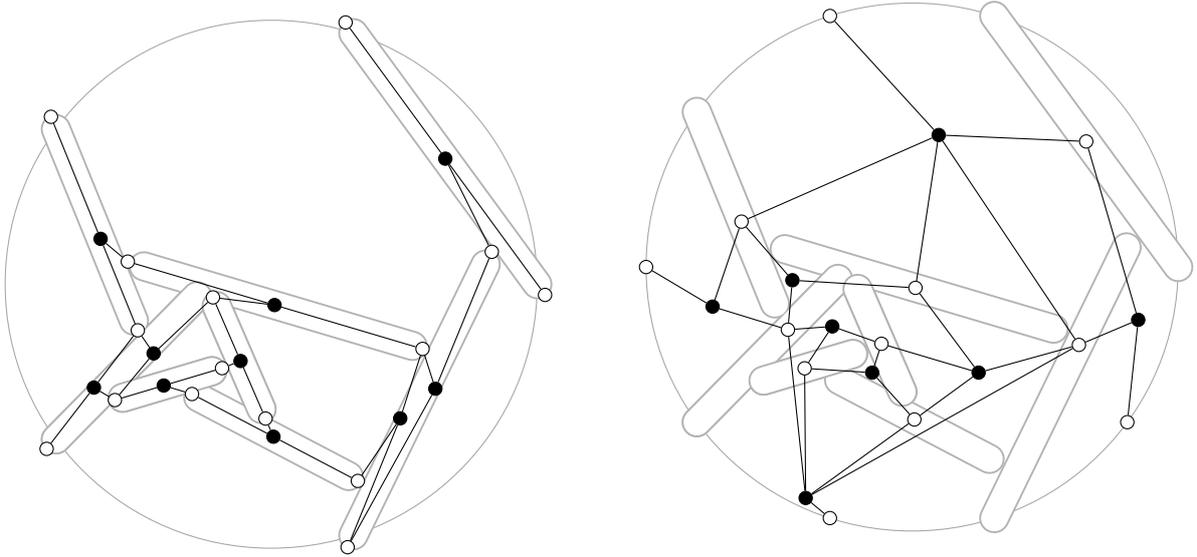

	
	\caption{On the left the VRC of a T-graph (gray). On the right $\tgrad$ that is also the section of that VRC on the left.}
	\label{fig:tgraphs}
\end{figure}

Next we want to compare the bipartite graph $\pb$ to the bipartite graph $\tgrad$ obtained by Kenyon and Sheffield.

\begin{definition}\label{def:tgraphsgd}
	Let $\tgra$ be a T-graph. Define the planar, bipartite graph $\tgrad$ (see also Figure \ref{fig:tgraphs}) by adding a white vertex $w_i$ for every line segment $L_i$ of the T-graph and a black vertex $b_f$ for every face $f$ of $\tgra$. Add an edge $e$ between a $w_i$ and $b_f$ in $\tgrad$ if $L_i$ is incident to $f$ in $\tgra$.
\end{definition}

Note that we swapped black and white in Definition \ref{def:tgraphsgd} in comparison to the definition by Kenyon and Sheffield, in order to be more consistent with our conventions.

\begin{lemma}
	Let $\tgra$ be a one-sided T-graph and let $T$ be the associated TCD map. Then there is a choice of $\sigma(\pb)$ such that, if we contract all 2-valent white vertices in both $\sigma(\pb)$ and $\tgrad$ the two graphs coincide.
\end{lemma}
\proof{
	Let us denote by $\tgra'$ the T-graph $\tgra$ after we removed all line segments that do not contain a vertex. By definition, $\sigma(\pb)$ has a white vertex for every line segment in $\tgra'$. All other white vertices of $\sigma(\pb)$ are two-valent and thus are removed by the contractions. Also by definition $\tgrad$ has a white vertex for every line segment of $\tgra$, and the white vertices for line segments in $\tgra$ but not in $\tgra'$ are removed by contraction. Thus the set of white vertices of $\sigma(\pb)$ and $\tgrad$ coincides. Moreover, after contraction in $\sigma(\pb)$ there is a black vertex for every face in $\tgra'$. On the other hand, $\tgrad$ has a black vertex for every face of $\tgra$. The contractions however join faces in $\tgra$ that are only seperated by line segments not in $\tgra'$. Thus the set of black vertices of $\sigma(\pb)$ and $\tgrad$ also coincide. It is clear that the edges also coincide and thus the claim is proven. \qed
}

\begin{remark}
	Note that there is an interesting phenomenon here. The reason we restrict our T-graph considerations to one-sided T-graphs is that we need $\pb$ to be planar. However, it is intuitive from a geometric point of view that one can define $\sigma(\pb)$ even if $\pb$ is not planar. And the case of T-graphs demonstrates that even if $\pb$ is not planar, $\sigma(\pb)$ can be planar. Thus it may well be that there are other examples of systems where it is not possible to associate a planar graph $\pb$ to the system, but there may exist a planar graph for one of the (possibly iterated) sections of the system.
\end{remark}

Kenyon and Sheffield defined a dimer model associated to a T-graph $\tgra$. Because we only introduce dimer models in Section \ref{sec:dimers}, we translate the definition of Kenyon and Sheffield directly to a cluster structure.

\begin{definition}\label{def:tgraphcluster}
	Let $\tgra$ be a T-graph and $\tgrad$ the graph as in Definition \ref{def:tgraphsgd}. By definition of $\tgrad$, each edge of $\tgrad$ corresponds to a pair $(f_e,L_e)$, where $f_e$ is a face of $\tgra$ and $L_e$ is a line segment of $\tgra$. Associate \emph{edge-weights} $\omega_e$ to the edges of $\tgrad$ such that for each edge $e = (f_e,L_e)$, $\omega_e$ equals the length of the part of the line segment $L_e$ on the boundary of $f_e$ in $\tgra$. The quiver $\qui_{\mathrm T}$ of the \emph{cluster structure $(\qui_{\mathrm T},U)$ of $\tgra$} is defined in the same way as the projective cluster structure is defined for $\pb$ in Definition \ref{def:projclusterstructure}. The variables $U_f$ are defined as the alternating ratios around faces of the edge-weights $\omega$, that is
	\begin{align}
			U_f = \prod_{i=1}^n \frac{\omega(b_i,w_i)}{\omega(w_i,b_{i+1})} ,
	\end{align}
	where $f\in F(\tgrad)$ is the face $(b_1,w_1,b_2,\dots,b_n,w_n)$ in counterclockwise order.
\end{definition}

\begin{theorem}
	Let $\tgra$ be a one-sided T-graph in $\R^2\subset \RP^2$ with the line $H$ at infinity and let $T$ be the associated TCD map. Then the cluster structure $(\qui_{\mathrm T},U)$ of $\tgra$ (see Definition \ref{def:tgraphcluster}) coincides with the affine cluster structure of $T$ with respect to $H$.
\end{theorem}
\proof{Note that because of the definition of $\tgra$, $T$ is 1-generic and $H$ is generic with respect to $T$. We can use Lemma \ref{lem:starratioviadistances} to define the variables $Y_v$ of the affine cluster structure for every vertex $v$ of $\tgra$. It is clear that the occurring distances correspond to the edge-weights $\omega$ of Definition \ref{def:tgraphcluster}. Moreover, we have established in Theorem \ref{th:affprojcluster} how the affine cluster structure of $T$ corresponds to the projective cluster structure of $\sigma_H(T)$. It remains to check that the signs are correct. Clearly all variables $U_f$ are positive by definition. On the other hand, all the star-ratios $Y_v$ are also positive, because the minus sign of the definition of $Y_v$ cancels with the unique minus sign coming from the oriented distance ratio along $L_v$ for every vertex $v$ of $\tgra$.\qed
}

For TCD maps we have discussed how to perform 2-2 moves in previous sections. We do not go into much detail, but 2-2 moves are also possible for T-graphs. Spider moves correspond to reparametrizations of the TCD map $T$ corresponding to a T-graph $\tgra$ and do not affect $\tgrad$. Resplits on the other hand correspond to mutations in the cluster structure of $\tgra$, a fact that has previously been observed by Kenyon and Sheffield and follows as a corollary from our results on the affine cluster structure of $T$. It is not difficult to check that $T$ after a resplit is still a T-graph.

By definition, the affine cluster variables of a TCD map associated to a T-graph only take positive values. It is tempting to think that maybe all TCD maps in $\RP^2$ with positive affine cluster variables are T-graphs, but this is not true. However, if one fixes the combinatorics and the values of the affine cluster variables then it is in many cases possible to find a corresponding T-graph. We begin with such a statement for the case that the corresponding TCD $\tcd$ is minimal, and then we also mention the (stronger) results by Kenyon and Sheffield.

\begin{figure}
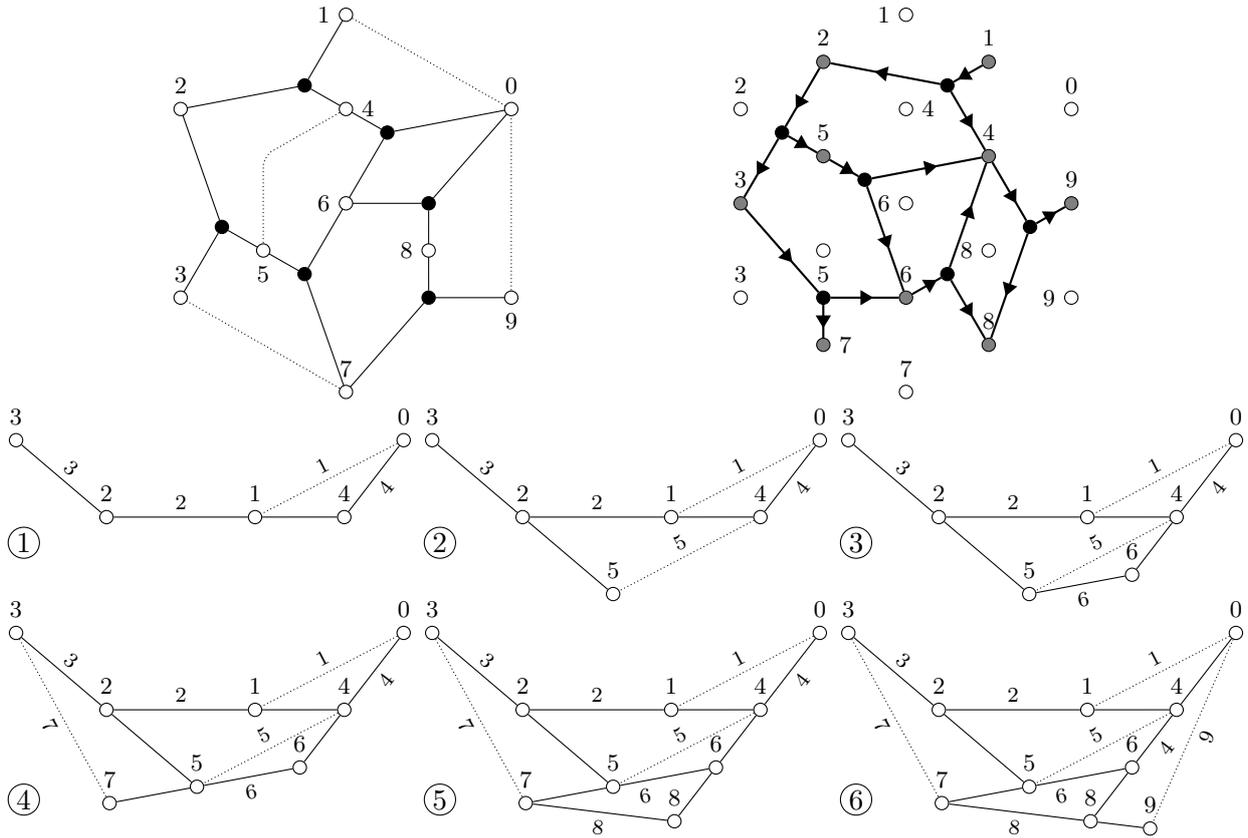


	\caption{A TCD $\tcd$ in the top left together with a distinguished diagonal and boundary lines (dotted). In the top right the li-orientation of the section $\sigma(\tcd)$. In the bottom we construct a T-graph from $\tcd$ step by step.}
	\label{fig:affineconstruction}
\end{figure}

\begin{theorem}\label{th:tgraphconstr}
	Let $\tcd$ be a minimal TCD and let $Y': W(\pb) \rightarrow \R_+$ be a positive function on the white vertices of $\pb$. Then a T-graph $T:\tcd \rightarrow \R^2$ exists such that the affine cluster variables $Y(T)$ satisfy $Y(T) = Y'$.
\end{theorem}
\proof{We give an explicit construction. Consider $\R^2\subset \RP^2$ and let $H$ be the line at infinity. Due to Theorem \ref{th:affprojcluster} we know that $Y_E(T) = X(\sigma_H(T))$. We solve the problem $X(\sigma_H(T)) = Y'$ first via Theorem \ref{th:positiveinvariants}. Then it remains to show that we can construct a T-graph $T$ from $\sigma_H(T)$. See Figure \ref{fig:affineconstruction} for an example. Let us first assume that $\tcd$ is sweepable and that we chose a labeling. Consider the li-orientation $\lio(\sigma(\pb))$ on the section. Each white vertex (gray vertex in the Figure) of the section corresponds to a line and each face to a point of the T-graph. Moreover, in the construction of $\sigma_H(T)$ via Theorem  \ref{th:positiveinvariants} we chose a point at infinity in $H$, which corresponds to an equivalence class of parallel lines in $\R^2$. Without loss of generality choose this direction to be the vertical direction in $\R^2$. The directions of the lines occurring in $T$ adhere to the lr-ordering $\lrp(\sigma(\pb))$. As a consequence, the boundary lines of $T$ above the sweep strand appear with increasing slope. We begin constructing $T$ by choosing the boundary points of $T$ above the sweep line such that consecutive points are on the lines predetermined by $\sigma_H(T)$ and such that they form the boundary of a convex polygon. Because of the aforementioned ordering of the slopes this polygon will be open towards the vertical direction. Now we iteratively place new lines of $T$ by descending through the li-orientation of the section. In each iteration, we choose a white vertex $w$ of $\sigma(\pb)$ that is maximal among those that we have not chosen previously. One of the adjacent faces to $w$ in $\sigma(\pb)$ corresponds to a point of $T$ that has already been placed, and thus the line of $T$ corresponding to $w$ is determined. The other adjacent face $f$ to $w$ does not correspond to an already placed point. However, $f$ has a unique other white vertex $w'$ that was already chosen. Thus the point in $T$ that corresponds to $f$ is now determined as the intersection of the two lines corresponding to $w$ and $w'$. Because of the ordering of the lines, the new point corresponding to $f$ is necessarily below the convex curve that bounds our construction. Therefore in the sweepable case we indeed iteratively construct a T-graph without intersecting segments. As every TCD $\tcd$ is the subdiagram of a sweepable TCD $\tcd'$, we can solve the problem first for the larger sweepable TCD $\tcd'$ and then restrict to $\tcd$.\qed
}

\begin{remark}
	Kenyon and Sheffield prove a stronger theorem than Theorem \ref{th:tgraphconstr}  \cite{kenyonsheffield}, as they show the existence of T-graphs even for certain non-minimal TCDs. They solve the gauge-problem first (as discussed in Section \ref{sec:fromprojinvariants}) and then use the invertability of the Kasteleyn matrix to construct a solution. We think the construction presented here has the benefit of involving only elementary steps and no matrix-calculus. However, in the case of T-graphs excluding non-minimal TCDs is quite a restriction as there are many natural examples of T-graphs that are non-minimal. On the other hand, not every non-minimal TCD can be realized as a T-graph. The non-realizable T-graphs have been characterized via so called $k$-cuts \cite{kenyonsheffield}. It would be interesting to understand the precise relation between $k$-cuts and minimality of TCDs. Moreover, it is an interesting question whether it is possible that the construction algorithm we gave using the li-orientation may be extendable to all T-graphs, by choosing a suitable generalization of the li-orientation. If such a generalization exists, then we expect that the corresponding construction algorithm in the non-minimal case will involve local degrees of freedom that cannot be a priori fixed by choices on the boundary.
\end{remark}

\section{Cluster variables and projective duality of TCD maps} \label{sec:projclusterduality}

In Section \ref{sec:projduality} we have presented a way to describe the geometric and combinatorial relations of TCD maps and projective duality. At this point we understand the associated quivers and can therefore make the observation in the next Lemma. Recall $\rho$ and $\Upsilon^-$ as defined in Definition \ref{def:reciprocalcluster} and Definition \ref{def:upsilon} respectively.

\begin{lemma}\label{lem:dualaffinequivers}
	Let $\tcd$ be a TCD and let $\sigma(\tcd)$ be a section and $\eta(\tcd)$ a line dual of $\tcd$ such that $\iota(\sigma(\tcd)) = \sigma(\eta(\tcd))$, see Lemma \ref{lem:iotaetasigma}. Let $\qui^\times_\sigma$ be the extended affine quiver of $\sigma(\tcd)$ and $\qui^\times_\eta$ the extended affine quiver of $\eta(\tcd)$. Then we claim that $\qui^\times_\eta = \Upsilon^-(\rho(\qui^\times_\sigma))$. Analogously, let $\qui_\sigma$ be the affine quiver of $\sigma(\tcd)$ and $\qui_\eta$ the affine quiver of $\eta(\tcd)$. Then we claim that $\qui_\eta = \rho(\qui_\sigma)$.

\end{lemma}
\proof{
	 We prove the claim by employing the methods of the proof of Lemma \ref{lem:iotaetasigma}, see also Figure \ref{fig:iotaetasigma}. In fact, if one gives the auxiliary graph $\mathbf B$ consistent orientation around faces and cancels opposite arrows, $B$ becomes either $\qui_\sigma$ or $\qui_\eta$ depending on the choice of consistent orientation. The black faces of $\sigma(\tcd)$ (resp. $\eta(\tcd)$) contain a TCD with endpoint matching $\enm{}{-1}$ ($\enm{}{1}$) while the white faces contain a TCD with endpoint matching $\enm{}{1}$ ($\enm{}{-1}$). Therefore due to Definition \ref{def:affinecluster} if we orient the white faces of $\mathbf B$ clockwise we obtain $\qui_\sigma^\times$, but if we instead orient the black faces of $\mathbf B$ clockwise we obtain $\qui_\eta^\times$ which proves the lemma. The claim about the non-extended quivers follows by restriction.\qed
}

\begin{figure}
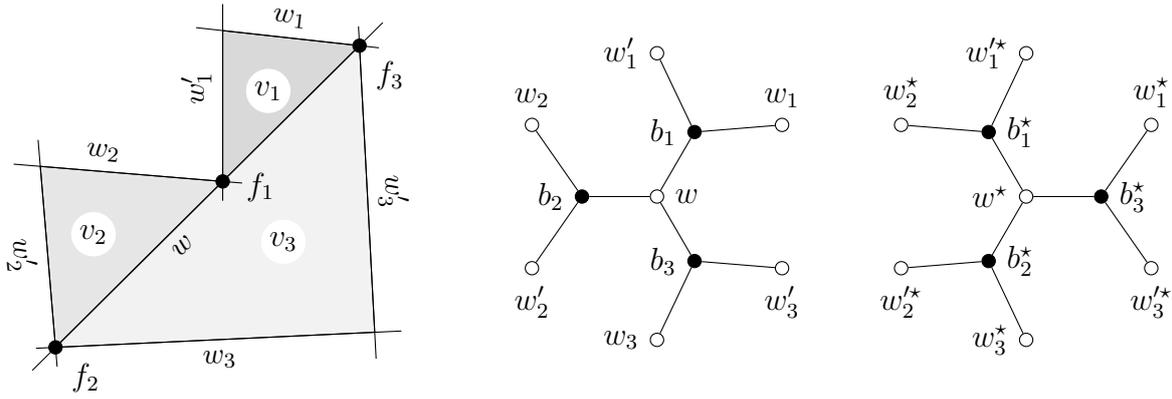


	\caption{An example for $n=m=3$ and $k=2$. On the left a drawing of a geometric configuration of $T_{k+1}$, the planes are labeled by the vertices of $\pb_{k-1}$, the lines by the vertices of $\pb_k$ and the points by vertices of $\pb_{k+1}$. In the center $\pb_k$ at $w$ and on the right $\pb^\star_{\bar k}$ at $w^\star$.}
	\label{fig:pduallabels}
\end{figure}

Now that we have this correspondence of the quivers, it seems natural to ask if there is also a correspondence of the cluster variables. We give an affirmative answer.
\begin{theorem} \label{th:dualcluster}
	Let $(T_k)_{1 \leq k \leq n}$, $(E_k)_{0 \leq k \leq n}$ be a flag of TCD maps and let $(T^\star_k)_{1 \leq k \leq n}$,  $(E^\star_k)_{0 \leq k \leq n}$ be the projective dual flag. Assume $\tcd$ is minimal and that all TCD maps $T_k,T^\star_k$ are flip generic. Then the identity
	\begin{equation}
		\aff_{E_{k-1}}(T_k) = \rho(\aff_{E_{n-k}^\star}(T_{n-k+1}^\star)),
	\end{equation}
	holds for $1 \leq k < n$.
\end{theorem}

\proof{
	For the proof let $\bar k = n-k+1$. On the level of combinatorics, we have that $\tcd_k = \sigma(\tcd_{k+1})$ and $\tcd^\star_{\bar k} = \eta(\tcd_{k+1})$ and therefore by Lemma \ref{lem:dualaffinequivers} the affine quivers of $\tcd_k$ and $\tcd^\star_{\bar k}$ coincide. We fix a white vertex $w$ of the affine quiver of $T_k$ and denote the corresponding white vertex of the affine quiver of $T^\star_{\bar k}$ by $w^\star$. We need to show that $Y_w Y_{w^\star} = 1$. Because we assume that $T_k$ and $T^\star_{\bar k}$ are flip-generic and the TCDs are minimal, we can apply spider moves in $T_k$ and $T^\star_{\bar k}$ until there is no face of degree four adjacent to $w$ and $w^\star$, without changing the affine cluster variables. Because the affine quivers agree, the degree $m$ of $w$ and $w^\star$ is the same. Let us label vertices such that the neighbourhoods are
	\begin{align}
		N(w) &= (b_1,\dots,b_m),& N(w^\star) &= (b^\star_{1},b^\star_{2},\dots,b^\star_{m}),\\
		N(b_i) &= (w,w_i,w'_i), & N(b^\star_{i}) &= (w^\star,w'^\star_i,w^\star_{i+1}),
	\end{align}
	see Figure \ref{fig:pduallabels}. We assume that in this labeling $w_i$ in $\pb_k$ corresponds to $w_i^\star$ in $\pb^\star_{\bar k}$ and accordingly for $w'_i$ and $w'^\star_i$ for $i\in \{1,2,\dots,m\}$. Due to Lemma \ref{lem:srviacr} we can express the star-ratio $Y_w$ via a product of cross-ratios, each cross-ratio being assigned to one of the black vertices $b_1,b_2,\dots,b_m$. Analogously, $Y_{w^\star}$ can be expressed via cross-ratios assigned to the black vertices $b^\star_1,b^\star_2,\dots,b^\star_m$. The plan of the proof is to decompose each cross-ratio into two contributions, and then show that the total contributions in $Y_w$ are inverse to those in $Y_w^\star$. We recall that due to Lemma \ref{lem:subspacedual} the subspace maps are related as 
	\begin{align}
		U_k(w) = (U^\star_{\bar k}(w^\star))^\perp,\quad U_k(w_i) = (U^\star_{\bar k}(w^\star_i))^\perp,\quad U_k(w'^\star_i) = (U^\star_{\bar k}(w'^\star_i))^\perp,
	\end{align}
	for all indices $i$. Let us denote by $v_i$ the white vertex in $\pb_{k-1}$ that replaces $b_i$ when taking the section and analogously $v^\star_i$ for the white vertex of $\pb_{\bar k -1}$ that replaces $b^\star_i$. With the same arguments as in the proof of Lemma \ref{lem:subspacedual}, we observe that
	\begin{align}
		U_{k-1}(v_i) &= \spa\{U_k(w), U_k(w_i), U_k(w'_i) \},\\ 
		U^\star_{\bar k-1}(v^\star_i) &= (U_k(w) \cap U_k(w'_i) \cap  U_k(w_{i+1}))^\perp,
	\end{align}
	for all indices $i$. Via Lemma \ref{lem:srviacr} we obtain
	\begin{align}
		Y_w &= \prod_{i=1}^m \cro(T_k(w),T_k(w_i),T_{k-1}(v_i),T_k(w'_i)),\\
		Y_{w^\star} &= \prod_{i=1}^m \cro(T^\star_{\bar k}(w^\star),T^\star_{\bar k}(w'^\star_i),T^\star_{\bar k-1}(v^\star_i),T^\star_{\bar k}(w^\star_{i+1})).
	\end{align}
	Note that each cross-ratio in the second equation is a cross-ratio of four points on a line in the dual space. In primal space, the four points correspond to four hyperplanes that intersect in a codimension 2 space. It is common in projective geometry to define the cross-ratio of four such hyperplanes via their representatives in dual space. Denote by $H_k(w) = \spa \{U_k(w), E_{k-2}\}$ the ``completion'' of the subspace maps to hyperplanes. Then using Definition \ref{def:dualflag} we can write
	\begin{align}
		Y_{w^\star} &= \prod_{i=1}^m \cro(H_k(w), H_k(w'_i), H_{k+1}(f_i), H_k(w_{i+1})), \label{eq:ydualhyper}
	\end{align}
	where $f_i$ is the white vertex of $\pb_{k+1}$ that when taking the section becomes the face in $\pb_k$ that has $w,w'_i, w_{i+1}$ on its boundary. Next we intend to reexpress $Y_w$ to coincide with the hyperplane formula for $Y_{w^\star}$. Let us choose an affine chart of $\CP^n$ such that $E_k$ is at infinity. Another well known property of the cross-ratio of four hyperplanes as above is that it is also the cross-ratio of the intersections of the four hyperplanes with a generic line. Let us denote by $\dot U_k(w_i)$ the space $U_k(w_i)$ translated such that it contains $U_{k+1}(f_i)$ and let $\dot H_k(w_i) = \spa \{\dot U_k(w_i),E_{k-2} \}$. Then
	\begin{align}
		Y_w &= \prod_{i=1}^m \cro(U_k(w), \dot U_k(w_i), \spa\{T_{k-1}(v_i) , U_{k+1}(f_i)\}, U_k(w'_i)),
	\end{align}
	because the four arguments are four hyperplanes as subspaces of $U_{k-1}(v_1)$ that all contain $U_{k+1}(f_i)$ and have intersections $T_k(w),T_k(w_i),T_{k-1}(v_i),T_k(w'_i)$ with the line $L_k(b_i)$, where $L_k$ is the line map (see Definition \ref{def:linemap}) of $T_k$. By considering the spans with $E_{k-2}$, we obtain that
	\begin{align}
		Y_w &= \prod_{i=1}^m \cro(H_k(w), \dot H_k(w_i),\spa\{ T_{k-1}(v_i) , E_{k-2} , U_{k+1}(f_1)\}, H_k(w'_i)).
	\end{align}
	On the other hand, $\spa\{T_{k-1}(v_i) , E_{k-2}\} = E_{k-1}$ and therefore
	\begin{align}
		Y_w &= \prod_{i=1}^m \cro(H_k(w), \dot H_k(w_i), H_{k+1}(f_1), H_k(w'_i)).
	\end{align}
	This expression almost coincides with the expression in Equation \eqref{eq:ydualhyper}, except that the ordering is different and that $\dot H_k(w_i)$ appears instead of $H_k(w_i)$. However, in the affine chart that we chose, $\dot H_k(w_i)$ and $H_k(w_i)$ belong to the same class of parallel hyperplanes, and the choice of representative does not matter for the cross-ratios. Moreover, the different ordering of the terms explains why $Y_w Y_{w^\star} = 1$ and not $Y_w = Y_{w^\star}$. Therefore the claim is proven. \qed

}

Of course, whenever both sides in the equation of Theorem \ref{th:dualcluster} do not involve a zero dimensional, the theorem can be restated for projective cluster structures as well.

\begin{corollary} \label{cor:dualcluster}
	Let $(T_k)_{1 \leq k \leq n}$, $(E_k)_{0 \leq k \leq n}$ be a flag of TCD maps and let $(T^\star_k)_{1 \leq k \leq n}$,  $(E^\star_k)_{0 \leq k \leq n}$ be the projective dual flag. Assume $\tcd$ is minimal and that all TCD maps $T_k,T^\star_k$ are flip generic. Then the identity
	\begin{equation}
		\pro(T_{k}) = \rho(\pro(T^\star_{n-k-1})),
	\end{equation}
	holds for $1 \leq k < n$.
\end{corollary}
\proof{Direct consequence of Theorem \ref{th:dualcluster} and Theorem \ref{th:affprojcluster}.\qed}

Moreover, Theorem \ref{th:dualcluster} does relate the ``zero-dimensional'' affine cluster structure to a projective cluster structure via
\begin{align}
	\aff_{E_0}(T_1) = \rho(\pro(T^\star_{n-1})).
\end{align}
Geometrically, this is saying that given the point $E_0 \in \CP^n$ we may consider $(n-2)$-spaces occurring in $T_n$ and extend them to $(n-1)$-spaces via $E_0$. The affine variables of $T_1$ then correspond to cross-ratios of these $(n-2)$-spaces. Note that in total from a TCD map we now obtain $(n+2)$ different cluster structures. The primal gives $(n+1)$ structures and the projective dual adds exactly one, namely $\pro(T_n^\star)$.

\begin{figure}
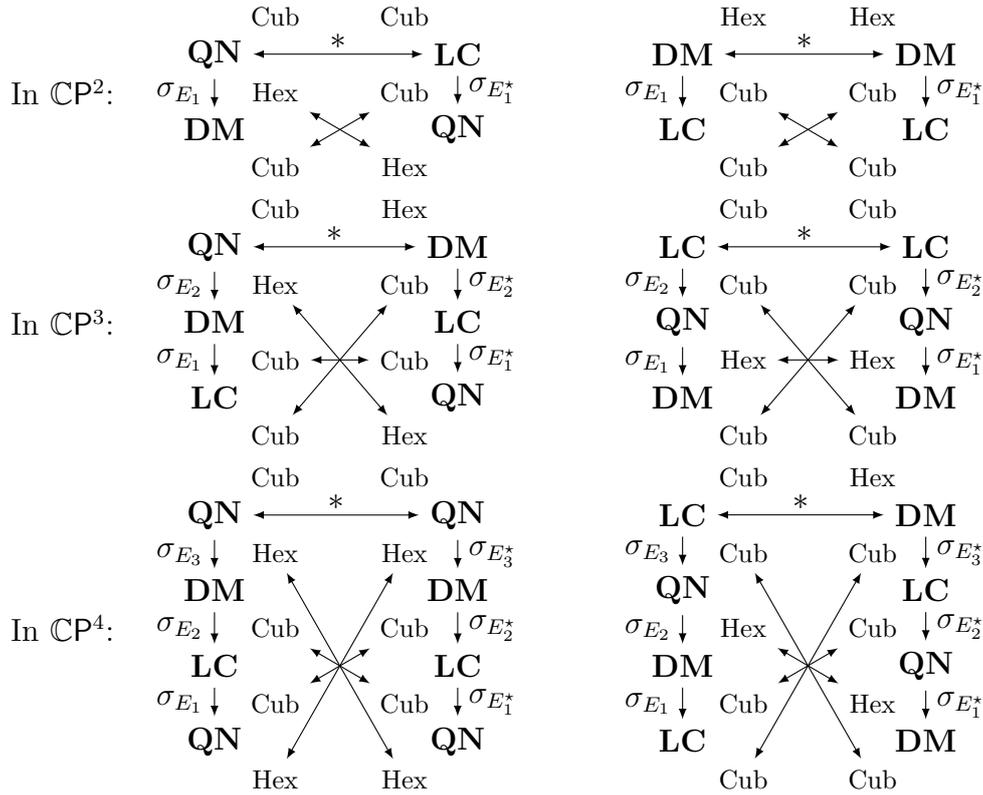


	\caption{We are considering $\Z^3$ combinatorics. The superscript of each map is its projective cluster structure while the subscript is its affine cluster structure. Two cluster structures that are connected by a diagonal arrow are reciprocal. We consider the two TCD maps in the top row to be the projective duals of each other.}
	\label{fig:qnclusterdualrelations}
\end{figure}

\begin{example}
	Let us consider Q-nets, Darboux maps and line complexes with $\Z^3$ combinatorics. If we represent the Cauchy data of these maps as stepped surface, then Example \ref{ex:ddgsteppedsections} shows that we can consider flags of TCD maps that consist solely of Q-nets, Darboux maps and line complexes. As a consequence of Example \ref{ex:tcdprojdual}, the projective dual flag then also consists of these three maps. The cluster structures of these flags therefore are either a cuboctahedral or a hexahedral cluster structure. Theorem \ref{th:dualcluster} then implies identifications between these cluster structures, see Figure \ref{fig:qnclusterdualrelations}. This will become of particular interest later on, as we will show that the resistor and Ising subvarieties of the cluster variables are in correspondence with certain geometric reductions of Q-nets, Darboux maps and line complexes. Theorem \ref{th:dualcluster} then immediately implies that the same reductions occur in the dual flag.
\end{example}

Looking at Figure \ref{fig:qnclusterdualrelations} it may appear a bit odd that all cluster structures in the primal flag have an equivalent in the dual flag, except for $\pro(T_n)$. We claim this is an artifact of our viewpoint. First of all, our examples are mostly from DDG. In these examples usually the dimension of the target space of the TCD maps is much smaller than the maximal dimension. Also, we studied taking a section $\sigma_E(T)$ extensively, but did not introduce or investigate a possible $\sigma_?^{-1}(T)$ operation, even though for example the results of Section \ref{sec:sweeps} imply that it is certainly possible on a combinatorial level and in Section \ref{sec:tcdextensions} we also discuss some geometric aspects. Let us therefore also consider an example that is ``maximal'' in both a geometric and combinatorial sense.

\begin{example}
	Consider a flag of TCD maps $(T_k)_{0 \leq k \leq n}$, $(E_k)_{0 \leq k \leq n}$ such that the endpoint matching of $\tcd_k$ is $\enm n{k-1}$ and such that $T_k$ extends its maximal dimension for $0\leq k \leq n$, see Figure \ref{fig:balanceddualflags} for an example of the combinatorics. Note that we included the trivial TCD map in dimension 0 that maps one white vertex to the one point in $\CP^0$. Then actually $\pro(T_n), \pro(T_{n-1})$ are trivial in the sense that the quivers have no vertices, because we recall that we only consider interior faces of of the bipartite graphs as quivers of the projective cluster structure. Analogously, $\aff(T_1), \aff(T_{0})$ are trivial as well. Therefore the identities $\aff_{E_{k-1}}(T_k) = \rho(\aff_{E_{n-k}^\star}(T_{n-k+1}^\star))$ and $\pro(T_{k}) = \rho(\pro(T^\star_{n-k-1}))$ hold for all indices such that $0 \leq k \leq n$.
\end{example}

\section{The perfect dual of a TCD map}\label{sec:perfdual}

We have previously encountered the TCD dual $\iota$ that reverses the orientations of all strands of a triple crossing diagram. We also noticed that the affine quiver of a TCD $\tcd$ is exactly the projective quiver of $\iota(\tcd)$ and vice versa. Thus one may wonder whether there are two TCD maps that have interchanged projective and affine cluster structures.

\begin{definition}
	Let $T,T'$ be two TCD maps defined on $\tcdp$ resp. $(\iota(\tcd))_{\bm{\circlearrowright}}$ to $\CP^n$ and let $E$ be a hyperplane that is generic with respect to both $T$ and $T'$. We say $T'$ is the \emph{perfect dual} of $T$ if
	\begin{equation}
		\pro(T) = \rho(\aff_E(T')) \quad\mbox{and}\quad \pro(T') = \rho(\aff_E(T))\label{eq:reverseclusters}.\qedhere
	\end{equation}	
\end{definition}

In the previous section we observed that $\iota$ plays a role in defining dual projective TCD maps. If we look at the $\CP^2$ example in Figure \ref{fig:qnclusterdualrelations}, we realize that the sections living in $E_1$ actually satisfy Equation \eqref{eq:reverseclusters}! The same phenomenon is also visible in $\CP^4$, in this case the sections in $E_2$ satisfy Equation \eqref{eq:reverseclusters}.

\begin{theorem}\label{th:perfdualinevenflags}
	Let $n=2m,\ m \in \N$ and let $(T_k)_{1\leq k \leq n}$ be a flag of TCD maps and $(T^*_k)_{1\leq k \leq n}$ its dual flag. Then $T_m$ is the perfect dual of $T^*_m$.
\end{theorem}
\proof{Immediate consequence of the definition of the perfect dual and Theorem \ref{th:dualcluster}. \qed}

Moreover, assume we are given a TCD map $T: \tcdp \rightarrow E_m = \CP^m$. Assume we can lift it to a TCD map $\hat T: \hat\tcdp \rightarrow \CP^{2m}$ such that $\sigma_{E_m}(\hat T_{2m}) = T$. Then we may extend $\hat T$ to a flag of TCD maps that involves $T$ and can construct the corresponding projective dual flag $\hat T^*$. Then Theorem \ref{th:perfdualinevenflags} states that $\sigma_{E_n^*}(\hat T^*)$ is actually a perfect dual of $T$. Thus a reasonable question is under which circumstances these lifts exist. We have discussed how to construct lifts in Section \ref{sec:tcdextensions} on extensions. The prerequisite in that case is the sweepability of the TCDs. We claim without proof that for Q-nets, Darboux maps and line complexes on rectangular subgraphs of $\Z^2$ these strands always exist and therefore so do the lifts and therefore so do perfect duals. For more general combinatorics, the answer to this question is unclear though. 

For the reader that is particularly interested in the perfect dual, we mention that in Theorem \ref{th:orthoharmperfdual} (on h-embeddings and orthodiagonal maps) another instance of the perfect dual occurs where $T' = \sigma(T)$.

In Section \ref{sec:pentagram} we show how the pentagram map is a TCD map on a torus. It is an interesting and in general open question whether there exist polygons that have a perfect dual when considered as TCD map. However, in the case of the pentagram map it is not possible to construct a lift to $\RP^4$ such that its section is the initial polygon again. Indeed, it seems not every polygon admits such a perfect dual. A short remark on a relation between the perfect dual and Glick's operator exists in work by Izosimov \cite{izosimov}.

\section{Focal nets and Laplace transforms}\label{sec:focalnets}
We have already looked at transformations of nets that involve the focal points in terms of Laplace-Darboux dynamics of $\Z^2$ Q-nets, see Section \ref{sec:laplacedarboux}. However, Laplace transformations also exist in the case of $\Z^3$ Q-nets. An introduction is given in the DDG book \cite{ddgbook}. In that context there is also the notion of a focal line complex of a Q-net. We now introduce a version of focal nets and Laplace transforms that treats Q-nets, Darboux maps and line complexes defined on $\Z^3$ on an equal footing. Because we work on $\Z^3$, we identify line complexes and line compounds, but we work with the definition of line compounds. That is, we consider the points of the line complex to be on the faces of $\Z^3$. Also note that we do not perform an analysis of necessary genericity conditions for focal and Laplace transforms, instead we generally assume that the maps involved are generic enough for the transforms to exist.

\begin{figure}
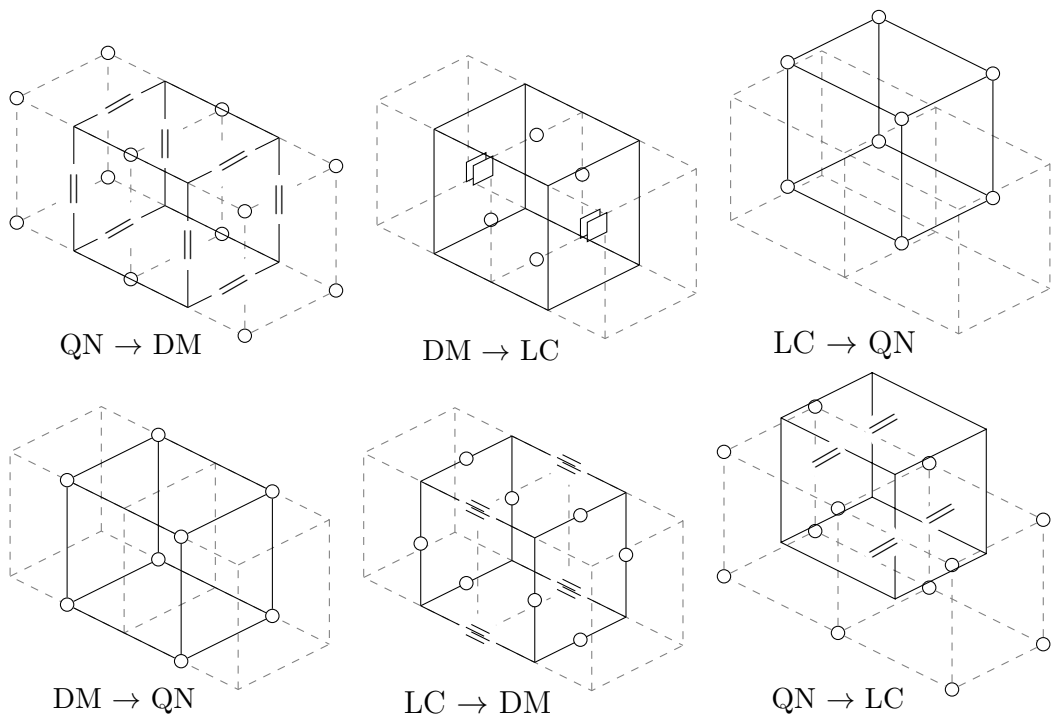

	
	\caption{The six cases of a focal net (solid) of a net of another type (dashed). White points represent points that coincide with points of the original map. Two (resp. 4) parallel lines represent the intersection of two (resp. 4) lines of the same type in $\Z^3$. Two parallel planes represent the intersection of two lines associated to two quads of the same type in $\Z^3$.}
	\label{fig:comboffocal}
\end{figure}

\begin{definition}\label{def:focalnet}
For each type of map -- Q-net, Darboux map and line complex -- we define a \emph{focal net} of each other type of map. Each focal net $f_i$ is taken in a specified direction $i\in\{1,2,3\}$. Positive direction $i$ versus negative direction $\bar i$ determines the type of focal net that we obtain.
	\begin{enumerate}
	\item The \emph{focal line complex $f_{\bar i}$ of a Q-net}  $V(\Z^3) \rightarrow \CP^n$ in direction $i$ consists of all lines associated to edges $(v,v_i)\in E(\Z^3)$. The intersection points of the focal line complex live naturally on the faces of $\Z^3 + \f12 (e_j+e_k)$ (see Figure \ref{fig:comboffocal}) and consist of all the points of the Q-net and all focal points $(v,v_i)\cap (v_j,v_{ij})$ and $(v,v_i)\cap (v_k,v_{ik})$.
	\item The \emph{focal Darboux map $f_{i}$ of a Q-net} $V(\Z^3) \rightarrow \CP^n$ in direction $i$ consists of all planes associated to faces $(v,v_j,v_k,v_{jk})\in F(\Z^3)$. The intersection points of this Darboux map live naturally on the edges of $\Z^3 + \f12 e_i$ and consist of all the points of the Q-net and all focal points $(v,v_j)\cap (v_i,v_{ij})$ and $(v,v_k)\cap (v_i,v_{ik})$.
	\item The \emph{focal Q-net $f_{\bar i}$ of a Darboux map} $E(\Z^3) \rightarrow \CP^n$ in direction $i$ consists of all points associated to edges $(v,v_i)\in E(\Z^3)$. The intersection points of the focal Q-net live naturally on the edges of $\Z^3 + \f12 e_i$.
	\item The \emph{focal line complex $f_{i}$ of a Darboux map} $E(\Z^3) \rightarrow \CP^n$ in direction $i$ consists of all lines associated to faces $(v,v_j,v_k,v_{jk})\in F(\Z^3)$. The intersection points of the focal congruence live naturally on the faces of $\Z^3 + \f12 (e_j+e_k)$. The intersection points consist of all the points of the Darboux map that live on edges of type $(v,v_j)$ and $(v,v_k)$ plus the intersections $(v,v_j,v_k,v_{jk})\cap(v_i,v_{ij},v_{ik},v_{ijk})$ of the lines associated to the faces of the Darboux map.
	\item The \emph{focal Q-net $f_{i}$ of a line complex} with intersection points $F(\Z^3) \rightarrow \CP^n$ in direction $i$ consists of all points associated to faces $(v,v_j,v_k,v_{jk})\in F(\Z^3)$. The points of the Q-net live naturally on the faces of $\Z^3 + \f12 (e_j+e_k)$. They are the intersection points of lines $\ell \cap \ell_i$. The planarity of the resulting quads in $\CP^3$ is a consequence of the coplanarity property of line complexes \cite{bobenkoschieflinecomplexes}. 
	\item The \emph{focal Darboux map $f_{\bar i}$ of a line complex} with intersection points $F(\Z^3) \rightarrow \CP^n$ in direction $i$ consists of all planes spanned by the two lines incident to $(v,v_j,v_k,v_{jk})\in F(\Z^3)$. The points of the Darboux map live naturally on the faces of $\Z^3 + \f12 e_i$ if we identify them with edges of ${\Z^3}^*$. The points consist of the intersection points of lines $\ell\cap \ell_j$ and $\ell\cap \ell_k$ as well as the points $f^{ij}f^{ij}_i \cap f^{ik}f^{ik}_i$.	\qedhere		
	\end{enumerate}
\end{definition}

We introduced the notation $f_i$ because we claim that taking focal nets of focal nets and iterations can be described by a commutative diagram that has the combinatorics of a quotient of $\Z^3$.

\begin{theorem}\label{th:focalrelations}
	For focal nets the equations
	\begin{align}
		f_{\bar i} \circ f_i &= \mbox{id},\\
		f_j \circ f_i &= f_i \circ f_j,\\
		f_3 \circ f_2 \circ f_1 &= \mbox{id},
	\end{align}
	hold.
\end{theorem}
\proof{These equations can be verified in each case by the definition of focal nets. \qed}

\begin{figure}
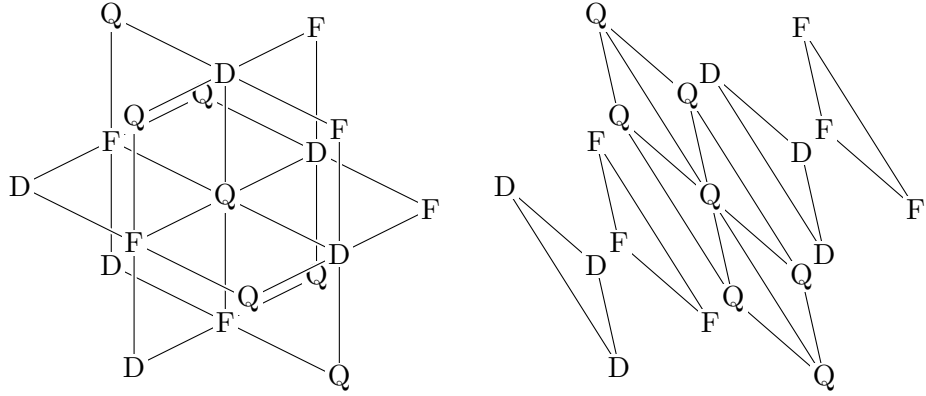

	
	\caption{Left: Part of the commuting diagram of focal nets. Right: The same part but with the edges that correspond to Laplace transforms. }
	\label{fig:focalmeta}
\end{figure}

This theorem essentially states that we can represent the commutative diagram of focal nets as a quotient of the $\Z^3$ lattice, see Figure \ref{fig:focalmeta}.

\begin{corollary}
	Consider the lattice $\lat=\Z^3/(1,1,1)$, that is $\Z^3$ modulo integer translations in direction $(1,1,1)$. We can assign
	\begin{itemize}
		\item a Q-net to each point $n=(n_1,n_2,n_3)\in \Z^3$ with $n_1+n_2+n_3 \in 3\Z$,
		\item a Darboux map to each point $n\in \Z^3$ with $n_1+n_2+n_3 \in 3\Z+1$,
		\item a line complex to each point $n\in \Z^3$ with $n_1+n_2+n_3 \in 3\Z+2$
	\end{itemize}
	such that for any $n\in L$ the maps $h$ at $n$ and $h'$ at $n+\mathbf{e}_k$ are related by 
	\begin{equation}
		h'= f_k(h)\quad \mbox{and}\quad h=f_{\bar k}(h').\qedhere \label{eq:focalgenerators}
	\end{equation}
\end{corollary}
\proof{We begin by assigning a Q-net to the origin of $\lat$. Then we assign maps to the remainder of $\lat$ via Equation \eqref{eq:focalgenerators}. This is consistent along one edge because of the first equation of Theorem \ref{th:focalrelations}. The consistency around a quad follows from the second equation of the same theorem. Together we have consistency on $\Z^3$. Lastly, the consistency also on the quotient $\lat$ follows from the third equation of Theorem \ref{th:focalrelations}.\qed
}

Thus even though we can generate focal nets in three directions, the generated lattice $\lat$ is just a 3-layer thick 2-dimensional lattice. Moreover, we can also define a general Laplace transform as follows.
\begin{definition}
	We define the \emph{Laplace transform} $\Delta_{ij} = f_j \circ f_{\bar i}$ for $i\neq j$.
\end{definition}
For Q-nets this definition coincides with the classical definition of Laplace transforms. Indeed, $f_{\bar i}$ translates to taking all the lines in one coordinate direction and $f_j \circ f_{\bar i}$ translates to intersecting lines of $f_{\bar i}$ adjacent in the other coordinate direction, which is the classical definition of the Laplace transform. The fact that $\Delta_{ij}\circ \Delta_{ji} = \mbox{id}$ (see Equation \eqref{eq:laplacetrafoinverse}) thus is simply a consequence of the fact that $\Delta_{ij}$ corresponds to going two steps on $\lat$ while $\Delta_{ji}$ corresponds to going those steps backwards. We observe that $\lat$ decomposes into three $A_2$ lattices $\lat_r = \{n\in \lat : n_1+n_2+n_3=r\}$. The Laplace transforms always take a map in $\lat_r$ to another map in $\lat_r$. Clearly two maps in $\lat_r$ are related by a sequence of Laplace transforms as the vectors $\{\mathbf{e_1}-\mathbf{e_2},\mathbf{e_2}-\mathbf{e_3}\}$ span $A_2$, see Figure \ref{fig:focalmeta}. The $A_2$ lattices $\lat_r$ are therefore the three disjoint commuting diagrams of the Laplace transforms for either Q-nets, Darboux maps or line complexes.

Now that we are considering Q-nets, Darboux maps and line complexes on equal footing with regards to focal nets and Laplace transforms as well, it is not so surprising that we can in fact capture focal nets via moves in TCD maps again.

\begin{figure}
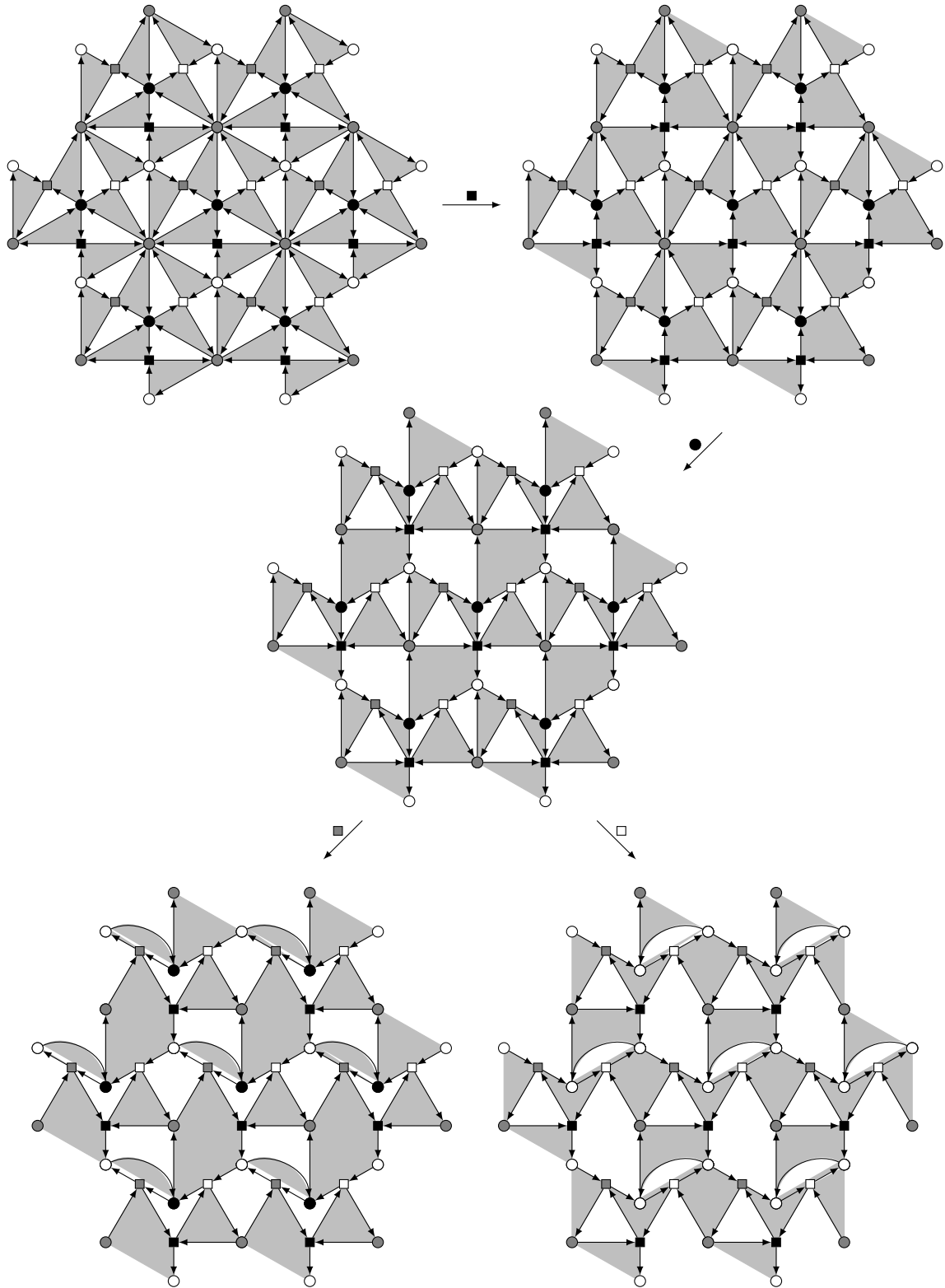


	
	\caption{The series of mutations that takes the affine quiver of a stepped surface of a $\Z^3$ Q-net $q$ (top left) to a line complex $f_{\bar i}(q)$ (bottom left) or Darboux map $f_i(q)$ (bottom right). Focal points in $q$ are drawn as squares, vertices as circles. The points on the boundary of gray regions are colinear.}
	\label{fig:focalmutations}
\end{figure}

\begin{theorem}\label{th:focalmutations}
	Let $h'=f_i(h)$ be two maps related by the focal transform $f_i$. Consider stepped surface Cauchy data for $h$ and $h'$, and let $T$ (resp. $T'$) be the two TCD maps that correspond to $h$ ($h'$). Then $T$ is related to $T'$ by a sequence of 2-2 moves.
\end{theorem}
\proof{We give the global sequences of 2-2 moves that take a Q-net to a Darboux map or line complex in Figure \ref{fig:focalmutations} via their affine quivers. Assume that $q$ is a Q-net and the TCD map $T$ includes as points the focal points in directions $(k,k+1)$ for cyclical indices $k\in\Z_3$. The sequence begins with mutating all focal points of type $(i-1,i)$ followed by mutation at all vertices of $\Z^3$ that have degree 3 in the stepped surface. If we want to obtain $f_i(q)$ then the third step is to mutate at all focal points $\fp^{i,i+1}$. If instead we want to obtain $f_{\bar i}(q)$ then the third step is to mutate at all focal points $\fp^{i+1,i-1}$. We leave it as an exercise to verify that the resulting points are exactly the defining points of the focal nets as in Definition \ref{def:focalnet}. The sequences producing focal transforms of Darboux maps and line complexes are obtained by composing the steps of the sequence given above in different orders.\qed}

Forgetting about geometry for a moment, an interesting consequence of Theorem \ref{th:focalmutations} is that it is possible to take a cuboctahedral quiver into a hexahedral quiver and vice versa via a sequence of mutations.

\chapter{$Z$-invariance and TCD maps}\label{cha:subvar} 
In this chapter we relate PDB quivers (see Definition \ref{def:planquiver}) to the \emph{dimer model} known from statistical physics. To this end we reiterate some of the relevant basics and results from the literature. In particular, we show how mutations in quivers correspond to $Z$-invariant moves in the dimer model. We then reiterate how the \emph{spanning tree model} \cite{kennelly,temperleybijection, gkdimers} as well as the \emph{Ising model} \cite{lenzising, dubedattrick, kenyonpemantle} can be identified with dimer models that have particular combinatorics and cluster variables in the \emph{resistor subvariety} and the \emph{Ising subvariety} respectively. In these subvarieties the $Z$-invariant moves reproduce known discrete integrable equations: The \emph{dKP, dBKP and dCKP equations}. There are many instances of objects of interest from DDG or discrete integrable systems that exhibit these equations as well. These objects turn out to be various reductions of Q-nets, Darboux maps or line complexes and in the literature the integrable equations are discovered by finding the right sequence of calculations and transformations on a case to case basis. The idea of associating discrete integrable equations to cluster subvarities will allow us instead to find the integrable equations in a \emph{canonical way} in these objects. In fact we can make precise statements and say whether the BKP resp. CKP structure belongs to the projective cluster structure of the object when considered as TCD map, to the projective cluster structure of the projective dual of the object, or to a section with some hyperplane of the object. When we restrict ourselves to the projective cluster structures, we show that in Q-nets we can recover the BKP equation precisely if the Q-net is a \emph{K{\oe}nigs net} \cite{bsmoutard, doliwatnets}. In Darboux maps we recover the CKP equation precisely in the case of \emph{Carnot maps}, which are maps discovered by Schief who called them \emph{CKP Darboux maps}. Moreover, for line complexes we find the BKP equation for \emph{Doliwa complexes}, which we introduce specifically for this purpose. However, Doliwa complexes are closely related to another definition of K{\oe}nigs nets due to Doliwa \cite{doliwakoenigs}, as we explain in Section \ref{sec:Doliwacomplex}. We also look at the cluster structures of sections of maps. For example, we show that \emph{Schief maps} \cite{schieflattice}, introduced as BKP maps, are a special case of affine Darboux maps. We also show that \emph{CQ-nets}, introduced by Doliwa as C-quadrilateral lattices \cite{doliwacqnet} are exactly the affine CKP Q-nets. Chapter \ref{cha:bilinear} will focus on special maps associated to bilinear forms, and these in turn can be related to many more canonical instances of the BKP and CKP equation, but we leave a discussion of these maps and relations to that chapter.

\section{$Z$-invariance in Statistical Physics}

In this section we introduce basic concepts that appear while studying discrete statistical mechanics models, but we have three particular models in mind: The dimer model, the spanning tree model and the Ising model.

Let $G$ be a finite planar graph embedded in a disk, let $E$ be its set of edges and let $w: E\rightarrow \R^+$ denote the \emph{edge weights}. The weight of some subset $R\subset E$ is then defined to be
\begin{align}
	w(R) = \prod_{e\in R} w(e).
\end{align}
Since we are interested in the statistics for arbitrary edge weights, we will generally consider the edge weights $w_e$ as formal variables. In order to obtain interesting statistics, we choose a set $\Omega \subset 2^E$ of subsets of edges that are the allowed configurations. The choice of $\Omega$ defines the statistical model. 
\begin{definition}
	Let $\Omega \subset 2^E$ be the \emph{sample space}. The \emph{partition function $Z$} is defined as the formal polynomial
	\begin{align}
		Z: \R^E\rightarrow \R,\quad Z = \sum_{R\in \Omega} w(R),
	\end{align}
	and the \emph{probability} of an element $R\in \Omega$ is defined as
	\begin{equation}
		\mathbb P_Z(R) = \frac{w(R)}{Z}\ . \label{eq:omegaprob}\qedhere
	\end{equation}
\end{definition}

The partition function $Z$ has the useful property that we can calculate the probability of an edge occurring via taking derivatives, that is
\begin{align}
	\mathbb P_Z(e\in R) = w_e\frac{\del\log Z}{\del w_e}\ .
\end{align}
Many other observables can be expressed via derivatives of the partition function. Among these are for example certain correlations like the probability that two fixed edges both appear. 

In the cases we consider in the thesis, the statistical models are not just given for one particular graph $G$, but for a whole class of graphs. To that end, the models are defined by definitions of the sample space $\Omega$ depending on the graph $G$. 

Before focusing on particular models, we discuss the key concept of \emph{$Z$-invariance}. This concept relates statistical physics models to discrete integrable systems.  

\begin{definition} Let $L$ and $\tilde L$ be finite, planar graphs together with an identification of the boundary vertices $\del L \leftrightarrow \del \tilde L$. Let $G$ be a planar graph such that $L$ is a subgraph and let $L^c=G\setminus L$. 
	\begin{enumerate}
	\item We can replace $L$ with $\tilde L$ in $G$ via the identification along $\del L$ to obtain the new graph $\tilde G$. We call this replacement a \emph{local change of combinatorics}.
	\item Let $f: \R^{E(L)} \rightarrow \R^{E(\tilde L)}$ be a function such that
	\begin{align}
		(f(\omega))_e &= \omega_e && \mbox{for all } e\in E(L^c),\\
		\del_{\omega_{e'}} (f(\omega))_e &= 0 && \mbox{for all } e\in E(\tilde L), e' \in E(L^c).
	\end{align}
	Then $f(\omega)$ defines new edge weights $\tilde \omega$ on $\tilde G$. The triple $(L,\tilde L, f)$ defines a \emph{local move}.
	\item We say a local move $(L,\tilde L, f)$ exhibits \emph{$Z$-invariance} if
	\begin{align}
		Z_{\tilde G} &= \lambda Z_G	
	\end{align}
	for some  $\lambda \in \R$. \qedhere
	\end{enumerate}
\end{definition}

We call the moves local because the function $f$ does not change edge weights outside of these subgraphs and the new edge weights in $\tilde L$ only depend on the edge weights of the old subgraph $L$. An important property of a local move is that it preserves the edge probabilities , that is
\begin{align}
	\mathbb P_{Z_G}(e\in R) = \mathbb P_{Z_{\tilde G}}(e\in R)
\end{align}
holds for any for $e \in E(L^c)$. In general, given the rule that determines $\Omega_G$, we are interested in finding $Z$-invariant local moves. This means that we are looking for two graphs $L$ and $\tilde L$ such that there is an accompanying function $f$ that makes the triple $(L,\tilde L,f)$ $Z$-invariant for all graphs $G$. If we fix $L$ and $\tilde L$ then $f$ has to satisfy a set of algebraic equations. In particular, we obtain an equation for every $R\in \Omega_G$. If we write $R^e = R\cap E(L^c)$, then the condition is that  
\begin{align}
	\sum_{R' \in \Omega_G : R' \cap E(L^c) = R^e} w(R') = \lambda\sum_{\tilde R' \in \Omega_{\tilde G} :  \tilde R' \cap E(L^c) = R^e} w(\tilde R')
\end{align}
for all $R\in \Omega_G$ and a $\lambda \in \R$. On the left and the right side the exact same edge weights on $L^c$ appear in every term, so that we can rewrite the equations as
\begin{align}
	\sum_{R' \in \Omega_G : R' \cap E(L^c) = R^e} w(R' \cap E(L)) = \lambda \sum_{\tilde R' \in \Omega_{\tilde G}:  \tilde R' \cap E(L^c) = R^e} w(\tilde R' \cap E(\tilde L)).
\end{align}

In practice however, many of the edge subsets $R$ induce the same equation. In fact there are only a finite number of equations in the particular $Z$-invariant moves that we will consider. 

\begin{remark}
	An interesting aspect of $Z$-invariance is, that in the statistical physics community the existence of these moves is sometimes considered as a manifestation of integrability. This makes sense because statistical models often involve edge-weights that are already factorizing in the 2D-data, that is the edge-weights are functions of parameters that are defined on 1D-data. Thus in these factorizing cases we are dealing with 2D-systems that exhibit 3D-consistency and are thus integrable also in the sense of discrete integrable systems. With non-factorizing edge weights on $G$, the existence of the moves themselves establishes a 3D-system, but not necessarily a 3D-system with 4D-consistency. Therefore for arbitrary edge weights the existence of $Z$-invariance is a priori not enough for discrete integrability. Surprisingly, in the examples we consider -- the dimer model, the spanning tree model and the Ising model -- it turns out that the models are actually 4D-consistent, and thus also discretely integrable.
\end{remark}

\begin{remark}
	We are not aware of any approach, in any of the statistical models, that classifies all the possible $Z$-invariant local moves in a systematic way. Thus the local moves that we present in the following sections are not derived by some algorithm. We therefore consider the existence of local moves as discoveries, and there is no known pattern, except to use trial and error for combinations of statistical models and reasonably simple local changes of combinatorics.
\end{remark}

\section{Dimers, \lowercase{d}KP and cluster algebras}\label{sec:dimers}
With the definitions of the previous section, all we need to define dimer statistics is the definition of the sample space depending on $G$.
\begin{definition}\label{def:dimer}
	The sample space of the \emph{dimer model} is
	\begin{align}
		\Omega_G = \{ R\subset E : \forall v\in V\ \exists! e \in R \mbox{ with } v\in e  \}.
	\end{align}
	Every set $R \in \Omega_G$ is called a \emph{dimer configuration}.
\end{definition}

In words, for each dimer configuration $R\in \Omega_G$ and for every vertex $v\in V$ there is exactly one edge $e\in E$ that contains $v$. In other words, every dimer configuration $R$ consists of a disjoint edge cover of the vertex set. Note that in the combinatorics community dimer configurations are called \emph{perfect matchings}. This model is well-defined on general graphs, but for our purposes it suffices to limit ourselves to planar, \emph{bipartite} graphs.

In this section we assume that all graphs $G$ that we consider are such that $\Omega_G \neq \emptyset$.  In general, bipartite graphs that admit dimer configurations are characterized by Hall's theorem. In Section \ref{sec:dimerinvariants} we generalize perfect matchings to almost perfect matchings, and show that whenever $G$ is the associated graph $\pb$ of some TCD, then $\pb$ admits an almost perfect matching.

In the dimer model, there exists a particular sort of \emph{gauge freedom}. Assume we fix $\rho \in \R^+, v\in V$ and replace the weight function $w$ by
\begin{align}
	w': E\rightarrow \R, \quad w'_e = \begin{cases}
		w_e & e \not\sim v,\\
		\rho w_e & e\sim v.
	\end{cases}
\end{align}
This operation scales the edge weights around a fixed vertex. Let us denote by $Z[w]$ and $Z[\omega']$ the partition function as a polynomial over the edge weights $\omega$ and $\omega'$ respectively. Given a dimer configuration, each vertex appears in exactly one of its edges, thus we observe that
\begin{align}
	Z[w'] = \rho Z[w],
\end{align}
so the probability (as in Equation \eqref{eq:omegaprob}) of each dimer configuration is unchanged. Hence there is a gauge or scaling freedom for each vertex that does not affect statistics. A complete set of invariants with respect to gauge are the so called \emph{face variables}
\begin{align}
	X : F \rightarrow \R^+, \ X_f = \frac{w(e_1)\cdot w(e_3)\dotsm w(e_{m-1})}{w(e_2)\cdot w(e_4)\dotsm w(e_{m})},\label{eq:dimerfacevars}
\end{align}
where $(e_k)_{0<k\leq m}$ is a counterclockwise oriented cycle bounding $f$ that starts at a black vertex. The face weights are gauge invariant because any scaling factor $\rho$ at a vertex does either not appear in $X_f$ or cancels in the fraction. If we include the outer face then the only constraint on the face variables is that
\begin{align}
	\prod_{f\in F} X_f = 1.
\end{align}
In fact for any set of face variables that satisfy the global constraint there are edge weights that realize the face variables. This was already proven in Remark \ref{rem:altconstruction}, where we called this ``solving the gauge problem''. Thus we think of the face variables as coordinates of possible dimer statistics on a given graph $G$.

\begin{figure}
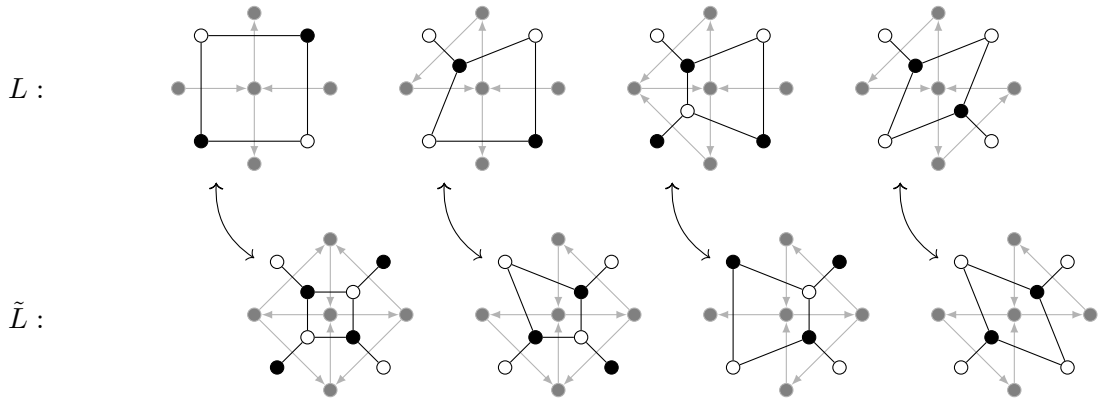


	
	\caption{Local changes of combinatorics in a bipartite graph $G$. The corresponding quiver mutations are drawn in gray.}
	\label{fig:localmovedimer}
\end{figure}

\begin{figure}
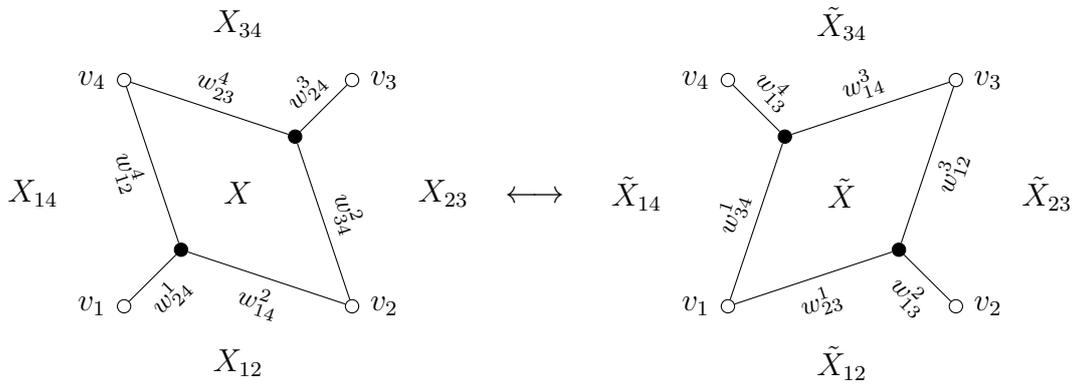

	%
	\caption{Labeling of edge weights in the spider move.}
	\label{fig:dimerweightsspider}
\end{figure}

Let us now shift our attention to the existence of $Z$-invariant moves. We have not chosen $X$ as the letter for the face variables by accident, indeed it turns out that these are the $X$-cluster variables for the quiver $\qui=G^*$, where the arrows are oriented counterclockwise around each black vertex. We present four pairs of subgraphs $L, \tilde L$. For each pair there is a function $f$ such that $(L,\tilde L,f)$ is a $Z$-invariant move. Also, for each pair there is a function $g$ such that $(\tilde  L,L,g)$ is a $Z$-invariant move. These subgraphs are listed in Figure \ref{fig:localmovedimer} and are exactly the duals of the graphs appearing in Figure \ref{fig:quivermut}, where we did list the allowed quiver mutations. Note that each local move is \emph{centered} at a face of degree 4, as depicted in Figure \ref{fig:localmovedimer}. We do the calculations for one example pair, the pair that is related by the spider move. The labeling is shown in Figure \ref{fig:dimerweightsspider}. The corresponding graphic equations, that are induced by comparing which boundary vertices of $L$ and $\tilde L$ are matched by the restricted dimer configuration, are
\begin{align}
	\tzspider11 &= \lambda \left[\ \tzspider[90]22 + \tzspider[90]33\ \right],\label{eq:dimergraphicsstart} \\
	\tzspider12 &= \lambda \ \tzspider[90]31, \qquad \tzspider21 = \lambda \ \tzspider[90]13,\\
	\tzspider13 &= \lambda \ \tzspider[90]12, \qquad \tzspider31 = \lambda \ \tzspider[90]21,\\
	\tzspider22 + \tzspider33 &= \lambda \ \tzspider[90]11,\label{eq:dimergraphics}
\end{align}
where $\lambda\in \R$. In terms of edge weights we obtain that
\begin{align}
	w_{24}^1 w_{24}^3 &= \lambda [w_{34}^1 w_{12}^3+w_{23}^1 w_{14}^3  ],\\
	w_{24}^1 w_{23}^4 &= \lambda w_{23}^1 w_{13}^4, \qquad w_{14}^2 w_{24}^3 = \lambda w_{13}^2 w_{14}^3,\\
	w_{24}^1 w_{34}^2 &= \lambda w_{34}^1 w_{13}^2, \qquad w_{24}^3 w_{23}^4 = \lambda w_{12}^3 w_{13}^4,\\
	w_{14}^2 w_{23}^4 + w_{34}^2 w_{12}^4 &= \lambda w_{13}^2 w_{13}^4.
\end{align}

As we have seen earlier, there is gauge freedom in the edge weights and thus it is not a surprise that there is no unique solution for the new edge weights. Instead it makes sense to look at how the face variables have to change so that the local moves are $Z$-invariant. It turns out that the corresponding equations for the face variables are the same for any of the pairs $(L,\tilde L)$.

\begin{theorem}
	Let $L,\tilde L$ be a pair of graphs of Figure \ref{fig:localmovedimer}. Let $f: \R^E \rightarrow \R^{\tilde E}$ be a function such that $(L,\tilde L,f)$ is a local move and the face weights satisfy 
	\begin{align}
		\tilde X &= X^{-1},\nonumber\\
		\tilde X_{12} X_{12}^{-1}  &= \tilde X_{34} X_{34}^{-1} = 1+X,\label{eq:dimerfacemove} \\
		\tilde X_{23} X_{23}^{-1}  &= \tilde X_{14} X_{14}^{-1} = (1+X^{-1})^{-1}.\nonumber
	\end{align}
	 Then $(L,\tilde L,f)$ is a $Z$-invariant move.
\end{theorem}
\proof{
	Straightforward calculation.\qed
}

There is actually a correspondence between a local move in the dimer model and a mutation in an associated cluster structure.

\begin{definition}\label{def:dimerclusterstructure}
	Let $G$ be a bipartite planar graph, on which we consider a dimer model with face weights $X_G: F(G) \rightarrow \C$. Let $\qui$ be the PDB quiver that is the oriented graph dual of $G$ and let $X: V(\qui) \rightarrow \C$ be cluster variables on $\qui$ such that $X(v) = X_G(v^*)$ for every vertex $v$ of $\qui$. We call $(\qui,X)$ the \emph{cluster structure of the dimer model} on $G$. 
	A mutation at a vertex $v \in V(\qui)$ of degree 4 corresponds to the local move centered at $v^* \in F(G)$.
\end{definition}

This is well-defined because the combinatorics of $G$ are dual to $\qui$ and stay dual to $\qui$ as we observe in Figure \ref{fig:localmovedimer}, and because Equations $\eqref{eq:dimerfacemove}$ correspond to the equations for the mutations of cluster variables as in Definition \ref{def:mutationx}. The converse is also possible: to any cluster structure on a PDB quiver $\qui$ (see Definition \ref{def:planquiver}) one can assign a dimer model on $\qui^*$ with face weights that are equal to the cluster variables.

Let us now show one way in which the dKP equation appears in the dimer model. Assume the edge weights factorize into face- and vertex-potentials, that is
\begin{align}
	w_e = \frac{\theta_v\theta_{v'}}{\tau_f\tau_{f'}} \quad \mbox{where } e^*=(f,f') \mbox{ and } e=(v,v').
\end{align}
If we mutate at face $f$, the edge weights still factor onto faces and vertices as before. The vertex potentials do not affect the face weights and can be chosen arbitrarily. The new face potentials  $\tilde \tau$ only differ from the old potentials at $f$. If $f_1,f_2,f_3,f_4$ are the neighbouring faces of $f$ then the potentials satisfy
\begin{align}
	\tau_f \tilde \tau_f = \tau_{f_1} \tau_{f_3} + \tau_{f_2} \tau_{f_4}.
\end{align}
This is the \emph{dKP equation} \cite{hirotaequation}. It is also the mutation rule for the cluster variables $\tau$, see Definition~\ref{def:mutationtau}.

\begin{remark}\label{rem:dimerextras}
	We have given a very brief introduction into the dimer model, a more encompassing introduction can be found in Kenyon's lecture notes \cite{kenyondimerintro}. Let us mention a few more references that we found particularly interesting from an integrability or geometry perspective.
	\begin{itemize}
		\item The dimer model received considerable attention when Kasteleyn \cite{kasteleyn} as well as Temperley and Fisher \cite{tfdimers} were able to give exact solutions in the planar case using determinants. 
		\item The first mention of $Z$-invariant moves was by Kuperberg \cite{kuperberg}. It is interesting to note that $Z$-invariant moves on the dimer model can be interpreted as the Desnanot–Jacobi identity on determinants of certain matrices \cite{zeilberger}.
		\item $Z$-invariant moves were applied successfully by Propp \cite{propp} to generate uniformly distributed dimer configurations.
		\item Goncharov and Kenyon \cite{gkdimers} discovered the so called dimer integrable systems that give Hamiltonians for dimer models on tori in terms of partition functions.
		\item There has been much investigation of the dimer model on isoradial graphs initiated by an article of Kenyon \cite{kenyonisoradial}. Since then a lot of results have been discovered, see \cite{btisosurvey} for a survey and additional references. We also want to highlight very interesting recent developments \cite{bctelliptic} that relate isoradial graphs, dimers and elliptic curves.\qedhere
	\end{itemize}
	  
\end{remark}

\section{Spanning trees, \lowercase{d}BKP and the cuboctahedral quiver} \label{sub:spanningtrees}\label{sec:spanningtrees}

In this section we look at the statistics of spanning trees, the $Z$-invariant move and the relation to cluster structures. The first to look at spanning trees was Kirchhoff \cite{kirchhofftrue}, although he did not think about them in a probabilistic manner.

\begin{definition}
	The sample space of the \emph{spanning tree model} is
	\begin{align}		
		\Omega_G = \{ R\subset E : R \mbox{ covers } V, \ R \mbox{ is cycle free}, R \mbox{ is connected}  \}.
	\end{align}
	Each $R\in \Omega_G$ is called a \emph{spanning tree}.
\end{definition}

\begin{figure}
	\begin{tikzpicture}[baseline={([yshift=-.7ex]current bounding box.center)},scale=2] 
	\tikzstyle{bvert}=[draw,circle,fill=black,minimum size=3pt,inner sep=0pt]
	\tikzstyle{gvert}=[draw,circle,fill=white,minimum size=3pt,inner sep=0pt]
	
	\node[bqgvert,label=left:$v_1$] (v1)  at (-.87,-.5) {};
	\node[bqgvert,label=right:$v_2$] (v2) at (.87,-.5) {};
	\node[bqgvert,label=left:$v_3$] (v3) at (0,1) {};
	
	\node[wqgvert,label=left:$v_{12}$] (v12)  at (0,-1) {};
	\node[wqgvert,label=right:$v_{23}$] (v23) at (.87,.5) {};
	\node[wqgvert,label=left:$v_{13}$] (v13) at (-.87,.5) {};
	
	\node[wqgvert,label=left:$v$] (v) at (0,0) {};
	
	\draw[-] (v2) -- (v3) node[midway, left]{$c^{23}$};
	\draw[-] (v1) -- (v3) node[midway, right]{$c^{13}$}; 
	\draw[-] (v1) -- (v2) node[midway, above]{$c^{12}$};

	\draw[-,gray] (v) -- (v1);
	\draw[-,gray] (v) -- (v2);
	\draw[-,gray] (v) -- (v3);
	\draw[-,gray] (v1) -- (v12) -- (v2) -- (v23) -- (v3) -- (v13) -- (v1);
	
	\end{tikzpicture}\hspace{1cm}
	\begin{tikzpicture}[baseline={([yshift=-.7ex]current bounding box.center)},scale=2] 
	\tikzstyle{bvert}=[draw,circle,fill=black,minimum size=3pt,inner sep=0pt]
	\tikzstyle{gvert}=[draw,circle,fill=white,minimum size=3pt,inner sep=0pt]
	
	\node[bqgvert,label=left:$v_1$] (v1)  at (-.87,-.5) {};
	\node[bqgvert,label=right:$v_2$] (v2) at (.87,-.5) {};
	\node[bqgvert,label=left:$v_3$] (v3) at (0,1) {};
	
	\node[wqgvert,label=left:$v_{12}$] (v12)  at (0,-1) {};
	\node[wqgvert,label=right:$v_{23}$] (v23) at (.87,.5) {};
	\node[wqgvert,label=left:$v_{13}$] (v13) at (-.87,.5) {};
	
	\node[bqgvert,label=left:$v_{123}$] (v123) at (0,0) {};
	
	\draw[-,gray] (v123) -- (v13);
	\draw[-,gray] (v123) -- (v12); 
	\draw[-,gray] (v123) -- (v23);
	
	\draw[-] (v123) -- (v1) node[midway, below]{$c^{23}_1$};
	\draw[-] (v123) -- (v2) node[midway, above]{$c^{13}_2$};
	\draw[-] (v123) -- (v3) node[midway, left]{$c^{12}_3$};
	\draw[-,gray] (v1) -- (v12) -- (v2) -- (v23) -- (v3) -- (v13) -- (v1);
	
	\end{tikzpicture}
	\caption{Labeling conventions for the front (left) and back (right) of a cube for spanning trees.}
	\label{fig:sptreecube}
\end{figure}
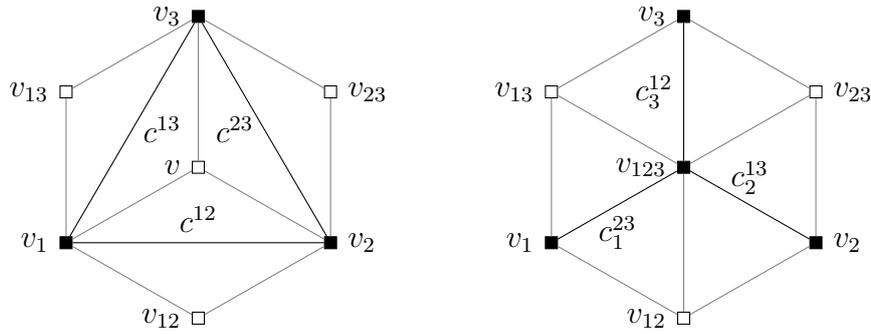

We can replace a triangle with a star in the graph $G$ or vice versa. As we intend to find a $Z$-invariant move, this leads to the graphic equations
\begin{align}
	\tztriangle000 &= \lambda\left[\tzstar100 + \tzstar010 + \tzstar001\right], \\
	\tztriangle100 &= \lambda\left[\tzstar011\right] , \qquad \tztriangle010 = \lambda\left[\tzstar101\right] , \qquad \tztriangle001 = \lambda\left[\tzstar110\right],\\		
	\tztriangle110 + \tztriangle101 + \tztriangle011 &= \lambda\left[\tzstar111	\right].
\end{align}
We identify the four vertices that appear in the star with the four black vertices of a cube. Subsequently we identify the six edges of the star and the triangle with the six quads of the cube, see Figure \ref{fig:sptreecube}. With this identification the translation of the graphic equations above into algebra is
\begin{align}
	1 &= \lambda [c^{12} + c^{23} + c^{13}],\\
	c^{12}_3 &= \lambda [c^{23} c^{13}], \qquad c^{23}_1 = \lambda [c^{12} c^{13}], \qquad c^{13}_2 = \lambda [c^{12} c^{23}],\\
	c^{12}_3 c^{23}_1 + c^{13}_2 c^{23}_1 + c^{12}_3 c^{13}_2  &= \lambda [c^{12} c^{23} c^{13}].
\end{align}
From these equations one calculates that the following theorem holds.

\begin{theorem}\label{th:startriangle}
	Let $L$ be the triangle and $\tilde L$ the star.  Let $f$ be the function that maps the three edge weights $c^{kl}$ to 
	\begin{align}
		c^{kl}_j  = \frac{1}{c^{kl}}\frac{c^{12}c^{23}c^{13}}{c^{12}+c^{23}+c^{13}} \quad  \mbox{for } \{j,k,l\} = \{1,2,3\}
	\end{align}
	and let $g$ be the function that maps the edge weight $c^{kl}_j$ to
	\begin{align}
		c^{kl} =\frac{1}{c^{kl}_j} (c^{12}_3c^{23}_1 + c^{12}_3c^{13}_2 + c^{23}_1c^{13}_2) \quad  \mbox{for } \{j,k,l\} = \{1,2,3\}.
	\end{align}
	Then the triples $(L,\tilde L, f)$ and $(\tilde L,L,g)$ are $Z$-invariant moves. These moves are called the \emph{star-triangle moves} \cite{kennelly}.
\end{theorem}

Temperley's bijection \cite{temperleybijection} relates spanning trees on $G$ to dimer configurations on a refined graph $G^D$. As we focus on the cluster formulation on the dual of a dimer graph, we explain the relevant construction directly on the level of the quiver \cite{gkdimers}. 

\begin{figure}
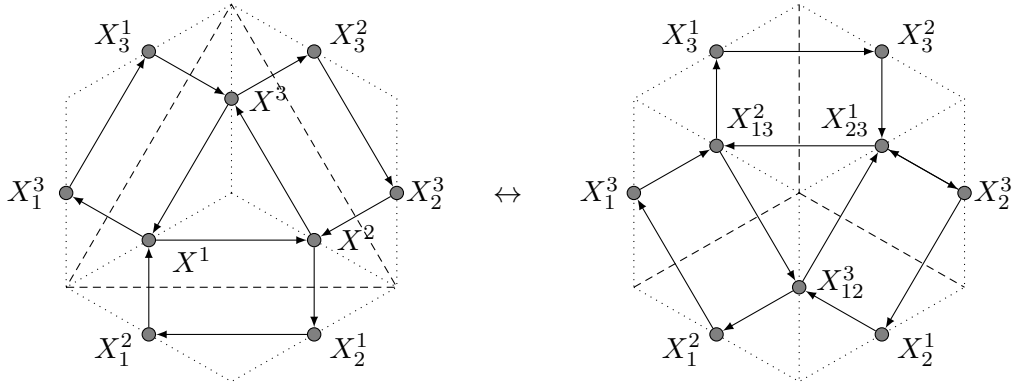
 
	\centering\small

	\caption{Before and after the cube-flip: the graph $G$ for the spanning trees (dashed), the quad-graph $\qg_G$ (dotted) and the corresponding cuboctahedral quiver $\qui$ (arrows and gray vertices).}
	\label{fig:spantreequiver}
\end{figure}

\begin{definition}\label{def:resistorclusterstructure}
	Let $G$ be a planar graph, $\qg_G$ the associated quad-graph (Definition \ref{def:graphtoqg}) and let $\qui$ be the cuboctahedral quiver (Definition \ref{def:cuboquiver}) of $\qg_G$, see also Figure \ref{fig:spantreequiver}. Let $c: E(G) \rightarrow \C$ be a function on the edges of $G$ (called \emph{conductances}), and define a function $X: V(\qui) \rightarrow \C$ by
	\begin{align}
		X_{(e,e')} = \frac{c_e}{c_{e'}}, \label{eq:clusterconductances}
	\end{align}
	where $e,e' \in E(G)$ and $(e,e')$ is a corner of $G$ (thus a vertex of $\qui$) such that $e'$ appears after $e$ in counterclockwise order at the common vertex of $e$ and $e'$. Then we call $(\qui,X)$ the \emph{resistor cluster structure} of $G$.
\end{definition}

In particular, we observe in Figure \ref{fig:spantreequiver} that triangle and star of $G$ are replaced by the two fundamental domains of the cuboctahedral quiver (see Section \ref{sus:cuboctahedral}). Note that we use shift notation (see page \pageref{pag:shiftnot}) in Figure \ref{fig:spantreequiver}, thus $c_e$ is replaced by $c^{ij}$ whenever $e = (v_i, v_j)$.

\begin{lemma}
	Let $c$ be the edge weights for a spanning tree model on the graph $G$ and $\tilde c$ the edge weights on $\tilde G$ after a star-triangle move on $G$. Let $(\qui,X), (\tilde \qui, \tilde X)$ be the resistor cluster structures of $G$ and $\tilde G$ before and after the star-triangle move respectively. Then $(\tilde \qui, \tilde X)$ is obtained from $(\qui,X)$ via the cuboctahedral flip, which is a sequence of four mutations (see Definition \ref{def:cuboflip}).
\end{lemma}
\proof{A small calculation using the fact that $X^1X^2X^3=1$ in the resistor subvariety.\qed}

Not every cluster structure that has a cuboctahedral quiver is a resistor cluster structure, but there is a simple criterion. 

\begin{definition}\label{def:resistorsub}
	Let $(Q,X)$ be a cluster structure with cuboctahedral quiver of the quad-graph $\qg$. Then the cluster variables $X$ are in the \emph{resistor subvariety} \cite{gkdimers} if for every non-boundary vertex $v$ of $\qg$
	\begin{align}
		\prod_{(e,e')\sim v} X_{(e,e')} = 1
	\end{align}
	holds, where the product is over corners $(e,e')$ adjacent to $v$.
\end{definition}

The resistor cluster structure is in the resistor subvariety because the cluster variables on the quad edges are quotients of the conductances $c_e$ that live on the quads of the quad-graph. Conversely, if $X$ are variables in the resistor subvariety then we can integrate them to conductances on the quads. The integration constant that we can pick results in a simultaneous scaling of all conductances. This simultaneous scaling does not affect the statistics of the spanning tree model and is considered a gauge freedom.

\begin{remark}\label{rem:linresistorproblem}
	The name ``resistor network subvariety'' is due to the fact that the rule for the replacements of conductances originates in the study of the star-triangle move on resistor networks. In particular, from Kirchhoff's circuit laws \cite{kirchhofftrue} and Ohm's law \cite{ohm} one finds that the voltage solves a discrete Laplace equation. Now we again want to replace a triangle with a in the graph $G$, and assign new conductances (independent of $I$ and $U$) such that the solutions for current $I$ and voltage $U$ only have to be changed at the involved edges. This results in the three equations
	\begin{align}
		(U_1 - U_2) c^{12} + (U_1 - U_3) c^{13} &= (U_1 - U_{123}) c^{23}_1,\\
		(U_2 - U_1) c^{12} + (U_2 - U_3) c^{23} &= (U_2 - U_{123}) c^{13}_2,\\
		(U_3 - U_1) c^{13} + (U_3 - U_2) c^{23} &= (U_3 - U_{123}) c^{12}_3.
	\end{align}
	
	These relations are satisfied if and only if the conductances satisfy the conductance equations that appear in the star-triangle move, see Theorem \ref{th:startriangle}. Note that it is possible to consider the voltages also as functions $U: V\rightarrow \CP^2$ to a 2-dimensional space and still the conductances from the star-triangle move satisfy the equations above. But this is not generally possible in higher dimensions. In fact, the condition for such a change of conductances to exist is exactly that all voltages corresponding to the black vertices of a cube are mapped to a common plane. This condition will manifest in K{\oe}nigs nets, see Section \ref{sec:konigs}. 
\end{remark}

Let us investigate in what sense spanning trees and the resistor subvariety relate to the discrete BKP equation. 

\begin{definition}\label{def:bkp}
	A map $f: \Z^3\rightarrow \C$ satisfies the \emph{discrete BKP equation} \cite{miwabkp} if on each cube of $\Z^3$ it satisfies
	\begin{equation}
		f f_{123} + f_1f_{23} + f_2f_{13} + f_3f_{12} = 0.\label{eq:bkp}\qedhere
	\end{equation}
\end{definition}
	
We define conductances $c: F(\Z^3) \rightarrow \C$ on the quads of $\Z^3$ from a function $f:\Z^3\rightarrow \C$ by
\begin{align}
	c^{kl} = \begin{cases}
	\frac{ff_{kl}}{f_kf_l} & \mbox{$f$ is at a white vertex},\\
	\frac{f_kf_l}{ff_{kl}} & \mbox{$f$ is at a black vertex}.
	\end{cases}\label{eq:conductancefactorization}
\end{align}
If the conductances are defined in this way, a straightforward calculation \cite{gkdimers} reveals that the conductances around a cube satisfy
\begin{align}\label{eq:condrecurrence}
	c^{kl}_j c^{kl} = \frac{c^{12}c^{23}c^{13}}{c^{12}+c^{23}+c^{13}} = c^{12}_3c^{23}_1 + c^{12}_3c^{13}_2 + c^{23}_1c^{13}_2,
\end{align}
for $\{j,k,l\} = \{1,2,3\}$. Thus the dBKP equation together with the factorization \eqref{eq:conductancefactorization} implies the star-triangle relations from Theorem \ref{th:startriangle}. Conversely, any conductances on $F(\Z^2)$ or a rectangular subgraph of $F(\Z^2)$ can be factorized as in \eqref{eq:conductancefactorization}. The factorization can be extended to all of $\Z^3$ if and only if the conductances satisfy the star-triangle relations in every cube.

\begin{remark}
	In enumerative combinatorics the dBKP equation is also known as the \emph{cube recurrence} \cite{cscube}. Under that alias the 4D-consistency of the discrete BKP equation was shown by Henriques and Speyer \cite{hscube}.
\end{remark}

\section{The Ising model, \lowercase{d}CKP and the hexahedral quiver} \label{sub:ising}\label{sec:ising}

The Ising model assigns a spin $s_v$ (spin up or down) to each vertex $v$ of a graph $G$. We call such an assignment the spin state $s\in 2^V$. The weight $\omega(s)$ of a spin state is 
\begin{align}
	\omega(s) = \prod_{v\in V(G)}\prod_{v'\sim v} \begin{cases}
		\omega_e & \mbox{if } s_v \neq s_{v'},\\
		1 & \mbox{if } s_v = s_{v'},
	\end{cases}
\end{align}
where $\omega_e$ are the edge weights. The partition function of the Ising model is thus
\begin{align}
	Z_{\mathcal I} = \sum_{s\in 2^V} \omega(s).
\end{align}
Now assume that $G$ is planar as well. Then each spin state $s$ can be identified with a subset of the dual edges $R_s\subset E^*$, where an edge $e^*$ with $e=(v,v')$ is in $R_s$ if and only if $s_v \neq s_{v'}$. Every subset $R_s$ consists of a union of cycles, called a contour of $G^*$. Also $R_s=R_s'$ if and only if $s=s'$ or if $s'$ has all spins reversed when compared to $s$. Therefore the partition function of the Ising model can also be written as 
\begin{align}
	Z_{\mathcal I} = \sum_{s\in 2^V} \prod_{e^*\in R_s} \omega_e.
\end{align}
This allows us to define the Ising model via edge subsets on a graph $G$, where we think of the spin states as instead living on the dual $G^*$.
\begin{definition}
	The sample space of the \emph{Ising model} is
	\begin{align}
		\Omega_{G} &= \{ R \subset E : \forall v\in V : |\{e\in R : v\in e\}| \in 2\Z \}.
	\end{align}
	Each $R\in \Omega_{G}$ is called a \emph{contour}.
\end{definition}

As in the case of the dimer and the spanning tree model, we are looking for a $Z$-invariant local move. As in the case of spanning trees, we can replace a triangle with a star in the graph $G$ or vice versa. This leads to the graphic equations
\begin{align}
	\tztriangle000 + \tztriangle111 &= \lambda\left[\tzstar000\right] ,\\
	\tztriangle100 + \tztriangle011 = \lambda\left[\tzstar011\right] , \qquad \tztriangle010 + \tztriangle101 &= \lambda\left[\tzstar101\right] , \qquad \tztriangle001 + \tztriangle110 = \lambda\left[\tzstar110\right].
\end{align}
We can translate this into the four algebraic equations
\begin{align}
	1+w^{12}w^{13}w^{23} &= \lambda, \\
	w^{jk}+w^{ik}w^{ij} &= \lambda w^{ik}_j  w^{ij}_k
\end{align}
for $\{i,j,k\} = \{1,2,3\}$. These equations can be solved which leads to the following theorem.
\begin{theorem}\label{th:startriangleising}
	Let $L$ be the triangle and $\tilde L$ the star.  Let $f$ be the function that maps the three edge weights $w^{ij}$ to 
	\begin{align}
		w^{ij}_k = \sqrt{\frac{(w^{jk}+w^{ik}w^{ij})(w^{ik}+w^{jk}w^{ij})}{(w^{ij}+w^{ik}w^{jk})(1+w^{12}w^{13}w^{23})}  }  
	\end{align}
	for $\{i,j,k\} = \{1,2,3\}$ and let $g$ be the function that maps the edge weight $w^{ij}_k$ to
	\begin{align}
		w^{ij} &= \frac{(w^{ij}_kw^{ik}_j)^2 + (w^{ij}_kw^{jk}_i)^2 - (w^{ik}_jw^{jk}_i)^2 - 1 +\sqrt{S^{+++}S^{+--}S^{-+-}S^{--+}}}{2 ((w^{ij}_k)^2-1) w^{ik}_j w^{jk}_i}			    
	\end{align}
	for $\{i,j,k\} = \{1,2,3\}$ where
	\begin{align}
		S^{\pm_1\pm_2\pm_3} &= 1\pm_1w^{ij}_kw^{ik}_j\pm_2w^{jk}_iw^{ij}_k\pm_3w^{ik}_jw^{jk}_i.
	\end{align}
	Then the triples $(L,\tilde L, f)$ and $(\tilde L,L,g)$ are $Z$-invariant moves. These moves are called the \emph{star-triangle moves of the Ising model} \cite{kashaev}.
\end{theorem}

To find the $\tau$ cluster variables, we need to perform a transformation of variables:
\begin{align}
	\frac14\left(w^{ij}-\frac1{w^{ij}}\right)^2= \frac{\tau \tau_{ij}}{\tau_i\tau_j}.\label{eq:kashaevtransform}
\end{align}
In order to define propagation, that is $\tau_{123}$ from the other variables we need to make a choice of sign for the square root. We choose the positive root, but other choices are possible and were investigated by Leaf \cite{leafkashaev}. The $\tau$ variables satisfy
\begin{align}\mbox{\small 
	$(\tau\tau_{123}-\tau_1\tau_{23}-\tau_2\tau_{13}-\tau_3\tau_{12})^2=4(\tau_1\tau_2\tau_{23}\tau_{13}+\tau_2\tau_3\tau_{13}\tau_{12}+\tau_3\tau_1\tau_{12}\tau_{23}-\tau_1\tau_2\tau_3\tau_{123}-\tau\tau_{12}\tau_{23}\tau_{13})$.}\label{eq:ckp}
\end{align}
This equation is called the \emph{discrete CKP equation} or also the \emph{Kashaev recurrence} \cite{kashaev}. In order to identify these $\tau$ variables as cluster variables of a quiver, let us first define the corresponding quiver.

\begin{definition}\label{def:isingcluster}
	Let $G$ be a planar graph. The quiver $\qui$ of the \emph{Ising cluster structure} is the hexahedral quiver of the quad-graph $\qg_G$ of $G$, see Definition \ref{def:hexaquiver}. The $\tau$ cluster variables of the Ising cluster structure at the vertices of the quad-graph are determined by Equation \eqref{eq:kashaevtransform}. The $\tau$ cluster variables at the quad-centers of the quad-graph are determined by
	\begin{equation}
		\tau^{kl} = \sqrt{\tau \tau_{kl} + \tau_k \tau_l} \label{eq:kashaevred},
	\end{equation}
	where $\tau^{kl}$ is the variable at the center of the quad and $\tau,\tau_k,\tau_{kl},\tau_l$ are the variables at the vertices.
\end{definition}

We observe that triangle and star are replaced by the two fundamental domains of the hexahedral quiver (see \circled1 and \circled6 of Figure \ref{fig:hexahedralquiver}). Kenyon and Pemantle \cite{kenyonpemantle} showed that in fact, the Kashaev recurrence is a reduction of the recurrence induced by the hexahedral flip. This discovery can be translated into the following lemma.

\begin{lemma}
	Let $\omega$ be the edge weights for the Ising model on a graph $G$ and $\tilde \omega$ the edge weights on $\tilde G$ after a star-triangle move on $G$. Let $(\qui,X), (\tilde \qui, \tilde X)$ be the Ising cluster structures before and after the star-triangle move. Then $(\tilde \qui, \tilde X)$ is obtained from $(\qui,X)$ via the hexahedral flip, which is a sequence of six mutations (see Definition \ref{def:hexaflip}).
\end{lemma}
\proof{
	Direct calculation.\qed
}

Not every cluster structure that has a hexahedral quiver is an Ising cluster structure. Moreover, we are more interested in the $X$ variables than the $\tau$ variables, because the $X$ variables are determined uniquely from TCD maps. Fortunately, there is an algebraic criterion to determine whether the $X$ variables of a hexahedral quiver are the Ising cluster structure of some graph. 

\begin{definition}\label{def:isingsub}
	Let $\qg$ be a quad-graph, and let $(\qui,X)$ be a cluster structure with hexahedral quiver $\qui = \qui_\qg$. Then we say the cluster variables $X$ are in the \emph{Ising subvariety} \cite{kenyonpemantle} if the equations
	\begin{align}
		(X_v)^2\prod_{f\sim v} (1+X^{f}) &= 1, \quad  \mbox{ for every black non-boundary vertex $v$ of } \qg,\\  \label{eq:isingsubvariety}
		(X_v)^2\prod_{f\sim v} (1+(X^{f})^{-1})^{-1} &= 1, \quad \mbox{ for every white non-boundary vertex $v$ of } \qg
	\end{align}
	hold, where $X^f$ is the cluster variable in the center of quad $f$ and the product is over all quads adjacent to $v$.
\end{definition}

\begin{lemma}
	The $X$ variables of the Ising cluster structure of a graph $G$ (Definition \ref{def:isingcluster}) are in the Ising subvariety. Vice versa, for any $X$ variables in the Ising subvariety there are corresponding $\tau$ variables of some graph $G$.
\end{lemma}
\proof{
	Assume the $\tau$ variables are the variables of an Ising cluster structure which live on the vertices of $\qg$. We label the variables such that $\tau$ is at a black vertex and so is $X$. We denote the cluster variables in the centers of the incident quads by $\tau^{k,k+1}$ resp. $X^{k,k+1}$. Then we observe that by Equation \eqref{eq:tauandx} the $X$ variables satisfy
\begin{align}
	X = \prod_k \frac{\tau_k}{\tau^{k,k+1}} &= \prod_k\frac{\tau_k}{\sqrt{\tau_k \tau_{k+1} + \tau \tau_{k,k+1}  }},\\
	1+X^{k,k+1} &= \frac{\tau_k \tau_{k+1} + \tau \tau_{k,k+1} }{\tau_k  \tau_{k,k+1}} \label{eq:isingxtaucalc}.
\end{align}

Therefore Equation \eqref{eq:isingsubvariety} is satisfied. The analogous argument works if we are at a white vertex. Now we prove the converse, that is we can begin with $X$-variables in the Ising subvariety and recover $\tau$ variables that define an Ising model. We follow the standard idea to propagate data in 2D-systems on quad-graphs by using the li-orientation as explained in Section \ref{sec:sweepsqg}. Given some initial values for the $\tau$-variables on the minimal elements of the labeled quad-graph, we can propagate to the remaining $\tau$-variables on the quad vertices via the $X$-variables in the quad centers. The $\tau$ variables in the quad centers are then determined by Equation \eqref{eq:kashaevred}. Therefore all $\tau$ variables are determined by the $X$ variables in the quad centers. Because of Equations \eqref{eq:isingxtaucalc}, the so determined $\tau$-variables are also consistent with the $X$-variables at the quad vertices, and the proof is finished.\qed
}

\begin{remark}
	Let us mention a few observations that we think add to the understanding of the Ising model and its relations to TCD maps, but that we cannot expand upon in this thesis because it would require a disproportionate amount of explanations of combinatorics and algebra.
	\begin{itemize}
		\item There is a linear problem accompanying the Kashaev equation similar to the resistor equations that accompany the spanning tree model. It involves spinors and is central to the study of s-embeddings and isoradial embeddings of the Ising model  \cite{csuniversality}.
		\item The $C^{(1)}_2$ loop model is a slight generalization of the Ising model. The Ising subvariety in fact captures more data than is necessary to define the statistics of the Ising model, but the data precisely defines the data of the $C^{(1)}_2$ loop model \cite{melottikashaev}.
		\item A relation between the Kashaev recurrence and M-systems which occur in the study of line complexes has been found by Bobenko and Schief \cite{bobenkoschieflinecomplexes, bobenkoschiefcirclecomplexes}. We will investigate a similar relation in Section \ref{sec:linearlinecomplexes}.
		\item For further interesting results on Ising models defined on isoradial graphs, see for example \cite{btcriticalisoradialising,btrisoradialising,tilierezdirac}. \qedhere
	\end{itemize}
\end{remark}

\section{Almost perfect matchings and TCD maps}\label{sec:dimerinvariants}

In Section \ref{sec:dimers} we have introduced the dimer model on planar bipartite graphs in order to establish relations to TCD maps. Recall that the graph $\pb$ associated to each TCD $\tcd$ is a planar bipartite graph. Moreover, with the VRC of each TCD map come edge-weights defined on the edges of $\pb$, and we also introduced local moves for VRCs, see Section \ref{sec:vrc}. Therefore $\pb$ is a natural candidate to study the dimer model on. 

At this point, we would like to state that the projective cluster structure of TCD maps coincides with the cluster structure of the dimer model on $\pb$. However, there is an important detail that requires attention first. Defining the dimer model on a graph $\pb$ is only interesting if there are actually any dimer configurations on $\pb$. A necessary condition for $\pb$ to admit dimer configurations is that the number of white vertices equals the number of black vertices. However, according to Theorem \ref{th:maxdim} the difference between white and black vertices is the maximal dimension of $\tcd$ plus one. But the maximal dimension of a minimal TCD $\tcd$ is always at least zero. Thus for the definition of the dimer model on $\pb$ we need a slight generalization of perfect matchings.

\begin{definition}\label{def:almostperfect}
	Let $G$ be a planar, bipartite graph such that all boundary vertices are white and such that $|W| \geq |B|$. An \emph{almost perfect matching} is a subset of edges that matches every vertex at most once, matches every interior vertex, and matches all but $|W|-|B|$ of the exterior vertices. If $|W|>|B|$ we define the sample space $\Omega_G$ of the dimer model on $G$ to be the set of almost perfect matchings of $G$.
\end{definition}

In the case that $G$ admits a perfect matching we have $|W|=|B|$ and Definition \ref{def:almostperfect} agrees with Definition \ref{def:dimer}. It is not a priori clear whether every graph $\pb$ coming from a TCD admits an almost perfect matching. 

\begin{figure}
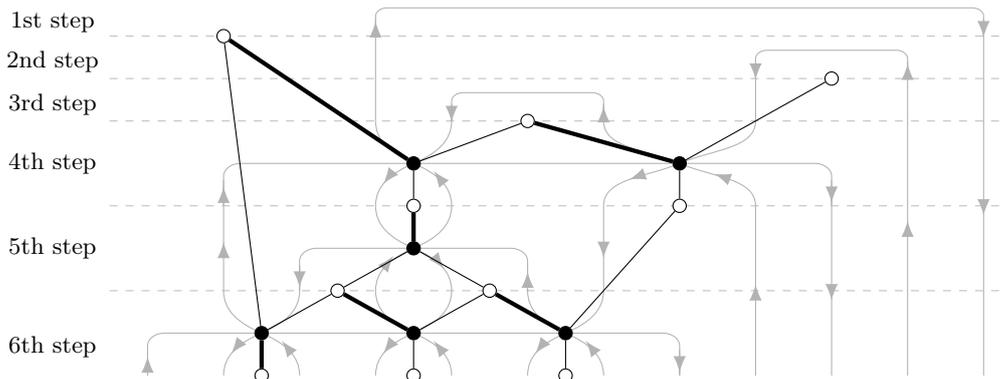
\hspace{-0.5cm}
\vspace{-4mm}
\caption{An almost perfect matching (bold) on the VRC of a standard diagram.}
\label{fig:standarddimer}
\end{figure}

\begin{lemma}
	Every bipartite graph $\pb$ associated to a TCD admits an almost perfect matching.
\end{lemma}
\proof{
	First of all, one observes that if $\pb$ admits an almost perfect matching then so does any $\tilde \pb$ related to $\pb$ via a spider move or resplit. The reason for this is that we can locally substitute the almost perfect matching. The graphic Equations \eqref{eq:dimergraphicsstart} - \eqref{eq:dimergraphics} illustrate this for the spider move and it is easy to verify the same for the resplit. Recall that in Definition \ref{def:standardtcd} we constructed the standard diagram for a given endpoint matching. We follow the idea of the proof of Theorem \ref{th:maxdim}, where we constructed a TCD map of maximal dimension for a given endpoint matching, see also Figure \ref{fig:standardtcdmap}. In that proof we went backwards through the steps of the construction of a standard diagram. We go backwards through these steps constructing the standard diagram here again, but we also select edges for an almost perfect matching of $\pb$ in the process, see Figure \ref{fig:standarddimer} for an example construction. In the first step we add one white vertex to the outer face, it is unmatched as there are no edges. Whenever we add a right moving strand, we add as many white as black vertices below that strand. Let $w_1,\dots,w_{k+1}$ be the white vertices above the strand, $v_1,\dots,v_{k}$ be the white vertices below the strand, and $b_1,\dots,b_k$ be the black vertices on the strand. If $w_i$ for $i\leq k$ is not already matched we add the edge $(w_i,b_i)$ to the matching, else we add the edge $(b_i,v_i)$. We observe that the number of unmatched vertices is preserved and that only boundary vertices are unmatched. Whenever we add a left moving strand, we add one more white vertex below the strand than we add black vertices on the strand. Let $w_1,\dots,w_{k-1}$ be the white vertices above the strand, $v_1,\dots,v_{k}$ be the white vertices below the strand, and $b_1,\dots,b_{k-1}$ be the black vertices on the strand. If $w_i$ for $i< k$ is not already matched we add the edge $(w_i,b_i)$ to the matching, else we add the edge $(b_i,v_i)$. We observe that the number of unmatched vertices is increased by one and that only boundary vertices are unmatched. Therefore we can always construct a TCD with the same endpoint matching as any given TCD that admits an almost perfect matching. Because this property is preserved under resplits and spider moves the claim is proven.\qed
}

Apart from the combinatorics, there is another important aspect that we have to address before associating a dimer model to a TCD map. Note that Definition \ref{def:projclusterstructure}, of the $X$-cluster variables of a TCD map contains a sign factor that does not appear in Equation \eqref{eq:dimerfacevars}, where we introduce the $X$-cluster variables of the dimer model. Hence, we cannot directly identify the edge-weights that appear in the VRC of a TCD map with the edge-weights of the dimer model. Instead we need an additional phase factor.

\begin{definition}
	Let $G$ be a bipartite planar graph. A Kasteleyn phase $\varphi:E\rightarrow \C$ is a function defined on the edges of $G$ such that
	\begin{align}
		\prod_{e\sim f} \varphi_e = (-1)^{\frac12\deg f +1},
	\end{align}
	for every face $f$ of $G$.
\end{definition}

It is almost obvious that every bipartite planar graph admits a Kasteleyn phase, for example by the same arguments as we used to solve the gauge problem in Remark \ref{rem:altconstruction}. Moreover, one can also define the Kasteleyn phase as taking only values in $\{\pm1\}$ and these phases also exist for any bipartite planar graph \cite{kasteleyn, tfdimers}. However, sometimes (see for example Section \ref{sec:cptemb}) it is useful to work with a Kasteleyn phase with values in $\C$.

\begin{definition}\label{def:tcddimers}
	Let $T: \tcdp \rightarrow \CP^n$ be a TCD map, let $\pb$ be the associated bipartite graph and let $\varphi$ be a Kasteleyn phase of $\pb$. The \emph{dimer model} associated to a TCD map $T$ is the dimer model on $\pb$, such that the sample space consists of the set of almost perfect matchings as explained in Definition \ref{def:almostperfect}. The edge-weights $\omega$ of the dimer model are determined by the edge-weights $\lambda$ of the VRC and the Kasteleyn phase $\phi$ as $\omega_e = \varphi^{-1}_e \lambda_e$ for all $e\in E$.
\end{definition}

In this definition, it is possible and will generically be the case that the edge-weights of the dimer model are not positive reals. Thus, a statistical interpretation is not available. However, the algebraic properties, especially the invariance under local moves of the partition function and the transformation behavior of the face variables are unchanged. On the other hand, there are TCD maps where the edge-weights of the associated dimer model are positive reals, see Chapter \ref{cha:cpone} for notable examples. A TCD map can be equipped with a VRC with positive edge-weights if and only if all the $X$-variables are positive. Moreover, one can of course associate a dimer model not only to $T$ itself but to each of the sections of $T$. In fact, the examples of Chapter~\ref{cha:cpone} are all with respect to the first section of $T$. It would definitely be an interesting undertaking to understand all TCD maps with positive $X$ variables, either of $T$ itself or with respect to one of its sections.  

\begin{theorem}\label{th:projdimercluster}
	Let $T$ be a TCD map with vector-relation configuration $\vrc$. The projective cluster structure for $T$ from Definition \ref{def:projclusterstructure} is the same as the cluster structure of the dimer model of $T$ as in Definition \ref{def:tcddimers}.
\end{theorem}
\proof{With the introduction of the Kasteleyn phase, this holds by definition.\qed}

Note that historically \cite{kasteleyn}, the Kasteleyn phase was used to introduce the Kasteleyn matrix $K$ with columns and rows corresponding to the black and white vertices of a planar bipartite graph $G$ and entries that correspond to $\varphi_e \omega_e$. Surprisingly, this allowed to write the partition function $Z$ simply as $Z=|\det K|$. It is also worth noting that one can see the lifts of TCD maps as elements of the cokernel of $K$. We do not pursue this approach here, but it has proven useful in the context of TCD maps defined on the torus, see for example joint work with George and Ramassamy \cite{agrcrdyn}. 

The main result and intent of this section is Theorem \ref{th:projdimercluster}, that is the identification of the TCD map cluster structure with the dimer cluster structure. However, the dimer model can also be used to associate \emph{global invariants} to TCD maps. We have not done a complete analysis of this subject, and thus cannot give a complete exposition either. It is the belief of the author that most of the necessary tools and results have already been obtained in a different context, in the context of Postnikov theory \cite{postgrass} (see also \cite{scattamp} for a less technical exposition). Still, we want to give a short outline of the possibilities.

Consider a TCD map $T$ and its VRC. Fix a subset $I$ of the boundary vertices of a graph $\pb$ such that $|I| = |W| - |B|$. Let $Z_I$ be the partition function associated to $\pb$ restricted to the almost perfect matchings that match $I$. By the results of Section \ref{sec:dimers}, $Z_I$ is invariant under choosing different homogeneous lifts of interior vertices of $\pb$, as this corresponds to scalings of edge-weights adjacent to vertices. Moreover, under local moves the partition function $Z_I$ is only scaled by a factor. Thus any quantity $\frac{Z_I}{Z_{I'}}$ is an invariant of $T$ under local moves, as well as rescalings of interior vertices. Moreover, consider the case where we have two families of sets $I_1,I_2,\dots, I_m$ and $I'_1,I'_2,\dots, I'_m$ such that $\cup_{i=1}^m I_i = \cup_{i=1}^m I'_i$. Then in fact the quantity
\begin{align}
	\frac{\prod_{i=1}^m Z_{I_i}}{ \prod{i=1}^m Z_{I'_i}}\label{eq:partitioninvariants}
\end{align}
is invariant under different lifts, local moves and projective transformations. Thus we obtain true global projective invariants of $T$. It is not a priori clear for which sets $I$ the restricted partition function $Z_I$ is actually non-trivial, that is for which $I$ the set of almost perfect matchings that match $I$ is non-empty. We strongly suspect however that for TCD maps defined on balanced TCD, that is where the endpoint matching is a cyclic shift, all $Z_I$ are non-trivial. Moreover, not all of the $Z_I$ are independent. Instead, they can be interpreted as $\tau$-cluster variables. Certain sets of true global invariants as in Equation \eqref{eq:partitioninvariants} can be interpreted as $X$-cluster variables. An additional question is whether, and if so how it is possible to reconstruct TCD maps from the true global invariants. We also believe this to be true, and we also believe the answers only need to be translated from Postnikov theory.

Finally, let us mention that there are other invariants that one can assign to TCD maps in particular affine gauges. We will focus on examples where the vertices of $\pb$ can be interpreted to live on the lattice $A_3$ in upcoming work with Melotti and de Tilière \cite{amdtdskp}.

\section{Geometry of TCD maps in the resistor and Ising subvarieties}

We have shown that every TCD map features a dimer model, see Definition \ref{def:tcddimers}. Every dimer model gives rise to a cluster structure. If the dimer model is associated to a TCD map, then the cluster structure of the dimer model is the projective cluster structure of the TCD map. Moreover, we have identified that if the cluster variables of a dimer model are in the resistor subvariety (Definition \ref{def:resistorsub}), they give rise to a spanning tree model. Analogously, if the cluster variables are in the Ising subvariety (Definition \ref{def:isingsub}) they give rise to an Ising model. Therefore it is reasonable to ask: For which TCD maps do the cluster variables take values in the resistor subvariety or the Ising subvariety? To answer this question one first has to confirm that the combinatorics are correct. Both the resistor and the Ising subvariety are associated to quad-graphs, albeit with different quivers. 
We show in the following sections that
\begin{enumerate}
	\item Q-nets are K{\oe}nigs nets \cite{ddgbook} if and only if their projective cluster variables are in the resistor subvariety (Section \ref{sec:konigs}),	
	\item Darboux maps are Carnot maps \cite{schieflattice} if and only if their projective cluster variables are in the Ising subvariety (Section \ref{sec:darbouxckp}),	
	\item line compounds are Doliwa compounds if and only if their projective cluster variables are in the resistor subvariety (Section \ref{sec:Doliwacomplex}).
\end{enumerate}

On the other hand, we have shown that the resistor subvariety is always accompanied by the dBKP equation while the Ising subvariety is always accompanied by the dCKP equation. The occurrences of dKP, dBKP and dCKP equations provide a link to the study of \emph{discrete integrable systems}. These equations have repeatedly been discovered to accompany objects in DDG as well, although the necessary calculations are done on a case by case basis.

\begin{definition}\label{def:bkpckp}
	Let $T$ be a TCD map. We say that $T$ 
	\begin{enumerate}
		\item is \emph{projective BKP} if the projective quiver is cuboctahedral and the projective cluster variables of $T$ are in the resistor subvariety,
		\item is \emph{projective CKP} if the projective quiver is hexahedral and the projective cluster variables of $T$ are in the Ising subvariety.
	\end{enumerate}
	We also say that $T$ is \emph{affine BKP (resp. CKP)} with respect to some hyperplane if the affine quiver is cuboctahedral (hexahedral) and the affine cluster variables are in the Resistor (resp. Ising) subvariety.
\end{definition}

With this definition, the claim that Q-nets are K{\oe}nigs nets if and only if their projective cluster variables are in the resistor subvariety is equivalent to claiming that a Q-net is projective BKP if and only if it is a K{\oe}nigs net. Analogously, we can rephrase the other two claims from the beginning of this section to state that a Darboux map is projective CKP if and only if it is a Carnot map and a line compound is projective BKP if and only if it is a Doliwa compound.

\begin{remark}
	Note that in the literature, the reductions of Q-nets, Darboux maps and line complexes (like K{\oe}nigs nets etc.) are defined by viewing these maps as defined on $\Z^N$. The definitions for the Resistor and Ising subvarieties on the other hand are given for general quad-graphs. As before, we can view (minimal) quad-graphs as Cauchy-data for the maps defined on $\Z^N$. We will show that the conditions associated to reductions of the maps on $\Z^N$ indeed translate to subvarieties for the cluster variables of the Cauchy-data.	
\end{remark}

We have already listed the equivalences between subvarieties of the projective cluster structures and geometric reductions. Let us list further results that we are going to show in the following.

\begin{enumerate}
	\item A CQ-net (also called C-quadrilateral lattice) is affine CKP with respect to the plane at infinity, see Section \ref{sec:darbouxckp}.
	\item A linear line complexes is affine CKP with respect to a line that is isotropic with respect to the associated anti-symmetric bilinear form, see Section \ref{sec:linearlinecomplexes}.
	\item An S-graph is affine CKP with respect to any line through the distinguished point at infinity, see Section \ref{sec:sgraphs}.
	\item An s-embeddings is affine CKP with respect to the point at infinity, see Section \ref{sec:sembeddings}.
	\item A Schief map is affine BKP with respect to the plane at infinity, see Section \ref{sec:pardm}.
	\item A reciprocal figure is projective BKP and affine BKP with respect to a line if we first take a lift to $\CP^3$, see Section \ref{sec:reciprocal}.
	\item A Q-net inscribed in a conic in $\CP^2$ is affine BKP with respect to any point of the conic, see Section \ref{sec:inscribedinconics}.
	\item A Q-net inscribed in a quadric in $\CP^3$ is affine BKP with respect to a generator of the quadric, see Section \ref{sec:quadricqnet}.
	\item Darboux maps tangent to a quadric are affine BKP with respect to a generator of the quadric, see Section \ref{sec:tangentdarboux}.
	\item The Plücker lift of an A-net is affine BKP with respect to an isotropic plane of the Plücker quadric, see Section \ref{sec:anets}.
	\item The Plücker lift of a Cox-lattice is affine BKP with respect to an isotropic line of the Plücker quadric, and affine BKP with respect to specific 3-spaces, see Section \ref{sub:cox}.
	\item An h-embedding is affine BKP with respect to the point at infinity, see Section \ref{sec:harmonicemb}.	
\end{enumerate}

\section{K{\oe}nigs nets, BKP and the Resistor subvariety}\label{sec:konigs}

The K{\oe}nigs nets we consider in this section are due to Bobenko and Suris \cite{bsmoutard} and independently Doliwa \cite{doliwatnets}, who called them B-quadrilateral lattices. There are multiple geometric and algebraic characterizations of K{\oe}nigs nets on lattices, see the DDG book \cite[Section 2.3]{ddgbook}. We give a definition directly on quad-graphs, which is a generalization of the geometric characterization on $\Z^2$ \cite{bsmoutard} as well as a reformulation of the geometric characterization on $\Z^N$ \cite{doliwatnets}.

\begin{definition}\label{def:koenigs}
	Let $\qg$ be a quad-graph and $q: V(\qg)\rightarrow \CP^n$ a Q-net. At any vertex $v\in V(\qg)$ of degree $d_v$, there are $d_v$ neighbouring vertices $v_1,v_2,\dots,v_{d_v}$ and $d_v$ diagonally neighbouring vertices $v_{12},v_{23},\dots,v_{d_v1}$. We call $q$ a \emph{K{\oe}nigs net} if at every vertex $v\in V(\qg)$ 
	\begin{align}
		\dim \spa \{q(v),q(v_{1}),q(v_{2}),\dots,q(v_{d_v})\} = d_v, \label{eq:koenigsgenericity}
	\end{align}
	and 
	\begin{equation}
		\dim \spa \{q(v_{12}),q(v_{23}),\dots,q(v_{d_v1})\} = d_v-1, \label{eq:koenigscondition}
	\end{equation}
	are satisfied.
\end{definition}
Note that Equation \eqref{eq:koenigsgenericity} is a genericity assumption on the Q-net $q$ while Equation \eqref{eq:koenigscondition} characterizes K{\oe}nigs nets.

\begin{figure}
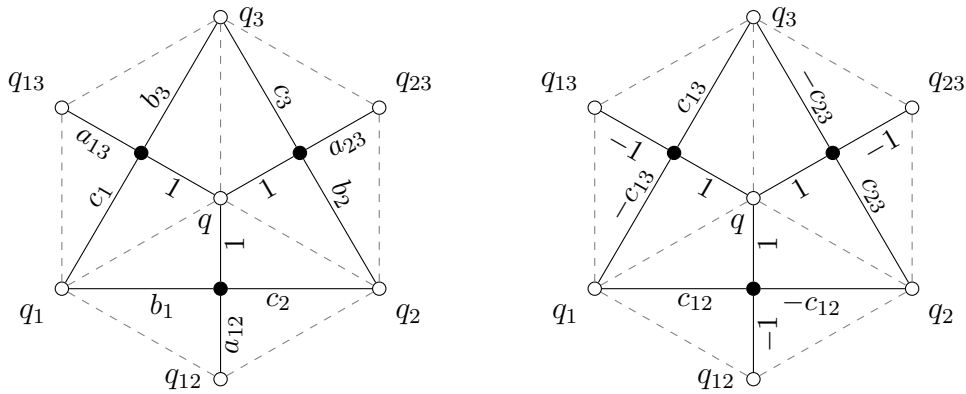

	
	\caption{Labels and edge-weights for the VRC of three adjacent quads in a Q-net cube on the left, for a T-net on the right.}
	\label{fig:konigslabels}
\end{figure}

\begin{theorem}\label{th:konigsbkp}
	A Q-net $q: V(\qg)\rightarrow \CP^n$ that satisfies the genericity condition \eqref{eq:koenigsgenericity} is projective BKP if and only if it is a K{\oe}nigs net.
\end{theorem}
\proof{
	We have to check that Equation \eqref{eq:koenigscondition} is satisfied for any $v\in V(\qg)$. Consider the case that $d_v=3$ first. Let us look at Figure \ref{fig:konigslabels}, where we scaled at the black vertices such that the edge weights incident to $v$ are 1. Therefore we obtain the three relations for the homogeneous lifts
	\begin{align}
		\hat q+b_1\hat q_1+a_{12}\hat q_{12} + c_2 \hat q_2 &= 0,\\
		\hat q+b_2\hat q_2+a_{23}\hat q_{23} + c_3 \hat q_3 &= 0,\\
		\hat q+b_3\hat q_3+a_{13}\hat q_{13} + c_1 \hat q_1 &= 0.
	\end{align} 
	We can solve the second equation for $q_2$ and the third for $q_3$ and then insert those two expressions into the first equation to obtain that
	\begin{align}
		\hat q\left(1-\frac{c_2}{b_2} + \frac{c_2}{b_2}\frac{c_3}{b_3}\right) + a_{12}\hat q_{12} - a_{23}\hat q_{23} + a_{13}\hat q_{13} +\hat q_1\left(b_1 + \frac{c_2}{b_2}\frac{c_3}{b_3}c_1\right) = 0.
	\end{align}
	If the coefficient of $\hat q$ vanishes then $q_1,q_{12},q_{23},q_{13}$ are coplanar, from which one can deduce that the whole cube is planar in contradiction to the genericity assumption. Therefore, $q,q_{12},q_{23},q_{13}$ are coplanar if and only if the coefficient of $\hat q_1$ vanishes, that is if
	\begin{align}
		\frac{c_1c_2c_3}{b_1b_2b_3} = -1
	\end{align}
	holds. Due to Definition \ref{def:bkpckp} a Q-net is BKP if the projective cluster variables (Definition~\ref{def:projclusterstructure}) are in the resistor subvariety, which in turn implies by Definition \ref{def:resistorsub} that the cluster variables multiply to 1 around each vertex of the quad-graph. At $v$ the three incident cluster variables are
	\begin{align}
		X^k &= -\frac{c_k}{b_k}.
	\end{align}
	That $X^1X^2X^3$ multiply to 1 is therefore equivalent to the coplanarity of $q,q_{12},q_{23},q_{13}$. If $d_v\neq 3$ the proof proceeds in the same manner.\qed
}

Without the relation to the resistor subvariety the result of Theorem \ref{th:konigsbkp} is already known, in the following sense: A Q-net is a K{\oe}nigs net if and only if its Laplace invariants multiply to one \cite[Exercise 2.21]{ddgbook}. Here, the Laplace invariants (see Definition \ref{def:laplaceinvariant}) are known from DDG and are exactly the projective cluster variables of the Q-net (see Lemma \ref{lem:laplaceinvprojvar}).

Note that the characterization of K{\oe}nigs nets in Theorem \ref{th:konigsbkp} is via projective invariants and is thus invariant under projections. We can thus drop the genericity assumption and consider Q-nets to be K{\oe}nigs nets if and only if the projective cluster variables are in the resistor subvariety.

The following theorem is well known \cite{doliwatnets}, but it is now also an easy corollary of Theorem \ref{th:konigsbkp}.
\begin{corollary}
	K{\oe}nigs nets are consistent.
\end{corollary}
\proof{We have to check that if a map $q: V(\qg)\rightarrow \CP^n$ is a K{\oe}nigs net then it is also a K{\oe}nigs net after a cube flip. However, we know that the resistor subvariety is preserved by cube flips and thus the claim is an immediate consequence of Theorem \ref{th:konigsbkp}.\qed}

In Theorem \ref{th:konigsbkp} we have established that the projective cluster variables of a K{\oe}nigs net are in the resistor subvariety. We have also given a result by Goncharov and Kenyon in Section \ref{sub:spanningtrees} that show how the cluster variables relate to the discrete BKP equation (Definition \ref{def:bkp}). We give a short outline now how these two results coincide in the case of K{\oe}nigs nets.

It was shown by Doliwa \cite{doliwatnets} that the VRC of K{\oe}nigs nets admit particular edge weights such that the equation for the lifts in each quad $(\hat q, \hat q_i,\hat q_j, \hat q_{ij})$ becomes
\begin{align}
	\hat q-\hat q_{ij} = c_{ij} (\hat q_i - \hat q_j),\label{eq:moutard}
\end{align}
with $c_{ij} \in \C$. This equation is also called the \emph{Moutard equation} \cite{miwabkp} and nets satisfying the Moutard equation in $\C^n$ -- not $\CP^n$ -- are known as \emph{T-nets} \cite{bsmoutard}. Note that in this representation the cluster variables (resp. the Laplace invariants) become
\begin{align}
	X^i = -\frac{c_{i,i-1}}{c_{i,i+1}}.
\end{align}
This equation is the same as in Definition \ref{def:resistorclusterstructure} for the cluster variables of the resistor cluster structure. Thus we can identify the coefficients $c_{ij}$ with conductances in a spanning tree model. The factorization that relates the conductances to the BKP equation (see Section \ref{sub:spanningtrees}) due to Goncharov and Kenyon \cite{gkdimers} is the same as the factorization that relates the coefficients of the T-net to the BKP equation due to Doliwa \cite{doliwatnets}. It would also be interesting if one could reproduce Doliwa's result of the existence of edge-weights \eqref{eq:moutard} for T-nets by using the li-orientation of quad-graphs. However, as we do not consider $T$-nets in the remainder, we do not consider this question here.

\section{Carnot maps, CKP, and the Ising subvariety}\label{sec:darbouxckp}

Schief introduced a consistent reduction of Darboux maps \cite{schieflattice}, that he related to a multi-ratio equation and the CKP equation. We give a slightly different characterization that is more advantageous for considering quad-graphs.

\begin{figure}
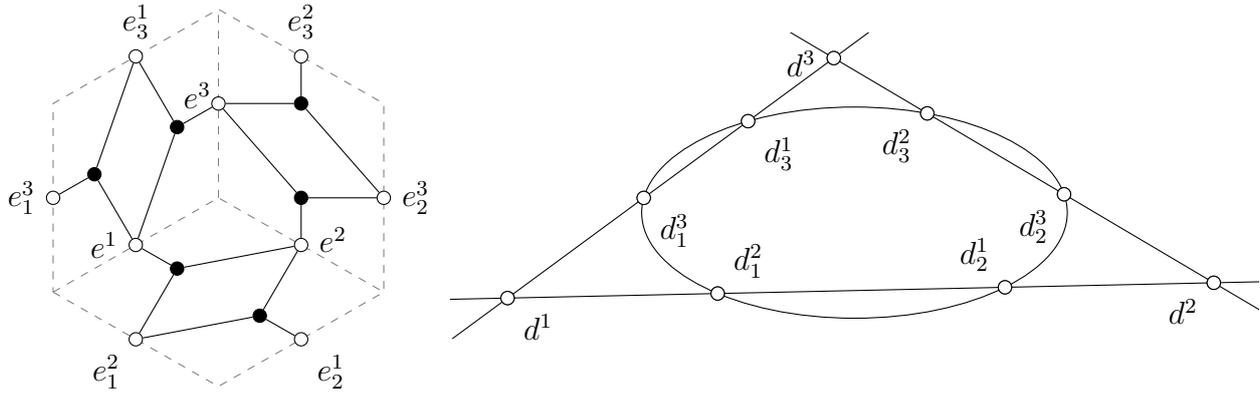


	\caption{The bipartite graph $\pb$ of a Darboux map and the geometry of a Carnot map.}
	\label{fig:carnot}
\end{figure}

\begin{definition}\label{def:carnot}
	Let $\qg$ be a quad-graph and $d: E(\qg) \rightarrow \CP^n$ a Darboux map. At any vertex $v\in V(\qg)$ of degree $d_v$ denote by $e^1,e^2,\dots e^{d_v}$ the adjacent edges and by $e^1_2,e^2_3,\dots,e^{d_v}_1$ and $e^1_{d_v}, e^2_1,\dots, e^{d_v}_{d_v-1}$ the other edges of the quads adjacent to $v$. We call $d$ a \emph{Carnot map} if around every vertex $v$ of $\qg$ 
	\begin{align}
		\dim \spa \{d(e^1),d(e^{2}),\dots,d(e^{d_v})\} = d_v - 1, \label{eq:carnotgenericity}
	\end{align}
	and if $d(e^1_2),d(e^2_3),\dots,d(e^{d_v}_1),d(e^1_{d_v}), d(e^2_1),\dots, d(e^{d_v}_{d_v-1})$ are on a common quadric.
\end{definition}

See Figure \ref{fig:carnot} for an example labeling and geometry. Note that Equation \eqref{eq:carnotgenericity} is a genericity assumption on $d$, while we consider the existence of the common quadric the characterizing property of Carnot maps. Schief calls these maps \emph{CKP lattices}. Unfortunately, from our point of view there are different maps resp.~lattices that feature the CKP equation, which is why we opted for an alternative naming. It will become obvious why we chose the name Carnot map in a moment. Also note that Schief's characterization is the special case for $d_v=3$, in which case the quadric is a conic. On $\Z^n$ Schief's characterization implies the characterization above for minimal quad-graphs as subgraphs of $\Z^n$.

The points in Definition \ref{def:carnot} are the intersections of the polygon $(d^1,d^2,\dots,d^n)$ with the quadric. This allows us to characterize CKP Darboux maps via a well known theorem of Carnot \cite{carnot}.

\begin{theorem}[Carnot's Theorem]\label{th:carnot}
	Let $d^1,d^2,\dots d^{n}$ be $n$ points that span $\CP^{n-1}$. For each $i\in \{1,2,\dots, n\}$ let $d_i^{i+1},d_{i+1}^i$ be two points on the line $d^id^{i+1}$ but not $d^i$ or $d^{i+1}$. Then the points $2n$ $d_i^{i+1},d_{i+1}^i$ are on a quadric not containing any of the points $d^i$ if and only if
	\begin{align}
		\mr(d^1,d^2_1, d^2,d^3_2, \dots,  d^n,d^1_n) \mr(d^1,d^1_2,d^2,d^2_3,\dots, d^n,d^n_1) = 1 \label{eq:darbouxckp}
	\end{align}
	holds.
\end{theorem}
\proof{
	Let us fix a homogeneous lift $\hat d^i$ of each $d^i$ and the other lifts via
	\begin{align}
		\hat d^i_{i+1} = u_{i} \hat d^i + v_{i+1} \hat d^{i+1} \quad\mbox{and}\quad \hat d_i^{i+1} = s_{i} \hat d^i + t_{i+1} \hat d^{i+1}. \label{eq:carnotlifts}
	\end{align}
	Let us consider the quadric as given by a bilinear form $b$ and write $b_{ij} = b(\hat d^i, \hat d^j)$. Then we obtain the equations
	\begin{align}
		u_iv_{i+1}^{-1} b_{ii} + u_i^{-1}v_{i+1} b_{i+1,i+1} + 2 b_{i,i+1} &= 0,\\
		s_it_{i+1}^{-1} b_{ii} + s_i^{-1}t_{i+1} b_{i+1,i+1} + 2 b_{i,i+1} &= 0.
	\end{align}
	We can eliminate $b_{i,i+1}$ to obtain 
	\begin{align}
		b_{i+1,i+1} = \frac{s_it_{i+1}^{-1}- u_iv_{i+1}^{-1}}{u_i^{-1}v_{i+1}-s_i^{-1}t_{i+1}}b_{ii} = \frac{s_iu_i}{t_{i+1}v_{i+1}} b_{ii}.
	\end{align}
	If we combine this equation for all the indices $i$ we obtain the consistency condition
	\begin{align}
		\prod_{i=1}^n \frac{s_iu_i}{t_{i+1}v_{i+1}} = 1.\label{eq:carnotconsistency}
	\end{align}
	On the other hand, we can use the same arguments as in Lemma \ref{lem:projclusterviadistances} to express the multi-ratios in the statement of the theorem via the coefficients in Equation \eqref{eq:carnotlifts}. As a consequence, the left hand side of Equation \eqref{eq:carnotconsistency} is nothing but the product of the multi-ratios in the statement of the theorem, which concludes the proof.\qed
}

After an initial choice of $b_{11}$, the equations appearing in the proof also determine the coefficients $b_{i,i+1}$. Thus, in $\CP^2$ the quadric (a conic) is completely determined, but in $\CP^n$ for $n>2$ this is not the case. 

\begin{corollary}[\cite{schieflattice}]\label{cor:carnotmr}
	A Darboux map $d: E(\qg)\rightarrow \CP^n$ that satisfies the genericity condition \eqref{eq:carnotgenericity} is a Carnot map if and only if Equation \eqref{eq:darbouxckp} holds around every vertex of the quad-graph.
\end{corollary}
Note that the dimensional genericity assumptions in Carnot's theorem is the genericity assumption of Definition \ref{def:carnot}, and the additional assumption in Carnot's theorem that certain points do not coincide is the general 0-genericity assumption on TCD maps.

\begin{theorem}\label{th:carnotckp}
	A Darboux map $d: E(\qg)\rightarrow \CP^n$ that satisfies the genericity condition \eqref{eq:carnotgenericity} is a Carnot map if and only if is projective CKP.
\end{theorem}
\proof{
	We compare Corollary \ref{cor:carnotmr} to Definition \ref{def:isingsub} of the Ising subvariety. We give the proof for a black vertex $v$ of the quad-graph, for white vertices the proof works analogously. The projective cluster variables of a Darboux map around $v$ are
	\begin{align}
		X^{i,i+1} &= -\cro(d^{i+1},d^i,d^{i+1}_i,d^i_{i+1}),\\
		X &= \mr(d^1,d^2_1,d^2,d^3_2,\dots,d^n,d^1_n),
	\end{align}
	where $X$ are the variables associated to the vertices of $\qg$ and $X^{i,i+1}$ are the variables associated to the faces of $\qg$. Due to the symmetries of the cross-ratio we can rearrange so that
	\begin{align}
		1 + X^{i,i+1} = \cro(d^{i+1},d^{i+1}_i,d^i,d^i_{i+1}) = \cro(d^i,d^i_{i+1},d^{i+1},d^{i+1}_i).
	\end{align}
	Now we recall that cluster variables are in the Ising subvariety if and only if
	\begin{align}
		(X)^2\prod_{i=1}^n(1+X^{i,i+1}) = 1,
	\end{align}
	at every vertex of the quad-graph, see Definition \ref{def:isingsub}. If we insert the expressions for $X$ and $(1+X^{i,i+1})$ from above and cancel terms then the last equation is equivalent to Carnot's Theorem.\qed	
}

The consistency of Darboux maps was shown when they were introduced \cite{schieflattice}. In our setup the consistency follows as an almost trivial corollary.

\begin{corollary}\label{cor:carnotconsistent}
	Carnot maps are consistent.
\end{corollary}
\proof{We have to check that if a map $d: E(\qg)\rightarrow \CP^n$ is a Carnot map then it is also a Carnot map after a cube flip. However, we know that the Ising subvariety is preserved by cube flips and thus the claim is an immediate consequence of Theorem \ref{th:carnotckp}.\qed}

\section{CQ-nets}\label{sec:cqnets}

Let us relate Carnot maps to another type of map that appeared in the literature. Doliwa introduced \emph{C-quadrilateral lattices} \cite{doliwacqnet}, a reduction of Q-nets that features the CKP equation in an affine gauge. In order to keep wording consistent in the thesis, we will call these nets CQ-nets. Doliwa gave a definition of CQ-nets via a property on 3-cubes in $\Z^N$. It is not obvious, how to immediately generalize this definition to quad-graphs, for reasons we explain below. Therefore we give Doliwa's definition and then after explaining the relation to Carnot maps, we explain how to understand CQ-nets on quad-graphs. Note that until then we also do not explain the precise genericity assumptions for Carnot maps.

\begin{definition}\label{def:cqnet}
	Let $H$ be a hyperplane in $\CP^n$ and $q: V(\Z^N) \rightarrow \CP^n$ be a Q-net. In every cube, consider the three lines $\ell^1,\ell^2,\ell^3$ that are the intersections of the planes of opposite faces, that is
	\begin{align}
		\ell^k = \spa(q,q_i,q_j) \cap \spa(q_k,q_{ik},q_{jk})
	\end{align}
	for any pairwise different indices $i,j,k\in \{1,2,\dots,n \}$. Then $q$ is called a \emph{CQ-net with respect to $H$} if the three intersections of the lines $\ell^k$ with $H$ are colinear.
\end{definition}

Doliwa himself gave a relation to Carnot maps via the \emph{Darboux system} \cite{doliwacqnet}. We do not discuss the Darboux system in this thesis. However, we show that there is an alternative relation between Carnot maps and CQ-nets. In the end, our relation to Carnot maps is quite possibly the same as Doliwa's, but with a considerable shortcut that makes it almost trivial.

\begin{lemma}\label{lem:cqandcarnot}
	A Q-net $q: V(\Z^3) \rightarrow \CP^n$ is a CQ-net with respect to $H$ if and only if $\sigma_H(q)$ is a Carnot map.
\end{lemma}
\proof{
	Let us denote $d = \sigma_H(q)$, which is a Darboux map. It is Carnot if the points of the hexagon $d^1,d^2_1,d^3_{12},d^1_{23},d^2_3,d^3$ are on a conic. By Pascal's theorem \cite{pascal} this is the case if and only if the intersection points of opposite sides of the hexagon are contained in a line. The six sides of the hexagon are the intersections of the six planes of every 3-cube of $q$ with $H$.\qed
}

Therefore the definition of a CQ-net on a quad-graph is clear, it is a 1-generic Q-net $q:  V(\Z^3) \rightarrow \CP^n$ such that $H$ is generic and $\sigma_H(q)$ is a Carnot map. The reason that one cannot immediately generalize Doliwa's definition of CQ-nets to quad-graphs, is that unlike for Menelaus' theorem and Carnot's theorem, there is to the best of our knowledge no obvious generalization of Pascal's theorem to larger polygons. Let us close this section with two obvious but relevant corollaries.

\begin{corollary}
	A Q-net is affine CKP if and only if it is a CQ-net.
\end{corollary} 
\proof{Because the section of a Q-net is a Darboux map and a Darboux map is projective CKP if and only if it is CKP, see Theorem \ref{th:carnotckp}.\qed}

\begin{corollary}
	CQ-nets are consistent.
\end{corollary}
\proof{Immediate from Corollary \ref{cor:carnotconsistent}.\qed}

\section{Doliwa compounds, BKP and the resistor subvariety}\label{sec:Doliwacomplex}

In this section we introduce a reduction of line compounds (see Section \ref{sec:lcco}) that corresponds to the resistor subvariety. We begin with some useful observations on line compounds. Recall that a line compound $l$ associates a $(d_v-2)$-space $S_v$ to every vertex $v\in V(\qg)$ with degree $d_v$, see Definition \ref{def:colc}. Let $f,f'\in F(\qg)$ be two faces adjacent to $v$ and denote by $S_v^{f,f'}$ the subspace of $S_v$ that is spanned by all points $l(f'')$ such that $f''$ is adjacent to $v$ but is neither $f$ nor $f'$.

\begin{definition}
	Let $l:F(\qg) \rightarrow \CP^n$ be a line compound. Then the \emph{focal points} of $l$ are the map $\fp:\vec{E}(\qg) \rightarrow \CP^n$ such that
	\begin{align}
		\fp(e) = l(f)l(f') \cap S_{v'}^{f,f'},
	\end{align}
	where $e=(v,v')$ and $f,f'$ are the two faces adjacent to $e$.
\end{definition}
Note that by $\vec{E}(\qg)$ we mean the oriented edges of $\qg$, that is 
\begin{align}
	\vec{E}(\qg) = \{(v,v') : \{v,v'\} \in E(\qg)\}.
\end{align}
This definition of focal points is analogous to the auxiliary points defined in \cite{vrc} to calculate the cluster variables. The introduction of focal points enables us to also introduce Laplace invariants in a very similar fashion as in the case of Q-nets (see Definition \ref{def:laplaceinvariant}).

\begin{definition}
	Let $l:F(\qg) \rightarrow \CP^n$ be a a line compound. The \emph{Laplace invariant} is a map $\Lambda: E(\qg) \rightarrow \C\setminus\{0\}$ such that
		\begin{align}
			\Lambda_e = -\cro(l(f),\fp(v'),l(f'),\fp(v)),
		\end{align}
		for every edge $e=(v,v')\in E(\qg)$, where $e^*=(f,f')$  and $v,f,v',f'$ appear in counterclockwise order in $\qg$.
\end{definition}

We now give a definition of Doliwa compounds on quad-graphs in terms of the Laplace invariants. We also give an interpretation in terms of a lattice multi-ratio equation below, and we will also discuss a geometric interpretation later on.

\begin{definition}\label{def:doliwacomplex}
	A \emph{Doliwa compound} is a line compound $l: F(\qg) \rightarrow \CP^n$ such that around every interior vertex of of $\qg$ the Laplace invariants multiply to 1.
\end{definition}
To best of the author's knowledge Doliwa compounds have not been explicitly defined in the literature. However, they appear implicitly in relation to K{\oe}nigs nets, as well as to another version of K{\oe}nigs nets defined by Doliwa. We will elaborate on these relations in this section. Let us first understand the relation to the resistor subvariety.

\begin{lemma}\label{lem:Doliwalaplacecluster}
	Let $l: F(\qg) \rightarrow \CP^n$ be a flip-generic line compound. For each non-boundary edge $e\in E(\qg)$ there is a unique face $f_e\in F(\pb)$ and vice versa. The Laplace invariants of $l$ and the projective cluster variables coincide, that is
	\begin{align}
		\Lambda_e = X_{f_e}
	\end{align}
	for all $e\in E(\qg)$.
\end{lemma}
\proof{
	Fix an edge $e\in E(\qg)$. We can always perform resplits at the white vertices of $\pb$ that correspond to the adjacent vertices of $\qg$ such that the face $f_e$ has degree four, because $l$ is flip-generic. Recall that resplits do not change the projective invariants $X_{f_e}$. Then the claim at $e$ is equivalent to Lemma \ref{lem:projclusterviadistances}.\qed
}

\begin{figure}
	\small
	\begin{tikzpicture}[scale=2.5] 
	\node[bvert] (v1)  at (-.87,-.5) {};
	\node[bvert] (v2) at (.87,-.5) {};
	\node[bvert] (v3) at (0,1) {};
	\node[bvert] (v12)  at (0,-1) {};
	\node[bvert] (v23) at (.87,.5) {};
	\node[bvert] (v13) at (-.87,.5) {};
	\node[bvert] (v) at (0,0) {};
	\draw[line width=2pt, gray,dashed]
		(v1) -- (v12) -- (v2) -- (v23) -- (v3) -- (v13) -- (v1)
		(v) edge (v12) edge (v13) edge (v23)
	;		
	\node[wvert,label={below:$l^{12}$}] (e1) at ($(v)!.5!(v1)$) {};
	\node[wvert,label={right:$\fp^2$}] (f12) at ($(v)!.5!(v12)$) {};
	\node[wvert,label={above right:$l^{23}$}] (e2) at ($(v)!.5!(v2)$) {};
	\node[wvert,label={above:$\fp^3$}] (f23) at ($(v)!.5!(v23)$) {};
	\node[wvert,label={left:$l^{13}$}] (e3) at ($(v)!.5!(v3)$) {};	
	\node[wvert,label={below left:$\fp^1$}] (f13) at ($(v)!.5!(v13)$) {};	
	\coordinate (e12) at ($(v2)!.5!(v12)$) {};
	\coordinate (e13) at ($(v3)!.5!(v13)$) {};
	\coordinate (e21) at ($(v1)!.5!(v12)$) {};
	\coordinate (e23) at ($(v3)!.5!(v23)$) {};
	\coordinate (e31) at ($(v1)!.5!(v13)$) {};
	\coordinate (e32) at ($(v2)!.5!(v23)$) {};
	
	\node[bvert] (c12) at ($(v)!.5!(f12)$) {};
	\node[bvert] (c23) at ($(v)!.5!(f23)$) {};
	\node[bvert] (c13) at ($(v)!.5!(f13)$) {};

	\draw[-] (v) -- (e1) -- (v1) (v) -- (e2) -- (v2) (v) -- (e3) -- (v3);
	\draw[-] (v12) -- (f12) -- (c12) -- (e1) (c12) -- (e2);
	\draw[-] (v23) -- (f23) -- (c23) -- (e2) (c23) -- (e3);
	\draw[-] (v13) -- (f13) -- (c13) -- (e3) (c13) -- (e1);

	\end{tikzpicture}
	\caption{Labeling and the graph $\pb$ for three quads (dashed, thick) of a line compound.}
	\label{fig:Doliwacomplexlabels}
\end{figure}
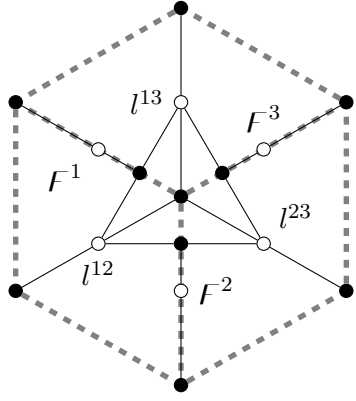

\begin{theorem}\label{th:projcomplex}
	A line compound is projective BKP if and only if it is a Doliwa compound.
\end{theorem}
\proof{
	Immediate consequence of the Definition \ref{def:resistorsub} of the resistor subvariety, Definition \ref{def:doliwacomplex} of Doliwa compounds and Lemma \ref{lem:Doliwalaplacecluster}.\qed
}

\begin{corollary}\label{cor:Doliwaconsistent}
	Doliwa compounds are consistent.
\end{corollary}
\proof{We have to check that if a map $l: E(\qg)\rightarrow \CP^n$ is a Doliwa compound then it is also a Doliwa compound after a cube flip. However, we know that the Resistor subvariety is preserved by cube flips and thus the claim is an immediate consequence of Theorem \ref{th:projcomplex}.\qed}

Therefore we can also think of Doliwa compounds as defined on the $\Z^N$ lattice. There is a simple lattice characterization of Doliwa compounds. Recall that in a line compound every cube $C(\Z^N)$ corresponds to a line in $\CP^n$.

\begin{lemma}\label{lem:Doliwacomplexlattice}
	Let $N\in \N, N\geq 3$ and let $l:F(\Z^N) \rightarrow \CP^n$ be a Doliwa compound. Then for every 3-cube of $\Z^N$ the six intersection points are in involution, that is
	\begin{align}
		\mr(l^{12},l^{23},l^{13},l^{12}_3,l^{23}_1,l^{13}_2) = \mr(l^{12},l^{23},l^{13},\fp^3,\fp^1,\fp^2) = -1,
	\end{align}
	holds, see Figure \ref{fig:Doliwacomplexlabels} for the labeling.
\end{lemma}
\proof{
	Consider Figure \ref{fig:Doliwacomplexlabels}. Note that due to the definition of line compounds and the definition of focal points, the three focal points before a cube-flip correspond to the three points of the line compound after the cube-flip, that is
	\begin{align}
		l^{23}_1 = \fp^1, \quad l^{13}_2 = \fp^2, \quad l^{12}_3 = \fp^3.
	\end{align}
	The product of the three cross-ratios associated to the three faces of $\pb$ have to multiply to $-1$ due to the definition of Doliwa compounds. At the same time the product of the three cross-ratios is also the multi-ratio in the claim, which finishes the proof. \qed
}

Of course, on $\Z^3$ there is no difference between line complexes and line compounds, thus we can also define Doliwa complexes on $\Z^3$ via Lemma \ref{lem:Doliwacomplexlattice}. However, by definition only Doliwa compounds are multi-dimensionally consistent. Unlike for Doliwa compounds, it is not clear if there is a version of multi-dimensional consistency for Doliwa complexes. A natural candidate is that Doliwa complexes on $\Z^N$ are line complexes on $\Z^N$ such that on every line the $2N$ intersection points are in involution, but we did not investigate whether such maps exist.

Let us turn to a geometric interpretation of Doliwa compounds.

\begin{theorem}\label{th:Doliwageometry}
	Let $l:F(\qg)\rightarrow \CP^n$ be a Doliwa compound. Let $v\in V(\qg)$ be an interior vertex with neighbours $v_1,v_2,\dots,v_{d_v}$  and let $f_1,f_2,\dots, f_{d_v}$ be the faces adjacent to $v$ such that $(v,v_i)=(f_i,f_{i+1})^*$ for every cyclic index $i$. Assume $l(f_1), l(f_2),\dots, l(f_{d_v})$ span a $(d_v-2)$-space. If $d_v\geq 4$ and $d_v\in 2\Z$ then the $2d_v$ points
	\begin{align}
		\fp(v,v_1),\fp(v,v_2),\dots,\fp(v,v_{d_v}) \quad \mbox{ and }\quad \fp(v_1,v),\fp(v_2,v),\dots,\fp(v_{d_v},v)\label{eq:focalmr}
	\end{align}
	are on a common quadric. If $d_v\geq 5$ and $d_v\in 2\Z+1$ then the two lines $l(f_1)l(f_2)$ and $\fp(v,v_1)\fp(v,v_2)$ intersect the $(d_v-3)$-space spanned by $\fp(v,v_3),\fp(v,v_4),\dots,\fp(v,v_{d_v})$ in a point.
\end{theorem}
\proof{
	Consider the even case first. We view the $d_v$ points $l(f_i)$ as the vertices of a $d_v$-gon $\mathfrak P$, and the $2d_v$ focal points as marked points on the sides of $\mathfrak P$. First we prove that
	\begin{align}
		\mr(l(f_1),\fp(v_1,v),l(f_2),\fp(v_2,v),\dots,\fp(v_{d_v},v)) = 1.
	\end{align}
	Assume the vertices of $\mathfrak P$ span the maximal possible dimension, that is a $(d_v-2)$-space. Inside this space choose homogeneous lifts for the involved points and denote the lifts with a hat. Then we can write
	\begin{align}
		\hat \fp(v_i,v) = \det(\hat l(f_{i+1}),\hat l(f_{i+2}),\dots,\hat l(f_{i-1})) \hat l(f_i) - (-1)^{d_v}\det(\hat l(f_{i+2}),\hat l(f_{i+3}),\dots,\hat l(f_{i})) \hat l(f_{i+1}).\label{eq:genfpdet}
	\end{align}
	The multi-ratio in Equation \eqref{eq:focalmr} can be expressed by the alternating ratio of the coefficients of Equation \eqref{eq:genfpdet}. In this expression all determinants cancel and what remains is a sign factor $(-1)^{d_v+d_v^2}$ which equals $1$. On the other hand, we prove that
	\begin{align}
		\mr(l(f_1),\fp(v,v_1),l(f_2),\fp(v,v_2),\dots,\fp(v,v_{d_v})) = (-1)^{d_v}. \label{eq:mrDoliwadecorated}
	\end{align}
	The multi-ratio on the left-hand side appears as a alternating ratio of edge-weights appearing in $\pb$ around $v$. It equals $(-1)^{d_v}X_{f_1}X_{f_2}\cdots X_{f_{d_v}}$, which is $(-1)^{d_v}$ due to Theorem \ref{th:projcomplex}. The claim of the theorem then follows from Carnot's theorem (see Theorem \ref{th:carnot}), because the two multi-ratios above multiply to 1 in the case of even $d_v$. Note that in Carnot's theorem the points of $\mathfrak P$ span maximal dimension. That assumption is not necessary for the multi-ratio condition to hold if the points are on a quadric, as is obvious from the proof of Carnot's theorem. On the other hand, a quadric in $\CP^{d_v-2}$ through $d_v$ points always exists, though it is not necessarily unique. Let us turn to the case of odd $d_v$. In this case the corresponding multi-ratio as in Equation \eqref{eq:mrDoliwadecorated} is $-1$. Therefore this is an instance of Menelaus' theorem (see Theorem \ref{th:genmenelaus}), but in dimension $d_v-2$ instead of $d_v-1$. However, this means that the polygon $l(f_1),l(f_2),\dots, l(f_{d_v})$ with marked points $\fp(v,v_1),\fp(v,v_2),\dots,\fp(v,v_{d_v})$ is the projection of a Menelaus configuration in dimension $d_v - 1$. The coincidence of intersection points in the claim of the theorem is a direct consequence of the fact that in a lift of the polygon do dimension $d_v - 1$ the points  $\fp(v,v_1),\fp(v,v_2),\dots,\fp(v,v_{d_v})$ only span a space of dimension $d_v-2$. \qed
}

Note that the focal points $\fp(v_i,v)$ in Theorem \ref{th:Doliwageometry} are points appearing as interior white vertices in the TCD map piece $T_v$ that represents the space $S_v$. If we would require the quadric to pass through all points appearing as interior white vertices in $T_v$, the quadric would be uniquely determined. This is because there are exactly $\binom{d_v}{2}-d_v$ interior white vertices. Together with the $d_v$ points $\fp(v,v_i)$ we are therefore requiring the quadric to pass through $\binom{d_v}{2}$ points. In general, a quadric is determined by $\binom{d_v}{2}-1$ points, which fits perfectly with the one equation we pose on the points.

In principle, the observation of Theorem \ref{th:Doliwageometry} is characterizing of Doliwa compounds, except for the case of interior vertices of degree $3$. In that case the points locally only span a 1-dimensional space, and therefore no incidence geometric characterization is possible.

\begin{figure}
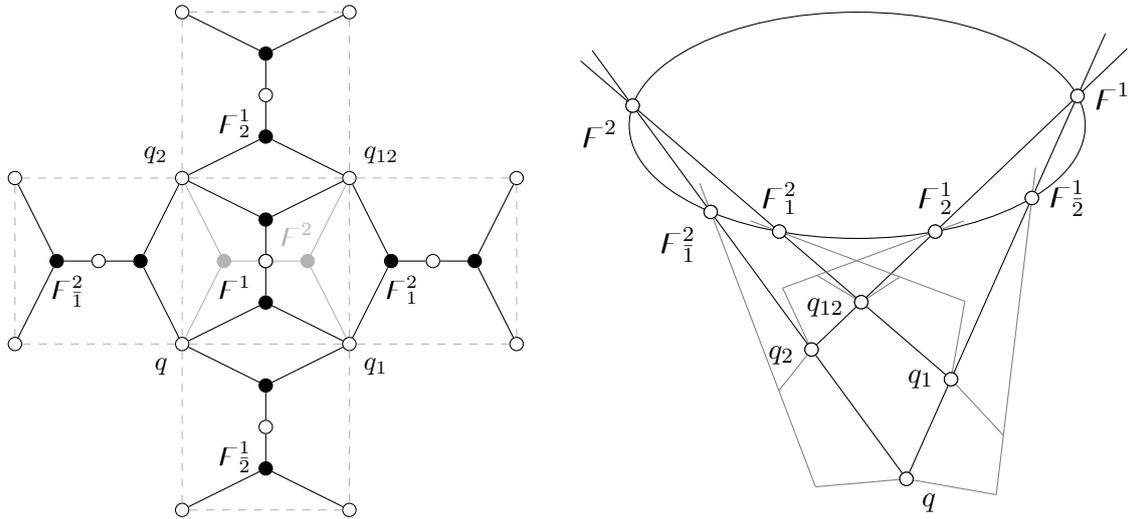


	\caption{Right: the conic that appears in Doliwa's K{\oe}nigs nets. Left: the VRC in which the projective cluster variables around $f^1$ multiply to 1. Note that the resplit does not affect the projective cluster variables.}
	\label{fig:doliwakonigsnet}
\end{figure}

The geometric characterization in the even case above is motivated by a relation between Doliwa compounds and a reduction of Q-nets that is due to Doliwa \cite{doliwakoenigs}.
\begin{definition}\label{def:Doliwaqnet}
	Let $q: \Z^2 \rightarrow \CP^n$ be a Q-net and let $\fp^1,\fp^2$ denote the focal points of $q$ in the two lattice directions. 
	Then we call $q$ a \emph{Doliwa-K{\oe}nigs net} if the six points
	\begin{align}		
		\fp^1_{\bar2},\fp^1,\fp^1_2,\fp^2_{\bar1},\fp^2,\fp^2_1
	\end{align}
	are on a conic, see Figure \ref{fig:doliwakonigsnet}.
\end{definition}

\begin{lemma}
	A $\Z^2$ slice of the focal Q-net of a Doliwa compound defined on $\Z^N$ is a Doliwa-K{\oe}nigs net.
\end{lemma}
\proof{
	We have discussed focal nets in Section \ref{sec:focalnets}, in particular how to take the focal Q-net of a line complex. The factorization property of the resistor subvariety carries over to $\Z^2$ slices of a Doliwa complex. The characterization in Definition \ref{def:Doliwaqnet} is equivalent to the geometric characterization of Doliwa complexes in Theorem \ref{th:Doliwageometry}.\qed
}

The fact that certain Laplace invariants multiply to 1 in Doliwa-K{\oe}nigs nets is already mentioned in \cite[Exercise 2.21]{ddgbook}.

\section{Diagonal intersection nets}

In a sense, the consequence of the last section is that Doliwa compounds are a -- and possibly the -- integrable extension of Doliwa nets to $\Z^N$. There is one more relation that has been pointed out in the DDG book for the $\Z^2$ case \cite{ddgbook}. We give it here in more generality for $\Z^N$, which gives further support to the consideration of Doliwa compounds defined on line compounds.
\begin{definition}
	The \emph{diagonal intersection net} $m: F(\qg) \rightarrow \CP^n$ of a Q-net $q: V(\qg) \rightarrow \CP^n$ consists of the intersection points of the diagonals of every quad.
\end{definition}

\begin{lemma}\label{lem:diagkoenigs}
	The diagonal intersection net $m: F(\qg) \rightarrow \CP^n$ of a K{\oe}nigs net $q: V(\qg) \rightarrow \CP^n$ is a Doliwa compound.
\end{lemma}
\proof{We give a proof based on the lattice interpretation of the maps. Due to Definition \ref{def:koenigs}, in every 3-cube of a K{\oe}nigs net all black vertices are in a plane $H$ and all white vertices are in a plane $H'$. Therefore all intersections of diagonals in a 3-cube are on the line $L=H\cap H'$. Consequently, the diagonal intersection net is a line compound. Moreover, the six points of $m$ on $L$ are the intersections of all six lines defined by the four points $q,q_{12},q_{23},q_{13}$. This is a complete quadrangle and it is a well known that the intersection of the lines of a complete quadrangle with a line are six points with multi-ratio $-1$. This proves the claim due to Lemma \ref{lem:Doliwacomplexlattice} that characterizes Doliwa compounds via multi-ratios.\qed
}

Note that the converse of this construction is also possible, see the thesis of Steinmeier \cite{steinmeierkoenigs}. Moreover, if we include some genericity assumptions the construction of a diagonal intersection net is a special case of a section (sections were introduced in Section \ref{sec:sections}). 
\begin{lemma}\label{lem:generickoenigsdiagonal}
	Let $U= [0,1,\dots,k]^3  \subset \Z^3 $ and let $q: U \rightarrow \CP^{3k}$ be a K{\oe}nigs net such that $q$ spans $\CP^{3k}$. Then there is a projective subspace $H\subset \CP^{3k}$ of codimension 2 such that the diagonal intersection map of $q$ is the section $\sigma_H(q)$.
\end{lemma}
\proof{First of all let us observe that the points along the coordinate axes are Cauchy-data for the Q-net. They span at most $\CP^{3k}$, and by the assumptions of the Lemma they do span all of $\CP^{3k}$. In each cube we have a line $L$ defined as intersection of the two planes $E,E'$ as in the proof of Lemma \ref{lem:diagkoenigs}. Now define $H$ as the span of all lines $L$. Moreover, if we just look at the span $H'$ of all lines $L$ in the cubes along the coordinate axes, then in fact $H'=H$. This is because in any quad of cubes the line in the fourth cube intersects two of the other lines already. But the dimension of $H'$ is $(3k-2)$ and so is the dimension of $H$. Thus in general position $m=\sigma_H(q)$.\qed
}

We have limited the statement of the lemma to cubical domains in order to avoid technical arguments in the proof. In principal however we claim that one can consider any quad-graph and replace the dimension $3k$ with the maximal dimension of the TCD that is associated to the quad-graph.

\section{Schief maps}\label{sec:pardm}\label{sec:schiefmaps}

There is also a notion of a BKP map introduced by Schief \cite{schieflattice}, which are Darboux maps that carry a BKP structure. We call these maps Schief maps instead, because from our point of view the name BKP map is not concise. In particular, in higher projective dimensions a Darboux map can be BKP in many ways. We do show that Schief maps are affine BKP, but they are only a special case of affine BKP Darboux maps.

\begin{definition}\label{def:paralleldm}
	A Darboux map $d: E(\qg) \rightarrow \CP^n$ is a \emph{Schief map} with respect to a hyperplane $H$ if in some affine chart $\C^n$ with $H$ at infinity every quad is mapped to a parallelogram.
\end{definition}
Note that a parallelogram is an affine notion and thus it is not surprising that the next statement is affine as well.

\begin{theorem}\label{th:schiefbkp}
	A Schief  $d: E(\qg) \rightarrow \CP^n$ map is affine BKP with respect to $H$.
\end{theorem}
\proof{Consider the affine cluster variables, that is the star-ratios of $d$. The star-ratios for a Darboux map sit at the edges of the quad-graph. Each star-ratio is the quotient of the two dilations of the incident parallelograms. Thus if we take the product of the star-ratios around a vertex of the quad-graph the result is 1. Therefore the affine cluster variables are in the Resistor subvariety and the Schief map is affine BKP.\qed}

Note however that not every Darboux map that is affine BKP is also a Schief map. The relation to Doliwa complexes is as follows.

\begin{corollary}
	Let $d: E(\qg) \rightarrow \CP^n$ be a 1-generic Schief map such that $H$ is generic. Then the section $\sigma_H(d)$ is a Doliwa complex.
\end{corollary}
\proof{In Section \ref{sec:sectioncluster} we proved that the projective cluster structure of the section $\sigma(T)$ with hyperplane $H$ is the affine cluster structure of $T$ with respect to $H$. Moreover the section of a Darboux map is a line complex. A Schief map is affine BKP with respect to infinity. As a consequence, the section with infinity is a line compound that is projective BKP and therefore a Doliwa complex.\qed}

\begin{remark}
	In the specific generic setup of Lemma \ref{lem:generickoenigsdiagonal} any section with a hyperplane $H'\supset H$ is a Schief map. We leave the proof of this statement as an exercise.
\end{remark}

\section{Reciprocal figures}\label{sec:reciprocal}

King and Schief \cite{kstetraoctacubo} have also encountered maps in $\C$ that are defined on faces of $\Z^3$ and satisfy lattice equation which states that the multi-ratio of six points equals -1. They interpret these as Darboux maps in $\CP^1$ that satisfy this additional 6-point equation and show that they are accompanied by a BKP equation. However, we have mentioned in Section \ref{sub:cponelc} that Darboux maps and line complexes in $\CP^1$ coincide. We intend to argue now that it is more natural to consider these maps as line complexes. In fact we show that these maps are projective BKP and thus are Doliwa compounds in $\CP^1$. We limit ourselves to considering $\Z^3$ instead of arbitrary quad-graphs to simplify the exposition. Because we restrict ourselves to $\Z^3$ we consider line complexes and line compounds to coincide in this section.

\begin{definition}\label{def:reciprocal}
	A \emph{reciprocal figure} map $r: V(\Z^3) \rightarrow \R^2$ is a Q-net such that in every quad the two diagonals are parallel. Let $r=(r^x,r^y)$ and define the slope $\phi^{ij} \in \R$ in a quad as 
	\begin{equation}
		\phi^{ij} = \frac{r^y_{ij}-r^y}{r^x_{ij}-r^x}.\qedhere
	\end{equation}
\end{definition}

\begin{theorem}[\cite{kstetraoctacubo}]\label{th:reciprocal}
	The slopes $\phi$ of a reciprocal figure map are a Darboux map $\phi: E(\Z^3) \rightarrow \R \subset \RP^1$ and satisfy the multi-ratio equation
	\begin{equation}
		\mr(\phi^{ij},\phi^{jk},\phi^{ik},\phi^{ij}_k,\phi^{jk}_i,\phi^{ik}_j) = -1.
	\end{equation}
	for $\{i,j,k\}=\{1,2,3\}$.
\end{theorem}

A proof that uses explicit calculation can be found in the article by King and Schief \cite{kstetraoctacubo}. Note that the slopes are defined in Definition \ref{def:reciprocal} on the faces of quads, while the corresponding Darboux map in Theorem \ref{th:reciprocal} are defined on edges of $\Z^3$. This is possible because we can think of the quad-graph in Definition \ref{def:reciprocal} as living on $(\Z^3)^*$.

We proceed to incorporate reciprocal figure maps into our framework and then the theorem will follow from observations that relate K{\oe}nigs nets and Doliwa complexes that we made before. In the following, if we view $\R^2$ as part of $\RP^2$ then parallel diagonals means that each pair of diagonals intersects on a given line $L$ (at infinity).
\begin{definition}\label{def:reciprocallift}
	Consider V($\Z^3$) as bicolored into white and black vertices. Define the lift $\hat r: V(\Z^3) \rightarrow \RP^3$ of a reciprocal figure map $r = (r_x,r_y)$ by
	\begin{align}
		\hat r &= [r_x,r_y,1,1] \mbox{ for black vertices,}\\
		\hat r &= [r_x,r_y,0,1] \mbox{ for white vertices.}\qedhere
	\end{align}
\end{definition}

Note that the brackets $[\cdot]$ denote the projectivization. We consider $\hat r$ as the homogeneous lift of a lift of $\hat r$ in one of two hyperplanes of $\RP^3$, one hyperplane for the white vertices and one hyperplane for the black vertices. Projecting $\hat r$ onto its first two coordinates is clearly a projection onto $\R^2$ that takes $\hat r$ back to $r$.

\begin{lemma}
	The lift $\hat r$ of a reciprocal figure is a K{\oe}nigs net.
\end{lemma}
\proof{First of all, because the diagonals $rr_{ij}$ and $r_ir_j$ are parallel by definition so are the homogeneous lifts of $\hat r\hat r_{ij}$ and of $\hat r_i\hat r_j$ as chosen in Definition  \ref{def:reciprocallift}. Therefore $\hat r$ is a Q-net. Moreover, by the definition of the lift all black vertices of every cube of $\hat r$ are contained in a plane, and so are all white vertices. Therefore $\hat r$ is K{\oe}nigs.\qed
}

As we discussed in Section \ref{sec:konigs}, a K{\oe}nigs net is projective BKP and it is straightforward to show that this BKP structure coincides with the first BKP structure that King and Schief discussed in their article \cite{kstetraoctacubo}. Next we give an interpretation for the other BKP structure that they found in terms of slopes (already indicated in Theorem \ref{th:reciprocal}).

\begin{theorem}\label{th:recDoliwacomplex}
	Consider the line $L = [1,\psi,0,0], \psi\in \R \cup \infty$ and the section $\sigma_L(\hat r)$. Then $\sigma_L(\hat r)$ is a Doliwa complex and $\psi$ are the slopes $\phi$ from Definition \ref{def:reciprocal}.
\end{theorem}
\proof{Because $L$ has codimension 2 the section $\sigma_L(\hat r)$ is a line complex. Moreover, $L$ is the intersection of the plane that contains all white vertices of $\hat r$ with the plane that contains all black vertices of $\hat r$. The intersection of the plane spanned by quad $(r,r_i,r_j,r_{ij})$ is 
\begin{align}
	r- r_{ij} = [r^x_{ij} -r_x, r^y_{ij}-r_y,0,0].
\end{align}
If we bring this expression into the affine parametrization $[1,\psi,0,0]$ of $L$ we see that indeed $\psi = \phi$. Additionally, the section with $L$ is also the diagonal intersection map of $\hat r$ and therefore a Doliwa complex by Lemma \ref{lem:diagkoenigs}.\qed
}

In conclusion, reciprocal figures are equipped with two BKP structures. One structure because one can lift reciprocal figures to a K{\oe}nigs net, which is projective BKP. The other structure is because the section of the lift with the line $L$ from Theorem \ref{th:recDoliwacomplex} is a Doliwa complex, which is also projective BKP.


\chapter{Reductions of TCD maps via bilinear forms}\label{cha:bilinear}
\label{sec:tcdbilinear}

In this section we study TCD maps that satisfy not only linear constraints, as they do by definition, but also additional quadratic constraints. Before we give a small refresher of quadrics and null-polarities, let us motivate this section with a lemma and a theorem, both found by Doliwa \cite{doliwaqnetsinquadrics}.

\begin{lemma}\label{lem:qnetquadcube}
	If seven points of a Q-net cube are contained in a quadric $B$, then so is the eighth.
\end{lemma}
\begin{theorem}\label{th:qnetquadric}
	If the Cauchy data of a Q-net is inscribed in a quadric, then so is the entire Q-net.
\end{theorem}
Note that Doliwa formulated Theorem \ref{th:qnetquadric} for $\RP^n$, but the theorem also holds for $\CP^n$ without any modifications to the proof. The proof relies on the theorem for the associated point, which can be found in the DDG book \cite{ddgbook}, also for $\RP^n$. The proof also holds without modification for $\CP^n$.

We emphasize the importance of Theorem \ref{th:qnetquadric} with a quote of Doliwa \cite{doliwaqnetsinquadrics}.
\begin{corollary}
	``The above result can be obviously generalized to quadrilateral lattices in spaces obtained by intersection of many quadric hypersurfaces. Since the spaces of constant curvature, Grassmann manifolds and Segré or Veronese varieties can be realized in this way, the above results can be applied, in principle, to construct integrable lattices in such spaces as well.''
\end{corollary}

Another source of power of these two theorems stems from the fact that many classical geometries (Möbius, Laguerre, Lie, Plücker) can be formulated via a quadric in a projective space and the set of projective transformations that preserve that quadric. Especially in discrete surface theory, this formulation enables us to study discretizations that are invariant with respect to the right group of transformations \cite{bsorganizing}. We list three particular examples, where advantage of this formulation was taken.
\begin{enumerate}
	\item Circular nets in $\RP^n$ are objects of Möbius geometry. They correspond to Q-nets in $\S^n \subset \RP^{n+1}$. From the projective point of view, $\S^n$ is just some quadric with signature $(n+1,1)$ \cite{bobenkocircular, cdscircular}. We discuss the case of circular nets in a sphere in $\RP^3$ in Section \ref{sub:inscribedqnets} and will discuss more general circular nets in Section \ref{sub:circularq}.
	\item Conical nets in $\RP^3$ are objects of Laguerre geometry. They correspond to Q-nets in the Blaschke cylinder, which in the 3D case is the degenerate quadric of signature $(3,1,1)$ \cite{pwconical}.
	\item A-nets in $\RP^3$ are objects of projective geometry. However, we can lift the lines of an A-net to the Plücker quadric, which is a quadric in $\RP^5$ with signature $(3,3)$. Therefore we also view A-nets as objects of Plücker geometry. More precisely, A-nets correspond to isotropic line complexes in the Plücker quadric \cite{doliwaanetspluecker}. We will discuss A-nets in more detail in Section \ref{sub:anets}.
\end{enumerate}

Thus circular nets, conical nets and A-nets are all cases of TCD maps that are in non-general position with respect to a quadric. Motivated by this we show in this chapter that certain Q-nets, Darboux maps and line complexes that satisfy additional quadratic constraints have an associated cluster structure with variables in the resistor or Ising subvarieties. This in turn relates the maps to the BKP or CKP equation. 

Let us continue with a short refresher on quadrics and null-polarities. We call a map 
\begin{align}
	m: \C^n \times \C^n \rightarrow \C
\end{align}
a \emph{bilinear form} if it is linear in both arguments. The \emph{polar space} to a point $r\in \C^n$ is the set
\begin{align}
	\{x\in \C^n : c(r,x) = 0 \}.
\end{align}
A bilinear form $c$ is \emph{degenerate} if there is an $x\in \C^n$ such that its polar space is the whole space, and non-degenerate otherwise. We call a bilinear form \emph{symmetric} (in which case we denote it by $b$) if
\begin{align}
	b(x,y) = b(y,x),
\end{align}
for all $x,y\in \C^n$. We call a bilinear form \emph{anti-symmetric} (in which case we denote it by $\omega$) if
\begin{align}
	\omega(x,y) = -\omega(y,x),
\end{align}
for all $x,y\in \C^n$. In both cases we note that if $b(x,y) = 0$ (resp. $\omega(x,y) = 0$) then also $b(y,x) = 0$ (resp. $\omega(y,x) = 0$). Thus $x$ being in the polar space of $y$ is a symmetric relation in both cases. We say $x$ is \emph{polar} to $y$ if $x$ is in the polar space of $y$. Moreover, if $x,x' \in \C^n$ are lifts of the same projective point $\hat x \in \CP^{n-1}$, that is there is $\lambda \in \C$ such that $x'=\lambda x$, then $b(x,y)=0$ (resp. $\omega(x,y)=0$) is equivalent to $b(x',y)=0$ (resp. $\omega(x',y)=0$) for all $y\in \C^n$. Therefore polar spaces are a projective notion. Therefore we abuse notation to some extent and also write expressions $b(\hat x,\hat y)=0$ for points $\hat x,\hat y\in \CP^{n-1}$, as this expression is well-defined independent of the choice of lifts. We call a projective space $I\subset \CP^{n-1}$ an \emph{isotropic space} if $b(\hat x,\hat y)=0$ (resp. $\omega(\hat x,\hat y)=0$) for all $\hat x,\hat y\in I$. 

For a bilinear form $b$ we call the set of points 
\begin{align}
	B=\{\hat x\in \CP^{n-1} : b(\hat x,\hat x)=0\}
\end{align}
the \emph{quadric} $B$ associated to $b$. The term ``quadric'' comes from the fact that $b(x,x)=0$ is a quadratic equation for the components of $x\in \C^n$. We say $B$ is degenerate if $b$ is degenerate. If two polar points $\hat x,\hat y$ lie in a quadric $B$ then by linearity also all points on the line $\hat x \hat y$ are in $B$. For points $\hat x$ on $B$ we call the polar spaces \emph{tangent spaces}.

For an anti-symmetric form $\omega$, it is true for all points $\hat x\in \CP^{n-1}$ that $\omega(\hat x,\hat x) = 0$. Therefore there is no analogue to a quadric in this case. The map $x\mapsto \omega(x, \cdot)$ is a map from the primal space $\C^n$ to the dual space $(\C^n)^*$ and such a map is called a \emph{polarity}. Because $\omega(\hat x, \hat x) = 0$ for all $\hat x\in \CP^{n-1}$, $\omega$ is called a null-polarity. Moreover, if $\omega(\hat x,\hat y)=0$ for two different points in $\CP^{n-1}$, then $\omega(\hat x,\hat z) = \omega(\hat y,\hat z) = 0$ for all points $\hat z$ on the line $\hat x \hat y$. Thus even though no quadric exists for anti-symmetric forms, isotropic lines are still distinguished. In particular, for every point $\hat x\in \CP^{n-1}$ all lines through $\hat x$ that are also contained in the polar plane of $\hat x$ are isotropic.

\section{Q-nets in $\CP^2$ inscribed in conics}\label{sec:inscribedinconics}

\begin{definition}
	Let $A,A'$ be two different points in $\CP^2$ and let $\qg$ be a quad-graph. We say a Q-net $q: V(\qg)\rightarrow \CP^2$ is \emph{2-conical} if the line $AA'$ is generic (in the sense of Definition \ref{def:generichyperplane}) and if for every quad of $\qg$ the six points $A,A',q,q_1,q_2,q_{12}$ are on a non-degenerate conic.
\end{definition}

Note that a conic is a quadric in $\CP^2$. A natural question is whether cube flips of Q-nets preserve 2-conical nets. To answer this, let us look at a generalization of Miquel's theorem \cite{miquel}.

\begin{theorem}[The complex projective theorem of Miquel]\label{th:projmiquel}
	Let $A,A'$ be two different points in $\CP^2$. Assume the three sets of points
	\begin{align}
		\{A,A',q,q_1,q_2,q_{12}\},\quad \{A,A',q,q_1,q_3,q_{13}\} \mbox{  and  } \{A,A',q,q_2,q_3,q_{23}\}
	\end{align}
	are on a non-degenerate common conic each and that none of the occurring points is on the line $AA'$. Then there is a unique point $q_{123}$ in $\CP^2$ such that
	\begin{align}
		\{A,A',q_3,q_{13},q_{23},q_{123}\},\quad \{A,A',q_2,q_{12},q_{23},q_{123}\} \mbox{  and  } \{A,A',q_1,q_{12},q_{13},q_{123}\}
	\end{align}
	are on a common conic each.
\end{theorem}
\proof{
	Consider $\CP^2$ as a hyperplane $H$ in $\CP^3$ and choose two lines $\ell_A,\ell_{A'}$ through $A$ resp. $A'$ that intersect in an additional point $C$ not in $H$. Now, choose a non-degenerate quadric $B$ (and a corresponding bilinear form $b$) such that $\ell_A,\ell_{A'}$ are isotropic. Define the lifts $\hat q,\hat q_1, \hat q_2, \hat q_3, \hat q_{12}, \hat q_{13}, \hat q_{23}, \hat q_{123} \in B \setminus \{C\}$ of the corresponding points such that the central projection from $C$ onto $H$ recovers the original points $q,q_1,q_2,q_3, q_{12},  q_{13},  q_{23},  q_{123}$. Now we claim that the lifts $\hat q,\hat q_1, \hat q_2, \hat q_{12}$ are on a plane if and only if their projections and $A$ and $A'$ are on a common conic. The points $\hat q,\hat q_1, \hat q_2$ determine a plane $E$ and a conic $\hat c = E\cap B$. We distinguish two cases. In the first case $E$ contains $C$. As a consequence $q,q_1,q_2$ are on a line. The projections are only on a common conic if $q_{12}$ is on the degenerate conic $q_1q_2 \cup AA'$ and therefore by assumption on the line $q_1q_2$. This is equivalent to $\hat q_{12}$ being in $E$. In the second case $E$ does not contain $C$. We consider the conic $c$ that is the projection of $\hat c$ onto $E$. The conic $c$ has to contain $A,A'$ because $E$ intersects the two lines $\ell_A, \ell_{A'}$ in two points in $\hat c$, because the two lines are isotropic. The only points in $B$ projected onto $c$ are the points in $\hat c$ or in $\ell_A,\ell_{A'}$ because the quadric is non-degenerate. Therefore $\hat q_{12}$ is on $\hat c$ and thus in $E$ if and only if $q_{12}$ is on $c$. Therefore in this lift to $\CP^3$, the statement of the theorem is equivalent to Lemma \ref{lem:qnetquadcube}, which guarantees that there is a unique point $\hat q_{123}$. The projection of $\hat q_{123}$ to $H$ is the unique point $q_{123}$.\qed
}

The classical Miquel's theorem refers to the special case where we identify $\R^2$ with an affine chart  of $\RP^2$. Circles are conics that contain the two points with homogeneous coordinates $(\pm i,1,0)$ in the complexification $\CP^2$ of $\RP^2$ \cite[$\mathrm I$ and $\mathrm J$ in Section 1.6]{rgbook}. These two points are called \emph{imaginäre Kreispunkte} in the German literature, and for the lack of an English expression we stick with `imaginäre Kreispunkte'. The statement of Miquel's theorem is for Q-nets such that the four points of each quad are contained in a circle. This is the special case where $A,A'$ are the imaginäre Kreispunkte. In this special case one can choose the quadric $B$ to be a sphere and the projection is the standard stereographic projection.

A consequence of Theorem \ref{th:projmiquel} is that if the Cauchy data of a Q-net in $\CP^2$ are 2-conical then so is the whole Q-net. Now that we have established that 2-conical Q-nets are a consistent reduction, we give our main theorem on them.

\begin{theorem}\label{th:miquelbkp}
	Let $q: V(\qg) \rightarrow \CP^2$ be a 1-generic and 2-conical Q-net. Then the section $\sigma_{AA'}(q)$ is an affine BKP Darboux map with respect to $A$ and also with respect to $A'$. 
\end{theorem}

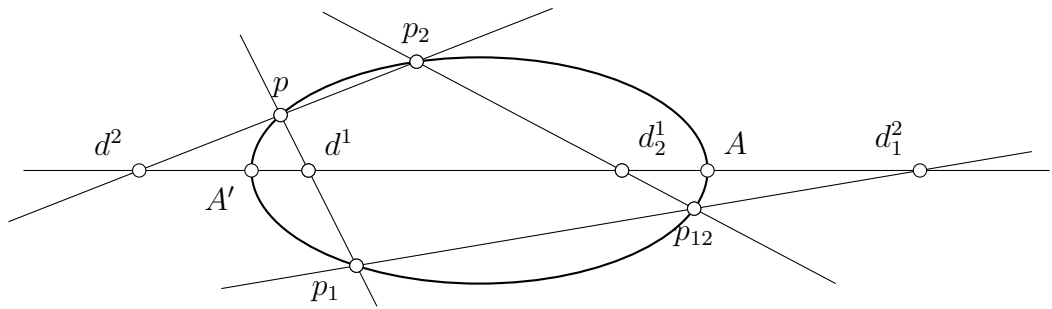
\begin{figure}
	\begin{tikzpicture}[scale=.75]
			\draw[black,thick,name path=conic] (0,0) circle[x radius = 4, y radius = 2];
			\coordinate (f1) at (-5,4);
			\coordinate (l1) at (-8,0);
			\coordinate (l2) at (10,0);
			\draw[name path=l, -] (l1) -- (l2);
			\path [name intersections={of = conic and l}];
			\node[wvert,label=above right:$A$] (a1)  at (intersection-1) {};
			\node[wvert,label=below left:$A'$] (a2)  at (intersection-2) {};
			
			\node[wvert,label=above right:$d^1$] (d1)  at (-3,0) {};
			\node[wvert,label=above right:$d^1_{2}$] (d12)  at (2.5,0) {};			
			\draw[name path=k1, -] ($(f1)!.4!(d1)$) -- ($(d1)!-.6!(f1)$);
			\draw[name path=k12, -] ($(f1)!.3!(d12)$) -- ($(d12)!-.5!(f1)$);
			
			\path [name intersections={of = conic and k1}];
			\node[wvert,label=above:$p$] (p) at (intersection-1) {};
			\node[wvert,label=below left:$p_1$] (p1)  at (intersection-2) {};
			\path [name intersections={of = conic and k12}];
			\node[wvert,label=above:$p_2$] (p2)  at (intersection-1) {};
			\node[wvert,label=below:$p_{12}$] (p12)  at (intersection-2) {};
			
			\draw[name path=k2, -] ($(p)!-2!(p2)$) -- ($(p)!2!(p2)$);
			\draw[name path=k21, -] ($(p1)!-.4!(p12)$) -- ($(p1)!2!(p12)$);

			\path [name intersections={of = k2 and l}];
			\node[wvert,label=above left:$d^2$] (d2)  at (intersection-1) {};
			\path [name intersections={of = k21 and l}];
			\node[wvert,label=above left:$d^2_1$] (d21)  at (intersection-1) {};

			\foreach \nodo in {(p),(p2),(d1),(p1),(d12),(p12)}{			
				\node[wvert] at \nodo {};
			};
	\end{tikzpicture}
	\caption{Desargues' involution theorem}
	\label{fig:desarguesinvolution}
\end{figure}

In order to prove this theorem, we need another classical theorem, see the book of Coxeter \cite[Section 9.3]{coxeterprojective} for a reference and a proof.

\begin{theorem}[Desargues' involution thereom]\label{th:desarguesinvolution}
	Let $q,q_1,q_2,q_{12},A,A'\in \CP^2$ be six different points on a conic $c$ and let $\ell$ be the line $AA'$. Denote by $d^l_k$ the intersection point of the line $q_kq_{kl}$ with $\ell$ (see Figure \ref{fig:desarguesinvolution}). Then the multi-ratio equation
	\begin{align}
	\mr(d^1,d^2,A,d^1_{2},d^2_{1},A')=-1
	\end{align}
	holds.
\end{theorem}

\begin{figure}
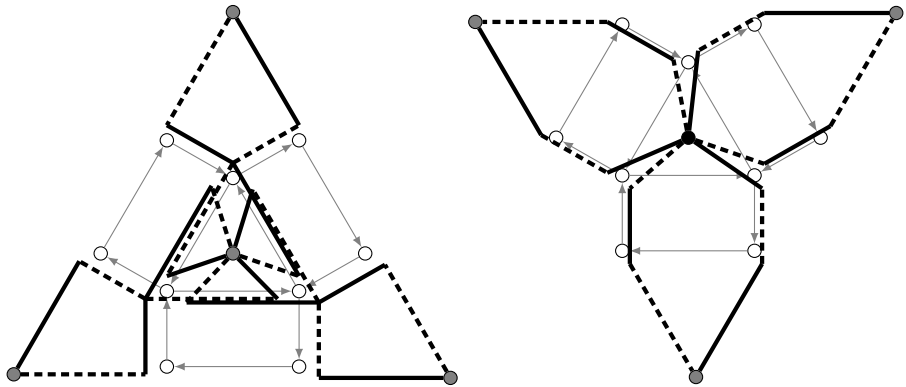

	
	\caption{The white points are the points of the section of a 2-conical Q-net. On the left the required factorization property for affine BKP with respect to $A$ (gray point). On the right three instances of Desargues' involution theorem with respect to $A$ (gray) and $A'$ (black).}
	\label{fig:miqueldesargues}
\end{figure}

\proof[Proof of Theorem \ref{th:miquelbkp}]{
	We visualize the proof in Figure \ref{fig:miqueldesargues} for three adjacent quads, but the proof readily generalizes to $m, m>3$ quads. Let $d=\sigma_{AA'}(q)$ be the Darboux map that is the section of $q$. The affine BKP requirement for the star-ratios is (we use Lemma \ref{lem:srviacr}), that
	\begin{align}
		\prod_{k=1}^m \frac{\cro(d^k,d^{k+1}_k,A,d^{k-1}_k)}{\cro(d^k,d^{k-1},A,d^{k+1})} = 1
	\end{align}
	holds. On the other hand, $m$ applications of Desargues' involution theorem yield that
	\begin{align}
		(-1)^m\prod_{k=1}^m \mr(d^1,d^{k+1}_k,A,d^k_{k+1},d^{k+1},A') = 1
	\end{align}
	holds. After canceling terms in the last two expressions they both yield
	\begin{align}
		\mr(d^k,d^{k+1}_k,A,d^k_{k+1},d^{k+1},d^{k+2}_{k+1},A,d^{k+1}_{k+2},\dots,A,d^m_{k}) = 1
	\end{align}
	and are therefore equivalent.\qed
}

\begin{remark}
	As explained, if we choose the two distinguished point of a 2-conical net to be the imaginäre Kreispunkte, then all conics in the 2-conical net are actually circles. In Section \ref{sec:cptemb} we will call these maps \emph{circle patterns}, but we will consider them as maps to $\CP^1$ instead of $\RP^2$. We will also find two associated cluster structures, but with different combinatorics. Moreover, for circle patterns considered in $\RP^2$ the natural move is the cube-flip and consists of a sequence of mutations in the cluster structure. In $\CP^1$ on the other hand we will do another move that corresponds to a single mutation, and is more like an octahedral flip. Also, we have just shown that the affine cluster structure for circle patterns in $\RP^2$ is naturally in the resistor subvariety. The same is not true for the cluster structure obtained when considering circle patterns in $\CP^1$.
\end{remark}

The naive converse of Theorem \ref{th:miquelbkp} is not true, not every doubly affine BKP Darboux map is the section of a 2-conical net. With some additional assumptions, there is a lift though.

\begin{theorem}\label{th:liftofdoublebkp}
	Let $A,A'$ be two points in $\CP^1$. Let $d: E(\qg) \rightarrow \CP^1$ be a Darboux map such that in every quad
	\begin{align}
		\mr(d^1,d^2,A,d^1_{2},d^2_{1},A')=-1\label{eq:lifttotwoconical}
	\end{align}
	holds. Then there is a Q-net $q: V(\qg) \rightarrow \CP^2$ such that $q$ is 2-conical with respect to the two points $A,A'$ and such that $\sigma_{AA'}(q) = d$.
\end{theorem}
\proof{We discussed in Section \ref{sec:tcdextensions} how to construct lifts of TCD maps. Assume we know three points $q,q_1,q_2$ of a quad of the lift. Then these three points together with $A,A'$ determine the conic $c$ on which $q_{12}$ has to lie. Because we also know $d^1_2$ the point $q_{12}$ is determined because $q_{12}$ has to be in the intersection $c\cap d^1_2q_2$ and not be $q_2$. But as $q_{12}$ is now determined, the point $d^2_1$ cannot be arbitrary but is also uniquely determined by the lift. Moreover, $d^2_1$ is also uniquely determined by Equation \eqref{eq:lifttotwoconical}. By Desargues' involution theorem (see Theorem \ref{th:desarguesinvolution}) these are the same points. Thus the lift is consistent on a per-quad basis and we can propagate through the quads by following a li-orientation of $\qg$ as explained in Section \ref{sec:sweepsqg}. \qed
}

\begin{corollary}
	Theorem \ref{th:miquelbkp} also holds in the case of a Q-net $q: V(\qg)\rightarrow \CP^2$ such that all its points are contained in one conic $c$. Conversely, let $d: E(\qg) \rightarrow \CP^1$ be a Darboux map as in Theorem \ref{th:liftofdoublebkp}. Then a Q-net lift of $d$ can also be found such that the points of the lift are all contained in one conic.
\end{corollary}
The proof of the corollary proceeds without modification to the proofs of the two theorems referenced. This introduces an intriguing construction.

\begin{corollary}
	Let $q: V(\qg) \rightarrow \CP^2$ be a 1-generic (Definition \ref{def:gentcdmap}) and 2-conical net with respect to two different points $A,A'$. Let $d = \sigma_{AA'}(q)$ be the section that is a Darboux map. Choose a quad $(v,v_1,v_2,v_{12})$ in $\qg$ and fix the conic $c$ through the points of that quad and $A,A'$. Then there is a unique Q-net $q': V(\qg) \rightarrow \CP^2$ such that 
	\begin{enumerate}
		\item the points of $q'$ are contained in $c$, that is $q'(V(\qg)) \subset c$,
		\item the sections of $q$ and $q'$ coincide, that is $\sigma_{AA'}(q') = d$,
		\item $q,q'$ coincide on $v$ and $v_1$. \qedhere
	\end{enumerate}
\end{corollary}

As an example, think of a circular net $q$ in $\R^2\subset \RP^2 \subset \CP^2$. In this case $A,A'$ are the imaginäre Kreispunkte at infinity. Thus from the corollary follows that there is a Q-net $q'$ contained in a circle with edges parallel to those of $q$. By construction, $q'$ has the same star-ratios as $q$.

In the case of a Q-net $q$ contained in one conic $c$, we can also consider stereographic projection from a cogeneric point $S\in c$ to the line $AA'$. Denote the stereographic projection of $q$ by $q^\infty$ and by $S^\infty$ the intersection with $AA'$ of the tangent to $c$ in $S$. Then we claim that
\begin{align}
	\mr(S^\infty,q^\infty,A,d^1,q^\infty_1,A')=-1
\end{align}
holds for every edge of the Q-net. In fact, this is a limit case of Desargues' involution theorem in which $S$ corresponds to two coinciding points. But this means that if we know the section $d=\sigma_{AA'}(q)$ and a single point of the projection $q^\infty$, then we know the whole projection. In the general case, finding the projection from the section is a 2D system (see Section \ref{sec:tcdextensions}). Unlike in the general case, here we are only left with a 1D system.

There is also a connection to Schief maps (see Definition \ref{def:paralleldm}), which appear in the particular limit case $A\rightarrow A'$, where we take the limit by moving $A$ towards $A'$ along a fixed line $L$.

\begin{definition}
	Let $A\in L$ be a point on a line in $\CP^2$. We say a Q-net $q:\Z^3\rightarrow \CP^2$ is \emph{contact-conical} if $L$ is generic and for every quad of $q$ the five points $A,q,q_1,q_2,q_{12}$ are on a conic with tangent $L$ in $A$.
\end{definition}

\begin{theorem}\label{th:contactconicalparalleldm}
	The section $\sigma_L(q)$ of a 1-generic contact-conical net with respect to $A\in L$ is a Schief map in an affine chart of $L$ with $A$ at infinity.
\end{theorem}
\proof{
	In a 2-conical net the points in the section of every quad satisfy the equation
	\begin{align}
		\mr(d^1,d^2,A,d^1_{2},d^2_{1},A')=-1.
	\end{align}
	Choosing the limit $A,A'\rightarrow \infty$ in an affine chart of $L$ we obtain that
	\begin{align}
		d^1 - d^2 = d^2_1 - d^1_2
	\end{align}
	holds.\qed
}

Again, there is a converse theorem.
\begin{theorem}\label{th:liftofparalleldm}
	Let $A$ be a point on the line $L$ in $\CP^2$. Let $d: E(\qg) \rightarrow L$ be a Darboux map such that in every quad
	\begin{align}
		\mr(d^1,d^2,A,d^1_{2},d^2_{1},A)=-1
	\end{align}
	holds. Then there is a Q-net $q: \qg \rightarrow \CP^2$ such that $q$ is contact-conical with respect to $A\in L$ and such that $\sigma_L(q) = d$.
\end{theorem}
\proof{Analogous to the proof of Theorem \ref{th:liftofdoublebkp}. \qed
}

Again, we may consider the special case where all points of $q$ are contained in a conic, but this time we take the section with respect to a tangent $L$ to $c$ in point $A$. We consider the multi-ratio equation
\begin{align}
	\mr(S^\infty,q^\infty,A,d^1,q^\infty_1,A')=-1
\end{align}
after the limit $A,A'\rightarrow \infty$, which yields
\begin{align}
	S^\infty + d^1 = q^\infty + q^\infty_1.
\end{align}

For one, this gives the 1D system a particularly simple shape. On the other hand, when Schief introduced Schief maps (as BKP maps \cite{schieflattice}) he also introduced an auxiliary map $m: V(\qg) \rightarrow \CP^1$ such that
\begin{align}
	d^k = \frac{m+m_k}{2}
\end{align}
for any direction $m$. In our setting, we regard $d$ as the section of a contact-conical Q-net $q$ with the tangent $L$. In this setting the auxiliary map $m$ is (up to a global translation and scaling) in fact the stereographic projection of the Q-net $q$. Choosing a different projection point $S$ clearly only translates the projection, as it only alters $S^\infty$.

\begin{remark}
	Consider a 2-conical Q-net $q$ with respect to $A,A'$. Then Theorem \ref{th:miquelbkp} states that $\sigma_{AA'}(q)$ is affine BKP with respect to both $A$ and $A'$. However, due to Theorem \ref{th:sectionsectioncluster}, we could consider any line $J$ through $A$ and still obtain that $\sigma_J(q)$ is affine BKP with respect to $A$. The analogous statement holds for contact-conical Q-nets. However, we have not considered whether in these cases there are constructions to obtain a lift from the section again. It is possible that a more thorough analysis of these cases could lead to understanding lifts to particular Q-nets for general affine BKP Darboux maps.
\end{remark}

\section{Q-nets in $\CP^3$ inscribed in a quadric}\label{sub:inscribedqnets}\label{sec:quadricqnet}

\begin{theorem}\label{th:inscribedqnetisbkp}
	Let $B$ be a non-degenerate quadric in $\CP^3$ and let $I$ be an isotropic line in $B$. Let $q: V(\qg) \rightarrow B$ be a 2-generic Q-net with points in $B$. Assume $I$ is generic. Then the section $\sigma_I(q)$ is a Doliwa complex, that is a projective BKP line complex.
\end{theorem}
\proof{Choose an additional point $K$ on the line $I$. There is a second isotropic line $I'$ through $K$. Together $I$ and $I'$ span a tangent plane $E$. Choose another plane $H$ in general position to $B,K,I$ and $I'$. Label the intersections $L = E\cap H, A=L\cap I$ and $A'= L\cap I'$. The proof consists of showing consecutive equivalences of the following four propositions:
\begin{enumerate}
	\item $\sigma_I(q)$ is projective BKP.
	\item $\sigma_E(q)$ is affine BKP with respect to $I$.
	\item $\pi_{K\rightarrow L}(\sigma_E(q))$ is affine BKP with respect to $A$.
	\item $\sigma_L(\pi_{K\rightarrow H}(q))$ is affine BKP with respect to $A'$.
\end{enumerate}
Proposition (1) is the claim of the theorem we intend to prove. Proposition (2) is equivalent to (1) because $\sigma_I(q) = \sigma_I(\sigma_E(q))$ and thus the projective cluster variables of $\sigma_I(q)$ are the affine cluster variables of $\sigma_E(q)$ with respect to $I$, see Theorem \ref{th:affprojcluster}. Proposition (3) is equivalent to (2) because the central projection $\pi_{K\rightarrow L}$ projects from the point $K$ which is on $I$. Therefore this projection preserves the affine cluster variables, see Theorem \ref{th:projectionandaffine}. Proposition (4) is equivalent to (3) because we claim that $\pi_{K\rightarrow L}(\sigma_E(q)) = \sigma_L(\pi_{K\rightarrow H}(q))$. Every edge $e$ of $\qg$ corresponds to a line $\ell_e$ of $q$. Both maps $\pi_{K\rightarrow L}(\sigma_E(q))$ and $\sigma_L(\pi_{K\rightarrow H}(q))$ map edges of $\qg$ to points in $L$, in fact they both map $e$ to $L\cap (\ell_e \cup K)$. Finally, (4) follows from Theorem \ref{th:miquelbkp}, which states that 2-conical nets are BKP with respect to $A,A'$. The fact that the stereographic projection $\pi_{K\rightarrow H}(q)$ is a 2-conical net follows from the arguments of the proof of Theorem \ref{th:projmiquel}.\qed
}

Unlike in the case of conics in $\CP^2$, we get a converse theorem without additional assumptions on the Doliwa complex.

\begin{theorem}\label{th:converseqnetquadricbkp}
	Let $B$ be a non-degenerate quadric in $\CP^3$ and let $I$ be an isotropic line in $B$. Let $l: F(\qg) \rightarrow I$ be a Doliwa complex. Then there is a lift, that is a Q-net $q: V(\qg) \rightarrow B$ such that $\sigma_I(q) = l$.
\end{theorem}
\proof{The planes of a Q-net determine the points of the Q-net, as every point is the intersection of at least three planes. The map $l$ determines one point in every plane of the Q-net. Let us consider one cube of the Q-net. We begin by choosing the plane through $l^{12}$. Next we choose the planes through $l^{23},l^{13},l^{23}_1,l^{13}_2$. After these choices there is a unique plane $U$ that completes the Q-net cube with points in $B$. By Theorem \ref{th:inscribedqnetisbkp} we know that 
\begin{align}
	\mr(l^{12},l^{23},l^{13},U\cap I, l^{23}_1,l^{13}_2) = -1
\end{align}
holds. But this is exactly the equation that is also satisfied if we substitute $U\cap I$ with $l^{12}_3$ and therefore $U$ contains $l^{12}_3$ indeed.\qed
}

We also mention a slight generalization.

\begin{corollary}\label{co:quadfamilyqnets}
	Let $I$ be a line in $\CP^3$ and let $q: V(\qg) \rightarrow B$ be a Q-net such that the points of every cube of $q$ is contained in a quadric that also contains $I$. Then the section $\sigma_I(q)$ is a Doliwa complex, that is a projective BKP line complex.	
\end{corollary}

The statement of the corollary follows from carefully considering the proof of Theorem \ref{th:inscribedqnetisbkp}. It is only based on considerations in one cube. Therefore it is not important that the quadric $B$ is the same for all cubes, it is only important that $I$ is a generator for all those quadrics. In the following, we return to the case of one fixed quadric $B$ for all cubes in order to explore what happens if $B$ is degenerate. However, the reader is invited to think about the generalizations in the spirit of Corollary \ref{co:quadfamilyqnets} by himself.

\begin{theorem}
	Assume $B$ is a cone (and thus a degenerate quadric) such that $E$ is a tangent plane that contains the isotropic line $I$. Let $q:V(\qg)\rightarrow B$ be a Q-net with points in $B$. Then the section $d=\sigma_E(q)$ is a Schief map with respect to $I$.
\end{theorem}
\proof{In every plane $U$ of a quad of $q$ we are actually in the situation of Theorem \ref{th:contactconicalparalleldm} (on sections of contact-conical nets). The intersection $U\cap B$ is a conic $C$ and $T = U \cap E$ is a tangent to $C$ in $P = I\cap U$. Because of Theorem \ref{th:contactconicalparalleldm} $\sigma_I(q)$ is a parallelogram in $U$ with respect to $P$ and thus also with respect to $I$.\qed}

We can even consider the case where the quadric $B$ degenerates to a pair of planes $B_+,B_-$ that intersect in the line $I$. However, we cannot distribute the points of the Q-net onto the two planes in an arbitrary manner. Denote by $V_+(\qg)$ the black and by $V_-(\qg)$ the white vertices of $\qg$.

\begin{theorem}\label{th:bkpdegquadric}
	Let $B_\pm$ be two planes in $\CP^3$ that intersect in the line $I$. Let $q: V(\qg) \rightarrow \CP^3$ be a Q-net such that
	\begin{align}
		q(V_+(\qg)) \subset B_+ \quad \mbox{ and }\quad q(V_-(\qg)) \subset B_-.
	\end{align}
	Then the section $\sigma_I(q)$ is a Doliwa complex.
\end{theorem}
\proof{The Q-net $q$ is a particularly simple case of a K{\oe}nigs net, Definition \ref{def:koenigs} is satisfied by the assumptions of the theorem. The only difference is that here, the two planes are actually the same for all three cubes. We also note that $I$ intersects all planes of the quads of $q$ in their diagonal intersection point. Thus $\sigma_I(q)$ is a special case of a Doliwa complex due to Lemma \ref{lem:diagkoenigs}. \qed   }

The only case of degenerate quadric that does not seem to lead to a theorem is the case of a double plane. It is possible that there is a sort of particular Q-net in $\CP^2$ that arises as limit of the cases mentioned above. The last theorem makes it plausible that the limits are in fact K{\oe}nigs nets, and the section is replaced with the diagonal intersection net. However, we do not carry out an investigation here.

In this section, the core statement was Theorem \ref{th:inscribedqnetisbkp} which associates a BKP structure to Q-nets with points in a non-degenerate quadric $B$ in $\CP^3$. We looked at two possible variations. On the one hand, we generalized quadrics to certain pencils of quadrics. On the other hand, we also looked at quadrics of varying degree of degeneracy (the rank of the associated bilinear form). An interesting further question is what happens if we look at $\CP^n$ with $n>3$. However, we postpone a discussion of that to Section \ref{sec:isotropicmaps}.

\section{S-graphs} \label{sec:sgraphs}

\begin{figure}
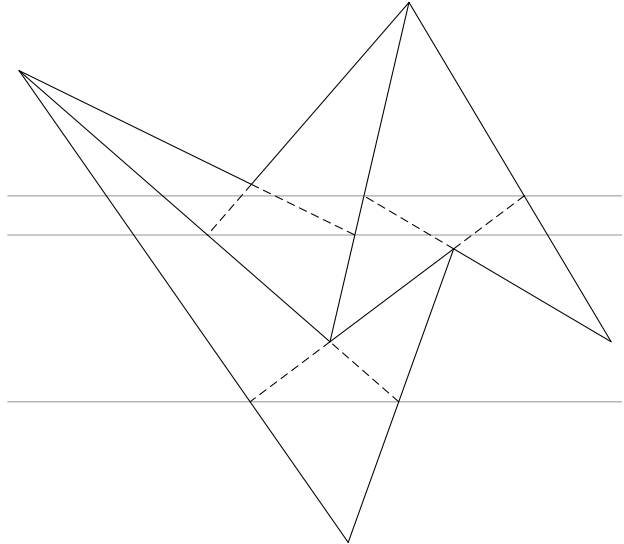


	\vspace{-4mm}
	\caption{Three quads in an S-graph.}
	\label{fig:sgraphs}
\end{figure}

\begin{definition}
	A Q-net $q: V(\qg) \rightarrow \CP^2$ is an \emph{S-graph} if there is a point $P \in \CP^2$ with the property that for every quad the two focal points and $P$ are colinear.
\end{definition}
In an affine chart with $P$ at infinity, this means that all lines that are spans of the two focal points of a quad, are parallel (see Figure \ref{fig:sgraphs}). This definition is a variation of the definition given by Chelkak \cite{chelkaksgraphs}, which is closely related to origami maps and s-embeddings.

\begin{theorem}\label{th:sgraphcarnot}
	Let $q:V(\qg) \rightarrow \CP^2$ be a 1-generic S-graph with respect to a point $P$ and let $H$ be a generic line that contains $P$. Then $q$ is affine CKP with respect to $H$, or equivalently $\sigma_H(q)$ is a Carnot map.
\end{theorem}
\proof{Consider $\CP^2$ as a plane $U \subset \CP^3$ and let $\hat q: V(\qg) \rightarrow \CP^3$ be a Q-net and $N\in \CP^3\setminus U$ a point such that $q$ is the projection $\proj{N}{U}(\hat q)$. Let $K$ be the line through $N$ and $P$ and let $E$ be the plane spanned by $K$ and $H$. Then we will show that $\sigma_E(\hat q)$ is a Carnot map. As a result (and with Theorem \ref{th:projectionandaffine}) $q$ is affine CKP with respect to $E\cap U = H$. Introduce the lines
\begin{align}
	\ell^{ij} &= \spa\{\fp^{ik},\fp^{ik}_j,\fp^{jk},\fp^{jk}_i\},\\
	h^{ij} &= \spa\{\fp^{ij},\fp^{ji}\},\\
	h^{ij}_k &= \spa\{\fp^{ij}_k,\fp^{ji}_k\},
\end{align}
for all triples of different indices $i,j,k$, where $\fp$ are the focal points as defined in Definition \ref{def:qgqnet}. Denote by $\hat \ell$ and $\hat h$ the corresponding lifts, that is the corresponding lines for $\hat q$. We observe that any one of the six lines $\hat h^{ij}$ intersects all of the three lines $\hat \ell^{kl}$ associated to a cube of $\Z^3$. Thus the lines $\hat h$ and $\hat \ell$ are generators of different type of a quadric $B \subset \CP^3$. As $K$ contains $P$ it intersects the lines $h^{ij}$. As we project from $N$ which is on $K$ the line $K$ also intersects the lines $\hat h^{ij}$, and is therefore itself a generator of the same type as the lines $\hat \ell^{ij}$. Thus $E$ is a tangent plane to the quadric $B$. As a result, $E$ also contains a second generator $G$ besides $K$ that is of the other type, and therefore intersects the lines $\hat \ell^{ij}$. But now we are in the setup of CQ-nets, see Definition \ref{def:cqnet} and Lemma \ref{lem:cqandcarnot}. By that Lemma we know that $\hat q$ is affine CKP with respect to any plane that contains $G$, and therefore also $E$. As $q$ is the projection of $\hat q$ from $N$, $q$ is affine CKP with respect to $E\cap U = H$.\qed
}

The proof of Theorem \ref{th:sgraphcarnot} is local, as we are arguing on a per-cube basis. Therefore, the converse is not necessarily true, not every projection from affine infinity of a CQ-net must be an S-graph. 

S-graphs will both help us to show how linear line complexes relate to CKP in Section \ref{sec:linearlinecomplexes} and also to show that s-embeddings are CKP in Section \ref{sec:sembeddings}.

\section{Linear line complexes in $\CP^3$}\label{sec:linearlinecomplexes}
In this subsection we study line complexes in relation to bilinear forms that are anti-symmetric instead of symmetric. A good introduction to the relation between linear line complexes and anti-symmetric bilinear forms can be found in Chapter 15 of Semple-Kneebone \cite{semplekneebone}. A general treatment of line geometry including a chapter on linear line complexes has been authored by Pottmann and Wallner \cite{pwlinegeometry}. With respect to integrability, linear line complexes were studied \cite{bobenkoschieflinecomplexes} as line complexes such that the lifts of all its lines to the Plücker quadric live in a codimension 1 space of $\CP^5$. Here we investigate the same discrete linear line complexes as defined by Bobenko and Schief \cite{bobenkoschieflinecomplexes}, but we view them in terms of anti-symmetric bilinear forms. In particular, one can choose an anti-symmetric bilinear form $\omega: \C^4 \times \C^4 \rightarrow \C$ on the homogeneous coordinates $\C^4$ of $\CP^3$. We assume throughout this subsection that $\omega$ is non-degenerate. For two points $p,p'\in\CP^3$ we will slightly abuse notation and write $\omega(p,p') = 0$ where we insert some homogeneous representatives of $p$ and $p'$ into $\omega$. The equality holds or fails independently of the choice of representatives. Let us gather a few basic properties of $\omega$. Given $p\in \CP^3$, the set of points $p'$ such that $\omega(p,p') = 0$ is clearly a plane, which we call the polar plane $p^\perp$ of $p$. Moreover, for two points $p,p'\in\RP^3$ the following is true due to the anti-symmetry of $\omega$:
\begin{align}
	\qquad\omega(p,p) &= 0\\ 
	\qquad\omega(p,p') &= 0 & &\eq & \omega(p',p) &= 0\\
	\qquad\omega(p,p') &= 0 & &\Rightarrow & \omega (p, \mu p + \mu' p') &= 0 \quad \forall \mu,\mu'\in \C
\end{align}
Ergo, every point is in its own polar plane. Call a line $\ell$ isotropic if for all $p,p'\in \ell$ holds $\omega(p,p')  = 0$. Given $p$ let $\ell$ be any line in the plane $p^\perp$ through $p$, then $\ell$ is isotropic. Thus there is a 1-parameter family of isotropic lines through every point and vice versa. There is a 3-parameter family of isotropic lines in total. The set of these lines is called a linear complex.

\begin{definition}
	Let $\omega$ be an anti-symmetric form as defined above and let $l: F(\Z^3) \rightarrow \CP^3$ be a line complex. Then $l$ is a \emph{linear line complex} if all lines of $l$ are isotropic with respect to $\omega$.
\end{definition}

\begin{lemma}\label{lem:sgraphandllc}
	Let $l: E(\qg) \rightarrow \CP^3$ be a 1-generic linear line complex with respect to an anti-symmetric form $\omega$. Let $P \in \CP^3$ be a point and let $P^\perp$ be the polar plane to $P$ and assume that $P^\perp$ is generic. Then $\sigma_{P^\perp}(l)$ is an S-graph with respect to $P$.
\end{lemma}
\proof{
	Let $q = \sigma_{P^\perp}(l)$ be the Q-net that arises by taking the section and consider a quad of $q$. Fix a homogeneous lift of $P$ and of every intersection point of the line complex. We abuse notation a bit and assume in the following linear equations that all the points are their respective lifts. Also abbreviate $\omega^i = \omega(l^i, P)$ et cetera. Then we can calculate that
	\begin{align}
		q &= \omega^2 l^1 - \omega^1 l^2,\\
		q_1 &= \omega^{21} l^1 - \omega^1 l^{21},\\
		q_2 &= \omega^{12} l^2 - \omega^2 l^{12},\\
		q_{12} &= \omega^{21} l^{12} - \omega^{12} l^{21}.
	\end{align}
	By definition, these points are on the corresponding lines of $l$. They are also in the plane $P^\perp$, as
	\begin{align}
		\omega(q,P) = \omega(q_1,P) = \omega(q_2,P) = \omega(q_{12},P) = 0 .
	\end{align}
	Moreover, we can form the following linear combinations of vertices of $q$ to express the focal points
	\begin{align}
		\fp^{12} &= \frac{1}{\omega^1}\left(\frac{1}{\omega^2}q - \frac{1}{\omega^{21}}q_1 \right) = \frac{1}{\omega^{12}}\left(\frac{1}{\omega^2}q_2 + \frac{1}{\omega^{21}}q_{12} \right) = \frac{1}{\omega^{21}}l^{21} - \frac{1}{\omega^{2}}l^2,\\
		\fp^{21} &= \frac{1}{\omega^2}\left(\frac{1}{\omega^1}q + \frac{1}{\omega^{12}}q_2 \right) = \frac{1}{\omega^{21}}\left(\frac{1}{\omega^1}q_1 - \frac{1}{\omega^{12}}q_{12} \right) = \frac{1}{\omega^{1}}l^{1} - \frac{1}{\omega^{12}}l^{12}.
	\end{align}
	But because $l$ are the lines of a linear line complex, its lines are isotropic and therefore the two points $l^2,l^{21}$ are both polar to the two points $l^1,l^{12}$. Hence, $\fp^{12}$ is polar to $\fp^{21}$. In the plane $P^\perp$ all isotropic lines pass through $P$. Therefore $\fp^{12},\fp^{21}$ and $P$ are on a line. \qed
}

\begin{theorem}\label{th:lincomplexckp}
	Let $\omega$ define a null-polarity in $\CP^3$ and let $I$ be an isotropic line of $\omega$. Let $l: E(\qg)\rightarrow \CP^3$ be a 2-generic linear line complex with respect to $\omega$. Assume $I$ is generic. Then the section $\sigma_I(l)$ is a Carnot map $d$, that is a Darboux map that is projective CKP.
\end{theorem}
\proof{Let $E$ be a plane that contains $I$. Lemma \ref{lem:sgraphandllc} states that $\sigma_E(l)$ is an S-graph with respect to $E^\perp$. Additionally, Theorem \ref{th:sgraphcarnot} states that $\sigma_{E^\perp}\circ\sigma_E(l) = \sigma_{E^\perp}(l)$ is a Carnot map.\qed
}

\begin{theorem}\label{th:converselinearlcckp}
	Let $\omega$ define a null-polarity in $\CP^3$ and let $I$ be an isotropic line of $\omega$. Let $d: E(\Z^3) \rightarrow I$ be a Carnot map. Then there is a lift, that is a linear line complex $l: F(\Z^3)\rightarrow \CP^3$ such that $\sigma_I(l) = d$.
\end{theorem}
\proof{
	Note that the planes that we intersect with $I$ to obtain $\sigma_I(l)$ are associated to the edges of $\qg$ and are the polar planes of the corresponding intersection points of the line complex. As a result of the anti-symmetry of $\omega$, each intersection point of $l$ has to lie in the polar plane of the corresponding Darboux map point $d$ on $I$. To begin with, we place $l^1$ on $d^{1\perp}$. Then we place $l^{2}$ and $l^{2}_1$ both on the designated polar planes and in the polar plane of $l^1$. Moreover $l^{1}_2$ is now in the intersection of three planes and therefore determined. The points $l^{3}$ and $l^3_1$ have to be on $l$ and $l_1$ as well as the designated polar plane and are therefore also determined, and $l^1_3$ is at the intersection of three planes. The point $l^3_2$ is also on a line and a plane and thus determined. Finally $l^2_3$ is on line $l_3$ and on $l^{3\perp}_2$, but should also be on $d^{2\perp}_3$ and is therefore over-determined. As a consequence, $d^2_3$ is not arbitrary but uniquely determined and by Theorem \ref{th:lincomplexckp} is the unique point that satisfies the CKP multi-ratio condition.\qed
}

A CKP structure for linear line complexes has previously been found by Bobenko and Schief \cite{bobenkoschiefcirclecomplexes}. In particular, they show that the CKP equation is satisfied by some $\tau$ variables that are minors of an M-system that they use to describe propagation of the line complex via Plücker coordinates. It would be interesting to understand the relations and to check if these two CKP structures coincide in some sense.

\section{Darboux maps tangent to a quadric}\label{sec:tangentdarboux}

We do not go into detail in this section, but mention the next theorem for completeness.
l

\begin{theorem}
	Let $B$ be a non-degenerate quadric in $\CP^3$ and let $I$ be an isotropic line in $B$. Let $d: E(\Z^3) \rightarrow B$ be a Darboux map such that all its planes (associated to cubes of $\Z^3$) are tangential to $B$. Then the section $\sigma_I(d)$ is a K{\oe}nigs net, that is a projective BKP Q-net.
\end{theorem}
\proof{Consider the dual $d^*$ of $d$, which is a Q-net (see Figure \ref{fig:qnclusterdualrelations}). The points of $d^*$ are the duals to the tangent planes occurring in $d$, and are thus contained in $B^*$. The dual $I^*$ of $I$ is a generator of $B^*$. Recall Theorem \ref{th:dualcluster} that identifies certain cluster structures of primal and dual TCD maps (again, see Figure \ref{fig:qnclusterdualrelations}). Therefore the projective cluster structure $\sigma_{I^*}(d^*)$ coincides with the projective cluster structure of $\sigma_I(d)$. However, Theorem \ref{th:inscribedqnetisbkp} tells us that $\sigma_{I^*}(d^*)$ is projective BKP and therefore so is $\sigma_I(d)$. Therefore $\sigma_I(d)$ is a K{\oe}nigs net.\qed
}

We observed in Theorem \ref{th:converseqnetquadricbkp} that any Doliwa complex (BKP line complex) in $\CP^1$ can be lifted to a Q-net with points in a quadric in $\CP^3$. In Theorem \ref{th:converselinearlcckp} we observed that any Carnot map (CKP Darboux map) in $\CP^1$ can be lifted to a linear line complex in $\CP^3$. However, in the case of a K{\oe}nigs net in $\CP^1$ we were not able to prove the analogous theorem and therefore leave it as a question.
l

\begin{question}
	Let $B$ be a non-degenerate quadric in $\CP^3$ and let $I$ be an isotropic line in $B$. Let $q: V(\Z^3) \rightarrow I$ be a K{\oe}nigs net. Is there a lift, that is a Darboux map $d: E(\Z^3) \rightarrow \CP^3$ such that all its planes are tangential to $B$ and such that $\sigma_I(d) = q$?
\end{question}

\section{Quadrirational Yang-Baxter Darboux maps}\label{sec:quadrirationalybmaps}

In this section we study maps that are a reduction of Darboux maps, which were introduced by Adler, Bobenko and Suris \cite{absyangbaxter}. We show that these maps fit naturally into the setup of TCD maps and show that certain sections of these maps are Doliwa compounds.

\begin{definition}\label{def:quadriybmap}
	Fix four different points $P_1,P_2,P_3,P_4\in \CP^2$, such that no three of them are colinear. Consider the pencil $\mathcal P$ of conics that pass through all four points $P_1,P_2,P_3,P_4$. Let $\qg$ be a quad-graph with $n$ strips. Fix $n$ different conics $C_1,C_2,\dots,C_n \in \mathcal P$. A Darboux map $d:E(\qg) \rightarrow \CP^2$ is called a \emph{quadrirational Yang-Baxter Darboux map} if for every edge $e\in E(\qg)$ crossing strip $i$ the imaged $d(e)$ is in conic $C_i$.
\end{definition}

Note that \cite{absyangbaxter} also considers the cases in which some of the points $P_1,P_2,P_3,P_4$ coincide. In those cases the pencil is defined by fixing higher order contact with some additionally chosen tangents. We keep this section brief by restricting ourselves to the most generic case of pencil.

In fact, it is clear that the image points $d(e^i),d(e^j)$ of two adjacent edges in the same quad of a quadrirational Yang-Baxter Darboux map determine the image points $d(e^i_j),d(e^j_i)$ of the other two edges in that quad, because these are the other intersection points of the line $d(e^i)d(e^j)$ with the conics $C_i,C_j$. Thus from the perspective of discrete consistency the quadrirational Yang-Baxter Darboux maps are 2-dimensional systems of Yang-Baxter type, as discussed briefly at the end of Section \ref{sec:sweepsqg}. 

Note that Darboux maps \cite{schieflattice} and quadrirational Yang-Baxter Darboux maps \cite{absyangbaxter} appeared in the literature around the same time, and the observation that the maps introduced by Adler, Bobenko and Suris are actually Darboux maps is already a new observation, albeit almost trivial. Adler, Bobenko and Suris also proved that quadrirational Yang-Baxter Darboux maps are 3-dimensionally consistent, which translates to being consistent reductions of Darboux maps.

\begin{theorem}
	Let $d:E(\qg) \rightarrow \CP^2$ be a quadrirational Yang-Baxter Darboux map with respect to points $P_1,P_2,P_3,P_4 \in \CP^2$. Fix two different indices $i,j\in \{1,2,3,4\}$ and let $H$ be the line $P_iP_j$. Assume that $d$ is 1-generic and that $H$ is generic with respect to $d$. Then $l=\sigma_H(d)$ is a Doliwa complex.
\end{theorem}
\proof{
	We have to show that in any 3-cube of $l$ the six intersection points of $l$ satisfy the multi-ratio equation in Lemma \ref{lem:Doliwacomplexlattice}. Choose a cube and denote the six intersection points of that cube by $l^{12},l^{23},l^{13},l^{12}_3,l^{23}_1,l^{13}_2$. Note that these points are the intersections of lines occurring in $d$ with $H$. Moreover, the four points $d^1,d^1_2,d^1_{23},d^1_3$ form a quad inscribed in the conic $C_i$. Thus Desargues' involution theorem (see Theorem \ref{th:desarguesinvolution}) applies, which states that in this case
	\begin{align}
		\mr(l^{12},l^{13},P_i,l^{12}_3,l^{13}_{2},P_j)=-1.
	\end{align}
	With the same reasoning we obtain that
	\begin{align}
		\mr(l^{12},l^{23},P_i,l^{12}_3,l^{23}_{1},P_j)=-1
	\end{align}
	holds as well. By combining the two previous multi-ratio equations we obtain that indeed
	\begin{align}
		\mr(l^{12},l^{23},l^{13},l^{12}_3,l^{23}_1,l^{13}_2) = -1,
	\end{align}		
	which is the requirement of Lemma \ref{lem:Doliwacomplexlattice} and thus concludes the proof.\qed
}

As a consequence, quadrirational Yang-Baxter Darboux maps are accompanied by the BKP equation.

\begin{remark}
	Note that we could also consider the other cases of quadrirational Yang-Baxter Darboux maps, where some of the base points of the pencil ($P_1,P_2,P_3,P_4$ above) coincide. This would require more per-case arguments, but Desargues' involution theorem also exists for tangent sections. We already took advantage of these tangent versions of Desargues' involution theorem in the study of contact-conical nets in Section \ref{sec:inscribedinconics}.
\end{remark}

We described quadrirational Yang-Baxter maps as quadrirational Yang-Baxter Darboux maps and found sections that are special line complexes. A natural follow up question is whether we can also understand the Darboux maps as the section of certain Q-nets. Assume there is a Q-net, such that the lines associated to the edges of the Q-net are generators of quadrics in a pencil, and generators associated to the same strip in the quad-graph belong to the same quadric. Then we can take the section with any (generic) plane to obtain a Darboux map with points in a pencil of conics as desired. Fortunately, it turns out that this is indeed possible. To show that these Q-nets exist, we use some facts from the theory of the curves that are intersections of quadrics. We give a very brief introduction but refer the reader to \cite{bstcheckerboardlaguerre, blptlaguerre} for more details.

Consider the following parametrized curve 
\begin{align}
	\gamma_k: \C \rightarrow \CP^3, \quad z \mapsto [\mathrm{cn}_k(z), \mathrm{sn}_k(z), \mathrm{dn}_k(z), 1].
\end{align}
Here we denoted by $\mathrm{cn}_k, \mathrm{sn}_k, \mathrm{dn}_k$ the \emph{Jacobi elliptic functions} of modulus $k$. Define the two quadrics
\begin{align}
	B_1 &= \{[x_1,x_2,x_3,x_4] \in \CP^3 : x_1^2+x_2^2=x_4^2 \}, \\  B_2 &= \{[x_1,x_2,x_3,x_4] \in \CP^3 : x_3^2+k^2x_2^2=x_4^2 \}.
\end{align}
Then due to the elementary relations between the Jacobi elliptic functions, the image of $\gamma_k$ coincides with $B_1\cap B_2$. Thus it makes sense to call $\gamma_k$ a quadric intersection curve. Moreover, every quadric in the pencil $\mathcal B$ spanned by $B_1$ and $B_2$ also contains $B_1\cap B_2$. In fact, in the generic case and up to projective transformations, there is one pencil of quadrics for each elliptic modulus $k$, and each pencil is completely determined by the corresponding quadric intersection curve.

Let us now look at two very handy properties of the elliptic parametrization of quadric intersection curves \cite[Proposition 8.4]{blptlaguerre}. First, the intersection of a quadric intersection curve with a generic plane consists of four different points. Four points $\gamma_k(z_1), \gamma_k(z_2), \gamma_k(z_3), \gamma_k(z_4)$ on a quadric intersection curve are coplanar if and only if $z_1+z_2+z_3+z_4 = 0$ (modulo the periods of the Jacobi elliptic functions). Secondly, every generator of a quadric in $\mathcal B$ intersects $\gamma_k$ twice. For every quadric $B(t) \in \mathcal B$ there is a $s_t\in \C$ such that for any $z\in \C$ the two other intersection points of the two generators of $B(t)$ through $\gamma_k(z)$ are $\gamma_k(s_t - z)$ and $\gamma_k(-s_t - z)$.

Together, these two properties allow us to construct certain planar quads. Associate the four points $\gamma_k(z_0),\gamma_k(t_1-z_0),\gamma_k(t_2-t_1+z_0),\gamma_k(-t_2-z_0)$ to the four vertices $q,q_1,q_{12},q_2$ of a quad. Clearly, the four points are coplanar. Moreover, opposite edge-lines of the quad correspond to generators of the same quadric. Therefore, if we intersect the edge-lines of the quad with a plane $H\subset \CP^3$, we obtain four points on a line, two of them are on the conic $B_{t_1}\cap H$ and two of them are on the conic $B_{t_2}\cap H$.

\begin{definition}
	Fix a generic pencil of quadrics $\mathcal B$ in $\CP^3$ and the corresponding intersection curve $\gamma_k$. Let $\qg$ be a quad-graph. We call a Q-net $q: V(\qg) \rightarrow \gamma(\C)$ an \emph{alternating generator Q-net} if in every quad the lines associated with the opposite edges are generators of the same quadric in $\mathcal B$.
\end{definition}

Note that in each quad the opposite generators are automatically of different type, unless the quad collapses to a point.

\begin{lemma}
	Let $q$ be a 1-generic alternating generator Q-net with respect to a generic pencil $\mathcal B$ and let $H$ be a plane that is generic with respect to $\mathcal B$ and $q$. Then $\sigma_H(q)$ is a quadrirational Yang-Baxter Darboux map.
\end{lemma}
\proof{
	The intersection of a pencil $\mathcal B$ of quadrics with a plane $H$ is a pencil $\mathcal P$ of conics. Because in every quad opposite edges of $q$ are generators of the same quadric $B(t)$, the intersections are in the same conic $C(t) \in \mathcal P$. Therefore along a strip in $\qg$ all points of $\sigma_H(q)$ belong to the same conic and that conic is in the pencil $\mathcal P$. The intersection of a generic pencil $\mathcal B$ with a generic hyperplane $H$ always consists of four different points and therefore the conditions of Definition \ref{def:quadriybmap} are satisfied.\qed
}

\begin{remark}
	Note that it is possible to drop the genericity conditions on $H$ and $\mathcal P$. In that case the section $\sigma_H(q)$ may be one of the non-generic cases of quadrirational Yang-Baxter Darboux maps considered in \cite{absyangbaxter}. Also note that in the study of so called checkerboard incircular nets \cite{bstcheckerboardlaguerre}, quadric intersection curves and generators of pencils are also used to generate these nets. The way the generators are combined to generate checkerboard incircular nets however is not the same as here. It would be interesting to understand how checkerboard incircular nets relate to quadrirational Yang-Baxter Darboux maps.
\end{remark}

\section{A-nets, isotropic line complexes in $\CP^5$}\label{sub:anets}\label{sec:anets}

A-nets were introduced by Sauer \cite{saueranet} and an introduction from the DDG perspective is given in the book \cite{ddgbook}. This section requires some basic knowledge in Plücker (line) geometry \cite{pwlinegeometry}. In Plücker geometry every line of $\RP^3$ is represented by a point in a quadric $\plu$ with signature $(3,3)$ in $\RP^5$, the quadric $\plu$ is called the Plücker quadric. The concept generalizes readily to $\CP^3$. Both points and planes of $\CP^3$ correspond to isotropic planes of the Plücker quadric. The points in these isotropic planes correspond to all lines through a point or alternatively all lines in a plane. 

\begin{definition}
	An \emph{A-net} is a map $a: V(\qg)\rightarrow \CP^3$ such that for each vertex $v\in V(\Z^3)$ the vertex $v$ and all its neighbours are mapped to a common plane. We denote by $a^*: V(\qg)\rightarrow \mbox{Gr}(2,\CP^3)$ the set of those planes.
\end{definition}

It is not immediately obvious from the definition, but A-nets naturally live in $\CP^3$. The coplanarity restrictions make sense only in more than two dimensions, and one can check that the points of an A-net span at most a three dimensional space.

Because there is a plane $a^*(v)$ and a point $a(v)$ associated to every vertex of $\qg$, we can think of the pair $(a(v),a^*(v))$ as a so called \emph{contact element}. From the point of view of Plücker geometry a contact element consists of all the lines in $a^*(v)$ through $a(v)$. In the Plücker quadric a contact element corresponds to an isotropic line. Moreover to an edge of an A-net belongs a line in $\CP^3$, which is the span of the two incident points, or alternatively the intersection of the two incident planes. 

\begin{definition}
	The \emph{Plücker lift} of an A-net $a$ is the map $\hat a: E(\qg) \rightarrow \plu$, such that $\hat a(v,v')$ is the lift of the line through $a(v)$ and $a(v')$.
\end{definition}

\begin{definition}
	Let $B \subset \CP^n$ be a quadric. We call a line complex $l: E(\Z^3) \rightarrow \CP^n$ an \emph{isotropic} line complex if all its lines associated to vertices of $\qg$ are isotropic with respect to $B$.
\end{definition}

Isotropic line complexes have been considered in the DDG book \cite[Section 3.3]{ddgbook}. The following observation is due to Doliwa \cite{doliwaanetspluecker}.

\begin{lemma}
	The Plücker lift $\hat a$ of an A-net $a$ is an isotropic line complex with respect to the Plücker quadric $\plu$.
\end{lemma}
\proof{
	Every line of $a$ is lifted to a point in the Plücker quadric. On each edge $(v,v')$ of $\qg$ the A-net defines a line through $a(v)$ and $a(v')$. All such lines associated to edges adjacent to a fixed vertex $v$ of $\qg$ belong to the contact element $(a,a^*)$. Therefore the lifts of those lines are contained in an isotropic line.\qed
}

Thus, A-nets are one more object from DDG that can be framed as a TCD map. Because the cube flip exists in the Plücker lift $\hat a$, it also exists for the A-net $a$. The Plücker lift is a line complex and thus A-nets are 4D consistent. While these last two statements are not new results, they follow almost trivially from the interpretation of the Plücker lift as line complex.

Also note that the representation of A-nets via the Plücker lift involves a line complex that satisfies some additional quadratic constraints. As we have seen in previous sections, TCD maps with additional quadratic constraints are often related to the BKP or CKP equation. Thus, one might suspect that there is an occurrence of the BKP or CKP equation also in A-nets, and indeed there is.

\begin{theorem}\label{th:cpfiveisolc}
	Let $B$ be a non-degenerate quadric in $\CP^5$ and let $I$ be an isotropic plane of $B$. Let $l:E(\qg) \rightarrow B$ be a 2-generic isotropic line complex. Assume $I$ is generic. Then $\sigma_I(l)$ is a Doliwa compound.
\end{theorem}
\proof{
	Choose a 4-space $H$ that contains $I$ and consider the section $q = \sigma_H(l)$. This is a Q-net with points in the quadric $H\cap B$, because every point of $q$ is the intersection of an isotropic line of $l$ with $H$. The quadric $H\cap B$ contains $I$ and thus contains an isotropic plane. Therefore $H\cap B$ is degenerate, as non-degenerate quadrics in $\CP^4$ only contain isotropic spaces of dimension strictly less than two. To see that $\sigma_I(q) = \sigma_I(l)$ is a line complex it suffices to look at a cube of $q$. The points of the cube in $q$ span a 3-space $F\subset H$. Thus, in this cube $\sigma_I(q)$ coincides with $\sigma_{I\cap F}(q)$, which is the section of a cube with a generator of the quadric $F\cap B$. Therefore, due to Theorem \ref{th:inscribedqnetisbkp}, the BKP condition holds in every cube of $\sigma_I(q)$ and therefore everywhere.\qed
}

By specifying $B$ to be the Plücker quadric, we immediately obtain a corollary concerning A-nets.

\begin{corollary}
	Let $a$ be an A-net such that its Plücker lift $\hat a$ is 2-generic, and let $E$ be an isotropic plane of the Plücker quadric that is generic with respect to $\hat a$. Then $\sigma_E(\hat a)$ is a Doliwa compound, that is a projective BKP line compound.
\end{corollary}

Interestingly, the DDG book \cite[Exercise 2.29]{ddgbook} contains a result, where a BKP structure appears in the context of A-nets. According to the DDG book, this is a discrete version of a result of K{\oe}nigs himself \cite{koenigsanets}. Note that we do not mention the full genericity assumptions in the next remainder of the section with respect to projections of A-nets, as we do not think they are of particular interest and would require explanations of how our standard genericity criteria translate from the Plücker lift to $\CP^3$. 

\begin{lemma}\label{lem:anetprojkoenigs}
	The projection $\proj PE:\CP^3 \rightarrow E$ of an A-net $a$ onto a plane $E$ from a point $P\notin E$ is a K{\oe}nigs net.
\end{lemma}

It is not immediately clear how this fits into our setup or generalizes to any isotropic line complexes in quadrics in $\CP^5$. However, we show that the two BKP structures are in fact the same ones (up to reciprocity). 

\begin{figure}
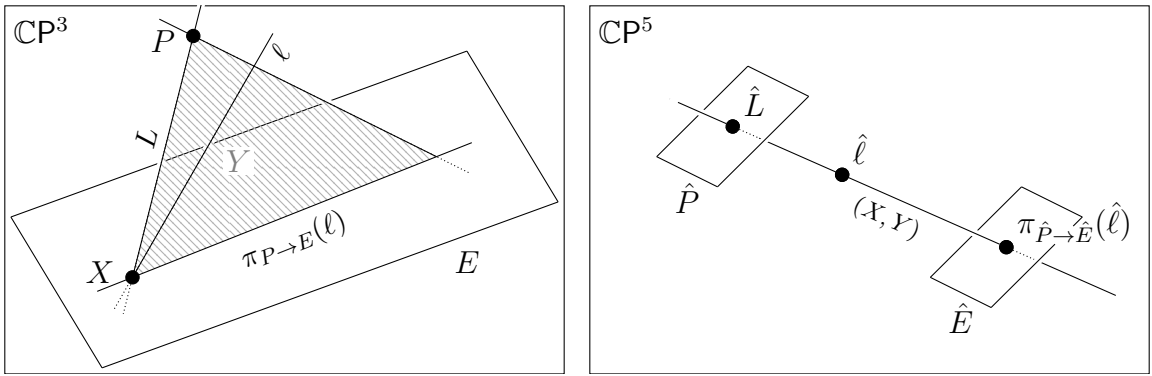


	\caption{Projection $\proj PE$ of the line $\ell$ from $P$ to $E$ on the left and the lift to the Plücker quadric on the right, with the corresponding projection $\proj{\hat P}{\hat E}$ in the lift.}
	\label{fig:pluliftvisual}
\end{figure}

\begin{figure}
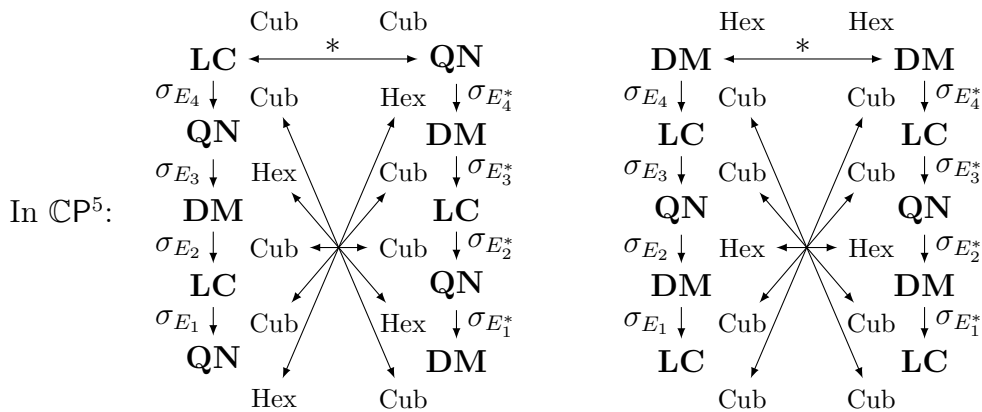

	%
	\vspace{-2mm}
	\caption{The two cluster structure duality diagrams in $\CP^5$. The left one relates line complexes and Q-nets.}
	\label{fig:cpfiveduality}
\end{figure}

\begin{theorem}
	Let $P$ be a point in $\CP^3$ and $E$ a plane not containing $P$. Let $a$ be an A-net. Denote the corresponding Plücker lifts by $\hat P,\hat E,\hat a$. Then the projective cluster structure of the Q-net $\proj{P}{E}(a)$ is reciprocal to the projective cluster structure of $\sigma_{\hat E}(\hat a)$.
	
\end{theorem}
\proof{
	We claim that the projection $\proj PE$ actually corresponds to the projection $\proj{\hat P}{\hat E}: \CP^5\rightarrow \hat E$, where $\hat P,\hat E$ are the Plücker lifts of $P$ and $E$, see Figure \ref{fig:pluliftvisual}. The correspondence is that for any line $\ell \subset \CP^3$ the lift of $\proj PE(\ell)$ equals $\proj{\hat P}{\hat E}(\hat \ell)$. Note that $\proj{\hat P}{\hat E}$ is a central projection from the 2-space $\hat P$ to the disjoint 2-space $\hat E$ in $\CP^5$. Moreover, every point of $\proj{\hat P}{\hat E}(\hat a)$ is the lift of the corresponding line of $\proj PE (a)$ and every line of $\proj{\hat P}{\hat E}(\hat a)$ is the lift of the contact element $(\pi(a),E)$. Thus $\proj{\hat P}{\hat E}(\hat a)$ is actually the projective dual of the K{\oe}nigs net $\proj PE(a)$. On the other hand, by our discussion in Section \ref{sec:projclusterduality} the projective dual of $\proj{\hat P}{\hat E}(\hat a)$ is $\sigma_{\hat E^\perp}(\hat a^\star)$ (see Figure \ref{fig:cpfiveduality}), where $\hat E^\perp=\hat E$ because it is an isotropic plane of $\plu$. Thus $\proj PE(a)$ and $\sigma_{\hat E^\perp}(\hat a^\star)$ are projectively equivalent and carry the same projective cluster structure. Moreover, by Theorem \ref{th:dualcluster} the cluster structure of $\proj{\hat P}{\hat E}(\hat a)$ is reciprocal to the cluster structure of $\sigma_{\hat E}(\hat a^\star)$. \qed
}

\begin{figure}
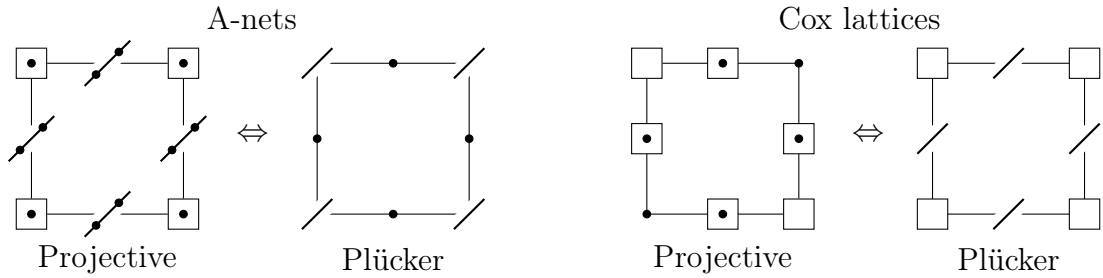


\caption{A quad in an A-net and in a Cox lattice, as well as the Plücker lifts. Dots represent points, squares represent planes, and a square with a dot represents a contact element.}
\label{fig:plulift}
\end{figure}

Before proceeding to the next section, let us note that there is also an analogue to Theorem \ref{th:cpfiveisolc} in $\CP^4$, even though we do not know of a relation to common examples.

\begin{theorem}
	Let $B$ be a non-degenerate quadric in $\CP^4$ and let $I$ be an isotropic line of $B$. Let $l:V(\qg) \rightarrow B$ be a 2-generic isotropic line complex and assume that $I$ is generic. Then $\sigma_I(l)$ is a Doliwa complex.
\end{theorem}
\proof{Choose a plane $E$ containing $I$. Then $\sigma_E(l)$ is a Q-net with points in the quadric $E\cap B$. Theorem \ref{th:inscribedqnetisbkp} proves the claim.\qed}

\section{Cox lattices, isotropic Darboux maps in $\CP^5$}\label{sub:cox}\label{sec:cox}

Cox lattices were introduced and studied by King and Schief \cite{kscox}. One of their findings is that Cox lattices are accompanied by a BKP structure. Continuing our approach, we cast Cox lattices as TCD maps, investigate, and show that there are two canonical appearances of the BKP equation in Cox lattices. In this section we restrict ourselves to maps defined on $\Z^3$ or ${\Z^3}^*$. The approach could be generalized to $\Z^n$ with $n > 3$ by introducing a new type of map defined on quad-graphs that have line complexes as sections.

\begin{definition}\label{def:cox}
	A \emph{Cox lattice} $c$ is a map from $V({\Z^3}^*)$ such that the image of each even (resp. odd) vertex is a point (plane) in $\CP^3$ and the image of every odd vertex is contained in the planes of the adjacent even vertices.
\end{definition}

We define Cox lattices on the dual lattice because it will simplify the identification with Darboux maps. As in the case of A-nets in Section \ref{sub:anets}, we make use of Plücker lifts. We identify every point and plane of a Cox lattice $c$ with an isotropic plane in the Plücker quadric. Every edge in the Cox lattice defines a contact element and thus an isotropic line in the Plücker quadric. To every face in the Cox lattice belongs the line spanned by the two points at even vertices of the quad. In the lift, every such line is a point on the Plücker quadric (see Figure \ref{fig:plulift}).

\begin{definition}
	Let $c: V({\Z^3}^*) \rightarrow \CP^3$ be a Cox lattice. The corresponding \emph{Plücker lift} is a map $\hat c: E(\Z^{3}) \rightarrow \plu$, where $\hat c((v,v'))$ is the Plücker lift of the line of $c$ associated to the quad $(v,v')^*$ in ${\Z^3}^*$.
\end{definition}

\begin{definition}
	Let $B$ be a quadric in $\CP^n$. A Darboux map $d: E(\Z^3) \rightarrow \CP^n$ is an \emph{isotropic Darboux map} with respect to $B$ if all its planes assigned to the cubes of $\Z^3$ are isotropic.
\end{definition}
We have stated the consistency of Q-nets in quadrics in Lemma \ref{lem:qnetquadcube} and mentioned the consistency of isotropic line complexes in Section \ref{sec:anets}. It is the same for isotropic Darboux maps: If the three adjacent planes $E_{\bar1},E_{\bar2},E_{\bar3}$ of a plane $E$ are isotropic then so is $E$ itself, because it intersects each $E_{\bar k}$ in an isotropic line. Thus the isotropic property of planes propagates from Cauchy data through $\Z^3$ in isotropic Darboux maps.

\begin{lemma}
	The Plücker lift $\hat c$ of a Cox lattice $c$ is an isotropic Darboux map.
\end{lemma}
\proof{The Plücker lift $\hat c$ assigns points in $\CP^5$ to edges of $\Z^3$. We argue that the points assigned by $\hat c$ to edges $e^1,e^2,e^{1}_2,e^2_1$ of a quad $u$ of $\Z^3$ are on a line in $\CP^5$. Note again that the dual to a quad in $\Z^3$ is an edge in ${\Z^3}^*$ and vice versa. Thus to the four edges $e^1,e^2,e^{1}_2,e^2_1$ of the quad $u$ in $\Z^3$ belong four quads ${e^1}^*,{e^2}^*,{e^{1}_2}^*,{e^2_1}^*$ around the edge $u^*$ in ${\Z^3}^*$. Thus $e^1,e^2,e^{1}_2,e^2_1$ are mapped by $\hat c$ to the lifts of the lines in $\CP^3$ associated by $c$ to the four quads ${e^1}^*,{e^2}^*,{e^{1}_2}^*,{e^2_1}^*$ adjacent to the common edge $u^*$ in ${\Z^3}^*$. If $u^*=(v,v')$ for $v,v'\in V(\Z^3)$ then the four lines in $\CP^3$ are part of the contact element defined by $c(v)$ and $c(v')$. Therefore $\hat c$ is a Darboux map. Additionally, the points in $\CP^5$ belonging to a cube of $\Z^3$ are on an isotropic plane, because the corresponding eight lines in $\CP^3$ of the Cox lattice $c$ are either all contained in a common plane, or all through a common point.\qed}

\begin{theorem}\label{th:coxDoliwa}
	Let $B$ be a non-degenerate quadric in $\CP^5$, let $I$ be an isotropic line of $B$ and let $d:E(\Z^3)\rightarrow \CP^5$ be a 2-generic isotropic Darboux map with respect to $B$. Assume $I$ is generic. Then $\sigma_I(d)$ is a Doliwa compound.
\end{theorem}
\proof{
	Consider any 3-space $H$ in general position to $B$ that contains $I$. Then the section $\sigma_H(d)$ is a Q-net with points contained in the quadric $B\cap H$. Because $\sigma_I(\sigma_H(d)) = \sigma_I(d)$ the claim follows from Theorem \ref{th:inscribedqnetisbkp}.\qed
}

As Cox lattices are special cases of isotropic Darboux maps, we have managed to associate a BKP structure to a Cox lattice via a section of their lift $\hat c$ to the Plücker quadric. However, we can take a more specific intermediate section that adds a second BKP structure to Cox lattices.

\begin{theorem}\label{th:coxkoenigs}
	Let $c$ be a Cox lattice with Plücker lift $\hat c$. Let $P$ be a point and $E$ a plane in $\CP^3$, such that $P\in E$ and such that $P$ is not contained in any of the planes of $c$ and $E$ contains none of the points of $c$. Let $H$ be the span of $\{\hat E, \hat P\}$ in $\CP^5$. Then $\sigma_H(\hat c)$ is a K{\oe}nigs net.
\end{theorem}
\proof{
	Note that $\hat E$ and $\hat P$ are two planes that intersect in a line and thus $H$ is a 3-dimensional space. Moreover, the intersection $H\cap B$ is $\hat E \cup \hat P$. By construction, the section $\sigma_H(d)$ is a Q-net $q$. The points of $q$ are contained in the two planes $\hat E$ and $\hat P$, but we need to verify that the odd points are in $\hat E$ and the even points are in $\hat P$ or vice versa. Let $v\in {\Z^3}^*$ and assume $c(v)$ is a point. Then $\sigma_H(\hat c)(v)$ is the lift of a line through $c(v)$ that is either contained in $E$ or containing $P$. However, by assumption $c(v) \notin E$ and therefore $\sigma_H(\hat c)(v)$ is in $\hat P$. An analogous argument convinces us that if $c(v)$ is a plane, then  $\sigma_H(\hat c)(v)$ is in $\hat E$. Thus, we are satisfying the assumptions of Definition \ref{def:koenigs}.\qed
}

Theorem \ref{th:coxkoenigs} could also be given for general isotropic Darboux maps in $\CP^5$, but the genericity assumption becomes somewhat opaque in that case.

Another way to understand Theorem \ref{th:coxkoenigs} is via projective duality of TCD maps. The dual of a Darboux map in $\CP^5$ is again a Darboux map (see Figure \ref{fig:cpfiveduality}). And indeed, the dual of an isotropic Darboux map with respect to $B$ is again an isotropic Darboux map with respect to $B^*$, because isotropic planes remain isotropic planes. Moreover, a generator $I$ of the Plücker quadric corresponds to a contact element $(P,E)$ in $\CP^3$. The dual $I^*$ of $I$ is the span of $\{\hat P, \hat E\}$. The fact that $\sigma_H(d)$ is a K{\oe}nigs net can also be deduced from our results on projective duality that state that the projective cluster structure of $\sigma_{I^*}(d)$ coincides with the projective cluster structure of $\sigma_{I}(d^*)$.

The BKP structure introduced by King and Schief \cite{kstetraoctacubo} uses homogeneous coordinates of $\CP^3$ and $\CP^{3*}$ plus a Hodge star operation or equivalently a bilinear form that identifies $\CP^3$ with $\CP^{3*}$. It is not clear whether their BKP structure corresponds to one of the two BKP structures we presented here.

\begin{remark}
	There is an interesting relation between A-nets and Cox lattices that is (without referring to Cox lattices) mentioned in the DDG book \cite[Section 2.4]{ddgbook}: Two layers of a Cox lattice yield one layer of an A-net. Above (in the $\Z^3$ sense) every point (resp.~plane) in the Cox lattice is a plane (point), together they define a contact element, and thus an A-net. As we have now also introduced the Plücker lift for Cox lattices, we can translate that observation to the lift. In the Plücker lift, we are actually ``collapsing'' two layers of $\hat c$. That is, assume we have a point $c(i,j,2k)$ and the plane $c(i,j,2k+1)$ above it. They define the contact element $(c(i,j,2k), c(i,j,2k+1))$. In the lift the contact element is represented by
	\begin{align}
		\hat c(i,j,2k) \cap \hat c(i,j,2k+1).
	\end{align}
	We are thus intersecting planes of the Darboux map $\hat{c}$ to obtain lines of the line complex $\hat{a}$. In fact, the relation is that
	\begin{align}
		\hat a = f_3(\hat c),
	\end{align}
	that is $\hat a$ is a focal line complex of the Darboux map $\hat c$ (see Section \ref{sec:focalnets} and Figure \ref{fig:comboffocal}). Thus Cox lattices and A-nets are actually related by focal transforms.
\end{remark}

\section{Anti-fundamental line-circle complexes}\label{sec:antifundamental}

In this section we make some observations on anti-fundamental line-circle complexes, which were introduced \cite{bobenkoschiefcirclecomplexes} as a reduction of linear line complexes. We use \emph{Laguerre geometry} \cite{laguerre} without giving an introduction. A classical treatment exists by Blaschke \cite{blaschkelaguerre}. As a modern (and English) treatment we recommend the book by Bobenko, Lutz, Pottmann and Techter \cite{blptlaguerre}.

\begin{definition}\label{def:aflc}
	An \emph{anti-fundamental line-circle complex} $x$ is a map from $V(\Z^3)$ such that the image of each even resp.~odd vertex is an oriented line resp.~an oriented circle in $\R^2$, and the images of adjacent vertices are in oriented contact.
\end{definition}

The \emph{Blaschke lift} is a lift of oriented lines in $\R^2$ to a cone $B$ in $\RP^3$, see \cite{blptlaguerre}. Each point in $B$ except for the apex corresponds to a unique oriented line in $\R^2$. Every plane that intersects $B$ in more than a point and not a line corresponds to all oriented lines in oriented contact with an oriented circle. This identification defines the Blaschke lift of oriented circles to planes in $\RP^3$.

\begin{definition}\label{def:aflclift}
	Let $x: V(\Z^3)\rightarrow \R^2$ be an anti-fundamental line-circle complex. The corresponding \emph{Blaschke lift} is a map $\hat x: V(\Z^{3}) \rightarrow B$, where $\hat x(v)$ is the Blaschke lift of the oriented line $x(v)$ if $v$ is even and the oriented circle $x(v)$ if $v$ is odd.
\end{definition}

\begin{theorem}
	The Blaschke lift $\hat x$ of an anti-fundamental line-circle complex $x$ is a Cox-lattice.
\end{theorem}
\proof{By comparing Definition \ref{def:aflclift} of the Blaschke lift to Definition \ref{def:cox} of Cox lattices the claim is obvious.\qed}

As an immediate corollary of the identification of anti-fundamental line-circle complexes with Cox lattices, we immediately obtain that anti-fundamental line-circle complexes are accompanied by two instances of the BKP equation, see Theorem \ref{th:coxDoliwa} and Theorem \ref{th:coxkoenigs} in Section \ref{sec:cox} on Cox lattices.

\section{Sequences of isotropic maps}\label{sec:isotropicmaps}

In all three examples of Q-nets, Darboux maps and line complexes defined on $\Z^3$ there is a partial order of projective subspaces appearing. The subspaces are first the points of the TCD map, then the lines, planes, 3-spaces, 4-spaces et cetera. They correspond to the points in sections of increasing codimension. We say a Q-net, Darboux map or line complex is $k$-isotropic if all of its $k$-spaces are isotropic (with respect to some given quadric $B$).

\begin{theorem}
	There are three sequences of isotropic maps that possess sections that are projective BKP.
	\begin{enumerate}
		\item Let $B$ be a non-degenerate quadric in $\CP^{2n+3}$ and let $I$ be an isotropic $(n+1)$-space. Let $\phi$ be an $n$-isotropic map to $\CP^{2n+3}$ that is a
		\begin{align}
			\begin{cases}
				\mbox{Q-net} & \mbox{if } n = 3k\\
				\mbox{Darboux map} & \mbox{if } n = 3k + 2\\
				\mbox{line complex} & \mbox{if } n = 3k + 1 
			\end{cases}
		\end{align}
		for some $k\in \N$. Then $\sigma_I(\phi)$ is projective BKP.
		
		\item Let $B$ be a non-degenerate quadric in $\CP^{2n+3}$ and let $I$ be an isotropic $n$-space. Let $\phi$ be an $(n+1)$-isotropic map to $\CP^{2n+3}$ that is a
		\begin{align}
			\begin{cases}
				\mbox{Q-net} & \mbox{if } n = 3k+2\\
				\mbox{Darboux map} & \mbox{if } n = 3k + 1\\
				\mbox{line complex} & \mbox{if } n = 3k
			\end{cases}
		\end{align}
		for some $k\in \N$. Then $\sigma_I(\phi)$ is projective BKP.
		
		\item Let $B$ be a non-degenerate quadric in $\CP^{2n+2}$ and let $I$ be an isotropic $n$-space. Let $\phi$ be an $n$-isotropic map to $\CP^{2n+2}$ that is a
		\begin{align}
			\begin{cases}
				\mbox{Q-net} & \mbox{if } n = 3k\\
				\mbox{Darboux map} & \mbox{if } n = 3k + 2\\
				\mbox{line complex} & \mbox{if } n = 3k+1
			\end{cases}
		\end{align}
		for some $k\in \N$. Then $\sigma_I(\phi)$ is projective BKP.\qedhere
	\end{enumerate}
\end{theorem}
\proof{
	For case (1) consider an $(n+3)$-space $E$ that contains $I$. Then the section $\sigma_E(\phi)$ is of codimension $n$. Therefore it is a Q-net with points in the quadric $B\cap E$. The quadric $B\cap E$ still contains the isotropic $(n+1)$-space $I$ and is thus degenerate for $n>0$. The span of the points of any cube of $\sigma_E(\phi)$ is 3-dimensional in an ambient $(n+3)$-space and thus intersects the $(n+1)$-space $I$ in a line $J$ which is also isotropic. Hence the claim follows from Theorem \ref{th:inscribedqnetisbkp}. 
	
	For case (2) consider an $(n+2)$-space $E$ that contains $I$. Then the section $\sigma_E(\phi)$ is of codimension $(n+1)$. Therefore it is a Q-net with points in the quadric $B\cap E$. The quadric $B\cap E$ still contains the isotropic $n$-space $I$ and is thus degenerate for $n>1$. The span of the points of any cube of $\sigma_E(\phi)$ is 3-dimensional in an ambient $(n+2)$-space and thus intersects the $n$-space $I$ in a line $J$ which is also isotropic. Hence the claim follows from Theorem \ref{th:inscribedqnetisbkp}.
	
	For case (3) consider an $(n+2)$-space $E$ that contains $I$. Then the section $\sigma_E(\phi)$ is of codimension $n$. Therefore it is a Q-net with points in the quadric $B\cap E$. The quadric $B\cap E$ still contains the isotropic $n$-space $I$ and is thus degenerate for $n>1$. The span of the points of any cube of $\sigma_E(\phi)$ is 3-dimensional in an ambient $(n+2)$-space and thus intersects the $n$-space $I$ in a line $J$ which is also isotropic. Hence the claim follows from Theorem \ref{th:inscribedqnetisbkp}\qed
}

\begin{question}
	Let $B$ be a quadric in $\CP^n, n>3$. Let $q: \Z^3\rightarrow B$ be a Q-net with points in $B$, such that the points of $q$ do span all of $\CP^n$. Are the cluster variables of $q$, or a section of $q$ in a BKP or CKP subvariety? Or are they in some other subvariety? What about the analogous questions for line complexes, Darboux maps and also for null-polarities instead of quadrics?
\end{question}

\chapter{Circular Q-nets, Poisson structure and quantization}\label{cha:poisson}

The purpose of this chapter is to relate results \cite{bmscircular} of Bazhanov, Mangazeev and Sergeev (BMS) on circular Q-nets to the TCD map framework. In particular we show that the affine cluster variables as defined in Definition \ref{def:affinecluster} can be expressed as a (simple) algebraic function of the BMS variables. Via this expression we show that both the Poisson and the quantum structure introduced by BMS are isomorphic to the canonical Poisson structure and quantum structure on cluster algebras.

The article of BMS \cite{bmscircular} does not give definitions of the Poisson and quantum structure they introduce, but instead focuses on a Poisson bracket as well as a Lie bracket and its invariance under the cube-flip. Therefore a rigid comparison of results is difficult. Thus we forego a complete introduction of the canonical Poisson and quantum structures and refer the reader to the article and book by Gekhtman, Shapiro and Vainshtein \cite{gsvpoissonpaper, gsvpoissonbook} for cluster Poisson structures and the article by Berenstein and Zelevinsky \cite{bzquantumcluster} for quantum cluster algebras. Instead, we will also restrict ourselves to comparing Poisson and Lie brackets.

\section{Poisson algebra, quantization}\label{sec:poisson}

\begin{definition}
	Let $A,A'$ be associative, commutative algebras. We call a bilinear map $\{\cdot,\cdot\}: A\times A\rightarrow A$ a \emph{Poisson bracket} if
	\begin{enumerate}
		\item $\{f,g\} = -\{g,f\}$ (anti-symmetry),
		\item $\{f,gh\} = g\{f,h\} + \{f,g\}h$ (Leibniz rule),
		\item $\{f,\{g,h\}\} + \{g,\{h,f\}\} + \{h,\{f,g\}\} = 0$ (Jacobi identity),
	\end{enumerate}
	for all $f,g,h\in A$. Moreover, a map $\phi: A \rightarrow A'$ preserves the Poisson bracket if
	\begin{align}
		\phi(\{f ,g\}_A) = \{\phi(f),\phi(g)\}_{A'}
	\end{align}
	for all $f,g\in A$. We also call $\phi$ a \emph{Poisson map}.
\end{definition}

Recall that it is possible to capture the arrows of a quiver with the bilinear, anti-symmetric form  $\nu$, see Equation \eqref{eq:nubil}. 

\begin{definition}\label{def:clusterpoisson}
	Let $\qui$ be a quiver with cluster variables $X_{1},X_2,\dots,X_m$. Consider the algebra $A_X$ that is the field of fractions of $\C[X_{1},X_2,\dots,X_m]$. The canonical Poisson bracket $\{\cdot,\cdot\}_\qui$ is generated by the relations
	\begin{align}
		\{X_v,X_{v'}\}_Q = \nu_{vv'}X_vX_{v'},
	\end{align}
	for all vertices $v,v'$ of $\qui$.
\end{definition}

Canonical cluster Poisson brackets are characterized by the following property, see \cite{gsvpoissonbook}.

\begin{lemma}\label{lem:clusterpoissonmap}
	The canonical cluster Poisson bracket is invariant under mutation.
\end{lemma}

Let us sketch the meaning and the proof of Lemma \ref{lem:clusterpoissonmap}. Consider vertices $v,v'\in V(\qui)$ and assume there is an arrow from $v$ to $v'$. Let $\tilde \qui$ be $\qui$ after a mutation at $v$. Then $\tilde X_v=X_v^{-1}$ and $\tilde X_{v'}=X_{v'}(1+X_v)$. On the one hand, we can calculate
\begin{align}
	\{\tilde X_v,\tilde X_{v'}\}_{\tilde \qui} = \tilde\nu_{vv'}\tilde X_v \tilde X_{v'} = -\tilde X_v \tilde X_{v'}.
\end{align}
On the other hand, we can calculate
\begin{align}
	\{\tilde X_v,\tilde X_{v'}\}_{\qui} &= \{X_v^{-1},X_{v'}(1+X_v)\}_{\qui} = -\frac{1}{X_v^2} \{X_v,X_{v'}\}_{\qui}(1+X_v)\\
	&=-\nu_{vv'}\frac{X_{v'}}{X_v} (1+X_v)=-\tilde X_v \tilde X_{v'}.
\end{align}

Similar calculations can be carried out for the other possibilities of local configurations in the quiver, and show that indeed the Poisson bracket is invariant under mutation.

Apart from the Poisson algebra there is also a definition for a canonical quantum structure \cite{bzquantumcluster} associated to each cluster algebra. A rigorous exposition is beyond the scope of this thesis. One of the obstacles is that in the non-commutative case, it is not so straightforward to define the space of rational functions. We give a formal definition of an algebra that allows us to demonstrate the commutation relations, which in turn will allow us to compare these to commutation relations previously obtained \cite{bmscircular} in the case of circular Q-nets.

\begin{definition}\label{def:clusterquantum}
	Let $\qui$ be a quiver with cluster variables $X_{1},X_2,\dots,X_m$. Then we define the algebra $\mathcal A_X$ to be the skew-field of fractions of the associative, non-commutative algebra over $\C$ that is generated by the elements 
	\begin{align}
		\{q\} \cup \{X_v : v\in Q\},
	\end{align}
	subject to the \emph{commutation relations}
	\begin{align}
		q^{-\nu_{vv'}}X_vX_{v'} &= q^{-\nu_{v'v}}X_{v'}X_{v},
	\end{align}
	for all quiver vertices $v,v'$, and with $q$ in the center.

\end{definition}

It is not necessary to consider the whole skew-field of fractions to proceed in the following and there are also some technical details to consider for the existence of the skew-field, but we refer the reader on both matters to \cite{bzquantumcluster}.

It is not a priori clear what the effect of a mutation is on the cluster variables in the quantum case, as the variables are now non-commuting. In the setup of quantum cluster algebras \cite{bzquantumcluster} though mutation can be defined via Equation \eqref{eq:mutationx} with the variables in the order shown in Equation \eqref{eq:mutationx}. Let us give an example calculation to illustrate the invariance of the commutation relations under mutation. We make the same assumptions as in the Poisson example, that is we consider vertices $v,v'\in V(\qui)$ and assume there is an arrow from $v$ to $v'$. This time we calculate, while considering the order of variables, that on the one hand
\begin{align}
	\tilde X_v\tilde X_{v'} = q^{2\tilde \nu_{vv'}} \tilde X_{v'} \tilde X_v = q^{-2} \tilde X_{v'} \tilde X_v.
\end{align}
On the other hand
\begin{align}
	\tilde X_v\tilde X_{v'} &= X_v^{-1} X_{v'}(1+X_v) = q^{-2\nu_{vv'}}X_{v'} (1+X_v)X_v^{-1} = q^{-2} \tilde X_{v'} \tilde X_v.
\end{align}

\begin{remark}
	This was a very brief introduction to cluster Poisson and quantum structures. We will only encounter them in the next section on circular Q-nets. But in the literature, both structures feature heavily in dimer cluster integrable systems \cite{gkdimers} and in higher Teichmüller theory \cite{fghighertm}. The Poisson structure is also used to relate results on the Liouville-integrability of cross-ratio dynamics \cite{afitcrossratio} to dimer integrable systems in \cite{agrcrdyn}. It is also suspected that there is a (direct) relation to the Poisson structure of the pentagram map \cite{ostpentagram} (see also \cite{gdnotes}).
\end{remark}

\section{Circular Q-nets} \label{sub:circularq}\label{sec:circularnets}

\begin{definition}
	A \emph{circular Q-net} $q:V(\qg) \rightarrow \R^3$ is a Q-net such that each quad is inscribed in a circle.
\end{definition}

Circular Q-nets were introduced by Bobenko and independently by Cie\'sli\'nski, Doliwa and Santini \cite{cdscircular} in relation to triply orthogonal coordinate systems \cite{bobenkocircular}. See also the bibliographical notes of \cite[Chapter 3]{ddgbook} for further references.

One can also consider a projective constraint on Q-nets in $\CP^3$ that can be specialized to circular Q-nets as follows: Fix a plane $E$ and a non-degenerate conic $C\subset E$. Consider Q-nets such that for each cube, the whole cube is inscribed in a quadric that contains the conic $C$. A circular Q-net is just the case where $E$ is the plane at infinity and $C$ is the conic
\begin{align}
	C = \{(x_1,x_2,x_3,0)\in \C^4\ | \ x_1^2+x_2^2+x_3^3 = 0\}.
\end{align}
If we allow complex projective transformations, we can always normalize any plane to be the plane at infinity and any non-degenerate conic in that plane to be the conic in the equation above. Thus, the projective perspective does not add a significant amount of freedom and we confine ourselves to consider circular Q-nets. The projective perspective however motivates us to study the section of a circular Q-net with the plane at infinity. The aim of the current section is to relate the affine cluster variables with respect to the plane at infinity to variables introduced by Bazhanov, Mangazeev and Sergeev \cite{bmscircular} (we abbreviate this triplet of authors by BMS). BMS investigated the propagation of these variables in a circular Q-net and introduced a Poisson structure as well as a variational principle. Moreover, they found a quantization that they can relate to multiple known quantum structures. We will show the BMS variables are potentials for the affine cluster structure and how BMS's Poisson and quantum structure induce the canonical cluster Poisson algebra and cluster quantum structure (see Section \ref{sec:poisson}) on the affine cluster structure of the Q-net.

\begin{figure}[]
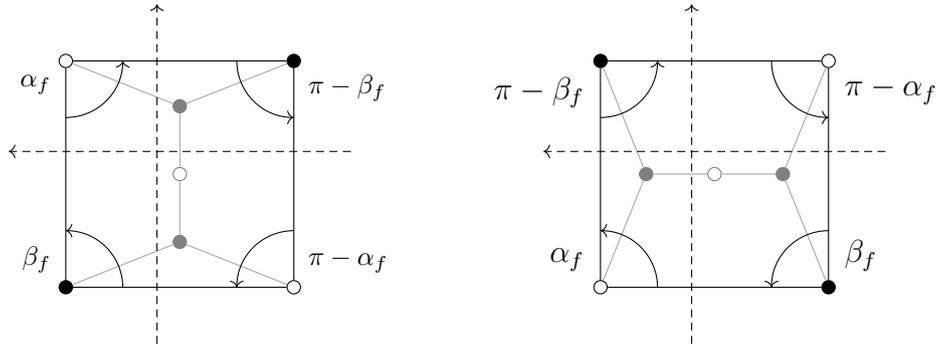
\centering


	
	\caption{Conventions for the angles in a quad $f$ of a circular Q-net. }
	\label{fig:quadstripesconvention}
\end{figure}

A useful observation for circular Q-nets is that opposing angles in each quad sum to $\pi$, which is an elementary consequence of the inscribed angle theorem. We can therefore use two variables $\alpha_f,\beta_f$ in each quad $f$ to parametrize the four angles by $\alpha_f, \beta_f, \pi - \alpha_f$ and $\pi - \beta_f$. We copy the approach of BMS of viewing a Q-net as a map on a quad-graph and using oriented strips of quad-graphs, see also Section \ref{sec:quadgraphs}. In each quad, we assign the angles $\alpha_f$, $\beta_f$, $\pi-\alpha_f$ and $\pi -\beta_f$ according to Figure \ref{fig:quadstripesconvention}. Note that our convention is slightly different from the convention of BMS, who orient the strips in a particular way (as in pseudoline arrangements, see \cite{felsnerbook}). We assume a fixed but arbitrary orientation is chosen for each strip. Note that in counterclockwise direction the angles always appear as $\beta_f,\alpha_f,\pi-\beta_f,\pi-\alpha_f$, and the $\alpha_f$ angles are always located at black vertices of $\qg$ while the $\beta_f$ angles are always located at white vertices of $\qg$.

\begin{definition}
	The \emph{BMS variables} $a_f,\bar a_f,k_f$ for a quad $f$ are
	\begin{align}
		k_f = \frac{\sin \beta_f}{\sin \alpha_f},\quad a_f = \frac{\sin(\alpha_f+\beta_f)}{\sin(\alpha_f)},\quad \bar a_f = \frac{\sin(\alpha_f-\beta_f)}{\sin(\alpha_f)},
	\end{align}
	subject to the constraint $a \bar a=1-k^2$. The \emph{BMS Poisson structure} is generated by the relations
	\begin{align}
		\{k_f,a_{f'}\}_{\mathrm{BMS}} = \delta_{ff'} k_fa_{f'},\quad \{k_f,\bar a_{f'}\}_{\mathrm{BMS}} = -\delta_{ff'} k_f\bar a_{f'},\quad \{a_f,\bar a_{f'}\}_{\mathrm{BMS}} = 2\delta_{ff'} k_f^2,
	\end{align}
	for any two quads $f,f'$ of $\qg$.
\end{definition}

The BMS Poisson structure is called ultra-local as the only non-trivial relations are among variables assigned to one quad. Consider the effect of reversing the orientation of a strip on the angles in a quad $f$ that is crossed by the strip. There are two cases, either the angles remain unchanged or all angles are replaced by $\pi$ minus the respective angle. Therefore, the BMS variables of $f$ transform as $k_f\mapsto k_f, a_f\mapsto \pm a_f$ and $\bar a_f\mapsto \pm \bar a_f$. So, the Poisson relations are preserved under the change of strip orientation.

Now, we can relate the BMS variables and the BMS Poisson structure to the affine cluster variables $Y$ of Q-nets. Recall that the affine cluster variables of a Q-net fall into three groups, variables $Y_w$ situated at white vertices of $\qg$, variables $Y_b$ situated at black vertices of $\qg$ and variables $Y_f$ associated to the quads of $\qg$.

\begin{lemma}\label{lem:circularfac}
	The affine cluster variables $Y$ with respect to the plane at infinity are related to the BMS variables by the formulas
	\begin{align}
		Y_{f} = -k_f^2,\qquad Y_w = (-1)^{d_w} \prod_{f\sim w} \bar a_fk_f^{-1},\qquad Y_b = (-1)^{d_b} \prod_{f\sim b} \bar a_f^{-1},\label{eq:srbms}
	\end{align}
	where $w,b$ are white and black vertices of the quad-graph $\qg$ and $d_w,d_b$ are the respective degrees of the vertices $w,b$. 
\end{lemma}
\proof{According to Lemma \ref{lem:srviaangles} the variable $Y_{f}$ at quad $f$ is
	\begin{align}
		Y_{f} &= -\frac{\sin(\beta_f)^2}{\sin(\alpha_f)^2} = -k_f^2,
	\end{align}
	independently of how the strips cross the quad. Moreover, the multiplicative contribution of a quad at a white vertex is
	\begin{align}
		\pm \frac{\sin(\alpha_f - \beta_f)}{\sin(\beta_f)} = \pm \bar a_f k_f^{-1},
	\end{align}
	where the sign depends on how the strips cross the quad. However, we note that reorienting a strip will change the sign at two quads incident to a given white vertex. The reversal of the orientation therefore leaves the total sign unchanged, and we can assume that the strips are oriented cyclically around the black vertex. In this case the claim follows. Finally, at a black vertex, the contribution of each quad is
	\begin{align}
		\pm \frac{\sin(\alpha_f)}{\sin(\alpha_f - \beta_f)} = \pm \bar a_f^{-1},
	\end{align}
	where the same arguments for the signs hold as before.\qed
}

\begin{lemma}
	The map that expresses the $Y$ variables in terms of the BMS variables is  a Poisson map up to a rescaling by a factor 2.
\end{lemma}
\proof{Direct computation. For example, if we express two cluster variables $Y_{f},Y_w$ by the BMS variables and then compute the BMS Poisson bracket we obtain
	\begin{align}
		\{Y_{f}, Y_b\}_{\mathrm{BMS}} = \{k_f^2, \prod_{f'\sim b} \bar a_{f'}^{-1}\}_{\mathrm{BMS}} = -2k_f\bar a_f^{-2}\{k_f,a_f\}_{\mathrm{BMS}}\prod_{\substack{f'\sim b\\ f'\neq f}} \bar a_{f'}^{-1}= 2 Y_{f}Y_b.
	\end{align}
	We can do these calculations for all relevant combinations to obtain the relations
	\begin{align}
		\{Y_{f},Y_w\}_{\mathrm{BMS}} = 2 \delta_{fw} Y_{f} Y_w,\quad \{Y_{f},Y_b\}_{\mathrm{BMS}} = -2\delta_{fb} Y_fY_b,\quad \{Y_w,Y_b\}_{\mathrm{BMS}} = 2\delta_{wb} Y_w Y_b.
	\end{align}
	Up to the time-rescaling factor 2, these relations coincide with the canonical Poisson relations
	\begin{align}
		\{Y_{f},Y_b\} = \delta_{fb} Y_{f} Y_b,\quad \{Y_{f},Y_w\} = -\delta_{fb} Y_fY_w,\quad \{Y_b,Y_w\} = \delta_{bw} Y_b Y_w,	
	\end{align}
	introduced in Section \ref{sec:poisson}, in the case of the affine cluster structure.\qed
}

It may already have been a surprise that it is possible to define an ultra-local Poisson bracket on the BMS variables. But even somewhat more astonishingly, the BMS Poisson algebra can be reduced to the following brackets on the angles:
\begin{align}
	\{\alpha_f,\beta_{f'}\} = \delta_{f,f'}, \quad \{\alpha_f,\alpha_{f'}\} = 0, \quad \{\beta_f,\beta_{f'}\} = 0.
\end{align}
BMS employ these to introduce the canonical quantization
\begin{align}
	[\alpha_f,\beta_{f'}] =\hbar \delta_{f,f'} , \quad [\alpha_f,\alpha_{f'}] = 0, \quad [\beta_f,\beta_{f'}] = 0.
\end{align}
From this Lie bracket one can derive the following commutation rule for the variables-turned-operators $\bar a$ and $k$:
\begin{align}
	k_f \bar a_f = q \bar a_f k_f,
\end{align}

where $q$ is a quantization constant. The Lie-brackets on the angle-variables also imply that operators of different quads commute. Assume the factorization formula \eqref{eq:srbms} holds for the $Y$-variables-turned-operators in the order given. Because of the ultra-locality of the commutation rules the ordering of the product over quads does not matter, only that $\bar a_f$ precedes $k^{-1}_f$. Let us do an example calculation to find the commutation rule for the $Y$ operators for $f\sim w$:
\begin{align}
	Y_{f} Y_w = k_f^2   \prod_{f'\sim w} \bar a_{f'}k_{f'}^{-1} = q^2 \prod_{f'\sim w}\bar  a_{f'}k_{f'}^{-1}  k_f^2 = q^2 Y_{f} Y_w.
\end{align}
In fact, one readily verifies that cluster operators belonging to variables that are not adjacent in the quiver commute, and if they are adjacent then the relations are
\begin{align}
	Y_{f} Y_w = q^2 Y_{f} Y_w, \quad Y_{f} Y_b = q^{-2} Y_{f} Y_b,\quad Y_w Y_b = q^2 Y_w Y_b.
\end{align}
These coincide precisely with the canonical commutation rules of cluster algebras, see Definition \ref{def:clusterquantum}.

BMS also obtained some results for the case of general Q-nets \cite{bmscircular}. In the general case, the BMS variables become $A_f,B_f,C_f,D_f$ and are neatly represented in the matrix
\begin{align}
	\ve{A_f & B_f\\ C_f & D_f} = \ve{\frac{\sin (\gamma_f)}{\sin (\delta_f)} & \frac{\sin (\delta_f + \sin \beta_f)}{\sin (\delta_f)}\\ \frac{\sin (\delta_f + \gamma_f)}{\sin (\delta_f)} & \frac{\sin (\beta_f)}{\sin (\delta_f)}},
\end{align}
where $\alpha_f,\beta_f,\gamma_f,\delta_f$ are the angles of quad $f$. These variables are subject to the constraint $(AD-BC)(DB-AC)=(AB-CD)$. It is not hard to see that for general Q-nets, one can obtain the affine cluster variables as certain products of the general BMS variables in the same manner as for the case of circular Q-nets, for example $Y_f = - A_fD_f$. However, as BMS do not give any Poisson algebra structure in the general case, there is not much to prove for us here and we skip a detailed exposition. It would be interesting to understand whether there is a Poisson algebra for the variables $A_f,B_f,C_f,D_f$ that is Poisson equivalent to the canonical affine cluster Poisson algebra. This would possibly allow to find similar results for general TCD maps.

One may also wonder whether the affine cluster variables of circular Q-nets are in some subvariety, possibly induced by the factorization as in Lemma \ref{lem:circularfac}. However, one can show that for any Q-net with affine variables $Y_f,Y_w,Y_b$, there exist BMS-variables $k_f,a_f,\bar a_f$ such that the factorization of Lemma \ref{lem:circularfac} holds. To see this, choose a square root of every variable $Y_f$ to determine the variables $k_f$ such that $-k_f^2=Y_f$. Determining the $\bar a_f$ variables then becomes a 2D-system. Let us for the moment also consider the frozen $Y_w, Y_b$ variables at boundary vertices, see the comment below Definition \ref{def:affinecluster}. We can follow the li-orientation for quads (see Section \ref{sec:sweepsqg}). In each step we add a quad we also close off a vertex, which therefore uniquely defines $\bar a$. Therefore, for general Q-nets it is always possible to ultra-localize the Poisson variables via this semi-local transformation of variables. It would be interesting to know if similar transformations exist in the general literature. It is also not clear whether the 2D system explained above is actually 3D integrable. This would most likely involve a careful choice of roots for $k_f$, and one has to keep in mind that the $Y$-variables also change in a cube flip.

There is a curious similarity of the factorization in Lemma \ref{lem:circularfac} and the definition of the Ising subvariety (see Definition \ref{def:isingsub}). Indeed, one can combine the two formulas and deduce that a circular Q-net has affine cluster variables in the Ising subvariety if and only if
\begin{align}
	\prod_{f\sim v} a_f\bar a_f^{-1} = \pm1, \label{eq:circularising}
\end{align}
for every interior vertex of $\qg$. We do not specify the sign of the product, because here the sign does depend on the signs of the $a_f$ variables, which change under reorientation of strips. It would be interesting to understand the geometric meaning of condition \eqref{eq:circularising}. Note that a Q-net has affine cluster variables in the Ising subvariety if and only if it is a CQ-net, see Section \ref{sec:cqnets}. Therefore Equation \eqref{eq:circularising} yields a way to recognize circular Q-nets that are CQ-nets (or vice versa) via an angle condition.


\chapter{Embeddings from statistical physics as TCD maps in $\CP^1$}\label{cha:cpone}
\chaptermark{Embeddings from stat. physics as TCD maps in $\CP^1$}

\section{Circle patterns and t-embeddings}\label{sec:cptemb}

\begin{definition}
	Let $G$ be a bipartite planar graph. A \emph{circle pattern} is a map $c: V(G)\rightarrow \CP^1$ such that the image of each face of $G$ is concyclic.	
\end{definition}
The predecessor of circle patterns were circle packings in the spirit of discretizing conformal maps. This approach is originally due to Koebe \cite{koebecp} and Thurston \cite{thurstoncp}. A good introduction to circle packings is by Stephenson \cite{stephensoncp}. The generalization to circle patterns with arbitrary intersection angles was started by Thurston \cite[Section 13.7]{thurstonlecturenotes} and Rivin \cite{rivin}.

Note that there are similarities between circle patterns and circular nets (Section \ref{sec:circularnets}) as well as special cases of 2-conical nets (Section \ref{sec:inscribedinconics}). However, circular nets that live in $\R^n$ with $n>2$ cannot generally be identified with circle patterns in $\CP^1$. Special 2-conical nets in $\RP^2$ with the `imaginäre Kreispunkte' as distinguished points can be identified with circle patterns, if there are no points at the line at infinity. However, the canonical cluster variables of 2-conical nets in $\RP^2$ are different from those we will introduce for circle patterns in $\CP^1$. Moreover, the local moves that we investigate are different.

Recall that we defined PDB quivers (planar, dualy bipartite quivers) in Definition \ref{def:planquiver}, as oriented duals of bipartite graphs.

\begin{definition}\label{def:tembedding}
	Let $\qui$ be a PDB quiver. A map $t: V(\qui)\rightarrow \C$ is a \emph{t-embedding} if for every vertex $v$ of $\qui$ with cyclically ordered neighbours $v_1,v_2,\dots,v_{2m}$ the consecutive angles satisfy
	\begin{equation}
		\sum_{k=1}^{m}\arg \left( \frac{t(v_{2k-1})-t(v)}{t(v_{2k})-t(v)} \right) \in \pi \Z.\qedhere
	\end{equation}	
	 
\end{definition}

Note that in the definitions of both circle patterns and t-embeddings we do not require any embeddedness properties. The term t-embeddings was coined in \cite{clrtembeddings}, however t-embeddings have also previously and independently appeared under the name of Coulomb gauge \cite{kenyonlam}. Moreover, in the DDG community t-embeddings are the planar case of conical nets, which has also been studied separately \cite{muellerconical}. We use the name t-embedding instead of Coulomb gauge because we focus on the centers, and instead of conical net to avoid confusion with several competing definitions of conical nets in higher dimensions.

\begin{figure}
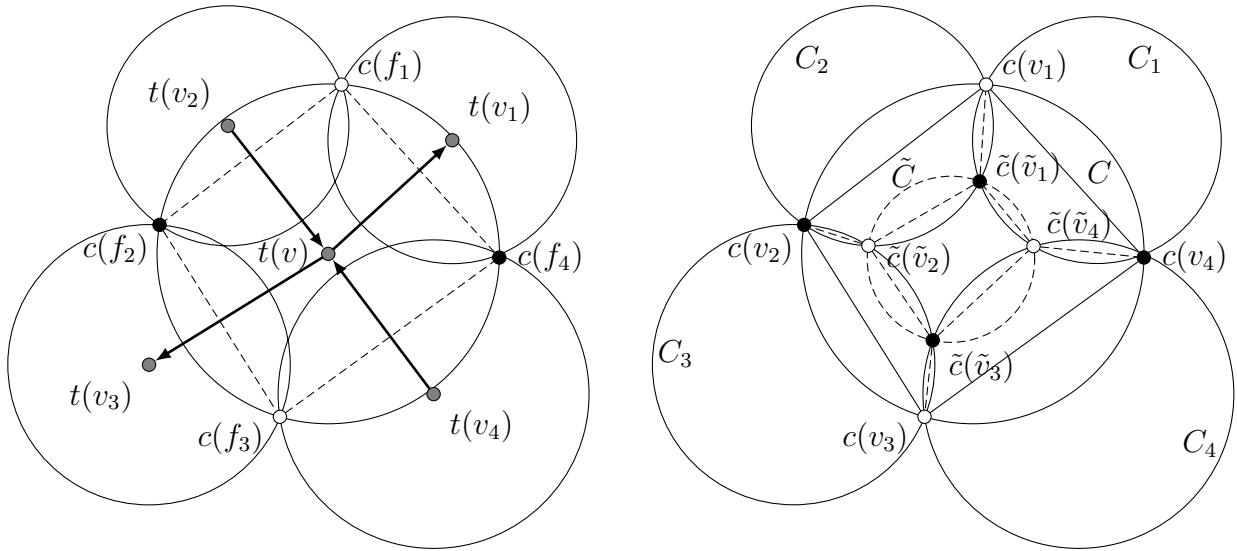

	\vspace{-1mm}

	\vspace{-1mm}
	\caption{Left: Labeling of t-embedding defined on $\qui$ (arrows), and circle pattern (circles) defined on $G$ (dashed) as in Lemma \ref{lem:cptotemb}. Right: Miquel's move, an example of $L$ (solid) and $\tilde L$ (dashed) as in Definition \ref{def:miquelmove}.}
	\label{fig:cptembmiquel}
\end{figure}

\begin{lemma}\label{lem:cptotemb}
	Let $G$ be a bipartite planar graph and let $\qui$ be the PDB quiver dual to $G$. Consider $\C$ as an affine chart of $\CP^1$ and let $c:V(G)\rightarrow \C$ be a circle pattern in that chart. Then the map $t: V(\qui) \rightarrow \C$ that consists of the circle centers of $c$ is a t-embedding.
\end{lemma}
\proof{
	Consider a vertex $v$ of $\qui$ with neighbours $v_1,v_2,\dots, v_{m}$. Because $\qui$ is a PDB quiver we have that $m\in 2\Z$. Let $f_1,f_2,\dots,f_m$ denote the faces of $\qui$ adjacent to $v$, such that $f_i$ is adjacent to $v,v_i$ and $v_{i+1}$ for all $i$, see Figure \ref{fig:cptembmiquel}. The intersection point $c(f_{i+1})$ is the reflection of $c(f_i)$ about the line $t(v)t(v_{i+1})$, because this line is the perpendicular bisector of $c(f_{i+1})$ and $c(f_i)$. Applying two consecutive reflections about the lines $t(v)t(v_{i+1})$ and $t(v)t(v_{i+2})$ maps $c(f_i)$ to $c(f_{i+2})$. Two consecutive reflections are a rotation around $t(v)$ about twice the angle between $t(v_{i+1})$ and $t(v_{i+2})$ with respect to $t(v)$. Concatenating all $m$ reflections yields a rotation that maps $c(f_1)$ to $c(f_1)$, and thus the angle of the rotation has to be a multiple of $2\pi$. Because the total rotation angle is the sum of every second angle between the centers, the angle condition of Definition \ref{def:tembedding} is satisfied.\qed
}

Note that if one chooses a different affine chart one gets a different map that is a t-embedding in that chart. However, a t-embedding constructed in one chart is not a t-embedding in another chart which has a different point at infinity. Still, both maps belong to a larger class of maps. We can think of a circle center as being the point at infinity reflected about the circle. Instead of infinity, we can fix another point and then consider all reflections about all circles of this point. The resulting points are called the conformal centers. In this way one could introduce a projective generalization of t-embeddings, and to each circle pattern there would correspond a (complex) one-parameter family of such generalized t-embeddings. We do not pursue this approach further.

We can also go in the other direction. We can construct a circle pattern from a t-embedding. To make this precise, we use the concept of a discrete connection. Recall that we denote by $\vec E = \{(v,v') : \{v,v'\} \in E\}$ the set of oriented edges of $G$.
\begin{definition}
	Let $S$ be a set and let $G$ be a graph. A \emph{discrete connection} $\gamma:\vec E \rightarrow \mbox{End}(S)$ is a map such that
	\begin{align}
		\gamma(v,v') \circ \gamma(v',v) = \mbox{id}
	\end{align}
	for all oriented edges $(v,v')\in \vec E$. A discrete connection $\gamma$ is called \emph{flat}, if for any cycle $(v_1,v_2,\dots,v_m,v_1)$ of $G$
	\begin{equation}
		\gamma(v_1,v_2) \circ \gamma(v_2,v_3) \circ \dots \circ \gamma(v_m,v_1) = \mbox{id}.\qedhere
	\end{equation}
\end{definition}

As we do not use non-discrete connections we refer to discrete connections simply as connections. 
We denote by $\mbox{Aff}(\C)$ the affine transformations of $\C$, which is a subgroup of $\mbox{End}(\C)$.

\begin{definition}
	Let $t:V(\qui)\rightarrow \C$ be a t-embedding and let $G$ be a graph such that $G^*=\qui$. We call $\varrho: \vec E(G) \rightarrow \mbox{Aff}(\C)$, such that $\varrho(v,v')$ is the reflection about the line $t(v^*)t(v'^*)$ the \emph{t-embedding connection}.
\end{definition}

\begin{lemma}\label{lem:uisflat}
	For any t-embedding $t:V(\qui)\rightarrow \C$ the t-embedding connection $\varrho: \vec E(G) \rightarrow \mbox{Aff}(\C)$ is flat.
\end{lemma}
\proof{
	The condition for a discrete connection is satisfied, because for any adjacent $v,v' \in V(G)$ the maps $\varrho(v,v') = \varrho(v',v)$ are both the same reflection, and therefore $\varrho(v,v') \circ \varrho(v',v)$ is the identity. It suffices to check flatness of $\varrho$ around every face of $G$, that is every vertex of $\qui = G^*$. But around a face $f$ of $G$ an even number of compositions of $\varrho$ correspond to a rotation around $t(f^*)$. The t-embedding angle condition in Definition \ref{def:tembedding} around $f^*$ is precisely the condition that the composition of $\varrho$ around $f$ is the identity.\qed
}

\begin{lemma}\label{lem:cpfromtemb}
	Consider a t-embedding $t: V(\qui)\rightarrow \C$ and its t-embedding connection $\varrho$, and identify the copy of the target space $\C$ of the t-embedding with the copy of $\C$ in $\mbox{Aff}(\C)$ which is the target space of the connection $\varrho$. Denote by $G$ the graph dual to $\qui$. Choose an initial point $c_0\in \C$ and a vertex $v_0$ of $G$. For every vertex $v$ of $G$ fix a path $v_0,v_1,\dots,v_m,v$. Define a map $c:V(G) \rightarrow \C$ by
	\begin{align}
		c(v) = \varrho(v_m,v) \circ \varrho(v_{m-1},v_m) \circ \dots \circ \varrho(v_0,v_1) (c_0).
	\end{align}
	Then $c$ is a circle pattern.
\end{lemma}
\proof{
	First of all, the definition is well-defined because $\varrho$ is a flat connection, see Lemma \ref{lem:uisflat}. Moreover, for a given face $f$ of $G$ with vertices $v'_1,v'_2,\dots,v'_{2m}$ on its boundary, we see that we get any $c(v'_k)$ from $c(v'_1)$ by a composition of reflections about lines that pass through $t(f^*)$. Thus the distance from $c(v'_k)$ to $t(f^*)$ is independent of $k$ and the points $c(v'_1)$ to $c(v'_{2m})$ are situated on a circle with center $t(f^*)$. \qed
}

As a consequence, there is a complex one-parameter family of circle patterns associated to any given t-embedding. Note that the \emph{origami map} \cite{clrtembeddings} is a map that is important in the study of t-embeddings and is also constructed from the t-embedding connection, but the construction is different from the construction in Lemma \ref{lem:cpfromtemb}.

Next we consider how to locally change combinatorics in circle patterns and t-embeddings. Again, it turns out that the key ingredient is a classic theorem of incidence geometry \cite{miquel} (see Figure \ref{fig:cptembmiquel}).

\begin{theorem}[Miquel's theorem]
	Let $C_1,C_2,C_3,C_4$ be four circles in $\CP^1$ such that each circle $C_k$ intersects $C_{k+1}$ in one or two points. Let $I = \bigcup_{k=1}^4 (C_k\cap C_{k+1})$ be the set of intersection points of consecutive circles. Assume there is a circle $C$ that intersects every circle $C_k$, and intersects them only in $I$. Then there is a second circle $\tilde C$ with the same property, such that $I$ is covered by $C \cup \tilde C$.
\end{theorem}

It is possible that the set $I$ in the theorem consists only of four points, in this case $C$ and $\tilde C$ coincide. We have changed the usual phrasing of Miquel's theorem to suit our purposes, by focusing on the simultaneous existence of $C$ and $\tilde C$, instead of the existence of the eighth intersection point in a cube. We employ Miquel's theorem to define a local change of combinatorics in circle patterns.

\begin{definition}\label{def:miquelmove}
	Let $c: V(G)\rightarrow \CP^1$ be a circle pattern and let $f$ be a face of degree four (a quad) in $G$. Consider a local change of combinatorics centered at $f$, as defined in the dimer model. See Figure \ref{fig:localmovedimer} for the four possible pairs of local configurations $L, \tilde L$, before and after the move respectively. Let $v_1,v_2,v_3$ and $v_4$ be the four vertices of $L\subset G$ adjacent to $f$ before the move, and let $\tilde v_1,\tilde v_2,\tilde v_3$ and $\tilde v_4$ be the four vertices of $\tilde L\subset \tilde G$ adjacent to $f$ after the move. Label the four faces of $G$ adjacent to $f$ by $f_1,f_2,f_3$ and $f_4$ such that each $f_i$ is adjacent to $v_i$ and $v_{i-1}$. The circle pattern $c$ associates five circumcircles with the faces $f,f_1,f_2,f_3,f_4$ that we identify with the circles $C,C_1,C_2,C_3,C_4$ of Miquel's theorem. Let $\tilde C$ be the sixth circle of Miquel's theorem. Define the circle pattern $\tilde c: V(\tilde G) \rightarrow \CP^1$, such that $\tilde c(\tilde v_i)$ is the intersection $C_i \cap C_{i+1} \cap \tilde C$, and such that $\tilde c$ agrees with $c$ everywhere else. The circle pattern $\tilde c: V(\tilde G) \rightarrow \CP^1$ is the result of the \emph{Miquel move at $f$} applied to $c$.
\end{definition}

Surprisingly, the effect of the Miquel move on circle centers does not depend on the actual circles, but only on the circle centers themselves.

\begin{theorem}\label{th:tmiquel}
	Let $t(v_1),t(v_2),t(v_3),t(v_4)$ be the centers of the circles $C_1,C_2,C_3,C_4$ in Miquel's theorem and let $t(v),\tilde t(v)$ be the centers of $C$ and $\tilde C$ in the same theorem. Then the centers satisfy the dSKP equation
	\begin{equation}
		\mr(t(v),t(v_1),t(v_2),\tilde t(v),t(v_3),t(v_4)) = -1.\qedhere
	\end{equation}
\end{theorem}
\proof{See \cite[Lemma 6.3]{amiquel} and \cite[Section 5.1]{kenyonlam}.\qed}

As a consequence of the theorem, we can define a local change of combinatorics for the t-embedding as well.

\begin{definition}\label{def:tembmut}
	Let $t: V(\qui)\rightarrow \C$ be a t-embedding and $v$ a vertex of $\qui$ of degree 4. Consider the quiver $\tilde \qui$ that appears after mutating at $v$, see Figure \ref{fig:localmovedimer}. Define the \emph{t-embedding after mutation at $v$} as the t-embedding $\tilde{t}: V(\tilde \qui)\rightarrow \C$, which agrees with $t$ everywhere except on $v$, where we define the image $\tilde t(v)$ as in Theorem \ref{th:tmiquel}.
\end{definition}

The next corollary states how the mutation in t-embeddings and the Miquel move in the corresponding circle patterns relate.

\begin{corollary}
	Consider $\C$ as an affine chart of $\CP^1$ and let $\qui$ be the PDB quiver that is dual to $G$. Let $v$ be a vertex of degree four of $\qui$ and $v^*$ the corresponding quad in $G$.
	\begin{enumerate}
		\item Let $c: G\rightarrow \CP^1$ be a circle pattern and let $t$ be the t-embedding consisting of the centers of $c$. Let $\tilde t$ (respectively $\tilde c$) be the t-embedding (circle pattern) after a mutation at vertex $v$ (face $v^*$). Then $\tilde t$ consists of the centers of $\tilde c$.
		\item Let $t: \qui\rightarrow \C$ be a t-embedding and let $c$ be a circle pattern constructed via the discrete t-embedding connection as in Lemma \ref{lem:cpfromtemb} with initial point $c_0$. Let $\tilde t$ (respectively $\tilde c$) be the t-embedding (circle pattern) after a mutation at vertex $v$ ($v^*$). Then the circle pattern constructed from $\tilde t$ with initial point $c_0$ coincides with $\tilde c$.\qedhere
	\end{enumerate}	
\end{corollary}

We have associated a quiver and also local changes of combinatorics to circle patterns and t-embeddings. The next goal is to show that t-embeddings are in fact special cases of TCD maps. Recall that the faces of a PDB quiver consist of two partitions: the faces that are oriented clockwise and those that are oriented counterclockwise.

\begin{lemma}
	Let $\qui$ be a PDB quiver. There is a TCD $\tcd$ that has $\qui$ as its affine quiver. 
\end{lemma}
\proof{We construct $\tcd$ by constructing the corresponding graph $\pb$. Add a white vertex $w_v$ to $\pb$ for every vertex $v$ of $\qui$. Into each clockwise face $f$ of degree $d_f$ glue a piece of graph $G_f$ that corresponds to a TCD with endpoint matching $\enm {d_f}{1}$ (see Definition \ref{def:endpointmatching}), such that the white boundary vertices of $G_f$ are glued to the white vertices $w_v$ that correspond to boundary vertices $v$ of $f$. The graph pieces $G_f$ do not contain interior white vertices and thus the white vertices of $\pb$ are in bijection with the vertices of $\qui$, as it should be. The affine quiver of a $G_f$ graph consists of a single clockwise oriented $d_f$-gon, which concludes the proof.  \qed 
}

Note that there is an alternate construction for $\pb$, by triangulating every clockwise face of $\qui$ and then replacing each such triangle with a black vertex and three white neighbours (the building block of TCDs). We observe in the construction that two TCDs that yield the same affine quiver only differ in a sequence of spider moves.

\begin{lemma}
	 Let $\tcd$ be a TCD with affine quiver $\qui$. Let $t: V(\qui)\rightarrow \C$ be a t-embedding and let $T:\tcdp \rightarrow \C$ be a TCD map that maps each white verticex $w_v$ of $\pb$ to the image $t(v)$. Then any resplit in $T$ at $w_v$ corresponds to a mutation in $t$ at $v$. 
\end{lemma}
\proof{The statement is true because Theorem \ref{th:tmiquel} and Equation \eqref{eq:resplitmr} coincide and determine the move uniquely.\qed}

As a consequence of the lemma, from now on we consider a t-embedding to be a TCD map as well and write $t: \tcdp \rightarrow \C$, as $\tcdp$ is in bijection with $V(\qui)$. 

\begin{remark}
	Let us add the following observation. Given the geometry of a t-embedding defined on $\qui$, we might as well define it on $\rho(\qui)$, the quiver with all arrows reversed, without affecting any of the previous statements in this section. The combinatorics of the two TCDs of the TCD maps corresponding to $\qui$ and $\rho(\qui)$ however are quite different. Of course, this is only possible in this case because all points of a t-embedding are projectively colinear in $\CP^1$ to begin with, and the black vertices of $\pb$ do not capture any incidence geometry.
\end{remark}

\begin{lemma}
	The cluster structure associated to a t-embedding $t:\tcdp \rightarrow \C$ in \cite{amiquel,kenyonlam} is the affine cluster of $t$.
\end{lemma}
\proof{The definition of the cluster structure of a t-embedding in \cite{amiquel,kenyonlam} is identical to the expression of the affine cluster structure via oriented differences as in Lemma \ref{lem:starratioviadistances}.\qed }

Note that a t-embedding, viewed as a TCD map actually possesses two cluster structures, an affine one and a projective one. However, for probabilistic purposes the affine cluster structure is of more interest because of the following lemma.

\begin{lemma}\label{lem:tembreal}
	The affine cluster variables of a t-embedding are real. Conversely, a dSKP lattice with real affine cluster variables is a t-embedding. Additionally, if the t-embedding is \emph{proper} \cite{clrtembeddings}, that is the faces are bounded by straight edges, convex and disjoint, then the affine cluster variables are real and positive.
\end{lemma}
\proof{The affine cluster variables of a TCD map $T:\tcdp \rightarrow \C$ are real if and only if the angle condition of Definition \ref{def:tembedding} is satisfied. Moreover, if the t-embedding is proper then the angle sum has to be $\pi$ which implies that the affine cluster variable is positive.\qed}

If the circle pattern is defined on a graph $G$ with $\Z^2$ combinatorics, then the corresponding dual quiver $Q$ also has $\Z^2$ combinatorics. Moreover, due to the $\Z^2$ combinatorics the circles can be partitioned into even and odd circles. A circle corresponding to a quad $((i,j),(i+1,j),(i+1,j+1),(i,j+1))$ is even if $(i+j)\in 2\Z$ and odd if not. For $\Z^2$ combinatorics it is possible to define global dynamics on circle patterns as well as t-embeddings.

\begin{definition}
	Consider a circle pattern $c^{(0)}: V(\Z^2) \rightarrow \CP^1$. Obtain the circle pattern $c^{(2k+1)}$ from the circle pattern $c^{(2k)}$ by performing the Miquel move at every odd circle, and the circle pattern $c^{(2k)}$ from the circle pattern $c^{(2k-1)}$ by performing the Miquel move at every even circle. This iteration defines \emph{Miquel dynamics}.
\end{definition}

Note that all circle patterns $C^{(k)}$ that are constructed in iterations of Miquel dynamics are still $\Z^2$ circle patterns. Miquel dynamics go back to an idea of Kenyon, and were first investigated by Ramassamy \cite{ramassamymiquel}. A first integrability result for Miquel dynamics was found by Ramassamy and Glutsyuk \cite{grmiquel}. Because of Theorem \ref{th:tmiquel} we see that on the level of t-embeddings, Miquel dynamics are a special case of dSKP lattices. However, up to this point, it was still unclear whether there is a cluster structure that can be associated to the circle patterns themselves, that is a cluster structure that can be read off from the intersection points without the need for an affine chart to determine circle centers. To show that such a cluster structure indeed exists, we begin with another classic theorem of incidence geometry, namely Clifford's four circle theorem \cite{clifford} and a relation to the dSKP equation. Let us use the notation $\bigcirc(a,b,c)$ for the unique circle through the threee different points $a,b,c\in \CP^1$.

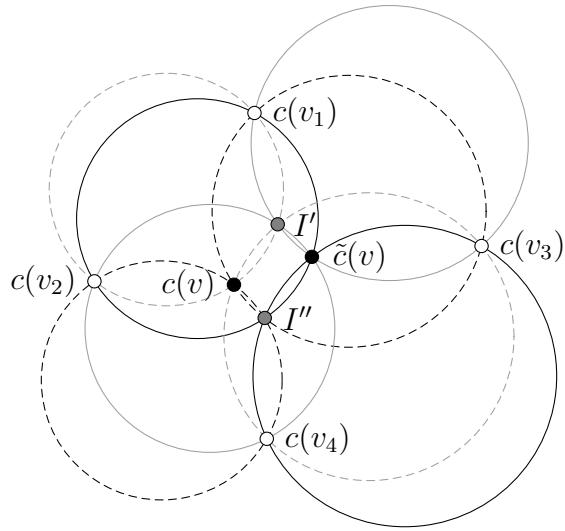
\begin{figure}
	\vspace{-1mm}
		\begin{tikzpicture}[line cap=round,line join=round,x=1.0cm,y=1.0cm,scale=0.9]
		\definecolor{aqaqaq}{rgb}{0.6274509803921569,0.6274509803921569,0.6274509803921569}
		\definecolor{cqcqcq}{rgb}{0.7529411764705882,0.7529411764705882,0.7529411764705882}
		\clip(-9,-0.52) rectangle (1,7.3);
		\draw [] (-5.774122657580919,4.0830834752981255) circle (1.7662813119035186cm);
		\draw [] (-2.741560665362035,1.7641903131115462) circle (2.2199606038087625cm);
		\draw [densely dashed] (-6.296087954204693,1.700111360088006) circle (1.760499291339022cm);
		\draw [densely dashed] (-3.5565359768143083,4.191581916824743) circle (2.002969706989656cm);
		\draw [color=aqaqaq] (-5.591048667439164,2.4659298957126303) circle (1.8260038098055804cm);
		\draw [color=aqaqaq] (-2.963026706231454,5.196913946587538) circle (2.026017930572848cm);
		\draw [densely dashed,color=aqaqaq] (-3.266361967731181,2.3444283727732063) circle (2.1199668158357827cm);
		\draw [densely dashed,color=aqaqaq] (-6.228272554597654,4.511595878128432) circle (1.712583439481506cm);
		\node[wvert,label={right:$c(v_1)$}] (c1) at (-4.94,5.64) {};
		\node[wvert,label={left:$c(v_2)$}] (c2) at (-7.28,3.16) {};
		\node[wvert,label={right:$c(v_3)$}] (c3) at (-1.62,3.68) {};
		\node[wvert,label={right:$c(v_4)$}] (c4) at (-4.76,0.84) {};
		\node[bvert,label={left:$c(v)$}] (c) at (-5.24,3.11) {};
		\node[bvert,label={[xshift=1]right:$\tilde c(v)$}] (c) at (-4.1,3.52)  {};
		\node[qvert,label={right:$I''$}] (i) at (-4.79,2.62)  {};
		\node[qvert,label={[xshift=-1]right:$I'$}] (i) at (-4.60,4.00)  {};
	\end{tikzpicture}
	\vspace{-2mm}
	\caption{Clifford's four circle theorem. The construction using the gray circles yields the same second black point as the construction using the black circles.} 
	\label{fig:clifford}
\end{figure}

\begin{theorem}[Clifford's four circle theorem]\label{th:cliff}
	Consider five different points $c(v)$, $c(v_1)$, $c(v_2)$, $c(v_3)$, $c(v_4)$ not on a common circle  in $\CP^1$ (see Figure \ref{fig:clifford}). Let $I',I''$ be such that
	\begin{align}
		\{c(v),I'\} &= \bigcirc(c(v),c(v_1),c(v_2)) \cap \bigcirc(c(v),c(v_3),c(v_4)),\\
		\{c(v),I''\} &= \bigcirc(c(v),c(v_2),c(v_3)) \cap \bigcirc(c(v),c(v_4),c(v_1)).
	\end{align}
	Then there is a unique sixth point $\tilde c(v)$ such that 
	\begin{align}
		\{\tilde c(v),I'\} &= \bigcirc(I',c(v_2),c(v_3)) \cap \bigcirc(I',c(v_4),c(v_1)),\hspace{8mm}\\
		\{\tilde c(v),I''\} &= \bigcirc(I'',c(v_1),c(v_2)) \cap \bigcirc(I'',c(v_3),c(v_4)).\qedhere
	\end{align}
\end{theorem}

The relationship between the dSKP equation and Clifford's four circle theorem stated below is due to Konopelchenko and Schief \cite{ksclifford}.

\begin{lemma}\label{lem:clifdskp}
	Let $c(v),\tilde c(v)$ and $c(v_1),c(v_2),c(v_3),c(v_4)$ be the six points from Clifford's four circle theorem. Then
	\begin{align}
		\mr(c(v),c(v_1),c(v_2),\tilde c(v),c(v_3),c(v_4)) = -1.
	\end{align}
\end{lemma}

We proceed by introducing an analogue of t-embeddings that is closely related to circle patterns.

\begin{definition}\label{def:uemb}
	Let $u: \tcdp \rightarrow \CP^1$ be a TCD map. We call $u$ a u-embedding if all its projective cluster variables are real.
\end{definition}

\begin{definition}\label{def:cptouemb}
	Let $c^{(0)}$ be a $\Z^2$ circle pattern and let $c^{(k)}$ for $k\in \Z$ be the circle patterns constructed via Miquel dynamics. Denote the odd circles of $c^{(2k)}$ by the map 
	\begin{align}
		C^{(2k)}: V_-(\Z^2) \rightarrow \{\mbox{Circles of }  \CP^1 \},
	\end{align} and the even circles of $c^{(2k)}$ by the map 
	\begin{align}
		C^{(2k+1)}: V_+(\Z^2) \rightarrow \{\mbox{Circles of } \CP^1 \}.
	\end{align}
	Let $u^{(k)}$ be a map from $V(\Z^2)$ to $\CP^1$ such that
	\begin{align}
		u^{(k)}(i,j) &= C^{(k)}(i,j) \cap C^{(k)}(i-1,j-1) \cap C^{(k+1)} (i-1,j) \cap C^{(k+1)} (i,j-1),
	\end{align}
	if $i+j+k\in 2\Z$ and otherwise
	\begin{align}
		u^{(k)}(i,j) &= C^{(k+1)} (i,j) \cap C^{(k+1)} (i-1,j-1) \cap C^{(k+2)}(i-1,j) \cap C^{(k+2)}(i,j-1).
	\end{align}
	Define a TCD $\tcd^{(k)}$ such that the set of white vertices of the associated graph $\pb^{(k)}$ is identified with $V(\Z^2)$, and for all $i,j$ with $i+j+k \in 2\Z$ glue the TCD
	\begin{center}

	\end{center}
	into the quad with lower left corner $(i,j)$. Then for any $k\in \N$ the map $u^{(k)}: \tcdp^{(k)}\rightarrow \CP^1$ is the \emph{u-embedding associated to $c^{(k)}$}.
\end{definition}

\begin{figure}
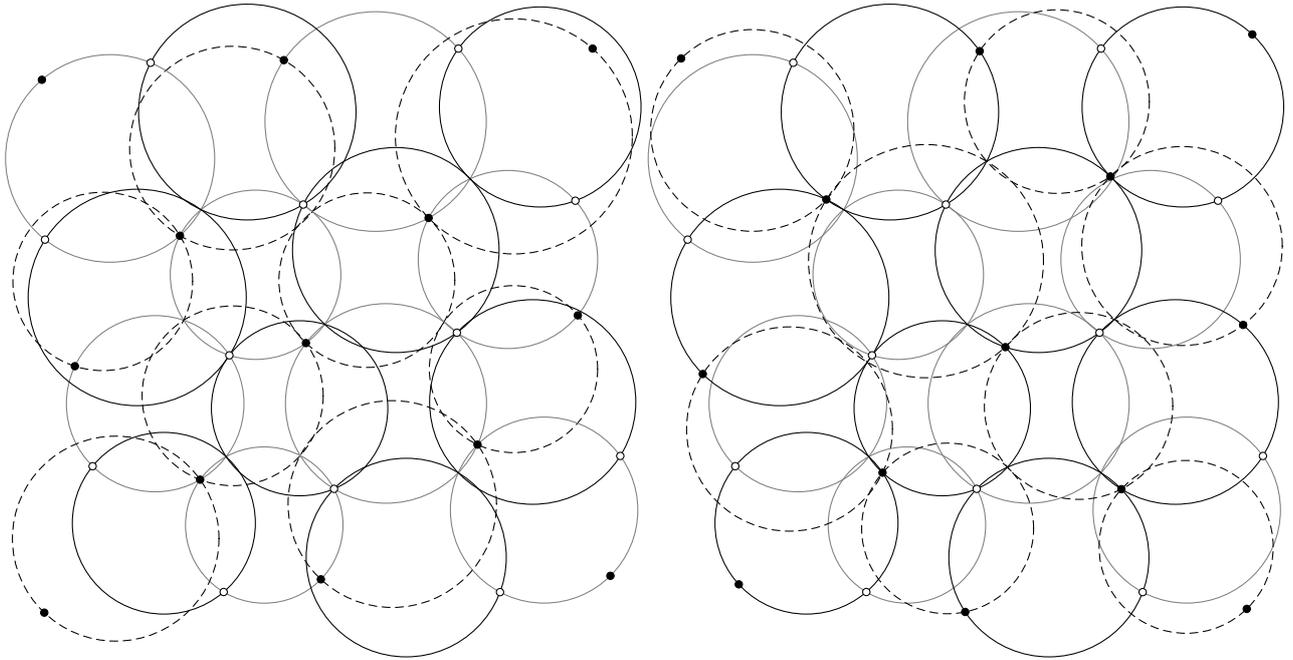

	
	%
	\caption{Circles related by Miquel dynamics, specifically the circles $C^{(k)}$ (dashed, only left), $C^{(k+1)}$ (gray), $C^{(k+2)}$ (black), $C^{(k+3)}$ (dashed, only right). Left: The even points ($i+j+k\in 2\Z$) of the u-embedding $u^{(k)}$ (white) as well as the odd points of $u^{(k)}$ (black). Right: The even points of $u^{(k+1)}$ (black) as well as the odd points of $u^{(k+1)}$ (white). Only the black points change from left to right.} 
	\label{fig:uembedding}
\end{figure}

See Figure \ref{fig:uembedding} for an example of $u^{(k)}$ and $u^{(k+1)}$ related by Miquel dynamics.

We will show in Lemma \ref{lem:cpanduemb} that the u-embeddings associated to circle patterns are indeed u-embeddings as defined in Definition \ref{def:uemb}. Moreover, there is in fact a second u-embedding that one could associate to $c$ by swapping the role of the even and the odd vertices in Definition \ref{def:cptouemb}. For our purposes it suffices to only consider the u-embeddings as defined by Definition \ref{def:cptouemb}. We also note that in Laplace-Darboux dynamics (see Figure \ref{fig:tcdlaplacedarboux}) we have already encountered the TCDs that correspond to $\Z^2$ combinatorics.

\begin{theorem}\label{th:uembandtwotwo}
	The associated u-embeddings $u^{(k)}$ and $u^{(k+1)}$ are related by a sequence of 2-2 moves.
\end{theorem}
\proof{
	Without loss of generality assume $k=0$. Then for $i+j\in 2\Z+1$ the points $u^{(k)}(i,j)$ and $u^{(k+1)}(i,j)$ agree. Moreover, the points 
	\begin{align}
		u^{(k)}(i,j),u^{(k)}(i-1,j),u^{(k)}(i,j+1) &\mbox{ are on } C^{(k+1)}_{i-1,j},\\
		u^{(k)}(i,j),u^{(k)}(i,j-1),u^{(k)}(i+1,j) &\mbox{ are on } C^{(k+1)}_{i,j-1},\\
		u^{(k+1)}(i,j),u^{(k)}(i+1,j),u^{(k)}(i,j+1) &\mbox{ are on } C^{(k+2)}_{i,j},\\
		u^{(k+1)}(i,j),u^{(k)}(i-1,j),u^{(k)}(i,j-1) &\mbox{ are on } C^{(k+2)}_{i-1,j-1}.
	\end{align}	
	Due to the definition of the associated u-embedding, the four circles in the list above intersect in a common point of the circle pattern $c^{(k+1)}$. We identify this configuration with the configuration of Clifford's four circle theorem. The common point of the four circles is $I'$; $u^{(k)}(i,j)$ and $u^{(k+1)}(i,j)$ are $c(v)$ and $\tilde c(v)$ respectively; the points $u^{(k)}(i-1,j),u^{(k)}(i,j+1),u^{(k)}(i+1,j)$ and $u^{(k)}(i,j-1)$ are $c(v_1)$, $c(v_2)$, $c(v_3)$ and $c(v_4)$ respectively. Due to Lemma \ref{lem:clifdskp} we know that the six involved points of $u^{(k)}$ and $u^{(k+1)}$ satisfy the dSKP equation. Thus the points are related by a resplit, see Lemma \ref{lem:resplitmr}. The first step to relate $u^{(k)}$ and $u^{(k+1)}$ as TCD maps is therefore to perform a resplit in $\pb^{(k)}$ at every even vertex of $\Z^2$. After a resplit in $\pb^{(k)}$ and the corresponding 2-2 move in $\tcd^{(k)}$, the resulting TCD is not yet $\tcd^{(k+1)}$. But it is, if we follow the resplits up by performing a spider move at every even face of $\Z^2$. This sequence relates $u^{(k)}$ to $u^{(k+1)}$ as TCD maps.\qed
}

\begin{lemma}\label{lem:cpanduemb}
 	Let $u^{(k)}$ be the maps from Definition \ref{def:cptouemb}. Every $u^{(k)}$ is a u-embedding.
\end{lemma}
\proof{
	Without loss of generality assume $k=0$.  We want to show that the projective cluster variables are real. In the odd faces of $\Z^2$ the bipartite graph $\pb$ look like the bipartite graph occurring in a spider move, thus the projective cluster variable in the center of an odd face is a cross-ratio of four points. As these four points are on a circle, the cross-ratio is real. In an even quad of $\Z^2$, the corresponding face of $\pb^{(0)}$ consists of eight vertices. However, two resplits (which do not change projective cluster variables) reduce the 8-point multi-ratio to a cross-ratio. These are exactly the resplits as discussed in the proof of Theorem \ref{th:uembandtwotwo}. Therefore we know that the four points occurring in the cross-ratio are actually four points on a circle in $u^{(1)}$, and thus the cross-ratio is also real.\qed
}

We have now found not one but two real cluster structures associated to Miquel dynamics. The t-embedding structure is invariant under affine transformations while the u-embedding structure is invariant under projective transformations. Note that in order to give the u-embeddings a probabilistic interpretation,  we have to consider u-embeddings where the cluster variables are also positive. The positivity corresponds to a non-intersecting ordering of the points of the u-embeddings on the circles $C^{(k)}(i,j)$. However, it is not immediately clear what the corresponding requirement for the original circle pattern is.

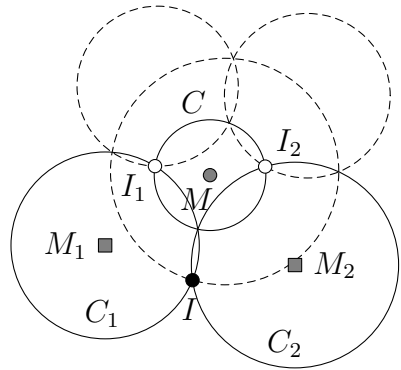
\begin{figure}
	\begin{tikzpicture}[line cap=round,line join=round,x=1.0cm,y=1.0cm,scale=0.4]
		\definecolor{ffffff}{rgb}{1.,1.,1.}
		\definecolor{yqyqyq}{rgb}{0.5019607843137255,0.5019607843137255,0.5019607843137255}
		\clip(-6.086666666666668,-8.393333333333333) rectangle (6.82,3.7590476190476165);
		\draw [] (0.54,-1.92) circle (1.835755975068582cm);
		\draw [densely dashed] (1.02,-1.8) circle (3.7544640096823403cm);
		\draw [densely dashed] (3.7398964420999006,0.7132147407080767) circle (2.7267367536296905cm);
		\draw [densely dashed] (-1.183183415059178,1.0297182347505474) circle (2.6520785446719892cm);
		\draw [] (-2.9108098938557383,-4.248228691322084) circle (3.09695017303653cm);
		\draw [] (3.334424142242881,-4.899354674051365) circle (3.408815738136784cm);
		\node[qvert,label={[yshift=2,xshift=-5]below:$M$}] at (0.54,-1.92) {}; 
		\node[wvert,label={[xshift=3,yshift=3]below left:$I_1$}] at (-1.2712256358537388,-1.6208985188498408) {};
		\node[wvert,label={[xshift=-2,yshift=-2]above right:$I_2$}] at (2.353456673813424,-1.6347371524338472) {};
		\node[bvert,label={[xshift=-1]below:$I$}] at (-0.03704106121086381,-5.4025913166655775) {};
		\node[qvert,rectangle,label=right:$M_2$] at (3.334424142242881,-4.899354674051365) {}; 
		\node[qvert,rectangle,label=left:$M_1$] at (-2.9108098938557383,-4.248228691322084) {}; 
		\node (c1) at (-3,-6.5) {$C_1$};
		\node (c2) at (3,-7.5) {$C_2$};
		\node (c) at (0,0.5) {$C$};
	\end{tikzpicture}
	\vspace{-3mm}
	\caption{The configuration of Lemma \ref{lem:centerintersecdksp}.} 
	\label{fig:centerintersecdskp}
\end{figure}

We proceed to show further relations between t- and u-embeddings.

\begin{lemma}\label{lem:centerintersecdksp}
	Consider three different circles $C$, $C_1$ and $C_2$ in an affine chart $\C$ of $\CP^1$ that intersect in a common point $I$, see Figure \ref{fig:centerintersecdskp}. Let the intersection sets between circles be
	\begin{align}
		C_1 \cap C = \{I,I_1\}, \qquad C_2 \cap C = \{I,I_2\}, \qquad C_1 \cap C_2 = \{I,I'\}.
	\end{align}
	Moreover, let $M,M_1,M_2$ be the centers of $C,C_1,C_2$. Then the dSKP equation
	\begin{align}
		\mr(I_1,M_1,I',M_2,I_2,M) = - 1
	\end{align}
	holds.
\end{lemma}
\proof{It is clear that the absolute value of the multi-ratio is 1, as we have
	\begin{align}
		|\mr(I_1,M_1,I',M_2,I_2,M)| = \frac{r_1r_2r}{r_1r_2r},
	\end{align}
	where $r_1,r_2,r$ are the radii of the circles $C_1,C_2,C$. Moreover, the dihedral angles at $I$ between the circles $C,C_1,C_2$ add up to $2\pi$. But the sum of the dihedral angles is also $\pi$ minus the angles appearing in the argument of the multi-ratio, thus proving that the argument of the multi-ratio is $\pi$.\qed
}

\begin{figure}
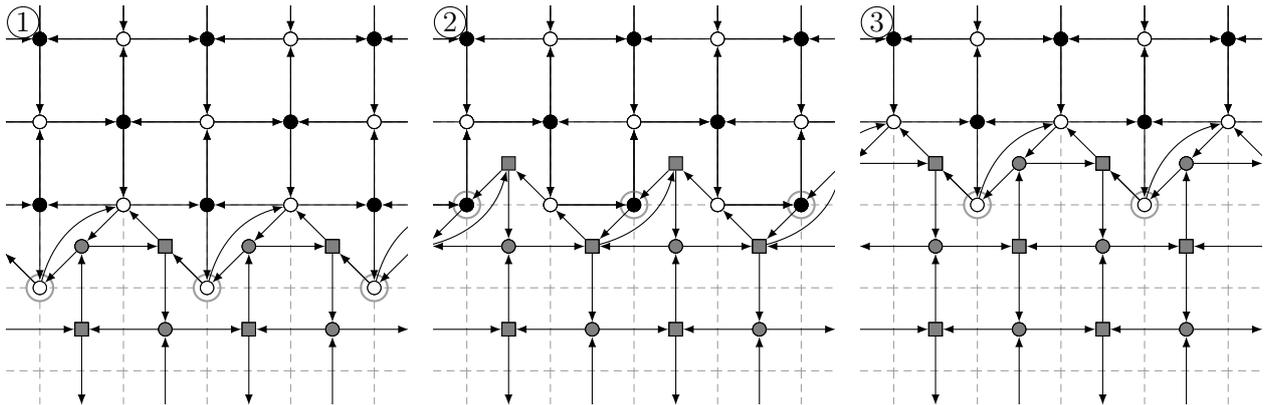


	\caption{Three steps in the sequence of mutations that turns a u-embedding into a t-embedding. White and black points correspond to intersection points of the u-embedding, while the gray nodes correspond to circle centers.}
	\label{fig:cpextension}
\end{figure}

\begin{theorem}\label{th:utextension}
	Let $c:V(\Z^2) \rightarrow \C$ be a circle pattern in an affine chart of $\CP^1$ and consider the t-embedding and the u-embedding associated to it. Then the t-embedding is an extension (see Definition \ref{def:extension}) of the u-embedding, and vice versa.
\end{theorem}
\proof{We show how the t-embedding data propagates through the u-embedding data in Figure \ref{fig:cpextension}, where we drew three iterations of the affine quiver. Two steps of mutations propagate the t-embedding data by one step in $\Z^2$. That mutations take the points of the u-embedding to the correct points of the t-embedding is a direct consequence of Lemma \ref{lem:centerintersecdksp}. Compare Figures \ref{fig:centerintersecdskp} and Figure \ref{fig:cpextension}, where we chose the same node-symbols to represent corresponding centers and points. This sequence of mutations moves a strand through $\Z^2$ and therefore corresponds combinatorially to a sweep and geometrically to an extension. Because the space is $\CP^1$, there is no geometric way to distinguish whether the u-embedding is the extension of the t-embedding or vice versa. This is reflected in the fact that we could just as well reverse all arrows in Figure \ref{fig:cpextension} and this would still be the affine quiver for a TCD map with the same set of white vertices that are mapped to the same points in $\CP^1$ but with a different TCD. Therefore we can consider the t-embedding as an extension of the u-embedding or vice versa.\qed 
}

Surprisingly, Theorem \ref{th:utextension} shows that associated u- and t-embeddings are related by a sequence of mutations. The combinatorics of all the iterations of a t-embedding live on an $A_3$ lattice, and the same is true for u-embeddings. Theorem \ref{th:utextension} implies that one can view the t-embedding, the associated u-embeddings and all their evolutions together as living on two $A_3$-slices of $A_4$. In this sense one could view the t-embedding as a discrete Bäcklund transformation of the u-embedding or vice versa. There is even a parameter involved that could be considered a spectral parameter. If we go from t-embedding to u-embedding, then we have to go from t-embedding to circle pattern, in which case we can choose one initial point arbitrarily. On the other hand, if we go from u-embedding to t-embedding, then in the step from circle pattern to t-embedding there is the choice of affine chart. This choice corresponds to the choice of the point at infinity, which in this direction could be considered to be the spectral parameter. As we have these different choices, we can look at an alternating sequence of t- and u-embeddings which all together form a map from $A_4$ that satisfies the dSKP equation.

There is one more approach to relate t- and u-embeddings. Recall that flags of TCD maps were defined in Definition \ref{def:tcdflag}.

\begin{lemma}
	Let $(T_2,T_1), (E_2,E_1,E_0)$ be a flag of TCD maps in $\CP^2$ and let $(T^\star_2,T^\star_1)$,  $(E^\star_2,E^\star_1,E^\star_0)$ be the dual flag of TCD maps in $(\CP^2)^*$. Choose an affine chart of $E_1$ such that $E_0$ is at infinity and an affine chart of $E_1^\star$ such that $E_0^\star$ is at infinity. Then
	\begin{enumerate}
		\item if $T_1$ is a u-embedding, $T^\star_1$ is a t-embedding in $E_1^\star$,
		\item and if $T_1$ is a t-embedding, $T^\star_1$ is a u-embedding in $E_1^\star$.\qedhere
	\end{enumerate}
\end{lemma}
\proof{
	We have characterized u-embeddings as maps that have real projective cluster variables (Definition \ref{def:uemb}) and t-embeddings as maps that have real affine cluster variables (Lemma \ref{lem:tembreal}). Moreover, we have shown in Theorem \ref{th:dualcluster} that
	\begin{align}
		\aff_{E_0}(T_1) = \pro(T_1^\star), \qquad \mbox{ } \qquad \pro(T_1) = \aff_{E_0^\star}(T_1^\star).
	\end{align}
	Thus the lemma is proven.\qed
}

Therefore, starting from a t- or u-embedding we can lift, dualize, and then take a section to obtain the other type of embedding. Of course, there is the freedom of choosing one dimensional boundary data involved in choosing the lift. On the other hand, Theorem \ref{th:utextension} shows that if $T_1$ is a t-embedding and $T_1'$ an associated u-embedding then there is a map $T_2$ such that $T_1=\sigma(T_2)$ is the section and $T_1' = \pi(T_2)$ the projection (recall the results of Section \ref{sec:tcdextensions}). Therefore both the affine cluster variables of $T_1$ and $T_2$ are real, as the projective variables of $T_2$ coincide with those of the u-embedding $T_1'$. Then when dualizing, also $T_2^\star$ and $T_1^\star$ have real affine cluster variables. Then the appropriate projection to $E_1^\star$ of $T_2^\star$ will be a u-embedding that extends the t-embedding $T_1^\star$. It would be very interesting to understand the map $T_2$ and $T_2^\star$ in more detail, including how to find the circle patterns that are dual to each other via this construction.

\section{Harmonic embeddings, h-embeddings and orthodiagonal maps} \label{sec:harmonicemb}

Harmonic embeddings were first considered by Tutte \cite{tutteembedding} and are also known as \emph{Tutte embeddings}. The relationship between harmonic embeddings and t-embeddings was investigated by Kenyon, Lam, Ramassamy and Russkikh \cite[Section 6]{kenyonlam}. We give some of their results without proof before we proceed to new results.

\begin{definition}
	Let $\qg$ be a quad-graph. A map $h: E(\qg) \rightarrow \C$ is an \emph{h-embedding} if for every quad with edges $(e^1,e^{2},e^{1}_2,e^2_1)$ in cyclic order, the points $h(e^1),h(e^{2}),h(e^{1}_2),h(e^2_1)$ form a rectangle.
\end{definition}

The term ``h-embedding'' does not appear in the literature. We introduce it here as h-embeddings are t-embeddings derived from harmonic embeddings, as we will explain in Remark~\ref{rem:harmonicemb}.

Note that the definition of h-embeddings is affine in $\C$, because the multiplication with a complex number and the addition of a complex number map rectangles to rectangles. We now consider the cube-flip in h-embeddings.

\begin{lemma}
	Let $\qg,\tilde \qg$ be two quad-graphs such that $\tilde \qg$ is $\qg$ after a cube flip. Let $h: E(\qg) \rightarrow \C$ be an h-embedding. Then there is a unique h-embedding $\tilde h: E(\tilde \qg) \rightarrow \C$ that only differs from $h$ on the three edges that are in $\tilde \qg$ but not $\qg$. 
\end{lemma}
\proof{
	See \cite[Section 6.3]{kenyonlam}.\qed
}

The following observation is also due to \cite[Theorem 23]{kenyonlam}.

\begin{lemma}
	Let $\qg$ be a quad-graph and $\qui_\qg$ be the cuboctahedral quiver (see Definition \ref{def:cuboquiver}) of $\qg$. Then an h-embedding $h: E(\qg) \rightarrow \C$ is also a t-embedding $t_h: V(\qui_\qg) \rightarrow \C$. 
\end{lemma}
\proof{
	Due to the Definition of the cuboctahedral quiver, $E(\qg) = V(\qui_\qg)$. Moreover, at any vertex of $\qui_\qg$ the two opposite angles that belong to the rectangles of $h$ sum to $\pi$, therefore satisfying the requirement of Definition \ref{def:tembedding}.\qed
}

If we consider an h-embedding $h$ as a t-embedding, then to every mutation of $\qui_\qg$ corresponds a mutation of the t-embedding as defined in Definition \ref{def:tembmut}. 

\begin{figure}
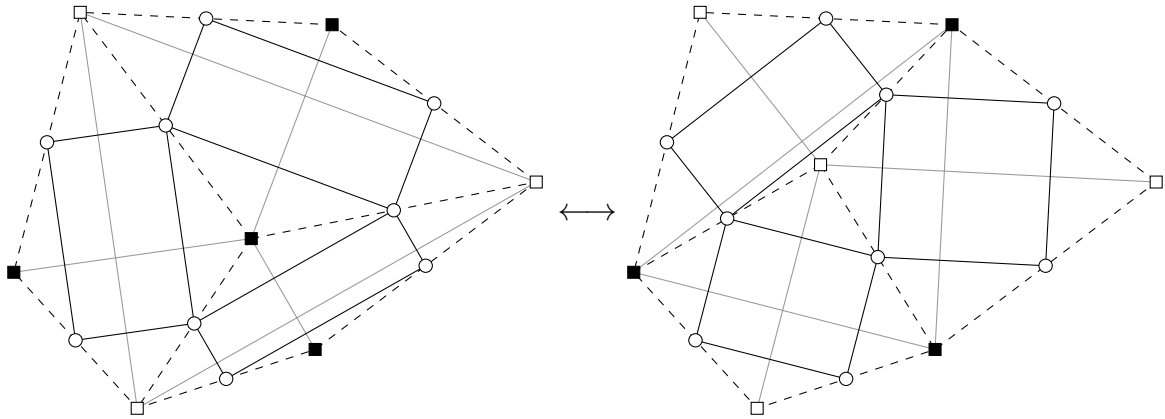


	\caption{The cube flip for an h-embedding (circle nodes, black lines) and a corresponding orthodiagonal map (square nodes and dashed lines).}
	\label{fig:harmandorthodiag}
\end{figure}

\begin{lemma}\label{lem:hembcuboctaflip}
	Let $t_h: V(\qui_\qg) \rightarrow \C$ be a h-embedding $h$ considered as a t-embedding. Apply the cuboctahedral flip (see Definition \ref{def:cuboflip}) to $\qg$ and the corresponding mutations to the t-embedding $t_h$. The result is the cube flip of the h-embedding (see also Figure \ref{fig:harmandorthodiag}).
\end{lemma}
\proof{See \cite[Section 6.3]{kenyonlam}.}

These lemmas together with our results from previous sections allow us to state a list of immediate corollaries.

\begin{corollary}
	Let $h: E(\qg) \rightarrow \C$ be an h-embedding. Then the following observations hold:
	\begin{enumerate}
		\item $h$ is a Darboux map.
		\item $h$ is a TCD map.
		\item The cube flip of $h$ corresponds to the cube flip of Darboux maps.
		\item $h$ is a Schief map.
		\item The affine cluster variables of $h$ are real.
		\item $h$ is affine BKP.		\qedhere
	\end{enumerate}
\end{corollary}
\proof{
	(1) Any map from the edges of a quad-graph to $\CP^1$ is a Darboux map. (2) A Darboux map is also a TCD map. (3) The cube flip of a Darboux map is defined as the sequence of mutations that belong to the cuboctahedral flip, thus the claim follows due to Lemma \ref{lem:hembcuboctaflip}. 
	(4) Rectangles are a special case of parallelograms, and parallelograms characterize Schief maps (Definition \ref{def:paralleldm}). (5) Because $h$ is also a t-embedding, the affine variables are real. (6) We have previously shown that Schief maps are affine BKP (Theorem \ref{th:schiefbkp}).\qed
}

The observations (1) to (4) are new. Observations (5) and (6) have already appeared in \cite[Secton 6]{kenyonlam}, we have only reformulated them in terms of the general definitions of cluster variables and cluster subvarieties for TCD maps.

Let us introduce another well known type of map, the orthodiagonal maps (see \cite{josefssonorthodiagonal}).

\begin{definition}
	Let $\qg$ be a quad-graph. An \emph{orthodiagonal map} is a map $o: V(\qg) \rightarrow \C$ such that the diagonals of every quad are orthogonal.
\end{definition}

\begin{definition}\label{def:horthocorrespondence}
 	Let $o: V(\qg) \rightarrow \C$ be an orthodiagonal map and $h: E(\qg) \rightarrow \C$ be a h-embedding. The maps $o$ and $h$ are \emph{in correspondence} if
	\begin{align}
		h(v,v') = \frac{o(v) + o(v')}{2},
	\end{align}
	for any edge $(v,v')$ of $\qg$.
\end{definition}

It is clear that for any given orthodiagonal map the corresponding h-embedding is determined by the definition. Conversely, given the h-embedding, we can choose one point of a corresponding orthodiagonal map arbitrarily and then the remaining points are determined. Indeed, there also exists a cube-flip for orthodiagonal maps (as shown in Figure \ref{fig:harmandorthodiag}), that is compatible with the cube-flip of h-embeddings (see \cite{kenyonlam}). 

\begin{remark}\label{rem:harmonicemb}
	For every quad-graph $\qg$ there are two dual graphs $G,G^*$, such that $\qg$ is the associated quad-graph $\qg_G$, see Definition \ref{def:graphtoqg}. Essentially, the vertices of $G$ (resp. $G^*$) consist of the white (black) vertices of $\qg$. Each edge of $G$ (and $G^*$) corresponds to a face of $\qg$. The restriction of an orthodiagonal map $o: V(\qg) \rightarrow \C$ to $G$ is a \emph{harmonic embedding} $U: V(G) \rightarrow \C$. The restriction to $G^*$ is a harmonic embedding $I: V(G^*) \rightarrow \C$ of the dual. The embeddings $U,I$ are called harmonic because they solve a discrete Laplace problem, see Remark \ref{rem:linresistorproblem} in the spanning tree section. The corresponding conductances $c_e$ for $U$ can be read off the orthodiagonal map as the ratio of diagonals
	\begin{align}
		c_e =  i\frac{o(v_1) - o(v_2)}{o(v) - o(v_{12})},
	\end{align}
	where $e$ corresponds to the edge $(v,v_{12})$ in $G$ and the quad $(v,v_1,v_{12},v_2)$ in $\qg$. The conductances for $I$ are $c_e^{-1}$.
\end{remark}

In Section \ref{sec:inscribedinconics} on Q-nets inscribed in conics we have seen that Schief maps arise as tangent sections of Q-nets inscribed in a conic in $\CP^2$. The case of h-embeddings is therefore a special case. Additionally, Definition \ref{def:horthocorrespondence} implies that the orthodiagonal map corresponds to the auxiliary map $m$ of Section \ref{sec:inscribedinconics}. There we showed that the auxiliary map is closely related to the stereographic projection of the corresponding Q-net, which has as section the Schief map. With this in mind, the next theorem is not so surprising.

\begin{figure}
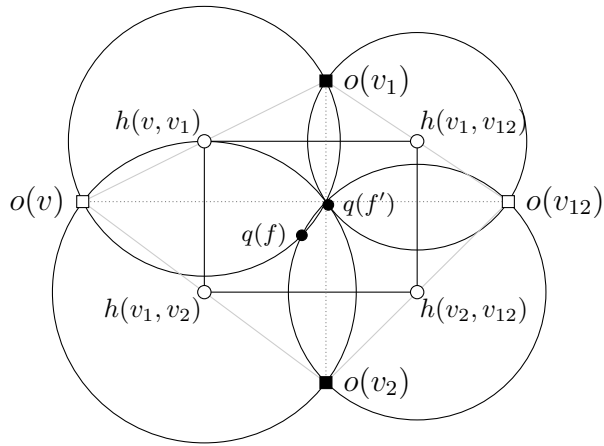

	\centering

	\caption{A rectangle of an h-embedding (white circle nodes) and the corresponding orthodiagonal quad (square nodes) including the two focal points (black circle nodes). The circles are the circles of the corresponding t-embedding.}
	\label{fig:harmandod}
\end{figure}

\begin{theorem}\label{th:odqnet}
	Let $o: V(\qg) \rightarrow \C \subset \CP^1$ be an orthodiagonal map. Extend $o$ to a map $q: V(\qg) \cup F(\qg) \rightarrow \C$ by requiring that $q$ agrees with $o$ on $V(\qg)$ and on every quad $f$ with white vertices $v,v_{12}$ and black vertices $v_1,v_2$ we require that
	\begin{align}
		q(f) &= \frac{q(v) q(v_2)-q(v_{12})q(v_1)}{q(v)+q(v_2)-q(v_{12})-q(v_1)}
	\end{align}
	holds, see also Figure \ref{fig:harmandod}. Then $q$ is a Q-net, which we call the \emph{corresponding orthodiagonal Q-net}, and the points assigned to the quads are the corresponding focal points. The other focal point in each quad $f$ as above is
	\begin{align}
		 q'(f) := \frac{q(v) q(v_1)-q(v_{12})q(v_2)}{q(v)+q(v_1)-q(v_{12})-q(v_2)}. 
	\end{align}
	Moreover, if $h$ is the h-embedding corresponding to $o$, then $q$ is an extension of $h$ in the sense of Definition \ref{def:extension}. 
\end{theorem}
\proof{Any map from $V(\qg) \cup F(\qg)$ to $\CP^1$ is a Q-net, as the incidence relations are trivially satisfied. However, the two focal points of a Q-net are related by a resplit and therefore also by a multi-ratio equation even in $\CP^1$. Indeed, a direct calculation shows that 
\begin{align}
	\mr(q(v),q(v_1),q(f),q(v_{12}),q(v_2),q(f')) = -1
\end{align}
holds, as required by Lemma \ref{lem:resplitmr}.
Another calculation shows that
\begin{align}
	\mr(q(f), h(v,v_1),q(v_1),h(v_1,v_2),q(v_{12}),h(v,v_2)) &= -1
\end{align}
is satisfied. Recall that an extension in $\CP^1$ is the projection of an extension in higher dimensions. In higher dimensions the last multi-ratio equation characterizes that the points $h(v,v_1), h(v_1,v_2)$ and $h(v,v_2)$ are the section of the triangle $q(f),q(v_1),q(v_{12})$ with a line, in fact this is Menelaus' theorem (see Theorem \ref{th:menelaus}). As the multi-ratio equation is preserved under projection, the claim is proven.\qed
}

There is also an alternative, more geometric, proof possible for the last theorem. To every quad $f$ in an orthodiagonal map $o$ we can assign the point of intersection of its two diagonals, let us call that point $D(f)$. Given an h-embedding $h$ corresponding to $o$, one can choose the circle pattern $c: V(\qg) \cup F(\qg) \rightarrow \C$ that has as intersection points the points of $o$ on vertices and the diagonal intersection points $d$ on faces of $\qg$. These circles are Thales circles, that is their centers are the centers of the segments of the orthodiagonal map. Therefore these circles pass through the diagonal intersection points, as the diagonals intersect orthogonally. Recall that in Definition \ref{def:cptouemb} we associate a u-embedding to a circle pattern with $\Z^2$ combinatorics. We argue now that $q$ is essentially the analogue to the u-embedding for cuboctahedral combinatorics. Indeed, the corresponding orthodiagonal Q-net $q$ consists of the intersection points of $c$ on $V(\qg)$ and of alternate intersection points of circles on $F(\qg)$ (see Figure \ref{fig:harmandod}). To see the latter, let $C: E(\qg) \rightarrow \{\mbox{Circles of } \CP^1\}$ be the circles of the circle pattern $c$, then one can verify by calculation that 
\begin{align}
	\{D(f),q(f)\} &= C(v,v_2) \cap C(v_1,v_{12}),\\
	\{D(f),q'(f)\} &= C(v,v_1) \cap C(v_2,v_{12}).
\end{align}
With this characterization the first multi-ratio equation in the proof of Theorem \ref{th:odqnet} is a consequence of Clifford's four circle theorem (Theorem \ref{th:cliff}), and the second multi-ratio equation is a consequence of Lemma \ref{lem:centerintersecdksp}. 

It turns out that for an orthodiagonal map $o$ the corresponding orthodiagonal Q-net $q$ is in fact a u-embedding, as a consequence of the following, stronger theorem.

\begin{theorem}\label{th:orthoharmperfdual}
	Let $\C$ be an affine chart of $\CP^1 = \C \cup \{P^\infty\}$. Let $o:V(\qg) \rightarrow \C$ be an orthodiagonal map and let $h,q$ be the corresponding h-embedding and the orthodiagonal Q-net respectively. Consider $h$ and $q$ as TCD maps, then
	\begin{equation}
		\aff_{P^\infty}(h) = \rho(\pro(q)) \qquad \mbox{ and } \qquad \aff_{P^\infty}(q) = \pro(h),
	\end{equation}
	where $\rho$ is the reciprocal cluster structure, see Definition \ref{def:reciprocalcluster}.
\end{theorem}
\proof{Looking at each quad of $\qg$ separately we verify that indeed the affine quiver of $h$ is the projective quiver of $q$ with arrows reversed and the affine quiver of $q$ is the projective quiver of $h$. Moreover, expressing the points of $h$ and the focal points of $q$ via the points of $o$ we verify by calculation that
\begin{align}
	\frac{h(v,v_2) - h(v,v_1)}{h(v_1,v_{12}) - h(v,v_1)} = \frac{o(v_1) - o(v_2)}{o(v) - o(v_{12})} = - \frac{q(v_1) - q(f)}{q(f) - q(v_2)}
\end{align}
holds in every quad of $\qg$. The ratio in the center term is nothing but the ratio of the diagonals of the orthodiagonal quads. The affine cluster variables of $h$ are composed of ratios as in the left most term (see Lemma \ref{lem:starratioviadistances}), while the projective cluster variables of $q$ are composed of ratios as in the right most term (see Lemma \ref{lem:projclusterviadistances}). This proves that $\aff_{P^\infty}(h) = \rho(\pro(q))$ is true. On the other hand, Theorem \ref{th:odqnet} states that $q$ is an extension of $h$ and therefore by Theorem \ref{th:affprojcluster} we obtain that $\aff_{P^\infty}(q) = \pro(h)$ is also true. \qed
}

\begin{corollary}
	Let $o$ be an orthodiagonal map and $q$ the corresponding orthodiagonal Q-net. Then $q$ is a u-embedding.
\end{corollary}
\proof{Because $\pro(q) = \aff(h)$ up to reciprocity, where $h$ is the corresponding h-embedding and the variables of $\aff(h)$ are real.\qed}

Note the close similarity between the concept of perfect duals (see Section \ref{sec:perfdual}) and Theorem \ref{th:orthoharmperfdual}. On the level of combinatorics though the difference is that in one of the equations the quiver is not reversed. On the level of geometry,  Theorem \ref{th:orthoharmperfdual} is also not quite an instance of Theorem \ref{th:perfdualinevenflags} of perfect duals of two maps $T,T'$, where $T$ lives in some primal projective space and $T'$ lives in the corresponding dual space. Instead both $h$ and $o$ live in the primal space.

\begin{remark}
	Let us sketch another curious property of h-embeddings. Assume for a moment that $h$ is \emph{proper}, that is all the images of the quads are oriented counterclockwise in $\C$. Then it turns out that the projective cluster variables of $h$ satisfy the following equations
	\begin{align}
		|X_v|^2\prod_{f\sim v} (1+X^{f}) &= 1, \qquad  \mbox{ for every black non-boundary vertex $v$ of } \qg,\\ 
		|X_v|^2\prod_{f\sim v} (1+(X^{f})^{-1})^{-1} &= 1, \qquad  \mbox{ for every white non-boundary vertex $v$ of } \qg,
	\end{align}
	where $X^f$ is the cluster variable in the center of quad $f$. These equations are almost the equations for the Ising subvariety (see Definition \ref{def:isingsub}), except that we take the absolute value squared instead of just the square of $X_v$ in the equations. The proof is also a simple calculation that gives these equations a geometric interpretation. Introduce variables for the length $l_f$ and width $b_f$ of every rectangle of $h$ and the associated conductance $c_f = \frac{b_f}{l_f}$. The projective cluster variable of $h$ in every quad is 
	\begin{align}
		X^f = \frac{b_f^2}{l_f^2} = c_f^2.
	\end{align}
	Moreover, for a black vertex $v$ of $\qg$ we have that
	\begin{align}
		|X_v|^2 &= |\mr(h(v,v_1),h(v_2,v_{12}),h(v,v_2),\dots,h(v_{1},v_{d_v1}))|^2\\ &= \prod_{q\sim v}^m\frac{b_f^2+l_f^2}{l_f^2} = \prod_{f\sim v}^m\left(1+X^f\right)^{-1}
	\end{align}
	holds. Analogously in the case of white vertices. So far, we have not shown that $X_v$ is completely determined by the conductances. However, trigonometry yields that the argument of $X_v$ is
	\begin{align}
		\arg(X_v) = \sum_{q \sim v}\left( \frac\pi2 - \arctan c_f\right) = \sum_{f \sim v} \left( \frac\pi2 - \arctan \sqrt{X^f} \right).
	\end{align}
	Thus indeed, $X_v$ is completely determined by the conductances. As the projective quiver of $h$ is hexahedral, and all positive real conductances can be realized as h-embedding, the recursion for conductances under a cube flip is not only described by a cuboctahedral cluster structure but \emph{also} by a hexahedral cluster structure. It is unclear, if there is any probabilistic use for the hexahedral approach, as the variables $X_v$ will generally not be real positive. On the other hand, one can consider a reduction of conductances such that the corresponding h-embeddings have $X_v$ real positive everywhere. Indeed, the extensively investigated case of \emph{isoradial embeddings} is a special case of such a reduction. We note that in isoradial embeddings the conductances are $\tan \frac{\theta}{2}$ \cite{kenyonisoradial}, while the Ising weights are $\tan^2 \frac{\theta}{2}$ \cite{btcriticalisoradialising}. In fact, the conductances are the affine cluster variables $Y_q$ of the h-embedding associated to the isoradial embedding and the Ising weights are the projective cluster variables $X_q$. On the other hand, we can view the isoradial embedding itself as a orthodiagonal Q-net $q$ and then due to Theorem \ref{th:orthoharmperfdual} the affine cluster variables of $q$ are the Ising weights and the projective cluster variables are the conductances.
\end{remark}

\section{S-embeddings}\label{sec:sembeddings}

\begin{definition}
	An \emph{s-embedding} is a map $s: V(\qg) \rightarrow \C$ such that each quad has an incircle.
\end{definition}
The notion of s-embedding is due to Chelkak \cite{chelkaksembeddings}, generalizing a notion of Smirnov \cite{smirnovtowards} that was introduced to obtain conformal invariance results of certain Ising observables \cite{csanalysisisoradial, csuniversality}. As in the case of h-embeddings, s-embeddings are an object of affine geometry (in $\C$). We call an s-embedding \emph{proper} if no two line segments associated to edges of $\qg$ intersect (except at the vertices).

\begin{lemma}[\cite{mrtalphaquads}]
	Let $\qg,\tilde \qg$ be two quad-graphs such that $\tilde \qg$ is $\qg$ after a cube flip. Let $s: E(\qg) \rightarrow \C$ be a proper s-embedding. Then there is a unique proper s-embedding $\tilde s: E(\tilde \qg) \rightarrow \C$ that only differs from $s$ on the three edges that are in $\tilde \qg$ but not $\qg$.
\end{lemma}

If $s$ is not proper then it is possible that there is no or more than one proper $\tilde s$, and even if $s$ is proper then there may be more than one (not necessarily proper) $\tilde s$. However, by relating $s$-embeddings to $t$-embeddings, there is a canonical cube-flip even in the non-proper case.

\begin{definition}
	Let $s: V(\qg) \rightarrow \C$ be an s-embedding. Define a map $t_s: V(\qg) \cup F(\qg) \rightarrow \C$ that agrees with $s$ on $V(\qg)$ and maps $F(\qg)$ to the corresponding incircle centers. We call $t_s$ the t-embedding corresponding to $s$, as justified by the next lemma.
\end{definition}

\begin{lemma}[\cite{kenyonlam}]
	Let $\qg$ be a quad-graph and $\qui_\qg$ be the hexahedral quiver (see Definition \ref{def:hexaquiver}) of $\qg$. Then the t-embedding corresponding to an s-embedding is indeed a t-embedding. 
\end{lemma}
\proof{Consider the star-ratio at an incircle center. Every angle is the sum of two angles that are $\frac{\pi}{2}$ minus half of one of the corner angles of the quad. Therefore the sum of the even angles and the sum of the odd angles are equal and thus both are equal to $\pi$. At a vertex of $\qg$ the two angles that come from the same quad are equal and therefore the star-ratio at vertices of $\qg$ are also real.\qed
}

\begin{figure}
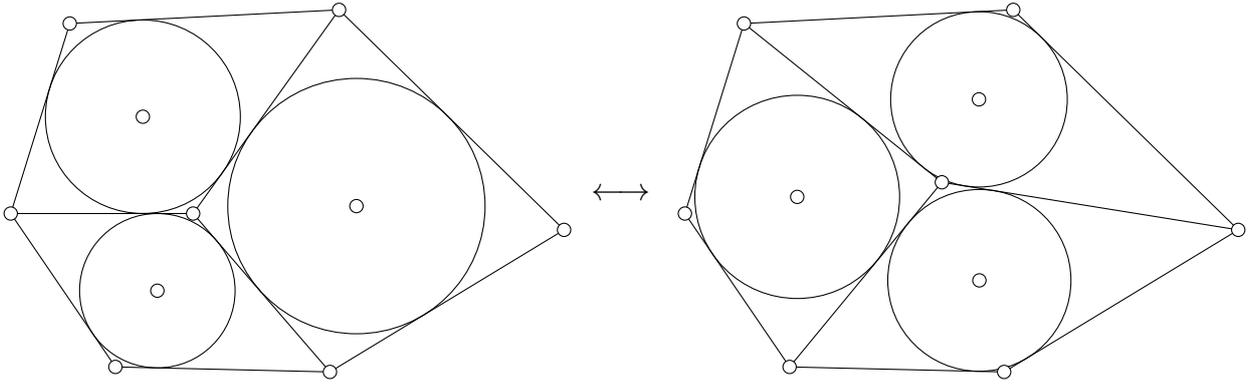

	
	\caption{The cube flip in an s-embedding.}
	\label{fig:sembcubeflip}
\end{figure}

\begin{lemma}
	Let $t_s: V(\qui_\qg) \rightarrow \C$ be a t-embedding corresponding to an s-embedding. Apply the hexahedral flip (see Definition \ref{def:hexaflip}) to $\qg$ and the corresponding mutations to the t-embedding $t_s$. The result is the cube flip of the corresponding s-embedding (see also Figure \ref{fig:sembcubeflip}).
\end{lemma}
\proof{See \cite[Remark 3.12]{mrtalphaquads}.\qed}

These lemmas together with our results from previous sections allow us to state a list of immediate corollaries.

\begin{corollary}
	Let $s$ be an s-embedding and $t_s$ the corresponding t-embedding.
	\begin{enumerate}
		\item $t_s$ is a Q-net.
		\item $t_s$ is a TCD map.
		\item The cube flip of $s$ corresponds to the cube flip of Q-nets.
		\item The cluster structure of s-embeddings \cite{kenyonlam} corresponds to the affine cluster structure of $t_s$ as a TCD map.
		\item The affine cluster variables of $t_s$ are real.\qedhere
	\end{enumerate}
\end{corollary}
\proof{Any map from the vertices and faces of a quad-graph to $\C \subset \CP^1$ is a Q-net (1). A Q-net is also a TCD map (2). The cube flip of a Q-net is defined as the sequence of mutations that belong to the hexahedral flip (3). The cluster structure of s-embeddings is the cluster structure of $t_s$, which is the affine cluster structure of the corresponding TCD map (4). Because $t_s$ is a t-embedding, the affine variables are real (5).\qed
}

The new result in this corollary is that an s-embedding can therefore be viewed as a projection of a Q-net in higher dimensions. In the case of h-embeddings the corresponding corollary also stated that h-embeddings are affine BKP. The proof was an immediate consequence of the fact that h-embeddings are a special case of Schief maps. We claim that s-embeddings are affine CKP, and a proof using direct calculations in $\C$ is possible \cite[Sectino 7]{kenyonlam}. We give an alternate proof that is based on incidence geometry in Corollary \ref{cor:sembckp}.

\begin{definition}
	Let $q: V(\qg) \rightarrow \CP^1$ be a Q-net. We call $q$ a \emph{fixed focal map} if every quad has the same two focal points.
\end{definition}
This definition is purely projective. We use shift notation in the following calculation to keep the equation readable. A calculation shows that the focal points of a quad in a fixed focal map have to be
\begin{align}
	\fp^{kl} = \frac{q q_{kl} - q_k q_l \pm \sqrt{(q-q_k)(q_k-q_{kl})(q_{kl}-q_l)(q_l-q)} }{q-q_k+q_{kl}-q_l},
\end{align}
as they are the fixed points of the involution that maps  $q\leftrightarrow q_{12}$ and $q_1\leftrightarrow q_2$. Let us now relate fixed-focal maps to s-embeddings.

\begin{theorem}\label{th:fixandincircular}
	Every s-embedding $s:V(\qg) \rightarrow \C \subset \CP^1$ is a fixed focal map. A fixed focal map is an s-embedding if its star-ratios are real.
\end{theorem}
\proof{
	We argue per quad. If a quad is incircular, then we know we can choose circles centered in the vertices that touch each other cyclically exactly in the touching points of the incircle with the quad edges. Thus the incircle is replaced with itself by Miquel dynamics. Therefore the center of the incircle is a fixed-point. 
	Conversely, consider a quad of a fixed-focal map with real star-ratio at the fixed focal point $\fp$. Let us denote the points $q,q_k,q_{kl},q_l$ by $p_1,p_2,p_3,p_4$ in order to identify them with the Miquel configuration. Then we know we can consider the four vertices $q_1,q_2,q_3,q_4$ of the quad and $\fp$ as centers of a Miquel configuration. Here, $\fp$ is also the center of the sixth circle in the Miquel configuration. In order to choose the circles, we need to pick one of the intersection points. Pick the point $P$ on the line (in $\C$) through $q_1,q_2$ closest to $\fp$ as the intersection point of the circles $C_1,C_2$ resp. $C$ centered at $q,q_1$ resp. $\fp$. Then the circles $C_1$ and $C_2$ are touching. However, $P$ is now also the intersection point of the sixth circle $\tilde C$, and therefore $\tilde C$ and $C$ coincide. As a consequencem every circle $C_i$ touches the circle $C_{i+1}$ and the quad therefore has an incircle.\qed
}

Fixed focal maps are a link that relates s-embeddings to S-graphs, which we introduced in Section \ref{sec:sgraphs}.

\begin{lemma}\label{lem:fixedfocalsgraph}
	The projection $\pi_{P\rightarrow L}: \CP^2 \rightarrow L$ of an S-graph $q$ from $P$ to a line $L$ that does not contain $P$ is a fixed focal map. Conversely, given a line $L \subset \CP^2$ and a point $P\notin L$ every Q-net $q:\Z^3 \rightarrow \CP^2$ such that $\pi_{P\rightarrow L}(q)$ is a fixed focal map is an S-graph.
\end{lemma}
\proof{Because of the colinearity of any two focal points of a quad with $P$ the focal points are mapped to the same point in $L$ by $\pi_{P\rightarrow L}$.\qed}

\begin{corollary}\label{cor:fixedfocalckp}
	Every fixed focal map $q: V(\qg) \rightarrow \CP^1$ is affine CKP with respect to any point $R\in \CP^1$.
\end{corollary}
\proof{Let $\CP^1 = L \subset \CP^2$ and choose a point $P\notin L$ and let $H$ be the line through $P$ and $R$. By Lemma \ref{lem:fixedfocalsgraph} there is an S-graph $\hat q: V(\qg) \rightarrow \CP^2$ with respect to $P$ such that $q$ is the projection of $\hat q$ from $P$. By Theorem \ref{th:sgraphcarnot} $\hat q$ is affine CKP with respect to $H$. Therefore $q$ is affine CKP with respect to $H \cap L = R$, due to our findings on cluster variables and projections (Theorem \ref{th:projectionandaffine}). \qed
}

\begin{corollary}\label{cor:sembckp}
	Every s-embedding $s:V(\qg) \rightarrow \C$ resp. the corresponding t-embedding $t_s$ is affine CKP with respect to any point in $\CP^1$, in particular also with respect to the point at infinity.
\end{corollary}
\proof{Direct consequence of the fact that an s-embedding is a fixed focal map (Theorem \ref{th:fixandincircular}) and Corollary \ref{cor:fixedfocalckp}.\qed}

In the case of t-embeddings defined on $\Z^2$ and h-embeddings we observed that they naturally give rise to u-embeddings. The same happens in the case of s-embeddings.

\begin{definition}
	Let $s:V(\qg) \rightarrow \C$ be an s-embedding and let $t_s: V(\qui_\qg) \rightarrow \C$ be the corresponding t-embedding, where $\qui_\qg$ is the hexahedral quiver. Fix two adjacent vertices $v_0,v_0'$ of $V(\qg)$. Define a circle pattern $c: V(\qui_\qg^*) \rightarrow \C$ with centers $t_s$ by choosing the initial intersection point of the circles centered at $t_s(v_0)$ and $t_s(v_0')$ on the line $s(v_0)s(v_0')$ (in $\C$). For any two adjacent vertices $v,v' \in V(\qg)$, the two corresponding circles only intersect in one point on the line $s(v)s(v')$. Therefore these points define a map $l_s: E(\qg) \rightarrow C$.
\end{definition}

The fact that the intersection points of the circles in $c$ for adjacent vertices $v,v'$ coincide, is a consequence of the fact that every quad in $s$ possesses an incircle. This property guarantees that
\begin{align}
	|s(v)-s(v_1)| - |s(v_1)-s(v_{12})| + |s(v_{12})-s(v_2)| - |s(v_2) - s(v)| = 0,
\end{align}
for every quad $(v,v_1,v_{12},v_2)$. This can be viewed as a consistency equation for 
\begin{align}
	|s(v)-s(v_i)| = r(v) + r(v_i),
\end{align}
where $r$ is a real function on the vertices of $\qg$. The function $r$ encodes where the points of $l_s$ are. More specifically, given a solution $r$ we can choose
\begin{align}
	|l_s(v,v_i) - s(v)| = r(v)
\end{align}
everywhere.

We have little doubt that the next statement can be proven by direct calculation, but have not found a geometric proof. Therefore we leave it as a conjecture.

\begin{conjecture}
	Given an s-embedding $s:V(\qg) \rightarrow \C$, the map $l_s: E(\qg) \rightarrow C$ is a u-embedding and a line complex that extends $t_s$, where we view $t_s$ as a Q-net.
\end{conjecture}


\chapter{Cluster $\tau$-variables}\label{cha:tau}

\section{Cluster $\tau$-variables for strongly generic TCD maps}\label{sec:tautcd}

We have mentioned $\tau$-variables briefly in our introduction to cluster variables, see Definition \ref{def:mutationtau}. However, we have not considered $\tau$-variables in relation to TCD maps. In the case of the $X$-variables we aimed to give a comprehensive treatment of the relations between the $X$-variables and various geometric operations and examples. In contrast, in the case of $\tau$-variables we only aim to give an introduction with the intention to enable future research. The three goals we want to accomplish with our definition of $\tau$-variables for TCD maps are that the variables
\begin{enumerate}
	\item satisfy the $\tau$-mutation rule (Definition \ref{def:mutationtau}) for spider-moves and are invariant under resplits,
	\item are compatible with the $X$-variables associated to TCD maps (via Equation \eqref{eq:tauandx}),
	\item reproduce the Fock-Goncharov $\tau$-variables \cite{fghighertm, fgtwoflags, gwebs} in the case of projective flag configurations.
\end{enumerate}

Recall that in Section \ref{sec:tcdconsistency} we assigned a label $\an'(f)$ to each face of a labeled TCD $\tcd$, such that $\an'(f)$ contains a strand $s$ if $s$ is to the left of $f$ in $\tcd$, see also Equation \eqref{eq:shiftamap}. In this section it is more practical to consider the complementary labeling. That is, we denote by $\ban(f)$ the maximal set of strands such that $f$ is to the right of every strand in $\ban(f)$. Because faces and white vertices of the graph $\pb$ correspond to faces of $\tcd$, we also write $\ban(w)$ and $\ban(f)$ for white vertices and faces of $\pb$ respectively.

We introduce $\tau$-variables only with fairly strong genericity assumptions, not because we know that these are necessary but because they are sufficient and because we want to avoid going into more technical detail than is necessary for an introduction.

\begin{definition}\label{def:stronglygeneric}
	Let $\tcd$ be a balanced TCD with $n$ strands and maximal dimension $k$, that is, the endpoint matching of $\tcd$ is $\enm {n}{k}$. To the right of each strand are $k$ white boundary vertices, and every white vertex $w$ is to the right of $k$ strands. We call a TCD map $T:\tcdp \rightarrow \CP^k$ \emph{strongly generic} if 
	\begin{enumerate}
		\item for every strand $i$, the images of the $k$ consecutive white boundary vertices to the right of the strand span a hyperplane $H_i$,
		\item the intersection of any set of $k$ distinct hyperplanes $H_{i_1},H_{i_2},\dots,H_{i_k}$ is a point, and
		\item for every white vertex $w$, the image $T(w)$ is the unique intersection of the $k$ hyperplanes $H_{j}$ for $j\in \ban(w)$.\qedhere
	\end{enumerate}
\end{definition}

The definition contains the three properties that we employ below. However, the second property actually already implies the third. This is because due to Lemma \ref{lem:rightofstrand}, everything to the right of strand $i$ is contained in $H_i$.

\begin{definition}\label{def:tcdtau}
	Let $T: \tcdp \rightarrow \CP^k$ be a strongly generic TCD map with $n$ strands. For $1\leq i \leq n$ we identify each hyperplane $H_i$ with a point $H^\perp_i$ in the dual space $(\CP^k)^*$ and fix a homogeneous lift of $\hat H^\perp_i$ in $\C^{k+1}$. Then the \emph{$\tau$-variables} of $T$ are the function $\tau: F(\pb) \rightarrow \C$, such that for every face $f\in F(\pb)$
	\begin{align}
		\tau_f = \det \left[\hat H^\perp_{j_1}, \hat H^\perp_{j_2}, \dots, \hat H^\perp_{j_{k+1}}  \right],
	\end{align}
	where $j_1, j_2, \dots, j_{k+1}$ are the strands in $\ban(f)$ sorted in increasing order.
\end{definition}

A different choice of homogeneous lifts for the hyperplanes $H_i$ leads to $\tau$-variables that are \emph{gauge related}. We use the term `gauge related' in the same way as for edge-weights: two choices of $\tau$-variables are gauge related if they yield the same $X$-variables. Moreover, the effect of rescaling the homogeneous lift $\hat H_i^\perp$ of a hyperplane $H_i$ for fixed $i$ by a factor $\mu$ is that all $\tau$-variables to the right of $i$ are also scaled by $\mu$. If we apply a projective transformation to the points of $T$ in $\CP^k$, then the effect on the dual space is also a projective transformation. Because the $\tau$-variables are defined as determinants of points in the dual space, the effect of a projective transformation is just a scaling of all $\tau$-variables, which is a gauge transformation.

\begin{lemma}
	The $\tau$-variables of Definition \ref{def:tcdtau} are invariant under the resplit and satisfy the mutation rules of Definition \ref{def:mutationtau}.
\end{lemma}
\proof{
	The invariance of every variable $\tau_f$ under a resplit is clear, because $\ban(f)$ is invariant for every $f\in F(\pb)$. Moreover, if we perform a spider move at face $f$ then every variable $\tau_{f'}$ for $f'\neq f$ is also invariant because $\ban(f')$ is invariant. It remains to verify that $\tau_f$ changes according to the mutation rule. Keeping Lemma \ref{lem:strandorder} in mind we assume without loss of generality that the four strands involved in the spider move are labeled $a,b,c,d$ with $a<b<c<d$, and assume that the face $f$ in $\tcd$ is bounded by $a$ and $c$ before the resplit. There are $k-1$ further strands $i_1,i_2,\dots i_{k-1}$, such that $\ban(f)$ is $\{a,c,i_1,\dots,i_{k-1}\}$. Let $M$ denote the $(k-1)\times(k+1)$-matrix with column vectors $\hat H^\perp_{i_1}, \hat H^\perp_{i_2}, \dots, \hat H^\perp_{i_{k-1}}$. We want to show that
	\begin{align}
		\det\left[\hat H^\perp_a, \hat H^\perp_c,M\right] \det\left[\hat H^\perp_b, \hat H^\perp_d, M \right]=&\phantom{+} \det\left[\hat H^\perp_a, \hat H^\perp_b,M \right]\det\left[\hat  H^\perp_c,\hat  H^\perp_d,M \right]\\&+\det\left[\hat H^\perp_a, \hat H^\perp_d, M \right]\det\left[\hat H^\perp_b, \hat H^\perp_c, M \right]\small\nonumber
	\end{align}
	holds. This is the mutation rule for $\tau_f$, except that we shuffled the column vectors $\hat H^\perp_a$, $\hat H^\perp_b, \hat H^\perp_c$, $\hat H^\perp_d$ to the left in each determinant, which only changes the overall sign of the equation. On the other hand, this is the simplest case of Plücker's determinant identities, which concludes the proof.\qed 
}

\begin{lemma}
	Let $T: \tcdp \rightarrow \CP^k$ be a maximally generic TCD map with $n$ strands. The $\tau$-variables of Definition \ref{def:tcdtau} are consistent with the projective $X$-cluster variables via Equation \eqref{eq:tauandx}. In other words, for every face $f\in F(\pb)$
	\begin{align}
		X_f = \frac{\prod_{(f',f)\in \vec E(\qui)}\tau_{f'} }{\prod_{(f,f')\in \vec  E(\qui)}\tau_{f'}},
	\end{align}
	where the products are over incoming and outgoing edges in the projective quiver $\qui$. 
\end{lemma}
\proof{
	Assume $f$ is a face of degree $m$, and let the strands bounding $f$ be denoted by $q_1,q_2,\dots,q_m$ in ascending order. For $1\leq i\leq m$ let $r_i$  denote the strand that shares a triple crossing with $q_{i}$ and $q_{i+1}$ at a black vertex on the boundary of $f$. Let $s_1,s_2,\dots,s_{k+1-m}$ be the strands that have $f$ to its right and are none of the aforementioned strands. Let $I=\{q_1,q_2,\dots,q_m, s_1, s_2, \dots, s_{k+1-m}\}$. We use Lemma \ref{lem:projclusterviadistances} to express $X_f$ via a multi-ratio. Therefore, we want to show the identity
	\begin{align}
		\mr\left(T_{I\setminus \{q_1\}}, T_{I \cup \{r_1\} \setminus \{q_1,q_2\} }, T_{I\setminus \{q_2\}}, T_{I \cup \{r_2\} \setminus \{q_2, q_3\}}, \dots, T_{I \cup \{r_m\} \setminus \{q_m, q_1\}}\right)
		&= \prod_{i=1}^m \frac{\det(M_i)}{\det(M'_i)},\label{eq:mrdualdetidentity}
	\end{align}
	where $M_i$ (resp. $M'_i$) is the matrix of column vectors 
	\begin{align}
		\left(\hat H^\perp_{r_i}, \hat H^\perp_{q_1}, \hat H^\perp_{q_2}, \dots, \hat H^\perp_{q_m},  \hat H^\perp_{s_1}, \hat H^\perp_{s_2},\dots, \hat H^\perp_{s_{k+1-m}}\right),
	\end{align} 
	with column $\hat H^\perp_{q_{i+1}}$ (resp. $\hat H^\perp_{q_{i}}$) deleted. Note that we have already reordered some terms in the determinant, but this reordering does not change the sign of each quotient appearing on the right-hand side. Both the left-hand and the right-hand side are invariant under choices of lifts and projective transformations. Thus we can assume that the $m$ vectors $\hat H^\perp_{q_i}$ together with the $k+1-m$ vectors $\hat H^\perp_{s_j}$ form a basis of the homogeneous coordinate space $\C^{k+1}$. The determinants in each quotient involve the column vectors associated to strands $s_i$. Therefore we can restrict ourselves to the space spanned by vectors $\hat H^\perp_{q_1}$ to $\hat H^\perp_{q_m}$, and replace vectors $\hat H^\perp_{r_i}$ for $1\leq i\leq m$ with their respective orthogonal projections to that space. Next, we forget about the strands $s_i$ and we replace $M_i$ (resp. $M'_i$) by the matrix 
	\begin{align}
		\left(\hat H^\perp_{r_i}, \hat H^\perp_{q_1}, \hat H^\perp_{q_2}, \dots, \hat H^\perp_{q_m}\right),
	\end{align} 
	with column $\hat H^\perp_{q_{i+1}}$ (resp. $\hat H^\perp_{q_{i}}$) deleted. For ease of notation, we write $T_{q_i}$ for $T_{I\setminus \{q_i\}}$ and $T_{r_i}$ for $T_{I\setminus \{q_i,q_{i+1}\}\cup \{r_i\} }$, for $1\leq i\leq m$. We choose lifts of $\hat H^\perp_{q_1}, \hat H^\perp_{q_2}, \dots, \hat H^\perp_{q_m}$ as a basis of $(\C^{m+1})^*$ and such that $\det[\hat H^\perp_{q_1}, \hat H^\perp_{q_2}, \dots, \hat H^\perp_{q_m}] = 1$. We also choose lifts $\hat T_{q_1}, \hat T_{q_2},\dots, \hat T_{q_m}$ as  basis for $\C^m$ and such that $\det[\hat T_{q_1}, \hat T_{q_2}, \dots, \hat T_{q_m}]=1$. For $1\leq i\leq m$ we introduce the functionals 
	\begin{align}
		\dot H^\perp_{q_i}: \C^m\rightarrow \C, \quad  \dot H^\perp_{q_i}(\cdot)&= \det[\hat T_{q_1}, \hat T_{q_2}, \dots,\xcancel{ \hat T_{q_i}},\dots \hat T_{q_m},\cdot],\\
		\dot T_{q_i}: (\C^m)^\perp\rightarrow \C,\quad \dot T_{q_i}(\cdot)&= \det[\hat H^\perp_{q_1}, \hat H^\perp_{q_2}, \dots, \xcancel{\hat H^\perp_{q_i}},\dots \hat H^\perp_{q_m},\cdot],
	\end{align}
	where the cross denotes that we skip this column. Any other vector $P$ in $\C^m$ or $(\C^m)^\perp$ also defines a functional $\dot P$ by linear superposition. The two bases formed by the vectors $\hat H^\perp_{q_i}$ and by the vectors $\hat T_{q_i}$ are dual in the sense that for any vector $\hat P\in \C^m$ (resp. $\hat Q\in (\C^m)^\perp$) holds
	\begin{align}
		\hat P = \sum_{i=1}^m \sigma_i \dot H^\perp_{q_i}(\hat P)\hat T_{q_i} \quad\mbox{ and }\quad
		\hat Q = \sum_{i=1}^m \sigma_i \dot T_{q_i}(\hat Q)\hat H^\perp_{q_i},
	\end{align}
	where $\sigma_i$ are in $\{\pm1 \}$. We readily check that, as a consequence any $\hat Q\in(\C^m)^*$ satisfies
	\begin{align}
		\dot Q(\hat T_i) = \dot T_i(\hat Q)
	\end{align}
	for all $1\leq i \leq m$. Let us now reformulate the left-hand and right-hand side of Equation \eqref{eq:mrdualdetidentity} in terms of such functionals. We recall that we can express a multi-ratio in terms of coefficients appearing in linear relations (as discussed in Section \ref{sec:projcluster}). The relations that appear on the boundary of the face that we are currently looking at are
	\begin{align}
		\hat T_{r_i} = \sigma_i \dot H^\perp_{r_i}(\hat T_{q_{i+1}})\hat T_{q_i} - \sigma_{i+1}\dot H^\perp_{r_i}(\hat T_{q_{i}}) \hat T_{q_{i+1}},
	\end{align}
	for $1\leq i \leq m$. By combining the sign factors the left-hand side of Equation \eqref{eq:mrdualdetidentity} becomes
	\begin{align}
		\prod_{i=1}^m \frac{\dot H^\perp_{r_i}(\hat T_{q_{i+1}})}{\dot H^\perp_{r_i}(\hat T_{q_{i}})}.
	\end{align}
	On the other hand, the right-hand side of Equation \eqref{eq:mrdualdetidentity} is nothing but
	\begin{align}
		\prod_{i=1}^m \frac{\dot T_{q_{i+1}}(\hat H^\perp_{r_i})}{\dot T_{q_{i}}(\hat H^\perp_{r_i})},
	\end{align}
	which concludes the proof.\qed
}

\begin{remark}
	It is not hard to see that for a TCD map $T$, the hyperplanes $H_s$ associated to strands, viewed as points in projective dual space are actually the boundary points of the projective dual TCD map $T^\star$ (as discussed in Section \ref{sec:projduality}). Indeed, $H_s = U_1(w_s)$, where $U$ is the subspace map (Definition \ref{def:subspacemap}) and $w_s$ is the unique white vertex in $\pb_1$ such that $\ban(w_s) = s$. Due to Lemma \ref{lem:subspacedual} $T^\star_k(w_s^\star) = H_s^\perp$.  Thus the points that we use to define the $\tau$-variables via determinants are actually boundary points of the projective dual TCD. It would be very interesting to further pursue these connections, as well as the interplay with taking sections. It would also be interesting to investigate whether there is a relation to the work of Muller and Speyer \cite{mstwist}, as some of the operations appearing in their work are reminiscent of projective duality. As a consequence, the $Z$-variables that we introduce in Section \ref{sec:dimerinvariants} hypothetically coincide with the $\tau$-variables of the dual TCD map, and the alternating ratios of $Z$-variables coincide with the $X$-variables of the dual TCD map.
\end{remark}

\begin{remark}
	Let us use the $\tau$-variables to introduce a particular gauge for the edge-weights of a strongly generic TCD map $T$. Even though we have not yet found any use for it, we think it is worth mentioning. Every edge $e\in E(\pb)$ is incident to two faces $f,f'$. The particular gauge for the edge-weights is $\lambda_e= \tau_f\tau_f'$. Why does this gauge always exist? We can use the alternative construction algorithm mentioned in Remark \ref{rem:altconstruction} to obtain a TCD map $T'$ with the given edge-weights, such that $T'$ also assumes its maximal dimension. By construction, $X(T') = X(T)$ and therefore $T$ and $T'$ are related by a projective transformation due to Theorem \ref{th:uniquefrominvariants}. Thus, also the points $H^\perp_i(T)$ and $H^\perp_i(T')$ are related by a projective transformation. Consequently, any choice $\tau(T')$ is gauge equivalent to $\tau(T)$, and we can indeed assume that $\lambda_e= \tau_f\tau_f'$ for each edge.
\end{remark}

We have satisfied our first two goals: to show that our definition of $\tau$-variables is compatible with our $X$-variables, and with mutation. Before we move on to discuss projective flag-configurations, let us discuss a bit how to extend the definition of $\tau$-variables to a more general setup. Let us first consider the case of a TCD map $T$ defined on a balanced TCD, but not realizing its maximal dimension. If $T$ is the projection $\pi(\hat T)$ of a strongly generic TCD map that does realize its maximal dimension, then we can set the $\tau$-variables of $T$ to the $\tau$-variables of $\hat T$. Because $X$-variables are invariant under projection and therefore $X(T) = X(\hat T)$, the copied $\tau$-variables on $T$ are compatible with the $X$-variables on $T$ as well. 

Another interesting question is how to extend the definition for $\tau$-variables to TCDs that are not balanced. Let us assume that the TCD attains maximal dimension. Then we can recall Lemma \ref{lem:subsweepable}, in the proof of which we showed that every TCD $\tcd$ can be iteratively expanded until it is a balanced TCD $\hat \tcd$. On the level of TCD maps, this corresponds to iteratively adding marked points on the boundary to obtain a map $\hat T: \hat\tcdp \rightarrow \CP^k$. This procedure does not change the $X$-variables of $T$ defined on $\tcdp$. In $T$, the points to the right of any strand span at most a hyperplane, spanning a hyperplane in the generic case. We can use these hyperplanes to define the $\tau$-variables for $\hat T$, and by restriction for $T$ as well. We will see an example of this procedure in Section \ref{sec:fgmoduli}.

Together, the two procedures -- reverse-projecting and adding points at the boundary -- should allow the introduction of $\tau$-variables for most TCD maps, although a more precise analysis would be interesting. Even with these procedures though, we are still only looking at TCD maps defined on minimal TCDs in the disc. In Section \ref{sec:fgmoduli} we will also look at an example where this procedure is possible in a non-minimal case, but a general setup in the non-minimal case is still unclear. Another very interesting question is the case of minimal TCDs defined not on a disc, but for example on the torus. A natural attempt would be to look at the corresponding TCD map defined on the universal cover of the torus. The maximal dimension of this TCD will be infinite, possibly one should try to look at limits of determinants on exhaustions of the torus.

\section{Projective flag configurations}\label{sec:fgmoduli}

In this section we want to cast the case of projective flags studied by Fock and Goncharov \cite{fghighertm} in terms of TCD maps. Because we limit ourselves to the case of projective flags, we introduce the objects in our own definitions, even though we stay close to follow-up work of Fock and Goncharov \cite{fgtwoflags}. Also, note that both the TCDs that we introduce and the corresponding bipartite graph already occur in work by Goncharov \cite{gwebs}. 

\begin{definition}\label{def:projkflagconfig}
	Let $\tg$ be a triangulation, and $k\in \N$. Let $\fl_k$ denote the set of projective $k$-flags $F_0 \subset F_1 \subset \dots \subset F_{k-1}$ of $\CP^k$, with $\dim F_i = i$ for $0\leq i < k$. A \emph{projective $k$-flag configuration} is a map $F: V(\tg) \rightarrow \fl_k$ such that for $0\leq i < k$ and any two different vertices $v,v'\in V(\tg)$ the dimensions satisfy $\dim (F_i(v) \cap F_i(v')) = i-1$.
\end{definition}

The assignment of $X$- and $\tau$-variables will only require that two flags intersect appropriately if the two corresponding vertices are incident in $\tg$. However, the definition above will avoid degeneracy problems when we consider edge-flips of the triangulation.

\begin{definition}
	A \emph{decorated projective $k$-flag configuration} is a projective $k$-flag configuration $F: V(\tg) \rightarrow \fl_k$ decorated by
	\begin{enumerate}
		\item a map $H^\perp: V(\tg) \rightarrow ((\CP^k)^*)^k$ such that for every $i \in \{0,1,\dots,k-1\}$ and every vertex $v\in V(\tg)$, the corresponding hyperplanes $H_j(v)$ satisfy $F_i(v) = \cap_{j=i}^{k-1} H_j(v)$,
		\item and a homogeneous lift $\hat H^\perp_j(v) \in (\C^{k+1})^*$ for every $i \in \{0,1,\dots,k-1\}$ and every vertex $v\in V(\tg)$. \qedhere
	\end{enumerate}

\end{definition}

\begin{definition}
	Let $\tg$ be a triangulation. The \emph{associated TCD $\tcd^k(\tg)$ of order $k$} and the corresponding graph $\pb^k(\tg)$ are obtained by gluing 
	\begin{center}
		
	\end{center}
	into each triangle of $\tg$ respectively.
\end{definition}

By definition, the TCD map $\tcd_\tg$ we associated to a triangulation in Section \ref{sec:tcdtriangulations} coincides with the associated TCD of order 1, that is, with $\tcd^1(\tg)$. As in the case of triangulations, $\tcd^k(\tg)$ is minimal if and only if $\tg$ has no interior vertices. Note that every strand turns clockwise around a vertex $v$ of $\tg$. Conversely, for each vertex $v$ of $\tg$, there are exactly $k$ strands of $\tcd^k(\tg)$ that turn clockwise around $v$. The $k$ strands around each vertex are nested. Let us denote the label of the most nested strand by $s_0(v)$, and the other strands by $s_i(v)$, for $0 \leq i < k$, in nesting order. Moreover, let us label the white vertices of the corresponding graph $\pb$ by strand labels as in Section \ref{sec:tcdconsistency}. We observe that every white vertex $w\in V(\pb)$ carries $k$ labels of strands from at most three vertices $v,v',v''$. If $w$ carries labels from three different vertices $v,v'$ and $v''$, then $(v,v',v'')$ is a triangle of $\tg$. If $w$ carries labels from two vertices $v$ and $v'$, then $(v,v')$ is an edge of $\tg$. Finally, if $w$ carries labels from only one vertex $v$, then $w$ coincides with $v$.

\begin{definition}\label{def:projflagconfigtcdmap}
	Let $F$ be a projective $k$-flag configuration on a triangulation $\tg$. The \emph{associated TCD map} $T_F:\tcd^k(\tg) \rightarrow \CP^k$ is the unique TCD map, such that if the strand $s_i(v)$ turns around $w$ clockwise then $T_F(w)$ is contained in $F_i(v)$ for all $i \in \{0,1,\dots,k-1\}$ and $v\in V(\tg)$.	
\end{definition}

It is straightforward to check that the associated TCD map above is indeed a TCD map. We have to check that the images of the white vertices of $\pb$ are uniquely defined, that the TCD map is 0-generic, and that the three points around black vertices of $\pb$ lie on a line. For example, if there are strands of three different vertices $v,v',v''$ turning around $w$, then 
\begin{align}
	T_F(w) = F_{s_{i}(v)}(v) \cap F_{s_{i'}(v')}(v') \cap F_{s_{i''}(v'')}(v''),
\end{align}
where $i,i',i''$ are the indices of the minimal strands in nesting order turning clockwise around $w$. Due to the combinatorics of $\tcd_\tg^k$, the sum $s_{i}(v) + s_{i'}(v') + s_{i''}(v'')$ equals $2k$, and thus the intersection of the three corresponding spaces $F_{s_{i}(v)}(v), F_{s_{i'}(v')}(v'), F_{s_{i''}(v'')}(v'')$ is indeed a point. The argument in the case of strands of two or one different vertex turning around $w$ works analogously. Moreover, for every black vertex $b$ of $\pb$ there are exactly three strands $s_{i}(v),s_{i'}(v'),s_{i''}(v'')$ passing through $b$. Let $w_1,w_2,w_3$ be the three white vertices adjacent to $b$. Then indeed 
\begin{align}
	T_F(w_1), T_F(w_2), T_F(w_3) \in F_{s_{i+1}(v)}(v) \cap F_{s_{i'+1}(v')}(v') \cap F_{s_{i''+1}(v'')}(v''),
\end{align}
where we employ the convention $F_k = \CP^k$. The right hand side is a line, because the sum $s_{i}(v) + s_{i'}(v') + s_{i''}(v'')$ equals $2k-2$. Finally, the 0-genericity of the $T_F$ follows from the genericity assumptions of Definition \ref{def:projkflagconfig}.

If we begin with a decorated projective $k$-flag, we can replace $F_i(v)$ in Definition \ref{def:projflagconfigtcdmap} with $H_i(v)$ to obtain an equivalent definition. Therefore, the equation
\begin{align}
	T_F(w) = \bigcap_{v\in V(\tg)}\left( \bigcap_{i\in I_v(w)} H_i(v)\right),
\end{align}
holds as well, where $I_v(w)$ is the set of nesting indices of strands that turn clockwise around $v$ and $w$. This equation shows that the points of $T_F$ can be obtained as unique intersections of hyperplanes associated to strands as required in Definition \ref{def:stronglygeneric}, even though the hyperplanes are not uniquely defined by the boundary points of $T_F$ but depend on the decoration of $T_F$.

\begin{figure}
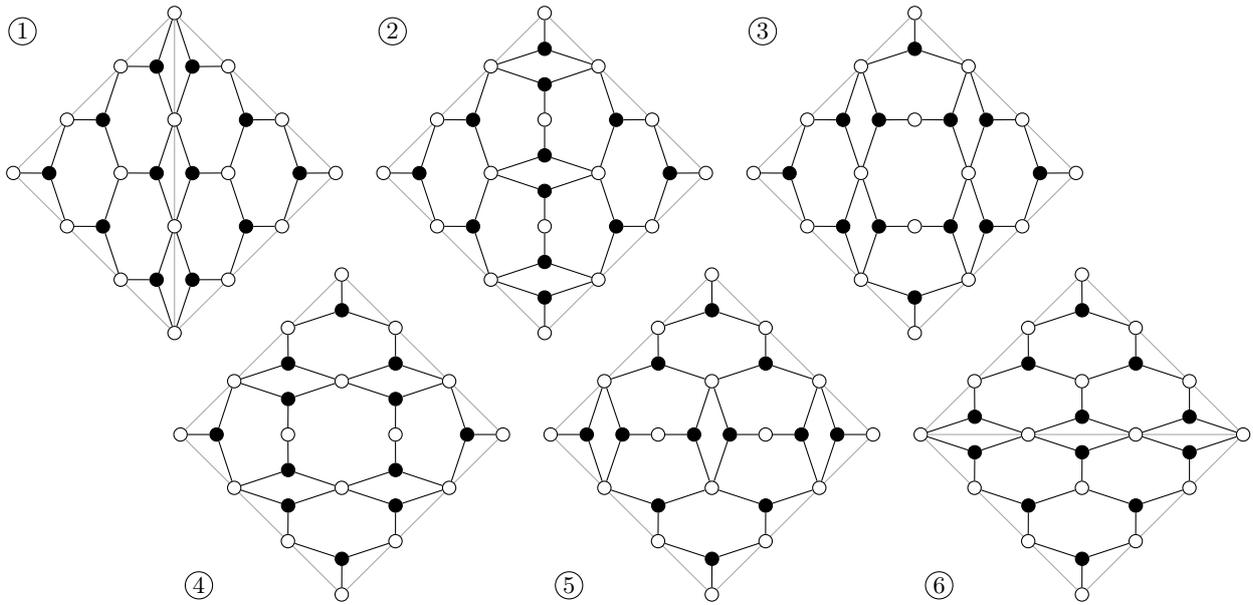

	\hspace{-22mm}
	
	\caption{The edge-flip in a projective 3-flag configuration as a sequence of fourteen 2-2 moves.}
	\label{fig:flagconfigflip}
\end{figure}

An edge-flip in the triangulation $\tg$ does not influence how we assign flags or decorations to projective $k$-flag configurations, but it does alter the values of the variables, and the way the variables are assigned to a given configuration. In fact, the sequence that corresponds to an edge-flip in terms of mutations of the projective quiver, can be found in \cite{fghighertm}, and the sequence in terms of 2-2 moves can be found in \cite{gwebs}. As an example, we show the edge-flip in a projective 3-flag configuration in Figure \ref{fig:flagconfigflip}.

In order to compare the variables introduced by Fock-Goncharov with our variables we need one more definition.

\begin{definition}
	Let $F$ be a projective $k$-flag configuration in $\CP^k$ on a triangulation $\tg$. The \emph{dual projective $k$-flag configuration} $F^\star$ in $(\CP^k)^*$ is obtained by projectively dualizing every $k$-flag $F(v)$ for each $v\in V(\tg)$.
\end{definition}

\begin{theorem}
	Let $F$ be a projective $k$-flag configuration on a triangulation $\tg$. Then the $\mathcal X$-cluster variables of $F$, as defined by Fock-Goncharov \cite{fghighertm}, are the $X$-variables of $T_{F^\star}$, and the $\mathcal X$-cluster variables of $F^\star$ are the $X$-variables of $T_{F}$.
\end{theorem}
\proof{
	Fock and Goncharov define the $X$-variables explicitly as cross-ratios and triple-ratios \cite[Section 9.4]{fghighertm}, that is multi-ratios of four or six points. Both coincide with the corresponding multi-ratios with which the $X$-variables can be expressed due to Lemma \ref{lem:projclusterviadistances}. \qed
}

\begin{conjecture}\label{con:flagtau}
	Consider a decorated projective $k$-flag configuration $F$ on a triangulation $\tg$. Let $T_F$ be the TCD map associated to $F$, and $H$ the hyperplanes of the decoration. Assign the homogeneous lift $H^\perp_i(v)$ to the $i$-th strand turning clockwise around $v$. Then the $\tau$-variables of $T_{F}$ are gauge equivalent to the $\mathcal A$-cluster variables of $F^\star$ as defined by Fock-Goncharov \cite{fghighertm}.
\end{conjecture}

We are very confident that Conjecture \ref{con:flagtau} holds. However, we were not able to find a direct and general definition of the $\mathcal A$-variables for projective $k$-flags in the literature. Therefore we do not give a proof of the conjecture.

It is natural to ask whether the $X$-variables determine a projective $k$-flag configuration $F$ up to projective transformations. This has been answered already by \cite{fghighertm}. Nevertheless, we give a short standard answer for triangulations without interior vertices in terms of TCD maps.

\begin{figure}
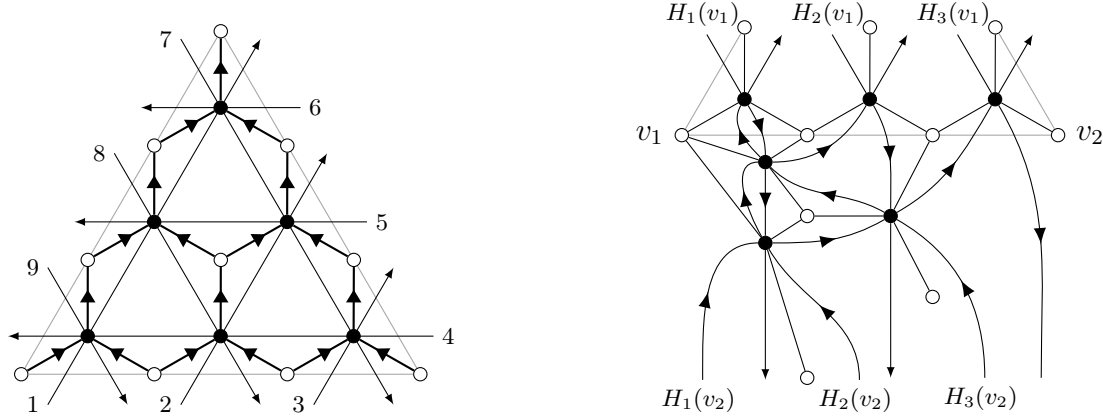


	\caption{Left: The li-orientation of a triangle in a projective 3-flag configuration, which propagates the geometric data via the $X$-variables. Right: Adding points on the boundary according to an upward sequence in a bridge-decomposition. The points determine the hyperplanes occurring in a decoration.}
	\label{fig:flagsstuff}
\end{figure}

Pick a boundary edge $e_0$. Then $e_0$ contains $k+1$ white vertices $w_1, w_2, \dots, w_{k+1}$ of $\pb_F$. On each line $T(w_i)T(w_{i+1})$, there is another marked point $T(w'_i)$ such that the white vertices $w_i,w_{i+1},w'_i$ are adjacent to a common black vertex of $\pb_F$. The $k$ marked points $w'_1,w'_2,\dots,w'_k$ span a hyperplane, which we denote by $J$. The choice of $k+1$ points and the hyperplane $J$ eliminates all the freedom given by projective transformations, because the only projective transformation that fixes $J$ and the $k+1$ points is the identity. Therefore we may choose the points $w_i$ and the hyperplane $J$ as initial data. We now claim that the $X$-variables determine all the points of $T_F$, and thus all the flags. This is particularly easy to see using the li-orientation (see Definition \ref{def:liorient}). In Figure \ref{fig:flagsstuff} we give an example for $k=3$.

In the case that the triangulation $\tg$ has internal vertices, we can look at a fundamental domain of the universal cover and propagate geometric data on this fundamental domain. As a result, certain $X$-variables ($2k$ per internal vertex) are not used for propagation. This leads to constraints on the $X$-variables. We claim that the constraints can be expressed via certain subvarieties given by comparatively simple equations for the $X$-variables, but we postpone a detailed discussion to feature research. Of course, one can also consider a triangulation $\tg$ with internal vertices and TCD maps that are only well defined up to monodromies around internal vertices, as in Teichmüller theory of punctured surfaces (that corresponds to $k=1$). In that case some, but not all, of the constraints on the $X$-variables are omitted, and we would be interested in a generalization of these constraints to $k>1$.

Consider again the case where $\tg$ has no internal vertices. Recall the discussion at the end of Section \ref{sec:tautcd}, where we outlined how to assign hyperplanes to strands for non-balanced TCDs $\tcd$. The idea is to consider $\tcd$ as a subdiagram of a balanced TCD $\tcd'$, so that $\tcd$ inherits the strand-hyperplanes from $\tcd'$. This procedure can be applied to the TCD maps of projective flag configurations as well. In particular, at every boundary edge of $\tg$, we can iteratively add white vertices to reverse the ordering of the ingoing strands. We show an example of this procedure in Figure \ref{fig:flagsstuff}. This added data at the boundary corresponds to a selection of the hyperplanes in the decoration.

\bibliographystyle{alpha}
\bibliography{references}

\newcommand{\etalchar}[1]{$^{#1}$}
\begin{thebibliography}{{K{\o}}92}

\bibitem[ABS03a]{absquads}
Vsevolod~E. Adler, Alexander~I. Bobenko, and Yuri~B. Suris.
\newblock Classification of integrable equations on quad-graphs. {T}he
  consistency approach.
\newblock {\em Comm. Math. Phys.}, 233(3):513--543, 2003.

\bibitem[ABS03b]{absyangbaxter}
Vsevolod~E. Adler, Alexander~I. Bobenko, and Yuri~B. Suris.
\newblock {Geometry of Yang-Baxter maps: Pencils of conics and quadrirational
  mappings}.
\newblock {\em Communications in Analysis and Geometry}, 12, 08 2003.

\bibitem[ABS09]{absgrassmannian}
Vsevolod~E. Adler, Alexander~I. Bobenko, and Yuri~B. Suris.
\newblock {Integrable Discrete Nets in Grassmannians}.
\newblock {\em Letters in Mathematical Physics}, 89:131--139, 2009.

\bibitem[ABS12]{absoctahedron}
Vsevolod~E. Adler, Alexander~I. Bobenko, and Yuri~B. Suris.
\newblock Classification of integrable discrete equations of octahedron type.
\newblock {\em Int. Math. Res. Not. IMRN}, pages 1822--1889, 2012.

\bibitem[Adl93]{adlerrecutting}
Vsevolod~E. Adler.
\newblock Recuttings of polygons.
\newblock {\em Functional Analysis and Its Applications}, 27(2):141--143, Apr
  1993.

\bibitem[Aff21]{amiquel}
Niklas~C. Affolter.
\newblock {Miquel dynamics, Clifford lattices and the Dimer model}.
\newblock {\em Letters in Mathematical Physics}, 111(3):61, May 2021.

\bibitem[AFIT20]{afitcrossratio}
Maxim Arnold, Dmitry Fuchs, Ivan Izmestiev, and Serge Tabachnikov.
\newblock {Cross-ratio Dynamics on Ideal Polygons}.
\newblock {\em International Mathematics Research Notices}, 12 2020.

\bibitem[AGPR19]{vrc}
Niklas~C. Affolter, Max Glick, Pavlo Pylyavskyy, and Sanjay Ramassamy.
\newblock Vector-relation configurations and plabic graphs, 2019.
\newblock Preprint, arXiv:1908.06959.

\bibitem[AGR21]{agrcrdyn}
Niklas~C. Affolter, Terrence George, and Sanjay Ramassamy.
\newblock Cross-ratio dynamics and the dimer cluster integrable system, 2021.
\newblock Preprint, arXiv:2108.12692.

\bibitem[AHBC{\etalchar{+}}14]{scattamp}
Nima Arkani-Hamed, Jacob~L. Bourjaily, Freddy Cachazo, Alexander~B. Goncharov,
  Alexander Postnikov, and Jaroslav Trnka.
\newblock {Scattering Amplitudes and the Positive Grassmannian}, 2014.
\newblock Preprint, arXiv:1212.5605.

\bibitem[AMdT22a]{amdtdskp}
Niklas~C. Affolter, Paul Melotti, and Béatrice de~Tilière.
\newblock {The Schwarzian octahedron recurrence (dSKP equation) I: explicit
  solutions}, 2022.
\newblock {Preprint, arXiv:2208.00239}.

\bibitem[AMdT22b]{amdtdskpb}
Niklas~C. Affolter, Paul Melotti, and Béatrice de~Tilière.
\newblock {The Schwarzian octahedron recurrence (dSKP equation) II: geometric
  systems}, 2022.
\newblock {Preprint, arXiv:2208.00244}.

\bibitem[Bax78]{baxterplanar}
Rodney~J. Baxter.
\newblock Solvable eight-vertex model on an arbitrary planar lattice.
\newblock {\em Philosophical Transactions of the Royal Society of London.
  Series A, Mathematical and Physical Sciences}, 289(1359):315--346, 1978.

\bibitem[BCdT20]{bctelliptic}
Cédric Boutillier, David Cimasoni, and Béatrice de~Tilière.
\newblock {Elliptic dimers on minimal graphs and genus 1 Harnack curves}, 2020.
\newblock Preprint, arXiv:2007.14699.

\bibitem[BdT11]{btcriticalisoradialising}
C{\'e}dric Boutillier and B{\'e}atrice de~Tili{\`e}re.
\newblock {The Critical Z-Invariant Ising Model via Dimers: Locality Property}.
\newblock {\em Communications in Mathematical Physics}, 301(2):473--516, Jan
  2011.

\bibitem[BdT12]{btisosurvey}
C{\'e}dric Boutillier and B{\'e}atrice de~Tili{\`e}re.
\newblock Statistical mechanics on isoradial graphs.
\newblock In {\em Probability in Complex Physical Systems}, pages 491--512,
  Berlin, Heidelberg, 2012. Springer Berlin Heidelberg.

\bibitem[BdTR19]{btrisoradialising}
C{\'e}dric Boutillier, B{\'e}atrice de~Tili{\`e}re, and Kilian Raschel.
\newblock {The Z-invariant Ising model via dimers}.
\newblock {\em Probability Theory and Related Fields}, 174(1):235--305, Jun
  2019.

\bibitem[BK95]{bkdarboux}
Leonid~V. Bogdanov and Boris~G. Konopelchenko.
\newblock {Lattice and q-difference Darboux-Zakharov-Manakov systems via delta
  -dressing method}.
\newblock {\em Journal of Physics A: Mathematical and General},
  28(5):L173--L178, mar 1995.

\bibitem[BK98a]{bkkphierarchy}
Leonid~V. Bogdanov and Boris~G. Konopelchenko.
\newblock {Analytic-bilinear approach to integrable hierarchies. I. Generalized
  KP hierarchy}.
\newblock {\em Journal of Mathematical Physics}, 39(9):4683--4700, 1998.

\bibitem[BK98b]{bkkphierarchymulti}
Leonid~V. Bogdanov and Boris~G. Konopelchenko.
\newblock {Analytic-bilinear approach to integrable hierarchies. II.
  Multicomponent KP and 2D Toda lattice hierarchies}.
\newblock {\em Journal of Mathematical Physics}, 39(9):4701--4728, 1998.

\bibitem[Bla40]{pascal}
Pascal Blaise.
\newblock Essay pour les coniques, 1640.

\bibitem[Bla10]{blaschkelaguerre}
Wilhelm Blaschke.
\newblock {Untersuchungen über die Geometrie der Speere in der Euklidischen
  Ebene}.
\newblock {\em Monatshefte f{\"u}r Mathematik und Physik}, 21(1):3--60, Dec
  1910.

\bibitem[BLPT21]{blptlaguerre}
Alexander~I. Bobenko, Carl~O.R. Lutz, H.~Pottmann, and Jan Techter.
\newblock {\em {Non-Euclidean Laguerre Geometry and Incircular Nets}}.
\newblock SpringerBriefs in Mathematics. Springer International Publishing,
  2021.

\bibitem[BMS05]{bmsconformal}
Alexander~I. {Bobenko}, Christian {Mercat}, and Yuri~B. {Suris}.
\newblock {Linear and nonlinear theories of discrete analytic functions.
  Integrable structure and isomonodromic Green's function}.
\newblock {\em Journal f\"ur die reine und angewandte Mathematik},
  2005:117--161, 2005.

\bibitem[BMS08]{bmscircular}
Vladimir~V. Bazhanov, Vladimir~V. Mangazeev, and Sergey~M. Sergeev.
\newblock Quantum geometry of three-dimensional lattices.
\newblock {\em Journal of Statistical Mechanics: Theory and Experiment},
  2008(07):P07004, 07 2008.

\bibitem[Bob99]{bobenkocircular}
Alexander~I. Bobenko.
\newblock {\em Discrete conformal maps and surfaces}, page 97–108.
\newblock London Mathematical Society Lecture Note Series. Cambridge University
  Press, 1999.

\bibitem[Boc16]{bocklandtdimerabc}
Raf Bocklandt.
\newblock {A dimer ABC}.
\newblock {\em Bulletin of the London Mathematical Society}, 48(3):387--451, 02
  2016.

\bibitem[BP96]{bpdisosurfaces}
Alexander~I. Bobenko and Ulrich Pinkall.
\newblock Discrete isothermic surfaces.
\newblock 1996(475):187--208, 1996.

\bibitem[Bro12]{broomheaddimer}
Nathan Broomhead.
\newblock {\em {Dimer Models and Calabi-Yau Algebras}}.
\newblock Memoirs of the American Mathematical Society. American Mathematical
  Society, 2012.

\bibitem[BS07a]{bsmoutard}
Alexander~I. Bobenko and Yuri~B. Suris.
\newblock {Discrete Koenigs Nets and Discrete Isothermic Surfaces}.
\newblock {\em International Mathematics Research Notices}, 09 2007.

\bibitem[BS07b]{bsorganizing}
Alexander~I. Bobenko and Yuri~B. Suris.
\newblock {On organizing principles of discrete differential geometry. Geometry
  of spheres}.
\newblock {\em Russian Mathematical Surveys}, 62(1):1--43, 02 2007.

\bibitem[BS08]{ddgbook}
Alexander~I. Bobenko and Yuri~B. Suris.
\newblock {\em Discrete Differential Geometry: Integrable Structure}, volume~98
  of {\em {Graduate Studies in Mathematics}}.
\newblock American Mathematical Society, 2008.

\bibitem[BS15]{bobenkoschieflinecomplexes}
Alexander~I. Bobenko and Wolfgang~K. Schief.
\newblock {Discrete line complexes and integrable evolution of minors}.
\newblock {\em Proceedings of the Royal Society A}, 471(2175), 03 2015.

\bibitem[BS16]{bobenkoschiefcirclecomplexes}
Alexander~I. Bobenko and Wolfgang~K. Schief.
\newblock {Circle Complexes and the Discrete CKP Equation}.
\newblock {\em International Mathematics Research Notices}, 2017(5):1504--1561,
  05 2016.

\bibitem[BSST16]{bsstconfocala}
Alexander~I. Bobenko, Wolfgang~K. Schief, Yuri~B. Suris, and Jan Techter.
\newblock {On a discretization of confocal quadrics. I. An integrable systems
  approach}.
\newblock {\em Journal of Integrable Systems}, 1(1), 08 2016.

\bibitem[BSST18]{bsstconfocalb}
Alexander~I. Bobenko, Wolfgang~K. Schief, Yuri~B. Suris, and Jan Techter.
\newblock {On a Discretization of Confocal Quadrics. A Geometric Approach to
  General Parametrizations}.
\newblock {\em International Mathematics Research Notices},
  2020(24):10180--10230, 12 2018.

\bibitem[BST20]{bstcheckerboardlaguerre}
Alexander~I. Bobenko, Wolfgang~K. Schief, and Jan Techter.
\newblock {Checkerboard incircular nets: Laguerre geometry and
  parametrisation}.
\newblock {\em Geometriae Dedicata}, 204(1):97--129, Feb 2020.

\bibitem[BW20]{bwtriple}
Alexey Balitskiy and Julian Wellman.
\newblock Flip cycles in plabic graphs.
\newblock {\em Selecta Mathematica}, 26(1):15, 2020.

\bibitem[BZ05]{bzquantumcluster}
Arkady Berenstein and Andrei Zelevinsky.
\newblock Quantum cluster algebras.
\newblock {\em Advances in Mathematics}, 195(2):405 -- 455, 2005.

\bibitem[Car06]{carnot}
Lazare N.~M. Carnot.
\newblock {Mémoire sur la relation qui existe entre les distances respectives
  de cinq points quelconques pris dans l'espace ; suivi d'un Essai sur la
  th\'eorie des transversals}.
\newblock {\em Paris}, 1806.

\bibitem[CDS97]{cdscircular}
Jan Cieśliński, Adam Doliwa, and Paolo~M. Santini.
\newblock The integrable discrete analogues of orthogonal coordinate systems
  are multi-dimensional circular lattices.
\newblock {\em Physics Letters A}, 235(5):480 -- 488, 1997.

\bibitem[Che18]{chelkaksembeddings}
Dmitry Chelkak.
\newblock Planar {I}sing model at criticality: state-of-the-art and
  perspectives.
\newblock In {\em Proceedings of the {I}nternational {C}ongress of
  {M}athematicians---{R}io de {J}aneiro 2018. {V}ol. {IV}. {I}nvited lectures},
  pages 2801--2828. World Sci. Publ., Hackensack, NJ, 2018.

\bibitem[Che20]{chelkaksgraphs}
Dmitry Chelkak.
\newblock Ising model and s-embeddings of planar graphs, 2020.
\newblock Preprint, arXiv:2006.14559.

\bibitem[CKP01]{ckpdimers}
Henry Cohn, Richard Kenyon, and James Propp.
\newblock A variational principle for domino tilings.
\newblock {\em J. Amer. Math. Soc.}, 14(2):297--346, 2001.

\bibitem[CLR20]{clrtembeddings}
Dmitry Chelkak, Benoît Laslier, and Marianna Russkikh.
\newblock Dimer model and holomorphic functions on t-embeddings of planar
  graphs, 2020.
\newblock Preprint, arXiv:2001.11871.

\bibitem[Cox03]{coxeterprojective}
Harold S.~M. Coxeter.
\newblock {\em Projective Geometry}.
\newblock Springer New York, second edition, 2003.

\bibitem[CR07]{crdimerspins}
David Cimasoni and Nicolai Reshetikhin.
\newblock Dimers on surface graphs and spin structures. i.
\newblock {\em Communications in Mathematical Physics}, 275:187--208, 10 2007.

\bibitem[CS04]{cscube}
Gabriel~D. Carroll and David~E Speyer.
\newblock The cube recurrence.
\newblock {\em Electron. J. Combin.}, 11(1):Research Paper 73, 31, 2004.

\bibitem[CS11]{csanalysisisoradial}
Dmitry Chelkak and Stanislav Smirnov.
\newblock Discrete complex analysis on isoradial graphs.
\newblock {\em Advances in Mathematics}, 228(3):1590 -- 1630, 2011.

\bibitem[CS12]{csuniversality}
Dmitry Chelkak and Stanislav Smirnov.
\newblock {Universality in the 2D Ising model and conformal invariance of
  fermionic observables}.
\newblock {\em Inventiones mathematicae}, 189(3):515--580, 2012.

\bibitem[CTS07]{clifford}
William~K. Clifford, Robert Tucker, and Henry J.~S. Smith.
\newblock {\em Mathematical Papers}.
\newblock AMS Chelsea Publishing Series. AMS Chelsea Publishing, 2007.

\bibitem[DCL19]{dlrandomcurrent}
Hugo Duminil-Copin and Marcin Lis.
\newblock On the double random current nesting field.
\newblock {\em Probability Theory and Related Fields}, 175(3):937--955, Dec
  2019.

\bibitem[DN91]{dndskp}
Irene~Ya. Dorfman and Frank~W. Nijhoff.
\newblock {On a (2+1)-dimensional version of the Krichever-Novikov equation}.
\newblock {\em Physics Letters A}, 157(2):107--112, 1991.

\bibitem[Dol97]{doliwalaplace}
Adam Doliwa.
\newblock {Geometric discretisation of the Toda system}.
\newblock {\em Physics Letters A}, 234(3):187 -- 192, 1997.

\bibitem[Dol99]{doliwaqnetsinquadrics}
Adam Doliwa.
\newblock Quadratic reductions of quadrilateral lattices.
\newblock {\em Journal of Geometry and Physics}, 30(2):169 -- 186, 1999.

\bibitem[Dol01]{doliwaanetspluecker}
Adam Doliwa.
\newblock {Discrete asymptotic nets and W-congruences in Plücker line
  geometry}.
\newblock {\em Journal of Geometry and Physics}, 39(1):9 -- 29, 2001.

\bibitem[Dol02]{doliwakoenigs}
Adam Doliwa.
\newblock {Geometric discretization of the Koenigs nets}.
\newblock {\em Journal of Mathematical Physics}, 44, 04 2002.

\bibitem[Dol07]{doliwatnets}
Adam Doliwa.
\newblock {The B-quadrilateral lattice, its transformations and the
  algebro-geometric construction}.
\newblock {\em J. Geom. Phys.}, 57:1171--1192, 2007.

\bibitem[Dol09]{doliwadesargues}
Adam Doliwa.
\newblock {Desargues maps and the Hirota–Miwa equation}.
\newblock {\em Proceedings of the Royal Society A: Mathematical, Physical and
  Engineering Sciences}, 466:1177 -- 1200, 2009.

\bibitem[Dol10a]{doliwadesarguesweyl}
Adam Doliwa.
\newblock {The affine Weyl group symmetry of Desargues maps and of the
  non-commutative Hirota-Miwa system}.
\newblock {\em Physics Letters A}, 375, 06 2010.

\bibitem[Dol10b]{doliwacqnet}
Adam Doliwa.
\newblock {The C-(symmetric) quadrilateral lattice, its transformations and the
  algebro-geometric construction}.
\newblock {\em Journal of Geometry and Physics}, 60(5):690 -- 707, 2010.

\bibitem[Dri92]{drinfeldyangbaxter}
Vladimir~G. Drinfeld.
\newblock On some unsolved problems in quantum group theory.
\newblock In {\em Quantum Groups}, pages 1--8, Berlin, Heidelberg, 1992.
  Springer Berlin Heidelberg.

\bibitem[DS97]{doliwasantiniqnet}
Adam Doliwa and Paolo~M. Santini.
\newblock Multidimensional quadrilateral lattices are integrable.
\newblock {\em Physics Letters A}, 233(4):365 -- 372, 1997.

\bibitem[DSM00]{dsmlinecongruence}
Adam Doliwa, Paolo~M. Santini, and Manuel Mañas.
\newblock Transformations of quadrilateral lattices.
\newblock {\em Journal of Mathematical Physics}, 41(2):944--990, 2000.

\bibitem[dT21]{tilierezdirac}
Béatrice de~Tilière.
\newblock {The Z-Dirac and massive Laplacian operators in the Z-invariant Ising
  model}.
\newblock {\em Electronic Journal of Probability}, 26(none):1 -- 86, 2021.

\bibitem[Dub11]{dubedattrick}
Julien Dubédat.
\newblock {Exact bosonization of the Ising model}, 2011.
\newblock Preprint, arXiv:1112.4399.

\bibitem[Fai22]{fairleythesis}
Alexander~Y. Fairley.
\newblock {In progress}, 2022.
\newblock PhD thesis.

\bibitem[Fel04]{felsnerbook}
Stefan Felsner.
\newblock {\em Geometric Graphs and Arrangements}.
\newblock Vieweg Verlag, 2004.

\bibitem[FG06]{fghighertm}
Vladimir~V. Fock and Alexander~B. Goncharov.
\newblock Moduli spaces of local systems and higher {T}eichm\"uller theory.
\newblock {\em Publications Math\'ematiques de l'IH\'ES}, 103:1--211, 2006.

\bibitem[FG07]{fgtwoflags}
Vladimir~V. Fock and Alexander~B. Goncharov.
\newblock Moduli spaces of convex projective structures on surfaces.
\newblock {\em Advances in Mathematics}, 208(1):249--273, 2007.

\bibitem[FM16]{fmloop}
Vladimir~V. Fock and Andrey Marshakov.
\newblock {\em Loop Groups, Clusters, Dimers and Integrable Systems}, pages
  1--65.
\newblock Springer International Publishing, Cham, 2016.

\bibitem[Foc15]{fockinverse}
Vladimir~V. Fock.
\newblock {Inverse spectral problem for GK integrable system}, 2015.
\newblock Preprint, arXiv:1503.00289.

\bibitem[FW01]{fwsignotopes}
Stefan Felsner and Helmut Weil.
\newblock Sweeps, arrangements and signotopes.
\newblock {\em Discrete Applied Mathematics}, 109(1):67--94, 2001.
\newblock 14th European Workshop on Computational Geometry.

\bibitem[FZ01]{fzclusteralgebra}
Sergey Fomin and Andrei Zelevinsky.
\newblock {Cluster algebras I: Foundations}.
\newblock {\em Journal of the American Mathematical Society}, 15, 05 2001.

\bibitem[FZ07]{fzcoefficients}
Sergey Fomin and Andrei Zelevinsky.
\newblock {Cluster algebras IV: Coefficients}.
\newblock {\em Compositio Mathematica}, 143(1):112–164, 2007.

\bibitem[Gal21]{galashincritical}
Pavel Galashin.
\newblock {Critical varieties in the Grassmannian}, 2021.
\newblock Preprint, arXiv:2102.13339.

\bibitem[GI19]{giclustermodular}
Terrence George and Giovanni Inchiostro.
\newblock Cluster modular groups of dimer models and networks, 2019.
\newblock Preprint, arXiv:1909.12896.

\bibitem[GK13]{gkdimers}
Alexander~B. Goncharov and Richard Kenyon.
\newblock Dimers and cluster integrable systems.
\newblock {\em Annales scientifiques de l'\'Ecole Normale Sup\'erieure},
  46(5):747--813, 2013.

\bibitem[Gli11]{glickpentagram}
Max Glick.
\newblock {The pentagram map and Y-patterns}.
\newblock {\em Advances in Mathematics}, 227(2):1019 -- 1045, 2011.

\bibitem[Gli15]{gdevron}
Max Glick.
\newblock {The Devron property}.
\newblock {\em Journal of Geometry and Physics}, 87:161--189, 2015.

\bibitem[Gon17]{gwebs}
Alexander~B. Goncharov.
\newblock {\em Ideal Webs, Moduli Spaces of Local Systems, and 3d Calabi--Yau
  Categories}, pages 31--97.
\newblock Springer International Publishing, Cham, 2017.

\bibitem[GP16]{gpymesh}
Max Glick and Pavlo Pylyavskyy.
\newblock {Y-meshes and generalized pentagram maps}.
\newblock {\em Proceedings of the London Mathematical Society},
  112(4):753--797, 03 2016.

\bibitem[GP20]{gpisingorthogonal}
Pavel Galashin and Pavlo Pylyavskyy.
\newblock {Ising model and the positive orthogonal Grassmannian}.
\newblock {\em Duke Math. J.}, 169(10):1877--1942, 2020.

\bibitem[GR17]{gdnotes}
Max Glick and Dylan Rupel.
\newblock {Introduction to Cluster Algebras}, 07 2017.
\newblock Notes for lecture at ASIDE conference 2016.

\bibitem[GR18]{grmiquel}
Alexey Glutsyuk and Sanjay Ramassamy.
\newblock {A first integrability result for Miquel dynamics}.
\newblock {\em Journal of Geometry and Physics}, 130:121 -- 129, 2018.

\bibitem[GSTV16]{gstvnetworks}
Michael Gekhtman, Michael Shapiro, Serge Tabachnikov, and Alek Vainshtein.
\newblock Integrable cluster dynamics of directed networks and pentagram maps.
\newblock {\em Advances in Mathematics}, 300:390--450, 2016.
\newblock Special volume honoring Andrei Zelevinsky.

\bibitem[GSV03]{gsvpoissonpaper}
Michael Gekhtman, Michael Shapiro, and Alek Vainshtein.
\newblock {Cluster Algebras and Poisson Geometry}.
\newblock {\em Moscow Mathematical Journal}, 3, 12 2003.

\bibitem[GSV10]{gsvpoissonbook}
Michael Gekhtman, Michael Shapiro, and Alek Vainshtein.
\newblock {\em {Cluster Algebras and Poisson Geometry}}.
\newblock Mathematical surveys and monographs. MPI, 2010.

\bibitem[Hir81]{hirotaequation}
Ryogo Hirota.
\newblock {Discrete Analogue of a Generalized Toda Equation}.
\newblock {\em Journal of the Physical Society of Japan}, 50(11):3785--3791,
  1981.

\bibitem[HS10]{hscube}
André Henriques and David~E Speyer.
\newblock The multidimensional cube recurrence.
\newblock {\em Adv. Math.}, 223(3):1107--1136, 2010.

\bibitem[Izo21a]{izosimovnetworks}
Anton Izosimov.
\newblock Dimers, networks, and cluster integrable systems, 2021.
\newblock Preprint, arXiv:2108.04975.

\bibitem[Izo21b]{izosimov}
Anton Izosimov.
\newblock {Intersecting the sides of a polygon}.
\newblock {\em Proc. Amer. Math. Soc.}, 2021.
\newblock To appear in.

\bibitem[Izo21c]{izosimovpoissonlie}
Anton Izosimov.
\newblock {Pentagram maps and refactorization in Poisson-Lie groups}, 2021.
\newblock Preprint, arXiv:1803.00726.

\bibitem[Izo22]{izosimovrecutting}
Anton Izosimov.
\newblock Polygon recutting as a cluster integrable system, 2022.
\newblock Preprint, arXiv:2201.12503.

\bibitem[Jos12]{josefssonorthodiagonal}
Martin Josefsson.
\newblock {Characterizations of Orthodiagonal Quadrilaterals}.
\newblock {\em Forum Geometricorum [electronic only]}, 12, 01 2012.

\bibitem[Kas61]{kasteleyn}
Pieter~W. Kasteleyn.
\newblock {The statistics of dimers on a lattice: I. The number of dimer
  arrangements on a quadratic lattice}.
\newblock {\em Physica}, 27(12):1209 -- 1225, 1961.

\bibitem[{Kas}96]{kashaev}
Rinat~M. {Kashaev}.
\newblock {On discrete three-dimensional equations associated with the local
  Yang-Baxter relation}.
\newblock {\em Letters in Mathematical Physics}, 38(4):389--397, Dec 1996.

\bibitem[Kel21]{kelsfcc}
Andrew~P. Kels.
\newblock Interaction-round-a-face and consistency-around-a-face-centered-cube.
\newblock {\em Journal of Mathematical Physics}, 62(3):033509, 2021.

\bibitem[Ken99]{kennelly}
Arthur~E. Kennelly.
\newblock The equivalence of triangles and three-pointed stars in conducting
  networks.
\newblock {\em Electrical world and engineer}, 34(12):413--414, 1899.

\bibitem[Ken02]{kenyonisoradial}
Richard Kenyon.
\newblock {The Laplacian and Dirac operators on critical planar graphs}.
\newblock {\em Inventiones mathematicae}, 150(2):409--439, 2002.

\bibitem[Ken03]{kenyondimerintro}
Richard Kenyon.
\newblock An introduction to the dimer model.
\newblock In {\em {Lecture notes from a minicourse given at the ICTP in May
  2002.}}, 2003.

\bibitem[Kir45]{kirchhoff}
Gustav~R. Kirchhoff.
\newblock {Ueber den Durchgang eines elektrischen Stromes durch eine Ebene,
  insbesondere durch eine kreisförmige}.
\newblock {\em Annalen der Physik}, 140(4):497--514, 1845.

\bibitem[Kir47]{kirchhofftrue}
Gustav~R. Kirchhoff.
\newblock {Ueber die Auflösung der Gleichungen, auf welche man bei der
  Untersuchung der linearen Vertheilung galvanischer Ströme geführt wird}.
\newblock {\em Annalen der Physik}, 148(12):497--508, 1847.

\bibitem[KLRR21]{kenyonlam}
Richard {Kenyon}, Wai~Yeung {Lam}, Sanjay {Ramassamy}, and Marianna {Russkikh}.
\newblock {Dimers and Circle patterns}.
\newblock {\em To appear in the Annales Scientifiques de l'ENS}, 10 2021.

\bibitem[{K{\o}}92]{koenigsanets}
Gabriel {K{\oe}nigs}.
\newblock {Sur les r\'eseaux plans \`a invariants \'egaux et les lignes
  asymptotiques.}
\newblock {\em {C. R. Acad. Sci., Paris}}, 114:55--57, 1892.

\bibitem[Koe36]{koebecp}
Paul Koebe.
\newblock {Kontaktprobleme der konformen Abbildung}.
\newblock {\em Ber. Sächs. Akad. Wiss. Leipzig, Math.-phys. Kl. 88}, pages
  141--164, 1936.

\bibitem[KP98]{kplelieuvre}
Boris~G. Konopelchenko and Ulrich Pinkall.
\newblock {Projective Generalizations of Lelieuvre's Formula}.
\newblock {\em Geometriae Dedicata}, 79:81--99, 1998.

\bibitem[KP16]{kenyonpemantle}
Richard Kenyon and Robin Pemantle.
\newblock {Double-dimers, the Ising model and the hexahedron recurrence}.
\newblock {\em Journal of Combinatorial Theory, Series A}, 137:27 -- 63, 2016.

\bibitem[KS01]{ksreciprocal}
Boris Konopelchenko and Wolfgang Schief.
\newblock {Reciprocal Figures, Graphical Statics, and Inversive Geometry of the
  Schwarzian BKP Hierarchy}.
\newblock {\em Studies in Applied Mathematics}, 109, 08 2001.

\bibitem[KS02]{ksclifford}
Boris~G. Konopelchenko and Wolfgang~K. Schief.
\newblock {Menelaus' theorem, Clifford configurations and inversive geometry of
  the Schwarzian {KP} hierarchy}.
\newblock {\em Journal of Physics A: Mathematical and General},
  35(29):6125--6144, 07 2002.

\bibitem[KS03]{kstetraoctacubo}
Alastair~D. King and Wolfgang~K. Schief.
\newblock Tetrahedra, octahedra and cubo-octahedra: integrable geometry of
  multi-ratios.
\newblock {\em Journal of Physics A: Mathematical and General}, 36(3):785--802,
  jan 2003.

\bibitem[KS04]{kenyonsheffield}
Richard Kenyon and Scott Sheffield.
\newblock Dimers, tilings and trees.
\newblock {\em Journal of Combinatorial Theory, Series B}, 92(2):295 -- 317,
  2004.
\newblock Special Issue Dedicated to Professor W.T. Tutte.

\bibitem[KS12]{ksconformaldesargues}
Alastair~D. King and Wolfgang~K. Schief.
\newblock {Clifford lattices and a conformal generalization of Desargues’
  theorem}.
\newblock {\em Journal of Geometry and Physics}, 62(5):1088 -- 1096, 2012.

\bibitem[KS14]{kscox}
Alastair~D. King and Wolfgang~K. Schief.
\newblock {Bianchi Hypercubes and a Geometric Unification of the Hirota and
  Miwa Equations}.
\newblock {\em International Mathematics Research Notices},
  2015(16):6842--6878, 09 2014.

\bibitem[Kup98]{kuperberg}
Greg Kuperberg.
\newblock {An exploration of the permanent-determinant method}.
\newblock {\em The electronic journal of combinatorics}, 5(1), 1998.

\bibitem[Lag85]{laguerre}
Edmond~N. Laguerre.
\newblock {\em {Recherches sur la Géométrie de Direction}}.
\newblock Gauthiers-Villars, 1885.

\bibitem[Lea19]{leafkashaev}
Alexander Leaf.
\newblock The {K}ashaev equation and related recurrences.
\newblock {\em SIGMA Symmetry Integrability Geom. Methods Appl.}, 15:Paper No.
  012, 64, 2019.

\bibitem[Len20]{lenzising}
Wilhelm Lenz.
\newblock {Beitrag zum Verständnis der magnetischen Erscheinungen in festen
  Körpern}.
\newblock {\em Z. Phys.}, 21:613--615, 1920.

\bibitem[LPW{\etalchar{+}}06]{pwconical}
Yang Liu, Helmut Pottmann, Johannes Wallner, Yongliang Yang, and Wenping Wang.
\newblock Geometric modeling with conical meshes and developable surfaces.
\newblock {\em ACM Transactions on Graphics}, 25(3):681--689, 7 2006.

\bibitem[Mel18]{melottikashaev}
Paul Melotti.
\newblock {The free-fermionic $C_2^{(1)}$ loop model, double dimers and
  Kashaev's recurrence}.
\newblock {\em Journal of Combinatorial Theory, Series A}, 158:407 -- 448,
  2018.

\bibitem[Miq38]{miquel}
Auguste Miquel.
\newblock {Th\'eor\`emes sur les intersections des cercles et des sph\`eres.}
\newblock {\em J. Math. Pures Appl.}, pages 517--522, 1838.

\bibitem[Miw82]{miwabkp}
Tetsuji Miwa.
\newblock On hirota's difference equations.
\newblock {\em Proc. Japan Acad. Ser. A Math. Sci.}, 58(1):9--12, 1982.

\bibitem[MRT20]{mrtalphaquads}
Paul Melotti, Sanjay Ramassamy, and Paul Thévenin.
\newblock {Cube moves for $s$-embeddings and $\alpha$-realizations}.
\newblock {\em To appear in Annales de l'Institut Henri Poincaré D,
  Combinatorics, Physics and their Interactions}, 2020.

\bibitem[MS16]{mstwist}
Greg Muller and David~E Speyer.
\newblock The twist for positroid varieties.
\newblock {\em Proceedings of the London Mathematical Society}, 115, 06 2016.

\bibitem[Mü15]{muellerconical}
Christian Müller.
\newblock Planar discrete isothermic nets of conical type.
\newblock {\em Beiträge zur Algebra und Geometrie / Contributions to Algebra
  and Geometry}, 57, 06 2015.

\bibitem[NCWQ84]{ncwqdskp}
Frank~W. Nijhoff, Hans~W. Capel, Gerlof~L. Wiersma, and G.~Reinout~W. Quispel.
\newblock Bäcklund transformations and three-dimensional lattice equations.
\newblock {\em Physics Letters A}, 105(6):267--272, 1984.

\bibitem[Ohm27]{ohm}
Georg~S. Ohm.
\newblock {\em Die galvanische Kette}.
\newblock T. H. Riemann, 1827.

\bibitem[OST10]{ostpentagram}
Valentin Ovsienko, Richard Schwartz, and Serge Tabachnikov.
\newblock {The Pentagram Map: A Discrete Integrable System}.
\newblock {\em Communications in Mathematical Physics}, 299(2):409--446, Oct
  2010.

\bibitem[Pen12]{pennerteichmuller}
Robert~C. Penner.
\newblock {\em Decorated {T}eichm\"{u}ller theory}.
\newblock QGM Master Class Series. European Mathematical Society (EMS),
  Z\"{u}rich, 2012.
\newblock With a foreword by Yuri I. Manin.

\bibitem[Pos06]{postgrass}
Alexander Postnikov.
\newblock {Total positivity, Grassmannians, and networks}, 2006.
\newblock Preprint, arXiv:0609764.

\bibitem[Pro03]{propp}
James Propp.
\newblock Generalized domino-shuffling.
\newblock {\em Theoretical Computer Science}, 303(2):267 -- 301, 2003.

\bibitem[PSBW21]{psbwamplituhedron}
Matteo Parisi, Melissa Sherman-Bennett, and Lauren Williams.
\newblock {The m=2 amplituhedron and the hypersimplex: signs, clusters,
  triangulations, Eulerian numbers}, 2021.
\newblock Preprint, arXiv:0609764.

\bibitem[PW01]{pwlinegeometry}
Helmut Pottmann and Johannes Wallner.
\newblock {\em Computational Line Geometry}.
\newblock Springer-Verlag, Berlin, Heidelberg, 2001.

\bibitem[Ram18]{ramassamymiquel}
Sanjay Ramassamy.
\newblock {Miquel Dynamics for Circle Patterns}.
\newblock {\em International Mathematics Research Notices}, Published online,
  2018.

\bibitem[RG11]{rgbook}
Jürgen Richter-Gebert.
\newblock {\em {Perspectives on Projective Geometry: A Guided Tour Through Real
  and Complex Geometry}}.
\newblock Springer Publishing Company, Incorporated, 1st edition, 2011.

\bibitem[RGST05]{rgsttropical}
Jürgen Richter-Gebert, Bernd Sturmfels, and Thorsten Theobald.
\newblock {\em First steps in tropical geometry}, page 289–317.
\newblock Idempotent Mathematics and Mathematical Physics, Proceedings Vienna
  2003. American Mathematical Society, 2005.

\bibitem[Riv94]{rivin}
Igor Rivin.
\newblock {Euclidean Structures on Simplicial Surfaces and Hyperbolic Volume}.
\newblock {\em Annals of Mathematics}, 139(3):553--580, 1994.

\bibitem[Sau33]{sauerqnet}
Robert Sauer.
\newblock {Weckelige Kurvennetze bei einer infinitesimalen
  Fl{\"a}chenverbiegung}.
\newblock {\em Mathematische Annalen}, 108(1):673--693, 1933.

\bibitem[Sau37]{saueranet}
Robert Sauer.
\newblock {\em {Projektive Liniengeometrie}}.
\newblock {G{\"o}schens Lehrb{\"u}cherei}. Walter de Gruyter, 1937.

\bibitem[Sch92]{schwartz}
Richard Schwartz.
\newblock The pentagram map.
\newblock {\em Experiment. Math.}, 1(1):71--81, 1992.

\bibitem[Sch97]{schramm}
Oded Schramm.
\newblock {Circle patterns with the combinatorics of the square grid}.
\newblock {\em Duke Mathematical Journal}, 86(2):347 -- 389, 1997.

\bibitem[Sch03]{schieflattice}
Wolfgang~K. Schief.
\newblock {Lattice Geometry of the Discrete Darboux, KP, BKP and CKP Equations.
  Menelaus’ and Carnot’s Theorems}.
\newblock {\em Journal of Nonlinear Mathematical Physics}, 10(sup2):194--208,
  2003.

\bibitem[Sch09]{schieftalk}
Wolfgang~K. Schief.
\newblock {Discrete Laplace-Darboux sequences, Menelaus' theorem and the
  pentagram map}.
\newblock {Talk at workshop: Geometric Aspects of Discrete and Ultra-discrete
  Integrable Systems}, 2009.

\bibitem[Sco06]{scottgrass}
Jeanne Scott.
\newblock Grassmannians and cluster algebras.
\newblock {\em Proceedings of the London Mathematical Society},
  92(2):345–380, 2006.

\bibitem[SK52]{semplekneebone}
John~G. Semple and Geoffrey~T. Kneebone.
\newblock {\em Algebraic projective geometry}.
\newblock Oxford, at the Clarendon Press, 1952.

\bibitem[Smi06]{smirnovtowards}
Stanislav Smirnov.
\newblock Towards conformal invariance of 2{D} lattice models.
\newblock In {\em International {C}ongress of {M}athematicians. {V}ol. {II}},
  pages 1421--1451. Eur. Math. Soc., Z\"{u}rich, 2006.

\bibitem[Spe07]{speyerdimers}
David~E Speyer.
\newblock Perfect matchings and the octahedron recurrence.
\newblock {\em Journal of Algebraic Combinatorics}, 25(3):309--348, 2007.

\bibitem[Ste05]{stephensoncp}
Kenneth Stephenson.
\newblock {\em Introduction to circle packing}.
\newblock Cambridge University Press, Cambridge, 2005.
\newblock The theory of discrete analytic functions.

\bibitem[Ste18]{steinmeierkoenigs}
Jannik Steinmeier.
\newblock {Two discretizations of Koenigs nets and their connection}, 2018.
\newblock Bachelor's thesis.

\bibitem[Tec21]{techterthesis}
Jan Techter.
\newblock Discrete confocal quadrics and checkerboard incircular nets, 2021.
\newblock PhD thesis.

\bibitem[Tem74]{temperleybijection}
Harold N.~V. Temperley.
\newblock {\em Enumeration of graphs on a large periodic lattice}, page
  155–160.
\newblock London Mathematical Society Lecture Note Series. Cambridge University
  Press, 1974.

\bibitem[TF61]{tfdimers}
Harold N.~V. Temperley and Michael~E. Fisher.
\newblock Dimer problem in statistical mechanics-an exact result.
\newblock {\em The Philosophical Magazine: A Journal of Theoretical
  Experimental and Applied Physics}, 6(68):1061--1063, 1961.

\bibitem[Thu79]{thurstonlecturenotes}
William~P. Thurston.
\newblock {The geometry and topology of three-manifolds}, 1979.
\newblock lecture notes.

\bibitem[Thu85]{thurstoncp}
William~P. Thurston.
\newblock {The finite Riemann mapping theorem}, 1985.
\newblock Invited talk at the International Symposium on the occasion of the
  proof of the Bierbach conjecture, Purdue University.

\bibitem[Thu17]{thurstontriple}
Dylan~P. Thurston.
\newblock From dominoes to hexagons.
\newblock In {\em Proceedings of the 2014 Maui and 2015 Qinhuangdao Conferences
  in Honour of Vaughan F.R. Jones’ 60th Birthday}, pages 399--414, Canberra
  AUS, 2017. Centre for Mathematics and its Applications, The Australian
  National University.

\bibitem[Tut63]{tutteembedding}
W.~T. Tutte.
\newblock {How to Draw a Graph}.
\newblock {\em Proceedings of the London Mathematical Society},
  s3-13(1):743--767, 1963.

\bibitem[Ves07]{veselovyangbaxter}
Alexander~P. Veselov.
\newblock Yang-baxter maps: Dynamical point of view.
\newblock In {\em Combinatorial aspect of integrable systems}, pages 145--167.
  Mathematical Society of Japan, 2007.

\bibitem[Wil12]{williamsintro}
Lauren~K. Williams.
\newblock Cluster algebras: An introduction.
\newblock {\em Bulletin of the American Mathematical Society}, 51, 12 2012.

\bibitem[Wil16]{williamsgrass}
Lauren~K. Williams.
\newblock A positive {G}rassmannian analogue of the permutohedron.
\newblock {\em Proc. Amer. Math. Soc.}, 144(6):2419--2436, 2016.

\bibitem[Zei97]{zeilberger}
Doron Zeilberger.
\newblock Dodgson's determinant-evaluation rule proved by two-timing men and
  women.
\newblock {\em Electron. J. Comb.}, 4, 1997.

\end{thebibliography}

\begin{figure}[bp]
	%
\end{figure}

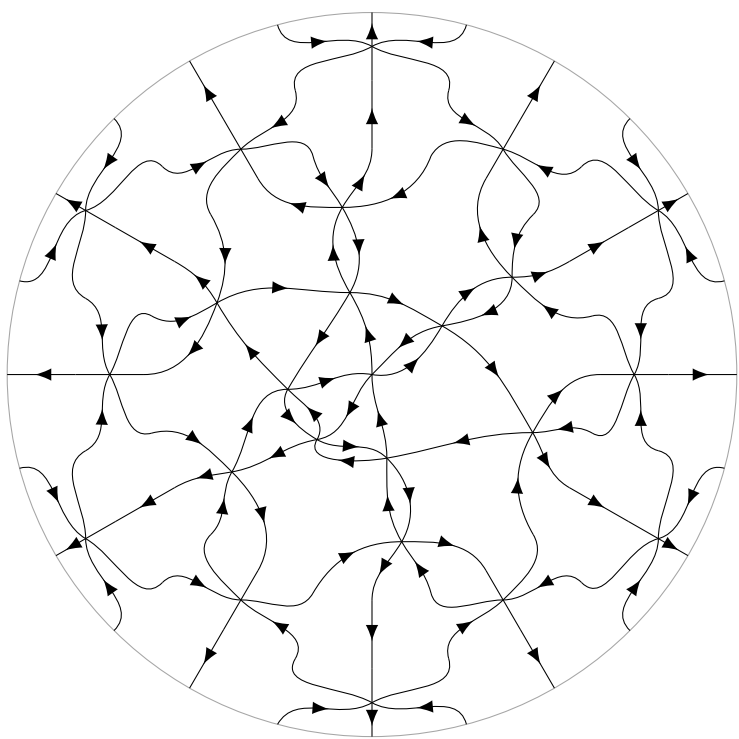
\end{document}